\documentclass[11pt,letterpaper]{article}
\RequirePackage[backend=biber,bibstyle=authoryear,citestyle=authoryear,
natbib=true,sorting=nyt,sortcites=true,maxnames=10,minnames=10,
hyperref=true,abbreviate=false,date=long,isbn=true,eprint=true,
uniquename=init,giveninits=true]{biblatex}
\usepackage{hyphenat,graphicx}
\usepackage[colorlinks,linktocpage,breaklinks]{hyperref}
\makeatletter
\def\@filecolor{blue}
\def\@linkcolor{blue}
\def\@citecolor{blue}
\def\@urlcolor{blue}
\let\@old@citep\citep
\let\@old@citet\citet
\let\@old@citeauthor\citeauthor
\def\citep{\@old@citep*}
\def\citet{\@old@citet*}
\def\citeauthor{\@old@citeauthor*}
\let\cite\citep
\makeatother
\usepackage[top=1.5in,bottom=1.2in,left=1.15in,right=1.15in,letterpaper]%
  {geometry}
\makeatletter
\def\ps@headings{%
    \let\@mkboth\@gobbletwo
    \def\@oddhead{\hss\scshape\shorttitle\hss\reset@font\rmfamily\thepage}
    \def\@evenhead{\reset@font\rmfamily\thepage\hss\scshape\shortauthors\hss}
    \let\@oddfoot\@empty\let\@evenfoot\@empty}
\@twosidetrue
\makeatother
\usepackage[neveradjust]{paralist}
\makeatletter
\@definecounter{@ADL@savenum}
\newcommand\savenum{\setcounter{@ADL@savenum}%
  {\the\@nameuse{c@\@listctr}}}
\newcommand\resumenum{\setcounter{\@listctr}{\arabic{@ADL@savenum}}}
\makeatother
\usepackage[reqno,centertags]{amsmath}
\usepackage{amssymb,xspace,twoopt,xy,mathrsfs}
\xyoption{arrow}
\xyoption{matrix}
\xyoption{curve}
\usepackage[mathcal]{eucal}
\makeatletter
\DeclareFontFamily{OMX}{MnSymbolE}{}
\DeclareSymbolFont{MnLargeSymbols}{OMX}{MnSymbolE}{m}{n}
\SetSymbolFont{MnLargeSymbols}{bold}{OMX}{MnSymbolE}{b}{n}
\DeclareFontShape{OMX}{MnSymbolE}{m}{n}{
    <-6>  MnSymbolE5
   <6-7>  MnSymbolE6
   <7-8>  MnSymbolE7
   <8-9>  MnSymbolE8
   <9-10> MnSymbolE9
  <10-12> MnSymbolE10
  <12->   MnSymbolE12
}{}
\DeclareFontShape{OMX}{MnSymbolE}{b}{n}{
    <-6>  MnSymbolE-Bold5
   <6-7>  MnSymbolE-Bold6
   <7-8>  MnSymbolE-Bold7
   <8-9>  MnSymbolE-Bold8
   <9-10> MnSymbolE-Bold9
  <10-12> MnSymbolE-Bold10
  <12->   MnSymbolE-Bold12
}{}
\let\llangle\@undefined
\let\rrangle\@undefined
\DeclareMathDelimiter{\llangle}{\mathopen}%
                     {MnLargeSymbols}{'164}{MnLargeSymbols}{'164}
\DeclareMathDelimiter{\rrangle}{\mathclose}%
                     {MnLargeSymbols}{'171}{MnLargeSymbols}{'171}
\DeclareRobustCommand\widecheck[1]{{\mathpalette\@widecheck{#1}}}
\def\@widecheck#1#2{%
    \setbox\z@\hbox{\m@th$#1#2$}%
    \setbox\tw@\hbox{\m@th$#1%
       \widehat{%
          \vrule\@width\z@\@height\ht\z@
          \vrule\@height\z@\@width\wd\z@}$}%
    \dp\tw@-\ht\z@
    \@tempdima\ht\z@ \advance\@tempdima2\ht\tw@ \divide\@tempdima\thr@@
    \setbox\tw@\hbox{%
       \raise\@tempdima\hbox{\scalebox{1}[-1]{\lower\@tempdima\box
\tw@}}}%
    {\ooalign{\box\tw@ \cr \box\z@}}}
\makeatother
\let\sinfty\infty
\input cyracc.def
\let\subsetneq\subset
\let\subset\subseteq

\newcommand\real{\mathbb{R}}
\newcommand\realp{\real_{>0}}

\newcommand\integer{\mathbb{Z}}
\newcommand\integerp{\integer_{>0}}
\newcommand\integernn{\integer_{\ge0}}

\newcommand\scirc{\raise1pt\hbox{$\,\scriptstyle\circ\,$}}
\newcommand\sscirc{\hbox{$\,\scriptscriptstyle\circ\,$}}
\newcommand\eqdef{\triangleq}

\newcommand\supscr[2]{#1^{\textup{#2}}}
\newcommand\ol[1]{\overline{#1}}
\newcommand\transpose[1]{#1^{T}}

\newcommand\disjointunion{\mathop{\overset{\sscirc}{\cup}}}
\newcommand\bigdisjointunion{\mathop{\overset{\sscirc}{\bigcup}}}
\newcommand\what[1]{\widehat{#1}\null}

\newdir{ >}{{}*!/-3\jot/\dir{>}}
\newcommand\slnorm{\lvert}
\newcommand\srnorm{\rvert}
\newcommand\snorm[1]{\slnorm #1\srnorm}
\newcommand\asnorm[1]{\left\slnorm #1\right\srnorm}
\newcommand\dlnorm{\lVert}
\newcommand\drnorm{\rVert}
\newcommand\dnorm[1]{\dlnorm #1\drnorm}
\newcommand\adnorm[1]{\left\dlnorm #1\right\drnorm}
\newcommand\setdef[2]{\{#1\;|\enspace#2\}}
\newcommand\asetdef[2]{\left\{#1\immediate\vphantom{#2}\;\right|
  \left.\immediate\vphantom{#1}\enspace#2\right\}}
\newcommand\ifam[1]{(#1)}
\newcommand\natpair[2]{\langle#1;#2\rangle}
\newcommand\ie{i.e.,}
\newcommand\eg{e.g.,}
\newcommand\cf{cf.}
\newcommand\etc{etc.}
\newcommand\resp{resp.}
\newcommand\mathupper[1]{\textup{#1}}
\renewcommand\d[1]{{\normalfont\textrm{d}}#1}
\newcommand\image{\operatorname{image}}
\newcommand\id{\operatorname{id}}
\newcommand\closure{\operatorname{cl}}
\newcommand\ver{\operatorname{ver}}
\newcommand\hor{\operatorname{hor}}
\newcommand\hol{\textup{hol}}
\newcommand\card{\mathupper{card}}
\newcommand\pr{\operatorname{pr}}
\newcommand\Sym{\operatorname{Sym}}
\newcommand\Ins{\mathupper{Ins}}
\newcommand\push{\mathupper{push}}
\newcommand\sB{\mathscr{B}}
\newcommand\sZ{\mathscr{Z}}
\let\algdual\stardual 
\let\dual\stardual 

\newcommand\deriv[2]{\frac{{\normalfont\mathupper{d}}#1}
  {{\normalfont\mathupper{d}}#2}}
\newcommand\derivatzero[2]{\frac{{\normalfont\mathupper{d}}#1}
  {{\normalfont\mathupper{d}}#2}\Big|_{#2=0}}

\newcommand\pderiv[2]{\frac{\partial#1}{\partial#2}}
\makeatletter
\newcommand\linder{\@ifnextchar[{\@ADL@rlinder}{\@ADL@linder}}
\def\@ADL@rlinder[#1]#2{\boldsymbol{D}^{#1}#2}
\newcommand\@ADL@linder[1]{\boldsymbol{D}#1}
\newcommand\plinder{\@ifnextchar[{\@ADL@rplinder}{\@ADL@plinder}}
\def\@ADL@rplinder[#1]#2#3{\boldsymbol{D}^{#1}_{#2}#3}
\newcommand\@ADL@plinder[2]{\boldsymbol{D}_{#1}#2}
\makeatother
\newcommand\map[3]{#1\colon#2\rightarrow#3}

\newcommand\mapdef[5]{\begin{aligned}
  #1\colon&\begin{aligned}[t]#2\end{aligned}\rightarrow
  \begin{aligned}[t]#3\end{aligned}\\&\begin{aligned}[t]#4\end{aligned}
  \mapsto\begin{aligned}[t]#5\end{aligned}\end{aligned}}
\newcommand\mapschar{\mathupper{C}}
\newcommand\C{\mapschar}
\renewcommand\c{\mathupper{c}}
\makeatletter
\newcommand\lin{\@ifnextchar[{\@ADL@klinmap}{\@ADL@linmap}}
\def\@ADL@klinmap[#1]#2#3{\def\@tempa{#3}%
  \ifx\@tempa\@empty\mathupper{L}^{#1}(#2;#2)\else%
  \mathupper{L}^{#1}(#2;#3)\fi}
\def\@ADL@linmap#1#2{\def\@tempa{#2}\ifx\@tempa\@empty\mathupper{L}(#1;#1)
  \else\mathupper{L}(#1;#2)\fi}
\makeatother
\newcommand\End{\mathupper{End}}
\newcommand\Hom{\mathupper{Hom}}
\newcommand\ann{\mathupper{ann}}
\newcommand\vect[1]{\boldsymbol{#1}}
\newcommand\alg[1]{\mathsf{#1}}
\newcommand\metric{\mathbb{G}}
\newcommand\hmetric{\mathbb{H}}
\newcommand\symmgroup[1]{\mathfrak{S}_{#1}}
\makeatletter
\newcommand\Symalg{\@ifstar{\@ADL@symalgs}{\@ADL@symalgns}}
\newcommand\@ADL@symalgns{\@ifnextchar[{\@ADL@symalg}{\@ADL@@symalg}}
\def\@ADL@symalg[#1]#2{\@ADL@tsymalgsym^{#1}(#2)}
\newcommand\@ADL@@symalg[1]{\@ADL@tsymalgsym(#1)}
\newcommand\@ADL@symalgs{\@ifnextchar[{\@ADL@ksymalgs}{\@ADL@noksymalgs}}
\def\@ADL@ksymalgs[#1]#2{\@ADL@symalgsym^{#1}(#2)}
\newcommand\@ADL@noksymalgs[1]{\@ADL@symalgsym(#1)}
\newcommand\Symalgsymbol[1]{\def\@ADL@symalgsym{#1}}
\Symalgsymbol{\mathupper{S}}
\newcommand\tensor{\@ifstar{\@ADL@tensors}{\@ADL@tensorns}}
\newcommand\@ADL@tensorns{\@ifnextchar[{\@ADL@rstensor}{\@ADL@tensor}}
\def\@ADL@rstensor[#1,#2]#3{\@ifnextchar[{\@ADL@rstensorbase{#1}{#2}{#3}}
  {\@ADL@rstensornobase{#1}{#2}{#3}}}
\def\@ADL@rstensorbase#1#2#3[#4]{\@ADL@tensorsym^{#1}_{#2}(#3)_{#4}}
\newcommand\@ADL@rstensornobase[3]{\@ADL@tensorsym^{#1}_{#2}(#3)}
\newcommand\@ADL@tensor[1]{\@ifnextchar[{\@ADL@tensorbase{#1}}
  {\@ADL@tensornobase{#1}}}
\def\@ADL@tensorbase#1[#2]{\@ADL@tensorsym(#1)_{#2}}
\newcommand\@ADL@tensornobase[1]{\@ADL@tensorsym(#1)}
\newcommand\@ADL@tensors{\@ifnextchar[{\@ADL@ktensors}{\@ADL@noktensors}}
\def\@ADL@ktensors[#1]#2{\@ADL@tensorsym^{#1}(#2)}
\newcommand\@ADL@noktensors[1]{\@ADL@tensorsym(#1)}
\newcommand\tensorsymbol[1]{\def\@ADL@tensorsym{#1}}
\tensorsymbol{\mathupper{T}}
\makeatother
\newcommand\ts[1]{\mathcal{#1}}
\newcommand\nbhd[1]{\mathcal{#1}}
\makeatletter
\newcommand\oball{\@ifnextchar[{\@ADL@oballarg}{\@ADL@oballnoarg}}
\def\@ADL@oballarg[#1]#2#3{\mathsf{B}_{#1}(#2,#3)}
\newcommand\@ADL@oballnoarg[2]{\mathsf{B}(#1,#2)}
\newcommand\cball{\@ifnextchar[{\@ADL@cballarg}{\@ADL@cballnoarg}}
\def\@ADL@cballarg[#1]#2#3{\ol{\mathsf{B}}_{#1}(#2,#3)}
\newcommand\@ADL@cballnoarg[2]{\ol{\mathsf{B}}(#1,#2)}
\makeatother
\newcommandtwoopt\sections[3][\sinfty][\null]{\Gamma^{#1}_{#2}(#3)}

\newcommand\man[1]{\mathsf{#1}}
\makeatletter
\newcommand\tb{\@ifnextchar[{\@ADL@tbarg}{\@ADL@tb}}
\def\@ADL@tbarg[#1]#2{\man{T}_{#1}#2}
\newcommand\@ADL@tb[1]{\man{T}#1}
\newcommand\tbproj[1]{\pi_{\tb{#1}}}
\newcommand\ctb{\@ifnextchar[{\@ADL@ctbarg}{\@ADL@ctb}}
\def\@ADL@ctbarg[#1]#2{\man{T}^*_{#1}#2}
\newcommand\@ADL@ctb[1]{\man{T}^*#1}
\newcommand\ctbproj[1]{\pi_{\ctb{#1}}}
\newcommand\nb{\@ifnextchar[{\@ADL@nbarg}{\@ADL@nb}}
\def\@ADL@nbarg[#1]#2{\man{N}_{#1}#2}
\newcommand\@ADL@nb[1]{\man{N}#1}
\newcommand\vb{\@ifnextchar[{\@ADL@vbarg}{\@ADL@vb}}
\def\@ADL@vbarg[#1]#2{\man{V}_{#1}#2}
\newcommand\@ADL@vb[1]{\man{V}#1}
\newcommand\cvb{\@ifnextchar[{\@ADL@cvbarg}{\@ADL@cvb}}
\def\@ADL@cvbarg[#1]#2{\man{V}^*_{#1}#2}
\newcommand\@ADL@cvb[1]{\man{V}^*#1}
\newcommand\hb{\@ifnextchar[{\@ADL@hbarg}{\@ADL@hb}}
\def\@ADL@hbarg[#1]#2{\man{H}_{#1}#2}
\newcommand\@ADL@hb[1]{\man{H}#1}
\newcommand\chb{\@ifnextchar[{\@ADL@chbarg}{\@ADL@chb}}
\def\@ADL@chbarg[#1]#2{\man{H}^*_{#1}#2}
\newcommand\@ADL@chb[1]{\man{H}^*#1}
\newcommand\func{\@ifnextchar[{\@ADL@crfuncs}{\@ADL@cinftyfuncs}}
\def\@ADL@crfuncs[#1]#2{\mapschar^{#1}(#2)}
\newcommand\@ADL@cinftyfuncs[1]{\mapschar^\infty(#1)}
\newcommand\sfunc{\@ifnextchar[{\@ADL@crsfuncs}{\@ADL@cinftysfuncs}}
\def\@ADL@crsfuncs[#1]#2{\mathscr{C}^{#1}_{#2}}
\newcommand\@ADL@cinftysfuncs[1]{\mathscr{C}^\sinfty_{#1}}
\newcommand\mappings{\@ifnextchar[{\@ADL@crmappings}{\@ADL@cinftymappings}}
\def\@ADL@crmappings[#1]#2#3{\mapschar^{#1}(#2;#3)}
\newcommand\@ADL@cinftymappings[2]{\mapschar^\infty(#1;#2)}
\newcommand\vbmappings{\@ifnextchar[{\@ADL@crvbmaps}{\@ADL@cinftyvbmaps}}
\def\@ADL@crvbmaps[#1]#2#3{\mathupper{VB}^{#1}(#2;#3)}
\newcommand\@ADL@cinftyvbmaps[2]{\mathupper{VB}^\sinfty(#1;#2)}
\newcommand\linfunc{\@ifnextchar[{\@ADL@crlinfuncs}{\@ADL@cinftylinfuncs}}
\def\@ADL@crlinfuncs[#1]#2{\mathupper{Lin}^{#1}(#2)}
\newcommand\@ADL@cinftylinfuncs[1]{\mathupper{Lin}^\sinfty(#1)}
\newcommand\gsections{\@ifnextchar[{\@ADL@rgsections}{\@ADL@gsections}}
\def\@ADL@rgsections[#1]#2#3{\mathscr{G}^{#1}_{#2,#3}}
\newcommand\@ADL@gsections[2]{\mathscr{G}^\infty_{#1,#2}}
\newcommand\ssections{\@ifnextchar[{\@ADL@rssections}{\@ADL@ssections}}
\def\@ADL@rssections[#1]#2{\mathscr{G}^{#1}_{#2}}
\newcommand\@ADL@ssections[1]{\mathscr{G}^\infty_{#1}}
\newcommand\tf{\@ifnextchar[{\@ADL@tfarg}{\@ADL@tf}}
\def\@ADL@tfarg[#1]#2{\@ifnextchar[{\@ADL@@tfargr{#1}{#2}}
  {\@ADL@@tfarg{#1}{#2}}}
\def\@ADL@@tfargr#1#2[#3]{T^{#3}_{#1}#2}
\newcommand\@ADL@@tfarg[2]{T_{#1}#2}
\newcommand\@ADL@tf[1]{\@ifnextchar[{\@ADL@@tfr{#1}}{\@ADL@@tf{#1}}}
\def\@ADL@@tfr#1[#2]{T^{#2}#1}
\newcommand\@ADL@@tf[1]{T#1}
\newcommand\ctf{\@ifnextchar[{\@ADL@ctfarg}{\@ADL@ctf}}
\def\@ADL@ctfarg[#1]#2{T^*_{#1}#2}
\newcommand\@ADL@ctf[1]{T^*#1}
\newcommand\lieder[2]{\def\@tempa{#2}\ifx\@tempa\@empty%
  \boldsymbol{\mathscr{L}}_{#1}\else\boldsymbol{\mathscr{L}}_{#1}#2\fi}
\newcommand\hl{\@ifstar{\@hl@star}{\@hl@nostar}}
\newcommand\@hl@star[1]{\supscr{#1}{h,$*$}}
\newcommand\@hl@nostar[1]{\supscr{#1}{h}}
\makeatother
\newcommand\ve[1]{\supscr{#1}{e}}
\newcommand\vl[1]{\supscr{#1}{v}}
\newcommand\tlift[1]{{#1}^T}
\newcommand\verlift{\operatorname{vlft}}
\newcommand\horlift{\operatorname{hlft}}
\makeatletter
\def\flow{\@ifstar{\@flow@star}{\@flow@nostar}}
\def\@flow@star#1{\Phi^{#1}}
\def\@flow@nostar#1#2{\@ifnextchar[{\@tflow{#1}{#2}}{\@flow{#1}{#2}}}
\def\@tflow#1#2[#3]{\Phi^{#1}_{#2,#3}}
\def\@flow#1#2{\Phi^{#1}_{#2}}
\def\jet{\@ifnextchar[{\@ADL@jetbase}{\@ADL@jetnobase}}
\def\@ADL@jetbase[#1]#2#3{\man{J}_{#1}^{#2}#3}
\def\@ADL@jetnobase#1#2{\man{J}^{#1}#2}
\def\jetalg{\@ifnextchar[{\@ADL@jetalgbase}{\@ADL@jetalgnobase}}
\def\@ADL@jetalgbase[#1]#2#3{\man{T}^{*#2}_{#1}#3}
\def\@ADL@jetalgnobase#1#2{\man{T}^{*#1}#2}
\def\Hjetalg{\@ifnextchar[{\@ADL@Hjetalgbase}{\@ADL@Hjetalgnobase}}
\def\@ADL@Hjetalgbase[#1]#2#3{\man{H}^{*#2}_{#1}#3}
\def\@ADL@Hjetalgnobase#1#2{\man{H}^{*#1}#2}
\def\Pjetalg{\@ifnextchar[{\@ADL@Pjetalgbase}{\@ADL@Pjetalgnobase}}
\def\@ADL@Pjetalgbase[#1]#2#3{\man{P}^{*#2}_{#1}#3}
\def\@ADL@Pjetalgnobase#1#2{\man{P}^{*#1}#2}
\def\Vjetalg{\@ifnextchar[{\@ADL@Vjetalgbase}{\@ADL@Vjetalgnobase}}
\def\@ADL@Vjetalgbase[#1]#2#3{\man{V}^{*#2}_{#1}#3}
\def\@ADL@Vjetalgnobase#1#2{\man{V}^{*#1}#2}
\def\Fjetalg{\@ifnextchar[{\@ADL@Fjetalgbase}{\@ADL@Fjetalgnobase}}
\def\@ADL@Fjetalgbase[#1]#2#3{\man{F}^{*#2}_{#1}#3}
\def\@ADL@Fjetalgnobase#1#2{\man{F}^{*#1}#2}
\def\Ljetalg{\@ifnextchar[{\@ADL@Ljetalgbase}{\@ADL@Ljetalgnobase}}
\def\@ADL@Ljetalgbase[#1]#2#3{\man{L}^{*#2}_{#1}#3}
\def\@ADL@Ljetalgnobase#1#2{\man{L}^{*#1}#2}
\def\Djetalg{\@ifnextchar[{\@ADL@Djetalgbase}{\@ADL@Djetalgnobase}}
\def\@ADL@Djetalgbase[#1]#2#3{\man{D}^{*#2}_{#1}#3}
\def\@ADL@Djetalgnobase#1#2{\man{D}^{*#1}#2}
\def\Cjetalg{\@ifnextchar[{\@ADL@Cjetalgbase}{\@ADL@Cjetalgnobase}}
\def\@ADL@Cjetalgbase[#1]#2#3{\man{C}^{*#2}_{#1}#3}
\def\@ADL@Cjetalgnobase#1#2{\man{C}^{*#1}#2}
\makeatother
\makeatletter
\def\interval{\@ifnextchar({\@ADL@openleftint}{\@ADL@closedleftint}}
\def\@ADL@openleftint(#1,#2{(#1,#2%
  \@ifnextchar){\@ADL@openrightint}{\@ADL@closedrightint}}
\def\@ADL@closedleftint[#1,#2{[#1,#2%
  \@ifnextchar){\@ADL@openrightint}{\@ADL@closedrightint}}
\def\@ADL@openrightint){)}
\def\@ADL@closedrightint]{]}
\makeatother
\usepackage{theorem}
\newtheorem{theorem}{Theorem}[section]
\newtheorem{lemma}[theorem]{Lemma}
\newtheorem{corollary}[theorem]{Corollary}
\newtheorem{prooflemma}{Lemma}[theorem]
\newtheorem{proofsublemma}{Sublemma}[theorem]

{\theorembodyfont{\normalfont\rmfamily}
\newtheorem{definition}[theorem]{Definition}
\newtheorem{remark}[theorem]{Remark}}
\makeatletter
\newcommand\pushright{\protect\@ADL@pushright}
\newcommand\@ADL@pushright[1]{{\ifvmode\null\hfill{#1}\par\else\ifmmode%
  \@ADLmaths@pushright{\hbox{#1}}\else\ifinner\@ADLhbox@pushright{#1}%
  \else\@ADLparag@pushright{#1}\fi\fi\fi}}
\newcommand\@ADLmaths@pushright[1]{{\ifinner\@ADLhbox@pushright{#1}\else%
  \tag*{$#1$}\fi}}
\newcommand\@ADLparag@pushright[1]{{\parfillskip=0pt\widowpenalty=10000%
  \displaywidowpenalty=10000\finalhyphendemerits=0\@ADLhbox@pushright#1\par}}
\newcommand\@ADLhbox@pushright{\unskip\nobreak\hfil\penalty50\hskip.2em%
  \null\hfill\hfill}
\newenvironment{proof}{\trivlist\item[\hskip\labelsep\textit{Proof:}\/]%
  \@ADLsave@set@qed\xspace\normalfont\rmfamily}
  {\qed\@ADLrestore@qed\endtrivlist}
\newif\if@ADL@qed\@ADL@qedfalse
\newcommand\qed{\protect\@ADL@qed{$\blacksquare$}}
\newcommand\@ADL@qed[1]{\if@ADL@qed\global\@ADL@qedfalse%
  \pushright{#1}\else\ifhmode\ifinner\else\par\fi\fi\fi}
\newcommand\@ADLrestore@qed{\global\let\if@ADL@qed\@ADLsaved@ifqed}
\newcommand\@ADLsave@set@qed{\let\@ADLsaved@ifqed
  \if@ADL@qed\global\@ADL@qedtrue}
\newcommand\eqqed{\tag*{\qed}}
\newenvironment{subproof}{\trivlist\item[\hskip\labelsep\textit{Proof:}\/]%
  \@ADLsave@set@subqed\normalfont\rmfamily}
  {\subqed\@ADLrestore@subqed\endtrivlist}
\newif\if@ADL@subqed\@ADL@subqedfalse
\newcommand\subqed{\protect\@ADL@subqed{$\blacktriangledown$}}
\newcommand\@ADL@subqed[1]{\if@ADL@subqed\global\@ADL@subqedfalse%
  \pushright{#1}\else\ifhmode\ifinner\else\par\fi\fi\fi}
\newcommand\@ADLrestore@subqed{\global\let\if@ADL@subqed\@ADLsaved@ifsubqed}
\newcommand\@ADLsave@set@subqed{\let\@ADLsaved@ifsubqed
  \if@ADL@subqed\global\@ADL@subqedtrue}

\newif\if@ADL@oprocend\@ADL@oprocendfalse
\newcommand\oprocend{\@ADLsave@set@oprocend
  \protect\@ADL@oprocend{$\bullet$}\@ADLrestore@oprocend}
\newcommand\@ADL@oprocend[1]{\if@ADL@oprocend\global\@ADL@oprocendfalse%
  \pushright{#1}\else\ifhmode\ifinner\else\par\fi\fi\fi}
\newcommand\@ADLrestore@oprocend{\global
  \let\if@ADL@oprocend\@ADLsaved@ifoprocend}
\newcommand\@ADLsave@set@oprocend{\let\@ADLsaved@ifoprocend\if@ADL@oprocend%
  \global\@ADL@oprocendtrue}
\newcommand\eqoprocend{\tag*{\oprocend}}
\makeatother
\newenvironment{keywords}{\quote\small\textbf{Keywords.}}{\endquote}
\newenvironment{AMS}{\quote\small\textbf{AMS Subject Classifications (2010).}}
   {\endquote}
\newcommand\defn[1]{{\normalfont\bfseries\emph{\mathversion{bold}#1}}}
\makeatletter
\@namedef{procnonum*}{\trivlist \@topsep \theorempreskipamount
  \@topsepadd \theorempostskipamount
  \@ifnextchar[{\@ADL@xprocnonumstar}{\@ADL@yprocnonumstar}}
\@namedef{endprocnonum*}{\endtrivlist}
\def\@ADL@xprocnonumstar[#1]{\item[\hskip \labelsep{\theorem@headerfont #1}]
  \normalfont\rmfamily}
\def\@ADL@yprocnonumstar{\item[] \normalfont\rmfamily}
\makeatother
\numberwithin{equation}{section}
\newcommand\pldblref[2]{\mbox{\ref{#1}(\ref{#2})}}
\title{Geometric analysis on real analytic manifolds\thanks{Research
supported in part by a grant from the Natural Sciences and Engineering
Research Council of Canada}}
\author{Andrew D.\ Lewis\thanks{Professor, Department of Mathematics and
Statistics, Queen's University, Kingston, ON K7L 3N6, Canada,
email:~\texttt{andrew.lewis@queensu.ca}}}
\date{2022/02/14}
\newcommand\shorttitle{Geometric analysis on real analytic manifolds}
\newcommand\shortauthors{A.\ D.\ Lewis}
\begin{document}
\maketitle

\begin{abstract}
The continuity, in a suitable topology, of algebraic and geometric operations
on real analytic manifolds and vector bundles is proved.  This is carried out
using recently arrived at seminorms for the real analytic topology.  A new
characterisation of the topology of the space of real analytic mappings
between manifolds is also developed.  To characterise these topologies,
geometric decompositions of various jet bundles are given by use of
connections.  These decompositions are then used to characterise many of the
standard operations from differential geometry: algebraic operations, tensor
evaluation, various lifts of tensor fields, compositions of mappings,~\etc
Apart from the main results, numerous techniques are developed that will
facilitate the performing of analysis on real analytic manifolds.
\end{abstract}

\tableofcontents

\section{Introduction}

Almost all operations/operators in differential geometry are formed by
combining a few essential operations such as composition, prolongation,
tensor evaluation, and/or some sort of lifting process.  Typically, these
operations are tacitly regarded as being continuous in some sense.  In smooth
differential geometry, ``in some sense'' usually means with respect to the
smooth compact-open topology (the topology of uniform convergence of
derivatives on compact sets) or, if one is interested in differential
topology, a topology like the Whitney topology.  An accounting of these sorts
of topologies can be found in~\cite{MWH:76,PWM:80}\@.  It is pretty easy to
convince oneself of the continuity of standard operations in the smooth
compact-open topology; one only needs to put suitable bounds on finitely many
derivatives on compact sets.  Thus there is some justification in not working
this out carefully in the smooth case.  However, if one is interested in
continuity in the real analytic category, it is not very easy to convince
oneself about the continuity of geometric operations.  Indeed, the more one
thinks about this, the harder the problem becomes.

A barrier right at the start is that the appropriate topology for real
analytic functions (functions, for simplicity) is not so easily envisaged.
While a suitable real analytic topology has been around since at least the
work of \citet{AM:66}\@\textemdash{}who provided two descriptions of such a
topology, and showed that they agree\textemdash{}there has not been a
``user-friendly'' description of the real analytic topology,~\ie~a
description using seminorms, until quite recently.  Some useful initial
formulae are provided by \citet{JM:84}\@, and seminorms are provided in the
lecture notes of \cite{PD:10}\@.  However, as far as we are aware, it is only
in the technical note of \citet{DV:13} that we see a proof of the suitability
of these seminorms.  These were adapted to the geometric setting for sections
of a real analytic vector bundle by \citet{SJ/ADL:14a}\@.  Part of this
development was a decomposition of jet bundles using connections.  The
initial developments of that monograph are the starting point for our
approach here.

Another complicating facet of the real analytic theory arises when one
considers lifts from the base space to the total space of a real analytic
vector bundle $\map{\pi_{\man{E}}}{\man{E}}{\man{M}}$\@,~\eg~vertical lift of
a section of $\man{E}$ or horizontal lift of a vector field on $\man{M}$\@.
The first of these operations requires no additional structure, but the
second requires a connection.  However, both require connections to study
their real analytic continuity, because one needs to provide bounds for the
jets on the codomain (\ie~on $\man{E}$) in terms of jets on the domain
(\ie~on $\man{M}$).  To provide seminorms, one also needs Riemannian and
vector bundle metrics, and all of this data has to fit together nicely to
provide the bounds required.  For instance, one has a natural Riemannian
metric on the manifold $\man{E}$ arising from~(1)~a Riemannian metric on
$\man{M}$\@, (2)~a fibre metric on $\man{E}$\@, (3)~an affine connection on
$\man{M}$\@, and (4)~a linear connection in $\man{E}$\@.  This structure
makes use of the resulting structure of
$\map{\pi_{\man{E}}}{\man{E}}{\man{M}}$ being a Riemannian submersion, and
using formulae of \citet{BO:68}\@.  The determination of a systematic means
to provide jet bundle estimates in this setting occupies us for a significant
portion of the paper.

In order to illustrate the nature of the difficulties one encounters, let us
consider a specific and illustrative instance of the sort of argument that
one must piece together to prove continuity in the real analytic case.
Suppose that we have a real analytic vector bundle
$\map{\pi_{\man{E}}}{\man{E}}{\man{M}}$ with $\nabla^{\pi_{\man{E}}}$ a real
analytic linear connection in $\man{E}$\@.  Let $X$ be a real analytic vector
field on $\man{M}$ which we horizontally lift to a real analytic vector field
$\hl{X}$ on $\man{E}$\@.  To assess the continuity of the map
$X\mapsto\hl{X}$ in the real analytic topology, one needs to compute jets of
$\hl{X}$ and relate these to jets of $X$\@.  Thus one needs to differentiate
$\hl{X}$ arbitrarily many times.  This differentiation must be done on
$\man{E}$\@, as this is the base on which $\hl{X}$ is defined.  Trying this
directly in local coordinates is, in principle, possible, but it is pretty
unlikely that one will be able to produce the refined estimates required in
this way.  Thus, in our approach, one needs an affine connection on $\man{E}$
(thinking of $\man{E}$ as just a manifold now).  One can now see that there
will be a complicated intermingling of the linear connection
$\nabla^{\pi_{\man{E}}}$\@, an affine connection $\nabla^{\man{M}}$ on
$\man{M}$ (to compute jets of $X$), and a fabricated affine connection
$\nabla^{\man{E}}$ on $\man{E}$\@.  This is only the beginning of the
difficulties one faces.  One also needs, not only formulae for the
derivatives of $\hl{X}$\@, but also recursive formulae relating how a
derivative of $\hl{X}$ of order, say, $k$ is related to the derivatives of
$X$ of orders $0,1,\dots,k$\@.  This recursive formulae is essential for
being able to obtain growth estimates for the derivatives needed to relate
the seminorms applied to $\hl{X}$ to those applied to $X$\@.  Moreover, since
the mapping $X\mapsto\hl{X}$ is injective, one might hope that the mapping is
not just continuous, but is an homeomorphism onto its image.  To prove this,
one now needs to get estimates for the jets of $X$ from formulae involving
the jets of $\hl{X}$\@.  Thus one needs estimates that go ``both ways.''  It
is also worth mentioning that the estimates one needs from these recursive
formulae are quite unforgiving, and so their form has to be very precisely
managed.  This requires extensive bookkeeping.  This bookkeeping occupies us
for a substantial portion of the paper.  This is contrasted with the smooth
case, where very coarse bounds suffice; we shall say a few words about this
contrast at illustrative places in the paper.

Another difficulty is that the use of connections to compute derivatives for
jets forces one to address the matter of whether the seminorms used for jets,
and derived from the use of connections, are actually not dependent on the
chosen connection.  Thus one must compare iterated covariant derivatives with
respect to different connections and show that these are related to one
another in such a way that the resulting real analytic topology is well
defined.  This, in itself, is a substantial undertaking.  It is done in an
\emph{ad hoc} way by \citet[Lemma~2.5]{SJ/ADL:14a}\@; here we do this in a
systematic and geometric way that offers many benefits towards the objectives
of this paper, apart from rendering more attractive the computations of
\citeauthor{SJ/ADL:14a}\@.

We mention that the idea of obtaining recursive formulae for derivatives is
given in a local setting by \citet{VT:97} during the course of the proof of
his Proposition~2.5, and can be applied to the mapping
$\func[\omega]{\man{N}}\ni f\mapsto\Phi^*f\in\func[\omega]{\man{M}}$ of
pull-back by a real analytic mapping
$\Phi\in\mappings[\omega]{\man{M}}{\man{N}}$\@.  We are able to extend the
ideas in \citeauthor{VT:97}' computations to general classes of geometric
operations.  For example, as we mention above, a local working out of the
estimates for the horizontal lift operation seems like it will be very
difficult.  However, once one \emph{does} get these things to work out, it is
relatively straightforward to prove the main results of the paper, which are
the continuity of the fundamental geometric operations mentioned in the first
paragraph.

One of the features of the paper is that almost all constructions are done
intrinsically.  While this may seem to unnecessarily complicate things, this
is not, in fact, so.  Even were one to work locally, there would still arise
two difficult problems that we overcome in our approach, but that still must
be overcome in a local approach: (1)~the difficulty of lifts as described in
detail above; (2)~the verification that the topologies do not depend on
various choices made (charts in the local calculations, and metrics and
connections in the intrinsic calculations).  Thus, while the intrinsic
calculations are sometimes complicated, they are only a little more
complicated than the necessarily already complicated local calculations.  And
we believe that the intrinsic approach is ultimately easier to use, once one
understands how to use it.  An objective of this paper is to do a lot of the
tedious hard work required to produce methods and results that are themselves
more or less straightforward.

As a side-benefit to our approach, we also are able to easily provide proofs
in the finitely differentiable and smooth cases.  We point out the relevant
places where modifications can be made to the real analytic proofs to give
the results in the finitely differentiable and smooth cases.

\subsection{Organisation of paper}

In Section~\ref{sec:ra-topology} we review the definition of the real
analytic topology and the geometric seminorms for this topology as
constructed in~\cite{SJ/ADL:14a}\@.

Section~\ref{sec:tensor-constructions} is the first of three sections,
forming the bulk of the paper in terms of words used, where we provide a host
of geometric constructions whose bearing on the main goal of the paper will
be difficult to glean on a first reading.  Some sketchy motivation for the
constructions of
Sections~\ref{sec:tensor-constructions}\@,~\ref{sec:tensor-derivatives}\@,
and~\ref{sec:lift-isomorphisms} is outlined above in our discussion of the
difficulties one will encounter trying to prove continuity of the horizontal
lift mapping $X\mapsto\hl{X}$\@.  In Section~\ref{sec:tensor-constructions}
we perform constructions with functions, vector fields, and tensors on the
total space of a vector bundle.  These form the basis for derivative
computations done in Section~\ref{sec:tensor-derivatives}\@.  Particularly,
in Section~\ref{subsec:Esubmersion} we give
$\map{\pi_{\man{E}}}{\man{E}}{\man{M}}$ the structure of a Riemannian
submersion, following \citet{BO:68}\@.  This allows us to relate, in a
natural way, constructions on $\man{E}$ with those on $\man{M}$\@.  In
Section~\ref{sec:lift-isomorphisms} we provide the crucial recursive formulae
that relate derivatives on $\man{E}$ with those on $\man{M}$\@.  We do this
for a few of the standard geometric lifts one has for a vector bundle with a
linear connection.  Some of these we do because they are intrinsically
interesting.  Some we do because they are required for our general approach,
even if one is not interested in them \emph{per se}\@.

In Section~\ref{sec:fibre-norms} we give fibre norms for various jet bundles
that are used to define seminorms corresponding to the geometric
constructions of interest.  In Section~\ref{sec:jet-estimates} we put all of
our work from
Sections~\ref{sec:tensor-constructions}--\ref{sec:fibre-norms} to use to
prove Lemma~\ref{lem:pissynabla}\@, the technical lemma which makes
everything work.  The lemma gives a very precise estimate for the fibre norms
of derivatives of coefficients that arise in the recursive constructions of
Section~\ref{sec:lift-isomorphisms}\@.  There is no wiggle room in the form
of the required estimate, and this is one of the reasons why the computations
of Sections~\ref{sec:tensor-constructions}--\ref{sec:lift-isomorphisms} are
so laboriously carried out; these computations need to be understood at a
high resolution.  Once we have these estimates, however, in
Section~\ref{sec:independent} we show that the fibre-norms for jet bundles
obtained in Section~\ref{sec:fibre-norms} behave in the proper way as to make
the topologies we construct independent of our choices of connections and
metrics.  This is stated as Lemma~\ref{lem:G1G2comp}\@.  The actual proving
of the independence of the topologies is carried out by proving in
Theorem~\ref{the:same-topology} that the topologies are each the same as a
topology described using local forms of the seminorms.  This device of using
a local description carries two benefits.
\begin{compactenum}
\item It provide the local description of the seminorms.  While our approach
is intrinsic as much as this is possible, sometimes in practice one must work
locally, and having the explicit local formulae for the fibre norms is
beneficial.
\item While we have tried to make our treatment intrinsic, there is a crucial
point where a \emph{local} estimate for the growth of derivatives becomes
unavoidable, resting as it does on the Cauchy estimates for holomorphic
functions.  In our proof of Theorem~\ref{the:same-topology} is where this
seemingly unavoidable local estimate is not avoided.
\end{compactenum}

In Section~\ref{sec:continuity} we prove continuity of some representative
and some important geometric constructions.  There is a long list of these
constructions and we only give representatives; we hope that the tools we
develop in the paper, and put to use in Section~\ref{sec:continuity}\@, will
make it easy for researchers down the road to prove some important results in
the real analytic setting where continuity is crucial.

We close the paper, in Section~\ref{sec:mappings}\@, by using the tools from
the previous parts of the paper to define a topology for the space of real
analytic mappings between manifolds.  This topology can be seen as an
adaptation to the real analytic setting of the topology\textemdash{}for
finitely differentiable and smooth mappings\textemdash{}of uniform
convergence on compacta of mappings and their derivatives.  We introduce for
this topology a class of semimetrics that play the r\^ole of the seminorms
used in earlier sections of the paper for space of sections of a vector
bundle.  As an important application, we prove the continuity of the
nonlinear operator
\begin{equation*}
\mappings[\omega]{\man{M}}{\man{N}}\ni\Phi\mapsto f\scirc\Phi
\in\func[\omega]{\man{M}}
\end{equation*}
for a fixed $f\in\func[\omega]{\man{N}}$\@,~\ie~of the so-called
``superposition operator.''

\subsection{Notation and background}\label{subsec:background}

We shall quickly review the notation we use.

\subsubsection*{Basic terminology and notation}

When $A$ is a subset of a set $X$\@, we write $A\subset X$\@.  If we wish to
exclude the possibility that $A=X$\@, we write $A\subsetneq X$\@.  For a
family of sets $\ifam{X_i}_{i\in I}$\@, we denote by $\prod_{i\in I}X_i$ the
product of these sets.  By $\map{\pr_j}{\prod_{i\in I}X_i}{X_j}$ we denote
the projection onto the $j$th factor.  The identity map on a set $X$ is
denoted by $\id_X$\@.

By $\integer$ we denote the set of integers.  We use the notation $\integerp$
and $\integernn$ to denote the subsets of positive and nonnegative integers.
By $\real$ we denote the sets of real numbers.  By $\realp$ we denote the
subset of positive real numbers.

\subsubsection*{Algebra and linear algebra}

By $\symmgroup{k}$ we denote the permutation group of $\{1,\dots,k\}$\@.  For
$k,l\in\integernn$\@, we denote by $\symmgroup{k,l}$ the subset of
$\symmgroup{k+l}$ consisting of permutations $\sigma$ satisfying
\begin{equation*}
\sigma(1)<\dots<\sigma(k),\quad\sigma(k+1)<\dots<\sigma(k+l).
\end{equation*}
We also denote by $\symmgroup{k|l}$ the subgroup of $\symmgroup{k+l}$ having
the form
\begin{equation*}
\begin{pmatrix}1&\cdots&k&k+1&\cdots&k+l\\
\sigma_1(1)&\cdots&\sigma_1(k)&k+\sigma_2(1)&\cdots&k+\sigma_2(l)\end{pmatrix}
\end{equation*}
for $\sigma_1\in\symmgroup{k}$ and $\sigma_2\in\symmgroup{l}$\@.  We note
that $\symmgroup{k|l}\backslash\symmgroup{k+l}\simeq\symmgroup{k,l}$\@, so
that (1)~if $\sigma\in\symmgroup{k+l}$\@, then
$\sigma=\sigma_1\scirc\sigma_2$ for $\sigma_1\in\symmgroup{k|l}$ and
$\sigma_2\in\symmgroup{k,l}$ and
(2)~$\card(\symmgroup{k,l})=\frac{(k+l)!}{k!l!}$\@.

We denote by $\real^n$ the $n$-fold Cartesian product of $\real$\@.  A point
in $\real^n$ will typically be denoted in a bold
font,~\eg~$\vect{x}=(x_1,\dots,x_n)$\@.  We denote the standard basis for
$\real^n$ by $\ifam{\vect{e}_1,\dots,\vect{e}_n}$\@.

For $\real$-vector spaces $\alg{U}$ and $\alg{V}$\@, we denote by
$\Hom_{\real}(\alg{U};\alg{V})$ the set of $\real$-linear mappings from
$\alg{U}$ to $\alg{V}$\@.  We denote
$\End_{\real}(\alg{V})=\Hom_{\real}(\alg{V};\alg{V})$\@.  We denote by
$\dual{\alg{V}}=\Hom_{\real}(\alg{V};\real)$ the algebraic dual.  If
$v\in\alg{V}$ and $\alpha\in\dual{\alg{V}}$\@, we will denote the evaluation
of $\alpha$ on $v$ at various points by $\alpha(v)$\@, $\alpha\cdot v$\@, or
$\natpair{\alpha}{v}$\@, whichever seems most pleasing to us at the moment.
If $A\in\Hom_{\real}(\alg{U};\alg{V})$\@, we denote by
$\dual{A}\in\Hom_{\real}(\dual{\alg{V}};\dual{\alg{U}})$ the dual of $A$\@.
If $S\subset\alg{V}$\@, then we denote by
\begin{equation*}
\ann(S)=\setdef{\alpha\in\dual{\alg{V}}}{\alpha(v)=0,\ v\in S}
\end{equation*}
the annihilator subspace.

For a $\real$-vector space $\alg{V}$\@, $\tensor*[k]{\alg{V}}$ is the
$k$-fold tensor product of $\alg{V}$ with itself.  For $r,s\in\integerp$\@,
we denote
\begin{equation*}
\tensor[r,s]{\alg{V}}=\underbrace{\alg{V}\otimes\dots\otimes
\alg{V}}_{r\ \textrm{times}}\otimes\underbrace{\dual{\alg{V}}\otimes\dots
\otimes\dual{\alg{V}}}_{s\ \textrm{times}}.
\end{equation*}
By $\Symalg*[k]{\alg{V}}$ we denote the $k$-fold symmetric tensor
product of $\alg{V}$ with itself, and we think of this as a subset of
$\tensor*[k]{\alg{V}}$\@.  For $A\in\Symalg*[k]{\alg{V}}$ and
$B\in\Symalg*[l]{\alg{V}}$\@, we define the symmetric tensor product of $A$
and $B$ to be
\begin{equation*}
A\odot B=\sum_{\sigma\in\symmgroup{k,l}}\sigma(A\otimes B).
\end{equation*}
We define $\map{\Sym_k}{\tensor*[k]{\alg{V}}}{\Symalg*[k]{\alg{V}}}$ by
\begin{equation*}
\Sym_k(v_1\otimes\dots\otimes v_k)=
\frac{1}{k!}\sum_{\sigma\in\symmgroup{k}}v_{\sigma(1)}\otimes\dots
\otimes v_{\sigma(k)}.
\end{equation*}
We note that we have the alternative formula
\begin{equation}\label{eq:altodot}
A\odot B=\frac{(k+l)!}{k!l!}\Sym_{k+l}(A\otimes B)
\end{equation}
for the product of $A\in\Symalg*[k]{\alg{V}}$ and
$B\in\Symalg*[l]{\alg{V}}$\@.  We recall that
\begin{equation}\label{eq:symalgdim}
\dim_{\real}(\Symalg*[k]{\alg{V}})=\binom{\dim_{\real}(\alg{V})+k-1}{k},
\end{equation}
when $\alg{V}$ is finite-dimensional.

For a $\real$-vector space $\alg{V}$\@, let us denote
\begin{equation*}
\tensor*[\le m]{\alg{V}}=\bigoplus_{j=0}^m\tensor*[j]{\alg{V}},\quad
\Symalg*[\le m]{\alg{V}}=\bigoplus_{j=0}^m\Symalg*[j]{\alg{V}},
\end{equation*}
and define
\begin{equation*}
\mapdef{\Sym_{\le m}}{\tensor*[\le m]{\alg{V}}}{\Symalg*[\le m]{\alg{V}}}
{(A_0,A_1,\dots,A_m)}{(A_0,\Sym_1(A_1),\dots,\Sym_m(A_m)).}
\end{equation*}

For $\real$-inner product spaces $(\alg{U},\metric_{\alg{U}})$ and
$(\alg{V},\metric_{\alg{V}})$\@, we denote the transpose of
$L\in\Hom_{\real}(\alg{U};\alg{V})$ as the linear map
$\transpose{L}\in\Hom_{\real}(\alg{V};\alg{U})$ defined by
\begin{equation*}
\metric_{\alg{V}}(L(u),v)=\metric_{\alg{U}}(u,\transpose{L}(v)),\qquad
u\in\alg{U},\ v\in\alg{V}.
\end{equation*}

\subsubsection*{Topology}

We shall not use any particular notation for the Euclidean norm for
$\real^n$\@, and so will just denote this norm by
\begin{equation*}
\dnorm{\vect{x}}=\left(\sum_{j=1}^n\snorm{x_j}^2\right)^{1/2}.
\end{equation*}
It is sometimes convenient to use other norms for $\real^n$\@, particularly
the $1$- and $\infty$-norms defined, as usual, by
\begin{equation*}
\dnorm{\vect{x}}_1=\sum_{j=1}^n\snorm{x_j},\quad
\dnorm{\vect{x}}_\infty=\sup\setdef{\snorm{x_j}}{j\in\{1,\dots,n\}}.
\end{equation*}
The following relationships between these norms are useful:
\begin{equation}\label{eq:12inftynorms}
\begin{gathered}
\dnorm{x}\le\dnorm{\vect{x}}_1\le\sqrt{n}\dnorm{\vect{x}},\quad
\dnorm{\vect{x}}_\infty\le\dnorm{\vect{x}}\le\sqrt{n}\dnorm{\vect{x}}_\infty,\\
\dnorm{\vect{x}}_\infty\le\dnorm{\vect{x}}_1\le n\dnorm{\vect{x}}_\infty.
\end{gathered}
\end{equation}
If we are using a norm whose definition is evident from context, we will
simply denote it by $\dnorm{\cdot}$\@, accepting that context will ensure
that there is no confusion.

For $\vect{x}\in\real^n$ and $r\in\realp$\@, we denote by
\begin{equation*}
\oball{r}{\vect{x}}=\setdef{\vect{y}\in\real^n}{\dnorm{\vect{y}-\vect{x}}<r}
\end{equation*}
and
\begin{equation*}
\cball{r}{\vect{x}}=\setdef{\vect{y}\in\real^n}
{\dnorm{\vect{y}-\vect{x}}\le r}
\end{equation*}
the \defn{open} and \defn{closed balls} of radius $r$ centred at
$\vect{x}$\@.  As with the notation for norms, we shall often use the
preceding notation for balls in settings different from $\real^n$\@, and
accept the abuse of notation.

\subsubsection*{Differential calculus}

If $\nbhd{U}\subset\real^n$ is open and if
$\map{\vect{\Phi}}{\nbhd{U}}{\real^m}$ is differentiable at
$\vect{x}\in\nbhd{U}$\@, we denote its derivative by
$\linder{\vect{\Phi}}(\vect{x})$\@.  Higher-order derivatives, when they
exist, are denoted by $\linder[k]{\vect{\Phi}}(\vect{x})$\@, $k$ being the
order of differentiation.  We recall that, if
$\map{\vect{\Phi}}{\nbhd{U}}{\real^m}$ is of class $\C^k$\@,
$k\in\integerp$\@, then $\linder[k]{\vect{\Phi}}(\vect{x})$ is symmetric.  We
shall sometimes find it convenient to use multi-index notation for
derivatives.  A \defn{multi-index} with length $n$ is an element of
$\integernn^n$\@,~\ie~an $n$-tuple $I=(i_1,\dots,i_n)$ of nonnegative
integers.  If $\map{\vect{\Phi}}{\nbhd{U}}{\real^m}$ is a smooth function,
then we denote
\begin{equation*}
\linder[I]{\vect{\Phi}}(\vect{x})=
\plinder[i_1]{1}{\cdots\plinder[i_n]{n}{\vect{\Phi}}}(\vect{x}).
\end{equation*}
We will use the symbol $\snorm{I}=i_1+\dots+i_n$ to denote the order of the derivative.  Another piece of multi-index notation we shall use is
\begin{equation*}
\vect{a}^I=a_1^{i_1}\cdots a_n^{i_n},
\end{equation*}
for $\vect{a}\in\real^n$ and $I\in\integernn^n$\@.  Also, we denote $I!=i_1!\cdots i_n!$\@.

\subsubsection*{Differential geometry}

We shall adopt the notation and conventions of smooth differential geometry
of~\cite{RA/JEM/TSR:88}\@.  We shall also make use of real analytic
differential geometry.  There are no useful textbook references dedicated to
real analytic differential geometry, but the book of \cite{KC/YE:12} contains
much of what we shall need.  Throughout the paper, unless otherwise stated,
manifolds are connected, second countable, Hausdorff manifolds.  The
assumption of connectedness can be dispensed with but is convenient as it
allows one to not have to worry about manifolds with components of different
dimensions and vector bundles with fibres of different dimensions.

We shall work with regularity classes $r\in\{\infty,\omega\}$\@, ``$\infty$''
meaning smooth, ``$\omega$'' meaning real analytic.  Sometimes we do not
require infinite differentiability, but will hypothesise it anyway.  Other
times we will precisely specify the regularity needed; but we will be a
little sloppy with this as (1)~it is not crucial to the purposes of this
paper and (2)~it is typically easy to know when infinite differentiability is
hypothesised but not required.

The tangent bundle of a manifold $\man{M}$ is denoted by
$\map{\tbproj{\man{M}}}{\tb{\man{M}}}{\man{M}}$ and the cotangent bundle by
$\map{\ctbproj{\man{M}}}{\ctb{\man{M}}}{\man{M}}$\@.

We denote by $\mappings[r]{\man{M}}{\man{N}}$ the set of mappings from a
manifold $\man{M}$ to a manifold $\man{N}$ of class $\C^r$\@.  When
$\man{N}=\real$\@, we denote by
$\func[r]{\man{M}}=\mappings[r]{\man{M}}{\real}$ the set of scalar-valued
functions of class $\C^r$\@.  For $\Phi\in\mappings[1]{\man{M}}{\man{N}}$\@,
$\map{\tf{\Phi}}{\tb{\man{M}}}{\tb{\man{N}}}$ denotes the derivative of
$\Phi$\@, and $\tf[x]{\Phi}=\tf{\Phi}|\tb[x]{\man{M}}$\@.  For
$f\in\func[r]{\man{M}}$\@, we denote by
$\d{f}\in\sections[r-1]{\ctb{\man{M}}}$ the \defn{differential} of $f$\@,
defined by
\begin{equation*}
\tf[x]{f}(v_x)=(f(x),\natpair{\d{f}(x)}{v_x},\qquad v_x\in\tb[x]{\man{M}}.
\end{equation*}
We denote by $\ctf[x]{\Phi}$ the dual of $\tf[x]{\Phi}$\@.  For a vector
field $X$ and a differentiable function $f$\@, $\lieder{X}{f}$ denotes the
Lie derivative of $f$ with respect to $X$\@.  We might also write
$Xf=\lieder{X}{f}$\@.  For differentiable vector fields $X$ and $Y$\@, we
denote by $[X,Y]$ the Lie bracket of these vector fields.  For
$X\in\sections[r]{\tb{\man{M}}}$\@, the flow of $X$ is denoted by
$\flow{X}{t}$\@, meaning that, for $x\in\man{M}$\@, we have
\begin{equation*}
\deriv{}{t}\flow{X}{t}(x)=X\scirc\flow{X}{t}(x),\quad\flow{X}{0}(x)=x.
\end{equation*}

The Lie derivative for vector fields extends to a derivation of the tensor
algebra for a manifold.  Specifically, for
$X\in\sections[\infty]{\tb{\man{M}}}$\@, we denote
\begin{equation*}
\lieder{X}{f}=\natpair{\d{f}}{X},\quad\lieder{X}{Y}=[X,Y],\qquad
f\in\func[\infty]{\man{M}},\ X\in\sections[\infty]{\tb{\man{M}}}.
\end{equation*}
For $\alpha\in\sections[\infty]{\ctb{\man{M}}}$\@, we can then define its Lie
derivative with respect to $X$ by
\begin{equation*}
\natpair{\lieder{X}{\alpha}}{Y}=\lieder{X}{\natpair{\alpha}{Y}}-
\natpair{\alpha}{\lieder{X}{Y}},\qquad Y\in\sections[\infty]{\tb{\man{M}}}.
\end{equation*}
The Lie derivative of a tensor field
$A\in\sections[\infty]{\tensor[r,s]{\tb{\man{M}}}}$ is then defined by
\begin{multline}\label{eq:willmore}
\lieder{X}{A}(\alpha^1,\dots,\alpha^r,X_1,\dots,X_s)=
\lieder{X}{(A(\alpha^1,\dots,\alpha^r,X_1,\dots,X_s))}\\
-\sum_{j=1}^rA(\alpha^1,\dots,\lieder{X}{\alpha^j},\dots,\alpha^r,X_1,\dots,X_s)
-\sum_{j=1}^sA(\alpha^1,\dots,\alpha^r,X_1,\dots,\lieder{X}{X_j},\dots,X_s).
\end{multline}
Of course, these constructions make sense for tensor fields and vector fields
that are less regular than smooth.

Let $\map{\pi_{\man{E}}}{\man{E}}{\man{M}}$ be a vector bundle of class
$\C^r$\@.  We shall sometimes denote the fibre over $x\in\man{M}$ by
$\man{E}_x$\@, noting that this has the structure of a $\real$-vector space.
If $A\subset\man{M}$\@, we denote $\man{E}|A=\pi_{\man{E}}^{-1}(A)$\@.  By
$\sections[r]{\man{E}}$ we denote the set of sections of $\man{E}$ of class
$\C^r$\@.  This space has the structure of a $\real$-vector space with the
vector space operations
\begin{equation*}
(\xi+\eta)(x)=\xi(x)+\eta(x),\quad(a\xi)(x)=a(\xi(x)),\qquad x\in\man{M},
\end{equation*}
and of a $\func[r]{\man{M}}$-module with the additional operation of
multiplication
\begin{equation*}
(f\xi)(x)=f(x)\xi(x),\qquad x\in\man{M},
\end{equation*}
for $f\in\func[r]{\man{M}}$\@, $\xi,\eta\in\sections[r]{\man{E}}$\@, and
$a\in\real$\@.  By $\ssections[r]{\man{E}}$ we denote the sheaf of
$\C^r$-sections of $\man{E}$\@.  Thus
\begin{equation*}
\ssections[r]{\man{E}}(\nbhd{U})=\sections[r]{\man{E}|\nbhd{U}}
\end{equation*}
when $\nbhd{U}\subset\man{M}$ is open.  By $\real^k_{\man{M}}$ we denote the
trivial bundle $\real^k_{\man{M}}=\man{M}\times\real^k$ with vector bundle
projection being projection onto the first factor.  The dual bundle
$\dual{\man{E}}$ of a vector bundle $\man{E}$ is the set of vector bundle
mappings from $\man{E}$ to $\real_{\man{M}}$ over $\id_{\man{M}}$\@.  We note
that there is a natural identification of $\sections[r]{\real_{\man{M}}}$
with $\func[r]{\man{M}}$\@.  Given a $\C^r$-vector bundle
$\map{\pi_{\man{E}}}{\man{E}}{\man{M}}$ and a mapping
$\Phi\in\mappings[r]{\man{N}}{\man{M}}$\@, we denote by
$\map{\Phi^*\pi_{\man{E}}}{\Phi^*\man{E}}{\man{N}}$ the pull-back bundle.
For $\C^r$-vector bundles $\map{\pi_{\man{E}}}{\man{E}}{\man{M}}$ and
$\map{\pi_{\man{F}}}{\man{F}}{\man{M}}$ over the same base, we denote by
$\vbmappings[r]{\man{E}}{\man{F}}$ the set of $\C^r$-vector bundle mappings
from $\man{E}$ to $\man{F}$ over $\id_{\man{M}}$\@.

\subsubsection*{Riemannian geometry and connections}

We shall make use of basic constructions from Riemannian geometry.  We also
work a great deal with connections, both affine connections and linear
connections in vector bundles.  We refer to~\cite{SK/KN:63a} as a standard
reference, and~\cite{IK/PWM/JS:93} is also a useful reference.

First suppose that $r\in\{\infty,\omega\}$\@.  A \defn{$\C^r$-fibre metric}
on a $\C^r$-vector bundle $\map{\pi_{\man{E}}}{\man{E}}{\man{M}}$ is
$\metric_{\pi_{\man{E}}}\in\sections[r]{\Symalg*[2]{\dual{\man{E}}}}$ such
that $\metric_{\pi_{\man{E}}}(x)$ is an inner product on $\man{E}_x$ for each
$x\in\man{M}$\@.  The associated norm on fibres we denote by
$\dnorm{\cdot}_{\metric}$\@.  In case $\man{E}$ is the tangent bundle of
$\man{M}$\@, then a fibre metric is a Riemannian metric, and we will use the
notation $\metric_{\man{M}}$ in this case.

A linear connection in a vector bundle
$\map{\pi_{\man{E}}}{\man{E}}{\man{M}}$ will be denoted by
$\nabla^{\pi_{\man{E}}}$\@.  In case $\man{E}$ is the tangent bundle of
$\man{M}$\@, then a linear connection is called an affine connection, and we
will denote it by $\nabla^{\man{M}}$\@.  A linear connection in a vector
bundle $\map{\pi_{\man{E}}}{\man{E}}{\man{M}}$ induces a splitting of the
short exact sequence
\begin{equation*}
\xymatrix{{0}\ar[r]&{\ker(\tf{\pi_{\man{E}}})}\ar[r]&{\tb{\man{E}}}
\ar[r]^{\tf{\pi_{\man{E}}}}&{\tb{\man{M}}}\ar[r]&{0}}
\end{equation*}
For $e\in\man{E}$\@, we thus have a splitting of the tangent space
$\tb[e]{\man{E}}\simeq\tb[\pi_{\man{E}}(e)]{\man{M}}\oplus
\man{E}_{\pi_{\man{E}}(e)}$\@.  The first component in this splitting we call
\defn{horizontal} and denote by $\hb[e]{\man{E}}$\@, and the second we call
\defn{vertical} and denote by $\vb[e]{\man{E}}$\@.  By $\hor$ and $\ver$ we
denote the projections onto the horizontal and vertical subspaces,
respectively.

We note that covariant differentiation with respect to a vector field $X$ of
sections of $\man{E}$\@, along with Lie differentiation of functions, gives
rise to covariant differentiation of tensors, just as we saw above for
$\lieder{X}{}$\@.  A little more generally, if we have vector bundles
$\map{\pi_{\man{E}}}{\man{E}}{\man{M}}$ and
$\map{\pi_{\man{F}}}{\man{F}}{\man{E}}$\@, and linear connections
$\nabla^{\pi_{\man{E}}}$ and $\nabla^{\pi_{\man{F}}}$\@, then we have a
connection in $\man{E}\otimes\man{F}$ denoted by
$\nabla^{\pi_{\man{E}}\otimes\pi_{\man{F}}}$ and defined by
\begin{equation*}
\nabla^{\pi_{\man{E}}\otimes\pi_{\man{F}}}(\xi\otimes\eta)=
(\nabla^{\pi_{\man{E}}}\xi)\otimes\eta+\xi\otimes(\nabla^{\pi_{\man{F}}}\eta).
\end{equation*}

\subsubsection*{Jet bundles}

We shall make extensive use of jet bundles of various sorts.  We can
recommend~\cite{DJS:89} and~\cite[\S{}12]{IK/PWM/JS:93} as useful references.

Let $\man{M}$ be a $\C^r$-manifold and let $m\in\integernn$\@.  For
$x\in\man{M}$ and $a\in\real$\@, by $\jet[(x,a)]{m}{(\man{M};\real)}$ we
denote the $m$-jets of functions at $x$ taking value $a$ at $x$\@.  For a
$\C^r$-function $f$ defined in a neighbourhood of $x$\@, we denote by
$j_mf(x)\in\jet[(x,f(x))]{m}{\man{M}}$ the $m$-ket of $f$\@.  Of particular
interest is the set $\jetalg[x]{m}{\man{M}}=\jet[(x,0)]{m}{(\man{M};\real)}$
of jets of functions taking the value $0$ at $x$\@.  This has the structure
of a $\real$-algebra with the algebra structure defined by the three
operations
\begin{equation*}
j_mf(x)+j_mg(x)=j_m(f+g)(x),\quad(j_mf(x))(j_mg(x))=j_m(fg)(x),\quad
a(j_mf(x))=j_m(af)(x),
\end{equation*}
for functions $f$ and $g$ and for $a\in\real$\@.  We denote
\begin{equation*}
\jetalg{m}{\man{M}}=\bigdisjointunion_{x\in\man{M}}\jetalg[x]{m}{\man{M}}.
\end{equation*}
For $m,l\in\integernn$ with $m\ge l$\@, we have projections
$\map{\rho^m_l}{\jetalg{m}{\man{M}}}{\jetalg{l}{\man{M}}}$\@.  Note that
$\jetalg{0}{\man{M}}\simeq\man{M}$ and that
$\jetalg{1}{\man{M}}\simeq\ctb{\man{M}}$\@.  We abbreviate
$\map{\rho_m\eqdef\rho^m_0}{\jetalg{m}{\man{M}}}{\man{M}}$ which has the
structure of a vector bundle.

Let $\map{\pi_{\man{E}}}{\man{E}}{\man{M}}$ be a $\C^r$-vector bundle.  For
$x\in\man{M}$ and $m\in\integernn$\@, $\jet[x]{m}{\man{E}}$ denotes the set
of $m$-jets of sections of $\man{E}$ at $x$\@.  For a $\C^r$-section $\xi$
defined in some neighbourhood of $x$\@, $j_m\xi(x)\in\jet[x]{m}{\man{E}}$
denotes the $m$-jet of $\xi$\@.  We denote by
$\jet{m}{\man{E}}=\disjointunion_{x\in\man{M}}\jet[x]{m}{\man{E}}$ the bundle
of $m$-jets.  For $m,l\in\integernn$ with $m\ge l$\@, we denote by
$\map{\pi^m_l}{\jet{m}{\man{E}}}{\jet{l}{\man{E}}}$ the projection.  Note
that $\jet{0}{\man{E}}\simeq\man{E}$\@.  We abbreviate
$\map{\pi_m\eqdef\pi_{\man{E}}\scirc\pi^m_0}{\jet{m}{\man{E}}}{\man{M}}$\@,
and note that $\jet{m}{\man{E}}$ has the structure of a vector bundle over
$\man{M}$\@, with addition and scalar multiplication defined by
\begin{equation*}
j_m\xi(x)+j_m\eta(x)=j_m(\xi+\eta)(x),\quad a(j_m\xi(x))=j_m(a\xi)(x)
\end{equation*}
for sections $\xi$ and $\eta$ and for $a\in\real$\@.  One can show that
\begin{equation}\label{eq:jet=jetalg}
\jet{m}{\man{E}}\simeq(\real_{\man{M}}\oplus\jetalg{m}{\man{M}})\otimes\man{E}.
\end{equation}

\section{The topology for sections of a real analytic vector
bundle}\label{sec:ra-topology}

In this section we shall provide a quick overview of the usual topology for
real analytic sections of a real analytic vector bundle, and will give three
descriptions of this topology, two due to \citet{AM:66} and one via seminorms
given by \citet{SJ/ADL:14a}\@, based on the note of \citet{DV:13}\@.

\subsection{Martineau's descriptions of the real analytic topology}

We shall give a brief characterisation of two topologies for the space
$\sections[\omega]{\man{E}}$ of real analytic sections of a real analytic
vector bundle $\map{\pi_{\man{E}}}{\man{E}}{\man{M}}$\@.  The original work
of \citet{AM:66} describes these topologies for the space of real analytic
functions, but it is evident that the same considerations apply to sections
of a general vector bundle.  Each description offers benefits in terms of
providing immediately some useful properties of the topology, although
showing that they agree is something of an undertaking, and we shall make
some comments in this direction.

Both characterisations rely on the fact that a real analytic vector bundle
$\map{\pi_{\man{E}}}{\man{E}}{\man{M}}$ can be complexified to an holomorphic
vector bundle $\map{\pi_{\ol{\man{E}}}}{\ol{\man{E}}}{\ol{\man{M}}}$\@,
following \citet{HW/FB:59}\@.  We denote by $\sections[\hol]{\ol{\man{E}}}$
the space of holomorphic sections of this vector bundle, which we equip with
its usual compact-open topology,~\ie~the topology of uniform convergence on
compact sets.  This renders $\sections[\hol]{\ol{\man{E}}}$ a Fr\'echet
space.  For a subset $A\subset\ol{\man{M}}$\@, we denote by
$\gsections[\hol]{A}{\ol{\man{E}}}$ the space of germs of holomorphic
sections of $\ol{\man{E}}$ about $A$\@.  The space
$\gsections[\hol]{A}{\ol{\man{E}}}$ has the direct limit topology over the
directed set of neighbourhoods of $A$\@.

In the first description of the topology of $\sections[\omega]{\man{E}}$\@,
we note that, if $\xi\in\sections[\omega]{\man{E}}$\@, then there is some
neighbourhood $\ol{\nbhd{U}}$ of $\man{M}$ in $\ol{\man{M}}$ to which $\xi$
admits a unique holomorphic extension
$\ol{\xi}\in\sections[\hol]{\ol{\man{E}}|\ol{\nbhd{U}}}$\@.  Thus we have a
mapping
\begin{equation*}
\sections[\omega]{\man{E}}\ni\xi\mapsto\ol{\xi}\in
\gsections[\hol]{\man{M}}{\ol{\man{E}}}.
\end{equation*}
This map is easily seen to be an isomorphism of vector spaces, and so equips
$\sections[\omega]{\man{E}}$ with the direct limit topology for the space of
germs of sections of $\ol{\man{E}}$ about $\man{M}\subset\ol{\man{M}}$\@.
This immediately shows that $\sections[\omega]{\man{E}}$ is
ultrabornological~\cite[Corollaries~13.1.4 and~13.1.5]{HJ:81}\@.

The other description of a locally convex topology first fixes a compact
subset $\nbhd{K}\subset\man{M}$\@.  We note, then, that $\nbhd{K}$ possesses
a countable collection of neighbourhoods in $\ol{\man{M}}$ that are cofinal
in the directed set of all neighbourhoods.  Thus the direct limit topology of
$\gsections[\hol]{\nbhd{K}}{\ol{\man{E}}}$ is that of a countable direct
limit of Fr\'echet spaces.  Indeed, by working instead with bounded sections,
one ensures that one has a countable direct limit of Banach spaces.  One can
additionally and importantly show that the linking maps for the direct limit
are compact, indeed nuclear.  Thus the topology of
$\gsections[\hol]{\nbhd{K}}{\ol{\man{E}}}$ inherits many nice properties: it
is a webbed, nuclear, Suslin space
by~\cite[Corollary~5.3.3]{HJ:81}\@,~\cite[Theorem~8.4]{AK/PWM:97}
and~\cite[Example~II.2(E)]{LS:74}\@, respectively.  We next note that we have
a natural mapping $\xi\mapsto[\ol{\xi}]_{\nbhd{K}}$ from
$\sections[\omega]{\man{E}}$ to $\gsections[\hol]{\nbhd{K}}{\ol{\man{E}}}$ by
taking the germ about $\nbhd{K}$ of an holomorphic extension.  Now, since
$\man{M}$ is second countable, it possesses a countable compact exhaustion
$\ifam{\nbhd{K}_j}_{j\in\integerp}$\@, and one can then reasonably easily see
that $\sections[\omega]{\man{E}}$ is the inverse limit (as a vector space) of
the inverse system $\gsections[\hol]{\nbhd{K}_j}{\ol{\man{E}}}$\@,
$j\in\integerp$\@, (with linking maps given by restriction).  We can then
give $\sections[\omega]{\man{E}}$ the inverse limit topology.  The resulting
topology is webbed, nuclear, and Suslin
by~\cite[Corollary~5.3.3]{HJ:81}\@,~\cite[Corollary~21.2.3]{HJ:81}\@,
and~\cite[Lemma~6.6.5(ii) and~(iii)]{VIB:07b}\@.  It is not, however,
metrisable as follows from~\cite[Theorem~10]{DV:10}\@.

One of the contributions of \citet{AM:66} is to show that the two preceding
topologies for $\sections[\omega]{\man{E}}$ agree.  \citeauthor{AM:66}'s
original proof was by showing that
$\cup_{j\in\integerp}\algdual{(\gsections[\hol]{\nbhd{K}_j}{\ol{\man{E}}})}$
is a dense subspace of the dual of $\sections[\omega]{\man{E}}$ equipped with
the direct limit topology, using earlier results in~\cite{AM:63} on analytic
functionals.  A modern approach, using homological methods, equates an
inverse limit being ultrabornological with the vanishing of
$\mathupper{Proj}^1$\@, where $\mathupper{Proj}$ is a functor on inverse
systems devised by \citet{VPP:68}\@.  In all cases, showing equality of the
two topologies is not straightforward.

The inverse limit description of the topology for
$\sections[\omega]{\man{E}}$ is the one that is most closely connected with
our approach here, since the seminorms we give are essentially for
$\gsections[\hol]{\nbhd{K}}{\ol{\man{E}}}$ for a compact subset
$\nbhd{K}\subset\man{M}\subset\ol{\man{M}}$\@.  It is to the description of
these seminorms that we now turn.

\subsection{Decompositions for jet bundles}\label{subsec:jetdecomp}

A prominent r\^ole in our characterisation of the topology for real analytic
sections is played by jets and a decomposition of jets using connections.
The reason for this is that the seminorms we define are given in terms of
infinite jets of real analytic sections.

Let $\map{\pi_{\man{E}}}{\man{E}}{\man{M}}$ be a smooth vector bundle.  We
suppose that we have a linear connection $\nabla^{\pi_{\man{E}}}$ on the
vector bundle $\man{E}$ and an affine connection $\nabla^{\man{M}}$ on
$\man{M}$\@.  We then have induced connections, that we also denote by
$\nabla^{\pi_{\man{E}}}$ and $\nabla^{\man{M}}$\@, in various tensor bundles
of $\man{E}$ and $\tb{\man{M}}$\@, respectively.  The connections
$\nabla^{\pi_{\man{E}}}$ and $\nabla^{\man{M}}$ extend naturally to
connections in various tensor products of $\tb{\man{M}}$ and $\man{E}$\@, all
of these being denoted by $\nabla^{\man{M},\pi_{\man{E}}}$\@.  Note that
\begin{equation}\label{eq:nabla(j)}
\nabla^{\man{M},\pi_{\man{E}},m}\xi
\eqdef\underbrace{\nabla^{\man{M},\pi_{\man{E}}}\cdots
(\nabla^{\man{M},\pi_{\man{E}}}}_{m-1~\textrm{times}}(\nabla^{\pi_{\man{E}}}\xi))
\in\sections[\infty]{\tensor*[m]{\ctb{\man{M}}}\otimes\man{E}}.
\end{equation}
Now, given $\xi\in\sections[\infty]{\man{E}}$ and $m\in\integernn$\@, we
define
\begin{equation*}
D^m_{\nabla^{\man{M}},\nabla^{\pi_{\man{E}}}}(\xi)=\Sym_m\otimes\id_{\man{E}}
(\nabla^{\man{M},\pi_{\man{E}},m}\xi)\in
\sections[\infty]{\Symalg*[m]{\ctb{\man{M}}}\otimes\man{E}},
\end{equation*}
We take the convention that $D^0_{\nabla^{\man{M}},\nabla^{\pi_{\man{E}}}}(\xi)=\xi$\@.

The following lemma is then key for our presentation, and is proved
in~\cite[Lemma~2.1]{SJ/ADL:14a} by means of induction and a diagram chase.
\begin{lemma}\label{lem:Jmdecomp}
The map
\begin{equation*}
\mapdef{S_{\nabla^{\man{M}},\nabla^{\pi_{\man{E}}}}^m}
{\jet{m}{\man{E}}}
{\bigoplus_{j=0}^m(\Symalg*[j]{\ctb{\man{M}}}\otimes\man{E})}
{j_m\xi(x)}{(\xi(x),D^1_{\nabla^{\man{M}},\nabla^{\pi_{\man{E}}}}(\xi)(x),
\dots,D^m_{\nabla^{\man{M}},\nabla^{\pi_{\man{E}}}}(\xi)(x))}
\end{equation*}
is an isomorphism of vector bundles, and, for each\/ $m\in\integerp$\@, the
diagram
\begin{equation*}
\xymatrix{{\jet{m+1}{\man{E}}}\ar[r]^(0.3)
{S_{\nabla^{\man{M}},\nabla^{\pi_{\man{E}}}}^{m+1}}\ar[d]_{\pi^{m+1}_m}&
{\oplus_{j=0}^{m+1}(\Symalg*[j]{\ctb{\man{M}}}\otimes\man{E})}
\ar[d]^{\pr^{m+1}_m}\\{\jet{m}{\man{E}}}
\ar[r]_(0.3){S_{\nabla^{\man{M}},\nabla^{\pi_{\man{E}}}}^m}&
{\oplus_{j=0}^m(\Symalg*[j]{\ctb{\man{M}}}\otimes\man{E})}}
\end{equation*}
commutes, where\/ $\pr^{m+1}_m$ is the obvious projection, stripping off the
last component of the direct sum.
\end{lemma}

There are a couple of special cases of interest.
\begin{compactenum}
\item Jets of functions fit into the framework of the lemma by using the
trivial line bundle $\real_{\man{M}}=\man{M}\times\real$\@.  The
identification of a function with a section of this bundle is specified by
$f\mapsto\xi_f$\@, with $\xi_f(x)=(x,f(x))$\@.  In this case, the bundle has
a canonical flat connection defined by $\nabla^{\pi_{\man{E}}}f=\d{f}$\@.
Therefore, the decomposition of Lemma~\ref{lem:Jmdecomp} is determined by an
affine connection $\nabla^{\man{M}}$ on $\man{M}$\@, and so we have a mapping
\begin{equation}\label{eq:Smnablafunc}
\mapdef{S^m_{\nabla^{\man{M}}}}{\jet{m}{(\man{M};\real)}}
{\bigoplus_{j=0}^m\Symalg*[j]{\ctb{\man{M}}}}{f(x)}
{(f(x),\d{f}(x),\dots,\Sym_m\scirc\nabla^{\man{M},m-1}\d{f}(x)).}
\end{equation}
This can be restricted to $\jetalg{m}{\man{M}}$ to give the mapping
\begin{equation}\label{eq:SmnablaTm}
\mapdef{S^m_{\nabla^{\man{M}}}}{\jetalg{m}{\man{M}}}
{\bigoplus_{j=1}^m\Symalg*[j]{\ctb{\man{M}}}}{f(x)}
{(\d{f}(x),\dots,\Sym_m\scirc\nabla^{\man{M},m-1}\d{f}(x)),}
\end{equation}
adopting a mild abuse of notation.  We recall that $\jetalg[x]{m}{\man{M}}$ is
an $\real$-algebra, and the induced $\real$-algebra structure on
$\oplus_{j=1}^m\Symalg*[j]{\ctb[x]{\man{M}}}$ is that of polynomial functions
that vanish at $0$ and with degree at most $m$\@.

\item Another special case is that of jets of vector fields.  In this case,
the vector bundle is $\map{\tbproj{\man{M}}}{\tb{\man{M}}}{\man{M}}$\@.  We
can make use of an affine connection $\nabla^{\man{M}}$ on $\man{M}$ to
provide everything we need to define the mapping
\begin{equation}\label{eq:Smnablavf}
\mapdef{S^m_{\nabla^{\man{M}}}}{\jet{m}{\tb{\man{M}}}}
{\bigoplus_{j=0}^m(\Symalg*[j]{\ctb{\man{M}}}\otimes\tb{\man{M}})}{X(x)}
{(X(x),\nabla^{\man{M}}X(x),\dots,\Sym_m\scirc\nabla^{\man{M},m}X(x)).}
\end{equation}
Of course, this applies equally well to jets of one-forms on $\man{M}$\@, or
any other sections of tensor bundles associated with the tangent bundle.

This case of vector fields is the setting of \citet{SJ/ADL:14a} in their
study of flows of time-varying vector fields.
\end{compactenum}

\subsection{Fibre norms for jet bundles of vector
bundles}\label{subsec:jet-norms}

Our discussion begins with general constructions for the fibres of jet
bundles.  Thus we let $r\in\{\infty,\omega\}$ and let
$\map{\pi_{\man{E}}}{\man{E}}{\man{M}}$ be a $\C^r$-vector bundle.  We shall
suppose that we have a $\C^r$-affine connection $\nabla^{\man{M}}$ on
$\man{M}$ and a $\C^r$-vector bundle connection $\nabla^{\pi_{\man{E}}}$ in
$\man{E}$\@, as in Section~\ref{subsec:jetdecomp}\@.  This allows us to give
the decomposition of $\jet{m}{\man{E}}$ as in Lemma~\ref{lem:Jmdecomp}\@.  By
additionally supposing that we have a $\C^r$-Riemannian metric
$\metric_{\man{M}}$ on $\man{M}$ and a $\C^r$-fibre metric
$\metric_{\pi_{\man{E}}}$ on $\man{E}$\@, we shall give a $\C^r$-fibre norm
on $\jet{m}{\man{E}}$\@.  Note that the existence of the metrics and
connections is ensured by~\cite[Lemma~2.4]{SJ/ADL:14a}\@.

The first step in making the construction is the following result concerning
inner products on tensor products.
\begin{lemma}\label{lem:inprodotimes}
Let\/ $\alg{U}$ and\/ $\alg{V}$ be finite-dimensional\/ $\real$-vector spaces
and let\/ $\metric$ and\/ $\hmetric$ be inner products on\/ $\alg{U}$ and\/
$\alg{V}$\@, respectively.  Then the element\/ $\metric\otimes\hmetric$ of\/
$\tensor*[2]{\dual{\alg{U}}\otimes\dual{\alg{V}}}$ defined by
\begin{equation*}
\metric\otimes\hmetric(u_1\otimes v_1,u_2\otimes v_2)=
\metric(u_1,u_2)\hmetric(v_1,v_2)
\end{equation*}
is an inner product on\/ $\alg{U}\otimes\alg{V}$\@.
\begin{proof}
Let $\ifam{e_1,\dots,e_m}$ and $\ifam{f_1,\dots,f_n}$ be orthonormal bases
for $\alg{U}$ and $\alg{V}$\@, respectively.  Then
\begin{equation}\label{eq:GHbasis}
\setdef{e_a\otimes f_j}{a\in\{1,\dots,m\},\ j\in\{1,\dots,n\}}
\end{equation}
is a basis for $\alg{U}\otimes\alg{V}$\@.  Note that
\begin{equation*}
\metric\otimes\hmetric(e_a\otimes f_j,e_b\otimes f_k)=
\metric(e_a,e_b)\hmetric(f_j,f_k)=\delta_{ab}\delta_{jk},
\end{equation*}
which shows that $\metric\otimes\hmetric$ is indeed an inner product,
as~\eqref{eq:GHbasis} is an orthonormal basis.
\end{proof}
\end{lemma}

With $\metric_{\pi_{\man{E}}}$ a fibre metric on $\man{E}$ and with
$\metric_{\man{M}}$ be a Riemannian metric on $\man{M}$ as above, let us
denote by $\metric_{\man{M}}^{-1}$ the associated fibre metric on
$\ctb{\man{M}}$ defined by
\begin{equation*}
\metric_{\man{M}}^{-1}(\alpha_x,\beta_x)=
\metric_{\man{M}}(\metric_{\man{M}}^\sharp(\alpha_x),
\metric_{\man{M}}^\sharp(\beta_x)).
\end{equation*}
In like manner, one has a fibre metric $\metric_{\pi_{\man{E}}}^{-1}$ on
$\dual{\man{E}}$\@.  Then, by induction using the preceding lemma, we have a
fibre metric in all tensor spaces associated with $\tb{\man{M}}$ and
$\man{E}$ and their tensor products.  We shall denote by
$\metric_{\man{M},\pi_{\man{E}}}$ any of these various fibre metrics.  In
particular, we have a fibre metric $\metric_{\man{M},\pi_{\man{E}}}$ on
$\tensor*[j]{\ctb{\man{M}}}\otimes\man{E}$ induced by
$\metric_{\man{M}}^{-1}$ and $\metric_{\pi_{\man{E}}}$\@.  By restriction,
this gives a fibre metric on $\Symalg*[j]{\ctb{\man{M}}}\otimes\man{E}$\@.
We can thus define a fibre metric $\metric_{\man{M},\pi_{\man{E}},m}$ on
$\jet{m}{\man{E}}$ given by
\begin{equation}\label{eq:olGpim}
\metric_{\man{M},\pi_{\man{E}},m}(j_m\xi(x),j_m\eta(x))=
\sum_{j=0}^m\metric_{\man{M},\pi_{\man{E}}}
\left(\frac{1}{j!}D^j_{\nabla^{\man{M}},\nabla^{\pi_{\man{E}}}}(\xi)(x),
\frac{1}{j!}D^j_{\nabla^{\man{M}},\nabla^{\pi_{\man{E}}}}(\eta)(x)\right).
\end{equation}
Associated to this inner product on fibres is the norm on fibres, which we
denote by $\dnorm{\cdot}_{\metric_{\man{M},\pi_{\man{E}},m}}$\@.  We shall
use these fibre norms continually in our descriptions of our various
topologies for real analytic vector
bundles,~\cf~Section~\ref{sec:fibre-norms}\@.  The appearance of the
factorials in the fibre metric~\eqref{eq:olGpim} appears superfluous at this
point.  However, it is essential in order for the real analytic topology
defined by our seminorms to be independent of the choices of
$\nabla^{\man{M}}$\@, $\nabla^{\pi_{\man{E}}}$\@, $\metric_{\man{M}}$\@, and
$\metric_{\pi_{\man{E}}}$\@,~\cf~Theorem~\ref{the:same-topology}\@.

The preceding constructions can be applied particularly to the tangent bundle
of the total space of a vector bundle
$\map{\pi_{\man{E}}}{\man{E}}{\man{M}}$\@.  Indeed, given a Riemannian metric
on $\man{M}$\@, a fibre metric on $\man{E}$\@, an affine connection on
$\man{M}$\@, and a vector bundle connection in $\man{E}$\@, the constructions
of Section~\ref{subsec:Esubmersion} give a Riemannian metric on $\man{E}$\@,
and this, along with its Levi-Civita connection, gives the data required to
define fibre norms for the jet bundles $\jet{m}{\tb{\man{E}}}$\@.  A
substantial amount of the work in the paper will be to consider lifts to
$\man{E}$ of objects on $\man{M}$\@.  The continuity of operations like this
requires us to relate the jet bundle decompositions of
$\jet{m}{\tb{\man{E}}}$ with those of $\jet{m}{\man{E}}$ and
$\jet{m}{\tb{\man{M}}}$\@.

\subsection{Seminorms for the real analytic topology}\label{subsec:seminorms}

In this section we provide explicit seminorms for Martineau's topologies for
$\sections[\omega]{\man{E}}$\@.  Throughout this section, we will work with a
vector bundle $\map{\pi_{\man{E}}}{\man{E}}{\man{M}}$ and the data
$\nabla^{\man{M}}$\@, $\nabla^{\pi_{\man{E}}}$\@, $\metric_{\man{M}}$\@, and
$\metric_{\pi_{\man{E}}}$ that define the fibre metrics for jet bundles as
per Section~\ref{subsec:jet-norms}\@.  To define seminorms for
$\sections[\omega]{\man{E}}$\@, let $\c_0(\integernn;\realp)$ denote the
space of sequences in $\realp$\@, indexed by $\integernn$\@, and converging
to zero.  We shall denote a typical element of $\c_0(\integernn;\realp)$ by
$\vect{a}=\ifam{a_j}_{j\in\integernn}$\@.  Now, for $\nbhd{K}\subset\man{M}$
and $\vect{a}\in\c_0(\integernn;\realp)$\@, we define a seminorm
$p_{\nbhd{K},\vect{a}}^\omega$ for $\sections[\omega]{\man{E}}$ by
\begin{equation*}
p_{\nbhd{K},\vect{a}}^\omega(\xi)=\sup
\setdef{a_0a_1\cdots a_m\dnorm{j_m\xi(x)}_{\metric_{\man{M},\pi_{\man{E}},m}}}
{x\in\nbhd{K},\ m\in\integernn}.
\end{equation*}
The family of seminorms $p_{\nbhd{K},\vect{a}}^\omega$\@,
$\nbhd{K}\subset\man{M}$ compact, $\vect{a}\in\c_0(\integernn;\realp)$\@,
defines a locally convex topology on $\sections[\omega]{\man{E}}$ that we
call the \defn{$\C^\omega$-topology}\@.  As we have mentioned above, this
topology is webbed, ultrabornological, nuclear, and Suslin, but is not
metrisable.

While we are in the process of defining seminorms, let us also define
seminorms for the set $\sections[\infty]{\man{E}}$ of smooth sections.  While
we are primarily interested in the difficult real analytic case in this
paper, it is useful and illustrative to, at times, make comparisons with the
smooth case.  In any case, the topology we consider for
$\sections[\infty]{\man{E}}$ is that of uniform convergence of derivatives on
compact sets.  A moment's thought will convince one that the appropriate
seminorms are
\begin{equation*}
p^\infty_{\nbhd{K},m}(\xi)=\sup
\setdef{\dnorm{j_m\xi(x)}_{\metric_{\man{M},\pi_{\man{E}},m}}}{x\in\nbhd{K}}
\end{equation*}
for $\nbhd{K}\subset\man{M}$ compact and for $m\in\integernn$\@.  These
seminorms define a Polish topology for $\sections[\infty]{\man{E}}$ called
the \defn{$\C^\infty$-topology}\@.  We note that, for the smooth topology,
the seminorms are defined for fixed order jets.  As we shall indicate as we
go along, it is this fact that leads to simplifications of the results in the
paper when applied to the smooth case.

The following lemma, providing bounds for real analytic sections, is a global
version of a well-known classical
result~\cite[\eg][Proposition~2.2.10]{SGK/HRP:02}\@.  We refer
to~\cite[Lemma~2.6]{SJ/ADL:14a} for a proof.
\begin{lemma}\label{lem:Comegachar}
Let\/ $\map{\pi_{\man{E}}}{\man{E}}{\man{M}}$ be a real analytic vector
bundle.  For\/ $\xi\in\sections[\infty]{\man{E}}$\@, the following statements
are equivalent:
\begin{compactenum}[(i)]
\item \label{pl:Comegachar1} $\xi\in\sections[\omega]{\man{E}}$\@;
\item \label{pl:Comegachar2} for\/ $\nbhd{K}\subset\man{M}$ compact, there
exists\/ $C,r\in\realp$ such that
\begin{equation*}
p^\infty_{\nbhd{K},m}(\xi)\le Cr^{-m},\qquad m\in\integernn.
\end{equation*}
\end{compactenum}
\end{lemma}

\section{Tensors on the total space of a vector bundle}\label{sec:tensor-constructions}

Many of the geometric constructions we undertake in the paper, and estimates
associated with these geometric constructions, involves tensors of various
sorts defined on the total space of a vector bundle.  In this section we
present the classes of such tensors as arise in our presentation.  We also
define a number of algebraic operations on these tensors.  Many of the
constructions we see here will seem, on an initial reading, disconnected from
the objectives of the paper.  However, the constructions are essential in
Section~\ref{sec:tensor-derivatives}\@.  This is not very encouraging,
however, since the constructions and results of
Section~\ref{sec:tensor-derivatives} themselves appear \emph{non sequitur} to
the objectives of the paper.  It is only in the later sections of the paper
that the relevance of all of these constructions will become apparent.  For
this reason, perhaps a good strategy would be to skip over this section and
the next in a first reading, coming back to them when they are subsequently
needed.

There is nothing particularly real analytic with the material in this
section, so the smooth and real analytic cases are considered side-by-side.

\subsection{Functions on vector bundles}\label{subsec:vbfunctions}

Among the geometric constructions we will consider are those associated to a
particular set of functions on a vector bundle.
\begin{definition}
Let $r\in\{\infty,\omega\}$ and let $\map{\pi_{\man{E}}}{\man{E}}{\man{M}}$
be a vector bundle of class $\C^r$\@.  A function $F\in\func[r]{\man{E}}$ is
\defn{fibre-linear} if, for each $x\in\man{M}$\@, $F|\man{E}_x$ is a linear
function.  We denote by $\linfunc[r]{\man{E}}$ the set of $\C^r$-fibre-linear
functions on $\man{E}$\@.\oprocend
\end{definition}

Let us give some elementary properties of the sets of fibre-linear functions.
\begin{lemma}
Let\/ $r\in\{\infty,\omega\}$ and let\/
$\map{\pi_{\man{E}}}{\man{E}}{\man{M}}$ be a vector bundle of class\/
$\C^r$\@.  Then the following statements hold:
\begin{compactenum}[(i)]
\item \label{pl:fibre-linear1} $\linfunc[r]{\man{E}}$ is a submodule of
the\/ $\func[r]{\man{M}}$-module\/ $\func[r]{\man{E}}$\@;

\item \label{pl:fibre-linear2} for\/ $F\in\linfunc[r]{\man{E}}$\@, there
exists\/ $\lambda_F\in\sections[r]{\dual{\man{E}}}$ such that
\begin{equation*}
F(e)=\natpair{\lambda_F\scirc\pi_{\man{E}}(e)}{e},\qquad e\in\man{E},
\end{equation*}
and, moreover, the map\/ $F\mapsto\lambda_F$ is an isomorphism of\/
$\func[r]{\man{M}}$-modules;
\end{compactenum}
\begin{proof}
\eqref{pl:fibre-linear1} Let $F\in\linfunc[r]{\man{E}}$ and
$f\in\func[r]{\man{M}}$\@.  Then
\begin{equation*}
f\cdot F(e)=(f\scirc\pi_{\man{E}}(e))F(e),
\end{equation*}
and so $f\cdot F$ is fibre-linear since a scalar multiple of a linear
function is a linear function.  Also, since the pointwise sum of linear
functions is a linear function, we conclude that $\linfunc[r]{\man{E}}$ is
indeed a submodule of $\func[r]{\man{E}}$\@.

\eqref{pl:fibre-linear2} This merely follows by definition of the dual bundle
$\dual{\man{E}}$\@.
\end{proof}
\end{lemma}

In a rather related manner, we can consider other classes of functions on
vector bundles.
\begin{definition}\label{def:func1formlift}
Let $r\in\{\infty,\omega\}$ and let $\map{\pi_{\man{E}}}{\man{E}}{\man{M}}$
be a $\C^r$-vector bundle.
\begin{compactenum}[(i)]
\item For $\lambda\in\sections[r]{\dual{\man{E}}}$\@, the \defn{vertical
evaluation} of $\lambda$ is $\ve{\lambda}\in\linfunc[r]{\man{E}}$ defined
by $\ve{\lambda}(e_x)=\natpair{\lambda(x)}{e_x}$\@.
\item For $f\in\func[r]{\man{M}}$\@, the \defn{horizontal lift} of $f$ is
the function $\hl{f}\in\func[r]{\man{E}}$ defined by
$\hl{f}=\pi_{\man{E}}^*f$\@.\oprocend
\end{compactenum}
\end{definition}

\subsection{Vector fields on vector bundles}\label{subsec:vbvfs}

Next we turn to vector fields on the total space of a vector bundle.  As with
our consideration of functions in the preceding section, we restrict
attention to vector fields that interact nicely with the vector bundle
structure.

We begin with the notion of the vertical lift of a section.
\begin{definition}
Let $r\in\{\infty,\omega\}$ and let $\map{\pi_{\man{E}}}{\man{E}}{\man{M}}$
be a vector bundle of class $\C^r$\@.
\begin{compactenum}[(i)]
\item For $e_x,e'_x\in\man{E}_x$\@, we define the \defn{vertical lift} of
$e'_x$ to $e_x$ to be
\begin{equation*}
\verlift(e_x,e'_x)=\derivatzero{}{t}(e_x+te'_x).
\end{equation*}
\item Given a section $\xi\in\sections[r]{\man{E}}$\@, we define the
\defn{vertical lift} of $\xi$ to $\man{E}$ to be the vector field
\begin{equation*}\eqoprocend
\vl{\xi}(e_x)=\verlift(e_x,\xi(x)).
\end{equation*}
\end{compactenum}
\end{definition}

Next we consider another sort of lift, this one requiring a connection
$\nabla^{\pi_{\man{E}}}$ in the vector bundle
$\map{\pi_{\man{E}}}{\man{E}}{\man{B}}$\@.  We let
$\vb{\man{E}}=\ker(\tf{\pi_{\man{E}}})$ be the vertical subbundle.  As
mentioned in Section~\ref{subsec:background}\@, the connection
$\nabla^{\pi_{\man{E}}}$ defines a complement $\hb{\man{E}}$ to
$\vb{\man{E}}$ called the horizontal subbundle.  We let
$\map{\ver,\hor}{\tb{\man{E}}}{\tb{\man{E}}}$ be the projections onto
$\vb{\man{E}}$ and $\hb{\man{E}}$\@.
\begin{definition}
Let $r\in\{\infty,\omega\}$\@, let $\map{\pi_{\man{E}}}{\man{E}}{\man{M}}$ be
a vector bundle of class $\C^r$\@, and let $\nabla^{\pi_{\man{E}}}$ be a
$\C^r$-connection in $\man{E}$\@.
\begin{compactenum}[(i)]
\item For $e_x\in\man{E}_x$ and $v_x\in\tb[x]{\man{M}}$\@, the
\defn{horizontal lift} of $v_x$ to $e_x$ is the unique vector
$\horlift(e_x,v_x)\in\hb[e_x]{\man{E}}$ satisfying
\begin{equation*}
\tf[e_x]{\pi_{\man{E}}}(\horlift(e_x,v_x))=v_x.
\end{equation*}
\item For $X\in\sections[r]{\tb{\man{M}}}$ on $\man{M}$\@, we denote by
$\hl{X}$ the \defn{horizontal lift} of $X$ to $\man{E}$\@, this being the
vector field $\hl{X}\in\sections[r]{\tb{\man{E}}}$ satisfying
\begin{equation*}\eqoprocend
\hl{X}(e_x)=\horlift(e_x,X(x)).
\end{equation*}
\end{compactenum}
\end{definition}

Next we provide formulae for differentiating various sorts of functions with
respect to various sorts of vector fields.
\begin{lemma}\label{lem:bundlefuncs}
Let\/ $r\in\{\infty,\omega\}$ and let\/
$\map{\pi_{\man{E}}}{\man{E}}{\man{M}}$ a vector bundle of class\/ $\C^r$\@.
Let\/ $f\in\func[r]{\man{M}}$\@,\/
$\lambda\in\sections[r]{\dual{\man{E}}}$\@,\/
$X\in\sections[r]{\tb{\man{M}}}$\@, and\/ $\xi\in\sections[r]{\man{E}}$\@,

Then the following statements hold:
\begin{compactenum}[(i)]
\item \label{pl:bundlefuncs1} $\lieder{\hl{X}}{\hl{f}}=\hl{(\lieder{X}{f})}$\@;
\item \label{pl:bundlefuncs2} $\lieder{\vl{\xi}}{\hl{f}}=0$\@;
\item \label{pl:bundlefuncs3}
$\lieder{\vl{\xi}}{\ve{\lambda}}=\hl{\natpair{\lambda}{\xi}}$\@;\savenum
\end{compactenum}
Additionally, let\/ $\nabla^{\pi_{\man{E}}}$ be a\/ $\C^r$-linear connection in\/
$\map{\pi_{\man{E}}}{\man{E}}{\man{M}}$\@.  Then
\begin{compactenum}[(i)]\resumenum
\item \label{pl:bundlefuncs4}
$\lieder{\hl{X}}{\ve{\lambda}}=\ve{(\nabla^{\pi_{\man{E}}}_X\lambda)}$\@.
\end{compactenum}
\begin{proof}
\eqref{pl:bundlefuncs1} We compute
\begin{align*}
\lieder{\hl{X}}{\hl{f}}(e)=&\;\natpair{\d{(\pi_{\man{E}}^*f)}(e)}{\hl{X}(e)}=
\natpair{\d{f}\scirc\pi_{\man{E}}(e)}{\tf[e]{\pi_{\man{E}}}(\hl{X}(e))}\\
=&\;\natpair{\d{f}\scirc\pi_{\man{E}}(e)}{X\scirc\pi_{\man{E}}(e)}=
\hl{(\lieder{X}{f})}(e).
\end{align*}

\eqref{pl:bundlefuncs2} Since $\hl{f}$ is constant on fibres of
$\pi_{\man{E}}$ and $\vl{\xi}$ is tangent to fibres, we have
\begin{equation*}
\hl{f}(e+t\xi\scirc\pi_{\man{E}}(e))=f(e).
\end{equation*}
Differentiating with respect to $t$ at $t=0$ gives the result.

\eqref{pl:bundlefuncs3} Here we compute
\begin{align*}
\lieder{\vl{\xi}}{\ve{\lambda}}(e)=&\;
\derivatzero{}{t}\natpair{\lambda(e+t\xi\scirc\pi_{\man{E}}(e))}
{e+t\xi\scirc\pi_{\man{E}}(e)}\\
=&\;\natpair{\lambda\scirc\pi_{\man{E}}(e)}{\xi\scirc\pi_{\man{E}}(e)}
=\hl{\natpair{\lambda}{\xi}},
\end{align*}
so completing the proof.

\eqref{pl:bundlefuncs4} Let $e\in\man{E}$ and let $t\mapsto\gamma(t)$ be the
integral curve for $X$ satisfying $\gamma(0)=\pi_{\man{E}}(e)$ and let
$t\mapsto\hl{\gamma}(t)$ be the integral curve for $\hl{X}$ satisfying
$\hl{\gamma}(0)=e$\@.  Then $t\mapsto\hl{\gamma}(t)$ is the parallel
translation of $e$ along $\gamma$\@, and as such we have
$\nabla^{\pi_{\man{E}}}_{\gamma'(t)}\hl{\gamma}(t)=0$\@.  Then
\begin{equation*}
\lieder{\hl{X}}{\ve{\lambda}}(e)=
\derivatzero{}{t}\natpair{\lambda\scirc\gamma(t)}{\hl{\gamma}(t)}=
\natpair{\nabla^{\pi_{\man{E}}}_X\lambda\scirc\pi_{\man{E}}(e)}{e},
\end{equation*}
as claimed.
\end{proof}
\end{lemma}

In Section~\ref{sec:tensor-derivatives} we shall have a great deal more to
say about differentiation of objects on the total space of a vector bundle
when one has more structure present than we use in the preceding result.

\subsection{Linear mappings on vector bundles}\label{subsec:vblinmaps}

Now we turn to an examination of linear maps associated to a vector bundle
$\map{\pi_{\man{E}}}{\man{E}}{\man{M}}$\@.  We shall consider vector bundle
mappings of two sorts:~(1)~with values in the trivial line bundle
$\real_{\man{M}}$\@; (2)~with values in $\man{E}$\@.  The first sort of
mappings are, of course, simply sections of the dual bundle, or linear
functions of the sort studied in Section~\ref{subsec:vbfunctions}\@.  Our
interest here is in lifting such objects to the total space.

First we work with one-forms.  If we have a connection
$\nabla^{\pi_{\man{E}}}$ in a vector bundle
$\map{\pi_{\man{E}}}{\man{E}}{\man{M}}$\@, then this gives us a splitting
$\tb{\man{E}}=\hb{\man{E}}\oplus\vb{\man{E}}$\@, and hence a splitting
$\ctb{\man{E}}=\chb{\man{E}}\oplus\cvb{\man{E}}$ with
\begin{equation*}
\chb{\man{E}}=\ann(\vb{\man{E}}),\quad\cvb{\man{E}}=\ann(\hb{\man{E}}).
\end{equation*}
Note that $\chb[e]{\man{E}}=\image(\ctf[e]{\pi_{\man{E}}})$\@.  
\begin{definition}\label{def:1formlift}
Let $r\in\{\infty,\omega\}$ and let $\map{\pi_{\man{E}}}{\man{E}}{\man{M}}$
be a $\C^r$-vector bundle.
\begin{compactenum}[(i)]
\item For $\alpha_x\in\ctb[x]{\man{M}}$ and $e_x\in\man{E}_x$\@, the
\defn{horizontal lift} of $\alpha_x$ to $e_x$ is
$\horlift(e_x,\alpha_x)=\ctf[e_x]{\pi_{\man{E}}}(\alpha_x)$\@.
\item The \defn{horizontal lift} of $\alpha\in\sections[r]{\ctb{\man{M}}}$ is
$\hl{\alpha}=\pi_{\man{E}}^*\alpha\in\sections[r]{\ctb{\man{E}}}$\@.\savenum
\end{compactenum}
Additionally, let $\nabla^{\pi_{\man{E}}}$ be a connection in $\man{E}$\@.
\begin{compactenum}[(i)]\resumenum
\item For $\lambda_x\in\dual{\man{E}}_x$ and $e_x\in\dual{\man{E}}_x$\@, then
the \defn{vertical lift} of $\lambda_x$ is the unique vector
$\verlift(e_x,\lambda_x)\in\cvb[e_x]{\man{E}}$ satisfying
\begin{equation*}
\natpair{\verlift(e_x,\lambda_x)}{\verlift(e_x,u_x)}=\natpair{\lambda_x}{u_x}
\end{equation*}
for every $u_x\in\man{E}_x$\@.
\item The \defn{vertical lift} of $\lambda\in\sections[r]{\dual{\man{E}}}$
is the one-form $\vl{\lambda}\in\sections[r]{\ctb{\man{E}}}$ satisfying
\begin{equation*}\eqoprocend
\vl{\lambda}(e_x)=\verlift(e_x,\lambda(x)).
\end{equation*}
\end{compactenum}
\end{definition}

We also have natural ways of lifting homomorphisms of vector bundles.
\begin{definition}\label{def:homlift}
Let $r\in\{\infty,\omega\}$\@, and let
$\map{\pi_{\man{E}}}{\man{E}}{\man{M}}$ and
$\map{\pi_{\man{F}}}{\man{F}}{\man{M}}$ be $\C^r$-vector bundles.  For
$L\in\sections[r]{\man{F}\otimes\dual{\man{E}}}$\@,
\begin{compactenum}[(i)]
\item the \defn{vertical evaluation} of $L$ is the section
$\ve{L}\in\sections[r]{\pi_{\man{E}}^*\man{F}}$ defined by
\begin{equation*}
\ve{L}(e_x)=(e_x,L(e_x)).
\end{equation*}
\savenum
\end{compactenum}
If, additionally, $\nabla^{\pi_{\man{E}}}$ is a connection in $\man{E}$\@,
\begin{compactenum}[(i)]\resumenum
\item the \defn{vertical lift} of $L$ is the vector bundle homomorphism
$\vl{L}\in\sections[r]{\pi_{\man{E}}^*\man{F}\otimes\ctb{\man{E}}}$ defined by
\begin{equation*}
\vl{L}(Z)=(e,L\scirc\ver(Z))
\end{equation*}
for\/ $Z\in\tb[e]{\man{E}}$\@, noting that $\ver(Z)\in\vb[e]{\man{E}}\simeq\man{E}_{\pi_{\man{E}}(e)}$\@.\oprocend
\end{compactenum}
\end{definition}

We shall be especially interested in two cases of the vector bundle
$\man{F}$\@.
\begin{compactenum}
\item $\man{F}=\real_{\man{M}}$\@: In this case,
$\man{F}\otimes\dual{\man{E}}\simeq\dual{\man{E}}$\@,
$\pi_{\man{E}}^*\man{F}\simeq\real_{\man{E}}$\@, and
$\pi_{\man{E}}^*\man{F}\otimes\ctb{\man{E}}\simeq\ctb{\man{E}}$\@.  One can
easily see that, if $\lambda\in\sections[r]{\dual{\man{E}}}$\@, then the
vertical evaluation as per Definition~\ref{def:homlift} agrees with that of
Definition~\ref{def:func1formlift}\@, and the vertical lift as per
Definition~\ref{def:homlift} agrees with that of
Definition~\ref{def:1formlift}\@.

\item $\man{F}=\man{E}$\@: In this case,
$\man{F}\otimes\dual{\man{E}}\simeq\tensor[1,1]{\man{E}}$\@,~\ie~the set of
endomorphisms of $\man{E}$\@.  We also have
$\pi_{\man{E}}^*\man{F}\simeq\vb{\man{E}}$~\cite[\S6.11]{IK/PWM/JS:93}\@.
Thus, for $L\in\sections[r]{\tensor[1,1]{\man{E}}}$\@, $\ve{L}$ is a
$\vb{\man{E}}$-valued vector field.  Also, $\vl{L}$ is a
$\vb{\man{E}}$-valued endomorphism of $\tb{\man{E}}$\@.
\end{compactenum}

Let us perform some analysis of the vertical evaluation and vertical
lift of an homomorphism.  First of all, for $e_1,e_2\in\man{E}_x$\@,
\begin{equation*}
\ve{L}(e_1+e_2)=(e_1,L(e_1))+(e_2,L(e_2))=
\ve{L}(e_1)+\ve{L}(e_2),
\end{equation*}
where addition is with respect to the vector bundle structure
\begin{equation*}
\xymatrix{{\man{E}}\ar[r]^{\ve{L}}\ar[d]_{\pi_{\man{E}}}&
{\pi_{\man{E}}^*\man{F}}\ar[d]_{\pi_{\man{E}}^*\pi_{\man{F}}}\\{\man{M}}\ar[r]_Z&
{\man{F}}}
\end{equation*}
where $Z$ is the zero section.  Thus $\ve{L}$ is a ``linear'' section over
$\man{E}$\@.  We define the vector bundle mapping
\begin{equation}\label{eq:PEFdef}
\mapdef{P_{\man{E},\man{F}}}{\pi_{\man{E}}^*\man{F}\otimes\cvb{\man{E}}}
{\pi_{\man{E}}^*\man{F}}{L_e}{L_e(e)}
\end{equation}
over $\id_{\man{E}}$\@, noting that
$e\in\man{E}_{\pi_{\man{E}}(e)}\simeq\vb[e]{\man{E}}$\@.  Then, given
$A\in\sections[r]{\pi_{\man{E}}^*\man{F}\otimes\cvb{\man{E}}}$\@,
$P_{\man{E},\man{F}}\scirc A$ is a section of $\pi_{\man{E}}^*\man{E}$\@.  Moreover,
$P_{\man{E},\man{F}}\scirc\vl{L}=\ve{L}$ for
$L\in\sections[r]{\man{F}\otimes\dual{\man{E}}}$\@.

We shall make use of these observations in Section~\ref{sec:lift-isomorphisms}\@.

Let us recast the preceding observations in a slightly different way.  To
start, note that, given $\lambda\in\sections[r]{\dual{\man{E}}}$ and
$\eta\in\sections[r]{\man{F}}$\@, we have
$\eta\otimes\lambda\in\sections[r]{\man{F}\otimes\dual{\man{E}}}$\@.  The
tensor product on the left can be thought of as being of
$\func[r]{\man{M}}$-modules.\footnote{This corresponds to the well-known
isomorphism
\begin{equation*}
\sections[r]{\man{E}}\otimes_{\func[r]{\man{M}}}\sections[r]{\man{F}}
\simeq\sections[r]{\man{E}\otimes\man{F}}
\end{equation*}
of $\func[r]{\man{M}}$-modules.  While this isomorphism is well-known, it is
commonly not correctly proved, as proofs are given that admit a direct
translation to the holomorphic setting, where the assertion is generally
false.  A correct proof in the smooth case is given by
\citet[Theorem~7.5.5]{LC:01}\@.  His proof makes use (without saying this
explicitly) of the Serre\textendash{}Swan Theorem.  Since the
Serre\textendash{}Swan Theorem is valid for vector bundles over smooth, real
analytic, and Stein manifolds (see~\cite[Theorem~6.5]{ADL:20c}),
\citeauthor{LC:01}'s proof applies in these cases.}  Moreover, such sections
of the bundle of endomorphisms locally generate the sections of the
homomorphism bundle.  Note that
\begin{equation*}
\ve{(\eta\otimes\lambda)}=\vl{\xi}\otimes\ve{\lambda},
\end{equation*}
as is directly verified.  In this case, since $\func[r]{\man{M}}$ is a
subring of $\func[r]{\man{E}}$ (by pull-back), we can regard the tensor
product as being of $\func[r]{\man{E}}$-modules.  Therefore,
\begin{equation*}
\ve{L}\in\sections[r]{\sections[r]{\pi_{\man{E}}^*\man{F}}\otimes
\linfunc[r]{\man{E}}}.
\end{equation*}
Since $\linfunc[r]{\man{E}}\subset\func[r]{\man{E}}$\@, the tensor product is
mere multiplication in this case.

A similar sort of analysis can be made for the vertical lift of an
homomorphism.  In this case, given $\lambda\in\sections[r]{\dual{E}}$ and
$\eta\in\sections[r]{\man{F}}$\@, we have
$\xi\otimes\lambda\in\sections[r]{\man{F}\otimes\dual{\man{E}}}$\@, as in
the preceding paragraph.  In this case, the vertical lift satisfies
\begin{equation*}
\vl{(\xi\otimes\lambda)}=\vl{\xi}\otimes\vl{\lambda}.
\end{equation*}

\subsection{Tensors fields on vector bundles}\label{subsec:vbtensors}

Next we discuss the extension of our lifts of functions, sections, and vector
fields to higher-order tensors.  The extension is to tensor powers of the
pull-back $\pi_{\man{E}}^*\ctb{\man{M}}$ of the cotangent bundle to the total space of
the vector bundle.  Other sorts of lifts are possible, especially in the
presence of a connection in the vector bundle.  We restrict ourselves to the
tensor powers of the pull-back of $\ctb{\man{M}}$ since our interest is in
jet bundles, and these tensor powers represent derivatives with respect to
the base.

We make the following definitions.
\begin{definition}\label{def:tensor-lifts}
Let $r\in\{\infty,\omega\}$\@, and let
$\map{\pi_{\man{E}}}{\man{E}}{\man{M}}$ and
$\map{\pi_{\man{F}}}{\man{F}}{\man{M}}$ be a $\C^r$-vector bundles.  Let
$k\in\integerp$\@.
\begin{compactenum}[(i)]
\item For $A\in\sections[r]{\tensor*[k]{\ctb{\man{M}}}}$\@, the
\defn{horizontal lift} of $A$ is
$\hl{A}\in\sections[r]{\tensor*[k]{\ctb{\man{E}}}}$ defined by
\begin{equation*}
\hl{A}(Z_1,\dots,Z_k)=
A(\tf[e]{\pi_{\man{E}}}(Z_1),\dots,\tf[e]{\pi_{\man{E}}}(Z_k))
\end{equation*}
for $Z_1,\dots,Z_k\in\tb[e]{\man{E}}$\@.\footnote{Of course, this is nothing
but the usual definition of pull-back, which we repeat for the sake of
symmetry.}
\item For $A\in\sections[r]{\tensor*[k]{\ctb{\man{M}}}\otimes\man{E}}$\@,
the \defn{vertical lift} of $A$ is
$\vl{A}\in\sections[r]{\tensor*[k]{\ctb{\man{E}}}\otimes\tb{\man{E}}}$
defined by
\begin{equation*}
\vl{A}(Z_1,\dots,Z_k)=
\verlift(e,A(\tf[e]{\pi_{\man{E}}}(Z_1),\dots,\tf[e]{\pi_{\man{E}}}(Z_k))),
\end{equation*}
for $Z_1,\dots,Z_k\in\tb[e]{\man{E}}$\@.
\item For $A\in\sections[r]{\tensor*[k]{\ctb{\man{M}}}\otimes\man{F}
\otimes\dual{\man{E}}}$\@, the \defn{vertical evaluation} of $A$ is
$\ve{A}\in\sections[r]{\tensor*[k]{\ctb{\man{E}}}\otimes
\pi_{\man{E}}^*\man{F}}$
defined by
\begin{equation*}
\ve{A}(Z_1,\dots,Z_k)=
(e,A(\tf[e]{\pi_{\man{E}}}(Z_1),\dots,\tf[e]{\pi_{\man{E}}}(Z_k))(e)),
\end{equation*}
for $Z_1,\dots,Z_k\in\tb[e]{\man{E}}$\@.\savenum
\end{compactenum}
Additionally, let $\nabla^{\pi_{\man{E}}}$ be a connection in $\man{E}$\@.
\begin{compactenum}[(i)]\resumenum
\item For
$A\in\sections[r]{\tensor*[k]{\ctb{\man{M}}}\otimes\tb{\man{M}}}$\@, the
\defn{horizontal lift} of $A$ is
$\hl{A}\in\sections[r]{\tensor*[k]{\ctb{\man{E}}}\otimes\tb{\man{E}}}$
defined by
\begin{equation*}
\hl{A}(Z_1,\dots,Z_k)=
\horlift(e,A(\tf[e]{\pi_{\man{E}}}(Z_1),\dots,\tf[e]{\pi_{\man{E}}}(Z_k)))
\end{equation*}
for $Z_1,\dots,Z_k\in\tb[e]{\man{E}}$\@.
\item For
$A\in\sections[r]{\tensor*[k]{\ctb{\man{M}}}\otimes\dual{\man{E}}}$\@, the
\defn{vertical lift} of $A$ is
$\vl{A}\in\sections[r]{\tensor*[k]{\ctb{\man{E}}}\otimes\ctb{\man{E}}}$
defined by
\begin{equation*}
\vl{A}(Z_1,\dots,Z_k)=
\verlift(e,A(\tf[e]{\pi_{\man{E}}}(Z_1),\dots,\tf[e]{\pi_{\man{E}}}(Z_k)))
\end{equation*}
for $Z_1,\dots,Z_k\in\tb[e]{\man{E}}$\@.
\item For $A\in\sections[r]{\tensor*[k]{\ctb{\man{M}}}\otimes\man{F}
\otimes\dual{\man{E}}}$\@, the \defn{vertical lift} of $A$ is
$\vl{A}\in\sections[r]{\tensor*[k]{\ctb{\man{E}}}\otimes
\pi_{\man{E}}^*\man{F}\otimes\ctb{\man{E}}}$
defined by
\begin{equation*}
\vl{A}(Z_1,\dots,Z_k)(Z)=
(e,A(\tf[e]{\pi_{\man{E}}}(Z_1),\dots,\tf[e]{\pi_{\man{E}}}(Z_k))(\ver(Z))),
\end{equation*}
for $Z_1,\dots,Z_k,Z\in\tb[e]{\man{E}}$\@.\oprocend
\end{compactenum}
\end{definition}

\subsection{Tensor contractions}\label{subsec:Ins}

In our differentiation results of Section~\ref{sec:tensor-derivatives}\@, we
shall make use of certain generalisations of the contraction operator on
tensors.  What we have is a sort of ``contraction and insertion'' operation.
We describe this here in the setting of linear algebra, since this is where
it most naturally resides.  The constructions can, of course, be extended to
vector bundles by performing the vector space constructions on fibres.

Let $\alg{V}$ be a finite-dimensional $\real$-vector space, let
$k\in\integerp$ and $l\in\integernn$\@, and let
$\alpha\in\tensor*[k]{\dual{\alg{V}}}$ and
$\beta\in\tensor*[l]{\dual{\alg{V}}}\otimes\alg{V}$\@.  For
$j\in\{1,\dots,k\}$\@, define the \defn{$j$th insertion of $\beta$ in
$\alpha$} by $\Ins_j(\alpha,\beta)\in\tensor*[k+l-1]{\dual{\alg{V}}}$ by
\begin{equation*}
\Ins_j(\alpha,\beta)(v_1,\dots,v_{k+l-1})=
\alpha(v_1,\dots,v_{j-1},\beta(v_j,v_{k+1},\dots,v_{k+l-1}),v_{j+1},\dots,v_k).
\end{equation*}
To be clear, when $l=0$ we have
\begin{equation*}
\Ins_j(\alpha,v)(v_1,\dots,v_{k-1})=
\alpha(v_1,\dots,v_{j-1},v,v_j,\dots,v_{k-1}).
\end{equation*}
We will also find it helpful to consider tensor contraction when one of the
arguments (the second is the one we care about) is fixed.  Thus let
$\beta\in\tensor*[l]{\dual{\alg{V}}}\otimes\alg{V}$ and define
$\Ins_{j,\beta}(\alpha)=\Ins_j(\alpha,\beta)$\@.

We shall also need notation for a specific sort of swapping of arguments of a
tensor.  Let $\alpha\in\tensor*[k]{\alg{V}}$ and let
$j_1,j_2\in\{1,\dots,k\}$\@.  We define
\begin{equation*}
\push_{j_1,j_2}\alpha(v_1,\dots,v_k)=\begin{cases}
\alpha(v_1,\dots,v_{j_1-1},v_{j_1+1},\dots,v_{j_2},v_{j_1},v_{j_2+1},\dots,v_k),&
j_1\le j_2,\\
\alpha(v_1,\dots,v_{j_2-1},v_{j_1},v_{j_2},\dots,v_{j_1-1},v_{j_1+1},\dots,v_k),&
j_1>j_2.\end{cases}
\end{equation*}
The idea is that $\push_{j_1,j_2}$ drops $v_{j_1}$ into the $j_2$-slot, and
shifts the arguments to make room for this.  The ``insertion'' and ``push''
mappings can be generalised in the obvious way to give $\Ins_j(A,\beta)$ and
$\push_{j_1,j_2}(A)$ for $A\in\tensor*[k]{\dual{\alg{V}}}\otimes\alg{U}$ and
$\beta\in(\tensor*[l]{\dual{\alg{V}}}\otimes\alg{V})\otimes\alg{U}$
(\resp~$A\in\alg{U}\otimes\tensor*[k]{\dual{\alg{V}}}$ and
$\beta\in\alg{U}\otimes(\tensor*[l]{\dual{\alg{V}}}\otimes\alg{V})$)\@, just
by acting on the first (\resp~second) component of the tensor product.

The final tensor construction we make is that of a linear tensor derivation.
Given $A\in\End_{\real}(\alg{V})$\@, we define a derivation $D_A$ of the
tensor algebra $\oplus_{r,s\in\integernn}\tensor[r,s]{\alg{V}}$ by $D_A(a)=0$
for $a\in\tensor[0,0]{\alg{V}}\simeq\real$\@, and $D_A(v)=A(v)$ for
$v\in\tensor[1,0]{\alg{V}}\simeq\alg{V}$\@.  It then follows that
$D_A(\alpha)=-\dual{A}(\alpha)$ for $\alpha\in\dual{\alg{V}}$\@.  More
generally, we have the following result which expresses a well-known
formula~\cite[\eg][\S3.4]{EN:67} in terms of our insertion operation.
\begin{lemma}
Let\/ $\alg{V}$ be a finite-dimensional\/ $\real$-vector space, let\/
$A\in\End_{\real}(\alg{V})$\@, let\/ $r,s\in\integerp$\@, and let\/
$T\in\tensor[r,s]{\alg{V}}$\@.  Then
\begin{equation*}
D_A(T)=\sum_{j=1}^r\Ins_j(T,\dual{A})-\sum_{j=1}^s\Ins_{r+j}(T,A).
\end{equation*}
\begin{proof}
It suffices to take $T=v_1\otimes\dots\otimes v_r\otimes\alpha^1\otimes
\dots\otimes\alpha^s$\@.  In this case we have
{\allowdisplaybreaks\begin{align*}
D_A(T)(\beta^1,&\dots,\beta^r,u_1,\dots,u_s)\\
=&\;\sum_{j=1}^rv_1\otimes\dots\otimes A(v_j)\otimes\dots v_r\otimes
\alpha^1\otimes\dots\otimes\alpha^s(\beta^1,\dots,\beta^r,u_1,\dots,u_s)\\
&\;-\sum_{j=1}^sv_1\otimes\dots\otimes v_r\otimes\alpha^1\otimes\dots
\otimes\dual{A}(\alpha^j)\otimes\dots\otimes\alpha^s(\beta^1,\dots,\beta^r,
u_1,\dots,u_s)\\
=&\;\sum_{j=1}^r\beta^1(v_1)\cdots\beta^j(A(v_j))\cdots\beta^r(v_r)
\alpha^1(u_1)\cdots\alpha^s(u_s)\\
&\;-\sum_{j=1}^s\beta^1(v_1)\cdots\beta^r(v_r)\alpha^1(u_1)\cdots
\dual{A}(\alpha^j)(u_j)\cdots\alpha^s(u_s)\\
=&\;\sum_{j=1}^r\beta^1(v_1)\cdots\dual{A}(\beta^j)(v_j)\cdots\beta^r(v_r)
\alpha^1(u_1)\cdots\alpha^s(u_s)\\
&\;-\sum_{j=1}^s\beta^1(v_1)\cdots\beta^r(v_r)\alpha^1(u_1)\cdots
\alpha^j(A(u_j))\cdots\alpha^s(u_s)\\
=&\;\sum_{j=1}^rT(\beta^1,\dots,\dual{A}(\beta^j),\dots,\beta^r,u_1,\dots,u_s)\\
&\;-\sum_{j=1}^sT(\beta^1,\dots,\beta^r,u_1,\dots,A(u_j),\dots,u_s).
\end{align*}}%
This is exactly the claimed formula.
\end{proof}
\end{lemma}

We shall make a minor extension of the preceding notion of a derivation
associated to an endomorphism.  Let $k,r,s\in\integerp$\@.  Here we let
$T\in\tensor[r,s]{\alg{V}}$ and $S\in\tensor[1,k]{\alg{V}}$\@.  For
$v_1,\dots,v_{k-1}\in\alg{V}$\@, we define
$S_{(v_1,\dots,v_{k-1})}\in\End_{\real}(\alg{V})$ by
\begin{equation*}
S_{(v_1,\dots,v_{k-1})}(v)=S(v,v_1,\dots,v_{k-1}).
\end{equation*}
Denote
$\dual{S}\in\tensor[1,k-1]{\alg{V}}\otimes\dual{\alg{V}}$ by
\begin{equation*}
\natpair{\dual{S}(\beta,v_1,\dots,v_{k-1})}{v}=
\natpair{\beta}{S(v,v_1,\dots,v_{k-1})}
\end{equation*}
so that
\begin{equation*}
\dual{S}(\beta,v_1,\dots,v_{k-1})=\dual{S}_{(v_1,\dots,v_{k-1})}(\beta).
\end{equation*}
We then define $D_S(T)\in\tensor[r,s+k-1]{\alg{V}}$ by
\begin{equation}\label{eq:DS-derivation}
D_S(T)(\beta^1,\dots,\beta^r,u_1,\dots,u_{s+k-1})=
D_{S_{(u_{s+1}\dots,u_{s+k-1})}}(T)(\beta^1,\dots,\beta^r,u_1,\dots,u_s).
\end{equation}

The following elementary lemma gives a simpler formula for the previous
constructions.
\begin{lemma}\label{lem:derinsertion2}
Let\/ $\alg{V}$ be a finite-dimensional\/ $\real$-vector space, let\/
$k\in\integerp$\@, let\/ $S\in\tensor[1,k]{\alg{V}}$\@, let\/
$r,s\in\integerp$\@, and let\/ $T\in\tensor[r,s]{\alg{V}}$\@.  Then
\begin{equation*}
D_S(T)=\sum_{j=1}^r\Ins_j(T,\dual{S})-\sum_{j=1}^s\Ins_{r+j}(T,S).
\end{equation*}
\begin{proof}
We have
\begin{align*}
D_S(T)(\beta^1,\dots,&\beta^r,u_1,\dots,u_{k+s-1})\\
=&\;\sum_{j=1}^r\Ins_j(T,\dual{S_{(u_{s+1},\dots,u_{s+k-1})}})
(\beta^1,\dots,\beta^r,u_1,\dots,u_s)\\
&\;-\sum_{j=1}^s\Ins_{r+j}(T,S_{(u_{s+1},\dots,u_{s+k-1})})
(\beta^1,\dots,\beta^r,u_1,\dots,u_s)\\
=&\;\sum_{j=1}^rT(\beta^1,\dots,\dual{S_{(u_{s+1},\dots,u_{s+k-1})}}(\beta_j),
\dots,\beta^r,u_1,\dots,u_s)\\
&\;-\sum_{j=1}^sT(\beta^1,\dots,\beta^r,u_1,\dots,
S_{(u_{s+1},\dots,u_{s+k-1})}(u_j),\dots,u_s)\\
=&\;\sum_{j=1}^r\Ins_j(T,\dual{S})(\beta^1,\dots,\beta^r,u_1,\dots,u_{s+k-1})\\
&\;-\sum_{j=1}^r\Ins_{r+j}(T,S)(\beta^1,\dots,\beta^r,u_1,\dots,u_{s+k-1}),
\end{align*}
as claimed.
\end{proof}
\end{lemma}

Let us summarise this in the cases of interest.  The cases of interest will
be two in number.  The first is when $S\in\tensor[1,2]{\alg{V}}$ and
$T=T_0\otimes v$ for $T_0\in\tensor*[k]{\dual{\alg{V}}}$ and $v\in\alg{V}$\@.
In this case the preceding lemma gives
\begin{equation}\label{eq:simpleAder1}
\begin{aligned}
D_S(T)(v_1,\dots,&v_{k+1},\beta)\\
=&\;\Ins_{k+1}(T_0\otimes v,S^*)(v_1,\dots,v_{k+1},\beta)
-\sum_{j=1}^k\Ins_j(T_0\otimes v,S)(v_1,\dots,v_{k+1},\beta)\\
=&\;T_0(v_1,\dots,v_k)\natpair{\beta}{S_{v_{k+1}}(v)}
-\natpair{\beta}{v}\sum_{j=1}^k\Ins_j(T_0,S)(v_1,\dots,v_{k+1})\\
=&\;T_0(v_1,\dots,v_k)\natpair{\beta}{S(v,v_{k+1})}
-\natpair{\beta}{v}\sum_{j=1}^k\Ins_j(T_0,S)(v_1,\dots,v_{k+1}).
\end{aligned}
\end{equation}
The second case we will consider is when $S\in\tensor[1,2]{\alg{V}}$ and
$T=T_0\otimes\alpha$ for $T_0\in\tensor*[k]{\dual{\alg{V}}}$ and
$\alpha\in\dual{\alg{V}}$\@.  In this case we have
\begin{equation}\label{eq:simpleAder2}
\begin{aligned}
D_S(T)(v_1,\dots,&v_{k+2})\\
=&\;-\Ins_{k+1}(T_0\otimes\alpha,S)(v_1,\dots,v_{k+2})-
\sum_{j=1}^k\Ins_j(T_0\otimes\alpha,S)(v_1,\dots,v_{k+2})\\
=&\;-T_0(v_1,\dots,v_k)\alpha(S(v_{k+1},v_{k+2}))-
\natpair{\alpha}{v_{k+2}}\sum_{j=1}^k\Ins_j(T_0,S)(v_1,\dots,v_{k+1}).
\end{aligned}
\end{equation}

\section{Differentiation of tensors on the total space of a vector
bundle}\label{sec:tensor-derivatives}

In this section we establish some technical results for differentiation via
connections of various objects\textemdash{}functions, vector fields,
tensors\textemdash{}on vector bundles.  These results will allow us to
intrinsically perform the many calculations required to determine the
recursive relations given in Section~\ref{sec:lift-isomorphisms} between jets
on $\man{M}$ and jets on $\man{E}$ for a vector bundle
$\map{\pi_{\man{E}}}{\man{E}}{\man{M}}$\@.  As with the constructions of the
preceding section, the results in this section might seem \emph{non sequitur}
to the objectives of the paper.  And, as with the results of the preceding
section, perhaps a good strategy is to hurdle over this section until the
results are subsequently needed.

As with the material in Section~\ref{sec:tensor-constructions}\@, there is
nothing in this section that really separates the real analytic case from the
smooth case, so the presentation treats both cases on an equal footing.  What
is true, however, is that the complications of the computations in this
section are most useful in the real analytic setting of the paper.  If one
only wants to prove the continuity results in Section~\ref{sec:continuity} in
the smooth case, then simpler computations would suffice.

\subsection{Vector bundles as Riemannian
submersions}\label{subsec:Esubmersion}

In this section we let $\map{\pi_{\man{E}}}{\man{E}}{\man{M}}$ be a vector
bundle with $\map{\tbproj{\man{E}}}{\tb{\man{E}}}{\man{E}}$ its tangent
bundle.  We shall construct on $\man{E}$ a Riemannian metric in a more or
less natural way, using a Riemannian metric on $\man{M}$\@, an affine
connection on $\man{M}$\@, a fibre metric on $\man{E}$\@, and a vector bundle
connection in $\man{E}$\@.  For the initial part of the construction, we do
not require the affine connection on $\man{M}$ to be the Levi-Civita
connection, but we will only work with the case when it is, since there are
useful formulae one can prove in this case.  Let us indicate how one builds
the Riemannian metric on $\man{E}$\@.

Let $r\in\{\infty,\omega\}$\@.  We let
$\map{\pi_{\man{E}}}{\man{E}}{\man{M}}$ be a vector bundle of class $\C^r$
and suppose that $\nabla^{\pi_{\man{E}}}$ is a vector bundle connection on
$\man{E}$\@, $\metric_{\man{M}}$ is a Riemannian metric on $\man{M}$\@, and
$\metric_{\pi_{\man{E}}}$ is a fibre metric for $\man{E}$ with all data of
class $\C^r$\@.  The total space $\man{E}$ can be equipped with a Riemannian
metric via a natural adaptation of the Sasaki metric for tangent
bundles~\cite{SS:58}\@.  To define the inner product, we use the splitting
determined by the connection to give the inner product on $\tb[e]{\man{E}}$
by
\begin{equation}\label{eq:metricE}
\metric_{\man{E}}(w_1,w_2)=
\metric_{\man{M}}(\hor(w_1),\hor(w_2))+
\metric_{\pi_{\man{E}}}(\ver(w_1),\ver(w_2)).
\end{equation}
This then turns $\man{E}$ into a Riemannian manifold.  We denote by
$\nabla^{\man{E}}$ the Levi-Civita connection associated with
$\metric_{\man{E}}$\@.  Since the connection giving the splitting is of class
$\C^r$ if $\nabla^{\pi_{\man{E}}}$ is of class $\C^r$\@, the Riemannian
metric $\metric_{\man{E}}$ and its Levi-Civita connection are of class $\C^r$
if $\metric_{\man{M}}$ and $\metric_{\pi_{\man{E}}}$ are of class $\C^r$\@.

When we are working in this setting of Riemannian metrics and Levi-Civita
connections on the total space of a vector bundle
$\map{\pi_{\man{E}}}{\man{E}}{\man{M}}$\@, we shall denote by
$\metric_{\man{M}}$ the Riemannian metric on $\man{M}$ and by
$\nabla^{\man{M}}$ its Levi-Civita connection.

We note that the choice of metric $\metric_{\man{E}}$ ensures that
$\map{\pi_{\man{E}}}{\man{E}}{\man{M}}$ is a Riemannian submersion if we
equip $\man{M}$ with its Riemannian metric $\metric_{\man{M}}$ used to build
$\metric_{\man{E}}$\@.  Moreover, the fibres of $\pi_{\man{E}}$ are totally
geodesic submanifolds.  There are a few constructions involving Riemannian
submersions that will be helpful for us, and we review these here.  Let us
introduce some notation apropos to this.  We do this in a general setting.
Thus let $(\man{F},\metric_{\man{F}})$ and $(\man{M},\metric_{\man{M}})$ be
Riemannian manifolds.  Let $\map{\pi}{\man{F}}{\man{M}}$ be a
\defn{Riemannian submersion}\@,~\ie~for each $y\in\man{F}$\@,
\begin{equation*}
\metric_{\man{M}}(\tf[y]{\pi}(u),\tf[y]{\pi}(v))=\metric_{\man{F}}(u,v)
\end{equation*}
for every $u,v\in\tb[y]{\man{F}}$ that are orthogonal to
$\ker(\tf[y]{\pi})$\@.  We let $\vb{\man{F}}=\ker(\tf{\pi})$ be the vertical
subbundle with $\hb{\man{F}}$ its $\metric_{\man{F}}$-orthogonal complement,
which we call the horizontal subbundle.  We let
$\map{\ver,\hor}{\tb{\man{F}}}{\tb{\man{F}}}$ be the projections onto
$\vb{\man{F}}$ and $\hb{\man{F}}$\@, just as we have done for vector bundles.
For a vector field $X$ on $\man{M}$\@, we denote by $\hl{X}$ the horizontal
lift of $X$ to $\man{F}$\@.  This is the unique $\hb{\man{F}}$-valued vector
field satisfying $\tf[y]{\pi}(\hl{X}(y))=X\scirc\pi(y)$ for each
$y\in\man{F}$\@.

Given a submanifold of $\man{S}$ of a Riemannian manifold
$(\man{M},\metric_{\man{M}})$\@, $\man{S}$ inherits the Riemannian metric
$\metric_{\man{S}}$ obtained by pulling back $\metric_{\man{M}}$ by the
inclusion $\map{\iota_{\man{S}}}{\man{S}}{\man{M}}$\@.  The submanifold
$\man{S}$ is \defn{totally geodesic} if every geodesic for
$(\man{S},\metric_{\man{S}})$ is also a geodesic for
$(\man{M},\metric_{\man{M}})$\@.

Following~\cite{BO:68}\@, for a $\C^r$-Riemannian submersion $\map{\pi}{\man{F}}{\man{N}}$\@, there are two
associated tensors that characterise the submersion.  Specifically, we define
\begin{equation*}
A_\pi,T_\pi\in
\sections[r]{\tensor*[2]{\ctb{\man{F}}}\otimes\tb{\man{F}}}
\end{equation*}
by
\begin{align}\label{eq:Api}
A_\pi(\xi,\eta)=&\;\ver(\nabla^{\man{F}}_{\hor(\xi)}\hor(\eta))+
\hor(\nabla^{\man{F}}_{\hor(\xi)}\ver(\eta)),\\\notag
T_\pi(\xi,\eta)=&\;\hor(\nabla^{\man{F}}_{\ver(\xi)}\ver(\eta))+
\ver(\nabla^{\man{F}}_{\ver(\xi)}\hor(\eta))
\end{align}
for $\xi,\eta\in\sections[1]{\tb{\man{F}}}$\@.  One can easily verify that
$A_\pi$ and $T_\pi$ are indeed tensors as claimed.  Since the fibres of $\pi$
are submanifolds, we can define the \defn{vertical covariant derivative} as
the projection of the covariant derivative onto sections:
\begin{equation*}
\nabla^{\ver}_UV=\ver(\nabla^{\man{F}}_UV)
\end{equation*}
for vertical vector fields $U$ and $V$\@.

With all this background, we have the following result with tells us how to
covariantly differentiate vector fields on the total space of a vector
bundle.
\begin{lemma}\label{lem:riemannsub}
Let\/ $r\in\{\infty,\omega\}$\@.  Let\/ $(\man{F},\metric_{\man{F}})$ and\/
$(\man{M},\metric_{\man{M}})$ be\/ $\C^r$-Riemannian manifolds with\/
$\nabla^{\man{F}}$ and\/ $\nabla^{\man{M}}$ the Levi-Civita connections.
Let\/ $\map{\pi}{\man{F}}{\man{M}}$ be a Riemannian submersion.  Let\/
$X,Y\in\sections[r]{\tb{\man{M}}}$ and let\/
$U,V\in\sections[r]{\tb{\man{F}}}$ be vertical vector fields.  Then the
following statements hold:
\begin{compactenum}[(i)]
\item \label{pl:riemannsub1}
$\hor(\nabla^{\man{F}}_{\hl{X}}\hl{Y})=\hl{(\nabla^{\man{M}}_XY)}$\@;
\item \label{pl:riemannsub2}
$A_\pi(\hl{X},\hl{Y})=-\frac{1}{2}\ver([\hl{X},\hl{Y}])$\@;
\item \label{pl:riemannsub3}
$\nabla^{\man{F}}_UV=\nabla^{\ver}_UV+T_\pi(U,V)$\@;
\item \label{pl:riemannsub4}
$\nabla^{\man{F}}_V\hl{X}=\hor(\nabla^{\man{F}}_V\hl{X})+T_\pi(V,\hl{X})$\@;
\item \label{pl:riemannsub5}
$\nabla^{\man{F}}_{\hl{X}}V=\ver(\nabla^{\man{F}}_{\hl{X}}V)+A_\pi(\hl{X},V)$\@;
\item \label{pl:riemannsub6}
$\nabla^{\man{F}}_{\hl{X}}\hl{Y}=\hl{(\nabla^{\man{M}}_XY)}+
A_\pi(\hl{X},\hl{Y})$\@.
\item \label{pl:riemannsub7}
$\metric_{\man{F}}(\nabla^{\man{F}}_V\hl{X},\hl{Y})=
-\frac{1}{2}\metric_{\man{F}}(\ver([\hl{X},\hl{Y}]),V)=
\metric_{\man{F}}(\nabla^{\man{F}}_V\hl{Y},\hl{X})$\@.\savenum
\end{compactenum}
Additionally, if the fibres of\/ $\pi$ are totally geodesic submanifolds of\/
$\man{F}$\@, then the following statements hold:
\begin{compactenum}[(i)]\resumenum
\item \label{pl:riemannsub8} $T_\pi=0$\@;
\item \label{pl:riemannsub9} $\nabla^{\ver}|\man{F}_x$ is the Levi-Civita
connection for the submanifold Riemannian metric on\/ $\man{F}_x$\@;
\item \label{pl:riemannsub10} $\ver(\nabla^{\man{F}}_{\hl{X}}V)=\ver([\hl{X},V])$\@;
\item \label{pl:riemannsub11} $\nabla^{\man{F}}_V\hl{X}$ is horizontal and $\nabla^{\man{F}}_V\hl{X}=A_\pi(\hl{X},V)$\@.\savenum
\end{compactenum}
Finally, if\/ $\man{F}=\man{E}$ is the total space of a vector bundle and
if\/ $\metric_{\man{E}}$ is the Riemannian metric on\/ $\man{E}$ defined
above, then the following additional statements hold for sections\/
$\xi,\eta\in\sections[r]{\man{E}}$\@:
\begin{compactenum}[(i)]\resumenum
\item \label{pl:riemannsub12} $\nabla^{\man{E}}_{\vl{\xi}}\vl{\eta}=0$\@;
\item \label{pl:riemannsub13}
$\ver(\nabla^{\man{E}}_{\hl{X}}\vl{\xi})=\vl{(\nabla^\pi_X\xi)}$\@.
\end{compactenum}
\begin{proof}
We use the Koszul formula for the Levi-Civita connection:
\begin{multline}\label{eq:koszul}
2\metric_{\man{F}}(\nabla^{\man{F}}_\xi\eta,\zeta)=
\lieder{\xi}{(\metric_{\man{F}}(\eta,\zeta))}+
\lieder{\eta}{(\metric_{\man{F}}(\xi,\zeta))}-
\lieder{\zeta}{(\metric(\xi,\eta))}\\
+\metric_{\man{F}}([\xi,\eta],\zeta)-
\metric_{\man{F}}([\xi,\zeta],\eta)-
\metric_{\man{F}}([\eta,\zeta],\xi)
\end{multline}
for vector fields $\xi$\@, $\eta$\@, and $\zeta$ on
$\man{F}$~\cite[Page~160]{SK/KN:63a}\@.  We shall also use the formulae
\begin{equation}\label{eq:levi-civita}
\lieder{\zeta}{(\metric_{\man{F}}(\xi,\eta))}=
\metric_{\man{F}}(\nabla^{\man{F}}_\zeta\xi,\eta)+
\metric_{\man{F}}(\xi,\nabla^{\man{F}}_\zeta\eta)
\end{equation}
(saying that the Levi-Civita connection is a metric connection) and
\begin{equation}\label{eq:torsion-free}
\nabla^{\man{F}}_\xi\eta-\nabla^{\man{F}}_\eta\xi=[\xi,\eta]
\end{equation}
(saying that the Levi-Civita connection is torsion-free).  Both of these
formulae are determinable from the Koszul formula.

Let us make some preliminary computations.  First, since $\hl{X}$ and $\hl{Y}$ are
$\pi$-related to $X$ and $Y$\@, we have that $[\hl{X},\hl{Y}]$ is $\pi$-related to
$[X,Y]$~\cite[Proposition~4.2.25]{RA/JEM/TSR:88}\@.  Thus
\begin{equation}\label{eq:horbrack}
\hor([\hl{X},\hl{Y}])=\hl{[X,Y]}.
\end{equation}
In like manner, since $V$ is $\pi$-related to the zero vector field and $\hl{X}$
is $\pi$-related to $X$\@, $[V,\hl{X}]$ is $\pi$-related to the zero vector
field.  That is,
\begin{equation}\label{eq:verhorbrack}
\hor([V,\hl{X}])=0.
\end{equation}
Next, if $f$ is a function on $\man{M}$\@, then
\begin{equation*}
\lieder{\hl{X}}{(\pi^*f)}=\natpair{\d{(\pi^*f)}}{\hl{X}}=
\natpair{\pi^*\d{f}}{\hl{X}},
\end{equation*}
from which we deduce
\begin{equation}\label{eq:horderiv}
\lieder{\hl{X}}{(\pi^*f)}(y)=\natpair{\d{f}\scirc\pi(y)}{X\scirc\pi(y)},
\qquad y\in\man{F}.
\end{equation}
We trivially have
\begin{equation*}
\lieder{V}{(\pi^*f)}=0.
\end{equation*}

\eqref{pl:riemannsub1} One can use~\eqref{eq:koszul} with $\xi=\hl{X}$\@,
$\eta=\hl{Y}$\@, and $\zeta=\hl{Z}$\@, and the formulae~\eqref{eq:horbrack}
and~\eqref{eq:horderiv} to give
\begin{equation*}
\metric_{\man{F}}(\nabla^{\man{F}}_{\hl{X}}\hl{Y},\hl{Z})=
\pi^*\metric_{\man{M}}(\nabla^{\man{M}}_XY,Z).
\end{equation*}
This shows that
\begin{equation}\label{eq:riemannsub1}
\hor(\nabla^{\man{F}}_{\hl{X}}\hl{Y})=\hl{(\nabla^{\man{M}}_XY)}.
\end{equation}

\eqref{pl:riemannsub2} Now we use~\eqref{eq:koszul} with $\xi=\hl{X}$\@,
$\eta=\hl{Y}$\@, and $\zeta=V$\@.  We immediately have that the first three
terms on the right in~\eqref{eq:koszul} are zero.
By~\eqref{eq:verhorbrack}\@, the last two terms on the right
in~\eqref{eq:koszul} are zero.  Thus we have
\begin{equation*}
2\metric_{\man{F}}(\nabla^{\man{F}}_{\hl{X}}\hl{Y},V)=
\metric_{\man{F}}(\hl{[X,Y]},V),
\end{equation*}
and so
\begin{equation*}
A_\pi(\hl{X},\hl{Y})=\ver(\nabla^{\man{F}}_{\hl{X}}\hl{Y})=
\frac{1}{2}\ver(\hl{[X,Y]}).
\end{equation*}

\eqref{pl:riemannsub3} We have
\begin{equation*}
\nabla^{\man{F}}_UV=\ver(\nabla^{\man{F}}_UV)+\hor(\nabla^{\man{F}}_UV)=
\nabla^{\ver}_UV+T_\pi(U,V),
\end{equation*}
as claimed.

\eqref{pl:riemannsub4} We have
\begin{equation*}
\nabla^{\man{F}}_V\hl{X}=\hor(\nabla^{\man{F}}_V\hl{X})+
\ver(\nabla^{\man{F}}_V\hl{X})=\hor(\nabla^{\man{F}}_V\hl{X})+T_\pi(V,\hl{X}),
\end{equation*}
as claimed.

\eqref{pl:riemannsub5} We have
\begin{equation*}
\nabla^{\man{F}}_{\hl{X}}V=\ver(\nabla^{\man{F}}_{\hl{X}}V)+
\hor(\nabla^{\man{F}}_{\hl{X}}V)=
\ver(\nabla^{\man{F}}_{\hl{X}}V)+A_\pi(\hl{X},V),
\end{equation*}
as claimed.

\eqref{pl:riemannsub6} We have
\begin{equation*}
\nabla^{\man{F}}_{\hl{X}}\hl{Y}=\hor(\nabla^{\man{F}}_{\hl{X}}\hl{Y})+
\ver(\nabla^{\man{F}}_{\hl{X}}\hl{Y})=
\hl{(\nabla^{\man{M}}_XY)}+A_\pi(\hl{X},\hl{Y}),
\end{equation*}
using part~\eqref{pl:riemannsub1}\@.

\eqref{pl:riemannsub7} This is a direct computation
using~\eqref{eq:torsion-free}\@,~\eqref{eq:verhorbrack}\@,%
~\eqref{eq:levi-civita}\@, and part~\eqref{pl:riemannsub1}\@:
\begin{align}\notag
\metric_{\man{F}}(\nabla^{\man{F}}_V\hl{X},\hl{Y})=&\;
\metric_{\man{F}}(\nabla^{\man{F}}_{\hl{X}}V,\hl{Y})+
\metric_{\man{F}}([V,\hl{X}],\hl{Y})\\\label{eq:riemannsub2}
=&\;-\metric_{\man{F}}(\nabla^{\man{F}}_{\hl{X}}\hl{Y},V)\\\notag
=&\;-\frac{1}{2}\metric_{\man{F}}([\hl{X},\hl{Y}],V)=
-\frac{1}{2}\metric_{\man{F}}(\ver([\hl{X},\hl{Y}]),V).
\end{align}

\eqref{pl:riemannsub8} and~\eqref{pl:riemannsub9} These are properties of
totally geodesic submanifolds, so we first prove the result for the following
situation.
\begin{proofsublemma}
Let\/ $(\man{M},\metric_{\man{M}})$ be a Riemannian manifold and let\/
$\man{S}\subset\man{M}$ be a submanifold.  We let\/
$\metric_{\man{S}}=\iota_{\man{S}}^*\metric$ be the induced Riemannian metric
on\/ $\man{S}$\@.  We let\/ $\nabla^{\man{M}}$ and\/ $\nabla^{\man{S}}$ be
the Levi-Civita connections.  Then\/ $\man{S}$ is totally geodesic if and
only if\/ $\nabla^{\man{M}}_XY$ is tangent to\/ $\man{S}$ whenever\/
$X,Y\in\sections[1]{\tb{\man{M}}}$ are tangent to\/ $\man{S}$\@.
\begin{subproof}
We let $\nb{\man{S}}\subset\tb{\man{M}}|\man{S}$ be the normal bundle.  We
define the second fundamental form for $\man{S}$ to be the section
$\Pi_{\man{S}}$ of $\tensor*[2]{\tb{\man{S}}}\otimes\nb{\man{S}}$ defined by
\begin{equation*}
\Pi_{\man{S}}(X,Y)=\pr_{\nb{\man{S}}}(\nabla^{\man{M}}_XY)
\end{equation*}
for vector fields $X$ and $Y$ on $\man{M}$ that are tangent to $\man{S}$\@,
where $\map{\pr_{\nb{\man{S}}}}{\tb{\man{M}}|\man{S}}{\nb{\man{S}}}$ is the
orthogonal projection onto $\nb{\man{S}}$\@.

We claim that $\Pi_{\man{S}}$ is symmetric.  Indeed, by
\eqref{eq:torsion-free} we have
\begin{equation*}
\Pi_{\man{S}}(X,Y)-\Pi_{\man{S}}(Y,X)=\pr_{\nb{\man{S}}}([X,Y])=0
\end{equation*}
since $[X,Y]$ is tangent to $\man{S}$ if $X$ and $Y$ are tangent to
$\man{S}$\@.

Next we claim that
$\pr_{\tb{\man{S}}}(\nabla^{\man{M}}_XY)=\nabla^{\man{S}}_XY$ for vector
fields $X$ and $Y$ that are tangent to $\man{S}$\@, where
$\map{\pr_{\tb{\man{S}}}}{\tb{\man{M}}|\man{S}}{\tb{\man{S}}}$ is the
orthogonal projection.  To prove this, we show that
\begin{equation*}
(X,Y)\mapsto\pr_{\tb{\man{S}}}(\nabla^{\man{M}}_XY),
\end{equation*}
when restricted to $\man{S}$\@, satisfies the defining
conditions~\eqref{eq:levi-civita} and~\eqref{eq:torsion-free} for the
Levi-Civita connection for $\metric_{\man{S}}$\@.  Indeed, because $[X,Y]$ is
tangent to $\man{S}$ whenever $X$ and $Y$ are tangent to $\man{S}$\@, we
determine that, when restricted to $\man{S}$\@,
\begin{equation*}
\pr_{\tb{\man{S}}}(\nabla^{\man{M}}_XY-\nabla^{\man{M}}_YX)=
\pr_{\tb{\man{S}}}([X,Y])=[X,Y]
\end{equation*}
for all vector fields $X$ and $Y$ tangent to $\man{S}$\@.  This shows that
$(X,Y)\mapsto\pr_{\tb{\man{S}}}(\nabla^{\man{M}}_XY)$
satisfies~\eqref{eq:torsion-free}\@.  Also, when we restrict to $\man{S}$\@,
we have
\begin{align*}
\lieder{Z}{(\metric_{\man{S}}(X,Y))}=&\;\lieder{Z}{(\metric_{\man{M}}(X,Y))}
=\metric_{\man{M}}(\nabla^{\man{M}}_ZX,Y)+
\metric_{\man{M}}(X,\nabla^{\man{S}}_ZY)\\
=&\;\metric_{\man{S}}(\pr_{\tb{\man{S}}}(\nabla^{\man{M}}_ZX),Y)+
\metric_{\man{S}}(X,\pr_{\tb{\man{S}}}(\nabla^{\man{M}}_ZY))
\end{align*}
for all vector fields $X$\@, $Y$\@, and $Z$ that are tangent to $\man{S}$\@.
This shows that $(X,Y)\mapsto\pr_{\tb{\man{S}}}(\nabla^{\man{M}}_XY)$ satisfies~\eqref{eq:levi-civita}\@.

Now we can prove the sublemma.  First suppose that $\man{S}$ is totally
geodesic.  Let $v_x\in\tb{\man{S}}$ and let $t\mapsto\gamma(t)$ be a geodesic
for $\nabla^{\man{S}}$ satisfying $\gamma'(0)=v_x$\@.  Then $\gamma$ is also
a geodesic for $\nabla^{\man{M}}$\@.  Thus
\begin{align*}
0=&\;\nabla^{\man{M}}_{\gamma'(t)}\gamma'(t)=
\nabla^{\man{S}}_{\gamma'(t)}\gamma'(t)\\
=&\;\pr_{\tb{\man{S}}}(\nabla^{\man{M}}_{\gamma'(t)}\gamma'(t))\\
=&\;\pr_{\tb{\man{S}}}(\nabla^{\man{M}}_{\gamma'(t)}\gamma'(t))+
\pr_{\nb{\man{S}}}(\nabla^{\man{M}}_{\gamma'(t)}\gamma'(t)),
\end{align*}
from which we conclude, evaluating at $t=0$\@, that
$\Pi_{\man{S}}(v_x,v_x)=0$\@.  Since $\Pi_{\man{S}}$ is symmetric,
$\Pi_{\man{S}}=0$\@.  Thus
\begin{equation*}
\nabla^{\man{S}}_XY=\pr_{\tb{\man{S}}}(\nabla^{\man{M}}_XY)=\nabla^{\man{M}}_XY
\end{equation*}
for vector fields $X$ and $Y$ on $\man{M}$ tangent to $\man{S}$\@.  The
converse, that $\man{S}$ is totally geodesic if
$\nabla^{\man{M}}_XY=\nabla^{\man{S}}_XY$ for all vector fields $X$ and $Y$
on $\man{M}$ tangent to $\man{S}$\@, is clear.
\end{subproof}
\end{proofsublemma}

Given the sublemma, let $x\in\man{M}$ and let $\man{S}=\pi^{-1}(x)$ be the
fibre.  As we showed in the proof of the sublemma, if $U$ and $V$ are
vertical vector fields (in particular, they are tangent to $\man{S}$), then
\begin{equation*}
\nabla^{\man{F}}_UV=\ver(\nabla^{\man{F}}_UV)+T_\pi(U,V)=\nabla^{\man{S}}_UV.
\end{equation*}
Matching vertical and horizontal parts on $\man{S}$ gives
\begin{equation*}
\nabla^{\ver}_UV=\nabla^{\man{S}}_UV,\quad T_\pi(U,V)=0,
\end{equation*}
as claimed.

\eqref{pl:riemannsub11} It follows immediately from
parts~\eqref{pl:riemannsub4} and~\eqref{pl:riemannsub8} that
$\nabla^{\man{F}}_V\hl{X}$ is horizontal.  We also have
\begin{equation*}
\nabla^{\man{F}}_V\hl{X}=\hor(\nabla^{\man{F}}_V\hl{X})=
\hor(\nabla^{\man{F}}_{\hl{X}}V)+\hor([V,\hl{X}])
\end{equation*}
by~\eqref{eq:torsion-free}\@.  By part~\eqref{pl:riemannsub5}\@, the first
term on the far right is $A_\pi(\hl{X},V)$ and, by~\eqref{eq:verhorbrack}\@,
the second term in the far right is zero.

\eqref{pl:riemannsub10} By~\eqref{eq:torsion-free}\@, we have
\begin{equation*}
\ver(\nabla^{\man{F}}_{\hl{X}}V)=
\ver(\nabla^{\man{F}}_V\hl{X})+\ver([\hl{X},V]).
\end{equation*}
By part~\eqref{pl:riemannsub11} the first term on the right is zero.

\eqref{pl:riemannsub12} We note here that the fibres of
$\map{\pi_{\man{E}}}{\man{E}}{\man{M}}$ are vector spaces and the restriction
of $\metric_{\man{E}}$ to $\man{E}_x$ is just the constant Riemannian metric
$\metric_{\pi_{\man{E}}}(x)$\@.  Thus covariant derivatives on fibres are
just ordinary derivatives.  Now, since vertical lifts restricted to fibres
are constant, their ordinary derivatives are zero, and this gives the
assertion.

\eqref{pl:riemannsub13} Here, by part~\eqref{pl:riemannsub10}\@, we have
\begin{equation*}
\ver(\nabla^{\man{E}}_{\hl{X}}\vl{\xi})=\ver([\hl{X},\vl{\xi}]).
\end{equation*}
By~\eqref{eq:verhorbrack}\@, $[\hl{X},\vl{\xi}]$ is vertical.
By~\cite[Proposition~4.2.34]{RA/JEM/TSR:88}\@, we have
\begin{equation*}
[\vl{\xi},\hl{X}]=\frac{1}{2}\left.\deriv{^2}{t^2}\right|_{t=0}
\flow{\hl{X}}{-t}\scirc
\flow{\vl{\xi}}{-t}\scirc\flow{\hl{X}}{t}\scirc\flow{\vl{\xi}}{t}(e).
\end{equation*}
Using the fact that $\flow{\vl{\xi}}{t}(e)=e+t\xi\scirc\pi_{\man{E}}(e)$ and
that $\flow{\hl{X}}{t}(e)$ is the parallel transport $t\mapsto\tau^\gamma_t$
along integral curve $\gamma$ for $X$ through $\pi_{\man{E}}(e)$\@, we
directly calculate
\begin{equation*}
\flow{\hl{X}}{-t}\scirc\flow{\vl{\xi}}{-t}\scirc\flow{\hl{X}}{t}\scirc
\flow{\vl{\xi}}{t}(e)=
e-t(\tau^\gamma_t(\xi\scirc\gamma(t))-\xi\scirc\gamma(0)),
\end{equation*}
and from this, using the relationship between parallel
transport~\cite[page~114]{SK/KN:63a} and covariant derivative, we have
$[\vl{\xi},\hl{X}]=-\vl{(\nabla^{\pi_{\man{E}}}_X\xi)}$\@.
\end{proof}
\end{lemma}

\subsection{Derivatives of tensor contractions}

In Section~\ref{subsec:Ins} we constructed a tensor contraction/insertion
operator.  Let us consider the derivative of this operation.
\begin{lemma}\label{lem:insder1}
Let\/ $r\in\{\infty,\omega\}$\@.  Let\/ $\map{\pi_{\man{E}}}{\man{E}}{\man{M}}$ a
vector bundle of class\/ $\C^r$\@, let\/ $\nabla^{\pi_{\man{E}}}$ be a\/ $\C^r$-vector
bundle connection in\/ $\man{E}$\@, let\/ $k,l\in\integerp$\@, let\/
$A\in\sections[r]{\tensor*[k]{\dual{\man{E}}}}$\@, and let\/
$S\in\sections[r]{\tensor*[l]{\dual{\man{E}}}\otimes\man{E}}$\@.  For\/
$j\in\{1,\dots,k\}$ we have
\begin{equation*}
\nabla^{\pi_{\man{E}}}(\Ins_j(A,S))=\Ins_j(\nabla^{\pi_{\man{E}}}A,S)+
\Ins_j(A,\nabla^{\pi_{\man{E}}}S).
\end{equation*}
\begin{proof}
We let $\xi_a\in\sections[r]{\man{E}}$\@, $a\in\{1,\dots,k+l-1\}$\@, and
$X\in\sections[r]{\tb{\man{M}}}$\@.  We calculate
{\allowdisplaybreaks\begin{align*}
\lieder{X}{(\Ins_j&(A,S)(\xi_1,\dots,\xi_{k+l-1}))}=
(\nabla^{\pi_{\man{E}}}_X\Ins_j(A,S))(\xi_1,\dots,\xi_{k+l-1})\\
&\;+\sum_{a=1}^{k+l-1}\Ins_j(A,S)(\xi_1,\dots,\nabla^{\pi_{\man{E}}}_X\xi_a,
\dots,\xi_{k+l-1})\\
=&\;(\nabla^{\pi_{\man{E}}}_X\Ins_j(A,S))(\xi_1,\dots,\xi_{k+l-1})\\
&\;+\sum_{a=1}^{j-1}A(\xi_1,\dots,\nabla^{\pi_{\man{E}}}_X\xi_a,\dots,\xi_{j-1},
S(\xi_j,\xi_{k+1},\dots,\xi_{k+l-1}),\xi_{j+1},\dots,\xi_k)\\
&\;+A(\xi_1,\dots,\xi_{j-1},S(\nabla^{\pi_{\man{E}}}_X\xi_j,\xi_{k+1},
\dots,\xi_{k+l-1}),\xi_{j+1},\dots,\xi_k)\\
&\;+\sum_{a=j+1}^kA(\xi_1,\dots,\xi_{j-1},S(\xi_j,\xi_{k+1},
\dots,\xi_{k+l-1}),\xi_{j+1},\dots,\nabla^{\pi_{\man{E}}}_X\xi_a,\dots,\xi_k)\\
&\;+\sum_{a=k+1}^{k+l-1}A(\xi_1,\dots,\xi_{j-1},
S(\xi_j,\xi_{k+1},\dots,\nabla^{\pi_{\man{E}}}_X\xi_a,\dots,\xi_{k+l-1}),
\xi_{j+1},\dots,\xi_k).
\end{align*}}%
We also calculate
{\allowdisplaybreaks\begin{align*}
\lieder{X}{(\Ins_j&(A,S)(\xi_1,\dots,\xi_{k+l-1}))}\\
=&\;\lieder{X}{(A(\xi_1,\dots,\xi_{j-1},S(\xi_j,\xi_{k+1},
\dots,\xi_{k+l-1}),\xi_{j+1},\dots,\xi_k))}\\
=&\;(\nabla^{\pi_{\man{E}}}_XA)(\xi_1,\dots,\xi_{j-1},S(\xi_j,\xi_{k+1},
\dots,\xi_{k+l-1}),\xi_{j+1},\dots,\xi_k)\\
&\;+\sum_{a=1}^{j-1}A(\xi_1,\dots,\nabla^{\pi_{\man{E}}}_X\xi_a,\dots,\xi_{j-1},
S(\xi_j,\xi_{k+1},\dots,\xi_{k+l-1}),\xi_{j+1},\dots,\xi_k)\\
&\;+A(\xi_1,\dots,\xi_{j-1},(\nabla^{\pi_{\man{E}}}_XS)(\xi_j,\xi_{k+1},
\dots,\xi_{k+l-1}),\xi_{j+1},\dots,\xi_k)\\
&\;+A(\xi_1,\dots,\xi_{j-1},S(\nabla^{\pi_{\man{E}}}_X\xi_j,\xi_{k+1},
\dots,\xi_{k+l-1}),\xi_{j+1},\dots,\xi_k)\\
&\;+\sum_{a=j+1}^kA(\xi_1,\dots,\xi_{j-1},S(\xi_j,\xi_{k+1},
\dots,\xi_{k+l-1}),\xi_{j+1},\dots,\nabla^{\pi_{\man{E}}}_X\xi_a,\dots,\xi_k)\\
&\;+\sum_{a=k+1}^{k+l-1}A(\xi_1,\dots,\xi_{j-1},
S(\xi_j,\xi_{k+1},\dots,\nabla^{\pi_{\man{E}}}_X\xi_a,\dots,\xi_{k+l-1}),
\xi_{j+1},\dots,\xi_k).
\end{align*}}%
Comparing the right-hand sides of the preceding calculations gives
\begin{align*}
(\nabla^{\pi_{\man{E}}}\Ins_j&(A,S))(\xi_1,\dots,\xi_{k+l-1},X)\\
=&\;(\nabla^{\pi_{\man{E}}}A)(\xi_1,\dots,\xi_{j-1},S(\xi_j,\xi_{k+1},
\dots,\xi_{k+l-1}),\xi_{j+1},\dots,\xi_k,\xi_{k+l-1},X)\\
&\;+A(\xi_1,\dots,\xi_{j-1},(\nabla^{\pi_{\man{E}}}S)(\xi_j,\xi_{k+1},
\dots,\xi_{k+l-1},X),\xi_{j+1},\dots,\xi_k)\\
=&\;\Ins_j(\nabla^{\pi_{\man{E}}}A,S)(\xi_1,\dots,\xi_{k+l-1},X)
+\Ins_j(A,\nabla^{\pi_{\man{E}}}S)(\xi_1,\dots,\xi_{k+l-1},X),
\end{align*}
and this gives the result.
\end{proof}
\end{lemma}

Using this result, we can easily compute the derivative for tensor insertion
with one of the arguments fixed.
\begin{lemma}\label{lem:insder2}
Let\/ $r\in\{\infty,\omega\}$\@.  Let\/
$\map{\pi_{\man{E}}}{\man{E}}{\man{M}}$ a vector bundle of class\/ $\C^r$\@,
let\/ $\nabla^{\pi_{\man{E}}}$ be a\/ $\C^r$-vector bundle connection in\/
$\man{E}$\@, let\/ $l\in\integerp$\@, and let\/
$S\in\sections[r]{\tensor*[l]{\dual{\man{E}}}\otimes\man{E}}$\@.  Then, for\/
$k\in\integerp$ and\/ $j\in\{1,\dots,k\}$\@,
\begin{equation*}
(\nabla^{\pi_{\man{E}}}\Ins_{S,j})(A)=\Ins_j(A,\nabla^{\pi_{\man{E}}}S).
\end{equation*}
\begin{proof}
We have
\begin{equation*}
\nabla^{\pi_{\man{E}}}(\Ins_{S,j}(A))=(\nabla^{\pi_{\man{E}}}\Ins_{S,j})(A)+
\Ins_{S,j}(\nabla^{\pi_{\man{E}}}A)
\end{equation*}
and
\begin{equation*}
\nabla^{\pi_{\man{E}}}(\Ins_j(A,S))=\Ins_j(\nabla^{\pi_{\man{E}}}A,S)+
\Ins_j(A,\nabla^{\pi_{\man{E}}}S).
\end{equation*}
Comparing the equations, noting that
$\Ins_{S,j}(\nabla^{\pi_{\man{E}}}A)=\Ins_j(\nabla^{\pi_{\man{E}}}A,S)$\@, the result follows.
\end{proof}
\end{lemma}

Related to tensor contraction is the evaluation of a vector bundle mapping.
We shall consider the derivative of this evaluation.  In stating the result,
we use a bit of tensor notation that we now introduce.  Let $\alg{V}$ be a
finite-dimensional $\real$-vector space and let
$A\in\tensor[1,k+1]{\dual{\alg{V}}}$ and $B\in\tensor[1,l]{\alg{V}}$\@.  We
then denote by $A(B)\in\tensor*[k+l]{\dual{\alg{V}}}$ the tensor defined by
\begin{equation}\label{eq:A(B)ins}
A(B)(v_1,\dots,v_k,v_{k+1},\dots,v_{k+l})=
A(v_1,\dots,v_k,B(v_{k+1},\dots,v_{k+l}))
\end{equation}
Thus $A(B)$ is shorthand for $\Ins_{k+1}(A,B)$\@.  With this notation, we
have the following result.
\begin{lemma}\label{lem:leibniz}
Let\/ $\map{\pi_{\man{E}}}{\man{E}}{\man{M}}$ and\/
$\map{\pi_{\man{F}}}{\man{F}}{\man{M}}$ be
smooth vector bundles and let\/ $\nabla^{\pi_{\man{E}}}$ and\/ $\nabla^{\pi_{\man{F}}}$ be smooth
vector bundle connections in\/ $\man{E}$ and\/ $\man{F}$\@, respectively.
Let\/ $\nabla^{\man{M}}$ be an affine connection on\/ $\man{M}$\@.  Let\/
$L\in\sections[\infty]{\man{F}\otimes\dual{\man{E}}}$\@.  Then
\begin{equation*}
D^k_{\nabla^{\man{M}},\nabla^{\pi_{\man{F}}}}(L\scirc\xi)=
\sum_{l=0}^k\binom{k}{l}\Sym_k\left(
D^l_{\nabla^{\man{M}},\nabla^{\pi_{\man{F}}\otimes\pi_{\man{E}}}}(L)
(D^{k-l}_{\nabla^{\man{M}},\nabla^{\pi_{\man{E}}}}(\xi))\right),
\end{equation*}
for\/ $\xi\in\sections[\infty]{\man{E}}$\@.
\begin{proof}
First we claim that
{\small\begin{multline}\label{eq:leibniz1}
\nabla^{\man{M},\pi_{\man{F}},k}(L\scirc\xi)(X_1,\dots,X_k)\\
=\sum_{l=0}^k\sum_{\sigma\in\symmgroup{l,k-l}}
(\nabla^{\man{M},\pi_{\man{F}}\otimes\pi_{\man{E}},l}
L(X_{\sigma(1)},\dots,X_{\sigma(l)}))
(\nabla^{\man{M},\pi_{\man{E}},k-l}\xi(X_{\sigma(l+1)},\dots,X_{\sigma(k)})),
\end{multline}}%
for\/ $\xi\in\sections[k]{\man{E}}$\@.  This clearly holds for $k=1$\@.
So suppose it true for $k\ge1$ and compute
{\small\allowdisplaybreaks\begin{align*}
\nabla^{\pi_{\man{F}}}(&\nabla^{\man{M},\pi_{\man{F}},k}
(L\scirc\xi)(X_1,\dots,X_k))(X_{k+1})\\
=&\;\nabla^{\man{M},\pi_{\man{F}},k+1}(L\scirc\xi)(X_1,\dots,X_k,X_{k+1})+
\sum_{j=1}^k\nabla^{\man{M},\pi_{\man{F}},k}(L\scirc\xi)(X_1,\dots,
\nabla^{\man{M}}_{X_{k+1}}X_j,\dots,X_k)\\
=&\;\nabla^{\man{M},\pi_{\man{F}},k+1}(L\scirc\xi)(X_1,\dots,X_k,X_{k+1})\\
&\;+\sum_{j=1}^m\sum_{l=0}^{j-1}\sum_{\sigma\in\symmgroup{l,k-l}}
(\nabla^{\man{M},\pi_{\man{F}}\otimes\pi_{\man{E}},l}
L(X_{\sigma(1)},\dots,X_{\sigma(l)}))\\
&\;\phantom{+\sum_{j=1}^m\sum_{l=0}^{j-1}\sum_{\sigma\in\symmgroup{l,k-l}}}
(\nabla^{\man{M},\pi_{\man{E}},k-l}\xi(X_{\sigma(l+1)},\dots,
\nabla^{\man{M}}_{X_{k+1}}X_j,\dots,X_{\sigma(k)}))\\
&\;+\sum_{j=1}^k\sum_{l=j}^k\sum_{\sigma\in\symmgroup{l,k-l}}
(\nabla^{\man{M},\pi_{\man{F}}\otimes\pi_{\man{E}},l}L(X_{\sigma(1)},\dots,
\nabla^{\man{M}}_{X_{k+1}}X_j,\dots,X_{\sigma(l)}))\\
&\;\phantom{+\sum_{j=1}^k\sum_{l=j}^k\sum_{\sigma\in\symmgroup{l,k-l}}}
(\nabla^{\man{M},\pi_{\man{E}},k-l}\xi(X_{\sigma(l+1)},\dots,X_{\sigma(k)})),
\end{align*}}%
using the induction hypothesis.  We also compute
{\small\allowdisplaybreaks\begin{align*}
\nabla^{\pi_{\man{F}}}(&\nabla^{\man{M},\pi_{\man{F}},k}
(L\scirc\xi)(X_1,\dots,X_k))(X_{k+1})\\
=&\;\sum_{l=0}^k\sum_{\sigma\in\symmgroup{l,k-l}}
(\nabla^{\man{M},\pi_{\man{F}}\otimes\pi_{\man{E}},l+1}L)
(X_{\sigma(1)},\dots,X_{\sigma(l)},X_{k+1})\\
&\;\phantom{\sum_{l=0}^k\sum_{\sigma\in\symmgroup{l,k-l}}}
(\nabla^{\man{M},\pi_{\man{E}},k-l}\xi(X_{\sigma(l+1)},\dots,X_{\sigma(k)}))\\
&\;+\sum_{l=0}^k\sum_{\sigma\in\symmgroup{l,k-l}}\sum_{j=1}^l
(\nabla^{\pi_{\man{F}}\otimes\pi_{\man{E}},l}L(X_{\sigma(l)},\dots,
\nabla^{\man{M}}_{X_{k+1}}X_{\sigma(j)},\dots,X_{\sigma(l)})\\
&\;\phantom{+\sum_{l=0}^k\sum_{\sigma\in\symmgroup{l,k-l}}\sum_{j=1}^l}
(\nabla^{\man{M},\pi_{\man{E}},k-l}\xi(X_{\sigma(l+1)},\dots,X_{\sigma(k)}))\\
&\;+\sum_{l=0}^k\sum_{\sigma\in\symmgroup{l,k-l}}
(\nabla^{\man{M},\pi_{\man{F}}\otimes\pi_{\man{E}},l}
L(X_{\sigma(1)},\dots,X_{\sigma(l)}))\\
&\;\phantom{+\sum_{l=0}^k\sum_{\sigma\in\symmgroup{l,k-l}}}
(\nabla^{\man{M},\pi_{\man{E}},k-l+1}\xi(X_{\sigma(l+1)},\dots,
X_{\sigma(k)},X_{k+1}))\\
&\;+\sum_{l=0}^k\sum_{\sigma\in\symmgroup{l,k-l}}\sum_{j=l+1}^k
(\nabla^{\man{M},\pi_{\man{F}}\otimes\pi_{\man{E}},l}
L(X_{\sigma(1)},\dots,X_{\sigma(l)}))\\
&\;\phantom{+\sum_{l=0}^k\sum_{\sigma\in\symmgroup{l,k-l}}\sum_{j=l+1}^k}
(\nabla^{\man{M},\pi_{\man{E}},k-l+1}\xi(X_{\sigma(l+1)},\dots,
\nabla^{\man{M}}_{X_{k+1}}X_{\sigma(j)},\dots,X_{\sigma(k)})).
\end{align*}}%
Comparing the preceding two equations gives
{\small\begin{align*}
\nabla^{\man{M},\pi_{\man{F}},k+1}&(L\scirc\xi)(X_1,\dots,X_k,X_{k+1})\\
=&\;\sum_{l=0}^k\sum_{\sigma\in\symmgroup{l,k-l}}
(\nabla^{\man{M},\pi_{\man{F}}\otimes\pi_{\man{E}},l+1}L)
(X_{\sigma(1)},\dots,X_{\sigma(l)},X_{k+1})\\
&\;\phantom{\sum_{l=0}^k\sum_{\sigma\in\symmgroup{l,k-l}}}
(\nabla^{\man{M},\pi_{\man{E}},k-l}\xi(X_{\sigma(l+1)},\dots,X_{\sigma(k)}))\\
&\;+\sum_{l=0}^k\sum_{\sigma\in\symmgroup{l,k-l}}
(\nabla^{\man{M},\pi_{\man{F}}\otimes\pi_{\man{E}},l}
L(X_{\sigma(1)},\dots,X_{\sigma(l)}))\\
&\;\phantom{+\sum_{l=0}^k\sum_{\sigma\in\symmgroup{l,k-l}}}
(\nabla^{\man{M},\pi_{\man{E}},k-l+1}\xi(X_{\sigma(l+1)},\dots,
X_{\sigma(k)},X_{k+1}))\\
=&\;\sum_{l=0}^{k+1}\sum_{\sigma\in\symmgroup{l,k+1-l}}
(\nabla^{\man{M},\pi_{\man{F}}\otimes\pi_{\man{E}},l}
L(X_{\sigma(1)},\dots,X_{\sigma(l)}))\\
&\;\phantom{\sum_{l=0}^{k+1}\sum_{\sigma\in\symmgroup{l,k+1-l}}}
(\nabla^{\man{M},\pi_{\man{E}},k+1-l}\xi(X_{\sigma(l+1)},\dots,X_{\sigma(k+1)})),
\end{align*}}%
giving~\eqref{eq:leibniz1}\@.

For $A\in\tensor*[k]{\dual{\alg{V}}}$ and $\sigma\in\symmgroup{k}$\@, we use
the notation
\begin{equation*}
\sigma(A)(v_1,\dots,v_k)=A(v_{\sigma(1)},\dots,v_{\sigma(k)}).
\end{equation*}
For $\sigma\in\symmgroup{k}$\@, write $\sigma=\sigma_1\scirc\sigma_2$ for
$\sigma_1\in\symmgroup{k,l}$ and $\sigma_2\in\symmgroup{k|l}$\@.  Now we
compute
{\small\begin{align*}
D^k_{\nabla^{\man{M}},\nabla^{\pi_{\man{F}}}}(L\scirc\xi)=&\;
\frac{1}{k!}\sum_{\sigma\in\symmgroup{k}}\sigma(\nabla^{\pi_{\man{F}},k}
(L\scirc\xi))\\
=&\;\frac{1}{k!}\sum_{l=0}^k\sum_{\sigma\in\symmgroup{k}}
\sum_{\sigma'\in\symmgroup{l,k-l}}
\sigma'\scirc\sigma(\nabla^{\pi_{\man{F}}\otimes\pi_{\man{E}},l}
L(\nabla^{\pi_{\man{E}},k-l}\xi))\\
=&\;\frac{1}{k!}
\sum_{l=0}^k\sum_{\sigma'\in\symmgroup{l,k-l}}\sum_{\sigma_1\in\symmgroup{l,k-l}}
\sum_{\sigma_2\in\symmgroup{l|k-l}}\sigma'\scirc\sigma_1\scirc\sigma_2
(\nabla^{\pi_{\man{F}}\otimes\pi_{\man{E}},l}L(\nabla^{\pi_{\man{E}},k-l}\xi))\\
=&\;\sum_{l=0}^k\sum_{\sigma'\in\symmgroup{l,k-l}}
\sum_{\sigma_1\in\symmgroup{l,k-l}}
\frac{l!(k-l)!}{k!}\sigma'\scirc\sigma_1(
D^l_{\nabla^{\man{M}},\nabla^{\pi_{\man{F}}\otimes\pi_{\man{E}}}}L
(D^{k-l}_{\nabla^{\man{M}},\nabla^{\pi_{\man{E}}}}(\xi)))\\
=&\;\sum_{l=0}^k\sum_{\sigma'\in\symmgroup{l,k-l}}
\frac{l!(k-l)!}{k!}\sigma'\scirc\left(\sum_{\sigma\in\symmgroup{k}}
\frac{k!}{l!(k-l)!}
\Sym_k(D^l_{\nabla^{\man{M}},\nabla^{\pi_{\man{F}}\otimes\pi_{\man{E}}}}L
(D^{k-l}_{\nabla^{\man{M}},\nabla^{\pi_{\man{E}}}}(\xi)))\right)\\
=&\;\sum_{l=0}^k\frac{k!}{l!(k-l)!}\left(
\Sym_k(D^l_{\nabla^{\man{M}},\nabla^{\pi_{\man{F}}\otimes\pi_{\man{E}}}}L
(D^{k-l}_{\nabla^{\man{M}},\nabla^{\pi_{\man{E}}}}(\xi)))\right),
\end{align*}}%
making reference to~\eqref{eq:altodot} in the penultimate step, and noting
that $\card(\symmgroup{l,k-l})=\frac{k!}{l!(k-l)!}$\@.  This is the desired
result.
\end{proof}
\end{lemma}

\subsection{Derivatives of tensors on the total space of a vector
bundle}\label{subsec:formsvf}

In Definition~\ref{def:tensor-lifts} we gave definitions for a variety of
lifts of tensor fields.  Here we give formulae for differentiating these.  We
shall make ongoing and detailed use of the formulae we develop in this
section, and decent notation is an integral part of arriving at useable
expressions.

Let $r\in\{\infty,\omega\}$ and let $\map{\pi_{\man{E}}}{\man{E}}{\man{M}}$
be a $\C^r$-vector bundle.  We consider a $\C^r$-affine connection
$\nabla^{\man{M}}$ on $\man{M}$ and a $\C^r$-vector bundle connection
$\nabla^{\pi_{\man{E}}}$ in $\man{E}$\@.  The connection $\nabla^{\man{M}}$
induces a covariant derivative for tensor fields
$A\in\sections[r]{\tensor[k,l]{\tb{\man{M}}}}$ on $\man{M}$\@,
$k,l\in\integernn$\@.  This covariant derivative we denote by
$\nabla^{\man{M}}$\@, dropping the particular $k$ and $l$\@.  Similarly, the
connection $\nabla^{\pi_{\man{E}}}$ induces a covariant derivative for
sections $B\in\sections[\infty]{\tensor[k,l]{\man{E}}}$ of the tensor bundles
associated with $\man{E}$\@, $k,l\in\integernn$\@.  This covariant derivative
we denote by $\nabla^{\pi_{\man{E}}}$\@, again dropping the particular $k$
and $l$\@.  We have already made use of these conventions,~\eg~in
Lemmata~\ref{lem:insder1} and~\ref{lem:insder2}\@.  We will also consider
differentiation of sections of
$\tensor[k_1,l_1]{\tb{\man{M}}}\otimes\tensor[k_2,l_2]{\man{E}}$\@.  Here we
denote the covariant derivative by $\nabla^{\man{M},\pi_{\man{E}}}$\@.  If we
have another $\C^r$-vector bundle $\map{\pi_{\man{F}}}{\man{F}}{\man{M}}$
with a $\C^r$-affine connection $\nabla^{\pi_{\man{F}}}$\@, then
$\nabla^{\pi_{\man{E}}}$ and $\nabla^{\pi_{\man{F}}}$ induce a covariant
derivative in $\tensor[k_1,l_1]{\man{E}}\otimes\tensor[k_2,l_2]{\man{F}}$\@,
and we denote this covariant derivative by
$\nabla^{\pi_{\man{F}}\otimes\pi_{\man{E}}}$\@.

Another construction we need in this section concerns pull-back bundles.  Let
$r\in\{\infty,\omega\}$\@, let $\man{M}$ and $\man{N}$ be $\C^r$-manifolds,
let $\map{\pi_{\man{F}}}{\man{F}}{\man{M}}$ be a $\C^r$-vector bundle, and
let $\Phi\in\mappings[r]{\man{N}}{\man{M}}$\@.  We then have the pull-back
bundle $\map{\Phi^*\pi_{\man{F}}}{\Phi^*\man{F}}{\man{N}}$\@, which is a
vector bundle over $\man{N}$\@.  Given a section $\eta$ of $\man{F}$\@, we
have a section $\Phi^*\eta$ of $\Phi^*\man{F}$ defined by
$\Phi^*\eta(y)=(y,\eta\scirc\Phi(y))$\@.  Given a $\C^r$-vector bundle
connection $\nabla^{\pi_{\man{F}}}$ in $\man{F}$\@, we can define a
$\C^r$-connection $\Phi^*\nabla^{\pi_{\man{F}}}$ in $\Phi^*\man{F}$ by
requiring that
\begin{equation*}
\Phi^*\nabla^{\pi_{\man{F}}}_Z\Phi^*\eta(y)=
\Phi^*(\nabla^{\pi_{\man{F}}}_{\tf[y]{\Phi}(Z)}\eta)
\end{equation*}
for a $\C^\infty$-section $\eta$ and for $Z\in\tb[y]{\man{N}}$\@.  Given an
affine connection $\nabla^{\man{N}}$ on $\man{N}$\@, we then have an affine
connection on $\tensor*[k]{\ctb{\man{N}}}\otimes\Phi^*\man{F}$ induced by
tensor product by $\nabla^{\man{N}}$ and $\Phi^*\nabla^{\pi_{\man{F}}}$\@.
This connection we denote by $\nabla^{\man{N},\Phi^*\pi_{\man{F}}}$\@,
consistent with our notation above.  If we additionally have an injection
$\map{\psi}{\Phi^*\man{F}}{\tb{\man{N}}}$\@, then we have
\begin{equation*}
\nabla^{\man{N}}_Z(\psi\scirc\Phi^*\eta)=
\psi\scirc(\Phi^*\nabla^{\pi_{\man{F}}}_Z\Phi^*\eta)+B_\psi(\Phi^*\eta,Z)
\end{equation*}
for some tensor $B_\psi\in\tensor[1,2]{\tb{\man{E}}}$\@.

A special case of the preceding paragraph is when $\Phi=\pi_{\man{E}}$ for a
vector bundle $\map{\pi_{\man{E}}}{\man{E}}{\man{M}}$ and
$\man{F}=\man{E}$\@.  In this case, $\pi_{\man{E}}^*\xi=\vl{\xi}$ and
$\pi_{\man{E}}^*\man{F}\simeq\vb{\man{E}}$ and so we indeed have a natural
inclusion of $\pi_{\man{E}}^*\man{F}$ in $\tb{\man{E}}$\@.  Moreover, by
Lemma~\pldblref{lem:riemannsub}{pl:riemannsub13}\@,
\begin{equation*}
\pi_{\man{E}}^*\nabla^{\pi_{\man{E}}}_Z\pi_{\man{E}}^*\xi=
\vl{(\nabla^{\pi_{\man{E}}}_{\tf{\pi_{\man{E}}}(Z)}\xi)},
\end{equation*}
and so
\begin{equation}\label{eq:pi*xi=xiv}
\nabla^{\man{E}}_Z\pi_{\man{E}}^*\xi=
\pi_{\man{E}}^*\nabla^{\pi_{\man{E}}}_Z\pi_{\man{E}}^*\xi+
A_{\pi_{\man{E}}}(Z,\vl{\xi}).
\end{equation}

With the preceding, we can give formulae for differentiating tensors on
vector bundles, rather mirroring what we did in Lemma~\ref{lem:bundlefuncs}
for functions.
\begin{lemma}\label{lem:bundleform}
Let\/ $r\in\{\infty,\omega\}$\@.  Let\/
$\map{\pi_{\man{E}}}{\man{E}}{\man{M}}$ and\/
$\map{\pi_{\man{F}}}{\man{F}}{\man{M}}$ be vector bundles of class $\C^r$\@.
Let\/ $\metric_{\man{M}}$ be a\/ $\C^r$-Riemannian metric on\/ $\man{M}$\@,
let\/ $\nabla^{\man{M}}$ be the Levi-Civita connection, let\/
$\metric_{\pi_{\man{E}}}$ be a\/ $\C^r$-fibre metric on\/ $\man{E}$\@, and
let\/ $\nabla^{\pi_{\man{E}}}$ be a\/ $\metric_{\pi_{\man{E}}}$-vector bundle
connection of class\/ $\C^r$ in\/ $\man{E}$\@.  Let\/
$\nabla^{\pi_{\man{F}}}$ be a\/ $\C^r$-vector bundle connection in\/
$\man{F}$\@.  Let\/ $\metric_{\man{E}}$ be the associated\/ $\C^r$-Riemannian
metric on\/ $\man{E}$ from~\eqref{eq:metricE}\@.  Define
\begin{equation*}
B_{\pi_{\man{E}}}=\push_{1,2}\Ins_1(\Ins_2(A_{\pi_{\man{E}}},\hor),\hor)+
\Ins_2(A_{\pi_{\man{E}}},\ver)+\push_{1,2}\Ins_2(A_{\pi_{\man{E}}},\ver),
\end{equation*}
where\/ $A_{\pi_{\man{E}}}$ is defined as in~\eqref{eq:Api}\@.

Then we have the following statements, recalling from~\eqref{eq:DS-derivation}
the derivation\/ $D_{B_{\pi_{\man{E}}}}$\@:
\begin{compactenum}[(i)]
\item \label{pl:bundleform1} for\/
$k\in\integerp$ and\/ $A\in\sections[r]{\tensor*[k]{\ctb{\man{M}}}}$\@,\/
we have
\begin{equation*}
\nabla^{\man{E}}(\hl{A})=\hl{(\nabla^{\man{M}}A)}+D_{B_{\pi_{\man{E}}}}(\hl{A});
\end{equation*}

\item \label{pl:bundleform2} for\/
$A\in\sections[r]{\tensor*[k]{\ctb{\man{M}}}\otimes\man{E}}$\@, we have
\begin{equation*}
\nabla^{\man{E}}(\vl{A})=
\vl{(\nabla^{\man{M},\pi_{\man{E}}}A)}+D_{B_{\pi_{\man{E}}}}(\vl{A});
\end{equation*}

\item \label{pl:bundleform3} for\/
$A\in\sections[r]{\tensor*[k]{\ctb{\man{M}}}\otimes\tb{\man{M}}}$\@, we
have
\begin{equation*}
\nabla^{\man{E}}(\hl{A})=\hl{(\nabla^{\man{M}}A)}+D_{B_{\pi_{\man{E}}}}(\hl{A});
\end{equation*}

\item \label{pl:bundleform4} for\/
$A\in\sections[r]{\tensor*[k]{\ctb{\man{M}}}\otimes\man{F}\otimes
\dual{\man{E}}}$\@, we have
\begin{equation*}
\nabla^{\man{E},\pi_{\man{F}}}(\vl{A})=
\vl{(\nabla^{\man{M},\pi_{\man{E}}\otimes\pi_{\man{F}}}A)}+
D_{B_{\pi_{\man{E}}}}(\vl{A});
\end{equation*}

\item \label{pl:bundleform5} for\/
$A\in\sections[r]{\tensor*[k]{\ctb{\man{M}}}\otimes
\tensor[1,1]{\man{E}}}$\@, we have
\begin{equation*}
\nabla^{\man{E}}(\vl{A})=
\vl{(\nabla^{\man{M},\pi_{\man{E}}}A)}+D_{B_{\pi_{\man{E}}}}(\vl{A});
\end{equation*}

\item \label{pl:bundleform6} for\/
$A\in\sections[r]{\tensor*[k]{\ctb{\man{M}}}\otimes\man{F}\otimes
\dual{\man{E}}}$\@, we have
\begin{equation*}
\nabla^{\man{E},\pi_{\man{F}}}(\ve{A})=
\ve{(\nabla^{\man{M},\pi_{\man{E}}\otimes\pi_{\man{F}}}A)}+
D_{B_{\pi_{\man{E}}}}(\ve{A})+\vl{A}.
\end{equation*}

\item \label{pl:bundleform7} for\/
$A\in\sections[r]{\tensor*[k]{\ctb{\man{M}}}\otimes
\tensor[1,1]{\man{E}}}$\@, we have
\begin{equation*}
\nabla^{\man{E}}(\ve{A})=\ve{(\nabla^{\man{M},\pi_{\man{E}}}A)}+
D_{B_{\pi_{\man{E}}}}(\ve{A})+\vl{A}.
\end{equation*}
\end{compactenum}
\begin{proof}
Before we begin the proof proper, let us justify a ``without loss of
generality'' argument that we will make for the last four parts of the proof.
The arguments all have to do with assuming that it is sufficient, when
working with differential operators on spaces of tensor products, to work
with pure tensor products.  Let us be a little specific about this.  Let
$\map{\pi_{\man{E}}}{\man{E}}{\man{M}}$\@,
$\map{\pi_{\man{F}}}{\man{F}}{\man{M}}$\@, and
$\map{\pi_{\man{G}}}{\man{G}}{\man{M}}$ be $\C^r$-vector bundles.  Suppose
that $\map{\Delta_1,\Delta_2}{\jet{m}{(\man{E}\otimes\man{F})}}{\alg{G}}$ are
linear differential operators of order $m$\@.  We wish to give conditions
under which $\Delta_1=\Delta_2$\@.  Of course, this is equivalent to giving
conditions under which, for a differential operator
$\map{\Delta}{\jet{m}{(\man{E}\otimes\man{F})}}{\man{G}}$\@, $\Delta=0$\@.
To do so, we claim that, without loss of generality, we can simply prove that
$\Delta(j_m(\xi\otimes\eta))=0$ for all $\xi\in\sections[r]{\man{E}}$ and
$\eta\in\sections[r]{\man{E}}$\@.

To prove this sufficiency, we state and prove a couple of sublemmata.  The
second is the one of interest to us, and the first is used to prove the
second.  Simpler versions of the first lemma are called Hadamard's Lemma, but
we could not find a reference to the form we require.
\begin{proofsublemma}\label{psublem:hadamard}
Let\/ $r\in\{\infty,\omega\}$\@.  Let\/ $\nbhd{U}\subset\real^n$ be a
neighbourhood of\/ $\vect{0}$\@, let\/ $\alg{S}\subset\real^n$ be the
subspace
\begin{equation*}
\alg{S}=\setdef{(x^1,\dots,x^n)\in\real^n}{x^1=\dots=x^s=0},
\end{equation*}
let\/ $k\in\integernn$\@, and let\/
$f\in\func[r]{\oball{\epsilon}{\vect{0}}}$ satisfy\/
$\linder[j]{f}(\vect{x})=0$ for all\/ $j\in\{0,1,\dots,k\}$ and\/
$\vect{x}\in\alg{S}\cap\nbhd{U}$\@.  Let\/
$\map{\pr_{\alg{S}}}{\real^n}{\alg{S}}$ be the natural projection onto the
first\/ $s$-components.  Then there exists a neighbourhood\/
$\nbhd{V}\subset\nbhd{U}$ of\/ $\vect{0}$\@,\/
$g_I\in\func[r]{\nbhd{V}}$\@,\/ $I\in\integerp^s$\@,\/ $\snorm{I}=k+1$\@,
such that
\begin{equation*}
f(\vect{x})=\sum_{\substack{I\in\integerp^s\\\snorm{I}=k+1}}^n
g_I(\vect{x})\pr_{\alg{S}}(\vect{x})^I,\qquad\vect{x}\in\nbhd{V}.
\end{equation*}
\begin{subproof}
We prove the sublemma by induction on $k$\@.  For $k=0$\@, the hypothesis is
that $f$ vanishes on $\man{S}\cap\nbhd{U}$\@.  Let $\nbhd{W}\subset\alg{S}$
be a neighbourhood of $\vect{0}$ and let $\epsilon\in\realp$ be such that
$\oball{\epsilon}{\vect{y}}\subset\nbhd{U}$ for all $\vect{x}\in\nbhd{W}$\@,
possibly after shrinking $\nbhd{W}$\@.  Let
\begin{equation*}
\nbhd{V}=\bigcup_{\vect{x}\in\nbhd{W}}\oball{\epsilon}{\vect{x}}.
\end{equation*}
Let $\vect{x}=(\vect{x}_1,\vect{x}_2)\in\nbhd{V}$ (with
$\vect{x}_1\in\alg{S}$) and define
\begin{equation*}
\mapdef{\gamma_{\vect{x}}}{\interval[0,1]}{\real}{t}
{f(\vect{x}_1,t\vect{x}_2).}
\end{equation*}
We calculate
\begin{align*}
f(\vect{x})=&\;f(\vect{x}_1,\vect{x}_2)=f(\vect{x}_1,\vect{x}_2)-
f(\vect{x}_1,\vect{0})\\
=&\;\gamma_{\vect{x}}(1)-\gamma_{\vect{x}}(0)=
\int_0^1\gamma'_{\vect{x}}(t)\,\d{t}\\
=&\;\sum_{j=1}^sx^j\pderiv{f}{x^j}((\vect{x}_1,t\vect{x}_2))\,\d{t}=
\sum_{j=1}^sx^jg_j(\vect{x}),
\end{align*}
where
\begin{equation*}
g_j(\vect{x})=\int_0^1\pderiv{f}{x^j}(\vect{x}_1,t\vect{x}_2)\,\d{t},
\qquad j\in\{1,\dots,s\}.
\end{equation*}
It remains to show that the functions $g_1,\dots,g_s$ are of class $\C^r$\@.
By standard theorems on interchanging derivatives and
integrals~\cite[Theorem~16.11]{JJ:05}\@, we can conclude that $g_1,\dots,g_m$
are smooth when $f$ is smooth.  If the data are holomorphic, swapping
integrals and derivatives allows us to conclude that $g_1,\dots,g_s$ are
holomorphic when $f$ is holomorphic, by verifying the
Cauchy\textendash{}Riemann equations.  In the real analytic case, we can
complexify to a complex neighbourhood of $\vect{0}$\@, and so conclude real
analyticity by holomorphicity of the complexification.

As a standin for a full proof by induction, let us see how the case $k=1$
follows from the case $k=0$\@.  The general inductive argument is the same,
only with more notation.

We note that, for $\vect{x}\in\nbhd{V}$\@, we have
\begin{equation*}
\pderiv{f}{x^k}(\vect{x})=\begin{cases}g_k(\vect{x})+
\sum_{j=1}^sx^j\pderiv{g_j}{x^k}(\vect{x}),&k\in\{1,\dots,s\},\\
\sum_{j=1}^sx^j\pderiv{g_j}{x^k}(\vect{x}),&k\in\{s+1,\dots,n\}.\end{cases}
\end{equation*}
Thus $\linder{f}(\vect{x})=0$ for $\vect{x}\in\nbhd{S}\cap\nbhd{V}$ if and
only if $g_1(\vect{x})=\dots=g_s(\vect{x})=0$\@.  Thus one can apply the
arguments from the first part of the proof to write
\begin{equation*}
g_k(\vect{x})=\sum_{j=1}^sx^jg_{kj}(\vect{x})
\end{equation*}
on a neighbourhood of $\vect{0}$\@.  Thus
\begin{equation*}
f(\vect{x})=\sum_{j,k=1}^sx^kx^jg_{kj}(\vect{x}),
\end{equation*}
giving the desired form of $f$ in this case.
\end{subproof}
\end{proofsublemma}

\begin{proofsublemma}\label{psublem:Sextend}
Let\/ $r\in\{\infty,\omega\}$\@, let\/
$\map{\pi_{\man{E}}}{\man{E}}{\man{M}}$ be a\/ $\C^r$-vector bundle, and
let\/ $\man{S}\subset\man{M}$ be a closed\/ $\C^r$-submanifold.  Let\/
$k\in\integernn$\@.  Let\/ $\nbhd{V}$ be a neighbourhood of\/ $\man{S}$ and
let\/ $\xi_{\man{S}}\in\sections[r]{\man{E}|\nbhd{V}}$\@.  Then there
exists\/ $\xi\in\sections[r]{\man{E}}$ such that\/
$j_k\xi(x)=j_k\xi_{\man{S}}(x)$\@.
\begin{subproof}
Let $\sZ^k_{\man{S}}$ be the sheaf of $\C^r$-sections of $\man{E}$ such that
whose $k$-jet vanishes on $\man{S}$\@.  We have the exact sequence
\begin{equation*}
\xymatrix{{0}\ar[r]&{\sZ^k_{\man{S}}}\ar[r]&{\ssections[r]{\man{E}}}
\ar[r]^(0.45){\Psi}&
{\ssections[r]{\man{E}}/\sZ^k_{\man{S}}}\ar[r]&{0}}
\end{equation*}
Note that the stalk of the quotient sheaf at $x\in\man{S}$ consists of germs
of sections whose $k$-jets agree on $\man{S}$\@.

Now, if $x\not\in\man{S}$\@, then there is a neighbourhood $\nbhd{U}$ of $x$
such that $\nbhd{V}\cap\man{S}=\emptyset$\@, and so
$\sZ^k_{\man{S}}(\nbhd{U})=\ssections[r]{\man{E}}(\nbhd{U})$\@.  That
$\sZ^k_{\man{S}}$ is locally finitely generated at $x$ then follows since
$\ssections[r]{\man{E}}$ is locally finitely generated.  If $x\in\man{S}$\@,
choose a submanifold chart $(\nbhd{U},\phi)$ for $\man{S}$ about $x$ so that
\begin{equation*}
\man{S}\cap\nbhd{U}=\setdef{y\in\nbhd{U}}{\phi(y)=
(0,\dots,0,x^{s+1},\dots,x^n)}.
\end{equation*}
Then the $k$-jet of a function $f$ on $\nbhd{U}$ vanishes on $\man{S}$ if and
only if it is a $\func[r]{\nbhd{U}}$-linear combination of polynomial
functions in $x^1,\dots,x^s$ of degree $k+1$\@; this follows by the previous
sublemma.  Thus, if $\xi_1,\dots,\xi_m$ is a local basis of sections of
$\man{E}$ about $x$\@, then the (finite) set of products of these sections
with the polynomial functions in $x^1,\dots,x^s$ of degree at least $k+1$
generates $\sections[r]{\man{E}|\nbhd{U}}$ as a $\func[r]{\nbhd{U}}$-module.
This shows that $\sZ^k_{\man{S}}$ is locally finitely generated about $x$\@.
This shows that $\sZ^k_{\man{S}}$ is coherent in the case $r=\omega$\@.

Cartan's Theorem~B~\cite[Proposition~6]{HC:57} shows, in the case
$r=\omega$\@, that $\Psi$ is surjective on global sections.  The case of
$r=\infty$ follows in a similar way, using the fact that positive cohomology
groups for sheave of modules of smooth functions vanish
(\cite[Proposition~3.11]{ROW:08}\@, along with
\cite[Examples~3.4(d,e)]{ROW:08} and~\cite[Proposition~3.5]{ROW:08}).  This
implies that there exists $\xi\in\sections[r]{\man{E}}$ such that, for each
$x\in\man{S}$\@, $[\xi]_x=[\xi_{\man{S}}]_x$\@.  This, however, means
precisely that $j_k\xi(x)=j_k\xi_{\man{S}}(x)$ for each $x\in\man{S}$\@.
\end{subproof}
\end{proofsublemma}

Now suppose that we have proved that $\Delta(j_m(\xi\otimes\eta))=0$ for all
$\xi\in\sections[r]{\man{E}}$ and $\eta\in\sections[r]{\man{E}}$\@.  Let
$x\in\man{X}$ and let $\alpha\in\jetalg[x]{m}{\man{M}}$\@.  Let
$u\in\man{E}_x$ and $v\in\man{F}_x$\@.  By the previous sublemma, there
exists $f\in\func[r]{\man{M}}$ such that $j_mf(x)=\alpha$\@.  Then, keeping
in mind the identification~\eqref{eq:jet=jetalg}\@,
\begin{equation*}
\Delta(\alpha\otimes(u\otimes v))=
\Delta(j_m(f(\xi\otimes\eta)))=\Delta(j_m((f\xi)\otimes\eta))=0.
\end{equation*}
Since every element of $\jet[x]{m}{\man{E}}$ is a finite linear combination
of terms of the form $\alpha\otimes(u\otimes v)$ for
$\alpha\in\jetalg[x]{m}{\man{M}}$\@, $u\in\man{E}_x$\@, and
$v\in\man{F}_x$\@, we conclude that $\Delta(j_mA)(x)=0$ for every $A\in\sections[r]{\man{E}\otimes\man{F}}$\@.

Now we proceed with the proof.

\eqref{pl:bundleform1} We have
\begin{equation*}
\lieder{Z_{k+1}}{(\hl{A}(Z_1,\dots,Z_k))}=
(\nabla^{\man{E}}_{Z_{k+1}}\hl{A})(Z_1,\dots,Z_k)
+\sum_{j=1}^k\hl{A}(Z_1,\dots,\nabla^{\man{E}}_{Z_{k+1}}Z_j,\dots,Z_k).
\end{equation*}
We consider four cases.
\begin{compactenum}
\item $Z_j=\hl{X_j}$\@, $j\in\{1,\dots,k+1\}$\@: Here we have
\begin{equation*}
\lieder{\hl{X_{k+1}}}{(\hl{A}(\hl{X_1},\dots,\hl{X_k}))}=
\hl{(\lieder{X_{k+1}}{(A(X_1,\dots,X_k))})}
\end{equation*}
(by Lemma~\pldblref{lem:bundlefuncs}{pl:bundlefuncs1}) and
\begin{equation*}
\hl{A}(\hl{X_1},\dots,\nabla^{\man{E}}_{\hl{X_{k+1}}}\hl{X_j},\dots,\hl{X_k})=
\hl{(A(X_1,\dots,\nabla^{\man{M}}_{X_{k+1}}X_j,\dots,X_k))}
\end{equation*}
(by Lemma~\pldblref{lem:riemannsub}{pl:riemannsub1}).  Thus we conclude that
\begin{equation*}
\nabla^{\man{E}}\hl{A}(\hl{X_1},\dots,\hl{X_{k+1}})=
\hl{((\nabla^{\man{M}}A)(X_1,\dots,X_{k+1}))}.
\end{equation*}

\item $Z_j=\hl{X_j}$\@, $j\in\{1,\dots,k\}$\@, $Z_{k+1}=\vl{\xi_{k+1}}$\@: Here
we calculate
\begin{equation*}
\lieder{\vl{\xi_{k+1}}}{(\hl{A}(\hl{X_1},\dots,\hl{X_k}))}=
\lieder{\vl{\xi_{k+1}}}{\hl{(A(X_1,\dots,X_k))}}=0
\end{equation*}
(using the definition of $\hl{A}$ and
Lemma~\pldblref{lem:bundlefuncs}{pl:bundlefuncs2}) and
\begin{equation*}
\hl{A}(\hl{X_1},\dots,\nabla^{\man{E}}_{\vl{\xi_{k+1}}}\hl{X_j},\dots,\hl{X_k})=
\hl{A}(\hl{X_1},\dots,A_{\pi_{\man{E}}}(\hl{X_j},\vl{\xi_{k+1}}),\dots,\hl{X_k})
\end{equation*}
(using Lemma~\pldblref{lem:riemannsub}{pl:riemannsub11}).  Thus we conclude
that
\begin{equation*}
\nabla^{\man{E}}\hl{A}(\hl{X_1},\dots,\hl{X_k},\vl{\xi_{k+1}})=
-\sum_{j=1}^k\hl{A}(\hl{X_1},\dots,A_{\pi_{\man{E}}}(\hl{X_j},
\vl{\xi_{k+1}}),\dots,\hl{X_k}).
\end{equation*}

\item $Z_j=\vl{\xi_j}$ for some $j\in\{1,\dots,k\}$\@,
$Z_{k+1}=\hl{X_{k+1}}$\@: We calculate
\begin{equation*}
\lieder{\hl{X_{k+1}}}{(\hl{A}(Z_1,\dots,\vl{\xi_j},Z_k))}=0
\end{equation*}
(by definition of $\hl{A}$) and
\begin{equation*}
\hl{A}(Z_1,\dots,\nabla^{\man{E}}_{\hl{X_{k+1}}}\vl{\xi_j},\dots,Z_k)=
\hl{A}(Z_1,\dots,A_{\pi_{\man{E}}}(\hl{X_{k+1}},\vl{\xi_j}),\dots,Z_k)
\end{equation*}
(by Lemma~\pldblref{lem:riemannsub}{pl:riemannsub5}).  Thus
\begin{equation*}
\nabla^{\man{E}}\hl{A}(Z_1,\dots,\vl{\xi_j},\dots,Z_k,\hl{X_{k+1}})=
-\hl{A}(Z_1,\dots,A_{\pi_{\man{E}}}(\hl{X_{k+1}},\vl{\xi_j}),\dots,Z_k).
\end{equation*}

\item $Z_j=\vl{\xi_j}$ for some $j\in\{1,\dots,k\}$\@,
$Z_{k+1}=\vl{\xi_{k+1}}$\@: We have
\begin{equation*}
\lieder{\vl{\xi_{k+1}}}{(\hl{A}(Z_1,\dots,\vl{\xi_j},\dots,Z_k))}=0
\end{equation*}
(by definition of $\hl{A}$) and
\begin{equation*}
\hl{A}(Z_1,\dots,\nabla^{\man{E}}_{\vl{\xi_{k+1}}}\vl{\xi_j},\dots,Z_k)=0
\end{equation*}
(by Lemma~\pldblref{lem:riemannsub}{pl:riemannsub3}).  Thus
\begin{equation*}
\nabla^{\man{E}}\hl{A}(Z_1,\dots,\vl{\xi_j},\dots,Z_k,\vl{\xi_{k+1}})=0.
\end{equation*}
\end{compactenum}
Putting this all together, and keeping in mind that $A_{\pi_{\man{E}}}$ is vertical when
both arguments are vertical, we have
\begin{multline*}
\nabla^{\man{E}}\hl{A}(Z_1,\dots,Z_{k+1})=
\hl{(\nabla^{\man{M}}A)}(Z_1,\dots,Z_{k+1})\\
-\sum_{j=1}^k\hl{A}(Z_1,\dots,A_{\pi_{\man{E}}}(\hor(Z_j),\ver(Z_{k+1})),
\dots,Z_k)\\
-\sum_{j=1}^k\hl{A}(Z_1,\dots,A_{\pi_{\man{E}}}(\hor(Z_{k+1}),\ver(Z_j)),
\dots,Z_k).
\end{multline*}
Now we note that
\begin{align*}
B_{\pi_{\man{E}}}(Z_j,Z_{k+1})=&\;A_{\pi_{\man{E}}}(\hor(Z_{k+1}),\hor(Z_j))+
A_{\pi_{\man{E}}}(Z_j,\ver(Z_{k+1}))+A_{\pi_{\man{E}}}(Z_{k+1},\ver(Z_j))\\
=&\;A_{\pi_{\man{E}}}(\hor(Z_j),\ver(Z_{k+1}))+
A_{\pi_{\man{E}}}(\hor(Z_{k+1}),\ver(Z_j))+\textrm{something vertical},
\end{align*}
using Lemma~\pldblref{lem:riemannsub}{pl:riemannsub2} and the definition of
$A_{\pi_{\man{E}}}$\@.  Thus
\begin{equation*}
\nabla^{\man{E}}\hl{A}=\hl{(\nabla^{\man{M}}A)}-
\sum_{j=1}^k\Ins_j(\hl{A},B_{\pi_{\man{E}}}),
\end{equation*}
which gives this part of the lemma by Lemma~\ref{lem:derinsertion2}\@.

\eqref{pl:bundleform2} First we compute, for
$Z\in\sections[r]{\tb{\man{E}}}$\@,
\begin{align*}
\nabla^{\man{E}}_Z\vl{\xi}=&\;\nabla^{\man{E}}_{\hor(Z)}\vl{\xi}+
\nabla^{\man{E}}_{\ver(Z)}\vl{\xi}=
\vl{(\nabla^{\man{M},\pi_{\man{E}}}_{\tf{\pi_{\man{E}}}(Z)}\xi)}+
A_{\pi_{\man{E}}}(\tf{\pi_{\man{E}}}(Z),\vl{\xi})\\
=&\;\vl{(\nabla^{\man{M},\pi_{\man{E}}}_{\tf{\pi_{\man{E}}}(Z)}\xi)}+
A_{\pi_{\man{E}}}(Z,\vl{\xi})
\end{align*}
using Lemma~\pldblref{lem:riemannsub}{pl:riemannsub3}\@,~\eqref{pl:riemannsub5}\@,
and~\eqref{pl:riemannsub13}\@, and the definition of $A_{\pi_{\man{E}}}$\@.
If we note that
\begin{equation*}
B_{\pi_{\man{E}}}(\vl{\xi},Z)=A_{\pi_{\man{E}}}(\hor(Z),\hor(\vl{\xi}))+
A_{\pi_{\man{E}}}(\vl{\xi},\ver(Z))+
A_{\pi_{\man{E}}}(Z,\ver(\vl{\xi}))=A_{\pi_{\man{E}}}(Z,\vl{\xi})
\end{equation*}
(using the definition of $A_{\pi_{\man{E}}}$), we have
\begin{equation*}
\nabla^{\man{E}}_Z\vl{\xi}=
\vl{(\nabla^{\man{M},\pi_{\man{E}}}_{\tf{\pi_{\man{E}}}(Z)}\xi)}+
B_{\pi_{\man{E}}}(\vl{\xi},Z).
\end{equation*}

Now, it suffices to prove this part of the lemma for
$A=\hl{A_0}\otimes\vl{\xi}$ for
$A_0\in\sections[r]{\tensor*[k]{\ctb{\man{M}}}}$ and
$\xi\in\sections[r]{\man{E}}$\@.  For $Z\in\sections[r]{\tb{\man{E}}}$\@, we
have
\begin{align*}
\nabla^{\man{E}}_Z(\hl{A_0}\otimes\vl{\xi})=&\;
(\nabla^{\man{E}}_Z\hl{A_0})\otimes\vl{\xi}+(\hl{A_0})\otimes
\nabla^{\man{E}}_Z\vl{\xi}\\
=&\;\hl{(\nabla^{\man{M}}_{\tf{\pi_{\man{E}}}(Z)}A_0)}\otimes\vl{\xi}-
\sum_{j=1}^k\Ins_j(\hl{A_0},B_{\pi_{\man{E}},Z})\otimes\vl{\xi}\\
&\;+\hl{A_0}\otimes\vl{(\nabla^{\pi_{\man{E}}}_{\tf{\pi_{\man{E}}}(Z)}\xi)}+
\hl{A_0}\otimes B_{\pi_{\man{E}}}(\vl{\xi},Z).
\end{align*}
We have
\begin{align*}
\hl{(\nabla^{\man{M}}_{\tf{\pi_{\man{E}}}(Z)}A_0)}\otimes\vl{\xi}+
\hl{A_0}\otimes\vl{(\nabla^{\pi_{\man{E}}}_{\tf{\pi_{\man{E}}}(Z)}\xi)}=&\;
\vl{((\nabla^{\man{M}}_{\tf{\pi_{\man{E}}(Z)}}A_0)\otimes\xi)}+
(\hl{A_0}\otimes\vl{(\nabla^{\pi_{\man{E}}}_{\tf{\pi_{\man{E}}(Z)}}\xi))}\\
=&\;\vl{(\nabla^{\man{M},\pi_{\man{E}}}_{\tf{\pi_{\man{E}}}(Z)}(A_0\otimes\xi))}.
\end{align*}
Thus, by~\eqref{eq:simpleAder1} and the first part of the lemma, we have
\begin{equation*}
\hl{A_0}\otimes
B_{\pi_{\man{E}}}(\vl{\xi},Z)-\sum_{j=1}^k\Ins_j(\hl{A_0},B_{\pi_{\man{E}},Z})
\otimes\vl{\xi}=D_{B_{\pi_{\man{E}},Z}}(\hl{A_0}\otimes\vl{\xi}).
\end{equation*}
Assembling the preceding three computations gives this part of the lemma.

\eqref{pl:bundleform3} First note that
\begin{equation*}
\nabla^{\man{E}}_Z\hl{X}=\nabla^{\man{E}}_{\hor(Z)}\hl{X}+
\nabla^{\man{E}}_{\ver(Z)}\hl{X}=\hl{(\nabla^{\man{M}}_{\tf{\pi_{\man{E}}}(Z)}X)}+
A_{\pi_{\man{E}}}(\hor(Z),\hl{X})+A_{\pi_{\man{E}}}(\hl{X},\ver(Z))
\end{equation*}
using Lemma~\pldblref{lem:riemannsub}{pl:riemannsub11}\@.  Now we have
\begin{align*}
B_{\pi_{\man{E}}}(\hl{X},Z)=&\;A_{\pi_{\man{E}}}(\hor(Z),\hl{X})+
A_{\pi_{\man{E}}}(\hl{X},\ver(Z))+A_{\pi_{\man{E}}}(Z,\ver(\hl{X}))\\
=&\;A_{\pi_{\man{E}}}(\hor(Z),\hl{X})+A_{\pi_{\man{E}}}(\hl{X},\ver(Z)).
\end{align*}
Thus we have
\begin{equation*}
\nabla^{\man{E}}_Z\hl{X}=\hl{(\nabla^{\man{M}}_{\tf{\pi_{\man{E}}}(Z)}X)}+
B_{\pi_{\man{E}}}(\hl{X},Z).
\end{equation*}

Now it suffices to prove this part of the lemma for $A=\hl{A_0}\otimes\hl{X}$
for $A_0\in\sections[r]{\tensor*[k]{\ctb{\man{M}}}}$ and
$X\in\sections[r]{\tb{\man{M}}}$\@.  In this case we calculate, for
$Z\in\sections[r]{\tb{\man{E}}}$\@,
\begin{align*}
\nabla^{\man{E}}_Z(\hl{A_0}\otimes\hl{X})=&\;
(\nabla^{\man{E}}_Z\hl{A_0})\otimes\hl{X}+
\hl{A_0}\otimes\nabla^{\man{E}}_Z\hl{X}\\
=&\;\hl{(\nabla^{\man{M}}_{\tf{\pi_{\man{E}}}(Z)}A_0)}\otimes\hl{X}-
\sum_{j=1}^k\Ins_j(\hl{A_0},B_{\pi_{\man{E}},Z})\otimes\hl{X}\\
&\;+\hl{A_0}\otimes\hl{(\nabla^{\man{M}}_{\tf{\pi_{\man{E}}}(Z)}X)}+
\hl{A_0}\otimes B_{\pi_{\man{E}}}(\hl{X},\ver(Z)).
\end{align*}
We have
\begin{equation*}
\hl{(\nabla^{\man{M}}_{\tf{\pi_{\man{E}}}(Z)}A_0)}\otimes\hl{X}
+\hl{A_0}\otimes\hl{(\nabla^{\man{M}}_{\tf{\pi_{\man{E}}}(Z)}X)}=
\hl{(\nabla^{\man{M}}_{\tf{\pi_{\man{E}}}(Z)}(A_0\otimes X))}.
\end{equation*}
We also have, by~\eqref{eq:simpleAder1} and the first part of the lemma,
\begin{equation*}
\hl{A_0}\otimes B_{\pi_{\man{E}}}(\hl{X},\ver(Z))-
\sum_{j=1}^k\Ins_j(\hl{A_0},B_{\pi_{\man{E}},Z})\otimes\hl{X}=
D_{(B_{\pi_{\man{E}}})_Z}(\hl{A_0}\otimes\hl{X}).
\end{equation*}
Putting the above computations together gives this part of the lemma.

\eqref{pl:bundleform4} First we need to compute
$\nabla^{\man{E}}\vl{\lambda}$\@.  We do this by using the formula
\begin{equation*}
\lieder{Z_1}{\natpair{\vl{\lambda}}{Z_2}}=
\natpair{\nabla^{\man{E}}_{Z_1}\vl{\lambda}}{Z_2}+
\natpair{\vl{\lambda}}{\nabla^{\man{E}}_{Z_1}Z_2}
\end{equation*}
in four cases.
\begin{compactenum}
\item $Z_1=\hl{X_1}$ and $Z_2=\hl{X_2}$\@: Here we have
\begin{equation*}
\lieder{\hl{X_1}}{\natpair{\vl{\lambda}}{\hl{X_2}}}=0
\end{equation*}
and
\begin{equation*}
\natpair{\vl{\lambda}}{\nabla^{\man{E}}_{\hl{X_1}}\hl{X_2}}=
\natpair{\vl{\lambda}}{A_{\pi_{\man{E}}}(\hl{X_1},\hl{X_2})}
\end{equation*}
(by Lemma~\pldblref{lem:riemannsub}{pl:riemannsub6}) giving
\begin{equation*}
\natpair{\nabla^{\man{E}}_{\hl{X_1}}\vl{\lambda}}{\hl{X_2}}=
-\natpair{\vl{\lambda}}{A_{\pi_{\man{E}}}(\hl{X_1},\hl{X_2})}=
\natpair{\vl{\lambda}}{A_{\pi_{\man{E}}}(\hl{X_2},\hl{X_1})}
\end{equation*}
(by Lemma~\pldblref{lem:riemannsub}{pl:riemannsub2}).  Thus we have
\begin{equation*}
\natpair{\nabla^{\man{E}}_{\hl{X}_1}\vl{\lambda}}{\hl{X_2}}=
\natpair{\dual{A}_{\pi_{\man{E}}}(\vl{\lambda},\hl{X_1})}{\hl{X_2}}.
\end{equation*}

\item $Z_1=\hl{X}$ and $Z_2=\vl{\xi}$\@: We compute
\begin{equation*}
\lieder{\hl{X}}{\natpair{\vl{\lambda}}{\vl{\xi}}}=
\hl{(\lieder{X}{\natpair{\lambda}{\xi}})}
\end{equation*}
(by Lemma~\pldblref{lem:bundlefuncs}{pl:bundlefuncs1}) and
\begin{equation*}
\natpair{\vl{\lambda}}{\nabla^{\man{E}}_{\hl{X}}\vl{\xi}}=
\natpair{\vl{\lambda}}{\vl{(\nabla^{\pi_{\man{E}}}_X\xi)}}=
\hl{\natpair{\lambda}{\nabla^{\pi_{\man{E}}}_X\xi}}
\end{equation*}
(by Lemma~\pldblref{lem:riemannsub}{pl:riemannsub13}).  Thus
\begin{equation*}
\natpair{\nabla^{\man{E}}_{\hl{X}}\vl{\lambda}}{\vl{\xi}}=
\hl{(\lieder{X}{\natpair{\lambda}{\xi}})}-
\hl{\natpair{\lambda}{\nabla^{\pi_{\man{E}}}_X\xi}}
\end{equation*}
or
\begin{equation*}
\natpair{\nabla^{\man{E}}_{\hl{X}}\vl{\lambda}}{\vl{\xi}}=
\natpair{\vl{(\nabla^{\pi_{\man{E}}}_X\lambda)}}{\vl{\xi}}.
\end{equation*}

\item $Z_1=\vl{\xi}$ and $Z_2=\hl{X}$\@: In this case we compute
\begin{equation*}
\lieder{\vl{\xi}}{\natpair{\vl{\lambda}}{\hl{X}}}=0
\end{equation*}
and
\begin{equation*}
\natpair{\vl{\lambda}}{\nabla^{\man{E}}_{\vl{\xi}}\hl{X}}=0
\end{equation*}
(by Lemma~\pldblref{lem:riemannsub}{pl:riemannsub11}) giving
\begin{equation*}
\natpair{\nabla^{\man{E}}_{\vl{\xi}}\vl{\lambda}}{\hl{X}}=0.
\end{equation*}

\item $Z_1=\vl{\xi_1}$ and $Z_2=\vl{\xi_2}$\@: We have
\begin{equation*}
\lieder{\vl{\xi_1}}{\natpair{\vl{\lambda}}{\vl{\xi_2}}}=
\hl{\lieder{\vl{\xi_1}}{\natpair{\lambda}{\xi_2}}}=0
\end{equation*}
(by Lemma~\pldblref{lem:bundlefuncs}{pl:bundlefuncs2}) and
\begin{equation*}
\natpair{\vl{\lambda}}{\nabla^{\man{E}}_{\vl{\xi_1}}\vl{\xi_2}}=0
\end{equation*}
(using Lemma~\pldblref{lem:riemannsub}{pl:riemannsub12}).  This gives
\begin{equation*}
\natpair{\nabla^{\man{E}}_{\vl{\xi_1}}\vl{\lambda}}{\vl{\xi_2}}=0.
\end{equation*}
\end{compactenum}
Putting the above together,
\begin{equation*}
\nabla^{\man{E}}_Z\vl{\lambda}=
\vl{(\nabla^{\pi_{\man{E}}}_{\tf{\pi_{\man{E}}}(Z)}\lambda)}+
\hor(\dual{A}_{\pi_{\man{E}}}(\vl{\lambda},\hor(Z))).
\end{equation*}
Now we note that
\begin{align*}
\natpair{\dual{B}_{\pi_{\man{E}}}(\vl{\lambda},Z_1)}{Z_2}=&\;
\natpair{\vl{\lambda}}{B_{\pi_{\man{E}}}(Z_2,Z_1)}\\
=&\;\natpair{\vl{\lambda}}{A_{\pi_{\man{E}}}(\hor(Z_1),\hor(Z_2))}+
\natpair{\vl{\lambda}}{A_{\pi_{\man{E}}}(Z_1,\ver(Z_2))}\\
&\;+\natpair{\vl{\lambda}}{A_{\pi_{\man{E}}}(Z_2,\ver(Z_1))}\\
=&\;-\natpair{\vl{\lambda}}{A_{\pi_{\man{E}}}(\hor(Z_2),\hor(Z_1))}\\
=&\;-\natpair{\dual{A}_{\pi_{\man{E}}}(\vl{\lambda},\hor(Z_1))}{\hor(Z_2)},
\end{align*}
using Lemma~\pldblref{lem:riemannsub}{pl:riemannsub11}\@.  Thus
\begin{equation*}
\nabla^{\man{E}}_Z\vl{\lambda}=
\vl{(\nabla^{\pi_{\man{E}}}_{\tf{\pi_{\man{E}}}(Z)}\lambda)}-
\dual{B}_{\pi_{\man{E}}}(\vl{\lambda},Z).
\end{equation*}

Now, it suffices to prove this part of the lemma for
$A=A_0\otimes\lambda\otimes\eta$ for
$A_0\in\sections[r]{\tensor*[k]{\ctb{\man{M}}}}$\@,
$\lambda\in\sections[r]{\dual{\man{E}}}$\@, and
$\eta\in\sections[r]{\man{F}}$\@.  Here we calculate, for
$Z\in\sections[r]{\tb{\man{E}}}$\@,
\begin{align*}
\nabla^{\man{E},\pi_{\man{F}}}_Z(\hl{A_0}&
\otimes\vl{\lambda}\otimes\pi_{\man{E}}^*\eta)=
(\nabla^{\man{E}}_Z\hl{A_0})\otimes\vl{\lambda}\otimes\pi_{\man{E}}^*\eta+
\hl{A_0}\otimes\nabla^{\man{E}}_Z\vl{\lambda}\otimes\eta+
\hl{A_0}\otimes\vl{\lambda}+
\pi_{\man{E}}^*\nabla^{\pi_{\man{F}}}_Z\pi_{\man{E}}^*\eta\\
=&\;\hl{(\nabla^{\man{M}}_{\tf{\pi_{\man{E}}}(Z)}A_0)}\otimes
\vl{\lambda}\otimes\pi_{\man{E}}^*\eta-
\sum_{j=1}^k\Ins_j(\hl{A_0},B_{\pi_{\man{E}},Z})\otimes
\vl{\lambda}\otimes\pi_{\man{E}}^*\eta\\
&\;+\hl{A_0}\otimes\vl{(\nabla^{\pi_{\man{E}}}_{\tf{\pi_{\man{E}}}(Z)}\lambda)}
\otimes\eta-\hl{A_0}\otimes\dual{B}_{\pi_{\man{E}}}(\vl{\lambda},Z)+
\hl{A_0}\otimes\vl{\lambda}\otimes
\pi_{\man{E}}^*(\nabla^{\pi_{\man{F}}}_{\tf{\pi_{\man{E}}}(Z)}\eta).
\end{align*}
We have
\begin{align*}
\hl{(\nabla^{\man{M}}_{\tf{\pi_{\man{E}}}(Z)}A_0)}&\otimes\vl{\lambda}+
\hl{A_0}\otimes\vl{(\nabla^{\pi_{\man{E}}}_{\tf{\pi_{\man{E}}}(Z)}\lambda)}+
\hl{A_0}\otimes\vl{\lambda}\otimes
\pi_{\man{E}}^*(\nabla^{\pi_{\man{F}}}_{\tf{\pi_{\man{E}}}(Z)}\eta)\\
=&\;\vl{(\nabla^{\man{M}}_{\tf{\pi_{\man{E}}}(Z)}A_0\otimes\lambda\otimes\eta)}+
\vl{(A_0\otimes\nabla^{\pi_{\man{E}}}_{\tf{\pi_{\man{E}}}(Z)}\lambda\otimes\eta)}+
\vl{(A_0\otimes\vl{\lambda}\otimes\nabla^{\pi_{\man{F}}}_{\tf{\pi_{\man{E}}}(Z)}
\eta)}\\
=&\;\vl{(\nabla^{\pi_{\man{E}}\otimes\pi_{\man{F}}}_{\tf{\pi_{\man{E}}}(Z)}
(A_0\otimes\lambda\otimes\eta))}
\end{align*}
and, by~\eqref{eq:simpleAder2} and the first part of the lemma,
\begin{equation*}
-\hl{A_0}\otimes\dual{B}_{\pi_{\man{E}}}(\vl{\lambda},Z)\otimes
\pi_{\man{E}}^*\eta-\sum_{j=1}^k\Ins_j(\hl{A_0},B_{\pi_{\man{E}},Z})\otimes
\vl{\lambda}\otimes\pi_{\man{E}}^*\eta=
D_{B_{\pi_{\man{E}},Z}}(\hl{A_0}\otimes\vl{\lambda}\otimes\pi_{\man{E}}^*\eta).
\end{equation*}
Assembling the preceding computations gives this part of the lemma.

\eqref{pl:bundleform5} This is a slight modification of the preceding part of
the proof, taking the formula~\eqref{eq:pi*xi=xiv} into account.

\eqref{pl:bundleform6} By Lemma~\pldblref{lem:bundlefuncs}{pl:bundlefuncs3}
and~\eqref{pl:bundlefuncs4} we have
\begin{equation*}
\nabla^{\man{E}}_Z\ve{\lambda}=\lieder{Z}{\ve{\lambda}}=
\ve{(\nabla^{\pi_{\man{E}}}_{\tf{\pi_{\man{E}}}(Z)}\lambda)}+
\natpair{\vl{\lambda}}{Z}.
\end{equation*}

With the constructions following Definition~\ref{def:homlift} in mind, we
work with $A=\hl{A_0}\otimes\ve{\lambda}\otimes\pi_{\man{E}}^*\eta$ for
$A_0\in\sections[r]{\tensor*[k]{\ctb{\man{M}}}}$\@,
$\lambda\in\sections[r]{\dual{\man{E}}}$\@, and
$\eta\in\sections[r]{\man{F}}$\@.  If we keep in mind that $\ve{\lambda}$ is
a function, then we can simply write
$A=A_0\otimes(\ve{\lambda}\pi_{\man{E}}^*\eta)$\@.  We now calculate
\begin{align*}
\nabla^{\man{E},\pi_{\man{F}}}_Z(\hl{A_0}\otimes&
\ve{\lambda}\otimes\pi_{\man{E}}^*\eta)=
(\nabla^{\man{E}}_Z\hl{A_0})\otimes\ve{\lambda}\otimes\pi_{\man{E}}^*\eta+
\hl{A_0}\otimes(\nabla^{\man{E}}_Z\ve{\lambda})\otimes\pi_{\man{E}}^*\eta+
\hl{A_0}\otimes\ve{\lambda}\otimes
(\pi_{\man{E}}^*\nabla^{\pi_{\man{F}}}_Z\pi_{\man{E}}^*\eta)\\
=&\;\hl{(\nabla^{\man{M}}_{\tf{\pi_{\man{E}}}(Z)}A_0)}\otimes
\ve{\lambda}\otimes\pi_{\man{E}}^*\eta-
\sum_{j=1}^k\Ins_j(\hl{A_0},B_{\pi_{\man{E}}})\otimes\ve{\lambda}
\otimes\pi_{\man{E}}^*\eta\\
&\;+\hl{A_0}\otimes\ve{(\nabla^{\pi_{\man{E}}}_{\tf{\pi_{\man{E}}}(Z)}\lambda)}
\otimes\pi_{\man{E}}^*\eta+
\hl{A_0}\otimes(\vl{\lambda}(Z))\otimes\pi_{\man{E}}^*\eta\\
&\;+\hl{A_0}\otimes\ve{\lambda}\otimes\pi_{\man{E}}^*
(\nabla^{\pi_{\man{F}}}_{\tf{\pi_{\man{E}}}(Z)}\eta)+
\hl{A_0}\otimes\ve{\lambda}\otimes B_{\pi_{\man{E}}}(\pi_{\man{E}}^*\eta,Z).
\end{align*}
We have
\begin{align*}
\hl{(\nabla^{\man{M}}_{\tf{\pi_{\man{E}}}(Z)}A_0)}&\otimes\ve{\lambda}
\otimes\pi_{\man{E}}^*\eta+
\hl{A_0}\otimes\ve{(\nabla^{\pi_{\man{E}}}_{\tf{\pi_{\man{E}}}(Z)}\lambda)}
\otimes\pi_{\man{E}}^*\eta+\hl{A_0}\otimes\ve{\lambda}\otimes\pi_{\man{E}}^*
(\nabla^{\pi_{\man{F}}}_{\tf{\pi_{\man{E}}}(Z)}\eta)\\
=&\;\ve{(\nabla^{\man{M}}_{\tf{\pi_{\man{E}}}(Z)}A_0\otimes\lambda\otimes\eta)}+
\ve{(A_0\otimes\nabla^{\pi_{\man{E}}}_{\tf{\pi_{\man{E}}}(Z)}\lambda\otimes\eta)}+
\ve{(A_0\otimes\lambda\otimes\nabla^{\pi_{\man{F}}}_{\tf{\pi_{\man{E}}}(Z)}\eta)}\\
=&\;\ve{(\nabla^{\man{M},\pi_{\man{E}}\otimes\pi_{\man{F}}}_{\tf{\pi_{\man{E}}}(Z)}
(A_0\otimes\lambda\otimes\eta))}.
\end{align*}
Next we note that
\begin{equation*}
\hl{A_0}\otimes\ve{\lambda}\otimes B_{\pi_{\man{E}}}(\pi_{\man{E}}^*\eta,Z)-
\sum_{j=1}^k\Ins_j(\hl{A_0},B_{\pi_{\man{E}}})\otimes\ve{\lambda}
\otimes\pi_{\man{E}}^*\eta=
D_{B_{\pi_{\man{E}}},Z}(A_0\otimes(\ve{\lambda}\pi_{\man{E}}^*\eta)),
\end{equation*}
keeping in mind that $\ve{\lambda}$ is a function, so the tensor products
with $\ve{\lambda}$ are just multiplication.  Again making reference to the
constructions following Definition~\ref{def:homlift}\@, we have
\begin{equation*}
\hl{A_0}\otimes\vl{\lambda}\otimes\pi_{\man{E}}^*\eta=
\vl{(A_0\otimes\lambda\otimes\eta)},
\end{equation*}
and the lemma follows by combining the preceding three formulae.

\eqref{pl:bundleform7} This is a slight modification of the preceding part of
the proof, taking the formula~\eqref{eq:pi*xi=xiv} into account.
\end{proof}
\end{lemma}

\subsection{Prolongation}\label{subsec:prolongation}

In our geometric setting, differentiation means ``prolongation'' by taking
jets.  In this section, we illustrate how our decompositions of
Section~\ref{subsec:jetdecomp} interact with prolongation.

Let $r\in\{\infty,\omega\}$\@.  We let
$\map{\pi_{\man{E}}}{\man{E}}{\man{M}}$ be a $\C^r$-vector bundle with
$\jet{m}{\man{E}}$\@, $m\in\integernn$\@, its jet bundles.  We suppose that
we have a $\C^r$-affine connection $\nabla^{\man{M}}$ on $\man{M}$ and a
$\C^r$-vector bundle connection $\nabla^{\pi_{\man{E}}}$ in $\man{E}$\@.
Because we have the decomposition
\begin{equation*}
\jet{m}{\man{E}}\simeq
\bigoplus_{j=0}^m(\Symalg*[j]{\ctb{\man{M}}}\otimes\man{E})
\end{equation*}
by Lemma~\ref{lem:Jmdecomp}\@, it follows that the vector bundle
$\jet{m}{\man{E}}$ has a $\C^r$-connection that we denote by
$\nabla^{\pi_m}$\@.  Explicitly,
\begin{equation*}
\nabla^{\pi_m}_Xj_m\xi=
(S^m_{\nabla^{\man{M}},\nabla^{\pi_{\man{E}}}})^{-1}(\nabla^{\pi_{\man{E}}}_X\xi,
\nabla^{\man{M},\pi_{\man{E}}}_XD^1_{\nabla^{\man{M}},\nabla^{\pi_{\man{E}}}}(\xi),
\dots,\nabla^{\man{M},\pi_{\man{E}}}_X
D^m_{\nabla^{\man{M}},\nabla^{\pi_{\man{E}}}}(\xi)).
\end{equation*}
Therefore, the construction of Lemma~\ref{lem:Jmdecomp} can be applied to
$\jet{m}{\man{E}}$\@, and all that remains to sort out is notation.

To this end, for $k,m\in\integernn$\@, let us denote by
$\nabla^{\man{M},\pi_m}$ the connection in
$\tensor*[k]{\ctb{\man{M}}}\oplus\jet{m}{\man{E}}$ induced, via tensor
product, by the connections $\nabla^{\man{M}}$ and $\nabla^{\pi_m}$\@.  Then,
for $\xi\in\sections[r]{\man{E}}$\@, denote
\begin{equation*}
\nabla^{\man{M},\pi_m,k}j_m\xi=\underbrace{\nabla^{\man{M},\pi_m}\cdots
(\nabla^{\man{M},\pi_m}}_{k-1~\textrm{times}}(\nabla^{\pi_m}j_m\xi))\in
\sections[r]{\tensor*[k]{\ctb{\man{M}}}\otimes\jet{m}{\man{E}}}.
\end{equation*}
We also denote
\begin{equation*}
D^k_{\nabla^{\man{M}},\nabla^{\pi_m}}(j_m\xi)=\Sym_k\otimes\id_{\jet{m}{\man{E}}}
(\nabla^{\man{M},\pi_m,k}j_m\xi)\in\sections[r]{\Symalg*[k]{\ctb{\man{M}}}
\otimes\jet{m}{\man{E}}}.
\end{equation*}

This can be refined further by explicitly decomposing $\jet{m}{\man{E}}$\@,
so let us provide the notation for making this refinement.  For
$m,k\in\integernn$ and for
$A\in\sections[r]{\tensor*[m]{\ctb{\man{M}}}\otimes\man{E}}$\@, we
denote
\begin{equation*}
\nabla^{\man{M},\pi_{\man{E}},k}A=
\underbrace{\nabla^{\man{M},\pi_{\man{E}}}\cdots
\nabla^{\man{M},\pi_{\man{E}}}}_{k~\textrm{times}}A
\in\sections[r]{\tensor*[m+k]{\ctb{\man{M}}}\otimes\man{E}}
\end{equation*}
and
\begin{equation*}
D^{k,m}_{\nabla^{\man{M}},\nabla^{\pi_{\man{E}}}}(A)=
\Sym_k\otimes\id_{\tensor*[m]{\ctb{\man{M}}}\otimes\man{E}}
(\nabla^{\man{M},\pi_{\man{E}},k}A)\in
\sections[r]{\Symalg*[k]{\ctb{\man{M}}}
\otimes\tensor*[m]{\ctb{\man{M}}}\otimes\man{E}}.
\end{equation*}
Note that, if
$A\in\sections[r]{\Symalg*[m]{\ctb{\man{M}}}\otimes\man{E}}$\@, then
\begin{equation*}
D^{k,m}_{\nabla^{\man{M}},\nabla^{\pi_{\man{E}}}}(A)\in
\sections[r]{\Symalg*[k]{\ctb{\man{M}}}
\otimes\Symalg*[m]{\ctb{\man{M}}}\otimes\man{E}}.
\end{equation*}

An immediate consequence of Lemma~\ref{lem:Jmdecomp} is then the following
result.
\begin{lemma}\label{lem:JkJmdecomp}
The maps
\begin{equation*}
\mapdef{S_{\nabla^{\man{M}},\nabla^{\pi_m}}^k}{\jet{k}{\jet{m}{\man{E}}}}
{\bigoplus_{j=0}^k(\Symalg*[j]{\ctb{\man{M}}}\otimes\jet{m}{\man{E}})}
{j_kj_m\xi(x)}{(j_m\xi(x),D^1_{\nabla^{\man{M}},\nabla^{\pi_m}}(j_m\xi)(x),
\dots,D^k_{\nabla^{\man{M}},\nabla^{\pi_m}}(j_m\xi)(x))}
\end{equation*}
and
\begin{equation*}
\map{S_{\nabla^{\man{M}},\nabla^{\pi_{\man{E}}}}^{k,m}}{\jet{k}{\jet{m}{\man{E}}}}
{\bigoplus_{j=0}^k\left(\Symalg*[j]{\ctb{\man{M}}}\otimes\left(
\bigoplus_{l=0}^m\Symalg*[l]{\ctb{\man{M}}}\otimes\man{E}\right)\right)}
\end{equation*}
defined by
\begin{multline*}
j_kj_m\xi(x)\mapsto((\xi(x),D^1_{\nabla^{\man{M}},\nabla^{\pi_{\man{E}}}}(\xi)(x),
\dots,D^m_{\nabla^{\man{M}},\nabla^{\pi_{\man{E}}}}(\xi)(x)),\\
(D^{1,0}_{\nabla^{\man{M}},\nabla^{\pi_{\man{E}}}}\xi(x),
D^{1,1}_{\nabla^{\man{M}},\nabla^{\pi_{\man{E}}}}\scirc D^1_{\nabla^{\man{M}},\nabla^{\pi_{\man{E}}}}(\xi)(x),\dots,
D^{1,m}_{\nabla^{\man{M}},\nabla^{\pi_{\man{E}}}}\scirc D^m_{\nabla^{\man{M}},\nabla^{\pi_{\man{E}}}}(\xi)(x)),\dots,\\
(D^{k,0}_{\nabla^{\man{M}},\nabla^{\pi_{\man{E}}}}\xi(x),
D^{k,1}_{\nabla^{\man{M}},\nabla^{\pi_{\man{E}}}}\scirc D^1_{\nabla^{\man{M}},\nabla^{\pi_{\man{E}}}}(\xi)(x),\dots,
D^{k,m}_{\nabla^{\man{M}},\nabla^{\pi_{\man{E}}}}\scirc D^m_{\nabla^{\man{M}},\nabla^{\pi_{\man{E}}}}(\xi)(x)))
\end{multline*}
are isomorphisms of vector bundles, and, for each\/ $k\in\integerp$\@, the
diagrams
\begin{equation*}
\xymatrix{{\jet{k+1}{\jet{m}{\man{E}}}}
\ar[r]^(0.35){S_{\nabla^{\man{M}},\nabla^{\pi_m}}^{k+1}}
\ar[d]_{(\pi_m)^{k+1}_k}&{\bigoplus_{j=0}^{k+1}(\Symalg*[j]{\ctb{\man{M}}}
\otimes\jet{m}{\man{E}})}\ar[d]^{\pr^{k+1}_k}\\
{\jet{k}{\jet{m}{\man{E}}}}\ar[r]_(0.35){S_{\nabla^{\man{M}},\nabla^{\pi_m}}^k}&
{\bigoplus_{j=0}^k(\Symalg*[j]{\ctb{\man{M}}}\otimes\jet{m}{\man{E})}}}
\end{equation*}
and
\begin{equation*}
\xymatrix{{\jet{k+1}{\jet{m}{\man{E}}}}
\ar[r]^(0.25){S_{\nabla^{\man{M}},\nabla^{\pi_{\man{E}}}}^{k+1,m}}
\ar[d]_{(\pi_m)^{k+1}_k}&{\bigoplus_{j=0}^{k+1}
\left(\Symalg*[j]{\ctb{\man{M}}}\otimes
\left(\bigoplus_{l=0}^m\Symalg*[l]{\ctb{\man{M}}}\otimes\man{E}\right)\right)}
\ar[d]^{\pr^{k+1}_k}\\
{\jet{k}{\jet{m}{\man{E}}}}\ar[r]_(0.25)
{S_{\nabla^{\man{M}},\nabla^{\pi_{\man{E}}}}^{k,m}}&
{\bigoplus_{j=0}^k\left(\Symalg*[j]{\ctb{\man{M}}}\otimes\left(
\bigoplus_{l=0}^m\Symalg*[l]{\ctb{\man{M}}}\otimes\man{E}\right)\right)}}
\end{equation*}
commute, where\/ $\pr^{k+1}_k$ are the obvious projections, stripping off the
last component of the direct sum.
\end{lemma}

Now we recall~\cite[Definition~6.2.25]{DJS:89} the inclusion, for
$k,m\in\integernn$\@,
\begin{equation*}
\mapdef{\pi_{k,m}}{\jet{k+m}{\man{E}}}{\jet{k}{\jet{m}{\man{E}}}}
{j_{m+k}\xi(x)}{j_kj_m\xi(x).}
\end{equation*}
Let us understand this mapping using our decompositions of jet bundles.  Note
that, for a $\real$-vector space $\alg{V}$ and for $r,s\in\integernn$\@, we
have an inclusion,
\begin{equation}\label{eq:Deltasym}
\map{\Delta_{r,s}}{\Symalg*[r+s]{\dual{\alg{V}}}}
{\Symalg*[r]{\dual{\alg{V}}}\otimes\Symalg*[s]{\dual{\alg{V}}}}.
\end{equation}
Let us give an explicit formula for these inclusions.
\begin{lemma}\label{lem:Deltasym}
For a finite-dimensional\/ $\real$-vector space\/ $\alg{V}$ and for\/
$r,s\in\integernn$\@,
\begin{equation*}
\Delta_{r,s}=(\Sym_r\otimes\Sym_s)\scirc\iota_{r,s},
\end{equation*}
where
\begin{equation*}
\map{\iota_{r,s}}{\Symalg*[r+s]{\dual{\alg{V}}}}
{\tensor*[r+s]{\dual{\alg{V}}}=\tensor*[r]{\dual{\alg{V}}}\otimes
\tensor*[s]{\dual{\alg{V}}}}
\end{equation*}
is the inclusion.
\begin{proof}
We note that
$\Sym_{r+s}\scirc\Delta_{r,s}=\id_{\Symalg*[r+s]{\dual{\alg{V}}}}$\@, simply
since $\Delta_{r,s}$ is the inclusion and $\Sym_{r+s}$ is orthogonal
projection onto
$\Symalg*[r+s]{\dual{\alg{V}}}\subset\tensor*[r+s]{\dual{\alg{V}}}$\@; see
Sublemma~\ref{psublem:Symorthproj} from the proof of
Lemma~\ref{lem:Symnorm}\@.  Thus it will suffice to show that
\begin{equation*}
\Sym_{r+s}\scirc(\Sym_r\otimes\Sym_s)\scirc\iota_{r,s}=
\id_{\Symalg*[r+s]{\dual{\alg{V}}}}.
\end{equation*}
For $A\in\Symalg*[r+s]{\dual{\alg{V}}}$ we have
\begin{align*}
\Sym_r\otimes\Sym_s(A)(v_1,\dots,v_s)=&\;
\frac{1}{r!s!}\sum_{\sigma_1\in\symmgroup{r}}\sum_{\sigma_2\in\symmgroup{s}}
A(v_{\sigma_1(1)},\dots,v_{\sigma_1(r)},v_{r+\sigma_2(1)},\dots,v_{r+\sigma_2(s)})\\
=&\;A(v_1,\dots,v_r,v_{r+1},\dots,v_{r+s}),
\end{align*}
since $A$ is symmetric, and so symmetric on the first $r$ and last $s$
entries.  Since $\Sym_{r+s}(A)=A$\@, our claim follows, and so does the lemma.
\end{proof}
\end{lemma}

Now, given a $\C^r$-vector bundle $\map{\pi_{\man{E}}}{\man{E}}{\man{M}}$\@,
the preceding mapping then induces by tensor product mappings
\begin{equation*}
\map{\Delta_{r,s}\otimes\id_{\man{E}}}
{\Symalg*[r+s]{\ctb{\man{M}}}\otimes\man{E}}
{\Symalg*[r]{\ctb{\man{M}}}\otimes\Symalg*[s]{\ctb{\man{M}}}\otimes\man{E}},
\end{equation*}
and so the mapping
\begin{equation*}
\map{\what{\Delta}_{k,m}\pi_{\man{E}}}{\bigoplus_{r=0}^{k+m}
\Symalg*[r]{\ctb{\man{M}}}\otimes\man{E}}
{\bigoplus_{j=0}^k\Symalg*[j]{\ctb{\man{M}}}\otimes\left(\bigoplus_{l=0}^m
\Symalg*[l]{\ctb{\man{M}}}\otimes\man{E}\right).}
\end{equation*}
Explicitly,
\begin{multline}\label{eq:Delathatrs}
\what{\Delta}_{k,m}\pi_{\man{E}}(A_0,\dots,A_{k+m})\\
=((\Delta_{0,0}(A_0),\Delta_{0,1}(A_1),\dots,\Delta_{0,m}(A_m)),
(\Delta_{1,1}(A_1),\Delta_{1,2}(A_2),\dots,\Delta_{1,m}(A_{m+1})),\\
\dots,(\Delta_{k,0}(A_k),\Delta_{k,1}(A_{k+1}),\dots,\Delta_{k,m}(A_{k+m}))),
\end{multline}
where $A_r\in\Symalg*[r]{\ctb{\man{E}}}\otimes\man{E}$\@,
$r\in\{0,1,\dots,k+m\}$\@, and where we abbreviate
$\Delta_{j,l}\otimes\id_{\man{E}}$ by $\Delta_{j,l}$ in an attempt to achieve
concision.

We now have the following result.
\begin{lemma}\label{lem:Jk+mdecomp}
For\/ $r\in\{\infty,\omega\}$ and a\/ $\C^r$-vector bundle\/
$\map{\pi_{\man{E}}}{\man{E}}{\man{M}}$ and for\/ $k,m\in\integernn$\@, the
following diagram commutes
\begin{equation*}
\xymatrix{\jet{k+m}{\man{E}}\ar[r]^{\pi_{k,m}}
\ar[d]_{S^{k+m}_{\nabla^{\man{M}},\nabla^{\pi_{\man{E}}}}}&
{\jet{k}{\jet{m}{\man{E}}}}
\ar[d]^{S^{k,m}_{\nabla^{\man{M}},\nabla^{\pi_{\man{E}}}}}\\
{\bigoplus_{r=0}^{k+m}\Symalg*[r]{\ctb{\man{M}}}\otimes\man{E}}
\ar[r]_(0.37){\what{\Delta}_{k,m}\pi_{\man{E}}}&
{\bigoplus_{j=0}^k\Symalg*[j]{\ctb{\man{M}}}\otimes\left(\bigoplus_{l=0}^m
\Symalg*[l]{\ctb{\man{M}}}\otimes\man{E}\right)}}
\end{equation*}
\begin{proof}
We note that, by definition of the symbols involved,
\begin{equation*}
\nabla^{\man{M},\pi_m,k}(\nabla^{\man{M},\pi,m}\xi)=
\nabla^{\man{M},\pi,m+k}\xi,\qquad k,m\in\integernn,\
\xi\in\sections[r]{\man{E}}.
\end{equation*}
By Lemma~\ref{lem:Deltasym}\@, we have
\begin{align*}
D^{k+m}_{\nabla^{\man{M}},\nabla^{\pi_{\man{E}}}}\xi=&\;
\Sym_{k+m}\otimes\id_{\man{E}}(\nabla^{\man{M},\pi_{\man{E}},k+m}\xi)\\
=&\;\Delta_{k,m}\otimes\id_{\man{E}}(\nabla^{\man{M},\pi_{\man{E}},k+m}\xi)\\
=&\;\Delta_{k,m}\otimes\id_{\man{E}}(\nabla^{\man{M},\pi_m,k}
(\nabla^{\man{M},\pi_{\man{E}},m}\xi))\\
=&\;\Sym_k\otimes\Sym_m\otimes\id_{\man{E}}
(\nabla^{\man{M},\pi_m,k}(\nabla^{\man{M},\pi_{\man{E}},m}\xi))\\
=&\;D^{k,m}_{\nabla^{\man{M}},\nabla^{\pi_{\man{E}}}}
(D^m_{\nabla^{\man{M}},\nabla^{\pi_{\man{E}}}}\xi).
\end{align*}
Using this observation, and the definition of the mappings
$S^{k+m}_{\nabla^{\man{M}},\nabla^{\pi_{\man{E}}}}$ and
$S^{k,m}_{\nabla^{\man{M}},\nabla^{\pi_{\man{E}}}}$\@, the lemma follows by a
straightforward computation involving mere notation.
\end{proof}
\end{lemma}

\section{Isomorphisms defined by lifts and pull-backs}\label{sec:lift-isomorphisms}

Let $r\in\{\infty,\omega\}$ and let $\map{\pi_{\man{E}}}{\man{E}}{\man{M}}$
be a $\C^r$-vector bundle.  In this section we carefully study isomorphisms
that arise from lifts of objects on $\man{M}$ to objects on $\man{E}$\@, of
the sorts introduced in
Sections~\ref{subsec:vbfunctions}\@,~\ref{subsec:vbvfs}\@,
and~\ref{subsec:vblinmaps}\@.  In particular, we shall see that jets of
geometric objects can be decomposed (as in Section~\ref{subsec:jetdecomp})
before or after lifting.  We wish here to relate these two sorts of
decompositions for all of the lifts we consider in this work.  This makes use
of our constructions of Section~\ref{sec:tensor-derivatives} to give explicit
decompositions for jets of certain sections of certain jet bundles on the
total space of a vector bundle.  Indeed, it is the results in this section
that provide the motivation for the rather intricate constructions of
Section~\ref{sec:tensor-derivatives}\@.  For these constructions, we
additionally suppose that we have a Riemannian metric $\metric_{\man{M}}$ on
$\man{M}$ and a fibre metric $\metric_{\pi_{\man{E}}}$ on $\man{E}$\@.  We
suppose that $\nabla^{\man{M}}$ is the Levi-Civita connection for
$\metric_{\man{M}}$\@.  This data gives rise to a Riemannian metric
$\metric_{\man{E}}$ on $\man{E}$ with its Levi-Civita connection
$\nabla^{\man{E}}$\@.  We break the discussion into eight cases,
corresponding to the seven parts of Lemma~\ref{lem:bundleform}\@, along with
a construction for pull-backs of functions.  The constructions, statements,
and proofs are somewhat repetitive, so we do not give proofs that are
essentially identical to previous proofs.  While the results are similar,
they are not the same, so we must go through all of the cases.  There is
probably a ``meta'' result here, but it would take a small journey in itself
to setup the framework for this.  For our purposes, we stick to a treatment
that is concrete at the cost of being dull.

In this section, as in the previous two sections, we shall state results on
an equal footing for the smooth and real analytic cases.  However, the
detailed recursion formulae we give are not really necessary if one wants to
prove the continuity results of Section~\ref{sec:continuity} in the smooth
case.  Thus one should really regard the results of this section as being
particular to the real analytic setting.

In any event, throughout this section we let $r\in\{\infty,\omega\}$\@.

\subsection{Isomorphisms for horizontal lifts of
functions}\label{subsec:Pjetisom}

Here we consider the horizontal lift mapping
\begin{equation*}
\func[r]{\man{M}}\ni f\mapsto\pi_{\man{E}}^*f\in\func[r]{\man{E}}.
\end{equation*}
We wish to relate the decomposition associated with the jets of $f$ to those
associated with the jets of $\pi_{\man{E}}^*f$\@.  Associated with this, let
us denote by $\Pjetalg{m}{\man{E}}$ the subbundle of
$\real_{\man{E}}\oplus\jetalg{m}{\man{E}}$ defined by
\begin{equation*}
\Pjetalg[e]{m}{\man{E}}=
\setdef{j_m(\pi_{\man{E}}^*f)(e)}{f\in\func[m]{\man{M}}}.
\end{equation*}
Following Lemma~\ref{lem:Jmdecomp}\@, our constructions have to do with
iterated covariant differentials.  The basis of all of our formulae will be a
formula for iterated covariant differentials of horizontal lifts of functions
on $\man{M}$\@.  Thus we let $f\in\func[\infty]{\man{M}}$ and consider
\begin{equation*}
\nabla^{\man{E},m}\pi_{\man{E}}^*f\eqdef\underbrace{\nabla^{\man{E}}\cdots
\nabla^{\man{E}}}_{m\ \textrm{times}}\pi_{\man{E}}^*f,\qquad m\in\integerp.
\end{equation*}

We state the first two lemmata that we will use.  We recall from
Lemma~\ref{lem:bundleform} the definition of $B_{\pi_{\man{E}}}$\@.
\begin{lemma}\label{lem:PiterateddersI}
Let\/ $r\in\{\infty,\omega\}$ and let\/
$\map{\pi_{\man{E}}}{\man{E}}{\man{M}}$ be a\/ $\C^r$-vector bundle, with the
data prescribed in Section~\ref{subsec:Esubmersion} to define the Riemannian
metric\/ $\metric_{\man{E}}$ on\/ $\man{E}$\@.  For\/ $m\in\integernn$\@,
there exist\/ $\C^r$-vector bundle mappings
\begin{equation*}
(A^m_s,\id_{\man{E}})\in
\vbmappings[r]{\tensor*[s]{\pi_{\man{E}}^*\ctb{\man{M}}}}
{\tensor*[m]{\ctb{\man{E}}}},\qquad s\in\{0,1,\dots,m\},
\end{equation*}
such that
\begin{equation*}
\nabla^{\man{E},m}\pi_{\man{E}}^*f=
\sum_{s=0}^mA^m_s(\pi_{\man{E}}^*\nabla^{\man{M},s}f)
\end{equation*}
for all\/ $f\in\func[m]{\man{M}}$\@.  Moreover, the vector bundle mappings\/
$A^m_0,A^m_1,\dots,A^m_m$ satisfy the recursion relations prescribed by
\begin{equation*}
A^0_0(\beta_0)=\beta_0,\enspace A^1_1(\beta_1)=\beta_1,\enspace A^1_0=0,
\end{equation*}
and
{\small\begin{align*}
A^{m+1}_{m+1}(\beta_{m+1})=&\;\beta_{m+1},\\
A^{m+1}_s(\beta_s)=&\;(\nabla^{\man{E}}A^m_s)(\beta_s)+
A^m_{s-1}\otimes\id_{\ctb{\man{E}}}(\beta_s)-
\sum_{j=1}^sA^m_s\otimes\id_{\ctb{\man{E}}}(\Ins_j(\beta_s,B_{\pi_{\man{E}}})),
\enspace s\in\{1,\dots,m\},\\
A^{m+1}_0(\beta_0)=&\;(\nabla^{\man{E}}A^m_0)(\beta_0),
\end{align*}}%
where\/ $\beta_s\in\tensor*[s]{\pi_{\man{E}}^*\ctb{\man{M}}}$\@,\/
$s\in\{0,1,\dots,m\}$\@.
\begin{proof}
The assertion clearly holds for the initial conditions of the recursion,
simply because
\begin{equation*}
\pi^*f=\pi^*f,\quad\d{(\pi^*f)}=\pi^*\d{f}+0f.
\end{equation*}
So suppose it true for $m\in\integerp$\@.  Thus
\begin{equation*}
\nabla^{\man{E},m}\pi_{\man{E}}^*f=
\sum_{s=0}^mA^m_s(\pi_{\man{E}}^*\nabla^{\man{M},s}f),
\end{equation*}
where the vector bundle mappings $A^a_s$\@, $a\in\{0,1,\dots,m\}$\@,
$s\in\{0,1,\dots,a\}$\@, satisfy the recursion relations from the statement
of the lemma.  Then
\begin{align*}
\nabla^{\man{E},m+1}\pi_{\man{E}}^*f=&\;
\sum_{s=0}^m(\nabla^{\man{E}}A^m_s)(\pi_{\man{E}}^*\nabla^{\man{M},s}f)+
\sum_{s=0}^mA^m_s\otimes\id_{\ctb{\man{E}}}
(\nabla^{\man{E}}\pi_{\man{E}}^*\nabla^{\man{M},s}f)\\
=&\;\sum_{s=0}^m(\nabla^{\man{E}}A^m_s)(\pi_{\man{E}}^*\nabla^{\man{M},s}f)+
\sum_{s=0}^mA^m_s\otimes
\id_{\ctb{\man{E}}}(\pi_{\man{E}}^*\nabla^{\man{M},s+1}f)\\
&\;-\sum_{s=0}^m\sum_{j=1}^sA^m_s\otimes\id_{\ctb{\man{E}}}
(\Ins_j(\pi_{\man{E}}^*\nabla^{\man{M},s}f,B_{\pi_{\man{E}}}))\\
=&\;\pi_{\man{E}}^*\nabla^{\man{M},m+1}f+
\left.\sum_{s=1}^m\right((\nabla^{\man{E}}A^m_s)
(\pi_{\man{E}}^*\nabla^{\man{M},s}f)
+A^m_{s-1}\otimes\id_{\ctb{\man{E}}}(\pi_{\man{E}}^*\nabla^{\man{M},s}f)\\
&\;-\left.\sum_{j=1}^sA^m_s\otimes\id_{\ctb{\man{E}}}
(\Ins_j(\pi_{\man{E}}^*\nabla^{\man{M},s}f,B_{\pi_{\man{E}}}))\right)
+(\nabla^{\man{E}}A^m_0)(\pi_{\man{E}}^*f)
\end{align*}
by Lemma~\pldblref{lem:bundleform}{pl:bundleform1}\@.  From this, the lemma
follows.
\end{proof}
\end{lemma}

We shall also need to ``invert'' the relationship of the preceding lemma.
\begin{lemma}\label{lem:PiterateddersII}
Let\/ $r\in\{\infty,\omega\}$ and let\/
$\map{\pi_{\man{E}}}{\man{E}}{\man{M}}$ be a\/ $\C^r$-vector bundle, with the
data prescribed in Section~\ref{subsec:Esubmersion} to define the Riemannian
metric\/ $\metric_{\man{E}}$ on\/ $\man{E}$\@.  For\/ $m\in\integernn$\@,
there exist\/ $\C^r$-vector bundle mappings
\begin{equation*}
(B^m_s,\id_{\man{E}})\in
\vbmappings[r]{\tensor*[s]{\ctb{\man{E}}}}
{\tensor*[m]{\pi_{\man{E}}^*\ctb{\man{M}}}},\qquad s\in\{0,1,\dots,m\},
\end{equation*}
such that
\begin{equation*}
\pi_{\man{E}}^*\nabla^{\man{M},m}f=
\sum_{s=0}^mB^m_s(\nabla^{\man{E},s}\pi_{\man{E}}^*f)
\end{equation*}
for all\/ $f\in\func[m]{\man{M}}$\@.  Moreover, the vector bundle mappings\/
$B^m_0,B^m_1,\dots,B^m_m$ satisfy the recursion relations prescribed by
\begin{equation*}
B^0_0(\alpha_0)=\alpha_0,\enspace B^1_1(\alpha_1)=\alpha_1,\enspace B^1_0=0,
\end{equation*}
and
\begin{align*}
B^{m+1}_{m+1}(\alpha_{m+1})=&\;\alpha_{m+1},\\
B^{m+1}_s(\alpha_s)=&\;(\nabla^{\man{E}}B^m_s)(\alpha_s)+
B^m_{s-1}\otimes\id_{\ctb{\man{E}}}(\alpha_s)+
\sum_{j=1}^m\Ins_j(B^m_s(\alpha_s),B_{\pi_{\man{E}}}),
\enspace s\in\{1,\dots,m\},\\
B^{m+1}_0(\alpha_0)=&\;(\nabla^{\man{E}}B^m_0)(\alpha_0)+
\sum_{j=1}^m\Ins_j(B^m_0(\alpha_0),B_{\pi_{\man{E}}}),
\end{align*}
where\/ $\alpha_s\in\tensor*[s]{\ctb{\man{E}}}$\@,\/
$s\in\{0,1,\dots,m\}$\@.
\begin{proof}
The assertion clearly holds for the initial conditions for the recursion because
\begin{equation*}
\pi^*f=\pi^*f,\quad\pi^*(\d{f})=\d{(\pi^*f)}+0f.
\end{equation*}
So suppose it true for $m\in\integerp$\@.  Thus
\begin{equation}\label{eq:PiteratedderII1}
\pi_{\man{E}}^*\nabla^{\man{M},m}f=
\sum_{s=0}^mB^m_s(\nabla^{\man{E},s}\pi_{\man{E}}^*f),
\end{equation}
where the vector bundle mappings $B^a_s$\@, $a\in\{0,1,\dots,m\}$\@,
$s\in\{0,1,\dots,a\}$\@, satisfy the recursion relations from the statement
of the lemma.  Then, by Lemma~\pldblref{lem:bundleform}{pl:bundleform1}\@, we
can work on the left-hand side of~\eqref{eq:PiteratedderII1} to give
\begin{align*}
\nabla^{\man{E}}\pi_{\man{E}}^*\nabla^{\man{M},m}f=&\;
\pi_{\man{E}}^*\nabla^{\man{M},m+1}f-\sum_{j=1}^m
\Ins_j(\pi_{\man{E}}^*\nabla^{\man{M},m}f,B_{\pi_{\man{E}}})\\
=&\;\pi_{\man{E}}^*\nabla^{\man{M},m+1}f-\sum_{s=0}^m\sum_{j=1}^m
\Ins_j(B^m_s(\nabla^{\man{E},s}\pi_{\man{E}}^*f),B_{\pi_{\man{E}}}).
\end{align*}
Working on the right-hand side of~\eqref{eq:PiteratedderII1} gives
\begin{equation*}
\nabla^{\man{E}}\pi_{\man{E}}^*\nabla^{\man{M},m}f=
\sum_{s=0}^m\nabla^{\man{E}}B^m_s(\nabla^{\man{E},s}\pi_{\man{E}}^*f)+
\sum_{s=0}^mB^m_s\otimes\id_{\ctb{\man{E}}}(\nabla^{\man{E},s+1}\pi_{\man{E}}^*f).
\end{equation*}
Combining the preceding two equations gives
\begin{align*}
\pi_{\man{E}}^*\nabla^{\man{M},m+1}f=&\;
\sum_{s=0}^m\nabla^{\man{E}}B^m_s(\nabla^{\man{E},s}\pi_{\man{E}}^*f)+
\sum_{s=0}^mB^m_s\otimes
\id_{\ctb{\man{E}}}(\nabla^{\man{E},s+1}\pi_{\man{E}}^*f)\\
&\;+\sum_{s=0}^m\sum_{j=1}^m\Ins_j
(B^m_s(\nabla^{\man{E},s}\pi_{\man{E}}^*f),B_{\pi_{\man{E}}})\\
=&\;\nabla^{\man{E},m+1}\pi_{\man{E}}^*f+
\left.\sum_{s=1}^m\right(\nabla^{\man{E}}B^m_s(\nabla^{\man{E},s}
\pi_{\man{E}}^*f)+B^m_{s-1}\otimes
\id_{\ctb{\man{E}}}(\nabla^{\man{E},s}\pi_{\man{E}}^*f)\\
&\;\left.+\sum_{j=1}^m\Ins_j(B^m_s(\nabla^{\man{E},s}\pi_{\man{E}}^*f),
B_{\pi_{\man{E}}})\right)+\nabla^{\man{E}}B^m_0(\pi_{\man{E}}^*f)+
\sum_{j=1}^m\Ins_j(B^m_0(\pi_{\man{E}}^*f),B_{\pi_{\man{E}}}),
\end{align*}
and the lemma follows from this.
\end{proof}
\end{lemma}

Next we turn to symmetrised versions of the preceding lemmata.  We show that
the preceding two lemmata induce corresponding mappings between symmetric
tensors.
\begin{lemma}\label{lem:PsymiterateddersI}
Let\/ $r\in\{\infty,\omega\}$ and let\/
$\map{\pi_{\man{E}}}{\man{E}}{\man{M}}$ be a\/ $\C^r$-vector bundle, with the
data prescribed in Section~\ref{subsec:Esubmersion} to define the Riemannian
metric\/ $\metric_{\man{E}}$ on\/ $\man{E}$\@.  For\/ $m\in\integernn$\@,
there exist\/ $\C^r$-vector bundle mappings
\begin{equation*}
(\what{A}^m_s,\id_{\man{E}})\in
\vbmappings[r]{\Symalg*[s]{\pi_{\man{E}}^*\ctb{\man{M}}}}
{\Symalg*[m]{\ctb{\man{E}}}},\qquad s\in\{0,1,\dots,m\},
\end{equation*}
such that
\begin{equation*}
\Sym_m\scirc\nabla^{\man{E},m}\pi_{\man{E}}^*f=
\sum_{s=0}^m\what{A}^m_s(\Sym_s\scirc\pi_{\man{E}}^*\nabla^{\man{M},s}f)
\end{equation*}
for all\/ $f\in\func[m]{\man{M}}$\@.
\begin{proof}
We define $\map{A^m}{\tensor*[\le m]{\pi_{\man{E}}^*\ctb{\man{M}}}}
{\tensor*[\le m]{\ctb{\man{E}}}}$ by
\begin{equation*}
A^m(\pi_{\man{E}}^*f,\pi_{\man{E}}^*\nabla^{\man{M}}f,\dots,
\pi_{\man{E}}^*\nabla^{\man{M},m}f)=
\left(A^0_0(\pi_{\man{E}}^*f),\sum_{s=0}^1
A^1_s(\pi_{\man{E}}^*\nabla^{\man{M},s}f),
\dots,\sum_{s=0}^mA^m_s(\pi_{\man{E}}^*\nabla^{\man{M},s}f)\right).
\end{equation*}
Let us organise the mappings we require into the following diagram:
\begin{equation}\label{eq:Psymmetrise1}
\xymatrix{{\tensor*[\le m]{\pi_{\man{E}}^*\ctb{\man{M}}}}
\ar[r]^{\Sym_{\le m}}\ar[d]_{A^m}&
{\Symalg*[\le m]{\pi_{\man{E}}^*\ctb{\man{M}}}}
\ar[r]^(0.45){S^m_{\nabla^{\man{M}}}}
\ar[d]^{\what{A}^m}&{\pi_{\man{E}}^*(\real_{\man{M}}\oplus\jetalg{m}{\man{M}})}
\ar[d]^{\id_{\real}\oplus j_m\pi_{\man{E}}}\\
{\tensor*[\le m]{\ctb{\man{E}}}}\ar[r]^{\Sym_{\le m}}&
{\Symalg*[\le m]{\ctb{\man{E}}}}\ar[r]^(0.5){S^m_{\nabla^{\man{E}}}}&
{\real_{\man{E}}\oplus\jetalg{m}{\man{E}}}}
\end{equation}
Here $\what{A}^m$ is defined so that the right square commutes.  We shall
show that the left square also commutes.  Indeed,
\begin{align*}
\what{A}^m\scirc\Sym_{\le m}(\pi_{\man{E}}^*f,\pi_{\man{E}}^*\nabla^{\man{M}}&f,\dots,
\pi_{\man{E}}^*\nabla^{\man{M},m}f)\\
=&\;(S^m_{\nabla^{\man{E}}})^{-1}\scirc(\id_{\real}\oplus j_m\pi_{\man{E}})\scirc
S^m_{\nabla^{\man{M}}}\scirc\Sym_{\le m}
(\pi_{\man{E}}^*f,\pi_{\man{E}}^*\nabla^{\man{M}}f,\dots,\pi_{\man{E}}^*
\nabla^{\man{M},m}f)\\
=&\;\Sym_{\le m}(\pi_{\man{E}}^*f,\nabla^{\man{E}}\pi_{\man{E}}^*f,\dots,
\nabla^{\man{E},m}\pi_{\man{E}}^*f)\\
=&\;\Sym_{\le m}\scirc A^m(\pi_{\man{E}}^*f,\pi_{\man{E}}^*\nabla^{\man{M}}f,
\dots,\pi_{\man{E}}^*\nabla^{\man{M},m}f).
\end{align*}
Thus the diagram~\eqref{eq:Psymmetrise1} commutes.  Now we have
\begin{multline*}
\what{A}^m\scirc\Sym_{\le m}(\pi_{\man{E}}^*f,\pi_{\man{E}}^*\nabla^{\man{M}}f,\dots,
\pi_{\man{E}}^*\nabla^{\man{M},m}f)\\
=\left(\Sym_1\scirc A^0_0(\pi_{\man{E}}^*f),
\sum_{s=0}^1\Sym_2\scirc A^1_s(\pi_{\man{E}}^*\nabla^{\man{M},s}f),
\dots,\sum_{s=0}^m\Sym_m\scirc A^m_s(\pi_{\man{E}}^*\nabla^{\man{M},s}f)\right).
\end{multline*}
Thus, if we define
\begin{equation}\label{eq:AhatA}
\what{A}^m_s(\Sym_s\scirc\pi_{\man{E}}^*\nabla^{\man{M},s}f)=
\Sym_m\scirc A^m_s(\pi_{\man{E}}^*\nabla^{\man{M},s}f),
\end{equation}
then we have
\begin{equation*}
\Sym_m\scirc\nabla^{\man{E},m}\pi_{\man{E}}^*f=
\sum_{s=0}^m\what{A}^m_s(\Sym_s\scirc\pi_{\man{E}}^*\nabla^{\man{M},s}f),
\end{equation*}
as desired.
\end{proof}
\end{lemma}

Next we consider the ``inverse'' of the preceding lemma.
\begin{lemma}\label{lem:PsymiterateddersII}
Let\/ $r\in\{\infty,\omega\}$ and let\/
$\map{\pi_{\man{E}}}{\man{E}}{\man{M}}$ be a\/ $\C^r$-vector bundle, with the
data prescribed in Section~\ref{subsec:Esubmersion} to define the Riemannian
metric\/ $\metric_{\man{E}}$ on\/ $\man{E}$\@.  For\/ $m\in\integernn$\@,
there exist\/ $\C^r$-vector bundle mappings
\begin{equation*}
(\what{B}^m_s,\id_{\man{E}})\in
\vbmappings[r]{\Symalg*[s]{\ctb{\man{E}}}}
{\Symalg*[m]{\pi_{\man{E}}^*\ctb{\man{M}}}},\qquad s\in\{0,1,\dots,m\},
\end{equation*}
such that
\begin{equation*}
\Sym_m\scirc\pi_{\man{E}}^*\nabla^{\man{M},m}f=
\sum_{s=0}^m\what{B}^m_s(\Sym_s\scirc\nabla^{\man{E},s}\pi_{\man{E}}^*f)
\end{equation*}
for all\/ $f\in\func[m]{\man{M}}$\@.
\begin{proof}
We define $\map{B^m}{\tensor*[\le m]{\ctb{\man{E}}}}
{\tensor*[\le m]{\pi_{\man{E}}^*\ctb{\man{M}}}}$ by requiring that
\begin{equation*}
B^m(\pi_{\man{E}}^*f,\dots,\nabla^{\man{E},m}\pi_{\man{E}}^*f)=
\left(B^0_0(\pi_{\man{E}}^*f),
\sum_{s=0}^1B^1_s(\nabla^{\man{E},s}\pi_{\man{E}}^*f),
\dots,\sum_{s=0}^mB^m_s(\nabla^{\man{E},m}\pi_{\man{E}}^*f)\right),
\end{equation*}
as in Lemma~\ref{lem:PiterateddersII}\@.  Note that the mapping
\begin{equation*}
\map{\id_{\real}\oplus j_m\pi_{\man{E}}}
{\pi_{\man{E}}^*(\real_{\man{M}}\oplus\jetalg{m}{\man{M}})}
{\Pjetalg{m}{\man{E}}}
\end{equation*}
is well-defined and a vector bundle isomorphism.  Let us organise the
mappings we require into the following diagram:
\begin{equation}\label{eq:Psymmetrise2}
\xymatrix{{\tensor*[\le m]{\ctb{\man{E}}}}
\ar[r]^{\Sym_{\le m}}\ar[d]_{B^m}&
{\Symalg*[\le m]{\ctb{\man{E}}}}\ar[r]^(0.55){S^m_{\nabla^{\man{E}}}}
\ar[d]^{\what{B}^m}&{\Pjetalg{m}{\man{E}}}
\ar@{<-}[d]^{\id_{\real}\oplus j_m\pi_{\man{E}}}\\
{\tensor*[\le m]{\pi_{\man{E}}^*\ctb{\man{M}}}}\ar[r]^{\Sym_{\le m}}&
{\Symalg*[\le m]{\pi_{\man{E}}^*\ctb{\man{M}}}}
\ar[r]^(0.45){S^m_{\nabla^{\man{M}}}}&
{\pi_{\man{E}}^*(\real_{\man{M}}\oplus\jetalg{m}{\man{M}})}}
\end{equation}
Here $\what{B}^m$ is defined so that the right square commutes.  We shall
show that the left square also commutes.  Indeed,
\begin{align*}
\what{B}^m\scirc\Sym_{\le m}(\pi_{\man{E}}^*f,\nabla^{\man{E}}&\pi_{\man{E}}^*f,\dots,
\nabla^{\man{E},m}\pi_{\man{E}}^*f)\\
=&\;(S^m_{\nabla^{\man{M}}})^{-1}\scirc
(\id_{\real}\oplus j_m\pi_{\man{E}})^{-1}\scirc
S^m_{\nabla^{\man{E}}}\scirc\Sym_{\le m}
(\pi_{\man{E}}^*f,\nabla^{\man{E}}\pi_{\man{E}}^*f,\dots,
\nabla^{\man{E},m}\pi_{\man{E}}^*f)\\
=&\;\Sym_{\le m}(\pi_{\man{E}}^*f,\pi_{\man{E}}^*\nabla^{\man{M}}f,\dots,
\pi_{\man{E}}^*\nabla^{\man{M},m}f)\\
=&\;\Sym_{\le m}\scirc B^m(\pi_{\man{E}}^*f,\nabla^{\man{E}}\pi_{\man{E}}^*f,
\dots,\nabla^{\man{E},m}\pi_{\man{E}}^*f).
\end{align*}
Thus the diagram~\eqref{eq:Psymmetrise2} commutes.  Thus, if we define
$\hat{B}^m_s$ so as to satisfy
\begin{equation*}
\what{B}^m_s(\Sym_s\scirc\nabla^{\man{E},s}\pi_{\man{E}}^*f)=
\Sym_m\scirc B^m_s(\nabla^{\man{E},s}\pi_{\man{E}}^*f),
\end{equation*}
then we have
\begin{equation*}
\Sym_m\scirc\pi_{\man{E}}^*\nabla^{\man{M},m}f=
\sum_{s=0}^m\what{B}^m_s(\Sym_s\scirc\nabla^{\man{E},s}\pi_{\man{E}}^*f),
\end{equation*}
as desired.
\end{proof}
\end{lemma}

The following lemma provides two decompositions of $\Pjetalg{m}{\man{E}}$\@,
one ``downstairs'' and one ``upstairs,'' and the relationship between them.
The assertion simply results from an examination of the preceding four
lemmata.
\begin{lemma}\label{lem:Pdecomp}
Let\/ $r\in\{\infty,\omega\}$ and let\/
$\map{\pi_{\man{E}}}{\man{E}}{\man{M}}$ be a\/ $\C^r$-vector bundle, with the
data prescribed in Section~\ref{subsec:Esubmersion} to define the Riemannian
metric\/ $\metric_{\man{E}}$ on\/ $\man{E}$\@.  Then there exist\/
$\C^r$-vector bundle mappings
\begin{equation*}
A^m_{\nabla^{\man{E}}}\in\vbmappings[r]{\Pjetalg{m}{\man{E}}}
{\Symalg*[\le m]{\pi_{\man{E}}^*\ctb{\man{M}}}},\quad
B^m_{\nabla^{\man{E}}}\in\vbmappings[r]{\Pjetalg{m}{\man{E}}}
{\Symalg*[\le m]{\ctb{\man{E}}}},
\end{equation*}
defined by
\begin{align*}
A^m_{\nabla^{\man{E}}}(j_m(\pi_{\man{E}}^*f)(e))=&\;
\Sym_{\le m}(\pi_{\man{E}}^*f(e),\pi_{\man{E}}^*\nabla^{\man{M}}f(e),\dots,
\pi_{\man{E}}^*\nabla^{\man{M},m}f(e)),\\
B^m_{\nabla^{\man{E}}}(j_m(\pi_{\man{E}}^*f)(e))=&\;
\Sym_{\le m}(\pi_{\man{E}}^*f(e),\nabla^{\man{E}}\pi_{\man{E}}^*f(e),\dots,
\nabla^{\man{E},m}\pi_{\man{E}}^*f(e)).
\end{align*}
Moreover,\/ $A^m_{\nabla^{\man{E}}}$ is an isomorphism,\/
$B^m_{\nabla^{\man{E}}}$ is injective, and
\begin{multline*}
B^m_{\nabla^{\man{E}}}\scirc(A^m_{\nabla^{\man{E}}})^{-1}\scirc
(\Sym_{\le m}(\pi_{\man{E}}^*f(e),\pi_{\man{E}}^*\nabla^{\man{M}}f(e),\dots,
\pi_{\man{E}}^*\nabla^{\man{M},m}f(e))\\
=\left(A^0_0(\pi_{\man{E}}^*f(e)),\sum_{s=0}^1\what{A}^1_s
(\Sym_s\scirc\pi_{\man{E}}^*\nabla^{\man{M},s}f(e)),\dots,
\sum_{s=0}^m\what{A}^m_s(\Sym_s\scirc
\pi_{\man{E}}^*\nabla^{\man{M},s}f(e))\right)
\end{multline*}
and
\begin{multline*}
A^m_{\nabla^{\man{E}}}\scirc(B^m_{\nabla^{\man{E}}})^{-1}\scirc
\Sym_{\le m}(\pi_{\man{E}}^*f(e),\nabla^{\man{E}}\pi_{\man{E}}^*f(e),\dots,
\nabla^{\man{E},m}\pi_{\man{E}}^*f(e))\\
=\left(B^0_0(\pi_{\man{E}}^*f(e)),
\sum_{s=0}^1\what{B}^1_s(\Sym_s\scirc\nabla^{\man{E},s}\pi_{\man{E}}^*f(e)),
\dots,\sum_{s=0}^m\what{B}^m_s(\Sym_s\scirc\nabla^{\man{E},s}
\pi_{\man{E}}^*f(e))\right),
\end{multline*}
where the vector bundle mappings\/ $\what{A}^m_s$ and\/ $\what{B}^m_s$\@,\/
$s\in\{0,1,\dots,m\}$\@, are as in Lemmata~\ref{lem:PsymiterateddersI}
and~\ref{lem:PsymiterateddersII}\@.
\end{lemma}

\subsection{Isomorphisms for vertical lifts of sections}

Next we consider vertical lifts of sections,~\ie~the mapping
\begin{equation*}
\sections[r]{\man{E}}\ni\xi\mapsto\vl{\xi}\in\sections[r]{\tb{\man{E}}}.
\end{equation*}
We wish to relate the decomposition of the jets of $\xi$ with those of
$\vl{\xi}$\@.  Associated with this, we denote
\begin{equation*}
\Vjetalg[e]{m}{\man{E}}=
\setdef{j_m\vl{\xi}(e)}{\xi\in\sections[m]{\man{E}}}.
\end{equation*}
By~\eqref{eq:jet=jetalg}\@, we have
\begin{equation*}
\Vjetalg[e]{m}{\man{E}}\simeq\Pjetalg[e]{m}{\man{E}}\otimes\vb[e]{\man{E}}.
\end{equation*}
As with the constructions of the preceding section, we wish to use
Lemma~\ref{lem:Jmdecomp} to provide a decomposition of
$\Vjetalg{m}{\man{E}}$\@, and to do so we need to understand the covariant
derivatives
\begin{equation*}
\nabla^{\man{E},m}\vl{\xi}\eqdef\underbrace{\nabla^{\man{E}}\cdots
\nabla^{\man{E}}}_{m\ \textrm{times}}\vl{\xi},\qquad m\in\integernn.
\end{equation*}
In our development, we shall use the notation used in the preceding section
in a slightly different, but similar, context.  This seems reasonable since
we have to do more or less the same thing six times, and using six different
pieces of notation will be excessively burdensome.

The first result we give is the following.
\begin{lemma}\label{lem:ViterateddersI}
Let\/ $r\in\{\infty,\omega\}$ and let\/
$\map{\pi_{\man{E}}}{\man{E}}{\man{M}}$ be a\/ $\C^r$-vector bundle, with the
data prescribed in Section~\ref{subsec:Esubmersion} to define the Riemannian
metric\/ $\metric_{\man{E}}$ on\/ $\man{E}$\@.  For\/ $m\in\integernn$\@,
there exist\/ $\C^r$-vector bundle mappings
\begin{equation*}
(A^m_s,\id_{\man{E}})\in
\vbmappings[r]{\tensor*[s]{\pi_{\man{E}}^*\ctb{\man{M}}}\otimes\vb{\man{E}}}
{\tensor*[m]{\ctb{\man{E}}}\otimes\vb{\man{E}}},\qquad
s\in\{0,1,\dots,m\},
\end{equation*}
such that
\begin{equation*}
\nabla^{\man{E},m}\vl{\xi}=
\sum_{s=0}^mA^m_s(\vl{(\nabla^{\man{M},\pi_{\man{E}},s}\xi)})
\end{equation*}
for all\/ $\xi\in\sections[m]{\man{E}}$\@.  Moreover, the vector bundle
mappings\/ $A^m_0,A^m_1,\dots,A^m_m$ satisfy the recursion relations
prescribed by\/ $A^0_0(\beta_0)=\beta_0$ and
\begin{align*}
A^{m+1}_{m+1}(\beta_{m+1})=&\;\beta_{m+1},\\
A^{m+1}_s(\beta_s)=&\;(\nabla^{\man{E}}A^m_s)(\beta_s)+
A^m_{s-1}\otimes\id_{\ctb{\man{E}}}(\beta_s)-
\sum_{j=1}^sA^m_s\otimes\id_{\ctb{\man{E}}}(\Ins_j(\beta_s,B_{\pi_{\man{E}}}))\\
&\;+A^m_s\otimes\id_{\ctb{\man{E}}}(\Ins_{s+1}(\beta_s,\dual{B}_{\pi_{\man{E}}})),
\enspace s\in\{1,\dots,m\},\\
A^{m+1}_0(\beta_0)=&\;(\nabla^{\man{E}}A^m_0)(\beta_0)+
A^m_0\otimes\id_{\ctb{\man{E}}}(\Ins_1(\beta_0,\dual{B}_{\pi_{\man{E}}})),
\end{align*}
where\/ $\beta_s\in\tensor*[s]{\pi_{\man{E}}^*\ctb{\man{M}}}\otimes\vb{\man{E}}$\@,\/
$s\in\{0,1,\dots,m+1\}$\@.
\begin{proof}
The assertion clearly holds for $m=0$\@, so suppose it true for
$m\in\integerp$\@.  Thus
\begin{equation*}
\nabla^{\man{E},m}\vl{\xi}=
\sum_{s=0}^mA^m_s(\vl{(\nabla^{\man{M},\pi_{\man{E}},s}\xi)}),
\end{equation*}
where the vector bundle mappings $A^a_s$\@, $a\in\{0,1,\dots,m\}$\@,
$s\in\{0,1,\dots,a\}$\@, satisfy the recursion relations from the statement
of the lemma.  Then
{\small\begin{align*}
\nabla^{\man{E},m+1}\vl{\xi}=&\;
\sum_{s=0}^m(\nabla^{\man{E}}A^m_s)(\vl{(\nabla^{\man{M},\pi_{\man{E}},s}\xi)})+
\sum_{s=0}^mA^m_s\otimes\id_{\ctb{\man{E}}}
(\nabla^{\man{E}}\vl{(\nabla^{\man{M},\pi_{\man{E}},s}\xi)})\\
=&\;\sum_{s=0}^m(\nabla^{\man{E}}A^m_s)
(\vl{(\nabla^{\man{M},\pi_{\man{E}},s}\xi)})+
\sum_{s=0}^mA^m_s\otimes\id_{\ctb{\man{E}}}
(\vl{(\nabla^{\man{M},\pi_{\man{E}},s+1}\xi)})\\
&\;-\sum_{s=1}^m\sum_{j=1}^sA^m_s\otimes\id_{\ctb{\man{E}}}
(\Ins_j(\vl{(\nabla^{\man{M},\pi_{\man{E}},s}\xi)},B_{\pi_{\man{E}}}))\\
&\;+\sum_{s=1}^mA^m_s\otimes\id_{\ctb{\man{E}}}
(\Ins_{s+1}(\vl{(\nabla^{\man{M},\pi_{\man{E}},s}\xi)},\dual{B}_{\pi_{\man{E}}}))
+A^m_0\otimes\id_{\ctb{\man{E}}}(\Ins_1(\vl{\xi},\dual{B}_{\pi_{\man{E}}}))\\
=&\;\vl{(\nabla^{\man{M},\pi_{\man{E}},m+1}\xi)}+
\left.\sum_{s=1}^m\right((\nabla^{\man{E}}A^m_s)
(\vl{(\nabla^{\man{M},\pi_{\man{E}},s}\xi)})
+A^m_{s-1}\otimes\id_{\ctb{\man{E}}}(\vl{(\nabla^{\man{M},\pi_{\man{E}},s}\xi)})\\
&\;-\sum_{j=1}^sA^m_s\otimes\id_{\ctb{\man{E}}}
(\Ins_j(\vl{(\nabla^{\man{M},\pi_{\man{E}},s}\xi)},B_{\pi_{\man{E}}}))
+\left.\vphantom{\sum_{j=1}^s}A^m_s\otimes\id_{\ctb{\man{E}}}
(\Ins_{s+1}(\vl{(\nabla^{\man{M},\pi_{\man{E}},s}\xi)},
\dual{B}_{\pi_{\man{E}}}))\right)\\
&\;+(\nabla^{\man{E}}A^m_0)(\vl{\xi})+
A^m_0\otimes\id_{\ctb{\man{E}}}(\Ins_1(\vl{\xi},\dual{B}_{\pi_{\man{E}}}))
\end{align*}}%
by Lemma~\pldblref{lem:bundleform}{pl:bundleform2}\@.  From this, the lemma
follows.
\end{proof}
\end{lemma}

Now we ``invert'' the constructions from the preceding lemma.
\begin{lemma}\label{lem:ViterateddersII}
Let\/ $r\in\{\infty,\omega\}$ and let\/
$\map{\pi_{\man{E}}}{\man{E}}{\man{M}}$ be a\/ $\C^r$-vector bundle, with the
data prescribed in Section~\ref{subsec:Esubmersion} to define the Riemannian
metric\/ $\metric_{\man{E}}$ on\/ $\man{E}$\@.  For\/ $m\in\integernn$\@,
there exist\/ $\C^r$-vector bundle mappings
\begin{equation*}
(B^m_s,\id_{\man{E}})\in
\vbmappings[r]{\tensor*[m]{\ctb{\man{E}}}\otimes\vb{\man{E}}}
{\tensor*[m]{\pi_{\man{E}}^*\ctb{\man{M}}}\otimes\vb{\man{E}}},\qquad
s\in\{0,1,\dots,m\},
\end{equation*}
such that
\begin{equation*}
\vl{(\nabla^{\man{M},\pi_{\man{E}},m}\xi)}=
\sum_{s=0}^mB^m_s(\nabla^{\man{E},s}\vl{\xi})
\end{equation*}
for all\/ $\xi\in\sections[m]{\man{E}}$\@.  Moreover, the vector bundle
mappings\/ $B^m_0,B^m_1,\dots,B^m_m$ satisfy the recursion relations
prescribed by\/ $B^0_0(\alpha_0)=\alpha_0$ and
\begin{align*}
B^{m+1}_{m+1}(\alpha_{m+1})=&\;\alpha_{m+1},\\
B^{m+1}_s(\alpha_s)=&\;(\nabla^{\man{E}}B^m_s)(\alpha_s)+
B^m_{s-1}\otimes\id_{\ctb{\man{E}}}(\alpha_s)+
\sum_{j=1}^m\Ins_j(B^m_s(\alpha_s),B_{\pi_{\man{E}}})\\
&\;-\Ins_{m+1}(B^m_s(\alpha_s),\dual{B}_{\pi_{\man{E}}}),\enspace s\in\{1,\dots,m\},\\
B^{m+1}_0(\alpha_0)=&\;(\nabla^{\man{E}}B^m_0)(\alpha_0)+
\sum_{j=1}^m\Ins_j(B^m_0(\alpha_0),B_{\pi_{\man{E}}})-
\Ins_{m+1}(B^m_0(\alpha_0),\dual{B}_{\pi_{\man{E}}}),
\end{align*}
where\/ $\alpha_s\in\tensor*[s]{\ctb{\man{E}}}\otimes\vb{\man{E}}$\@,\/
$s\in\{0,1,\dots,m+1\}$\@.
\begin{proof}
The assertion clearly holds for $m=0$\@, so suppose it true for
$m\in\integerp$\@.  Thus
\begin{equation}\label{eq:ViteratedderII1}
\vl{(\nabla^{\man{M},\pi_{\man{E}},m}\xi)}=
\sum_{s=0}^mB^m_s(\nabla^{\man{E},s}\vl{\xi}),
\end{equation}
where the vector bundle mappings $B^a_s$\@, $a\in\{0,1,\dots,m\}$\@,
$s\in\{0,1,\dots,a\}$\@, satisfy the recursion relations from the statement
of the lemma.  Then, by Lemma~\pldblref{lem:bundleform}{pl:bundleform2}\@, we
can work on the left-hand side of~\eqref{eq:ViteratedderII1} to give
{\small\begin{align*}
\nabla^{\man{E}}\vl{(\nabla^{\man{M},\pi_{\man{E}},m}\xi)}=&\;
\vl{(\nabla^{\man{M},\pi_{\man{E}},m+1}\xi)}-\sum_{j=1}^m
\Ins_j(\vl{(\nabla^{\man{M},\pi_{\man{E}},m}\xi)},B_{\pi_{\man{E}}})+
\Ins_{m+1}(\vl{(\nabla^{\man{M},\pi_{\man{E}},m}\xi)},\dual{B}_{\pi_{\man{E}}})\\
=&\;\vl{(\nabla^{\man{M},\pi_{\man{E}},m+1}\xi)}-\sum_{s=0}^m\sum_{j=1}^m
\Ins_j(B^m_s(\nabla^{\man{E},s}\vl{\xi}),B_{\pi_{\man{E}}})+
\sum_{s=0}^m\Ins_{m+1}(B^m_s(\nabla^{\man{E},s}\vl{\xi}),
\dual{B}_{\pi_{\man{E}}}).
\end{align*}}%
Working on the right-hand side of~\eqref{eq:ViteratedderII1} gives
\begin{equation*}
\nabla^{\man{E}}\vl{(\nabla^{\man{M},\pi_{\man{E}},m}\xi)}=
\sum_{s=0}^m\nabla^{\man{E}}B^m_s(\nabla^{\man{E},s}\vl{\xi})+
\sum_{s=0}^mB^m_s\otimes\id_{\ctb{\man{E}}}(\nabla^{\man{E},s+1}\vl{\xi}).
\end{equation*}
Combining the preceding two equations gives
\begin{align*}
\nabla^{\man{M},\pi_{\man{E}},m+1}\vl{\xi}=&\;
\sum_{s=0}^m\nabla^{\man{E}}B^m_s(\nabla^{\man{E},s}\vl{\xi})+
\sum_{s=0}^mB^m_s\otimes\id_{\ctb{\man{E}}}(\nabla^{\man{E},s+1}\vl{\xi})\\
&\;+\sum_{s=0}^m\sum_{j=1}^m\Ins_j(B^m_s(\nabla^{\man{E},s}\vl{\xi}),
B_{\pi_{\man{E}}})-\Ins_{m+1}(\vl{(\nabla^{\man{M},\pi_{\man{E}},m}\xi)},
\dual{B}_{\pi_{\man{E}}})\\
=&\;\nabla^{\man{E},m+1}\vl{\xi}+
\left.\sum_{s=1}^m\right(\nabla^{\man{E}}B^m_s(\nabla^{\man{E},s}\vl{\xi})+
B^m_{s-1}\otimes\id_{\ctb{\man{E}}}(\nabla^{\man{E},s}\vl{\xi})\\
&\;\left.+\sum_{j=1}^m\Ins_j(B^m_s(\nabla^{\man{E},s}\vl{\xi}),
B_{\pi_{\man{E}}})-\Ins_{m+1}(B^m_s(\nabla^{\man{E},s}\vl{\xi}),
\dual{B}_{\pi_{\man{E}}})\right)\\
&\;+\nabla^{\man{E}}B^m_0(\vl{\xi})+\sum_{j=1}^m\Ins_j(B^m_0(\vl{\xi}),
B_{\pi_{\man{E}}})-\Ins_{m+1}(B^m_0(\vl{\xi}),\dual{B}_{\pi_{\man{E}}}),
\end{align*}
and the lemma follows from this.
\end{proof}
\end{lemma}

Next we turn to symmetrised versions of the preceding lemmata.  We show that
the preceding two lemmata induce corresponding mappings between symmetric
tensors.
\begin{lemma}\label{lem:VsymiterateddersI}
Let\/ $r\in\{\infty,\omega\}$ and let\/
$\map{\pi_{\man{E}}}{\man{E}}{\man{M}}$ be a\/ $\C^r$-vector bundle, with the
data prescribed in Section~\ref{subsec:Esubmersion} to define the Riemannian
metric\/ $\metric_{\man{E}}$ on\/ $\man{E}$\@.  For\/ $m\in\integernn$\@,
there exist\/ $\C^r$-vector bundle mappings
\begin{equation*}
(\what{A}^m_s,\id_{\man{E}})\in
\vbmappings[r]{\Symalg*[s]{\pi_{\man{E}}^*\ctb{\man{M}}}\otimes\vb{\man{E}}}
{\Symalg*[m]{\ctb{\man{E}}}\otimes\vb{\man{E}}},\qquad s\in\{0,1,\dots,m\},
\end{equation*}
such that
\begin{equation*}
(\Sym_m\otimes\id_{\vb{\man{E}}})\scirc\nabla^{\man{E},m}\vl{\xi}=
\sum_{s=0}^m\what{A}^m_s((\Sym_s\otimes\id_{\vb{\man{E}}})
\scirc\vl{(\nabla^{\man{M},\pi_{\man{E}},s}\xi)})
\end{equation*}
for all\/ $\xi\in\sections[m]{\man{E}}$\@.
\begin{proof}
The proof follows very similarly to that of
Lemma~\ref{lem:PsymiterateddersI}\@, but taking the tensor product of
everything with $\vb{\man{E}}$\@.  We shall present the complete construction
here, but will not repeat it for similar proofs that follow.

We define $\map{A^m}{\tensor*[\le m]{\pi_{\man{E}}^*\ctb{\man{M}}}\otimes\vb{\man{E}}}
{\tensor*[\le m]{\ctb{\man{E}}}\otimes\vb{\man{E}}}$ by
\begin{multline*}
A^m(\vl{\xi},\vl{(\nabla^{\pi_{\man{E}}}\xi)},\dots,
\vl{(\nabla^{\man{M},\pi_{\man{E}},m}\xi)})\\
=\left(A^0_0(\vl{\xi}),\sum_{s=0}^1A^1_s
(\vl{(\nabla^{\man{M},\pi_{\man{E}},s}\xi)}),
\dots,\sum_{s=0}^mA^m_s(\vl{(\nabla^{\man{M},\pi_{\man{E}},s}\xi)})\right)
\end{multline*}
Let us organise the mappings we require into the following diagram:
\begin{equation}\label{eq:Vsymmetrise1}
\xymatrix{{\tensor*[\le m]{\pi_{\man{E}}^*\ctb{\man{M}}}\otimes\vb{\man{E}}}
\ar[rr]^{\Sym_{\le m}\otimes\id_{\vb{\man{E}}}}\ar[d]_{A^m}&&
{\Symalg*[\le m]{\pi_{\man{E}}^*\ctb{\man{M}}}\otimes\vb{\man{E}}}
\ar[rr]^(0.5){S^m_{\nabla^{\man{M}},\nabla^{\pi_{\man{E}}}}
\otimes\id_{\vb{\man{E}}}}\ar[d]^{\what{A}^m}&&
{\pi_{\man{E}}^*(\real_{\man{M}}\oplus\jetalg{m}{\man{M}})\otimes\vb{\man{E}}}
\ar[d]^{(\id_{\real}\oplus j_m\pi_{\man{E}})\otimes\id_{\vb{\man{E}}}}\\
{\tensor*[\le m]{\ctb{\man{E}}}\otimes\vb{\man{E}}}
\ar[rr]^{\Sym_{\le m}\otimes\id_{\vb{\man{E}}}}&&
{\Symalg*[\le m]{\ctb{\man{E}}}\otimes\vb{\man{E}}}
\ar[rr]^(0.5){S^m_{\nabla^{\man{E}}}\otimes\id_{\vb{\man{E}}}}&&
{(\real_{\man{M}}\oplus\jetalg{m}{\man{E}})\otimes\vb{\man{E}}}}
\end{equation}
Here $\what{A}^m$ is defined so that the right square commutes.  We shall
show that the left square also commutes.  Indeed,
\begin{align*}
\what{A}^m\scirc\Sym_{\le m}\otimes\id_{\vb{\man{E}}}&
(\vl{\xi},\vl{(\nabla^{\pi_{\man{E}}}\xi)},\dots,
\vl{(\nabla^{\man{M},\pi_{\man{E}},m}\xi)})\\
=&\;(S^m_{\nabla^{\man{E}}}\otimes\id_{\vb{\man{E}}})^{-1}\scirc
((\id_{\real}\oplus j_m\pi_{\man{E}})\otimes\id_{\vb{\man{E}}})\scirc
(S^m_{\nabla^{\man{M}},\nabla^{\pi_{\man{E}}}}\otimes\id_{\vb{\man{E}}})\\
&\;\scirc(\Sym_{\le m}\otimes\id_{\vb{\man{E}}})
(\vl{\xi},\vl{(\nabla^{\pi_{\man{E}}}\xi)},\dots,
\vl{(\nabla^{\man{M},\pi_{\man{E}},m}\xi)})\\
=&\;\Sym_{\le m}\otimes\id_{\vb{\man{E}}}
(\vl{\xi},\nabla^{\man{E}}\vl{\xi},\dots,\nabla^{\man{E},m}\vl{\xi})\\
=&\;(\Sym_{\le m}\otimes\id_{\vb{\man{E}}})\scirc
A^m(\vl{\xi},\vl{(\nabla^{\pi_{\man{E}}}\xi)},\dots,
\vl{(\nabla^{\man{M},\pi_{\man{E}},m}\xi)}).
\end{align*}
Thus the diagram~\eqref{eq:Vsymmetrise1} commutes.  Thus, if we define
\begin{equation*}
\what{A}^m_s((\Sym_s\otimes\id_{\vb{\man{E}}})\scirc
\vl{(\nabla^{\man{M},\pi_{\man{E}},s}\xi)})=
(\Sym_m\otimes\id_{\vb{\man{E}}})\scirc A^m_s(\vl{(\nabla^{\man{M},\pi_{\man{E}},s}\xi)}),
\end{equation*}
then we have
\begin{equation*}
(\Sym_m\otimes\id_{\vb{\man{E}}})\scirc\nabla^{\man{E},m}\vl{\xi}=
\sum_{s=0}^m\what{A}^m_s((\Sym_s\otimes\id_{\vb{\man{E}}})\scirc
\vl{(\nabla^{\man{M},\pi_{\man{E}},s}\xi)}),
\end{equation*}
as desired.
\end{proof}
\end{lemma}

The preceding lemma gives rise to an ``inverse,'' which we state in the
following lemma.
\begin{lemma}\label{lem:VsymiterateddersII}
Let\/ $r\in\{\infty,\omega\}$ and let\/
$\map{\pi_{\man{E}}}{\man{E}}{\man{M}}$ be a\/ $\C^r$-vector bundle, with the
data prescribed in Section~\ref{subsec:Esubmersion} to define the Riemannian
metric\/ $\metric_{\man{E}}$ on\/ $\man{E}$\@.  For\/ $m\in\integernn$\@,
there exist\/ $\C^r$-vector bundle mappings
\begin{equation*}
(\what{B}^m_s,\id_{\man{E}})\in
\vbmappings[r]{\Symalg*[s]{\ctb{\man{E}}}\otimes\vb{\man{E}}}
{\Symalg*[m]{\pi_{\man{E}}^*\ctb{\man{M}}}\otimes\vb{\man{E}}},\qquad
s\in\{0,1,\dots,m\},
\end{equation*}
such that
\begin{equation*}
(\Sym_m\otimes_{\id_{\vb{\man{E}}}})\scirc
\vl{(\nabla^{\man{M},\pi_{\man{E}},m}\xi)}=
\sum_{s=0}^m\what{B}^m_s((\Sym_s\otimes\id_{\vb{\man{E}}})\scirc
\nabla^{\man{E},s}\vl{\xi})
\end{equation*}
for all\/ $\xi\in\sections[m]{\man{E}}$\@.
\begin{proof}
This follows along the lines of Lemma~\ref{lem:PsymiterateddersII} in the
same manner as Lemma~\ref{lem:VsymiterateddersI} follows from
Lemma~\ref{lem:PsymiterateddersI}\@, by taking tensor products with
$\vb{\man{E}}$\@.
\end{proof}
\end{lemma}

We can put together the previous four lemmata into the following
decomposition result, which is to be regarded as the main result of this
section.
\begin{lemma}\label{lem:Vdecomp}
Let\/ $r\in\{\infty,\omega\}$ and let\/
$\map{\pi_{\man{E}}}{\man{E}}{\man{M}}$ be a\/ $\C^r$-vector bundle, with the
data prescribed in Section~\ref{subsec:Esubmersion} to define the Riemannian
metric\/ $\metric_{\man{E}}$ on\/ $\man{E}$\@.  Then there exist\/
$\C^r$-vector bundle mappings
\begin{equation*}
A^m_{\nabla^{\man{E}}}\in\vbmappings[r]{\Pjetalg{m}{\man{E}}\otimes\vb{\man{E}}}
{\Symalg*[\le m]{\pi_{\man{E}}^*\ctb{\man{M}}}\otimes\vb{\man{E}}},\quad
B^m_{\nabla^{\man{E}}}\in\vbmappings[r]{\Pjetalg{m}{\man{E}}\otimes\vb{\man{E}}}
{\Symalg*[\le m]{\ctb{\man{E}}}\otimes\vb{\man{E}}},
\end{equation*}
defined by
\begin{align*}
A^m_{\nabla^{\man{E}}}(j_m(\vl{\xi})(e))=&\;
\Sym_{\le m}\otimes\id_{\vb{\man{E}}}(\vl{\xi}(e),
\vl{(\nabla^{\pi_{\man{E}}}\xi)}(e),
\dots,\vl{(\nabla^{\man{M},\pi_{\man{E}},m}\xi)}(e)),\\
B^m_{\nabla^{\man{E}}}(j_m(\vl{\xi})(e))=&\;
\Sym_{\le m}\otimes\id_{\vb{\man{E}}}
(\vl{\xi}(e),\nabla^{\man{E}}\vl{\xi}(e),\dots,
\nabla^{\man{E},m}\vl{\xi}(e)).
\end{align*}
Moreover,\/ $A^m_{\nabla^{\man{E}}}$ is an isomorphism,\/
$B^m_{\nabla^{\man{E}}}$ is injective, and
{\small\begin{multline*}
B^m_{\nabla^{\man{E}}}\scirc(A^m_{\nabla^{\man{E}}})^{-1}\scirc
(\Sym_{\le m}\otimes\id_{\vb{\man{E}}})
(\vl{\xi}(e),\vl{(\nabla^{\pi_{\man{E}}}\xi)}(e),\dots,
\vl{(\nabla^{\man{M},\pi_{\man{E}},m}\xi)}(e))\\
=\left(\vl{\xi}(e),\sum_{s=0}^1\what{A}^1_s((\Sym_s\otimes\id_{\vb{\man{E}}})
\scirc\vl{(\nabla^{\man{M},\pi_{\man{E}},s}\xi)}(e)),\dots,
\sum_{s=0}^m\what{A}^m_s((\Sym_s\otimes\id_{\vb{\man{E}}})
\scirc\vl{(\nabla^{\man{M},\pi_{\man{E}},s}\xi)}(e))\right)
\end{multline*}}%
and
{\small\begin{multline*}
A^m_{\nabla^{\man{E}}}\scirc(B^m_{\nabla^{\man{E}}})^{-1}\scirc
(\Sym_{\le m}\otimes\id_{\vb{\man{E}}})
(\vl{\xi}(e),\nabla^{\man{E}}\vl{\xi}(e),\dots,\nabla^{\man{E},m}\vl{\xi}(e))\\
=\left(\vl{\xi}(e),\sum_{s=0}^1\what{B}^1_s((\Sym_s\otimes\id_{\vb{\man{E}}})
\scirc\nabla^{\man{E},s}\vl{\xi}(e)),\dots,
\sum_{s=0}^m\what{B}^m_s((\Sym_s\otimes\id_{\vb{\man{E}}})\scirc
\nabla^{\man{E},s}\vl{\xi}(e))\right),
\end{multline*}}%
where the vector bundle mappings\/ $\what{A}^m_s$ and\/ $\what{B}^m_s$\@,\/
$s\in\{0,1,\dots,m\}$\@, are as in Lemmata~\ref{lem:VsymiterateddersI}
and~\ref{lem:VsymiterateddersII}\@.
\end{lemma}

\subsection{Isomorphisms for horizontal lifts of vector fields}

Next we consider horizontal lifts of vector fields via the mapping
\begin{equation*}
\sections[r]{\tb{\man{M}}}\ni X\mapsto\hl{X}\in\sections[r]{\tb{\man{E}}}.
\end{equation*}
We wish to relate the decomposition of the jets of $X$ with the jets of
$\hl{X}$\@.  Associated with this, we denote
\begin{equation*}
\Hjetalg[e]{m}{\man{E}}=\setdef{j_m\hl{X}(e)}{X\in\sections[m]{\tb{\man{M}}}}.
\end{equation*}
By~\eqref{eq:jet=jetalg}\@, we have
\begin{equation*}
\Hjetalg[e]{m}{\man{E}}\simeq\Pjetalg[e]{m}{\man{E}}\otimes\hb[e]{\man{E}}.
\end{equation*}
As with the constructions of the preceding sections, we wish to use
Lemma~\ref{lem:Jmdecomp} to provide a decomposition of
$\Hjetalg{m}{\man{E}}$\@, and to do so we need to understand the covariant
derivatives
\begin{equation*}
\nabla^{\man{E},m}\hl{X}\eqdef\underbrace{\nabla^{\man{E}}\cdots
\nabla^{\man{E}}}_{m\ \textrm{times}}\hl{X},\qquad m\in\integernn.
\end{equation*}
In this section we omit proofs, since proofs follow along entirely similar
lines to those of the preceding section.

The first result we give is the following.
\begin{lemma}\label{lem:HiterateddersI}
Let\/ $r\in\{\infty,\omega\}$ and let\/
$\map{\pi_{\man{E}}}{\man{E}}{\man{M}}$ be a\/ $\C^r$-vector bundle, with the
data prescribed in Section~\ref{subsec:Esubmersion} to define the Riemannian
metric\/ $\metric_{\man{E}}$ on\/ $\man{E}$\@.  For\/ $m\in\integernn$\@,
there exist\/ $\C^r$-vector bundle mappings
\begin{equation*}
(A^m_s,\id_{\man{E}})\in
\vbmappings[r]{\tensor*[s]{\pi_{\man{E}}^*\ctb{\man{M}}}\otimes\hb{\man{E}}}
{\tensor*[m]{\ctb{\man{E}}}\otimes\hb{\man{E}}},\qquad
s\in\{0,1,\dots,m\},
\end{equation*}
such that
\begin{equation*}
\nabla^{\man{E},m}\hl{X}=\sum_{s=0}^mA^m_s(\hl{(\nabla^{\man{M},s}X)})
\end{equation*}
for all\/ $X\in\sections[m]{\tb{\man{M}}}$\@.  Moreover, the vector bundle
mappings\/ $A^m_0,A^m_1,\dots,A^m_m$ satisfy the recursion relations
prescribed by\/ $A^0_0(\beta_0)=\beta_0$ and
\begin{align*}
A^{m+1}_{m+1}(\beta_{m+1})=&\;\beta_{m+1},\\
A^{m+1}_s(\beta_s)=&\;(\nabla^{\man{E}}A^m_s)(\beta_s)+
A^m_{s-1}\otimes\id_{\ctb{\man{E}}}(\beta_s)-
\sum_{j=1}^sA^m_s\otimes\id_{\ctb{\man{E}}}(\Ins_j(\beta_s,B_{\pi_{\man{E}}}))\\
&\;+A^m_s\otimes\id_{\ctb{\man{E}}}(\Ins_{s+1}(\beta_s,\dual{B}_{\pi_{\man{E}}})),
\enspace s\in\{1,\dots,m\},\\
A^{m+1}_0(\beta_0)=&\;(\nabla^{\man{E}}A^m_0)(\beta_0)+
A^m_0\otimes\id_{\ctb{\man{E}}}(\Ins_1(\beta_0,\dual{B}_{\pi_{\man{E}}})),
\end{align*}
where\/ $\beta_s\in\tensor*[s]{\pi_{\man{E}}^*\ctb{\man{M}}}\otimes
\hb{\man{E}}$\@,\/ $s\in\{0,1,\dots,m+1\}$\@.
\begin{proof}
This follows in the same manner as Lemma~\ref{lem:ViterateddersI}\@, making
use of Lemma~\pldblref{lem:bundleform}{pl:bundleform3}\@.
\end{proof}
\end{lemma}

The following lemma ``inverts'' the relations from the preceding one.
\begin{lemma}\label{lem:HiterateddersII}
Let\/ $r\in\{\infty,\omega\}$ and let\/
$\map{\pi_{\man{E}}}{\man{E}}{\man{M}}$ be a\/ $\C^r$-vector bundle, with the
data prescribed in Section~\ref{subsec:Esubmersion} to define the Riemannian
metric\/ $\metric_{\man{E}}$ on\/ $\man{E}$\@.  For\/ $m\in\integernn$\@,
there exist\/ $\C^r$-vector bundle mappings
\begin{equation*}
(B^m_s,\id_{\man{E}})\in
\vbmappings[r]{\tensor*[s]{\ctb{\man{E}}}\otimes\hb{\man{E}}}
{\tensor*[m]{\pi_{\man{E}}^*\ctb{\man{M}}}\otimes\hb{\man{E}}},\qquad
s\in\{0,1,\dots,m\},
\end{equation*}
such that
\begin{equation*}
\hl{(\nabla^{\man{M},m}X)}=\sum_{s=0}^mB^m_s(\nabla^{\man{E},s}\hl{X})
\end{equation*}
for all\/ $X\in\sections[m]{\tb{\man{M}}}$\@.  Moreover, the vector bundle
mappings\/ $B^m_0,B^m_1,\dots,B^m_m$ satisfy the recursion relations
prescribed by\/ $B^0_0(\alpha_0)=\alpha_0$ and
\begin{align*}
B^{m+1}_{m+1}(\alpha_{m+1})=&\;\alpha_{m+1},\\
B^{m+1}_s(\alpha_s)=&\;(\nabla^{\man{E}}B^m_s)(\alpha_s)+
B^m_{s-1}\otimes\id_{\ctb{\man{E}}}(\alpha_s)+
\sum_{j=1}^m\Ins_j(B^m_s(\alpha_s),B_{\pi_{\man{E}}})\\
&\;-\Ins_{m+1}(B^m_s(\alpha_s),\dual{B}_{\pi_{\man{E}}}),\enspace s\in\{1,\dots,m\},\\
B^{m+1}_0(\alpha_0)=&\;(\nabla^{\man{E}}B^m_0)(\alpha_0)+
\sum_{j=1}^m\Ins_j(B^m_0(\alpha_0),B_{\pi_{\man{E}}})-
\Ins_{m+1}(B^m_0(\alpha_0),\dual{B}_{\pi_{\man{E}}}),
\end{align*}
where\/ $\alpha_s\in\tensor*[s]{\ctb{\man{E}}}\otimes\hb{\man{E}}$\@,\/
$s\in\{0,1,\dots,m+1\}$\@.
\begin{proof}
This follows in the same manner as Lemma~\ref{lem:ViterateddersII}\@, making
use of Lemma~\pldblref{lem:bundleform}{pl:bundleform3}\@.
\end{proof}
\end{lemma}

Now we can give the symmetrised versions of the preceding lemmata.
\begin{lemma}\label{lem:HsymiterateddersI}
Let\/ $r\in\{\infty,\omega\}$ and let\/
$\map{\pi_{\man{E}}}{\man{E}}{\man{M}}$ be a\/ $\C^r$-vector bundle, with the
data prescribed in Section~\ref{subsec:Esubmersion} to define the Riemannian
metric\/ $\metric_{\man{E}}$ on\/ $\man{E}$\@.  For\/ $m\in\integernn$\@,
there exist\/ $\C^r$-vector bundle mappings
\begin{equation*}
(\what{A}^m_s,\id_{\man{E}})\in
\vbmappings[r]{\Symalg*[s]{\pi_{\man{E}}^*\ctb{\man{M}}}\otimes\hb{\man{E}}}
{\Symalg*[m]{\ctb{\man{E}}}\otimes\hb{\man{E}}},\qquad s\in\{0,1,\dots,m\},
\end{equation*}
such that
\begin{equation*}
(\Sym_m\otimes\id_{\hb{\man{E}}})\scirc\nabla^{\man{E},m}\hl{X}=
\sum_{s=0}^m\what{A}^m_s((\Sym_s\otimes\id_{\hb{\man{E}}})
\scirc(\hl{\nabla^{\man{M},s}X)})
\end{equation*}
for all\/ $X\in\sections[m]{\tb{\man{M}}}$\@.
\begin{proof}
This follows along the lines of Lemma~\ref{lem:PsymiterateddersI} in the same
manner as Lemma~\ref{lem:VsymiterateddersI} follows from
Lemma~\ref{lem:PsymiterateddersI}\@, by taking tensor products with
$\hb{\man{E}}$\@.
\end{proof}
\end{lemma}

\begin{lemma}\label{lem:HsymiterateddersII}
Let\/ $r\in\{\infty,\omega\}$ and let\/
$\map{\pi_{\man{E}}}{\man{E}}{\man{M}}$ be a\/ $\C^r$-vector bundle, with the
data prescribed in Section~\ref{subsec:Esubmersion} to define the Riemannian
metric\/ $\metric_{\man{E}}$ on\/ $\man{E}$\@.  For\/ $m\in\integernn$\@,
there exist\/ $\C^r$-vector bundle mappings
\begin{equation*}
(\what{B}^m_s,\id_{\man{E}})\in
\vbmappings[r]{\Symalg*[s]{\ctb{\man{E}}}\otimes\hb{\man{E}}}
{\Symalg*[m]{\pi_{\man{E}}^*\ctb{\man{M}}}\otimes\hb{\man{E}}},\qquad
s\in\{0,1,\dots,m\},
\end{equation*}
such that
\begin{equation*}
(\Sym_m\otimes_{\id_{\hb{\man{E}}}})\scirc\hl{(\nabla^{\man{M},m}X)}=
\sum_{s=0}^m\what{B}^m_s((\Sym_s\otimes\id_{\hb{\man{E}}})\scirc
\nabla^{\man{E},s}\hl{X})
\end{equation*}
for all\/ $X\in\sections[m]{\tb{\man{M}}}$\@.
\begin{proof}
This follows along the lines of Lemma~\ref{lem:PsymiterateddersII} in the
same manner as Lemma~\ref{lem:VsymiterateddersI} follows from
Lemma~\ref{lem:PsymiterateddersI}\@, by taking tensor products with
$\hb{\man{E}}$\@.
\end{proof}
\end{lemma}

We can put together the previous four lemmata into the following
decomposition result, which is to be regarded as the main result of this
section.
\begin{lemma}\label{lem:Hdecomp}
Let\/ $r\in\{\infty,\omega\}$ and let\/
$\map{\pi_{\man{E}}}{\man{E}}{\man{M}}$ be a\/ $\C^r$-vector bundle, with the
data prescribed in Section~\ref{subsec:Esubmersion} to define the Riemannian
metric\/ $\metric_{\man{E}}$ on\/ $\man{E}$\@.  Then there exist\/
$\C^r$-vector bundle mappings
\begin{equation*}
A^m_{\nabla^{\man{E}}}\in\vbmappings[r]{\Pjetalg{m}{\man{E}}\otimes\hb{\man{E}}}
{\Symalg*[\le m]{\pi_{\man{E}}^*\ctb{\man{M}}}\otimes\hb{\man{E}}},\quad
B^m_{\nabla^{\man{E}}}\in\vbmappings[r]{\Pjetalg{m}{\man{E}}\otimes\hb{\man{E}}}
{\Symalg*[\le m]{\ctb{\man{E}}}\otimes\hb{\man{E}}},
\end{equation*}
defined by
\begin{align*}
A^m_{\nabla^{\man{E}}}(j_m(\hl{X})(e))=&\;
\Sym_{\le m}\otimes\id_{\hb{\man{E}}}(\hl{X}(e),\hl{(\nabla^{\man{M}}X)}(e),
\dots,\hl{(\nabla^{\man{M},m}X)}(e)),\\
B^m_{\nabla^{\man{E}}}(j_m(\hl{X})(e))=&\;
\Sym_{\le m}\otimes\id_{\hb{\man{E}}}(\hl{X}(e),\nabla^{\man{E}}\hl{X}(e),
\dots,\nabla^{\man{E},m}\hl{X}(e)).
\end{align*}
Moreover,\/ $A^m_{\nabla^{\man{E}}}$ is an isomorphism,\/
$B^m_{\nabla^{\man{E}}}$ is injective, and
{\small\begin{multline*}
B^m_{\nabla^{\man{E}}}\scirc(A^m_{\nabla^{\man{E}}})^{-1}\scirc
(\Sym_{\le m}\otimes_{\id_{\hb{\man{E}}}})(\hl{X}(e),\hl{(\nabla^{\man{M}}X)}(e),\dots,
\hl{(\nabla^{\man{M},m}X)}(e))\\
=\left(\hl{X}(e),\sum_{s=0}^1\what{A}^1_s((\Sym_s\otimes\id_{\hb{\man{E}}})
\scirc({\nabla^{\man{M},s}X)}(e)),\dots,
\sum_{s=0}^m\what{A}^m_s((\Sym_s\otimes\id_{\hb{\man{E}}})
\scirc\hl{(\nabla^{\man{M},s}X)}(e))\right)
\end{multline*}}%
and
{\small\begin{multline*}
A^m_{\nabla^{\man{E}}}\scirc(B^m_{\nabla^{\man{E}}})^{-1}\scirc
(\Sym_{\le m}\otimes\id_{\hb{\man{E}}})
(\hl{X}(e),\nabla^{\man{E}}\hl{X}(e),\dots,\nabla^{\man{E},m}\hl{X}(e))\\
=\left(\hl{X}(e),\sum_{s=0}^1\what{B}^1_s((\Sym_s\otimes\id_{\hb{\man{E}}})
\scirc\nabla^{\man{E},s}\hl{X}(e)),\dots,
\sum_{s=0}^m\what{B}^m_s((\Sym_s\otimes\id_{\hb{\man{E}}})\scirc
\nabla^{\man{E},s}\hl{X}(e))\right),
\end{multline*}}%
where the vector bundle mappings\/ $\what{A}^m_s$ and\/ $\what{B}^m_s$\@,\/
$s\in\{0,1,\dots,m\}$\@, are as in Lemmata~\ref{lem:HsymiterateddersI}
and~\ref{lem:HsymiterateddersII}\@.
\end{lemma}

\subsection{Isomorphisms for vertical lifts of dual sections}

Next we consider vertical lifts of sections of the dual bundle,~\ie~the
mapping defined by
\begin{equation*}
\sections[r]{\dual{\man{E}}}\ni\lambda\mapsto
\vl{\lambda}\in\sections[r]{\ctb{\man{E}}}.
\end{equation*}
Our objective is to relate the decomposition of the jets of $\lambda$ with
the decomposition of the jets of $\vl{\lambda}$\@.  To do this, we denote
\begin{equation*}
\Fjetalg[e]{m}{\man{E}}=\setdef{j_m\vl{\lambda}(e)}
{\lambda\in\sections[m]{\dual{\man{E}}}}.
\end{equation*}
By~\eqref{eq:jet=jetalg}\@, we have
\begin{equation*}
\Fjetalg[e]{m}{\man{E}}\simeq\Pjetalg[e]{m}{\man{E}}\otimes\cvb[e]{\man{E}}.
\end{equation*}
As with the constructions of the preceding sections, we wish to use
Lemma~\ref{lem:Jmdecomp} to provide a decomposition of
$\Fjetalg{m}{\man{E}}$\@, and to do so we need to understand the covariant
derivatives
\begin{equation*}
\nabla^{\man{E},m}\vl{\lambda}\eqdef\underbrace{\nabla^{\man{E}}\cdots
\nabla^{\man{E}}}_{m\ \textrm{times}}\vl{\lambda},\qquad m\in\integernn.
\end{equation*}
In this section we omit proofs, since proofs follow along entirely similar
lines to those of the preceding section.

The first result we give is the following.
\begin{lemma}\label{lem:V*iterateddersI}
Let\/ $r\in\{\infty,\omega\}$ and let\/
$\map{\pi_{\man{E}}}{\man{E}}{\man{M}}$ be a\/ $\C^r$-vector bundle, with the
data prescribed in Section~\ref{subsec:Esubmersion} to define the Riemannian
metric\/ $\metric_{\man{E}}$ on\/ $\man{E}$\@.  For\/ $m\in\integernn$\@,
there exist\/ $\C^r$-vector bundle mappings
\begin{equation*}
(A^m_s,\id_{\man{E}})\in
\vbmappings[r]{\tensor*[s]{\pi_{\man{E}}^*\ctb{\man{M}}}\otimes\cvb{\man{E}}}
{\tensor*[m]{\ctb{\man{E}}}\otimes\cvb{\man{E}}},\qquad
s\in\{0,1,\dots,m\},
\end{equation*}
such that
\begin{equation*}
\nabla^{\man{E},m}\vl{\lambda}=
\sum_{s=0}^mA^m_s(\vl{(\nabla^{\man{M},\pi_{\man{E}},s}\lambda)})
\end{equation*}
for all\/ $\lambda\in\sections[m]{\dual{\man{E}}}$\@.  Moreover, the vector
bundle mappings\/ $A^m_0,A^m_1,\dots,A^m_m$ satisfy the recursion relations
prescribed by\/ $A^0_0(\beta_0)=\beta_0$ and
{\small\begin{align*}
A^{m+1}_{m+1}(\beta_{m+1})=&\;\beta_{m+1},\\
A^{m+1}_s(\beta_s)=&\;(\nabla^{\man{E}}A^m_s)(\beta_s)+
A^m_{s-1}\otimes\id_{\ctb{\man{E}}}(\beta_s)-
\sum_{j=1}^sA^m_s\otimes\id_{\ctb{\man{E}}}(\Ins_j(\beta_s,B_{\pi_{\man{E}}})),
\enspace s\in\{1,\dots,m\},\\
A^{m+1}_0(\beta_0)=&\;(\nabla^{\man{E}}A^m_0)(\beta_0)-
A^m_0\otimes\id_{\ctb{\man{E}}}(\Ins_1(\beta_0,B_{\pi_{\man{E}}})),
\end{align*}}%
where\/ $\beta_s\in\tensor*[s]{\pi_{\man{E}}^*\ctb{\man{M}}}\otimes
\cvb{\man{E}}$\@,\/ $s\in\{0,1,\dots,m+1\}$\@.
\begin{proof}
This follows in the same manner as Lemma~\ref{lem:ViterateddersI}\@, making
use of Lemma~\pldblref{lem:bundleform}{pl:bundleform4}\@.
\end{proof}
\end{lemma}

The ``inverse'' of the preceding lemma is as follows.
\begin{lemma}\label{lem:V*iterateddersII}
Let\/ $r\in\{\infty,\omega\}$ and let\/
$\map{\pi_{\man{E}}}{\man{E}}{\man{M}}$ be a\/ $\C^r$-vector bundle, with the
data prescribed in Section~\ref{subsec:Esubmersion} to define the Riemannian
metric\/ $\metric_{\man{E}}$ on\/ $\man{E}$\@.  For\/ $m\in\integernn$\@,
there exist\/ $\C^r$-vector bundle mappings
\begin{equation*}
(B^m_s,\id_{\man{E}})\in
\vbmappings[r]{\tensor*[s]{\ctb{\man{E}}}\otimes\cvb{\man{E}}}
{\tensor*[m]{\pi_{\man{E}}^*\ctb{\man{M}}}\otimes\cvb{\man{E}}},\qquad s\in\{0,1,\dots,m\},
\end{equation*}
such that
\begin{equation*}
\vl{(\nabla^{\man{M},\pi_{\man{E}},m}\lambda)}=
\sum_{s=0}^mB^m_s(\nabla^{\man{E},s}\vl{\lambda})
\end{equation*}
for all\/ $\lambda\in\sections[m]{\dual{\man{E}}}$\@.  Moreover, the vector
bundle mappings\/ $B^m_0,B^m_1,\dots,B^m_m$ satisfy the recursion relations
prescribed by\/ $B^0_0(\alpha_0)=\alpha_0$ and
\begin{align*}
B^{m+1}_{m+1}(\alpha_{m+1})=&\;\alpha_{m+1},\\
B^{m+1}_s(\alpha_s)=&\;(\nabla^{\man{E}}B^m_s)(\alpha_s)+
B^m_{s-1}\otimes\id_{\ctb{\man{E}}}(\alpha_s)+
\sum_{j=1}^m\Ins_j(B^m_s(\alpha_s),B_{\pi_{\man{E}}}),\enspace s\in\{1,\dots,m\},\\
B^{m+1}_0(\alpha_0)=&\;(\nabla^{\man{E}}B^m_0)(\alpha_0)+
\sum_{j=1}^{m+1}\Ins_j(B^m_0(\alpha_0),B_{\pi_{\man{E}}}),
\end{align*}
where\/ $\alpha_s\in\tensor*[s]{\ctb{\man{E}}}\otimes\cvb{\man{E}}$\@,\/
$s\in\{0,1,\dots,m+1\}$\@.
\begin{proof}
This follows in the same manner as Lemma~\ref{lem:ViterateddersII}\@, making
use of Lemma~\pldblref{lem:bundleform}{pl:bundleform4}\@.
\end{proof}
\end{lemma}

Next we turn to symmetrised versions of the preceding lemmata.  We show that
the preceding two lemmata induce corresponding mappings between symmetric
tensors.
\begin{lemma}\label{lem:V*symiterateddersI}
Let\/ $r\in\{\infty,\omega\}$ and let\/
$\map{\pi_{\man{E}}}{\man{E}}{\man{M}}$ be a\/ $\C^r$-vector bundle, with the
data prescribed in Section~\ref{subsec:Esubmersion} to define the Riemannian
metric\/ $\metric_{\man{E}}$ on\/ $\man{E}$\@.  For\/ $m\in\integernn$\@,
there exist\/ $\C^r$-vector bundle mappings
\begin{equation*}
(\what{A}^m_s,\id_{\man{E}})\in
\vbmappings[r]{\Symalg*[s]{\pi_{\man{E}}^*\ctb{\man{M}}}\otimes\cvb{\man{E}}}
{\Symalg*[m]{\ctb{\man{E}}}\otimes\cvb{\man{E}}},\qquad s\in\{0,1,\dots,m\},
\end{equation*}
such that
\begin{equation*}
(\Sym_m\otimes\id_{\cvb{\man{E}}})\scirc\nabla^{\man{E},m}\vl{\lambda}=
\sum_{s=0}^m\what{A}^m_s((\Sym_s\otimes\id_{\cvb{\man{E}}})
\scirc\vl{(\nabla^{\man{M},\pi_{\man{E}},s}\lambda)})
\end{equation*}
for all\/ $\lambda\in\sections[m]{\dual{\man{E}}}$\@.
\begin{proof}
This follows along the lines of Lemma~\ref{lem:PsymiterateddersI} in the same
manner as Lemma~\ref{lem:VsymiterateddersI} follows from
Lemma~\ref{lem:PsymiterateddersI}\@, by taking tensor products with
$\cvb{\man{E}}$\@.
\end{proof}
\end{lemma}

The preceding lemma gives rise to an ``inverse,'' which we state in the
following lemma.
\begin{lemma}\label{lem:V*symiterateddersII}
Let\/ $r\in\{\infty,\omega\}$ and let\/
$\map{\pi_{\man{E}}}{\man{E}}{\man{M}}$ be a\/ $\C^r$-vector bundle, with the
data prescribed in Section~\ref{subsec:Esubmersion} to define the Riemannian
metric\/ $\metric_{\man{E}}$ on\/ $\man{E}$\@.  For\/ $m\in\integernn$\@,
there exist\/ $\C^r$-vector bundle mappings
\begin{equation*}
(\what{B}^m_s,\id_{\man{E}})\in
\vbmappings[r]{\Symalg*[s]{\ctb{\man{E}}}\otimes\cvb{\man{E}}}
{\Symalg*[m]{\pi_{\man{E}}^*\ctb{\man{M}}}\otimes\cvb{\man{E}}},\qquad s\in\{0,1,\dots,m\},
\end{equation*}
such that
\begin{equation*}
(\Sym_m\otimes_{\id_{\cvb{\man{E}}}})\scirc
\vl{(\nabla^{\man{M},\pi_{\man{E}},m}\lambda)}=
\sum_{s=0}^m\what{B}^m_s((\Sym_s\otimes\id_{\cvb{\man{E}}})\scirc
\nabla^{\man{E},s}\vl{\lambda})
\end{equation*}
for all\/ $\lambda\in\sections[m]{\dual{\man{E}}}$\@.
\begin{proof}
This follows along the lines of Lemma~\ref{lem:PsymiterateddersII} in the same
manner as Lemma~\ref{lem:VsymiterateddersI} follows from
Lemma~\ref{lem:PsymiterateddersI}\@, by taking tensor products with
$\cvb{\man{E}}$\@.
\end{proof}
\end{lemma}

We can put together the previous four lemmata into the following
decomposition result, which is to be regarded as the main result of this
section.
\begin{lemma}\label{lem:V*decomp}
Let\/ $r\in\{\infty,\omega\}$ and let\/
$\map{\pi_{\man{E}}}{\man{E}}{\man{M}}$ be a\/ $\C^r$-vector bundle, with the
data prescribed in Lemma~\ref{subsec:Esubmersion} to define the Riemannian
metric\/ $\metric_{\man{E}}$ on\/ $\man{E}$\@.  Then there exist\/
$\C^r$-vector bundle mappings
\begin{gather*}
A^m_{\nabla^{\man{E}}}\in\vbmappings[r]{\Pjetalg{m}{\man{E}}\otimes
\cvb{\man{E}}}{\Symalg*[\le m]{\pi_{\man{E}}^*\ctb{\man{M}}}\otimes
\cvb{\man{E}}},\\
B^m_{\nabla^{\man{E}}}\in\vbmappings[r]{\Pjetalg{m}{\man{E}}\otimes
\cvb{\man{E}}}{\Symalg*[\le m]{\ctb{\man{E}}}\otimes\cvb{\man{E}}},
\end{gather*}
defined by
\begin{align*}
A^m_{\nabla^{\man{E}}}(j_m(\vl{\lambda})(e))=&\;
\Sym_{\le m}\otimes\id_{\cvb{\man{E}}}(\vl{\lambda}(e),
\vl{(\nabla^{\pi_{\man{E}}}\lambda)}(e),
\dots,\vl{(\nabla^{\man{M},\pi_{\man{E}},m}\lambda)}(e)),\\
B^m_{\nabla^{\man{E}}}(j_m(\vl{\lambda})(e))=&\;
\Sym_{\le m}\otimes\id_{\cvb{\man{E}}}(\vl{\lambda}(e),
\nabla^{\man{E}}\vl{\lambda}(e),\dots,\nabla^{\man{E},m}\vl{\lambda}(e)).
\end{align*}
Moreover,\/ $A^m_{\nabla^{\man{E}}}$ is an isomorphism,\/
$B^m_{\nabla^{\man{E}}}$ is injective, and
{\small\begin{multline*}
B^m_{\nabla^{\man{E}}}\scirc(A^m_{\nabla^{\man{E}}})^{-1}\scirc
(\Sym_{\le m}\otimes\id_{\cvb{\man{E}}})
(\vl{\lambda}(e),\vl{(\nabla^{\pi_{\man{E}}}\lambda)}(e),\dots,
\vl{(\nabla^{\man{M},\pi_{\man{E}},m}\lambda)}(e))\\
=\left(\vl{\lambda}(e),\sum_{s=0}^1
\what{A}^1_s((\Sym_s\otimes\id_{\cvb{\man{E}}})
\scirc\vl{(\nabla^{\man{M},\pi_{\man{E}},s}\lambda)}(e)),\dots,
\sum_{s=0}^m\what{A}^m_s((\Sym_s\otimes\id_{\cvb{\man{E}}})
\scirc\vl{(\nabla^{\man{M},\pi_{\man{E}},s}\lambda)}(e))\right)
\end{multline*}}%
and
{\small\begin{multline*}
A^m_{\nabla^{\man{E}}}\scirc(B^m_{\nabla^{\man{E}}})^{-1}\scirc
(\Sym_{\le m}\otimes\id_{\cvb{\man{E}}})
(\vl{\lambda}(e),\nabla^{\man{E}}\vl{\lambda}(e),\dots,
\nabla^{\man{E},m}\vl{\lambda}(e))\\
=\left(\vl{\lambda}(e),\sum_{s=0}^1
\what{B}^1_s((\Sym_s\otimes\id_{\cvb{\man{E}}})
\scirc\nabla^{\man{E},s}\vl{\lambda}(e)),\dots,
\sum_{s=0}^m\what{B}^m_s((\Sym_s\otimes\id_{\cvb{\man{E}}})\scirc
\nabla^{\man{E},s}\vl{\lambda}(e))\right),
\end{multline*}}%
where the vector bundle mappings\/ $\what{A}^m_s$ and\/ $\what{B}^m_s$\@,\/
$s\in\{0,1,\dots,m\}$\@, are as in Lemmata~\ref{lem:V*symiterateddersI}
and~\ref{lem:V*symiterateddersII}\@.
\end{lemma}

\subsection{Isomorphisms for vertical lifts of endomorphisms}

Next we consider vertical lifts of endomorphisms defined by the mapping
\begin{equation*}
\sections[r]{\tensor[1,1]{\man{E}}}\ni L\mapsto
\vl{L}\in\sections[r]{\tensor[1,1]{\tb{\man{E}}}}.
\end{equation*}
We wish to relate the decomposition of the jets of $L$ with those of
$\vl{L}$\@.  Associated with this, we denote
\begin{equation*}
\Ljetalg[e]{m}{\man{E}}=\setdef{j_m\vl{L}(e)}
{L\in\sections[m]{\tensor[1,1]{\man{E}}}}.
\end{equation*}
By~\eqref{eq:jet=jetalg}\@, we have
\begin{equation*}
\Ljetalg[e]{m}{\man{E}}\simeq\Pjetalg[e]{m}{\man{E}}\otimes
\tensor[1,1]{\vb[e]{\man{E}}}.
\end{equation*}
As with the constructions of the preceding sections, we wish to use
Lemma~\ref{lem:Jmdecomp} to provide a decomposition of
$\Ljetalg{m}{\man{E}}$\@, and to do so we need to understand the covariant
derivatives
\begin{equation*}
\nabla^{\man{E},m}\vl{L}\eqdef\underbrace{\nabla^{\man{E}}\cdots
\nabla^{\man{E}}}_{m\ \textrm{times}}\vl{L},\qquad m\in\integernn.
\end{equation*}
In this section we omit proofs, since proofs follow along entirely similar
lines to those of the preceding section.

The first result we give is the following.
\begin{lemma}\label{lem:LiterateddersI}
Let\/ $r\in\{\infty,\omega\}$ and let\/
$\map{\pi_{\man{E}}}{\man{E}}{\man{M}}$ be a\/ $\C^r$-vector bundle, with the
data prescribed in Section~\ref{subsec:Esubmersion} to define the Riemannian
metric\/ $\metric_{\man{E}}$ on\/ $\man{E}$\@.  For\/ $m\in\integernn$\@,
there exist\/ $\C^r$-vector bundle mappings
\begin{equation*}
(A^m_s,\id_{\man{E}})\in
\vbmappings[r]{\tensor*[s]{\pi_{\man{E}}^*\ctb{\man{M}}}\otimes
\tensor[1,1]{\vb{\man{E}}}}{\tensor*[m]{\ctb{\man{E}}}\otimes
\tensor[1,1]{\vb{\man{E}}}},\qquad s\in\{0,1,\dots,m\},
\end{equation*}
such that
\begin{equation*}
\nabla^{\man{E},m}\vl{L}=
\sum_{s=0}^mA^m_s(\vl{(\nabla^{\man{M},\pi_{\man{E}},s}L)})
\end{equation*}
for all\/ $L\in\sections[m]{\tensor[1,1]{\man{E}}}$\@.  Moreover, the vector
bundle mappings\/ $A^m_0,A^m_1,\dots,A^m_m$ satisfy the recursion relations
prescribed by\/ $A^0_0(\beta_0)=\beta_0$ and
\begin{align*}
A^{m+1}_{m+1}(\beta_{m+1})=&\;\beta_{m+1},\\
A^{m+1}_s(\beta_s)=&\;(\nabla^{\man{E}}A^m_s)(\beta_s)+
A^m_{s-1}\otimes\id_{\ctb{\man{E}}}(\beta_s)-
\sum_{j=1}^sA^m_s\otimes\id_{\ctb{\man{E}}}(\Ins_j(\beta_s,B_{\pi_{\man{E}}}))\\
&\;+A^m_s\otimes\id_{\ctb{\man{E}}}
(\Ins_{s+1}(\beta_s,\dual{B}_{\pi_{\man{E}}})),
\enspace s\in\{1,\dots,m\},\\
A^{m+1}_0(\beta_0)=&\;(\nabla^{\man{E}}A^m_0)(\beta_0)-
A^m_0\otimes\id_{\ctb{\man{E}}}(\Ins_1(\beta_0,B_{\pi_{\man{E}}}))+
A^m_0\otimes\id_{\ctb{\man{E}}}(\Ins_2(\beta_0,\dual{B}_{\pi_{\man{E}}})),
\end{align*}
where\/ $\beta_s\in\tensor*[s]{\pi_{\man{E}}^*\ctb{\man{M}}}\otimes
\tensor[1,1]{\vb{\man{E}}}$\@,\/ $s\in\{0,1,\dots,m+1\}$\@.
\begin{proof}
This follows in the same manner as Lemma~\ref{lem:ViterateddersI}\@, making
use of Lemma~\pldblref{lem:bundleform}{pl:bundleform5}\@.
\end{proof}
\end{lemma}

The ``inverse'' of the preceding lemma is as follows.
\begin{lemma}\label{lem:LiterateddersII}
Let\/ $r\in\{\infty,\omega\}$ and let\/
$\map{\pi_{\man{E}}}{\man{E}}{\man{M}}$ be a\/ $\C^r$-vector bundle, with the
data prescribed in Section~\ref{subsec:Esubmersion} to define the Riemannian
metric\/ $\metric_{\man{E}}$ on\/ $\man{E}$\@.  For\/ $m\in\integernn$\@,
there exist\/ $\C^r$-vector bundle mappings
\begin{equation*}
(B^m_s,\id_{\man{E}})\in
\vbmappings[r]{\tensor*[s]{\ctb{\man{E}}}\otimes\tensor[1,1]{\vb{\man{E}}}}
{\tensor*[m]{\pi_{\man{E}}^*\ctb{\man{M}}}\otimes\tensor[1,1]{\vb{\man{E}}}},
\qquad s\in\{0,1,\dots,m\},
\end{equation*}
such that
\begin{equation*}
\vl{(\nabla^{\man{M},\pi_{\man{E}},m}L)}=
\sum_{s=0}^mB^m_s(\nabla^{\man{E},s}\vl{L})
\end{equation*}
for all\/ $L\in\sections[m]{\tensor[1,1]{\man{E}}}$\@.  Moreover, the vector
bundle mappings\/ $B^m_0,B^m_1,\dots,B^m_m$ satisfy the recursion relations
prescribed by\/ $B^0_0(\alpha_0)=\alpha_0$ and
\begin{align*}
B^{m+1}_{m+1}(\alpha_{m+1})=&\;\alpha_{m+1},\\
B^{m+1}_s(\alpha_s)=&\;(\nabla^{\man{E}}B^m_s)(\alpha_s)+
B^m_{s-1}\otimes\id_{\ctb{\man{E}}}(\alpha_s)+
\sum_{j=1}^m\Ins_j(B^m_s(\alpha_s),B_{\pi_{\man{E}}})\\
&\;-\Ins_{m+1}(B^m_s(\alpha_s),\dual{B}_{\pi_{\man{E}}}),\enspace s\in\{1,\dots,m\},\\
B^{m+1}_0(\alpha_0)=&\;(\nabla^{\man{E}}B^m_0)(\alpha_0)+
\sum_{j=1}^m\Ins_j(B^m_0(\alpha_0),B_{\pi_{\man{E}}})-
\Ins_{m+1}(B^m_0(\alpha_0),B_{\pi_{\man{E}}}),
\end{align*}
where\/ $\alpha_s\in\tensor*[s]{\ctb{\man{E}}}\otimes
\tensor[1,1]{\vb{\man{E}}}$\@,\/ $s\in\{0,1,\dots,m+1\}$\@.
\begin{proof}
This follows in the same manner as Lemma~\ref{lem:ViterateddersII}\@, making
use of Lemma~\pldblref{lem:bundleform}{pl:bundleform5}\@.
\end{proof}
\end{lemma}

Next we turn to symmetrised versions of the preceding lemmata.  We show that
the preceding two lemmata induce corresponding mappings between symmetric
tensors.
\begin{lemma}\label{lem:LsymiterateddersI}
Let\/ $r\in\{\infty,\omega\}$ and let\/
$\map{\pi_{\man{E}}}{\man{E}}{\man{M}}$ be a\/ $\C^r$-vector bundle, with the
data prescribed in Section~\ref{subsec:Esubmersion} to define the Riemannian
metric\/ $\metric_{\man{E}}$ on\/ $\man{E}$\@.  For\/ $m\in\integernn$\@,
there exist\/ $\C^r$-vector bundle mappings
\begin{equation*}
(\what{A}^m_s,\id_{\man{E}})\in
\vbmappings[r]{\Symalg*[s]{\pi_{\man{E}}^*\ctb{\man{M}}}\otimes
\tensor[1,1]{\vb{\man{E}}}}{\Symalg*[m]{\ctb{\man{E}}}\otimes
\tensor[1,1]{\vb{\man{E}}}},\qquad s\in\{0,1,\dots,m\},
\end{equation*}
such that
\begin{equation*}
(\Sym_m\otimes\id_{\tensor[1,1]{\vb{\man{E}}}})\scirc\nabla^{\man{E},m}\vl{L}=
\sum_{s=0}^m\what{A}^m_s((\Sym_s\otimes\id_{\tensor[1,1]{\vb{\man{E}}}})
\scirc\vl{(\nabla^{\man{M},\pi_{\man{E}},s}L)})
\end{equation*}
for all\/ $L\in\sections[m]{\tensor[1,1]{\man{E}}}$\@.
\begin{proof}
This follows along the lines of Lemma~\ref{lem:PsymiterateddersI} in the same
manner as Lemma~\ref{lem:VsymiterateddersI} follows from
Lemma~\ref{lem:PsymiterateddersI}\@, by taking tensor products with
$\tensor[1,1]{\vb{\man{E}}}$\@.
\end{proof}
\end{lemma}

The preceding lemma gives rise to an ``inverse,'' which we state in the
following lemma.
\begin{lemma}\label{lem:LsymiterateddersII}
Let\/ $r\in\{\infty,\omega\}$ and let\/
$\map{\pi_{\man{E}}}{\man{E}}{\man{M}}$ be a\/ $\C^r$-vector bundle, with the
data prescribed in Section~\ref{subsec:Esubmersion} to define the Riemannian
metric\/ $\metric_{\man{E}}$ on\/ $\man{E}$\@.  For\/ $m\in\integernn$\@,
there exist\/ $\C^r$-vector bundle mappings
\begin{equation*}
(\what{B}^m_s,\id_{\man{E}})\in
\vbmappings[r]{\Symalg*[s]{\ctb{\man{E}}}\otimes\tensor[1,1]{\vb{\man{E}}}}
{\Symalg*[m]{\pi_{\man{E}}^*\ctb{\man{M}}}\otimes\tensor[1,1]{\vb{\man{E}}}},
\qquad s\in\{0,1,\dots,m\},
\end{equation*}
such that
\begin{equation*}
(\Sym_m\otimes_{\id_{\tensor[1,1]{\vb{\man{E}}}}})\scirc
\vl{(\nabla^{\man{M},\pi_{\man{E}},m}L)}=
\sum_{s=0}^m\what{B}^m_s((\Sym_s\otimes\id_{\tensor[1,1]{\vb{\man{E}}}})\scirc
\nabla^{\man{E},s}\vl{L})
\end{equation*}
for all\/ $L\in\sections[m]{\tensor[1,1]{\man{E}}}$\@.
\begin{proof}
This follows along the lines of Lemma~\ref{lem:PsymiterateddersII} in the
same manner as Lemma~\ref{lem:VsymiterateddersI} follows from
Lemma~\ref{lem:PsymiterateddersI}\@, by taking tensor products with
$\tensor[1,1]{\vb{\man{E}}}$\@.
\end{proof}
\end{lemma}

We can put together the previous four lemmata into the following
decomposition result, which is to be regarded as the main result of this
section.
\begin{lemma}\label{lem:Ldecomp}
Let\/ $r\in\{\infty,\omega\}$ and let\/
$\map{\pi_{\man{E}}}{\man{E}}{\man{M}}$ be a\/ $\C^r$-vector bundle, with the
data prescribed in Section~\ref{subsec:Esubmersion} to define the Riemannian
metric\/ $\metric_{\man{E}}$ on\/ $\man{E}$\@.  Then there exist\/
$\C^r$-vector bundle mappings
\begin{gather*}
A^m_{\nabla^{\man{E}}}\in\vbmappings[r]{\Pjetalg{m}{\man{E}}\otimes
\tensor[1,1]{\vb{\man{E}}}}
{\Symalg*[\le m]{\pi_{\man{E}}^*\ctb{\man{M}}}\otimes
\tensor[1,1]{\vb{\man{E}}}},\\
B^m_{\nabla^{\man{E}}}\in\vbmappings[r]{\Pjetalg{m}{\man{E}}\otimes
\tensor[1,1]{\vb{\man{E}}}}
{\Symalg*[\le m]{\ctb{\man{E}}}\otimes\tensor[1,1]{\vb{\man{E}}}},
\end{gather*}
defined by
\begin{align*}
A^m_{\nabla^{\man{E}}}(j_m(\vl{L})(e))=&\;
\Sym_{\le m}\otimes\id_{\tensor[1,1]{\vb{\man{E}}}}(\vl{L}(e),
\vl{(\nabla^{\pi_{\man{E}}}L)}(e),\dots,
\vl{(\nabla^{\man{M},\pi_{\man{E}},m}L)}(e)),\\
B^m_{\nabla^{\man{E}}}(j_m(\vl{L})(e))=&\;
\Sym_{\le m}\otimes\id_{\tensor[1,1]{\vb{\man{E}}}}
(\vl{L}(e),\nabla^{\man{E}}\vl{L}(e),\dots,\nabla^{\man{E},m}\vl{L}(e)).
\end{align*}
Moreover,\/ $A^m_{\nabla^{\man{E}}}$ is an isomorphism,\/
$B^m_{\nabla^{\man{E}}}$ is injective, and
{\small\begin{multline*}
B^m_{\nabla^{\man{E}}}\scirc(A^m_{\nabla^{\man{E}}})^{-1}\scirc
(\Sym_{\le m}\otimes\id_{\tensor[1,1]{\vb{\man{E}}}})
(\vl{L}(e),\vl{(\nabla^{\pi_{\man{E}}}L)}(e),\dots,
\vl{(\nabla^{\man{M},\pi_{\man{E}},m}L)}(e))\\
=\left(\vl{L}(e),\sum_{s=0}^1\what{A}^1_s
((\Sym_s\otimes\id_{\tensor[1,1]{\vb{\man{E}}}})
\scirc\vl{(\nabla^{\man{M},\pi_{\man{E}},s}L)}(e)),\dots,\right.\\
\left.\sum_{s=0}^m\what{A}^m_s((\Sym_s\otimes\id_{\tensor[1,1]{\vb{\man{E}}}})
\scirc\vl{(\nabla^{\man{M},\pi_{\man{E}},s}L)}(e))\right)
\end{multline*}}%
and
{\small\begin{multline*}
A^m_{\nabla^{\man{E}}}\scirc(B^m_{\nabla^{\man{E}}})^{-1}\scirc
(\Sym_{\le m}\otimes\id_{\tensor[1,1]{\vb{\man{E}}}})
(\vl{L}(e),\nabla^{\man{E}}\vl{L}(e),\dots,\nabla^{\man{E},m}\vl{L}(e))\\
=\left(\vl{L}(e),\sum_{s=0}^1\what{B}^1_s((\Sym_s\otimes
\id_{\tensor[1,1]{\vb{\man{E}}}})\scirc\nabla^{\man{E},s}\vl{L}(e)),\dots,
\sum_{s=0}^m\what{B}^m_s((\Sym_s\otimes\id_{\tensor[1,1]{\vb{\man{E}}}})\scirc
\nabla^{\man{E},s}\vl{L}(e))\right),
\end{multline*}}%
where the vector bundle mappings\/ $\what{A}^m_s$ and\/ $\what{B}^m_s$\@,\/
$s\in\{0,1,\dots,m\}$\@, are as in Lemmata~\ref{lem:LsymiterateddersI}
and~\ref{lem:LsymiterateddersII}\@.
\end{lemma}

\subsection{Isomorphisms for vertical evaluations of dual sections}

Next we consider vertical evaluations of endomorphisms given by the mapping
\begin{equation*}
\sections[r]{\dual{\man{E}}}\ni\lambda\mapsto\ve{\lambda}\in\func[r]{\man{E}}.
\end{equation*}
To study the relationship between the decomposition of the jets of $\lambda$
with those of the jets of $\ve{\lambda}$\@, we denote
\begin{equation*}
\Djetalg[e]{m}{\man{E}}=\setdef{j_m\ve{\lambda}(e)}
{\lambda\in\sections[m]{\dual{\man{E}}}}.
\end{equation*}
By~\eqref{eq:jet=jetalg}\@, we have
\begin{equation*}
\Djetalg[e]{m}{\man{E}}\subset\Pjetalg[e]{m}{\man{E}}.
\end{equation*}
As we shall see, one can be a little more explicit about the nature of
$\Djetalg[e]{m}{\man{E}}$\@, and see that
\begin{equation*}
\Djetalg[e]{m}{\man{E}}\simeq
(\Pjetalg[e]{m}{\man{E}}\otimes\cvb{\man{E}})\oplus
(\Pjetalg[m-1]{\man{E}}\otimes\cvb{\man{E}}).
\end{equation*}
However, this sort of isomorphism is too cumbersome to make explicit.  As
with the constructions of the preceding sections, we wish to use
Lemma~\ref{lem:Jmdecomp} to provide a decomposition of
$\Djetalg{m}{\man{E}}$\@, and to do so we need to understand the covariant
derivatives
\begin{equation*}
\nabla^{\man{E},m}\ve{\lambda}\eqdef\underbrace{\nabla^{\man{E}}\cdots
\nabla^{\man{E}}}_{m\ \textrm{times}}\ve{\lambda},\qquad m\in\integernn.
\end{equation*}
The results in this section have a slightly different character than in the
preceding sections.  We will not give the complete proofs, but will note that
they are very similar to the complete proofs given in the next section.

Our first result is then the following.
\begin{lemma}\label{lem:DiterateddersI}
Let\/ $r\in\{\infty,\omega\}$ and let\/
$\map{\pi_{\man{E}}}{\man{E}}{\man{M}}$ be a\/ $\C^r$-vector bundle, with the
data prescribed in Section~\ref{subsec:Esubmersion} to define the Riemannian
metric\/ $\metric_{\man{E}}$ on\/ $\man{E}$\@.  For\/ $m\in\integernn$\@,
there exist\/ $\C^r$-vector bundle mappings
\begin{equation*}
(A^m_s,\id_{\man{E}})\in
\vbmappings[r]{\tensor*[s]{\pi_{\man{E}}^*\ctb{\man{M}}}}
{\tensor*[m]{\ctb{\man{E}}}},\qquad s\in\{0,1,\dots,m\},
\end{equation*}
and
\begin{equation*}
(C^m_s,\id_{\man{E}})\in
\vbmappings[r]{\tensor*[s]{\pi_{\man{E}}^*\ctb{\man{M}}}\otimes\cvb{\man{E}}}
{\tensor*[m-1]{\ctb{\man{E}}}\otimes\cvb{\man{E}}},\qquad
s\in\{0,1,\dots,m-1\},
\end{equation*}
such that
\begin{equation*}
\nabla^{\man{E},m}\ve{\lambda}=
\sum_{s=0}^mA^m_s(\ve{(\nabla^{\man{M},\pi_{\man{E}},s}\lambda)})+
\sum_{s=0}^{m-1}C^m_s(\vl{(\nabla^{\man{M},\pi_{\man{E}},s}\lambda)})%
\footnote{Here we regard $\cvb{\man{E}}$ as a subbundle of
$\ctb{\man{E}}$\@.}
\end{equation*}
for all\/ $\lambda\in\sections[m]{\dual{\man{E}}}$\@.  Moreover, the vector
bundle mappings\/ $A^m_0,A^m_1,\dots,A^m_m$ and\/
$C^m_0,C^m_1,\dots,C^m_{m-1}$ satisfy the recursion relations prescribed by
\begin{equation*}
A^0_0(\beta_0)=\beta_0,\enspace A^1_1(\beta_1)=\beta_1,\enspace
A^1_0(\beta_0)=\Ins_1(\beta_0,B_{\pi_{\man{E}}}),\enspace C^1_0(\gamma_0)=\gamma_0,
\end{equation*}
and, for\/ $m\ge2$\@,
\begin{align*}
A^{m+1}_{m+1}(\beta_{m+1})=&\;\beta_{m+1}\\
A^{m+1}_m(\beta_m)=&\;A^m_{m-1}\otimes\id_{\ctb{\man{E}}}(\beta_m)-
\sum_{j=1}^m\Ins_j(\beta_m,B_{\pi_{\man{E}}})\\
A^{m+1}_s(\beta_s)=&\;(\nabla^{\man{E}}A^m_s)(\beta_m)+
A^m_{s-1}\otimes\id_{\ctb{\man{E}}}(\beta_s)
-\sum_{j=1}^sA^m_s\otimes\id_{\ctb{\man{E}}}(\Ins_j(\beta_s,B_{\pi_{\man{E}}})),\\
&\;\enspace s\in\{1,\dots,m-1\},\\
A^{m+1}_0(\beta_0)=&\;(\nabla^{\man{E}}A^m_0)(\beta_0)
\end{align*}
and
\begin{align*}
C^{m+1}_m(\gamma_m)=&\;C^m_{m-1}\otimes\id_{\ctb{\man{E}}}(\gamma_m)+\gamma_m\\
C^{m+1}_s(\gamma_s)=&\;A^m_s\otimes\id_{\ctb{\man{E}}}(\gamma_s)+
(\nabla^{\man{E}}C^m_s)(\gamma_s)+C^m_{s-1}\otimes\id_{\ctb{\man{E}}}(\gamma_s)\\
&\;-\sum_{j=1}^{s+1}C^m_s\otimes
\id_{\ctb{\man{E}}}(\Ins_j(\gamma_s,B_{\pi_{\man{E}}})),
\enspace s\in\{1,\dots,m-1\},\\
C^{m+1}_0(\gamma_0)=&\;A^m_0\otimes\id_{\ctb{\man{E}}}(\gamma_0)+
(\nabla^{\man{E}}C^m_0)(\gamma_0)-C^m_0\otimes\id_{\ctb{\man{E}}}
(\Ins_1(\gamma_0,B_{\pi_{\man{E}}})),
\end{align*}
where\/ $\beta_s\in\tensor*[s]{\pi_{\man{E}}^*\ctb{\man{M}}}$\@,\/
$s\in\{0,1,\dots,m+1\}$\@, and\/
$\gamma_s\in\tensor*[s]{\pi_{\man{E}}^*\ctb{\man{M}}}\otimes
\cvb{\man{E}}$\@,\/ $s\in\{0,1,\dots,m-1\}$\@.
\begin{proof}
This follows in the same manner as Lemma~\ref{lem:CiterateddersI} below,
making use of Lemma~\pldblref{lem:bundleform}{pl:bundleform6}\@.
\end{proof}
\end{lemma}

Now we ``invert'' the constructions from the preceding lemma.
\begin{lemma}\label{lem:DiterateddersII}
Let\/ $r\in\{\infty,\omega\}$ and let\/
$\map{\pi_{\man{E}}}{\man{E}}{\man{M}}$ be a\/ $\C^r$-vector bundle, with the
data prescribed in Section~\ref{subsec:Esubmersion} to define the Riemannian
metric\/ $\metric_{\man{E}}$ on\/ $\man{E}$\@.  For\/ $m\in\integernn$\@,
there exist\/ $\C^r$-vector bundle mappings
\begin{equation*}
(B^m_s,\id_{\man{E}})\in
\vbmappings[r]{\tensor*[s]{\ctb{\man{E}}}}
{\tensor*[m]{\pi_{\man{E}}^*\ctb{\man{M}}}},\qquad s\in\{0,1,\dots,m\},
\end{equation*}
and
\begin{equation*}
(D^m_s,\id_{\man{E}})\in
\vbmappings[r]{\tensor*[s]{\ctb{\man{E}}}\otimes\cvb{\man{E}}}
{\tensor*[m-1]{\pi_{\man{E}}^*\ctb{\man{M}}}\otimes\cvb{\man{E}}},
\qquad s\in\{0,1,\dots,m-1\},
\end{equation*}
such that
\begin{equation*}
\ve{(\nabla^{\man{M},\pi_{\man{E}},m}\lambda)}=
\sum_{s=0}^mB^m_s(\nabla^{\man{E},s}\ve{\lambda})+
\sum_{s=0}^{m-1}D^m_s(\nabla^{\man{E},s}\vl{\lambda})
\end{equation*}
for all\/ $\lambda\in\sections[m]{\dual{\man{E}}}$\@.  Moreover, the
vector bundle mappings\/ $B^m_0,B^m_1,\dots,B^m_m$ and\/
$D^m_0,D^m_1,\dots,D^m_{m-1}$ satisfy the recursion relations prescribed by\/
$B^0_0(\alpha_0)=\alpha_0$\@, $D^1_0(\gamma_0)=\gamma_0$\@,
\begin{align*}
B^{m+1}_{m+1}(\alpha_{m+1})=&\;\alpha_{m+1}\\
B^{m+1}_m(\alpha_m)=&\;B^m_{m-1}\otimes\id_{\ctb{\man{E}}}(\alpha_m)+
\sum_{j=1}^m\Ins_j(\alpha_m,B_{\pi_{\man{E}}})-
\Ins_{m+1}(\alpha_m,\dual{B}_{\pi_{\man{E}}})\\
B^{m+1}_s=&\;(\nabla^{\man{E}}B^m_s)(\alpha_s)+
B^m_{s-1}\otimes\id_{\ctb{\man{E}}}(\alpha_s)
+\sum_{j=1}^m\Ins_j(B^m_s(\alpha_s),B_{\pi_{\man{E}}})\\
&\;-\Ins_{m+1}(B^m_s(\alpha_s),\dual{B}_{\pi_{\man{E}}}),\enspace s\in\{1,\dots,m-1\},\\
B^{m+1}_0(\alpha_0)=&\;(\nabla^{\man{E}}B^m_0)(\alpha_0)+
\sum_{j=1}^m\Ins_j(B^m_0(\alpha_0),B_{\pi_{\man{E}}})-
\Ins_{m+1}(B^m_0(\alpha_0),\dual{B}_{\pi_{\man{E}}})
\end{align*}
and
\begin{align*}
D^{m+1}_m(\gamma_m)=&\;D^m_{m-1}\otimes\id_{\ctb{\man{E}}}(\gamma_m)-\gamma_m\\
D^m_s(\gamma_s)=&\;(\nabla^{\man{E}}D^m_s)(\gamma_s)+
D^m_{s-1}\otimes\id_{\ctb{\man{E}}}(\gamma_s)-\ol{B}^m_s(\gamma_s),
\enspace s\in\{1,\dots,m-1\},\\
D^{m+1}_0=&\;(\nabla^{\man{E}}D^m_0)(\gamma_0)-\ol{B}^m_0(\gamma_0)
\end{align*}
for\/ $\alpha_s\in\tensor*[s]{\ctb{\man{E}}}$\@,\/ $s\in\{0,1,\dots,m\}$\@,
and\/ $\gamma_s\in\tensor*[s]{\ctb{\man{E}}}\otimes\cvb{\man{E}}$\@,\/
$s\in\{0,1,\dots,m-1\}$\@, and where
\begin{equation*}
(\ol{B}^m_s,\id_{\man{E}})\in
\vbmappings[r]{\tensor*[s]{\ctb{\man{E}}}\otimes\cvb{\man{E}}}
{\tensor*[m]{\pi_{\man{E}}^*\ctb{\man{M}}}\otimes\cvb{\man{E}}},\qquad
s\in\{0,1,\dots,m\},
\end{equation*}
are the vector bundle mappings from Lemma~\ref{lem:V*iterateddersII}\@.
\begin{proof}
This follows in the same manner as Lemma~\ref{lem:CiterateddersII} below,
making use of Lemma~\pldblref{lem:bundleform}{pl:bundleform6}\@.
\end{proof}
\end{lemma}

Next we turn to symmetrised versions of the preceding lemmata.  We show that
the preceding two lemmata induce corresponding mappings between symmetric
tensors.
\begin{lemma}\label{lem:DsymiterateddersI}
Let\/ $r\in\{\infty,\omega\}$ and let\/
$\map{\pi_{\man{E}}}{\man{E}}{\man{M}}$ be a\/ $\C^r$-vector bundle, with the
data prescribed in Section~\ref{subsec:Esubmersion} to define the Riemannian
metric\/ $\metric_{\man{E}}$ on\/ $\man{E}$\@.  For\/ $m\in\integernn$\@,
there exist\/ $\C^r$-vector bundle mappings
\begin{equation*}
(\what{A}^m_s,\id_{\man{E}})\in
\vbmappings[r]{\Symalg*[s]{\pi_{\man{E}}^*\ctb{\man{M}}}}
{\Symalg*[m]{\ctb{\man{E}}}},\qquad s\in\{0,1,\dots,m\},
\end{equation*}
and
\begin{equation*}
(\what{C}^m_s,\id_{\man{E}})\in
\vbmappings[r]{\Symalg*[s]{\pi_{\man{E}}^*\ctb{\man{M}}}\otimes\cvb{\man{E}}}
{\Symalg*[m]{\ctb{\man{E}}}},\qquad s\in\{0,1,\dots,m-1\},
\end{equation*}
such that
\begin{equation*}
\Sym_m\scirc\nabla^{\man{E},m}\ve{\lambda}
=\sum_{s=0}^m\what{A}^m_s(\Sym_s\scirc
\ve{(\nabla^{\man{M},\pi_{\man{E}},s}\lambda)})+
\sum_{s=0}^{m-1}\what{C}^m_s((\Sym_s\otimes\id_{\cvb{\man{E}}})
\scirc\vl{(\nabla^{\man{M},\pi_{\man{E}},s}\lambda)})
\end{equation*}
for all\/ $\lambda\in\sections[m]{\dual{\man{E}}}$\@.
\begin{proof}
The proof here follows along the lines of Lemma~\ref{lem:CsymiterateddersI}
below.
\end{proof}
\end{lemma}

The preceding lemma gives rise to an ``inverse,'' which we state in the
following lemma.
\begin{lemma}\label{lem:DsymiterateddersII}
Let\/ $r\in\{\infty,\omega\}$ and let\/
$\map{\pi_{\man{E}}}{\man{E}}{\man{M}}$ be a\/ $\C^r$-vector bundle, with the
data prescribed in Section~\ref{subsec:Esubmersion} to define the Riemannian
metric\/ $\metric_{\man{E}}$ on\/ $\man{E}$\@.  For\/ $m\in\integernn$\@,
there exist\/ $\C^r$-vector bundle mappings
\begin{equation*}
(\what{B}^m_s,\id_{\man{E}})\in
\vbmappings[r]{\Symalg*[s]{\ctb{\man{E}}}}
{\Symalg*[m]{\pi_{\man{E}}^*\ctb{\man{M}}}},\qquad s\in\{0,1,\dots,m\},
\end{equation*}
and
\begin{equation*}
(\what{D}^m_s,\id_{\man{E}})\in
\vbmappings[r]{\Symalg*[s]{\ctb{\man{E}}}\otimes\cvb{\man{E}}}
{\Symalg*[m]{\pi_{\man{E}}^*\ctb{\man{M}}}},\qquad s\in\{0,1,\dots,m-1\},
\end{equation*}
such that
\begin{equation*}
\Sym_m\scirc\ve{(\nabla^{\man{M},\pi_{\man{E}},m}\lambda)}
=\sum_{s=0}^m\what{B}^m_s(\Sym_s\scirc\nabla^{\man{E},s}\ve{\lambda})+
\sum_{s=0}^{m-1}\what{D}^m_s((\Sym_s\otimes\id_{\cvb{\man{E}}})\scirc
\nabla^{\man{E},s}\vl{\lambda})
\end{equation*}
for all\/ $\lambda\in\sections[m]{\dual{\man{E}}}$\@.
\begin{proof}
The proof here follows along the lines of Lemma~\ref{lem:CsymiterateddersI}
below.
\end{proof}
\end{lemma}

We can put together the previous four lemmata, along with
Lemma~\ref{lem:Ldecomp}\@, into the following decomposition result, which is
to be regarded as the main result of this section.
\begin{lemma}\label{lem:Ddecomp}
Let\/ $r\in\{\infty,\omega\}$ and let\/
$\map{\pi_{\man{E}}}{\man{E}}{\man{M}}$ be a\/ $\C^r$-vector bundle, with the
data prescribed in Section~\ref{subsec:Esubmersion} to define the Riemannian
metric\/ $\metric_{\man{E}}$ on\/ $\man{E}$\@.  Then there exist\/
$\C^r$-vector bundle mappings
\begin{equation*}
A^m_{\nabla^{\man{E}}}\in\vbmappings[r]{\Djetalg{m}{\man{E}}}
{\Symalg*[\le m]{\pi_{\man{E}}^*\ctb{\man{M}}}},\quad
B^m_{\nabla^{\man{E}}}\in\vbmappings[r]{\Djetalg{m}{\man{E}}}
{\Symalg*[\le m]{\ctb{\man{E}}}}
\end{equation*}
defined by
\begin{align*}
A^m_{\nabla^{\man{E}}}(j_m(\ve{\lambda})(e))=&\;
\Sym_{\le m}
(\ve{\lambda}(e),\ve{(\nabla^{\pi_{\man{E}}}\lambda)}(e),\dots,
\ve{(\nabla^{\man{M},\pi_{\man{E}},m}\lambda)}(e)),\\
B^m_{\nabla^{\man{E}}}(j_m(\ve{\lambda})(e))=&\;
\Sym_{\le m}(\ve{\lambda}(e),
\nabla^{\man{E}}\ve{\lambda}(e),\dots,\nabla^{\man{E},m}\ve{\lambda}(e)).
\end{align*}
Moreover,\/ $A^m_{\nabla^{\man{E}}}$ and\/ $B^m_{\nabla^{\man{E}}}$ are
injective, and
{\small\begin{multline*}
B^m_{\nabla^{\man{E}}}\scirc(A^m_{\nabla^{\man{E}}})^{-1}\scirc\Sym_{\le m}
(\ve{\lambda}(e),\ve{(\nabla^{\pi_{\man{E}}}\lambda)}(e),\dots,
\ve{(\nabla^{\man{M},\pi_{\man{E}},m}\lambda)}(e))\\
=\left(\ve{\lambda}(e),\sum_{s=0}^1\what{A}^1_s(\Sym_s
\scirc\ve{(\nabla^{\man{M},\pi_{\man{E}},s}\lambda)}(e)),\dots,
\sum_{s=0}^m\what{A}^m_s(\Sym_s
\scirc\ve{(\nabla^{\man{M},\pi_{\man{E}},s}\lambda)}(e))\right)\\
+\left(0,\vl{\lambda}(e),
\sum_{s=0}^1\what{C}^2_s((\Sym_s\otimes\id_{\cvb{\man{E}}})
\scirc\vl{(\nabla^{\man{M},\pi_{\man{E}},s}\lambda)}(e)),\dots,\right.\\
\left.\sum_{s=0}^{m-1}\what{C}^m_s((\Sym_s\otimes\id_{\cvb{\man{E}}})
\scirc\vl{(\nabla^{\man{M},\pi_{\man{E}},s}\lambda)}(e))\right)
\end{multline*}}%
and
{\small\begin{multline*}
A^m_{\nabla^{\man{E}}}\scirc(B^m_{\nabla^{\man{E}}})^{-1}\scirc\Sym_{\le m}
(\ve{\lambda}(e),\nabla^{\man{E}}\ve{\lambda}(e),\dots,
\nabla^{\man{E},m}\ve{\lambda}(e))\\
=\left(\ve{\lambda}(e),\sum_{s=0}^1\what{B}^1_s(\Sym_s
\scirc\nabla^{\man{E},s}\ve{\lambda}(e)),\dots,
\sum_{s=0}^m\what{B}^m_s(\Sym_s\scirc\nabla^{\man{E},s}
\ve{\lambda}(e))\right)\\
+\left(0,\vl{\lambda}(e),
\sum_{s=0}^1\what{D}^2_s((\Sym_s\otimes\id_{\cvb{\man{E}}})
\scirc\nabla^{\man{E},s}\vl{\lambda}(e)),\dots,
\sum_{s=0}^{m-1}\what{D}^m_s((\Sym_s\otimes\id_{\cvb{\man{E}}})\scirc
\nabla^{\man{E},s}\vl{\lambda}(e))\right),
\end{multline*}}%
where the vector bundle mappings\/ $\what{A}^m_s$ and\/ $\what{B}^m_s$\@,\/
$s\in\{0,1,\dots,m\}$\@, and\/ $\what{C}^m_s$ and\/ $\what{D}^m_s$\@,
$s\in\{0,1,\dots,m-1\}$\@, are as in Lemmata~\ref{lem:DsymiterateddersI}
and~\ref{lem:DsymiterateddersII}\@.
\end{lemma}

\subsection{Isomorphisms for vertical evaluations of
endomorphisms}\label{subsec:Cjetisom}

Next we consider vertical evaluations of endomorphisms,~\ie~the mapping given
by
\begin{equation*}
\sections[r]{\tensor[1,1]{\man{E}}}\ni L\mapsto
\ve{L}\in\sections[r]{\tb{\man{E}}}.
\end{equation*}
To study the relationship between the decomposition of jets of $L$ with those
of $\ve{L}$\@, we denote
\begin{equation*}
\Cjetalg[e]{m}{\man{E}}=\setdef{j_m\ve{L}(e)}
{L\in\sections[m]{\tensor[1,1]{\man{E}}}}.
\end{equation*}
By~\eqref{eq:jet=jetalg}\@, we have
\begin{equation*}
\Cjetalg[e]{m}{\man{E}}\subset\Pjetalg[e]{m}{\man{E}}\otimes\vb[e]{\man{E}}.
\end{equation*}
As we shall see, one can be a little more explicit about the nature of
$\Djetalg[e]{m}{\man{E}}$\@, and see that
\begin{equation*}
\Cjetalg[e]{m}{\man{E}}\simeq
(\Pjetalg[e]{m}{\man{E}}\otimes\cvb{\man{E}}\otimes\vb[e]{\man{E}})
\oplus(\Pjetalg[m-1]{\man{E}}\otimes\cvb{\man{E}}\otimes\vb[e]{\man{E}}).
\end{equation*}
However, this sort of isomorphism is too cumbersome to make explicit.  As
with the constructions of the preceding sections, we wish to use
Lemma~\ref{lem:Jmdecomp} to provide a decomposition of
$\Cjetalg{m}{\man{E}}$\@, and to do so we need to understand the covariant
derivatives
\begin{equation*}
\nabla^{\man{E},m}\ve{L}\eqdef\underbrace{\nabla^{\man{E}}\cdots
\nabla^{\man{E}}}_{m\ \textrm{times}}\ve{L},\qquad m\in\integernn.
\end{equation*}
The results in this section have a slightly different character than in the
preceding sections, so we provide complete proofs.

The first result we give is the following.
\begin{lemma}\label{lem:CiterateddersI}
Let\/ $r\in\{\infty,\omega\}$ and let\/
$\map{\pi_{\man{E}}}{\man{E}}{\man{M}}$ be a\/ $\C^r$-vector bundle, with the
data prescribed in Section~\ref{subsec:Esubmersion} to define the Riemannian
metric\/ $\metric_{\man{E}}$ on\/ $\man{E}$\@.  For\/ $m\in\integernn$\@,
there exist\/ $\C^r$-vector bundle mappings
\begin{equation*}
(A^m_s,\id_{\man{E}})\in
\vbmappings[r]{\tensor*[s]{\pi_{\man{E}}^*\ctb{\man{M}}}\otimes\vb{\man{E}}}
{\tensor*[m]{\ctb{\man{E}}}\otimes\vb{\man{E}}},\qquad
s\in\{0,1,\dots,m\},
\end{equation*}
and
\begin{equation*}s
(C^m_s,\id_{\man{E}})\in
\vbmappings[r]{\tensor*[s]{\pi_{\man{E}}^*\ctb{\man{M}}}\otimes
\tensor[1,1]{\vb{\man{E}}}}
{\tensor*[m-1]{\ctb{\man{E}}}\otimes\tensor[1,1]{\vb{\man{E}}}},\qquad
s\in\{0,1,\dots,m-1\},
\end{equation*}
such that
\begin{equation*}
\nabla^{\man{E},m}\ve{L}=
\sum_{s=0}^mA^m_s(\ve{(\nabla^{\man{M},\pi_{\man{E}},s}L)})+
\sum_{s=0}^{m-1}C^m_s(\vl{(\nabla^{\man{M},\pi_{\man{E}},s}L)})%
\footnote{Here we regard $\tensor[1,1]{\vb{\man{E}}}$ as a
subbundle of $\ctb{\man{E}}\otimes\vb{\man{E}}$ by the mapping
\begin{equation*}
\tensor[1,1]{\vb{\man{E}}}\ni A\mapsto
A\scirc\ver\in\ctb{\man{E}}\otimes\vb{\man{E}}.
\end{equation*}\label{fn:T11inT*}}
\end{equation*}
for all\/ $L\in\sections[m]{\tensor[1,1]{\man{E}}}$\@.  Moreover, the vector
bundle mappings\/ $A^m_0,A^m_1,\dots,A^m_m$ and\/
$C^m_0,C^m_1,\dots,C^m_{m-1}$ satisfy the recursion relations prescribed by
\begin{equation*}
A^0_0(\beta_0)=\beta_0,\enspace A^1_1(\beta_1)=\beta_1,\enspace
A^1_0(\beta_0)=\Ins_1(\beta_0,B_{\pi_{\man{E}}}),\enspace C^1_0(\gamma_0)=\gamma_0,
\end{equation*}
and, for\/ $m\ge2$\@,
{\small\begin{align*}
A^{m+1}_{m+1}(\beta_{m+1})=&\;\beta_{m+1}\\
A^{m+1}_m(\beta_m)=&\;A^m_{m-1}\otimes\id_{\ctb{\man{E}}}(\beta_m)-
\sum_{j=1}^m\Ins_j(\beta_m,B_{\pi_{\man{E}}})+
\Ins_{m+1}(\beta_m,\dual{B}_{\pi_{\man{E}}})\\
A^{m+1}_s(\beta_s)=&\;(\nabla^{\man{E}}A^m_s)(\beta_m)+
A^m_{s-1}\otimes\id_{\ctb{\man{E}}}(\beta_s)
-\sum_{j=1}^sA^m_s\otimes\id_{\ctb{\man{E}}}(\Ins_j(\beta_s,B_{\pi_{\man{E}}}))\\
&\;+A^m_s\otimes\id_{\ctb{\man{E}}}
(\Ins_{s+1}(\beta_s,\dual{B}_{\pi_{\man{E}}})),
\enspace s\in\{1,\dots,m-1\},\\
A^{m+1}_0(\beta_0)=&\;(\nabla^{\man{E}}A^m_0)(\beta_0)-
A^m_0\otimes\id_{\ctb{\man{E}}}(\Ins_1(\beta_0,\dual{B}_{\pi_{\man{E}}}))
\end{align*}}%
and
{\small\begin{align*}
C^{m+1}_m(\gamma_m)=&\;C^m_{m-1}\otimes\id_{\ctb{\man{E}}}(\gamma_m)+\gamma_m\\
C^{m+1}_s(\gamma_s)=&\;A^m_s\otimes\id_{\ctb{\man{E}}}(\gamma_s)+
(\nabla^{\man{E}}C^m_s)(\gamma_s)+C^m_{s-1}\otimes\id_{\ctb{\man{E}}}(\gamma_s)\\
&\;-\sum_{j=1}^{s+1}C^m_s\otimes\id_{\ctb{\man{E}}}
(\Ins_j(\gamma_s,B_{\pi_{\man{E}}}))+
C^m_s\otimes\id_{\ctb{\man{E}}}(\Ins_{s+1}(\gamma_s,\dual{B}_{\pi_{\man{E}}})),
\enspace s\in\{1,\dots,m-1\},\\
C^{m+1}_0(\gamma_0)=&\;A^m_0\otimes\id_{\ctb{\man{E}}}(\gamma_0)+
(\nabla^{\man{E}}C^m_0)(\gamma_0)-C^m_0\otimes\id_{\ctb{\man{E}}}
(\Ins_1(\gamma_0,B_{\pi_{\man{E}}}))\\
&\;+C^m_0\otimes\id_{\ctb{\man{E}}}(\Ins_2(\gamma_0,\dual{B}_{\pi_{\man{E}}})),
\end{align*}}%
where\/ $\beta_s\in\tensor*[s]{\pi_{\man{E}}^*\ctb{\man{M}}}\otimes\vb{\man{E}}$\@,\/
$s\in\{0,1,\dots,m\}$\@, and\/
$\gamma_s\in\tensor*[s]{\pi_{\man{E}}^*\ctb{\man{M}}}\otimes
\tensor[1,1]{\vb{\man{E}}}$\@,\/ $s\in\{0,1,\dots,m-1\}$\@. 
\begin{proof}
The assertion is clearly true for $m=0$ and, for $m=1$\@, we have
\begin{equation*}
\nabla^{\man{E}}\ve{L}=
\ve{(\nabla^{\pi_{\man{E}}}L)}+\Ins_1(L,B_{\pi_{\man{E}}})+\vl{L}
\end{equation*}
by Lemma~\pldblref{lem:bundleform}{pl:bundleform7}\@, which gives the result
for $m=1$\@.  Thus suppose the result true for $m\ge2$ so that
\begin{equation*}
\nabla^{\man{E},m}\ve{L}=
\sum_{s=0}^mA^m_s(\ve{(\nabla^{\man{M},\pi_{\man{E}},s}L)})+
\sum_{s=0}^{m-1}C^m_s(\vl{(\nabla^{\man{M},\pi_{\man{E}},s}L)})
\end{equation*}
for vector bundle mappings $A^m_s$ and $C^m_s$ satisfying the stated
recursion relations.  We then compute
{\allowdisplaybreaks\begin{align*}
\nabla^{\man{E},m+1}\ve{L}=&\;
\sum_{s=0}^m(\nabla^{\man{E}}A^m_s)(\ve{(\nabla^{\man{M},\pi_{\man{E}},s}L)})+
\sum_{s=0}^mA^m_s\otimes\id_{\ctb{\man{E}}}
(\nabla^{\man{E}}\ve{(\nabla^{\man{M},\pi_{\man{E}},s}L)})\\
&\;+\sum_{s=0}^{m-1}(\nabla^{\man{E}}C^m_s)
(\vl{(\nabla^{\man{M},\pi_{\man{E}},s}L)})+\sum_{s=0}^{m-1}C^m_s\otimes
\id_{\ctb{\man{E}}}(\nabla^{\man{E}}\vl{(\nabla^{\man{M},\pi_{\man{E}},s}L)})\\
=&\;\sum_{s=0}^m(\nabla^{\man{E}}A^m_s)(\ve{(\nabla^{\man{M},\pi_{\man{E}},s}L)})+
\sum_{s=0}^mA^m_s\otimes\id_{\ctb{\man{E}}}
(\ve{(\nabla^{\man{M},\pi_{\man{E}},s+1}L)})\\
&\;-\sum_{s=1}^m\sum_{j=1}^sA^m_s\otimes\id_{\ctb{\man{E}}}
(\Ins_j(\ve{(\nabla^{\man{M},\pi_{\man{E}},s}L)},B_{\pi_{\man{E}}}))\\
&\;+\sum_{s=0}^mA^m_s\otimes\id_{\ctb{\man{E}}}
(\Ins_{s+1}(\ve{(\nabla^{\man{M},\pi_{\man{E}},s}L)},\dual{B}_{\pi_{\man{E}}}))
+\sum_{s=0}^mA^m_s\otimes\id_{\ctb{\man{E}}}
(\vl{(\nabla^{\man{M},\pi_{\man{E}},s}L)})\\
&\;+\sum_{s=0}^{m-1}(\nabla^{\man{E}}C^m_s)
(\vl{(\nabla^{\man{M},\pi_{\man{E}},s}L)})+
\sum_{s=0}^{m-1}C^m_s\otimes\id_{\ctb{\man{E}}}
(\vl{(\nabla^{\man{M},\pi_{\man{E}},s+1}L)})\\
&\;-\sum_{s=1}^{m-1}\sum_{j=1}^sC^m_s\otimes\id_{\ctb{\man{E}}}
(\Ins_j(\vl{(\nabla^{\man{M},\pi_{\man{E}},s}L)},B_{\pi_{\man{E}}}))\\
&\;+\sum_{s=0}^{m-1}C^m_s\otimes\id_{\ctb{\man{E}}}
(\Ins_{s+1}(\vl{(\nabla^{\man{M},\pi_{\man{E}},s}L)},\dual{B}_{\pi_{\man{E}}}))\\
=&\;\ve{(\nabla^{\man{M},\pi_{\man{E}},m+1}L)}+
\left(A^m_{m-1}\otimes\id_{\ctb{\man{E}}}
(\ve{(\nabla^{\man{M},\pi_{\man{E}},m}L)})-
\sum_{j=1}^m\Ins_j(\ve{(\nabla^{\man{M},\pi_{\man{E}},m}L)},
B_{\pi_{\man{E}}})\right.\\
&\;\left.\vphantom{\sum_{s=0}^m}+
\Ins_{m+1}(\ve{(\nabla^{\man{M},\pi_{\man{E}},m}L)},\dual{B}_{\pi_{\man{E}}})+
\vl{(\nabla^{\man{M},\pi_{\man{E}},m}L)}+C^m_{m-1}\otimes\id_{\ctb{\man{E}}}
(\vl{(\nabla^{\man{M},\pi_{\man{E}},m}L)})\right)\\
&\;+\left(\sum_{s=1}^{m-1}(\nabla^{\man{E}}A^m_s)
(\ve{(\nabla^{\man{M},\pi_{\man{E}},s}L)})+
\sum_{s=1}^{m-1}A^m_{s-1}\otimes\id_{\ctb{\man{E}}}
(\ve{(\nabla^{\man{M},\pi_{\man{E}},s}L)})\right.\\
&\;-\sum_{s=1}^{m-1}\sum_{j=1}^sA^m_s\otimes\id_{\ctb{\man{E}}}
(\Ins_j(\ve{(\nabla^{\man{M},\pi_{\man{E}},s}L)}),B_{\pi_{\man{E}}})\\
&\;+\sum_{s=1}^{m-1}A^m_s\otimes\id_{\ctb{\man{E}}}
(\Ins_{s+1}(\ve{(\nabla^{\man{M},\pi_{\man{E}},s}L)}),\dual{B}_{\pi_{\man{E}}})
+\sum_{s=1}^{m-1}A^m_s\otimes\id_{\ctb{\man{E}}}
(\vl{(\nabla^{\man{M},\pi_{\man{E}},s}L)})\\
&\;+\sum_{s=1}^{m-1}(\nabla^{\man{E}}C^m_s)
(\vl{(\nabla^{\man{M},\pi_{\man{E}},s}L)})+
\sum_{s=1}^{m-1}C^m_{s-1}\otimes\id_{\ctb{\man{E}}}
(\vl{(\nabla^{\man{M},\pi_{\man{E}},s}L)})\\
&\;-\sum_{s=0}^{m-1}\sum_{j=1}^sC^m_s\otimes\id_{\ctb{\man{E}}}
(\Ins_j(\vl{(\nabla^{\man{M},\pi_{\man{E}},s}L)},B_{\pi_{\man{E}}}))\\
&\;\left.+\sum_{s=1}^{m-1}C^m_s\otimes\id_{\ctb{\man{E}}}
(\Ins_{s+1}(\vl{(\nabla^{\man{M},\pi_{\man{E}},s}L)},
\dual{B}_{\pi_{\man{E}}}))\right)\\
&+\;(\nabla^{\man{E}}A^m_0)(\ve{L})+A^m_0\otimes\id_{\ctb{\man{E}}}
(\Ins_1(\ve{L},\dual{B}_{\pi_{\man{E}}}))+A^m_0\otimes\id_{\ctb{\man{E}}}
(\vl{L})\\
&+\;(\nabla^{\man{E}}C^m_0)(\vl{L})-C^m_0\otimes
\id_{\ctb{\man{E}}}(\Ins_1(\vl{L},B_{\pi_{\man{E}}}))+
C^m_0\otimes\id_{\ctb{\man{E}}}(\Ins_2(\vl{L},\dual{B}_{\pi_{\man{E}}})).
\end{align*}}
From these calculations, the lemma follows.
\end{proof}
\end{lemma}

Now we ``invert'' the constructions from the preceding lemma.
\begin{lemma}\label{lem:CiterateddersII}
Let\/ $r\in\{\infty,\omega\}$ and let\/
$\map{\pi_{\man{E}}}{\man{E}}{\man{M}}$ be a\/ $\C^r$-vector bundle, with the
data prescribed in Section~\ref{subsec:Esubmersion} to define the Riemannian
metric\/ $\metric_{\man{E}}$ on\/ $\man{E}$\@.  For\/ $m\in\integernn$\@,
there exist\/ $\C^r$-vector bundle mappings
\begin{equation*}
(B^m_s,\id_{\man{E}})\in
\vbmappings[r]{\tensor*[s]{\ctb{\man{E}}}\otimes\vb{\man{E}}}
{\tensor*[m]{\pi_{\man{E}}^*\ctb{\man{M}}}\otimes\vb{\man{E}}},
\qquad s\in\{0,1,\dots,m\},
\end{equation*}
and
\begin{equation*}
(D^m_s,\id_{\man{E}})\in
\vbmappings[r]{\tensor*[s]{\ctb{\man{E}}}\otimes\tensor[1,1]{\vb{\man{E}}}}
{\tensor*[m-1]{\pi_{\man{E}}^*\ctb{\man{M}}}\otimes\tensor[1,1]{\vb{\man{E}}}},
\qquad s\in\{0,1,\dots,m-1\},
\end{equation*}
such that
\begin{equation*}
\ve{(\nabla^{\man{M},\pi_{\man{E}},m}L)}=
\sum_{s=0}^mB^m_s(\nabla^{\man{E},s}\ve{L})+
\sum_{s=0}^{m-1}D^m_s(\nabla^{\man{E},s}\vl{L})
\end{equation*}
for all\/ $L\in\sections[m]{\tensor[1,1]{\man{E}}}$\@.  Moreover, the vector
bundle mappings\/ $B^m_0,B^m_1,\dots,B^m_m$ and\/
$D^m_0,D^m_1,\dots,D^m_{m-1}$ satisfy the recursion relations prescribed by\/
$B^0_0(\alpha_0)=\alpha_0$\@, $D^1_0(\gamma_0)=\gamma_0$\@,
\begin{align*}
B^{m+1}_{m+1}(\alpha_{m+1})=&\;\alpha_{m+1}\\
B^{m+1}_m(\alpha_m)=&\;B^m_{m-1}\otimes\id_{\ctb{\man{E}}}(\alpha_m)+
\sum_{j=1}^m\Ins_j(\alpha_m,B_{\pi_{\man{E}}})-
\Ins_{m+1}(\alpha,\dual{B}_{\pi_{\man{E}}})\\
B^{m+1}_s=&\;(\nabla^{\man{E}}B^m_s)(\alpha_s)+
B^m_{s-1}\otimes\id_{\ctb{\man{E}}}(\alpha_s)
+\sum_{j=1}^m\Ins_j(B^m_s(\alpha_s),B_{\pi_{\man{E}}})\\
&\;-\Ins_{m+1}(B^m_s(\alpha_s),\dual{B}_{\pi_{\man{E}}}),\enspace s\in\{1,\dots,m-1\},\\
B^{m+1}_0(\alpha_0)=&\;(\nabla^{\man{E}}B^m_0)(\alpha_0)+
\sum_{j=1}^m\Ins_j(B^m_0(\alpha_0),B_{\pi_{\man{E}}})-
\Ins_{m+1}(B^m_0(\alpha_0),\dual{B}_{\pi_{\man{E}}})
\end{align*}
and
\begin{align*}
D^{m+1}_m(\gamma_m)=&\;D^m_{m-1}\otimes\id_{\ctb{\man{E}}}(\gamma_m)-\gamma_m\\
D^m_s(\gamma_s)=&\;(\nabla^{\man{E}}D^m_s)(\gamma_s)+
D^m_{s-1}\otimes\id_{\ctb{\man{E}}}(\gamma_s)-\ol{B}^m_s(\gamma_s),
\enspace s\in\{1,\dots,m-1\},\\
D^{m+1}_0=&\;(\nabla^{\man{E}}D^m_0)(\gamma_0)-\ol{B}^m_0(\gamma_0)
\end{align*}
for $\alpha_s\in\tensor*[s]{\ctb{\man{E}}\otimes\vb{\man{E}}}$\@,
$s\in\{0,1,\dots,m+1\}$\@, and
$\gamma_s\in\tensor*[s]{\ctb{\man{E}}}\otimes\tensor[1,1]{\vb{\man{E}}}$\@,
$s\in\{0,1,\dots,m\}$\@, and where
\begin{equation*}
(\ol{B}^m_s,\id_{\man{E}})\in
\vbmappings[r]{\tensor*[s]{\ctb{\man{E}}}\otimes\tensor[1,1]{\vb{\man{E}}}}
{\tensor*[m]{\pi_{\man{E}}^*\ctb{\man{M}}}\otimes\tensor[1,1]{\vb{\man{E}}}},
\qquad s\in\{0,1,\dots,m\},
\end{equation*}
are the vector bundle mappings from Lemma~\ref{lem:LiterateddersII}\@.
\begin{proof}
The assertion is clearly true for $m=0$\@, so suppose it true for
$m\in\integerp$\@.  Thus
\begin{equation}\label{eq:CiteratedderII1}
\ve{(\nabla^{\man{M},\pi_{\man{E}},m}L)}=
\sum_{s=0}^mB^m_s(\nabla^{\man{E},s}\ve{L})+
\sum_{s=0}^{m-1}D^m_s(\vl{(\nabla^{\man{E},s}L)}).
\end{equation}
Working on the left-hand side of this equation, using
Lemma~\pldblref{lem:bundleform}{pl:bundleform7}\@, we have
\begin{align*}
\nabla^{\man{E}}\ve{(\nabla^{\man{M},\pi_{\man{E}},m}L)}=&\;
\ve{(\nabla^{\man{M},\pi_{\man{E}},m+1}L)}-
\sum_{j=1}^m\Ins_j(\ve{(\nabla^{\man{M},\pi_{\man{E}},m}L)},B_{\pi_{\man{E}}})\\
&\;+\Ins_{m+1}(\ve{(\nabla^{\man{M},\pi_{\man{E}},m}L)},\dual{B}_{\pi_{\man{E}}})+
\vl{(\nabla^{\man{M},\pi_{\man{E}},m}L)}\\
=&\;\ve{(\nabla^{\man{M},\pi_{\man{E}},m+1}L)}-\sum_{s=0}^m\sum_{j=1}^m
\Ins_j(B^m_s(\nabla^{\man{E},s}\ve{L}),B_{\pi_{\man{E}}})\\
&\;+\sum_{s=0}^m\Ins_{m+1}(B^m_s(\nabla^{\man{E},s}\ve{L}),
\dual{B}_{\pi_{\man{E}}})+\sum_{s=0}^m\ol{B}^m_s(\nabla^{\man{E},s}\vl{L}).
\end{align*}
Working on the right-hand side of~\eqref{eq:CiteratedderII1}\@,
\begin{multline*}
\nabla^{\man{E}}\ve{(\nabla^{\man{M},\pi_{\man{E}},m}L)}=
\sum_{s=0}^m(\nabla^{\man{E}}B^m_s)(\nabla^{\man{E},s}\ve{L})+
\sum_{s=0}^mB^m_s\otimes\id_{\ctb{\man{E}}}(\nabla^{\man{E},s+1}\ve{L})\\
+\sum_{s=0}^{m-1}(\nabla^{\man{E}}D^m_s)(\nabla^{\man{E},s}\vl{L})+
\sum_{s=0}^{m-1}D^m_s\otimes\id_{\ctb{\man{E}}}(\nabla^{\man{E},s+1}\vl{L}).
\end{multline*}
Combining the preceding two computations,
{\small\allowdisplaybreaks\begin{align*}
\ve{(\nabla^{\man{M},\pi_{\man{E}},m+1}L)}=&\;
\sum_{s=0}^m(\nabla^{\man{E}}B^m_s)(\nabla^{\man{E},s}\ve{L})+
\sum_{s=0}^mB^m_s\otimes\id_{\ctb{\man{E}}}(\nabla^{\man{E},s+1}\ve{L})\\
&\;+\sum_{s=0}^{m-1}(\nabla^{\man{E}}D^m_s)(\nabla^{\man{E},s}\vl{L})+
\sum_{s=0}^{m-1}D^m_s\otimes\id_{\ctb{\man{E}}}(\nabla^{\man{E},s+1}\vl{L})\\
&\;+\sum_{s=0}^m\sum_{j=1}^m\Ins_j(B^m_s(\nabla^{\man{E},s}\ve{L}),
B_{\pi_{\man{E}}})-
\sum_{s=1}^m\Ins_{m+1}(B^m_s(\nabla^{\man{E},s}\ve{L}),\dual{B}_{\pi_{\man{E}}})\\
&\;-\sum_{s=0}^m\ol{B}^m_s(\nabla^{\man{E},s}\vl{L})\\
=&\;\nabla^{\man{E},m+1}\ve{L}+\left(B^m_{m-1}\otimes\id_{\ctb{\man{E}}}
(\nabla^{\man{E},m}\ve{L})+D^m_{m-1}\otimes\id_{\ctb{\man{E}}}
(\nabla^{\man{E},m}\vl{L})\vphantom{\sum_{s=0}^m}\right.\\
&\;\left.+\sum_{j=1}^m\Ins_j(\nabla^{\man{E},m}\ve{L},B_{\pi_{\man{E}}})
-\Ins_{m+1}(\nabla^{\man{E},m}\ve{L},\dual{B}_{\pi_{\man{E}}})-
(\nabla^{\man{E},m}\vl{L})\right)\\
&\;+\left(\sum_{s=1}^{m-1}(\nabla^{\man{E}}B^m_s)
(\nabla^{\man{E},s}\ve{L})\right.
+\sum_{s=1}^{m-1}B^m_{s-1}\otimes\id_{\ctb{\man{E}}}(\nabla^{\man{E},s}\ve{L})+
\sum_{s=1}^{m-1}(\nabla^{\man{E}}D^m_s)(\nabla^{\man{E},s}\vl{L})\\
&\;+\sum_{s=1}^{m-1}D^m_{s-1}\otimes\id_{\ctb{\man{E}}}(\nabla^{\man{E},s}\vl{L})+
\sum_{s=1}^{m-1}\sum_{j=1}^m\Ins_j(B^m_s(\nabla^{\man{E},s}\ve{L}),
B_{\pi_{\man{E}}})\\
&\;-\sum_{s=1}^{m-1}\Ins_{m+1}(B^m_s(\nabla^{\man{E},s}\ve{L}),
\dual{B}_{\pi_{\man{E}}})
\left.-\sum_{s=1}^{m-1}\ol{B}^m_s(\nabla^{\man{E},s}\vl{L})\right)\\
&\;+\left((\nabla^{\man{E}}B^m_0)(\ve{L})+(\nabla^{\man{E}}D^m_0)(\vl{L})+
\sum_{j=1}^m\Ins_j(B^m_0(\ve{L}),B_{\pi_{\man{E}}})\right.\\
&\;\left.\vphantom{\sum_{s=0}^m}-\Ins_{m+1}(B^m_0(\ve{L}),
\dual{B}_{\pi_{\man{E}}})-\ol{B}^m_0(\vl{L})\right).
\end{align*}}%
The lemma follows from these computations.
\end{proof}
\end{lemma}

Next we turn to symmetrised versions of the preceding lemmata.  We show that
the preceding two lemmata induce corresponding mappings between symmetric
tensors.
\begin{lemma}\label{lem:CsymiterateddersI}
Let\/ $r\in\{\infty,\omega\}$ and let\/
$\map{\pi_{\man{E}}}{\man{E}}{\man{M}}$ be a\/ $\C^r$-vector bundle, with the
data prescribed in Section~\ref{subsec:Esubmersion} to define the Riemannian
metric\/ $\metric_{\man{E}}$ on\/ $\man{E}$\@.  For\/ $m\in\integernn$\@,
there exist\/ $\C^r$-vector bundle mappings
\begin{equation*}
(\what{A}^m_s,\id_{\man{E}})\in
\vbmappings[r]{\Symalg*[s]{\pi_{\man{E}}^*\ctb{\man{M}}}\otimes\vb{\man{E}}}
{\Symalg*[m]{\ctb{\man{E}}}\otimes\vb{\man{E}}},\qquad s\in\{0,1,\dots,m\},
\end{equation*}
and
\begin{equation*}
(\what{C}^m_s,\id_{\man{E}})\in
\vbmappings[r]{\Symalg*[s]{\pi_{\man{E}}^*\ctb{\man{M}}}\otimes
\tensor[1,1]{\vb{\man{E}}}}
{\Symalg*[m]{\ctb{\man{E}}}\otimes\vb{\man{E}}},\qquad
s\in\{0,1,\dots,m-1\},
\end{equation*}
such that
\begin{multline*}
(\Sym_m\otimes\id_{\vb{\man{E}}})\scirc\nabla^{\man{E},m}\ve{L}\\
=\sum_{s=0}^m\what{A}^m_s((\Sym_s\otimes\id_{\vb{\man{E}}})\scirc
\ve{(\nabla^{\man{M},\pi_{\man{E}},s}L)}+
\sum_{s=0}^{m-1}\what{C}^m_s((\Sym_s\otimes\id_{\tensor[1,1]{\vb{\man{E}}}})
\scirc\vl{(\nabla^{\man{M},\pi_{\man{E}},s}L)})
\end{multline*}
for all\/ $L\in\sections[m]{\tensor[1,1]{\man{E}}}$\@.
\begin{proof}
Following along the lines of the proof of
Lemma~\ref{lem:VsymiterateddersI}\@, we define $\what{A}^m_s$ by requiring
that
\begin{equation*}
\what{A}^m_s((\Sym_s\otimes\id_{\vb{\man{E}}})\scirc
\ve{(\nabla^{\man{M},\pi_{\man{E}},s}L)})=
(\Sym_m\otimes\id_{\vb{\man{E}}})\scirc A^m_s(\ve{(\nabla^{\man{M},\pi_{\man{E}},s}L)}),
\end{equation*}
and $\what{C}^m_s$ by requiring that
\begin{equation*}
\what{C}^m_s((\Sym_s\otimes\id_{\tensor[1,1]{\vb{\man{E}}}})\scirc
\ve{(\nabla^{\man{M},\pi_{\man{E}},s}L)})=(\Sym_m\otimes\id_{\vb{\man{E}}})
\scirc C^m_s(\ve{(\nabla^{\man{M},\pi_{\man{E}},s}L)}).
\end{equation*}
That this definition of $\what{A}^m_s$ makes sense follows exactly as in the
proof of Lemma~\ref{lem:VsymiterateddersI}\@.  Let us see how the same
arguments also apply to the definition of $\what{C}^m_s$\@.

For $m\in\integerp$\@, we define
$\map{C^m}{\tensor*[\le m-1]{\pi_{\man{E}}^*\ctb{\man{M}}}\otimes
\tensor[1,1]{\vb{\man{E}}}}{\tensor*[\le
m]{\ctb{\man{E}}}\otimes\vb{\man{E}}}$ by
\begin{multline*}
C^m(\vl{L},\vl{(\nabla^{\pi_{\man{E}}}L)},\dots,
\vl{(\nabla^{\man{M},\pi_{\man{E}},m-1}L)})\\
=\left(C^1_0(\vl{L}),\sum_{s=0}^1C^2_s(\vl{(\nabla^{\man{M},\pi_{\man{E}},s}L)}),
\dots,\sum_{s=0}^{m-1}C^m_s(\vl{(\nabla^{\man{M},\pi_{\man{E}},s}L)})\right),
\end{multline*}
making the identification of $\tensor[1,1]{\vb{\man{E}}}$ with a subspace of
$\ctb{\man{E}}\otimes\vb{\man{E}}$ as in the footnote from
Lemma~\ref{fn:T11inT*}\@.  Note that we have a natural mapping
\begin{equation*}
\ctb{\man{E}}\otimes\jetalg{m-1}{\man{E}}\to\jetalg{m}{\man{E}}
\end{equation*}
\cf~\cite[Theorem~6.2.9]{DJS:89}\@.  This then induces a mapping
\begin{equation*}
\map{P_m}{(\real_{\man{E}}\oplus\jetalg{m-1}{\man{E}})\otimes
\tensor[1,1]{\vb{\man{E}}}}{(\real_{\man{E}}\otimes\jetalg{m}{\man{E}})
\otimes\vb{\man{E}}}.
\end{equation*}
Now define
\begin{equation*}
\map{\what{P}_m}{\pi_{\man{E}}^*(\real_{\man{M}}\oplus
\jetalg{m-1}{\man{M}})\otimes
\tensor[1,1]{\vb{\man{E}}}}{(\real_{\man{E}}\oplus\jetalg{m}{\man{E}})
\otimes\vb{\man{E}}}
\end{equation*}
by
\begin{equation*}
\what{P}_m=P_m\scirc((\id_{\real}\oplus j_{m-1}\pi_{\man{E}})
\otimes\id_{\tensor[1,1]{\vb{\man{E}}}}),
\end{equation*}
noting that
\begin{equation*}
\map{\id_{\real}\oplus j_{m-1}\pi_{\man{E}}}
{\pi_{\man{E}}^*(\real_{\man{M}}\jetalg{m-1}{\man{M}})}
{\real_{\man{E}}\oplus\jetalg{m-1}{\man{E}}}
\end{equation*}
is injective.  Also define
\begin{equation*}
\map{Q_m}{\Symalg*[\le m-1]{\ctb{\man{E}}}\otimes\tensor[1,1]{\vb{\man{E}}}}
{\Symalg*[m]{\ctb{\man{E}}}\otimes\vb{\man{E}}}
\end{equation*}
by
\begin{multline*}
Q_m(A_0\otimes\alpha_0\otimes u_0,\dots,A_{m-1}\otimes\alpha_{m-1}\otimes u_{m-1})\\
=(\Sym_1(A_0\otimes\alpha_0)\otimes u_0,\dots,
\Sym_m(A_{m-1}\otimes\alpha_{m-1})\otimes u_{m-1}).
\end{multline*}
Note that the diagram
\begin{equation*}
\xymatrix{{\Symalg*[\le m-1]{\ctb{\man{E}}}\otimes\tensor[1,1]{\vb{\man{E}}}}
\ar[rr]^(0.48){S^{m-1}_{\nabla^{\man{E}}}\otimes\id_{\tensor[1,1]{\vb{\man{E}}}}}
\ar[d]_{Q_m}&&
{(\real_{\man{E}}\oplus\jetalg{m-1}{\man{E}})
\otimes\tensor[1,1]{\vb{\man{E}}}}\ar[d]^{P_m}\\
{\Symalg*[\le m]{\ctb{\man{E}}}\otimes\vb{\man{E}}}
\ar[rr]_{S^m_{\nabla^{\man{E}}}\otimes\id_{\vb{\man{E}}}}&&
{(\real_{\man{E}}\oplus\jetalg{m}{\man{E}})\otimes\vb{\man{E}}}}
\end{equation*}
commutes.  We also define
\begin{equation*}
\what{Q}_m=Q_m\scirc(\pi^*_{m-1}\otimes\id_{\tensor[1,1]{\vb{\man{E}}}}),
\end{equation*}
where
\begin{equation*}
\map{\pi^*_{m-1}}{\Symalg*[\le m-1]{\pi_{\man{E}}^*\ctb{\man{M}}}}
{\Symalg*[\le m-1]{\ctb{\man{E}}}}
\end{equation*}
is the inclusion.  Note that the diagram
\begin{equation*}
\xymatrix{{\Symalg*[\le m-1]{\pi_{\man{E}}^*\ctb{\man{M}}}}\ar[r]^{\pi^*_{m-1}}
\ar[d]_{S^{m-1}_{\nabla^{\man{M}},\nabla^{\pi_{\man{E}}}}}&{\Symalg*[\le m-1]{\ctb{\man{E}}}}
\ar[d]^{S^{m-1}_{\nabla^{\man{E}}}}\\
{\pi_{\man{E}}^*(\real_{\man{M}}\oplus\jetalg{m-1}{\man{M}})}
\ar[r]_(0.52){\id_{\real}\oplus j_{m-1}\pi_{\man{E}}}&
{\real_{\man{E}}\oplus\jetalg{m-1}{\man{E}}}}
\end{equation*}
commutes.

Let us organise the mappings we require into the following diagram:
{\tiny\begin{equation}\label{eq:Csymmetrise1}
\xymatrix{{\tensor*[\le m-1]{\pi_{\man{E}}^*\ctb{\man{M}}}\otimes
\tensor[1,1]{\vb{\man{E}}}}
\ar[rr]^{\Sym_{\le m-1}\otimes\id_{\tensor[1,1]{\vb{\man{E}}}}}\ar[d]_{C^m}&&
{\Symalg*[\le m-1]{\pi_{\man{E}}^*\ctb{\man{M}}}\otimes
\tensor[1,1]{\vb{\man{E}}}}
\ar[rr]^(0.48){S^{m-1}_{\nabla^{\man{M}},\nabla^{\pi_{\man{E}}}}\otimes
\id_{\tensor[1,1]{\vb{\man{E}}}}}\ar[d]^{\what{C}^m}&&
{\pi_{\man{E}}^*(\real_{\man{M}}\oplus\jetalg{m-1}{\man{M}})\otimes
\tensor[1,1]{\vb{\man{E}}}}\ar[d]^{\what{P}_m}\\
{\tensor*[\le m]{\ctb{\man{E}}}\otimes\vb{\man{E}}}
\ar[rr]^{\Sym_{\le m}\otimes\id_{\vb{\man{E}}}}&&
{\Symalg*[\le m]{\ctb{\man{E}}}\otimes\vb{\man{E}}}
\ar[rr]^(0.5){S^m_{\nabla^{\man{E}}}\otimes\id_{\vb{\man{E}}}}&&
{(\real_{\man{E}}\oplus\jetalg{m}{\man{E}})\otimes\vb{\man{E}}}}
\end{equation}}%
Here $\what{C}^m$ is defined so that the right square commutes, which is
possible since the horizontal arrows in the right square are isomorphisms.
We shall show that the left square also commutes.  Indeed,
\begin{align*}
\what{C}^m\scirc(\Sym_{\le m-1}\otimes\id_{\tensor[1,1]{\vb{\man{E}}}})&
(\vl{L},\vl{(\nabla^{\pi_{\man{E}}}L)},\dots,
\vl{(\nabla^{\man{M},\pi_{\man{E}},m}L)})\\
=&\;(S^m_{\nabla^{\man{E}}}\otimes\id_{\vb{\man{E}}})^{-1}\scirc
\what{P}_m\scirc (S^{m-1}_{\nabla^{\man{M}},\nabla^{\pi_{\man{E}}}}\otimes
\id_{\tensor[1,1]{\vb{\man{E}}}})\\
&\;\scirc(\Sym_{\le m-1}\otimes\id_{\tensor[1,1]{\vb{\man{E}}}})
(\vl{L},\vl{(\nabla^{\pi_{\man{E}}}L)},\dots,
(\nabla^{\man{M},\pi_{\man{E}},m}\vl{L})\\
=&\;(\Sym_{\le m-1}\otimes\id_{\vb{\man{E}}})
(\vl{L},\nabla^{\man{E}}\vl{L},\dots,\nabla^{\man{E},m}\vl{L})\\
=&\;(\Sym_{\le m}\otimes\id_{\vb{\man{E}}})\scirc
C^m(\vl{L},\vl{(\nabla^{\pi_{\man{E}}}L)},\dots,
\vl{(\nabla^{\man{M},\pi_{\man{E}},m}L)}).
\end{align*}
Thus the diagram~\eqref{eq:Csymmetrise1} commutes.  Thus, if we define
\begin{equation*}
\what{C}^m_s((\Sym_s\otimes\id_{\tensor[1,1]{\vb{\man{E}}}})
\scirc\vl{(\nabla^{\man{M},\pi_{\man{E}},s}L)})=
(\Sym_m\otimes\id_{\vb{\man{E}}})\scirc C^m_s(\vl{(\nabla^{\man{M},\pi_{\man{E}},s}L)}),
\end{equation*}
then we have
\begin{equation*}
(\Sym_m\otimes\id_{\vb{\man{E}}})\scirc\nabla^{\man{E},m}\ve{L}=
\sum_{s=0}^m\what{C}^m_s((\Sym_s\otimes\id_{\tensor[1,1]{\vb{\man{E}}}})
\scirc\vl{(\nabla^{\man{M},\pi_{\man{E}},s}L)}),
\end{equation*}
as desired.
\end{proof}
\end{lemma}

The preceding lemma gives rise to an ``inverse,'' which we state in the
following lemma.
\begin{lemma}\label{lem:CsymiterateddersII}
Let\/ $r\in\{\infty,\omega\}$ and let\/
$\map{\pi_{\man{E}}}{\man{E}}{\man{M}}$ be a\/ $\C^r$-vector bundle, with the
data prescribed in Section~\ref{subsec:Esubmersion} to define the Riemannian
metric\/ $\metric_{\man{E}}$ on\/ $\man{E}$\@.  For\/ $m\in\integernn$\@,
there exist\/ $\C^r$-vector bundle mappings
\begin{equation*}
(\what{B}^m_s,\id_{\man{E}})\in
\vbmappings[r]{\Symalg*[s]{\ctb{\man{E}}}\otimes\vb{\man{E}}}
{\Symalg*[m]{\pi_{\man{E}}^*\ctb{\man{M}}}\otimes\vb{\man{E}}},\qquad
s\in\{0,1,\dots,m\},
\end{equation*}
and
\begin{equation*}
(\what{D}^m_s,\id_{\man{E}})\in
\vbmappings[r]{\Symalg*[s]{\ctb{\man{E}}}\otimes\tensor[1,1]{\vb{\man{E}}}}
{\Symalg*[m]{\pi_{\man{E}}^*\ctb{\man{M}}}\otimes\vb{\man{E}}},\qquad s\in\{0,1,\dots,m-1\},
\end{equation*}
such that
\begin{multline*}
(\Sym_m\otimes_{\id_{\vb{\man{E}}}})\scirc
\ve{(\nabla^{\man{M},\pi_{\man{E}},m}L)}\\
=\sum_{s=0}^m\what{B}^m_s((\Sym_s\otimes\id_{\vb{\man{E}}})\scirc
\nabla^{\man{E},s}\ve{L})+
\sum_{s=0}^{m-1}\what{D}^m_s((\Sym_s\otimes
\id_{\tensor[1,1]{\vb{\man{E}}}})\scirc\nabla^{\man{E},s}\vl{L})
\end{multline*}
for all\/ $L\in\sections[m]{\tensor[1,1]{\man{E}}}$\@.
\begin{proof}
Following along the lines of the proof of
Lemma~\ref{lem:VsymiterateddersI}\@, we define $\what{B}^m_s$ by requiring
that
\begin{equation*}
\what{B}^m_s((\Sym_s\otimes\id_{\vb{\man{E}}})\scirc\nabla^{\man{E},s}\ve{L}=
(\Sym_m\otimes\id_{\vb{\man{E}}})\scirc B^m_s(\nabla^{\man{E},s}\ve{L}),
\end{equation*}
and $\what{C}^m_s$ by requiring that
\begin{equation*}
\what{C}^m_s((\Sym_s\otimes\id_{\tensor[1,1]{\vb{\man{E}}}})\scirc
\nabla^{\man{E},s}\vl{L})=(\Sym_m\otimes\id_{\vb{\man{E}}})
\scirc C^m_s(\nabla^{\man{E},s}\vl{L}).
\end{equation*}
That these definitions make sense follows along the same lines as the proof
of Lemma~\ref{lem:CsymiterateddersI}\@.
\end{proof}
\end{lemma}

We can put together the previous four lemmata into the following
decomposition result, which is to be regarded as the main result of this
section.
\begin{lemma}\label{lem:Cdecomp}
Let\/ $r\in\{\infty,\omega\}$ and let\/
$\map{\pi_{\man{E}}}{\man{E}}{\man{M}}$ be a\/ $\C^r$-vector bundle, with the
data prescribed in Section~\ref{subsec:Esubmersion} to define the Riemannian
metric\/ $\metric_{\man{E}}$ on\/ $\man{E}$\@.  Then there exist\/
$\C^r$-vector bundle mappings
\begin{equation*}
A^m_{\nabla^{\man{E}}}\in\vbmappings[r]{\Cjetalg{m}{\man{E}}}
{\Symalg*[\le m]{\pi_{\man{E}}^*\ctb{\man{M}}}\otimes\vb{\man{E}}},\quad
B^m_{\nabla^{\man{E}}}\in\vbmappings[r]{\Cjetalg{m}{\man{E}}}
{\Symalg*[\le m]{\ctb{\man{E}}}\otimes\vb{\man{E}}}
\end{equation*}
defined by
\begin{align*}
A^m_{\nabla^{\man{E}}}(j_m(\ve{L})(e))=&\;
\Sym_{\le m}\otimes\id_{\vb{\man{E}}}
(\ve{L}(e),\ve{(\nabla^{\pi_{\man{E}}}L)}(e),\dots,
\ve{(\nabla^{\man{M},\pi_{\man{E}},m}L)}(e)),\\
B^m_{\nabla^{\man{E}}}(j_m(\ve{L})(e))=&\;
\Sym_{\le m}\otimes\id_{\vb{\man{E}}}
(\ve{L}(e),\nabla^{\man{E}}\ve{L}(e),\dots,\nabla^{\man{E},m}\ve{L}(e)).
\end{align*}
Moreover,\/ $A^m_{\nabla^{\man{E}}}$ and\/ $B^m_{\nabla^{\man{E}}}$ are
injective, and
{\small\begin{multline*}
B^m_{\nabla^{\man{E}}}\scirc(A^m_{\nabla^{\man{E}}})^{-1}\scirc
(\Sym_{\le m}\otimes\id_{\vb{\man{E}}})
(\ve{L}(e),\ve{(\nabla^{\pi_{\man{E}}}L)}(e),\dots,
\ve{(\nabla^{\man{M},\pi_{\man{E}},m}L)}(e))\\
=\left(\ve{L}(e),\sum_{s=0}^1\what{A}^1_s((\Sym_s\otimes\id_{\vb{\man{E}}})
\scirc\ve{(\nabla^{\man{M},\pi_{\man{E}},s}L)}(e)),\dots,\right.\\
\left.\sum_{s=0}^m\what{A}^m_s((\Sym_s\otimes\id_{\vb{\man{E}}})
\scirc\ve{(\nabla^{\man{M},\pi_{\man{E}},s}L)}(e))\right)\\
+\left(0,\vl{L}(e),\sum_{s=0}^1
\what{C}^2_s((\Sym_s\otimes\id_{\tensor[1,1]{\vb{\man{E}}}})
\scirc\vl{(\nabla^{\man{M},\pi_{\man{E}},s}L)}(e)),\dots,\right.\\
\left.\sum_{s=0}^{m-1}\what{C}^m_s((\Sym_s\otimes\id_{\tensor[1,1]{\vb{\man{E}}}})
\scirc\vl{(\nabla^{\man{M},\pi_{\man{E}},s}L)}(e))\right)
\end{multline*}}%
and
{\small\begin{multline*}
A^m_{\nabla^{\man{E}}}\scirc(B^m_{\nabla^{\man{E}}})^{-1}\scirc
(\Sym_{\le m}\otimes\id_{\vb{\man{E}}})
(\ve{L}(e),\nabla^{\man{E}}\ve{L}(e),\dots,\nabla^{\man{E},m}\ve{L}(e))\\
=\left(\ve{L}(e),\sum_{s=0}^1\what{B}^1_s((\Sym_s\otimes\id_{\vb{\man{E}}})
\scirc\nabla^{\man{E},s}\ve{L}(e)),\dots,
\sum_{s=0}^m\what{B}^m_s((\Sym_s\otimes\id_{\vb{\man{E}}})\scirc
\nabla^{\man{E},s}\ve{L}(e))\right)\\
+\left(0,\vl{L}(e),\sum_{s=0}^1\what{D}^2_s((\Sym_s\otimes
\id_{\tensor[1,1]{\vb{\man{E}}}})\scirc\nabla^{\man{E},s}\vl{L}(e)),
\dots,\right.\\
\left.\sum_{s=0}^{m-1}\what{D}^m_s((\Sym_s\otimes\id_{\tensor[1,1]{\vb{\man{E}}}})
\scirc\nabla^{\man{E},s}\vl{L}(e))\right),
\end{multline*}}%
where the vector bundle mappings\/ $\what{A}^m_s$ and\/ $\what{B}^m_s$\@,\/
$s\in\{0,1,\dots,m\}$\@, and\/ $\what{C}^m_s$ and\/ $\what{D}^m_s$\@,
$s\in\{0,1,\dots,m-1\}$\@, are as in Lemmata~\ref{lem:CsymiterateddersI}
and~\ref{lem:CsymiterateddersII}\@.
\end{lemma}

\subsection{Isomorphisms for pull-backs of functions}

Next we generalise the presentation of Section~\ref{subsec:Pjetisom} from the
pull-back of a vector bundle projection to the pull-back by a general
mapping.  The development here is a little different from the preceding
sections, so we first have a little bit of setting up to do.  For
$\C^r$-manifolds $\man{M}$ and $\man{N}$\@, and for
$\Phi\in\mappings[r]{\man{M}}{\man{N}}$\@, we consider the mapping
\begin{equation*}
\func[r]{\man{N}}\ni f\mapsto\Phi^*f\in\func[r]{\man{M}}.
\end{equation*}
We wish to compare the decomposition of jets of $f$ with those of
$\Phi^*f$\@, and to do so we consider the subbundle
$\jetalg[\Phi]{m}{\man{M}}$ of $\jetalg{m}{\man{M}}$ defined by
\begin{equation*}
\jetalg[\Phi,x]{m}{\man{M}}=\setdef{j_m(\Phi^*f)(x)}{f\in\func[m]{\man{N}}}.
\end{equation*}
Following Lemma~\ref{lem:Jmdecomp}\@, we shall give a formula for iterated
covariant differentials of pull-backs of functions on $\man{N}$\@.  To do
this, we let $\nabla^{\man{M}}$ and $\nabla^{\man{N}}$ be affine connections
on $\man{M}$ and $\man{N}$\@, respectively.  We note that we have the vector
bundle connection $\Phi^*\nabla^{\man{N}}$ in the vector bundle
$\Phi^*\tb{\man{N}}$ over $\man{M}$\@.  Explicitly,
\begin{equation*}
(\Phi^*\nabla^{\man{N}}_X\Phi^*Y)(x)=
(x,\nabla^{\man{N}}_{\tf[x]{\Phi}(X(x))}Y).
\end{equation*}
Following our usual mild notational abuse, we shall also denote by
$\Phi^*\nabla^{\man{N}}$ the connection in the dual bundle
$\dual{(\Phi^*\tb{\man{N}})}\simeq\Phi^*\ctb{\man{N}}$\@.  We have a natural
mapping
\begin{equation*}
\mapdef{\what{\Phi}}{\tb{\man{M}}}{\Phi^*\tb{\man{N}}}
{v_x}{(x,\tf[x]{\Phi}(v_x)).}
\end{equation*}
This mapping induces a mappings on sections which we denote by the same
symbol; thus we have the mapping
\begin{equation*}
\map{\what{\Phi}}{\sections[\infty]{\tb{\man{M}}}}
{\sections[\infty]{\Phi^*\tb{\man{N}}}}.
\end{equation*}

The following lemma gives an important tensor for our analysis.
\begin{lemma}\label{lem:APhidef}
Let\/ $r\in\{\infty,\omega\}$\@.  Let\/ $\man{M}$ and\/ $\man{N}$ be\/
$\C^r$-manifolds and let\/ $\nabla^{\man{M}}$ and\/ $\nabla^{\man{N}}$ be\/
$\C^r$-affine connections on\/ $\man{M}$ and\/ $\man{N}$\@, respectively.
Let\/ $\Phi\in\mappings[r]{\man{M}}{\man{N}}$\@.  Then there exists\/
$A_\Phi\in\sections[r]{\tensor*[2]{\ctb{\man{M}}}\otimes
\Phi^*\tb{\man{N}}}$ such that, for\/ $x\in\man{M}$\@,
\begin{equation*}
\what{\Phi}(\nabla^{\man{M}}_XY)(x)-
\Phi^*\nabla^{\man{N}}_X\what{\Phi}(Y)(x)=A_\Phi(X(x),Y(x))
\end{equation*}
for\/ $X,Y\in\sections[\infty]{\tb{\man{M}}}$\@.
\begin{proof}
Let $\map{K^{\man{M}}}{\tb{\tb{\man{M}}}}{\tb{\man{M}}}$ and
$\map{K^{\man{N}}}{\tb{\tb{\man{N}}}}{\tb{\man{N}}}$ be the connectors for
$\nabla^{\man{M}}$ and $\nabla^{\man{N}}$ so that
\begin{equation*}
\nabla^{\man{M}}_XY=K^{\man{M}}\scirc\tf{Y}\scirc X,\qquad
X,Y\in\sections[\infty]{\tb{\man{M}}},
\end{equation*}
and
\begin{equation*}
\nabla^{\man{N}}_UV=K^{\man{N}}\scirc\tf{V}\scirc U,\qquad
U,V\in\sections[\infty]{\tb{\man{N}}}.
\end{equation*}
We, moreover, have
\begin{equation*}
\what{\Phi}(\nabla^{\man{M}}_XY)=\tf{\Phi}\scirc
K^{\man{M}}\scirc\tf{Y}\scirc X,\qquad
X,Y\in\sections[\infty]{\tb{\man{M}}},
\end{equation*}
and
\begin{equation*}
\Phi^*\nabla^{\man{N}}_X\what{\Phi}(Y)=
K^{\man{N}}\scirc\tf{(\tf{\Phi}\scirc Y)}\scirc X,\qquad
X,Y\in\sections[\infty]{\tb{\man{M}}}
\end{equation*}
\cite[\S{}10.12]{PWM:08}\@.  In preparation to use these formulae, we have
the following results.
\begin{proofsublemma}\label{psublem:T(fxi)}
If\/ $\map{\pi_{\man{E}}}{\man{E}}{\man{M}}$ is a smooth vector bundle, if\/
$\xi\in\sections[\infty]{\man{E}}$\@, and if\/
$f\in\func[\infty]{\man{M}}$\@, then
\begin{equation*}
\tf[x]{(f\xi)}(v_x)=f(x)\tf[x]{\xi}(v_x)+\natpair{\d{f}(x)}{v_x}\vl{\xi}(x).
\end{equation*}
\begin{subproof}
Let $\nabla^{\pi_{\man{E}}}$ be a linear connection in the vector bundle
$\man{E}$ which gives the decomposition
$\tb{\man{E}}=\hb{\man{E}}\oplus\vb{\man{E}}$\@.  Let $\hor$ and $\ver$ be
the horizontal and vertical projections.  Let $v_x\in\tb[x]{\man{M}}$ and let
$\map{\gamma}{I}{\man{M}}$ be a smooth curve for which $\gamma'(0)=v_x$\@.
Denote $\Xi(t)=f\scirc\gamma(t)\xi\scirc\gamma(t)$ the corresponding curve in
$\man{E}$\@.  Then
\begin{equation*}
\hor(\Xi'(t))=\horlift(f\scirc\gamma(t)\xi\scirc\gamma(t),\gamma'(t)),\quad
\ver(\Xi'(t))=\verlift(f\scirc\gamma(t)\xi\scirc\gamma(t),
\nabla^{\pi_{\man{E}}}_{\gamma'(t)}\Xi(t)).
\end{equation*}
We now have
\begin{equation*}
\nabla^{\pi_{\man{E}}}_{\gamma'(t)}\Xi(t)=
f\scirc\gamma(t)\nabla^{\pi_{\man{E}}}_{\gamma'(t)}\xi\scirc\gamma(t)+
\natpair{\d{f}\scirc\gamma(t)}{\gamma'(t)}\xi\scirc\gamma(t).
\end{equation*}
Thus
\begin{align*}
\tf[x]{(f\xi)}(v_x)=&\;\derivatzero{}{t}f\scirc\gamma(t))\xi\scirc\gamma(t)\\
=&\;f(x)\horlift(f(x)\xi(x),v_x)+\verlift(f(x)\xi(x),
f(x)\nabla^{\pi_{\man{E}}}_{v_x}\xi+\natpair{\d{f}(x)}{v_x}\xi(x))\\
=&\;f(x)\Xi'(0)+\natpair{\d{f}(x)}{v_x}\xi(x)=f(x)\tf[x]{\xi}(v_x)+
\natpair{\d{f}(x)}{v_x}\xi(x)),
\end{align*}
as claimed.
\end{subproof}
\end{proofsublemma}

\begin{proofsublemma}\label{psublem:TTPhiXv}
If\/ $\man{M}$ and\/ $\man{N}$ are smooth manifolds, if\/
$\Phi\in\mappings[\infty]{\man{M}}{\man{N}}$\@, and if\/
$X\in\sections[\infty]{\tb{\man{M}}}$\@, then
\begin{equation*}
\tf{\tf{\Phi}}\scirc\vl{X}(v_x)=
\verlift(\tf[x]{\Phi}(v_x),\tf[x]{\Phi}(X(x))).
\end{equation*}
\begin{subproof}
We have
\begin{align*}
\tf{\tf{\Phi}}\scirc\vl{X}(v_x)=&\;\derivatzero{}{t}\tf[x]{\Phi}(v_x+tX(x))\\
=&\;\derivatzero{}{t}(\tf[x]{\Phi}(v_x)+t\tf[x]{\Phi}(X(x)))\\
=&\;\verlift(\tf[x]{\Phi}(v_x),\tf[x]{\Phi}(X(x))),
\end{align*}
as claimed.
\end{subproof}
\end{proofsublemma}

We now directly compute, using Sublemma~\ref{psublem:T(fxi)}\@,
\begin{align*}
\what{\Phi}(\nabla^{\man{M}}_XfY)(x)=&\;
\tf{\Phi}\scirc K^{\man{M}}\scirc\tf{(fY)}\scirc X(x)\\
=&\;f(x)\tf{\Phi}\scirc K^{\man{M}}\scirc\tf{Y}\scirc X(x)+
\natpair{\d{f}(x)}{X(x)}\tf{\Phi}\scirc K^{\man{M}}\scirc\vl{Y}\scirc X(x)\\
=&\;f(x)\what{\Phi}(\nabla^{\man{M}}_XY)(x)+
\natpair{\d{f}(x)}{X(x)}\tf{\Phi}\scirc X(x),
\end{align*}
noting that $K^{\man{M}}$ is a left-inverse for vertical lift.  We also
directly compute, using both of the sublemmata above,
\begin{align*}
\Phi^*\nabla^{\man{N}}_X\what{fY}(x)=&\;
K^{\man{N}}\scirc\tf{\tf{\Phi}}\scirc\tf{(fY)}\scirc X(x)\\
=&\;f(x)K^{\man{N}}\scirc\tf{(\tf{\Phi}\scirc Y)}\scirc X(x)+
\natpair{\d{f}(x)}{X(x)}K^{\man{N}}\scirc\tf{\tf{\Phi}}\scirc\vl{Y}
\scirc X(x)\\
=&\;f(x)\Phi^*\nabla^{\man{N}}_X\what{Y}(x)+
\natpair{\d{f}(x)}{X(x)}K^{\man{N}}(\verlift(\tf[x]{\Phi}(X(x)),
\tf[x]{\Phi}(X(x))))\\
=&\;f(x)\Phi^*\nabla^{\man{N}}_X\what{Y}(x)+
\natpair{\d{f}(x)}{X(x)}\tf{\Phi}\scirc Y\scirc X(x),
\end{align*}
again noting that $K^{\man{N}}$ is the left-inverse for the vertical lift.
Combining the preceding two computations gives the tensoriality of
\begin{equation*}
(X,Y)\mapsto\what{\Phi}(\nabla^{\man{M}}_XY)(x)-
\Phi^*\nabla^{\man{N}}_X\what{\Phi}(Y)(x),
\end{equation*}
and so gives
$A_\Phi\in\sections[r]{\tensor*[2]{\tb{\man{M}}}\otimes \Phi^*\tb{\man{N}}}$
satisfying the assertion of the lemma.
\end{proof}
\end{lemma}

Note that, if $A\in\sections[\infty]{\tensor*[k]{\ctb{\man{N}}}}$\@, then
$\Phi^*A$ denotes the pull-back of $A$ to
$\sections[\infty]{\tensor*[k]{\ctb{\man{M}}}}$ and also the section of the
tensor bundle $\tensor*[k]{\Phi^*\ctb{\man{N}}}$\@.  Let
$x\in\tb[x]{\man{M}}$\@, let $v_1,\dots,v_k\in\tb[x]{\man{M}}$\@, and denote
$u_j=\tb[x]{\Phi}(v_j)$\@, $j\in\{1,\dots,k\}$\@.  Note that
\begin{equation}\label{eq:Phi*Amultiple}
\begin{aligned}
\Phi^*A((x,u_1),\dots,(x,u_k))=&\;A(u_1,\dots,u_k)=
A(\tf[x]{\Phi}(v_1),\dots,\tf[x]{\Phi}(v_k))\\
=&\;\Phi^*A(v_1,\dots,v_k),
\end{aligned}
\end{equation}
where we are using the two interpretations of the symbol $\Phi^*A$\@.

With the above as background, we can now understand the iterated covariant derivatives
\begin{equation*}
\nabla^{\man{M},m}\Phi^*f=
\underbrace{\nabla^{\man{M}}\cdots\nabla^{\man{M}}}_{m\ \textrm{times}}
\Phi^*f,\qquad m\in\integerp,
\end{equation*}
and
\begin{equation*}
\nabla^{\man{N},m}f=
\underbrace{\nabla^{\man{N}}\cdots\nabla^{\man{N}}}_{m\
\textrm{times}}f,\qquad m\in\integerp,
\end{equation*}
for $f\in\func[\infty]{\man{N}}$\@.  The following lemma gives the first part
of this development, playing the r\^ole of Lemma~\ref{lem:bundleform} in this
case.
\begin{lemma}\label{lem:pullbackform}
Let\/ $r\in\{\infty,\omega\}$\@.  Let\/ $\man{M}$ and\/ $\man{N}$ be\/
$\C^r$-manifolds with\/ $\C^r$-affine connections\/ $\nabla^{\man{M}}$ and\/
$\nabla^{\man{N}}$\@, respectively.  Define\/ $B_\Phi=\push_{1,2}A_\Phi$
with\/ $A_\Phi$ as in Lemma~\ref{lem:APhidef}\@.  Then, for\/ $k\in\integerp$
and\/ $A\in\sections[r]{\tensor*[k]{\ctb{\man{N}}}}$\@,
\begin{equation*}
\nabla^{\man{M}}\Phi^*A=\Phi^*\nabla^{\man{N}}A+D_{B_\Phi}(\Phi^*A).
\end{equation*}
\begin{proof}
Let $x\in\man{M}$\@.  Let
$X_1,\dots,X_k\in\sections[\infty]{\tb{\man{M}}}$\@.  For
$X_{k+1}\in\sections[\infty]{\tb{\man{M}}}$\@, we have
\begin{multline*}
\lieder{X_{k+1}}(\Phi^*A(X_1,\dots,X_k))\\
=(\nabla^{\man{M}}_{X_{k+1}}\Phi^*A)(X_1,\dots,X_k)+
\sum_{j=1}^k\Phi^*A(X_1,\dots,\nabla^{\man{M}}_{X_{k+1}}X_j,\dots,X_k)
\end{multline*}
and
\begin{multline*}
\lieder{X_{k+1}}(\Phi^*A(\what{\Phi}(X_1),\dots,\what{\Phi}(X_k)))=
(\Phi^*\nabla^{\man{N}}_{X_{k+1}}\Phi^*A)
(\what{\Phi}(X_1),\dots,\what{\Phi}(X_k))\\
+\sum_{j=1}^k\Phi^*A(\what{\Phi}(X_1),\dots,
\Phi^*\nabla^{\man{N}}_{X_{k+1}}\what{\Phi}(X_j),\dots,\what{\Phi}(X_k)),
\end{multline*}
using the two interpretations of $\Phi^*A$\@.  By~\eqref{eq:Phi*Amultiple} we
have, in the above expressions,
\begin{equation*}
\Phi^*A(X_1,\dots,X_k)=\Phi^*A(\what{\Phi}(X_1),\dots,\what{\Phi}(X_k)).
\end{equation*}
By~\eqref{eq:Phi*Amultiple} again, we have
\begin{equation*}
\Phi^*A(X_1,\dots,\nabla^{\man{M}}_{X_{k+1}}X_j,\dots,X_k)=
\Phi^*A(\what{\Phi}(X_1),\dots,\what{\Phi}(\nabla^{\man{M}}_{X_{k+1}}X_j),
\dots,\what{\Phi}(X_k)).
\end{equation*}
Also note that
\begin{align*}
(\Phi^*\nabla^{\man{N}}_{X_{k+1}}\Phi^*A)(\what{\Phi}(X_1),\dots,
\what{\Phi}(X_k))(x)=&\;
(\Phi^*\nabla^{\man{N}}_{X_{k+1}}\Phi^*A)(\tf[x]{\Phi}(X_1(x)),\dots,
\tf[x]{\Phi}(X_k(x)))\\
=&\;\nabla^{\man{N}}_{\tf[x]{\Phi}(X_{k+1}(x))}A(\tf[x]{\Phi}(X_1(x)),\dots,
\tf[x]{\Phi}(X_{k+1}(x)))\\
=&\;\nabla^{\man{N}}A(\tf[x]{\Phi}(X_1(x)),\dots,\tf[x]{\Phi}(X_{k+1}(x)))\\
=&\;\Phi^*\nabla^{\man{N}}A(X_1,\dots,X_{k+1})(x).
\end{align*}
Combining the above gives
\begin{align*}
\nabla^{\man{M}}\Phi^*A&(X_1,\dots,X_{k+1})\\
=&\;\Phi^*\nabla^{\man{N}}A(X_1,\dots,X_{k+1})\\
&\;+\sum_{j=1}^k\Phi^*A(\what{\Phi}(X_1),\dots,\Phi^*\nabla^{\man{N}}_{X_{k+1}}
\what{\Phi}(X_j)-\what{\Phi}(\nabla^{\man{M}}_{X_{k+1}}X_j),\dots,
\what{\Phi}(X_k))\\
=&\;\Phi^*\nabla^{\man{N}}A(X_1,\dots,X_{k+1})-
\sum_{j=1}^k\Phi^*A(\what{\Phi}(X_1),\dots,A_\Phi(X_{k+1},X_j),
\dots,\what{\Phi}(X_k)).
\end{align*}
Thus
\begin{equation*}
\nabla^{\man{M}}\Phi^*A=\Phi^*\nabla^{\man{N}}A-
\sum_{j=1}^k\Ins_j(\Phi^*A,B_\Phi),
\end{equation*}
giving the result by Lemma~\ref{lem:derinsertion2}\@.
\end{proof}
\end{lemma}

We now have the following lemma, the first of two regarding iterated
covariant differentials.
\begin{lemma}\label{lem:pbiterateddersI}
Let\/ $r\in\{\infty,\omega\}$ and let\/ $\man{M}$ and\/ $\man{N}$ be\/
$\C^r$-manifolds with\/ $\C^r$-affine connections\/ $\nabla^{\man{M}}$ and\/
$\nabla^{\man{N}}$\@, respectively.  For\/ $m\in\integernn$\@, there exist\/
$\C^r$-vector bundle mappings
\begin{equation*}
(A^m_s,\id_{\man{M}})\in\vbmappings[r]{\tensor*[s]{\Phi^*\ctb{\man{N}}}}
{\tensor*[m]{\ctb{\man{M}}}},\qquad s\in\{0,1,\dots,m\},
\end{equation*}
such that
\begin{equation*}
\nabla^{\man{M},m}\Phi^*f=\sum_{s=0}^mA^m_s(\Phi^*\nabla^{\man{N},s}f)
\end{equation*}
for all\/ $f\in\func[m]{\man{N}}$\@.  Moreover, the vector bundle mappings\/
$A^m_0,A^m_1,\dots,A^m_m$ satisfy the recursion relations prescribed by
\begin{equation*}
A^0_0(\beta_0)=\beta_0,\enspace A^1_1(\beta_1)=\beta_1,\enspace A^1_0=0,
\end{equation*}
and
{\small\begin{align*}
A^{m+1}_{m+1}(\beta_{m+1})=&\;\beta_{m+1},\\
A^{m+1}_s(\beta_s)=&\;(\nabla^{\man{M}}A^m_s)(\beta_s)+
A^m_{s-1}\otimes\id_{\ctb{\man{M}}}(\beta_s)-
\sum_{j=1}^sA^m_s\otimes\id_{\ctb{\man{M}}}(\Ins_j(\beta_s,B_\Phi)),\
s\in\{1,\dots,m\},\\
A^{m+1}_0(\beta_0)=&\;(\nabla^{\man{M}}A^m_0)(\beta_0),
\end{align*}}%
where\/ $\beta_s\in\tensor*[s]{\Phi^*\ctb{\man{N}}}$\@,\/
$s\in\{0,1,\dots,m\}$\@.
\begin{proof}
The assertion clearly holds for the initial conditions of the recursion,
simply because
\begin{equation*}
\Phi^*f=\Phi^*f,\quad\d{(\Phi^*f)}=\Phi^*\d{f}+0f.
\end{equation*}
So suppose that it holds for $m\in\integerp$\@.  Thus
\begin{equation*}
\nabla^{\man{M},m}\Phi^*f=\sum_{s=0}^mA^m_s(\Phi^*\nabla^{\man{N},s}f),
\end{equation*}
where the vector bundle mappings $A^a_s$\@, $a\in\{0,1,\dots,m\}$\@,
$s\in\{0,1,\dots,a\}$\@, satisfy the stated recursion relations.  Then
\begin{align*}
\nabla^{\man{M},m+1}\Phi^*f=&\;
\sum_{s=0}^m(\nabla^{\man{M}}A^m_s)(\Phi^*\nabla^{\man{N},s}f)+
\sum_{s=0}^mA^m_s\otimes\id_{\ctb{\man{M}}}
(\nabla^{\man{M}}\Phi^*\nabla^{\man{N},s}f)\\
=&\;\sum_{s=0}^m(\nabla^{\man{M}}A^m_s)(\Phi^*\nabla^{\man{N},s}f)+
\sum_{s=0}^mA^m_s\otimes\id_{\ctb{\man{M}}}(\Phi^*\nabla^{\man{N},s+1}f)\\
&\;-\sum_{s=0}^m\sum_{j=1}^sA^m_s\otimes\id_{\ctb{\man{M}}}
\Ins_j(\Phi^*\nabla^{\man{N},s}f,B_\Phi)\\
=&\;\Phi^*\nabla^{\man{N},m+1}f+\sum_{s=1}^m\left(\vphantom{\sum_{j=1}^s}
(\nabla^{\man{M}}A^m_s)(\Phi^*\nabla^{\man{N},s}f)+
A^m_{s-1}\otimes\id_{\ctb{\man{M}}}(\Phi^*\nabla^{\man{N},s}f)\right.\\
&\;\left.-\sum_{j=1}^sA^m_s\otimes\id_{\ctb{\man{M}}}
(\Ins_j(\Phi^*\nabla^{\man{N},s}f,B_\Phi))\right)+
(\nabla^{\man{M}}A^m_0)(\Phi^*f)
\end{align*}
by Lemma~\ref{lem:pullbackform}\@.  From this the lemma follows.
\end{proof}
\end{lemma}

We shall also need to ``invert'' the relationship of the preceding lemma.
\begin{lemma}\label{lem:pbiterateddersII}
Let\/ $r\in\{\infty,\omega\}$ and let\/ $\man{M}$ and\/ $\man{N}$ be\/
$\C^r$-manifolds with\/ $\C^r$-affine connections\/ $\nabla^{\man{M}}$ and\/
$\nabla^{\man{N}}$\@, respectively.  For\/ $m\in\integernn$\@, there exist\/
$\C^r$-vector bundle mappings
\begin{equation*}
(B^m_s,\id_{\man{M}})\in
\vbmappings[r]{\tensor*[s]{\ctb{\man{M}}}}
{\tensor*[m]{\Phi^*\ctb{\man{N}}}},\qquad s\in\{0,1,\dots,m\},
\end{equation*}
such that
\begin{equation*}
\Phi^*\nabla^{\man{N},m}f=\sum_{s=0}^mB^m_s(\nabla^{\man{M},s}\Phi^*f)
\end{equation*}
for all\/ $f\in\func[m]{\man{N}}$\@.  Moreover, the vector bundle mappings\/
$B^m_0,B^m_1,\dots,B^m_m$ satisfy the recursion relations prescribed by
\begin{equation*}
B^0_0(\alpha_0)=\alpha_0,\enspace B^1_1(\alpha_1)=\alpha_1,\enspace B^1_0=0,
\end{equation*}
and
{\small\begin{align*}
B^{m+1}_{m+1}(\alpha_{m+1})=&\;\alpha_{m+1},\\
B^{m+1}_s(\alpha_s)=&\;(\nabla^{\man{M}}B^m_s)(\alpha_s)+
B^m_{s-1}\otimes\id_{\ctb{\man{M}}}(\alpha_s)+
\sum_{j=1}^m\Ins_j(B^m_s(\alpha_s),B_\Phi),
\enspace s\in\{1,\dots,m\},\\
B^{m+1}_0(\alpha_0)=&\;(\nabla^{\man{M}}B^m_0)(\alpha_0)+
\sum_{j=1}^m\Ins_j(B^m_0(\alpha_0),B_\Phi),
\end{align*}}%
where\/ $\alpha_s\in\tensor*[s]{\ctb{\man{M}}}$\@,\/
$s\in\{0,1,\dots,m\}$\@.
\begin{proof}
The assertion clearly holds for the initial conditions for the recursion
because
\begin{equation*}
\Phi^*f=\Phi^*f,\quad\Phi^*(\d{f})=\d{(\Phi^*f)}+0f.
\end{equation*}
So suppose it true for $m\in\integerp$\@.  Thus
\begin{equation}\label{eq:pbiteratedderII1}
\Phi^*\nabla^{\man{N},m}f=\sum_{s=0}^mB^m_s(\nabla^{\man{M},s}\Phi^*f),
\end{equation}
where the vector bundle mappings $B^a_s$\@, $a\in\{0,1,\dots,m\}$\@,
$s\in\{0,1,\dots,a\}$\@, satisfy the recursion relations from the statement
of the lemma.  Then, by Lemma~\ref{lem:pullbackform}\@, we can work on the
left-hand side of~\eqref{eq:pbiteratedderII1} to give
\begin{align*}
\nabla^{\man{M}}\Phi^*\nabla^{\man{N},m}f=&\;
\Phi^*\nabla^{\man{N},m+1}f-\sum_{j=1}^m
\Ins_j(\Phi^*\nabla^{\man{N},m}f,B_\Phi)\\
=&\;\Phi^*\nabla^{\man{N},m+1}f-\sum_{s=0}^m\sum_{j=1}^m
\Ins_j(B^m_s(\nabla^{\man{M},s}\Phi^*f),B_\Phi).
\end{align*}
Working on the right-hand side of~\eqref{eq:PiteratedderII1} gives
\begin{equation*}
\nabla^{\man{M}}\Phi^*\nabla^{\man{N},m}f=
\sum_{s=0}^m\nabla^{\man{M}}B^m_s(\nabla^{\man{M},s}\Phi^*f)+
\sum_{s=0}^mB^m_s\otimes\id_{\ctb{\man{M}}}(\nabla^{\man{M},s+1}\Phi^*f).
\end{equation*}
Combining the preceding two equations gives
\begin{align*}
\Phi^*\nabla^{\man{N},m+1}f=&\;
\sum_{s=0}^m\nabla^{\man{M}}B^m_s(\nabla^{\man{M},s}\Phi^*f)+
\sum_{s=0}^mB^m_s\otimes\id_{\ctb{\man{M}}}(\nabla^{\man{M},s+1}\Phi^*f)\\
&\;+\sum_{s=0}^m\sum_{j=1}^m\Ins_j(B^m_s(\nabla^{\man{M},s}\Phi^*f),B_\Phi)\\
=&\;\nabla^{\man{M},m+1}\Phi^*f+
\left.\sum_{s=1}^m\right(\nabla^{\man{M}}B^m_s(\nabla^{\man{M},s}\Phi^*f)+
B^m_{s-1}\otimes\id_{\ctb{\man{M}}}(\nabla^{\man{M},s}\Phi^*f)\\
&\;\left.+\sum_{j=1}^m\Ins_j(B^m_s(\nabla^{\man{M},s}\Phi^*f),B_\Phi)\right)+
\nabla^{\man{M}}B^m_0(\Phi^*f)+\sum_{j=1}^m\Ins_j(B^m_0(\Phi^*f),B_\Phi),
\end{align*}
and the lemma follows from this.
\end{proof}
\end{lemma}

With this data, we have the following result.
\begin{lemma}\label{lem:pbsymiterateddersI}
Let\/ $r\in\{\infty,\omega\}$ and let\/ $\man{M}$ and\/ $\man{N}$ be\/
$\C^r$-manifolds with\/ $\C^r$-affine connections\/ $\nabla^{\man{M}}$ and\/
$\nabla^{\man{N}}$\@, respectively.  For\/ $m\in\integernn$\@, there exist\/
$\C^r$-vector bundle mappings
\begin{equation*}
(\what{A}^m_s,\id_{\man{M}})\in
\vbmappings[r]{\Symalg*[s]{\Phi^*\ctb{\man{N}}}}
{\Symalg*[m]{\ctb{\man{M}}}},\qquad s\in\{0,1,\dots,m\},
\end{equation*}
such that
\begin{equation*}
\Sym_m\scirc\nabla^{\man{M},m}\Phi^*f=
\sum_{s=0}^m\what{A}^m_s(\Sym_s\scirc\Phi^*\nabla^{\man{N},s}f)
\end{equation*}
for all\/ $f\in\func[m]{\man{N}}$\@.
\begin{proof}
This follows from Lemma~\ref{lem:pbiterateddersI} in the same way as
Lemma~\ref{lem:PsymiterateddersI} follows from
Lemma~\ref{lem:PiterateddersI}\@.
\end{proof}
\end{lemma}

Next we consider the ``inverse'' of the preceding lemma.
\begin{lemma}\label{lem:pbsymiterateddersII}
Let\/ $r\in\{\infty,\omega\}$ and let\/ $\man{M}$ and\/ $\man{N}$ be\/
$\C^r$-manifolds with\/ $\C^r$-affine connections\/ $\nabla^{\man{M}}$ and\/
$\nabla^{\man{N}}$\@, respectively.  For\/ $m\in\integernn$\@, there exist\/
$\C^r$-vector bundle mappings
\begin{equation*}
(\what{B}^m_s,\id_{\man{M}})\in
\vbmappings[r]{\Symalg*[s]{\ctb{\man{M}}}}
{\Symalg*[m]{\Phi^*\ctb{\man{N}}}},\qquad s\in\{0,1,\dots,m\},
\end{equation*}
such that
\begin{equation*}
\Sym_m\scirc\Phi^*\nabla^{\man{N},m}f=
\sum_{s=0}^m\what{B}^m_s(\Sym_s\scirc\nabla^{\man{M},s}\Phi^*f)
\end{equation*}
for all\/ $f\in\func[m]{\man{N}}$\@.
\begin{proof}
This follows from Lemma~\ref{lem:pbiterateddersII} in the same way as
Lemma~\ref{lem:PsymiterateddersII} follows from
Lemma~\ref{lem:PiterateddersII}\@.
\end{proof}
\end{lemma}

The following lemma provides two decompositions of
$\jetalg[\Phi]{m}{\man{M}}$\@, one ``in the domain'' and one ``in the
codomain,'' and the relationship between them.  The assertion simply results
from an examination of the preceding four lemmata.
\begin{lemma}\label{lem:pbdecomp}
Let\/ $r\in\{\infty,\omega\}$ and let\/ $\man{M}$ and\/ $\man{N}$ be\/
$\C^r$-manifolds with\/ $\C^r$-affine connections\/ $\nabla^{\man{M}}$ and\/
$\nabla^{\man{N}}$\@, respectively.  Then there exist\/ $\C^r$-vector bundle
mappings
\begin{equation*}
A^m_{\nabla^{\man{M}},\nabla^{\man{N}}}\in\vbmappings[r]{\jetalg[\Phi]{m}{\man{M}}}
{\Symalg*[\le m]{\Phi^*\ctb{\man{N}}}},\quad
B^m_{\nabla^{\man{M}},\nabla^{\man{N}}}\in\vbmappings[r]{\jetalg[\Phi]{m}{\man{M}}}
{\Symalg*[\le m]{\ctb{\man{M}}}},
\end{equation*}
defined by
\begin{align*}
A^m_{\nabla^{\man{M}},\nabla^{\man{N}}}(j_m(\Phi^*f)(x))=&\;
\Sym_{\le m}(\Phi^*f(x),\Phi^*\nabla^{\man{N}}f(x),\dots,
\Phi^*\nabla^{\man{N},m}f(x)),\\
B^m_{\nabla^{\man{M}},\nabla^{\man{N}}}(j_m(\Phi^*f)(x))=&\;
\Sym_{\le m}(\Phi^*f(x),\nabla^{\man{M}}\Phi^*f(x),\dots,
\nabla^{\man{M},m}\Phi^*f(x)).
\end{align*}
Moreover,\/ $A^m_{\nabla^{\man{M}},\nabla^{\man{N}}}$ is an isomorphism,\/
$B^m_{\nabla^{\man{M}},\nabla^{\man{N}}}$ is injective, and
\begin{multline*}
B^m_{\nabla^{\man{M}},\nabla^{\man{N}}}\scirc
(A^m_{\nabla^{\man{M}},\nabla^{\man{N}}})^{-1}\scirc
(\Sym_{\le m}(\Phi^*f(e),\Phi^*\nabla^{\man{N}}f(x),\dots,
\Phi^*\nabla^{\man{N},m}f(x))\\
=\left(A^0_0(\Phi^*f(x)),\sum_{s=0}^1\what{A}^1_s
(\Sym_s\scirc\Phi^*\nabla^{\man{N},s}f(x)),\dots,
\sum_{s=0}^m\what{A}^m_s(\Sym_s\scirc\Phi^*\nabla^{\man{N},s}f(x))\right)
\end{multline*}
and
\begin{multline*}
A^m_{\nabla^{\man{M}},\nabla^{\man{N}}}\scirc
(B^m_{\nabla^{\man{M}},\nabla^{\man{N}}})^{-1}\scirc
\Sym_{\le m}(\Phi^*f(x),\nabla^{\man{M}}\Phi^*f(x),\dots,
\nabla^{\man{M},m}\Phi^*f(x))\\
=\left(B^0_0(\Phi^*f(x)),
\sum_{s=0}^1\what{B}^1_s(\Sym_s\scirc\nabla^{\man{M},s}\Phi^*f(x)),\dots,
\sum_{s=0}^m\what{B}^m_s(\Sym_s\scirc\nabla^{\man{M},s}\Phi^*f(x))\right),
\end{multline*}
where the vector bundle mappings\/ $\what{A}^m_s$ and\/ $\what{B}^m_s$\@,\/
$s\in\{0,1,\dots,m\}$\@, are as in Lemmata~\ref{lem:pbsymiterateddersI}
and~\ref{lem:pbsymiterateddersII}\@.
\end{lemma}

\section{Fibre norms for some useful jet bundles}\label{sec:fibre-norms}

In \ref{sec:lift-isomorphisms} we saw how to make decompositions for jets of
sections of vector bundles and jets of various lifts to the total space of a
vector bundle $\map{\pi_{\man{E}}}{\man{E}}{\man{M}}$\@, using the
Levi-Civita affine connection induced by a natural Riemannian metric on
$\man{E}$\@.  In this section we consider fibre norms for these jet bundles.
The fibre norm for the space of jets of sections of a vector bundle is
deduced in a natural way from a Riemannian metric on $\man{M}$ and a fibre
metric in $\map{\pi_{\man{E}}}{\man{E}}{\man{M}}$\@.  For fibre norms of
lifted objects, the story is more complicated.  Since the objects are lifted
from $\man{M}$\@, there are two natural fibre norms in each case, one coming
from the Riemannian metric on $\man{E}$\@, and the other coming from the
Riemannian metric on $\man{M}$ and the fibre metric on the vector bundle.

The setup is the following.  We let $r\in\{\infty,\omega\}$ and let
$\map{\pi_{\man{E}}}{\man{E}}{\man{M}}$ be a $\C^r$-vector bundle.  We
consider a Riemannian metric $\metric_{\man{M}}$ on $\man{M}$\@, a fibre
metric $\metric_{\pi_{\man{E}}}$ on $\man{E}$\@, the Levi-Civita connection
$\nabla^{\man{M}}$ on $\man{M}$\@, and a vector bundle connection
$\nabla^{\pi_{\man{E}}}$ in $\man{E}$\@, all being of class $\C^r$\@.  This
gives the Riemannian metric $\metric_{\man{E}}$ of~\eqref{eq:metricE} and the
associated Levi-Civita connection $\nabla^{\man{E}}$\@.  This data gives the
fibre metrics for all sorts of tensors defined on the total space
$\man{E}$\@.  We, however, are interested only in the lifted tensors such as
are described in Section~\ref{sec:tensor-constructions}\@.

The reader will definitely observe a certain repetitiveness to our
constructions in this section, rather similar to that seen in
Section~\ref{sec:lift-isomorphisms}\@.  However, the ideas here are important
and the notation is confusing, so we do not skip anything.

We treat the smooth and real analytic cases simultaneously in this section.
In the smooth case, the formulae we give are useful for applying the methods
of the paper to the setting of the paper in smooth category.

\subsection{Fibre norms for horizontal lifts of
functions}\label{subsec:Pjetnorm}

We let $r\in\{\infty,\omega\}$ and let
$\map{\pi_{\man{E}}}{\man{E}}{\man{M}}$ be a $\C^r$-vector bundle.  For
$f\in\func[m]{\man{M}}$\@, we have $\pi_{\man{E}}^*f\in\func[m]{\man{E}}$\@.
We can, therefore, think of the $m$-jet of $\pi_{\man{E}}^*f$ as being
characterised by $j_mf$\@, as well as by $j_m\pi_{\man{E}}^*f$\@, and of
comparing these two characterisations.  Thus we have the two fibre norms
\begin{equation*}
\dnorm{j_mf(x)}^2_{\metric_{\man{M},m}}=
\sum_{j=0}^m\frac{1}{(j!)^2}\dnorm{\nabla^{\man{M},j}f(x)}^2_{\metric_{\man{M}}}
\end{equation*}
and
\begin{equation}\label{eq:pi*fBnorm}
\dnorm{j_m\pi_{\man{E}}^*f(e)}^2_{\metric_{\man{E},m}}=
\sum_{j=0}^m\frac{1}{(j!)^2}
\dnorm{\nabla^{\man{E},j}\pi_{\man{E}}^*f(e)}^2_{\metric_{\man{E}}}.
\end{equation}
These fibre norms can be related by virtue of Lemma~\ref{lem:Pdecomp}\@.  To
do so, we make use of the following lemma.
\begin{lemma}\label{lem:Pnormbase}
\begin{equation*}
\dnorm{\pi_{\man{E}}^*\nabla^{\man{M},m}f(e)}_{\metric_{\man{E}}}=
\dnorm{\nabla^{\man{M},m}f(\pi_{\man{E}}(e))}_{\metric_{\man{M}}}.
\end{equation*}
\begin{proof}
We have the fibre metric $\metric_{\man{E}}^{-1}$ on $\ctb{\man{E}}$
associated with the Riemannian metric $\metric_{\man{E}}$\@.  The subbundles
$\chb{\man{E}}$ and $\cvb{\man{E}}$ are $\metric_{\man{E}}^{-1}$-orthogonal.
We note that
$\map{\ctf[e]{\pi_{\man{E}}}}
{\ctb[\pi_{\man{E}}(e)]{\man{M}}}{\chb[e]{\man{E}}}$ is an isometry.  Thus we
have the formula
\begin{equation*}
\dnorm{\pi_{\man{E}}^*B}_{\metric_{\man{E}}}=\dnorm{B}_{\metric_{\man{M}}},\qquad
B\in\sections[0]{\tensor*[m]{\ctb{\man{M}}}},
\end{equation*}
and the assertion of the lemma is merely a special case of this formula.
\end{proof}
\end{lemma}

We note that the fibre norm~\eqref{eq:pi*fBnorm} makes use of the vector
bundle mapping
\begin{equation*}
B^m_{\nabla^{\man{E}}}\in\vbmappings[r]{\Pjetalg{m}{\man{E}}}
{\Symalg*[\le m]{\ctb{\man{E}}}}
\end{equation*}
from Lemma~\ref{lem:Pdecomp}\@.  If instead we use the vector bundle mapping
\begin{equation*}
A^m_{\nabla^{\man{E}}}\in\vbmappings[r]{\Pjetalg{m}{\man{E}}}
{\Symalg*[\le m]{\pi_{\man{E}}^*\ctb{\man{M}}}}
\end{equation*}
from Lemma~\ref{lem:Pdecomp}\@, then we have the alternative fibre norm
\begin{equation*}
\dnorm{j_m\pi_{\man{E}}^*f(e)}^{\prime2}_{\metric_{\man{E},m}}=
\sum_{j=0}^m\frac{1}{(j!)^2}
\dnorm{\pi_{\man{E}}^*\nabla^{\man{M},j}f(e)}^2_{\metric_{\man{E}}}=
\sum_{j=0}^m\frac{1}{(j!)^2}
\dnorm{\nabla^{\man{M},j}f(\pi_{\man{E}}(e))}^2_{\metric_{\man{M}}}.
\end{equation*}
The relationship between the fibre norms
$\dnorm{\cdot}_{\metric_{\man{E},m}}$ and
$\dnorm{\cdot}'_{\metric_{\man{E},m}}$ can be phrased as, ``What is the
relationship between the jet of the lift and the lift of the jet?''  This is
a question we will phrase below for other sorts of lifts, and will address
comprehensively when we prove the continuity of the various lifting
operations in Section~\ref{subsec:liftcont}\@.

\subsection{Fibre norms for vertical lifts of sections}

We let $r\in\{\infty,\omega\}$ and let
$\map{\pi_{\man{E}}}{\man{E}}{\man{M}}$ be a $\C^r$-vector bundle.  For
$\xi\in\sections[m]{\man{E}}$\@, we have
$\vl{\xi}\in\sections[m]{\tb{\man{E}}}$\@.  We can, therefore, think of the
$m$-jet of $\vl{\xi}$ as being characterised by $j_m\xi$\@, as well as by
$j_m\vl{\xi}$\@, and of comparing these two characterisations.  Thus we have
the two fibre norms
\begin{equation*}
\dnorm{j_m\xi(x)}^2_{\metric_{\man{M},\pi_{\man{E}},m}}=
\sum_{j=0}^m\frac{1}{(j!)^2}
\dnorm{\nabla^{\man{M},\pi_{\man{E}},j}\xi(x)}^2_{\metric_{\man{M},\pi_{\man{E}}}}
\end{equation*}
and
\begin{equation}\label{eq:xivBnorm}
\dnorm{j_m\vl{\xi}(e)}^2_{\metric_{\man{E},m}}=\sum_{j=0}^m\frac{1}{(j!)^2}
\dnorm{\nabla^{\man{E},j}\vl{\xi}(e)}^2_{\metric_{\man{E}}}.
\end{equation}
These fibre norms can be related by virtue of Lemma~\ref{lem:Vdecomp}\@.  To
do so, we make use of the following lemma.
\begin{lemma}\label{lem:Vnormbase}
\begin{equation*}
\dnorm{\vl{(\nabla^{\man{M},\pi_{\man{E}},m}\xi)}(e)}_{\metric_{\man{E}}}=
\dnorm{\nabla^{\man{M},\pi_{\man{E}},m}
\xi(\pi_{\man{E}}(e))}_{\metric_{\man{M},\pi_{\man{E}}}}.
\end{equation*}
\begin{proof}
The subbundles $\hb{\man{E}}$ and $\vb{\man{E}}$ are
$\metric_{\man{E}}$-orthogonal and the subbundles $\chb{\man{E}}$ and
$\cvb{\man{E}}$ are $\metric_{\man{E}}^{-1}$-orthogonal.  We note that the
identification $\vb[e]{\man{E}}\simeq\man{E}_{\pi_{\man{E}}(e)}$ is an isometry and
that $\map{\ctf[e]{\pi_{\man{E}}}}{\ctb[\pi_{\man{E}}(e)]{\man{M}}}{
\chb[e]{\man{E}}}$ is an
isometry.  Thus we have the formula
\begin{equation*}
\dnorm{\vl{B}}_{\metric_{\man{E}}}=\dnorm{B}_{\metric_{\man{M},\pi_{\man{E}}}},
\qquad B\in\sections[0]{\tensor*[m]{\ctb{\man{M}}}\otimes\man{E}},
\end{equation*}
and the assertion of the lemma is merely a special case of this formula.
\end{proof}
\end{lemma}

We note that the fibre norm~\eqref{eq:xivBnorm} makes use of the vector
bundle mapping
\begin{equation*}
B^m_{\nabla^{\man{E}}}\in\vbmappings[r]{\Pjetalg{m}{\man{E}}\otimes
\vb{\man{E}}}{\Symalg*[\le m]{\ctb{\man{E}}}\otimes\vb{\man{E}}}
\end{equation*}
from Lemma~\ref{lem:Vdecomp}\@.  If instead we use the vector bundle mapping
\begin{equation*}
A^m_{\nabla^{\man{E}}}\in\vbmappings[r]{\Pjetalg{m}{\man{E}}\otimes
\vb{\man{E}}}{\Symalg*[\le m]{\pi_{\man{E}}^*\ctb{\man{M}}}\otimes\vb{\man{E}}}
\end{equation*}
from Lemma~\ref{lem:Vdecomp}\@, then we have the alternative fibre norm
\begin{equation*}
\dnorm{j_m\vl{\xi}(e)}^{\prime2}_{\metric_{\man{E},m}}=
\sum_{j=0}^m\frac{1}{(j!)^2}
\dnorm{\vl{(\nabla^{\man{M},\pi_{\man{E}},j}\xi)}(e)}^2_{\metric_{\man{E}}}=
\sum_{j=0}^m\frac{1}{(j!)^2}\dnorm{\nabla^{\man{M},\pi_{\man{E}},j}
\xi(\pi_{\man{E}}(e))}^2_{\metric_{\man{M},\pi_{\man{E}}}}.
\end{equation*}
Again, this points out the matter of the relationship between the jet of a
lift versus the lift of the jet, and this matter will be considered in detail
in the continuity results of Section~\ref{subsec:liftcont}\@.

\subsection{Fibre norms for horizontal lifts of vector fields}

We let $r\in\{\infty,\omega\}$ and let
$\map{\pi_{\man{E}}}{\man{E}}{\man{M}}$ be a $\C^r$-vector bundle.  For
$X\in\sections[m]{\tb{\man{M}}}$\@, we have
$\hl{X}\in\sections[m]{\tb{\man{E}}}$\@.  We can, therefore, think of the
$m$-jet of $\hl{X}$ as being characterised by $j_mX$\@, as well as by
$j_m\hl{X}$\@, and of comparing these two characterisations.  Thus we have
the two fibre norms
\begin{equation*}
\dnorm{j_mX(x)}^2_{\metric_{\man{M},m}}=
\sum_{j=0}^m\frac{1}{(j!)^2}\dnorm{\nabla^{\man{M},j}X(x)}^2_{\metric_{\man{M}}}
\end{equation*}
and
\begin{equation}\label{eq:XhBnorm}
\dnorm{j_m\hl{X}(e)}^2_{\metric_{\man{E},m}}=
\sum_{j=0}^m\frac{1}{(j!)^2}
\dnorm{\nabla^{\man{E},j}\hl{X}(e)}^2_{\metric_{\man{E}}}.
\end{equation}
These fibre norms can be related by virtue of Lemma~\ref{lem:Hdecomp}\@.  To
do so, we make use of the following lemma.
\begin{lemma}\label{lem:Hnormbase}
\begin{equation*}
\dnorm{\hl{(\nabla^{\man{M},m}X)}(e)}_{\metric_{\man{E}}}=
\dnorm{\nabla^{\man{M},m}X(\pi_{\man{E}}(e))}_{\metric_{\man{M}}}.
\end{equation*}
\begin{proof}
The subbundles $\hb{\man{E}}$ and $\vb{\man{E}}$ are
$\metric_{\man{E}}$-orthogonal.  We note that the identification
$\hb[e]{\man{E}}\simeq\tb[\pi_{\man{E}}(e)]{\man{M}}$ is an isometry and that
$\map{\ctf[e]{\pi_{\man{E}}}}{\ctb[\pi_{\man{E}}(e)]{\man{M}}}{\chb[e]{\man{E}}}$ is an isometry.
Thus we have the formula
\begin{equation*}
\dnorm{\hl{B}}_{\metric_{\man{E}}}=\dnorm{B}_{\metric_{\man{M}}},\qquad
B\in\sections[0]{\tensor*[m]{\ctb{\man{M}}}\otimes\tb{\man{M}}},
\end{equation*}
and the assertion of the lemma is merely a special case of this formula.
\end{proof}
\end{lemma}

We note that the fibre norm~\eqref{eq:XhBnorm} makes use of the vector
bundle mapping
\begin{equation*}
B^m_{\nabla^{\man{E}}}\in\vbmappings[r]{\Pjetalg{m}{\man{E}}\otimes
\hb{\man{E}}}{\Symalg*[\le m]{\ctb{\man{E}}}\otimes\hb{\man{E}}}
\end{equation*}
from Lemma~\ref{lem:Hdecomp}\@.  If instead we use the vector bundle mapping
\begin{equation*}
A^m_{\nabla^{\man{E}}}\in\vbmappings[r]{\Pjetalg{m}{\man{E}}
\otimes\hb{\man{E}}}{\Symalg*[\le m]{\pi_{\man{E}}^*\ctb{\man{M}}}
\otimes\hb{\man{E}}}
\end{equation*}
from Lemma~\ref{lem:Hdecomp}\@, then we have the alternative fibre norm
\begin{equation*}
\dnorm{j_m\hl{X}(e)}^{\prime2}_{\metric_{\man{E},m}}=
\sum_{j=0}^m\frac{1}{(j!)^2}
\dnorm{\hl{(\nabla^{\man{M},j}X)}(e)}^2_{\metric_{\man{E}}}=
\sum_{j=0}^m\frac{1}{(j!)^2}
\dnorm{\nabla^{\man{M},j}X(\pi_{\man{E}}(e))}^2_{\metric_{\man{M}}}.
\end{equation*}
Again, this points out the matter of the relationship between the jet of a
lift versus the lift of the jet, and this matter will be considered in detail
in the continuity results of Section~\ref{subsec:liftcont}\@.

\subsection{Fibre norms for vertical lifts of dual sections}

We let $r\in\{\infty,\omega\}$ and let
$\map{\pi_{\man{E}}}{\man{E}}{\man{M}}$ be a $\C^r$-vector bundle.  For
$\lambda\in\sections[m]{\dual{\man{E}}}$\@, we have
$\vl{\lambda}\in\sections[m]{\ctb{\man{E}}}$\@.  We can, therefore, think of
the $m$-jet of $\vl{\lambda}$ as being characterised by $j_m\lambda$\@, as
well as by $j_m\vl{\lambda}$\@, and of comparing these two characterisations.
Thus we have fibre norms
\begin{equation*}
\dnorm{j_m\lambda(x)}^2_{\metric_{\man{M},\pi_{\man{E}},m}}=
\sum_{j=0}^m\frac{1}{(j!)^2}\dnorm{\nabla^{\man{M},\pi_{\man{E}},j}
\lambda(x)}^2_{\metric_{\man{M},\pi_{\man{E}}}}
\end{equation*}
and
\begin{equation}\label{eq:lvBnorm}
\dnorm{j_m\vl{\lambda}(e)}^2_{\metric_{\man{E},m}}=\sum_{j=0}^m\frac{1}{(j!)^2}
\dnorm{\nabla^{\man{E},j}\vl{\lambda}(e)}^2_{\metric_{\man{E}}}.
\end{equation}
These fibre norms can be related by virtue of Lemma~\ref{lem:V*decomp}\@.  To
do so, we make use of the following lemma.
\begin{lemma}\label{lem:V*normbase}
\begin{equation*}
\dnorm{\vl{(\nabla^{\man{M},\pi_{\man{E}},m}\lambda)}(e)}_{\metric_{\man{E}}}=
\dnorm{\nabla^{\man{M},\pi_{\man{E}},m}
\lambda(\pi_{\man{E}}(e))}_{\metric_{\man{M},\pi_{\man{E}}}}.
\end{equation*}
\begin{proof}
The subbundles $\chb{\man{E}}$ and $\cvb{\man{E}}$ are
$\metric_{\man{E}}^{-1}$-orthogonal.  We note that the identification
$\cvb[e]{\man{E}}\simeq\dual{\man{E}}_{\pi_{\man{E}}(e)}$ is an isometry and
that $\map{\ctf[e]{\pi_{\man{E}}}}
{\ctb[\pi_{\man{E}}(e)]{\man{M}}}{\chb[e]{\man{E}}}$ is an isometry.  Thus we
have the formula
\begin{equation*}
\dnorm{\vl{B}}_{\metric_{\man{E}}}=\dnorm{B}_{\metric_{\man{M},\pi_{\man{E}}}},
\qquad B\in\sections[0]{\tensor*[m]{\ctb{\man{M}}}\otimes\dual{\man{E}}},
\end{equation*}
and the assertion of the lemma is merely a special case of this formula.
\end{proof}
\end{lemma}

We note that the fibre norm~\eqref{eq:lvBnorm} makes use of the vector
bundle mapping
\begin{equation*}
B^m_{\nabla^{\man{E}}}\in
\vbmappings[r]{\Pjetalg{m}{\man{E}}\otimes\cvb{\man{E}}}
{\Symalg*[\le m]{\ctb{\man{E}}}\otimes\cvb{\man{E}}}
\end{equation*}
from Lemma~\ref{lem:V*decomp}\@.  If instead we use the vector bundle mapping
\begin{equation*}
A^m_{\nabla^{\man{E}}}\in
\vbmappings[r]{\Pjetalg{m}{\man{E}}\otimes\cvb{\man{E}}}
{\Symalg*[\le m]{\pi_{\man{E}}^*\ctb{\man{M}}}\otimes\cvb{\man{E}}}
\end{equation*}
from Lemma~\ref{lem:V*decomp}\@, then we have the alternative fibre norm
\begin{equation*}
\dnorm{j_m\vl{\lambda}(e)}^{\prime2}_{\metric_{\man{E},m}}=
\sum_{j=0}^m\frac{1}{(j!)^2}
\dnorm{\vl{(\nabla^{\man{M},\pi_{\man{E}},j}\lambda)}(e)}^2_{\metric_{\man{E}}}=
\sum_{j=0}^m\frac{1}{(j!)^2}\dnorm{\nabla^{\man{M},\pi_{\man{E}},j}
\lambda(\pi_{\man{E}}(e))}^2_{\metric_{\man{M},\pi_{\man{E}}}}.
\end{equation*}
Again, this points out the matter of the relationship between the jet of a
lift versus the lift of the jet, and this matter will be considered in detail
in the continuity results of Section~\ref{subsec:liftcont}\@.

\subsection{Fibre norms for vertical lifts of endomorphisms}

We let $r\in\{\infty,\omega\}$ and let
$\map{\pi_{\man{E}}}{\man{E}}{\man{M}}$ be a $\C^r$-vector bundle.  For
$L\in\sections[m]{\tensor[1,1]{\man{E}}}$\@, we have
$\vl{L}\in\sections[m]{\tensor[1,1]{\man{E}}}$\@.  We can, therefore, think
of the $m$-jet of $\vl{L}$ as being characterised by $j_mL$\@, as well as by
$j_m\vl{L}$\@, and of comparing these two characterisations.  Thus we have
the two fibre norms
\begin{equation*}
\dnorm{j_mL(x)}^2_{\metric_{\man{M},\pi_{\man{E}},m}}=
\sum_{j=0}^m\frac{1}{(j!)^2}
\dnorm{\nabla^{\man{M},\pi_{\man{E}},j}L(x)}^2_{\metric_{\man{M},\pi_{\man{E}}}}
\end{equation*}
and
\begin{equation}\label{eq:LvBnorm}
\dnorm{j_m\vl{L}(e)}^2_{\metric_{\man{E},m}}=
\sum_{j=0}^m\frac{1}{(j!)^2}
\dnorm{\nabla^{\man{E},j}\vl{L}(e)}^2_{\metric_{\man{E}}}.
\end{equation}
These fibre norms can be related by virtue of Lemma~\ref{lem:Ldecomp}\@.  To
do so, we make use of the following lemma.
\begin{lemma}\label{lem:Lnormbase}
\begin{equation*}
\dnorm{\vl{(\nabla^{\man{M},\pi_{\man{E}},m}L)}(e)}_{\metric_{\man{E}}}=
\dnorm{\nabla^{\man{M},\pi_{\man{E}},m}
L(\pi_{\man{E}}(e))}_{\metric_{\man{M},\pi_{\man{E}}}}.
\end{equation*}
\begin{proof}
The subbundles $\chb{\man{E}}$ and $\cvb{\man{E}}$ are
$\metric_{\man{E}}^{-1}$-orthogonal.  We note that the identifications
$\vb[e]{\man{E}}\simeq\man{E}_{\pi_{\man{E}}(e)}$ and
$\cvb[e]{\man{E}}\simeq\dual{\man{E}}_{\pi_{\man{E}}(e)}$ are isometries, and
that $\map{\ctf[e]{\pi_{\man{E}}}}{\ctb[\pi_{\man{E}}(e)]{\man{M}}}
{\chb[e]{\man{E}}}$ is an isometry.  Thus we have the formula
\begin{equation*}
\dnorm{\vl{B}}_{\metric_{\man{E}}}=\dnorm{B}_{\metric_{\man{M},\pi_{\man{E}}}},
\qquad
B\in\sections[0]{\tensor*[m]{\ctb{\man{M}}}\otimes\tensor[1,1]{\man{E}}},
\end{equation*}
and the assertion of the lemma is merely a special case of this formula.
\end{proof}
\end{lemma}

We note that the fibre norm~\eqref{eq:LvBnorm} makes use of the vector
bundle mapping
\begin{equation*}
B^m_{\nabla^{\man{E}}}\in\vbmappings[r]{\Pjetalg{m}{\man{E}}\otimes
\tensor[1,1]{\vb{\man{E}}}}
{\Symalg*[\le m]{\ctb{\man{E}}}\otimes\tensor[1,1]{\vb{\man{E}}}}
\end{equation*}
from Lemma~\ref{lem:Ldecomp}\@.  If instead we use the vector bundle mapping
\begin{equation*}
A^m_{\nabla^{\man{E}}}\in\vbmappings[r]{\Pjetalg{m}{\man{E}}\otimes
\tensor[1,1]{\vb{\man{E}}}}
{\Symalg*[\le m]{\pi_{\man{E}}^*\ctb{\man{M}}}\otimes
\tensor[1,1]{\vb{\man{E}}}}
\end{equation*}
from Lemma~\ref{lem:Ldecomp}\@, then we have the alternative fibre norm
\begin{equation*}
\dnorm{j_m\vl{L}(e)}^{\prime2}_{\metric_{\man{E},m}}=
\sum_{j=0}^m\frac{1}{(j!)^2}
\dnorm{\vl{(\nabla^{\man{M},\pi_{\man{E}},j}L)}(e)}^2_{\metric_{\man{E}}}=
\sum_{j=0}^m\frac{1}{(j!)^2}\dnorm{\nabla^{\man{M},\pi_{\man{E}},j}
L(\pi_{\man{E}}(e))}^2_{\metric_{\man{M},\pi_{\man{E}}}}.
\end{equation*}
Again, this points out the matter of the relationship between the jet of a
lift versus the lift of the jet, and this matter will be considered in detail
in the continuity results of Section~\ref{subsec:liftcont}\@.

\subsection{Fibre norms for vertical evaluations of dual sections}

We let $r\in\{\infty,\omega\}$ and let
$\map{\pi_{\man{E}}}{\man{E}}{\man{M}}$ be a $\C^r$-vector bundle.  For
$\lambda\in\sections[m]{\dual{\man{E}}}$\@, we have
$\ve{\lambda}\in\func[m]{\man{E}}$\@.  We can, therefore, think of the
$m$-jet of $\ve{\lambda}$ as being characterised by $j_m\lambda$\@, as well
as by $j_m\ve{\lambda}$\@, and of comparing these two characterisations.

Thus we have the two fibre norms
\begin{equation*}
\dnorm{j_m\lambda(x)}^2_{\metric_{\man{M},\pi_{\man{E}},m}}=
\sum_{j=0}^m\frac{1}{(j!)^2}\dnorm{\nabla^{\man{M},\pi_{\man{E}},j}
\lambda(x)}^2_{\metric_{\man{M},\pi_{\man{E}}}}
\end{equation*}
and
\begin{equation}\label{eq:leBnorm}
\dnorm{j_m\ve{\lambda}(e)}^2_{\metric_{\man{E},m}}=\sum_{j=0}^m\frac{1}{(j!)^2}
\dnorm{\nabla^{\man{E},j}\ve{\lambda}(e)}^2_{\metric_{\man{E}}}.
\end{equation}
These fibre norms can be related by virtue of Lemma~\ref{lem:Ddecomp}\@.  To
do so, we make use of the following lemma.
\begin{lemma}\label{lem:Dnormbase}
\begin{equation*}
\dnorm{\ve{(\nabla^{\man{M},\pi_{\man{E}},m}\lambda)}(e)}_{\metric_{\man{E}}}=
\dnorm{\nabla^{\man{M},\pi_{\man{E}},m}
\lambda(\pi_{\man{E}}(e))(e)}_{\metric_{\man{M},\pi_{\man{E}}}}.
\end{equation*}
\begin{proof}
The subbundles $\chb{\man{E}}$ and $\cvb{\man{E}}$ are
$\metric_{\man{E}}^{-1}$-orthogonal.  We note that the identification
$\cvb[e]{\man{E}}\simeq\dual{\man{E}}_{\pi_{\man{E}}(e)}$ is an isometry, and
that $\map{\ctf[e]{\pi_{\man{E}}}}{\ctb[\pi_{\man{E}}(e)]{\man{M}}}
{\chb[e]{\man{E}}}$ is an isometry.  Thus we have the formula
\begin{equation*}
\dnorm{\ve{B}(e)}_{\metric_{\man{E}}}=
\dnorm{B(\pi_{\man{E}}(e))(e)}_{\metric_{\man{M},\pi_{\man{E}}}},\qquad
B\in\sections[0]{\tensor*[m]{\ctb{\man{M}}}\otimes\dual{\man{E}}},
\end{equation*}
and the assertion of the lemma is merely a special case of this formula.
\end{proof}
\end{lemma}

We note that the fibre norm~\eqref{eq:leBnorm} makes use of the vector
bundle mapping
\begin{equation*}
B^m_{\nabla^{\man{E}}}\in\vbmappings[r]{\Pjetalg{m}{\man{E}}}
{\Symalg*[\le m]{\ctb{\man{E}}}}
\end{equation*}
from Lemma~\ref{lem:Ddecomp}\@.  If instead we use the vector bundle mapping
\begin{equation*}
A^m_{\nabla^{\man{E}}}\in\vbmappings[r]{\Pjetalg{m}{\man{E}}}
{\Symalg*[\le m]{\pi_{\man{E}}^*\ctb{\man{M}}}}
\end{equation*}
from Lemma~\ref{lem:Ddecomp}\@, then we have the alternative fibre norm
\begin{equation*}
\dnorm{j_m\ve{\lambda}(e)}^{\prime2}_{\metric_{\man{E},m}}=
\sum_{j=0}^m\frac{1}{(j!)^2}
\dnorm{\ve{(\nabla^{\man{M},\pi_{\man{E}},j}\lambda)}(e)}^2_{\metric_{\man{E}}}=
\sum_{j=0}^m\frac{1}{(j!)^2}
\dnorm{\nabla^{\man{M},\pi_{\man{E}},j}
\lambda(\pi_{\man{E}}(e))(e)}^2_{\metric_{\man{M},\pi_{\man{E}}}}.
\end{equation*}
Again, this points out the matter of the relationship between the jet of a
lift versus the lift of the jet, and this matter will be considered in detail
in the continuity results of Section~\ref{subsec:liftcont}\@.

\subsection{Fibre norms for vertical evaluations of
endomorphisms}\label{subsec:Cjetnorm}

We let $r\in\{\infty,\omega\}$ and let
$\map{\pi_{\man{E}}}{\man{E}}{\man{M}}$ be a $\C^r$-vector bundle.  For
$L\in\sections[m]{\tensor[1,1]{\man{E}}}$\@, we have
$\ve{L}\in\sections[m]{\tb{\man{E}}}$\@.  We can, therefore, think
of the $m$-jet of $\ve{L}$ as being characterised by $j_mL$\@, as well as by
$j_m\ve{L}$\@, and of comparing these two characterisations.  Thus we have
the two fibre norms
\begin{equation*}
\dnorm{j_mL(x)}^2_{\metric_{\man{M},\pi_{\man{E}},m}}=
\sum_{j=0}^m\frac{1}{(j!)^2}
\dnorm{\nabla^{\man{M},\pi_{\man{E}},j}L(x)}^2_{\metric_{\man{M},\pi_{\man{E}}}}
\end{equation*}
and
\begin{equation}\label{eq:LeBnorm}
\dnorm{j_m\ve{L}(e)}^2_{\metric_{\man{E},m}}=
\sum_{j=0}^m\frac{1}{(j!)^2}
\dnorm{\nabla^{\man{E},j}\ve{L}(e)}^2_{\metric_{\man{E}}}.
\end{equation}
These fibre norms can be related by virtue of Lemma~\ref{lem:Cdecomp}\@.  To
do so, we make use of the following lemma.
\begin{lemma}\label{lem:Cnormbase}
\begin{equation*}
\dnorm{\ve{(\nabla^{\man{M},\pi_{\man{E}},m}L)}(e)}_{\metric_{\man{E}}}=
\dnorm{\nabla^{\man{M},\pi_{\man{E}},m}
L(\pi_{\man{E}}(e))(e)}_{\metric_{\man{M},\pi_{\man{E}}}}.
\end{equation*}
\begin{proof}
The subbundles $\chb{\man{E}}$ and $\cvb{\man{E}}$ are
$\metric_{\man{E}}^{-1}$-orthogonal.  We note that the identification
$\cvb[e]{\man{E}}\simeq\dual{\man{E}}_{\pi_{\man{E}}(e)}$ is an isometry and
that $\map{\ctf[e]{\pi_{\man{E}}}}{\ctb[\pi_{\man{E}}(e)]{\man{M}}}
{\chb[e]{\man{E}}}$ is an isometry.  Thus we have the formula
\begin{equation*}
\dnorm{\ve{B}(e)}_{\metric_{\man{E}}}=
\dnorm{B(\pi_{\man{E}}(e))(e)}_{\metric_{\man{M},\pi_{\man{E}}}},\qquad
B\in\sections[0]{\tensor*[m]{\ctb{\man{M}}}\otimes\tensor[1,1]{\man{E}}},
\end{equation*}
and the assertion of the lemma is merely a special case of this formula.
\end{proof}
\end{lemma}

We note that the fibre norm~\eqref{eq:LeBnorm} makes use of the vector
bundle mapping
\begin{equation*}
B^m_{\nabla^{\man{E}}}\in\vbmappings[r]{\Pjetalg{m}{\man{E}}\otimes
\vb{\man{E}}}{\Symalg*[\le m]{\ctb{\man{E}}}\otimes\vb{\man{E}}}
\end{equation*}
from Lemma~\ref{lem:Cdecomp}\@.  If instead we use the vector bundle mapping
\begin{equation*}
A^m_{\nabla^{\man{E}}}\in\vbmappings[r]{\Pjetalg{m}{\man{E}}\otimes
\vb{\man{E}}}{\Symalg*[\le m]{\pi_{\man{E}}^*\ctb{\man{M}}}\otimes\vb{\man{E}}}
\end{equation*}
from Lemma~\ref{lem:Cdecomp}\@, then we have the alternative fibre norm
\begin{equation*}
\dnorm{j_m\ve{L}(e)}^{\prime2}_{\metric_{\man{E},m}}=
\sum_{j=0}^m\frac{1}{(j!)^2}
\dnorm{\ve{(\nabla^{\man{M},\pi_{\man{E}},j}L)}(e)}^2_{\metric_{\man{E}}}=
\sum_{j=0}^m\frac{1}{(j!)^2}
\dnorm{\nabla^{\man{M},\pi_{\man{E}},j}L(e)}^2_{\metric_{\man{M},\pi_{\man{E}}}}.
\end{equation*}
Again, this points out the matter of the relationship between the jet of a
lift versus the lift of the jet, and this matter will be considered in detail
in the continuity results of Section~\ref{subsec:liftcont}\@.

\subsection{Fibre norms for pull-backs of functions}\label{subsec:pbjetnorm}

We let $r\in\{\infty,\omega\}$ and let $\man{M}$ and $\man{N}$ be
$\C^r$-manifolds, and let $\Phi\in\mappings[r]{\man{M}}{\man{N}}$\@.  For
$f\in\func[m]{\man{N}}$\@, we have $\Phi^*f\in\func[m]{\man{M}}$\@.  We can,
therefore, think of the $m$-jet of $\Phi^*f$ as being characterised by
$j_mf$\@, as well as by $j_m\Phi^*f$\@, and of comparing these two
characterisations.  Thus we have the two fibre norms
\begin{equation*}
\dnorm{j_mf(x)}^2_{\metric_{\man{N},m}}=
\sum_{j=0}^m\frac{1}{(j!)^2}\dnorm{\nabla^{\man{N},j}f(x)}^2_{\metric_{\man{N}}}
\end{equation*}
and
\begin{equation}\label{eq:Phi*fBnorm}
\dnorm{j_m\Phi^*f(e)}^2_{\metric_{\man{M},m}}=
\sum_{j=0}^m\frac{1}{(j!)^2}
\dnorm{\nabla^{\man{M},j}\Phi^*f(e)}^2_{\metric_{\man{M}}}.
\end{equation}
These fibre norms can be related by virtue of Lemma~\ref{lem:pbdecomp}\@.  To
make use of this relationship, we shall also need to relate the norms of the
terms in these expressions.  In the preceding sections, this was easy to do
since the Riemannian metric on $\man{E}$ was related in a specific way to the
Riemannian metric on $\man{M}$ and the fibre metric in $\man{E}$\@.  Here,
this is not so simple since, if we choose a Riemannian metric
$\metric_{\man{M}}$ on $\man{M}$ and a Riemannian metric $\metric_{\man{N}}$
on $\man{N}$\@, these will be have no useful relationship.  So, rather than
getting an equality between certain norms, the best we can achieve (and all
that we need) is a useful bound, and this is the content of the next lemma.
\begin{lemma}\label{lem:pbnormbase}
For a compact set\/ $\nbhd{K}\subset\man{M}$\@:
\begin{compactenum}[(i)]
\item \label{pl:pbnormbase1} there exists\/ $C\in\realp$ such that
\begin{equation*}
\dnorm{\Phi^*\nabla^{\man{N},m}f(x)}_{\metric_{\man{M}}}\le
C^m\dnorm{\nabla^{\man{N},m}f(\Phi(x))}_{\metric_{\man{N}}},\qquad
x\in\nbhd{K},\ m\in\integernn;
\end{equation*}
\item if\/ $\Phi$ is a submersion or an injective immersion, then\/ $C$ from
part~\eqref{pl:pbnormbase1} can be chosen so that it holds that
\begin{equation*}
\dnorm{\nabla^{\man{N},m}f(\Phi(x))}_{\metric_{\man{N}}}\le
C^m\dnorm{\Phi^*\nabla^{\man{N},m}f(x)}_{\metric_{\man{M}}},\qquad
x\in\nbhd{K},\ m\in\integernn.
\end{equation*}
\end{compactenum}
\begin{proof}
The essential part of the proof is the following linear algebraic sublemma.
\begin{proofsublemma}
Let\/ $(\alg{U},\metric_{\alg{U}})$ and\/ $(\alg{V},\metric_{\alg{V}})$ be
finite-dimensional\/ $\real$-inner product spaces and let\/
$\Phi\in\Hom_{\real}(\alg{U};\alg{V})$\@.  Then there exists\/ $C\in\realp$
such that
\begin{equation*}
\dnorm{\Phi^*A}_{\metric_{\alg{U}}}\le C^k\dnorm{A}_{\metric_{\alg{V}}}
\end{equation*}
for every\/ $A\in\tensor*[k]{\dual{\alg{V}}}$\@,\/ $k\in\integernn$\@.  If,
additionally,\/ $\Phi$ is a surjective submersion or an injective immersion,
then\/ $C$ can be chosen so that, additionally, it holds that
\begin{equation*}
\dnorm{A}_{\metric_{\alg{V}}}\le C^k\dnorm{\Phi^*A}_{\metric_{\alg{U}}}
\end{equation*}
for every\/ $A\in\tensor*[k]{\dual{\alg{V}}}$\@,\/ $k\in\integernn$\@.
\begin{subproof}
Let $\ifam{f_1,\dots,f_m}$ and $\ifam{e_1,\dots,e_n}$ be orthonormal bases
for $\alg{U}$ and $\alg{V}$ with dual bases $\ifam{f^1,\dots,f^m}$ and $\ifam{e^1,\dots,e^n}$\@.  Write
\begin{equation*}
A=\sum_{j_1,\dots,j_k=1}^nA_{j_1\cdots j_k}e^{j_1}\otimes\dots\otimes e^{j_k}
\end{equation*}
and
\begin{equation*}
\Phi=\sum_{j=1}^n\sum_{a=1}^m\Phi^j_ae_j\otimes f^a.
\end{equation*}
Then
\begin{equation*}
\Phi^*A=\sum_{j_1,\dots,j_k=1}^n\sum_{a_1,\dots,a_k=1}^m
\Phi^{j_1}_{a_1}\cdots\Phi^{j_k}_{a_k}A_{j_1\cdots j_k}
f^{a_1}\otimes\dots\otimes f^{a_k}.
\end{equation*}
Denote
\begin{equation*}
\dnorm{\Phi}_\infty=\max\asetdef{\asnorm{\Phi^j_a}}{a\in\{1,\dots,m\},
j\in\{1,\dots,n\}}.
\end{equation*}
We have
\begin{align*}
\dnorm{\Phi^*A}^2_{\metric_{\alg{U}}}=&\;
\sum_{a_1,\dots,a_k=1}^m\left(\sum_{j_1,\dots,j_k=1}^n
\Phi^{j_1}_{a_1}\cdots\Phi^{j_k}_{a_k}A_{j_1\cdots j_k}\right)^2\\
\le&\;\sum_{a_1,\dots,a_k=1}^m\left(\sum_{j_1,\dots,j_k=1}^n
\asnorm{\Phi^{j_1}_{a_1}\cdots\Phi^{j_k}_{a_k}A_{j_1\cdots j_k}}\right)^2\\
\le&\;\sum_{a_1,\dots,a_k=1}^m\left(\sum_{j_1,\dots,j_k=1}^n
\asnorm{\Phi^{j_1}_{a_1}\cdots\Phi^{j_k}_{a_k}}^2\right)
\left(\sum_{j_1,\dots,j_k=1}^n\asnorm{A_{j_1\cdots j_k}}^2\right)\\
\le&\;(nm\dnorm{\Phi}_\infty^2)^k\dnorm{A}^2_{\metric_{\alg{V}}}.
\end{align*}
The first part of the result follows by taking
$C=\sqrt{nm}\dnorm{\Phi}_\infty$\@.

If $\Phi$ is surjective, let
$\Psi\in\Hom_{\real}(\alg{V};\alg{U})$ be a right-inverse for $\Phi$\@.
Then, by the first part of the result, there exists $C\in\realp$ such that
\begin{equation*}
\dnorm{A}_{\metric_{\alg{V}}}=\dnorm{(\Phi\scirc\Psi)^*A}_{\metric_{\alg{V}}}
\le\dnorm{\Psi^*\Phi^*A}_{\metric_{\alg{V}}}\le C^k\dnorm{\Phi^*A}_{\metric_{\alg{U}}}
\end{equation*}
for every $A\in\tensor*[k]{\dual{\alg{V}}}$\@, $k\in\integernn$\@.

If $\Phi$ is injective, we choose the orthonormal basis
$\ifam{e_1,\dots,e_n}$ so that $\ifam{e_1,\dots,e_m}$ is a basis for
$\image(\Phi)$\@.  In this case we have
\begin{equation*}
\Phi=\sum_{a,b=1}^m\Phi^b_ae_b\otimes f^a,
\end{equation*}
where the $m\times m$ matrix with components $\Phi^b_a$\@,
$a,b\in\{1,\dots,m\}$\@, is invertible, and
\begin{equation*}
\Phi^*A=\sum_{b_1,\dots,b_k=1}^m\sum_{a_1,\dots,a_k=1}^m
\Phi^{b_1}_{a_1}\cdots\Phi^{b_k}_{a_k}A_{b_1\cdots b_k}
f^{a_1}\otimes\dots\otimes f^{a_k}.
\end{equation*}
Letting $\Psi^b_a$\@, $a,b\in\{1,\dots,m\}$\@, be defined by
\begin{equation*}
\Psi^c_a\Phi^b_c=\begin{cases}1,&a=b,\\0,&a\not=b,\end{cases}
\end{equation*}
we have
\begin{equation*}
A=\sum_{b_1,\dots,b_k=1}^m\sum_{a_1,\dots,a_k=1}^m
\Psi^{b_1}_{a_1}\cdots\Psi^{b_k}_{a_k}(\Phi^*A)_{b_1\cdots b_k}
e^{a_1}\otimes\dots\otimes e^{a_k},
\end{equation*}
and the conclusion in this case follows just as in the first part of the
proof.
\end{subproof}
\end{proofsublemma}

To prove the first part of the lemma, let $x\in\nbhd{K}$ and take
$C_x\in\realp$ as in the sublemma such that
\begin{equation*}
\dnorm{\Phi^*\nabla^{\man{N},m}f(x)}_{\metric_{\man{M}}}\le
C_x^m\dnorm{\nabla^{\man{N},m}f(\Phi(x))}_{\metric_{\man{N}}},\qquad
m\in\integernn.
\end{equation*}
By continuity, and noting the exact form of the constant $C$ from the
sublemma (\ie~depending on the size of the derivative of $\tf[x]{\Phi}$),
there exists a neighbourhood $\nbhd{U}_x$ of $x$ such that
\begin{equation*}
\dnorm{\Phi^*\nabla^{\man{N},m}f(y)}_{\metric_{\man{M}}}\le
(2C_x)^m\dnorm{\nabla^{\man{N},m}f(\Phi(y))}_{\metric_{\man{N}}},\qquad
y\in\nbhd{N}_x,\ m\in\integerp.
\end{equation*}
Then take $x_1,\dots,x_k\in\nbhd{K}$ such that
$\nbhd{K}\subset\cup_{j=1}^k\nbhd{U}_{x_j}$\@.  The first part of the lemma
then follows by taking
\begin{equation*}
C=\max\{2C_{x_1},\dots,2C_{x_k}\}.
\end{equation*}

The second part of the lemma follows, \emph{mutatis mutandis}\@, from the
second part of the sublemma since, locally, a surjective submersion and an
injective immersion can be made linear in an appropriate set of
coordinates~\cite[Theorems~3.5.2 and~3.5.7]{RA/JEM/TSR:88}\@.
\end{proof}
\end{lemma}

We note that the fibre norm~\eqref{eq:Phi*fBnorm} makes use of the vector
bundle mapping
\begin{equation*}
B^m_{\nabla^{\man{E}}}\in\vbmappings[r]{\jetalg[\Phi]{m}{\man{M}}}
{\Symalg*[\le m]{\ctb{\man{M}}}}
\end{equation*}
from Lemma~\ref{lem:pbdecomp}\@.  If instead we use the vector bundle mapping
\begin{equation*}
A^m_{\nabla^{\man{E}}}\in\vbmappings[r]{\jetalg[\Phi]{m}{\man{M}}}
{\Symalg*[\le m]{\Phi^*\ctb{\man{N}}}}
\end{equation*}
from Lemma~\ref{lem:pbdecomp}\@, then we have the alternative fibre norm
\begin{equation*}
\dnorm{j_m\Phi^*f(e)}^{\prime2}_{\metric_{\man{M},m}}=
\sum_{j=0}^m\frac{1}{(j!)^2}
\dnorm{\Phi^*\nabla^{\man{N},j}f(e)}^2_{\metric_{\man{M}}}.
\end{equation*}
The relationship between the fibre norms
$\dnorm{\cdot}_{\metric_{\man{M},m}}$ and
$\dnorm{\cdot}'_{\metric_{\man{M},m}}$ can be phrased as, ``What is the
relationship between the jet of the pull-back and the pull-back of the jet?''
This is a question we will phrase below for other sorts of lifts, and will
address comprehensively in the proof of continuity of pull-back in
Theorem~\ref{the:compcont}\@.

\section{Estimates related to jet bundle norms}\label{sec:jet-estimates}

In Section~\ref{sec:lift-isomorphisms} we gave formulae relating derivatives
of geometric objects to derivatives of their lifts, and vice versa.  In
Section-\ref{sec:fibre-norms} we defined fibre metrics associated with spaces
of lifted objects.  In each of the multitude of constructions, there arose
certain vector bundle mappings that satisfied recursion relations.  In order
to establish some important comparison results for different
characterisations of topologies, we will need some rather detailed technical
estimates concerning the growth of these recursively defined vector bundle
mappings in the real analytic case, and we develop these here.  As a part of
this, we establish a number of fairly simple, linear algebraic estimates.  It
is not the existence of these estimates that are of interest, but the form
they take.  As we shall see, for the real analytic topology, the dimensions
of various tensor spaces show up in ways that need to be bookkept.

The results in this section are important, but somewhat elaborate.  Moreover,
they apply specifically to the real analytic setting.  The algebraic
computations and estimates of Section~\ref{subsec:tensornorms}\@, when
applied in the smooth setting, do not require the very particular forms we
give here.

\subsection{Algebraic estimates}\label{subsec:tensornorms}

To work with the topologies we present in Section~\ref{subsec:seminorms}\@,
we will have to compute and estimate high-order derivatives of various sorts
of tensors.  In this section we collect the fairly elementary formulae we
shall need.  All norms on tensor products are those induced by an inner
product as in Lemma~\ref{lem:inprodotimes}\@.  For simplicity, therefore, we
shall often omit any particular symbols attached to ``$\dnorm{\cdot}$'' to
connote which norm we are talking about; all vector spaces have a unique norm
(given the data) that we shall use.

We start by giving the norm of the identity mapping on tensors.
\begin{lemma}\label{lem:idnorm}
If\/ $\alg{V}$ is a finite-dimensional\/ $\real$-vector space with inner
product\/ $\metric$\@, then\/
$\dnorm{\id_{\alg{V}}}=\sqrt{\dim_{\real}(\alg{V})}$\@.
\begin{proof}
Let $\ifam{e_1,\dots,e_n}$ be an orthonormal basis for $\alg{V}$ with dual
basis $\ifam{e^1,\dots,e^n}$ the dual basis.  Write
\begin{equation*}
\id_{\alg{V}}=\sum_{j=1}^n\sum_{k=1}^n\delta^k_je_k\otimes e^j.
\end{equation*}
We have
\begin{equation*}
\dnorm{A}^2=\sum_{j=1}^n\sum_{k=1}^n(\delta^k_j)^2=n,
\end{equation*}
as claimed.
\end{proof}
\end{lemma}

Next we consider the norm of the tensor product of linear maps.
\begin{lemma}\label{lem:otimesnorm}
Let\/ $\alg{U}$\@,\/ $\alg{V}$\@,\/ $\alg{W}$\@, and\/ $\alg{X}$ be
finite-dimensional\/ $\real$-vector spaces with inner products.  Then, for\/
$A\in\Hom_{\real}(\alg{U};\alg{V})$ and\/
$B\in\Hom_{\real}(\alg{W};\alg{X})$\@,
\begin{equation*}
\dnorm{A\otimes B}=\dnorm{A}\dnorm{B}.
\end{equation*}
\begin{proof}
Let $\ifam{e_1,\dots,e_n}$\@, $\ifam{f_1,\dots,f_m}$\@,
$\ifam{g_1,\dots,g_k}$\@, and $\ifam{h_1,\dots,h_l}$ be orthonormal bases for
$\alg{U}$\@, $\alg{V}$\@, $\alg{W}$\@, and $\alg{X}$\@, respectively.  Let $\ifam{e^1,\dots,e^n}$\@, $\ifam{f^1,\dots,f^m}$\@,
$\ifam{g^1,\dots,g^k}$\@, and $\ifam{h^1,\dots,h^l}$ be the dual bases.
Write
\begin{equation*}
A=\sum_{j=1}^n\sum_{a=1}^mA^a_jf_a\otimes e^j,\quad
B=\sum_{i=1}^k\sum_{b=1}^lB^b_ih_b\otimes g^i
\end{equation*}
so that
\begin{equation*}
A\otimes B=\sum_{j=1}^n\sum_{i=1}^k\sum_{a=1}^m\sum_{b=1}^lA^a_jB^b_i
(f_a\otimes h_b)\otimes(e^j\otimes g^i).
\end{equation*}
Then
\begin{align*}
\dnorm{A\otimes B}^2=&\;
\sum_{j=1}^n\sum_{i=1}^k\sum_{a=1}^m\sum_{b=1}^l(A^a_jB^b_i)^2\\
\le&\;\left(\sum_{j=1}^n\sum_{a=1}^m(A^a_j)^s\right)
\left(\sum_{i=1}^k\sum_{b=1}^l(B^b_i)^2\right)=\dnorm{A}^2\dnorm{B}^2,
\end{align*}
as claimed.
\end{proof}
\end{lemma}

Our next estimate concerns the relationship between norms of tensors
evaluated on arguments.
\begin{lemma}\label{lem:A(S)<AS}
Let\/ $\alg{U}$ and\/ $\alg{V}$ be finite-dimensional\/ $\real$-vector spaces
with inner products\/ $\metric$ and\/ $\hmetric$\@, respectively.  Then
\begin{equation*}
\dnorm{L(u)}\le\dnorm{L}\,\dnorm{u}
\end{equation*}
for all linear mappings\/ $L\in\Hom_{\real}(\alg{U};\alg{V})$ and for all\/
$u\in\alg{U}$\@.
\begin{proof}
Let $\ifam{f_1,\dots,f_m}$ and $\ifam{e_1,\dots,e_n}$ be an orthonormal basis
for $\alg{U}$ and $\alg{V}$\@.  For $L\in\Hom_{\real}(\alg{U};\alg{V})$\@,
write
\begin{equation*}
L=\sum_{a=1}^m\sum_{j=1}^nL^j_ae_j\otimes f^a.
\end{equation*}
Then we compute, using Cauchy\textendash{}Schwarz,
\begin{align*}
\dnorm{L(u)}^2=&\;\sum_{j=1}^n\left(\sum_{a=1}^mL_a^ju^a\right)^2
\le\sum_{j=1}^n\left(\sum_{a=1}^m\asnorm{L_a^ju^a}\right)^2\\
\le&\;\sum_{j=1}^n\left(\sum_{a=1}^m\asnorm{L_a^j}^2\right)
\left(\sum_{a=1}^m\snorm{u^a}^2\right)=\dnorm{L}^2\dnorm{u}^2,
\end{align*}
giving the lemma.
\end{proof}
\end{lemma}

We shall also make use of a sort of ``reverse inequality'' related to the
above.
\begin{lemma}\label{lem:opnormupper}
Let\/ $\alg{U}$ and\/ $\alg{V}$ be finite-dimensional\/ $\real$-vector spaces
with inner products\/ $\metric_{\alg{U}}$ and\/ $\metric_{\alg{V}}$\@.  For\/
$L\in\Hom_{\real}(\alg{U};\alg{V})$\@,
\begin{equation*}
\dnorm{L}\le\sqrt{\dim_{\real}(\alg{U})}
\sup\setdef{\dnorm{L(u)}}{\dnorm{u}=1}.
\end{equation*}
\begin{proof}
The result is true with equality and without the constant if one uses the
induced norm for $\Hom_{\real}(\alg{U};\alg{V})$\@, rather than the tensor
norm as we do here.  So the statement of the lemma is really about relating
the induced norm with the tensor norm.

The tensor norm, in the case of linear mappings as we have here, is really
the Frobenius norm, and as such it is computed as the $\ell^2$-norm of the
vector of the set of $\dim_{\real}(\alg{U})$ eigenvalues of
$\sqrt{\transpose{L}\scirc L}$\@.  On the other hand, the induced norm is the
$\ell^\infty$ norm of this same vector of eigenvalues of
$\sqrt{\transpose{L}\scirc L}$\@.  These interpretations can be found
in~\cite[page~7]{RB:97}\@.  For this reason, an application
of~\eqref{eq:12inftynorms} gives the result.
\end{proof}
\end{lemma}

Another tensor estimate we shall find useful concerns symmetrisation.
\begin{lemma}\label{lem:Symnorm}
Let\/ $\alg{V}$ be a finite-dimensional\/ $\real$-vector space and let\/
$\metric$ be an inner product on\/ $\alg{V}$\@.  Then
\begin{equation*}
\dnorm{\Sym_k(A)}\le\dnorm{A}
\end{equation*}
for every\/ $A\in\tensor*[k]{\dual{\alg{V}}}$ and\/ $k\in\integerp$\@.
\begin{proof}
The result follows from the following sublemma.
\begin{proofsublemma}\label{psublem:Symorthproj}
The map\/
$\map{\Sym_k}{\tensor*[k]{\dual{\alg{V}}}}{\Symalg*[k]{\dual{\alg{V}}}}$ is
the orthogonal projection.
\begin{subproof}
Let us simply denote by $\metric$ the inner product on
$\tensor*[k]{\dual{\alg{V}}}$\@, defined as in
Lemma~\ref{lem:inprodotimes}\@.  It suffices to show that
$\metric(A,S)=\metric(\Sym_k(A),S)$ for every
$A\in\tensor*[k]{\dual{\alg{V}}}$ and $S\in\Symalg*[k]{\dual{\alg{V}}}$\@.
It suffices to show that this is true as $A$ runs over a set of generators
for $\tensor*[k]{\dual{\alg{V}}}$ and $S$ runs over a set of generators for
$\Symalg*[k]{\dual{\alg{V}}}$\@.

Thus we let $\ifam{e_1,\dots,e_n}$ be an orthonormal basis for $\alg{V}$ with
dual basis $\ifam{e^1,\dots,e^n}$\@.  Then we have generators
\begin{equation*}
e^{a_1}\otimes\dots\otimes e^{a_k},\qquad a_1,\dots,a_k\in\{1,\dots,n\},
\end{equation*}
for $\tensor*[k]{\dual{\alg{V}}}$ and
\begin{equation*}
\Sym_k(e^{b_1}\otimes\dots\otimes e^{b_k}),
\qquad b_1,\dots,b_k\in\{1,\dots,n\},
\end{equation*}
for $\Symalg*[k]{\dual{\alg{V}}}$\@.  For
$a_1,\dots,a_k,b_1,\dots,b_k\in\{1,\dots,n\}$\@, we wish to show that the
inner product
\begin{multline*}
\metric(e^{a_1}\otimes\dots\otimes e^{a_k},
\Sym_k(e^{b_1}\otimes\dots\otimes e^{b_k}))\\
=\frac{1}{k!}\sum_{\sigma\in\symmgroup{k}}
\metric(e^{a_1}\otimes\dots\otimes e^{a_k},
e^{b_{\sigma(1)}}\otimes\dots\otimes e^{b_{\sigma(k)}})\\
=\frac{1}{k!}\sum_{\sigma\in\symmgroup{k}}
\metric(e^{a_1},e^{b_{\sigma(1)}})\cdots\metric(e^{a_k},e^{b_{\sigma(k)}})
\end{multline*}
is equal to
\begin{equation*}
\metric(\Sym_k(e^{a_1}\otimes\dots\otimes e^{a_k}),
\Sym_k(e^{b_1}\otimes\dots\otimes e^{b_k})).
\end{equation*}
Unless $\{a_1,\dots,a_k\}$ and $\{b_1,\dots,b_k\}$ agree as multisets, we
have
\begin{equation*}
0=\metric(e^{a_1}\otimes\dots\otimes e^{a_k},
\Sym_k(e^{b_1}\otimes\dots\otimes e^{b_k}))=
\metric(\Sym_k(e^{a_1}\otimes\dots\otimes e^{a_k}),
\Sym_k(e^{b_1}\otimes\dots\otimes e^{b_k})).
\end{equation*}
Thus we can suppose that $\{a_1,\dots,a_k\}$ and $\{b_1,\dots,b_k\}$ agree as
multisets.

In this case, since
\begin{equation*}
\Sym_k(e^{a_1}\otimes\dots\otimes e^{a_k})=
\Sym_k(e^{b_1}\otimes\dots\otimes e^{b_k}),
\end{equation*}
we can assume, without loss of generality, that $a_j=b_j$\@,
$j\in\{1,\dots,k\}$\@.  For $l\in\{1,\dots,n\}$\@, let
$k^{\vect{a}}_l\in\integernn$ be the number of occurrences of $l$ in the list
$\ifam{a_1,\dots,a_k}$\@.  Let $\symmgroup{k}^{\vect{a}}\subset\symmgroup{k}$
be those permutations $\sigma$ for which $a_j=a_{\sigma(j)}$\@,
$j\in\{1,\dots,k\}$\@.  Note that
$\card(\symmgroup{k}^{\vect{a}})=k^{\vect{a}}_1!\cdots k^{\vect{a}}_n!$ since
$\symmgroup{k}^{\vect{a}}$ consists of compositions of permutations that
permute all the $1$'s, all the $2$'s,~\etc, in the list
$\ifam{a_1,\dots,a_k}$\@.  With these bits of notation, we have
\begin{equation*}
e^{a_1}\otimes\dots\otimes e^{a_k}=
e^{a_{\sigma(1)}}\otimes\dots\otimes e^{a_{\sigma(k)}}\enspace\iff\enspace
\sigma\in\symmgroup{k}^{\vect{a}}.
\end{equation*}
Therefore,
\begin{equation*}
\metric(e^{a_1}\otimes\dots\otimes e^{a_k},
e^{a_{\sigma(1)}}\otimes\dots\otimes e^{a_{\sigma(k)}})=
\begin{cases}1,&\sigma\in\symmgroup{k}^{\vect{a}},\\0,&\textrm{otherwise}.
\end{cases}
\end{equation*}
We then have
\begin{multline*}
\metric(e^{a_1}\otimes\dots\otimes e^{a_k},
\Sym_k(e^{a_1}\otimes\dots\otimes e^{a_k}))\\=
\frac{k^{\vect{a}}_1!\cdots k^{\vect{a}}_n!}{k!}
\metric(e^{a_1}\otimes\dots\otimes e^{a_k},
e^{a_1}\otimes\dots\otimes e^{a_k})
=\frac{k^{\vect{a}}_1!\cdots k^{\vect{a}}_n!}{k!}.
\end{multline*}

Next we calculate
\begin{equation*}
\metric(\Sym_k(e^{a_1}\otimes\dots\otimes e^{a_k}),
\Sym_k(e^{a_1}\otimes\dots\otimes e^{a_k})).
\end{equation*}
Let $\sigma\in\symmgroup{k}$ and, for $l\in\{1,\dots,n\}$\@, let
$k^{\sigma(\vect{a})}_l\in\integernn$ be the number of occurrences of
$l$ in the list $\ifam{a_{\sigma(1)},\dots,a_{\sigma(k)}}$\@.  Let
$\symmgroup{k}^{\sigma(\vect{a})}\subset\symmgroup{k}$ be those permutations
$\sigma'$ for which $a_{\sigma(j)}=a_{\sigma'(j)}$\@, $j\in\{1,\dots,k\}$\@.
As above,
$\card(\symmgroup{k}^{\sigma(\vect{a})})= k^{\sigma(\vect{a})}_1!\cdots
k^{\sigma(\vect{a})}_n!$\@.  Also as above, we then have
\begin{equation*}
\metric(e^{\sigma(1)}\otimes\dots\otimes e^{\sigma(k)},
\Sym_k(e^{a_1}\otimes\dots\otimes e^{a_k}))=
\frac{k_1^{\sigma(\vect{a})}!\cdots k_n^{\sigma(\vect{a})}!}{k!}=
\frac{k^{\vect{a}}_1!\cdots k^{\vect{a}}_n!}{k!},
\end{equation*}
if $k^{\vect{a}}_1,\dots,k^{\vect{a}}_n$ are as in the preceding paragraph.
Therefore,
\begin{align*}
\metric(\Sym_k(e^{a_1}\otimes\dots\otimes e^{a_k}),&
\Sym_k(e^{a_1}\otimes\dots\otimes e^{a_k}))\\
=&\;\frac{1}{k!}\sum_{\sigma\in\symmgroup{k}}
\metric(e^{\sigma(1)}\otimes\dots\otimes e^{\sigma(k)},
\Sym_k(e^{a_1}\otimes\dots\otimes e^{a_k}))\\
=&\;\frac{1}{k!}\sum_{\sigma\in\symmgroup{k}}
\frac{k^{\vect{a}}_1!\cdots k^{\vect{a}}_n!}{k!}=
\frac{k^{\vect{a}}_1!\cdots k^{\vect{a}}_n!}{k!},
\end{align*}
and so we have
\begin{multline*}
\metric(e^{a_1}\otimes\dots\otimes e^{a_k},
\Sym_k(e^{b_1}\otimes\dots\otimes e^{b_k}))\\
=\metric(\Sym_k(e^{a_1}\otimes\dots\otimes e^{a_k}),
\Sym_k(e^{a_1}\otimes\dots\otimes e^{a_k})),
\end{multline*}
and the sublemma follows.
\end{subproof}
\end{proofsublemma}

Now, given $A\in\tensor*[k]{\dual{\alg{V}}}$\@, we write $A=\Sym_k(A)+A_1$
where $A_1$ is orthogonal to $\Symalg*[k]{\dual{\alg{V}}}$\@.  We then have
$\dnorm{A}^2=\dnorm{\Sym_k(A)}^2+\dnorm{A_1}^2$\@, from which the lemma
follows.
\end{proof}
\end{lemma}

The sublemma from the preceding lemma is proved, differently, by
\citet[page~124]{JWN:68}\@.

Let us also determine the norm of various insertion operators that we shall
use.  We shall use notation that is specific to the manner in which we shall
use these estimates, and this will seem unmotivated out of context.  Let
$\alg{U}$\@, $\alg{V}$\@, and $\alg{W}$ be finite-dimensional $\real$-vector
spaces, let $m,s,r\in\integerp$ and $a\in\{0,1,\dots,r\}$\@, let
$S\in\tensor[1,r-a+2]{\alg{U}}$\@, and let
\begin{equation*}
A\in\tensor[s+r-a+1,m+a+1]{\alg{U}}\otimes\alg{W}\otimes\dual{\alg{V}}.
\end{equation*}
We then have the mapping
\begin{equation*}
\map{I^1_{A,S,j}}{\tensor*[s]{\dual{\alg{U}}}\otimes\alg{V}}
{\tensor*[m+r+1]{\dual{\alg{U}}}\otimes\alg{W}}
\end{equation*}
defined by
\begin{equation*}
I^1_{A,S,j}(\beta)=A(\Ins_j(\beta,S)).
\end{equation*}
Here we implicitly use the isomorphism
\begin{equation*}
\map{\kappa}{\tensor[s,m]{\alg{U}}}
{\Hom_{\real}(\tensor*[s]{\dual{\alg{U}}};\tensor*[m]{\dual{\alg{U}}})},
\end{equation*}
for a finite-dimensional $\real$-vector space $\alg{U}$ and for $m,s\in\integernn$\@, via
\begin{equation*}
\kappa(v_1\otimes\dots\otimes v_s\otimes\alpha^1\otimes\dots\otimes\alpha^m)
(\beta^1\otimes\dots\otimes\beta^s)=
\natpair{\beta^1}{v_1}\cdots\natpair{\beta^s}{v_s}
\alpha^1\otimes\dots\otimes\alpha^m,
\end{equation*}
for $v_a\in\alg{U}$\@, $a\in\{1,\dots,s\}$\@, and
$\alpha^j,\beta^b\in\dual{\alg{U}}$\@, $b\in\{1,\dots,s\}$\@,
$j\in\{1,\dots,m\}$\@.  Thus, for additional finite-dimensional
$\real$-vector spaces $\alg{V}$ and $\alg{W}$\@, we have the identification
\begin{equation*}
\tensor[m,s]{\alg{U}}\otimes\alg{W}\otimes\dual{\alg{V}}\simeq
\Hom_{\real}(\tensor*[s]{\dual{\alg{U}}}\otimes\alg{V};
\tensor*[m]{\dual{\alg{U}}}\otimes\alg{W}).
\end{equation*}

We now have the following result.
\begin{lemma}\label{lem:I1AS}
With the preceding notation,
\begin{equation*}
\dnorm{I^1_{A,S,j}}\le\dnorm{A}\dnorm{S}.
\end{equation*}
\begin{proof}
Let $\ifam{f_1,\dots,f_m}$ be an orthonormal basis for $\alg{U}$ with dual
basis $\ifam{f^1,\dots,f^m}$\@.  Let $\ifam{e_1,\dots,e_n}$ be an orthonormal
basis for $\alg{V}$ with $\ifam{e^1,\dots,e^n}$ the dual basis.  Let
$\ifam{g_1,\dots,g_k}$ be an orthonormal basis for $\alg{W}$ with
$\ifam{g^1,\dots,g^k}$ the dual basis.  Let us write
\begin{equation*}
S=\sum_{a=1}^m\sum_{a_1,\dots,a_{r-a+2}}^mS^a_{a_1\cdots a_{r-a+2}}
f_a\otimes f^{a_1}\otimes\dots\otimes f^{a_{r-a+2}}
\end{equation*}
and
\begin{multline*}
A=\sum_{a_1,\dots,a_{s+r-a+1}=1}^m\sum_{b_1,\dots,b_{m+a+1}=1}^m
\sum_{\alpha=1}^k\sum_{l=1}^nA_{b_1\cdots b_{m+a+1}l}^{a_a\cdots a_{s+r-a+1}\alpha}\\
\times f^{b_1}\otimes\dots\otimes f^{b_{m+a+1}}\otimes
f_{a_1}\otimes\dots\otimes f_{a_{s+1}}\otimes g_\alpha\otimes e^l.
\end{multline*}
We then have, for $a_1,\dots,a_s\in\{1,\dots,m\}$\@,
$\alpha\in\{1,\dots,k\}$\@, and $l\in\{1,\dots,n\}$\@,
\begin{align*}
\Ins_{S,j}(f^{a_1}&\otimes\dots\otimes f^{a_s}\otimes
g_\alpha\otimes e^l)\\
=&\;\Ins_j(f^{a_1}\otimes\dots\otimes f^{a_j}\otimes\dots
\otimes f^{a_s}\otimes g_\alpha\otimes e^l,S)\\
=&\;\sum_{b_1,\dots,b_{r-a+2}=1}^mS^{a_j}_{b_1\cdots b_{r-a+2}}\\
&\;\times f^{a_1}\otimes\dots\otimes f^{a_{j-1}}
\otimes f^{b_1}\otimes\dots\otimes f^{b_{r-a+2}}\otimes f^{a_{j+1}}
\otimes\dots\otimes f^{a_s}\otimes g_\alpha\otimes e^l\\
=&\;\sum_{c_1,\dots,c_{j-1}=1}^m\sum_{c_{j+1},\dots,c_s=1}^m
\sum_{b_1,\dots,b_{r-a+2}=1}^n\sum_{\beta=1}^k\sum_{p=1}^n
S^{a_j}_{b_1\cdots b_{r-a+2}}
\delta^{a_1}_{c_1}\cdots\delta^{a_{j-1}}_{c_{j-1}}
\delta^{a_{j+1}}_{c_{j+1}}\cdots\delta^{a_s}_{c_s}\delta^\beta_\alpha\delta^l_p\\
&\;\times f^{c_1}\otimes\dots\otimes f^{c_{j-1}}\otimes f^{b_1}\otimes
\dots\otimes f^{b_{r-a+2}}\otimes f^{c_{j+1}}\otimes\dots\otimes f^{c_s}
\otimes g_\beta\otimes e^p.
\end{align*}
Thus
\begin{multline*}
I^1_{A,S,j}(f^{a_1}\otimes\dots\otimes f^{a_s}\otimes
g_\alpha\otimes e^l)\\
=\sum_{b_1,\dots,b_{r-a+2}=1}^m\sum_{d_1,\dots,d_{m+a+1}=1}^m
\sum_{\alpha=1}^k\sum_{l=1}^n
A^{a_1\cdots a_{j-1}b_1\cdots b_{r-a+2}a_{j+1}\cdots a_s\alpha}_{d_1\cdots d_{m+a+1}l}
S^{a_j}_{b_1\cdots b_{r-a+2}}\\
\times f^{d_1}\otimes\dots\otimes f^{d_{m+a+1}}\otimes g_\alpha\otimes e^l.
\end{multline*}
Then we calculate, using Cauchy\textendash{}Schwarz,
{\small\begin{align*}
\dlnorm&I^1_{A,S,j}\drnorm^2=
\sum_{a_1,\dots,a_s=1}^m\sum_{d_1,\dots,d_{m+a+1}=1}^m
\sum_{\alpha=1}^k\sum_{l=1}^n\left(\sum_{b_1,\dots,b_{r-a+2}=1}^m
A^{a_1\cdots a_{j-1}b_1\cdots b_{r-a+2}a_{j+1}\cdots a_s\alpha}_{d_1\cdots d_{m+a+1}l}
S^{a_j}_{b_1\cdots b_{r-a+2}}\right)^2\\
&\le\hspace{1em}\sum_{a_1,\dots,a_s=1}^m\sum_{d_1,\dots,d_{m+a+1}=1}^m
\sum_{\alpha=1}^k\sum_{l=1}^n\left(\sum_{b_1,\dots,b_{r-a+2}=1}^m\asnorm{
A^{a_1\cdots a_{j-1}b_1\cdots b_{r-a+2}a_{j+1}\cdots a_s\alpha}_{d_1\cdots d_{m+a+1}l}
S^{a_j}_{b_1\cdots b_{r-a+2}}}\right)^2\\
&\le\hspace{1em}\sum_{a_1,\dots,a_s=1}^m\sum_{d_1,\dots,d_{m+a+1}=1}^m
\sum_{\alpha=1}^k\sum_{l=1}^n\left(\sum_{b_1,\dots,b_{r-a+2}=1}^m\asnorm{
A^{a_1\cdots a_{j-1}b_1\cdots b_{r-a+2}a_{j+1}\cdots a_s}_{d_1\cdots d_{m+a+1}l}
}^2\right)\\
&\phantom{\le}\hspace{1em}\times
\left(\sum_{b_1,\dots,b_{r-a+2}=1}^m
\asnorm{S^{a_j}_{b_1\cdots b_{r-a+2}}}^2\right)\\
&\le\hspace{1em}\dnorm{A}^2\dnorm{S}^2,
\end{align*}}%
as claimed.
\end{proof}
\end{lemma}

Now we perform the same sort of estimate for a similar construction.  We take
$\alg{U}$\@, $\alg{V}$\@, and $\alg{W}$ as above, and $m$\@, $s$\@, $r$\@,
and $a$ as above.  We also still take $S\in\tensor[1,r-a+2]{\alg{U}}$\@, but
here we take
\begin{equation*}
B\in\tensor[s,m+a]{\alg{U}}\otimes\alg{W}\otimes\dual{\alg{V}}.
\end{equation*}
We then have the mapping
\begin{equation*}
\map{I^2_{B,S,j}}{\tensor*[s]{\dual{\alg{U}}}\otimes\alg{V}}
{\tensor*[m+r+1]{\dual{\alg{U}}}\otimes\alg{W}}
\end{equation*}
defined by
\begin{equation*}
I^2_{B,S,j}(\beta)=\Ins_j(B(\beta),S)
\end{equation*}
We now have the following result, whose proof follows from direct
computation, just as does Lemma~\ref{lem:I1AS}\@.
\begin{lemma}\label{lem:I2BS}
With the preceding notation,
\begin{equation*}
\dnorm{I^2_{B,S,j}}\le\dnorm{B}\dnorm{S}.
\end{equation*}
\end{lemma}

\subsection{Tensor field estimates}

We next turn to providing estimates for the tensors $A^m_s$\@, $B^m_s$\@,
$C^m_s$\@, and $D^m_s$\@, $m\in\integernn$\@, $s\in\{0,1,\dots,m\}$\@, that
appear in the lemmata from Section~\ref{sec:lift-isomorphisms}\@.  In this
section is where all of our seemingly pointless computations from
Sections~\ref{sec:tensor-constructions} and~\ref{sec:tensor-derivatives}\@,
and our only slightly less seemingly pointless constructions from
Sections~\ref{sec:lift-isomorphisms} and~\ref{sec:fibre-norms}\@, bear fruit.
We first develop a general estimate, and then show how this estimate can be
made to apply to all of the required tensors from
Section~\ref{sec:lift-isomorphisms}\@.

We work with real analytic vector bundles
$\map{\pi_{\man{E}}}{\man{E}}{\man{M}}$ and
$\map{\pi_{\man{F}}}{\man{F}}{\man{M}}$\@.  The r\^ole of
$\map{\pi_{\man{E}}}{\man{E}}{\man{M}}$ in this discussion and that in
Section~\ref{sec:lift-isomorphisms} is different.  One should think of
$\man{E}$ in Section~\ref{sec:lift-isomorphisms} as being played by $\man{M}$
here.  This is because the tensors in Section~\ref{sec:lift-isomorphisms} are
defined as having $\man{E}$ as their base space.  So here we rename this base
space as $\man{M}$\@.  As a consequence of this, one should think of (1)~the
r\^ole of $\man{M}$ in the lemma below as being played by $\man{E}$ in the
lemmata of Section~\ref{sec:lift-isomorphisms}\@, (2)~the r\^ole of
$\nabla^{\man{M}}$ in the lemma below as being played by $\nabla^{\man{E}}$
in the lemmata of Section~\ref{sec:lift-isomorphisms}\@, and (2)~the r\^ole
of $\nabla^{\pi_{\man{E}}}$ in the lemma below as being played by the induced
connection in an appropriate tensor bundle in the lemmata of
Section~\ref{sec:lift-isomorphisms}\@.  In our development here, we use the
symbol $\nabla^{\man{M},\pi_{\man{E}}}$ to denote the connection induced in
any of the myriad bundles formed by taking tensor products of
$\tb{\man{M}}$\@, $\ctb{\man{M}}$\@, $\man{E}$\@, and
$\dual{\man{E}}$\@,~\cf~the constructions at the beginning of
Section~\ref{subsec:jetdecomp}\@.

With this as backdrop, the main technical result we have is the following.
\begin{lemma}\label{lem:pissynabla}
Let\/ $\map{\pi_{\man{E}}}{\man{E}}{\man{M}}$ and\/
$\map{\pi_{\man{F}}}{\man{F}}{\man{M}}$ be real analytic vector bundles, let\/
$\nabla^{\man{M}}$ be a real analytic affine connection on\/ $\man{M}$\@,
let\/ $\nabla^{\pi_{\man{E}}}$ and\/ $\nabla^{\pi_{\man{F}}}$ be real
analytic vector bundle connections in\/ $\man{E}$ and\/ $\man{F}$\@,
respectively.  Let\/ $\metric_{\man{M}}$ be a real analytic Riemannian metric
on\/ $\man{M}$\@, and let\/ $\metric_{\pi_{\man{E}}}$ and\/
$\metric_{\pi_{\man{F}}}$ be real analytic fibre metrics for\/ $\man{E}$
and\/ $\man{F}$\@, respectively.  Suppose that we are given the following
data:
\begin{compactenum}[(i)]
\item $\phi_m\in\sections[\omega]{\tensor[m,m]{\tb{\man{M}}}\otimes
\man{F}\otimes\dual{\man{E}}}$\@,\/ $m\in\integernn$\@;
\item $\Phi^s_m\in\sections[\omega]{\End(\tensor[s,m+1]{\tb{\man{M}}}\otimes
\man{F}\otimes\dual{\man{F}})}$\@,\/ $m\in\integernn$\@,\/
$s\in\{0,1,\dots,m\}$\@;
\item $\Psi^s_{jm}\in\sections[\omega]{\Hom(\tensor[s,m]{\tb{\man{M}}}
\otimes\man{F}\otimes\dual{\man{E}};\tensor[s,m+1]{\tb{\man{M}}}\otimes
\man{F}\otimes\dual{\man{E}})}$\@,\/ $m\in\integernn$\@,\/
$s\in\{0,1,\dots,m\}$\@,\/ $j\in\{0,1,\dots,m\}$\@;
\item $\Lambda^s_m\in\sections[\omega]{\Hom(\tensor[s-1,m]{\tb{\man{M}}}
\otimes\man{F}\otimes\dual{\man{E}};\tensor[s,m+1]{\tb{\man{M}}}\otimes
\man{F}\otimes\dual{\man{E}})}$\@,\/ $m\in\integerp$\@,\/
$s\in\{1,\dots,m\}$\@;
\item $A^m_s\in\sections[\omega]{\tensor[s,m]{\tb{\man{M}}}\otimes
\man{F}\otimes\dual{\alg{E}}}$\@,\/ $m\in\integernn$\@,\/
$s\in\{0,1,\dots,m\}$\@,
\end{compactenum}
and that the data satisfies the recursion relations prescribed
by\/ $A^0_0=\phi_0$ and
\begin{xalignat*}{2}
A^{m+1}_{m+1}=&\;\Phi_m^{m+1}\scirc\phi_{m+1},&&m\in\integernn\\
A^{m+1}_s=&\;\Phi_m^s\scirc\nabla^{\man{M},\pi_{\man{E}}\otimes\pi_{\man{F}}}A^m_s+
\sum_{j=0}^m\Psi_{jm}^s\scirc A^m_s+\Lambda^s_m\scirc A^m_{s-1},&&m\in\integerp,\ s\in\{1,\dots,m\},\\
A^{m+1}_0=&\;\Phi^0_m\scirc\nabla^{\man{M},\pi_{\man{E}}\otimes\pi_{\man{F}}}A^m_0+
\sum_{j=0}^m\Psi^0_{jm}\scirc A^m_0,&&m\in\integernn.
\end{xalignat*}
Suppose that the data are such that, for each compact\/
$\nbhd{K}\subset\man{M}$\@, there exist\/ $C_1,\sigma_1\in\realp$ satisfying
\begin{compactenum}[(i)]
\item
$\dnorm{D^r_{\nabla^{\man{M}},\nabla^{\pi_{\man{E}}\otimes\pi_{\man{F}}}}
\phi_m(x)}_{\metric_{\man{M},\pi_{\man{E}}\otimes\pi_{\man{F}}}}\le
C_1\sigma_1^{-r}r!$\@,\/ $m,r\in\integernn$\@;

\item
$\dnorm{D^r_{\nabla^{\man{M}},\nabla^{\pi_{\man{F}}}}
\Phi_m^s(x)\scirc A}_{\metric_{\man{M},\pi_{\man{E}}\otimes\pi_{\man{F}}}}\le
C_1\sigma_1^{-r}r!
\dnorm{A}_{\metric_{\man{M},\pi_{\man{E}}\otimes\pi_{\man{F}}}}$\@,\/
$A\in\tensor[s,m+a]{\tb[x]{\man{M}}\otimes\man{F}_x
\otimes\dual{\man{E}}_x}$\@,\/ $m,r,a\in\integernn$\@,\/
$s\in\{0,1,\dots,m+1\}$\@;

\item $\dnorm{D^r_{\nabla^{\man{M}},\nabla^{\pi_{\man{E}}\otimes\pi_{\man{F}}}}
\Psi_{jm}^s(x)\scirc A}_{\metric_{\man{M},\pi_{\man{E}}\otimes\pi_{\man{F}}}}\le
C_1\sigma_1^{-r}r!
\dnorm{A}_{\metric_{\man{M},\pi_{\man{E}}\otimes\pi_{\man{F}}}}$\@,\/
$A\in\tensor[s,m+a]{\tb[x]{\man{M}}\otimes\man{F}_x\otimes
\dual{\man{E}}_x}$\@,\/ $m,r,a\in\integernn$\@,\/ $s\in\{0,1,\dots,m\}$\@,\/ $j\in\{0,1,\dots,m\}$\@;

\item $\dnorm{D^r_{\nabla^{\man{M}},\nabla^{\pi_{\man{E}}\otimes\pi_{\man{F}}}}
\Lambda_m^s(x)\scirc A}_{\metric_{\man{M},\pi_{\man{E}}\otimes\pi_{\man{F}}}}\le
C_1\sigma_1^{-r}r!
\dnorm{A}_{\metric_{\man{M},\pi_{\man{E}}\otimes\pi_{\man{F}}}}$\@,\/
$A\in\tensor[s-1,m+a]{\tb[x]{\man{M}}\otimes\man{F}_x
\otimes\dual{\man{E}}_x}$\@,\/ $m,r,a\in\integernn$\@,\/ $s\in\{0,1,\dots,m\}$\@.
\end{compactenum}
for\/ $x\in\nbhd{K}$\@.

Then, for\/ $\nbhd{K}\subset\man{M}$ compact, there exist\/
$C,\sigma,\rho\in\realp$ such that
\begin{equation*}
\dnorm{D^r_{\nabla^{\man{M}},\nabla^{\pi_{\man{E}}\otimes\pi_{\man{F}}}}
A^m_s(x)}_{\metric_{\man{M},\pi_{\man{E}}\otimes\pi_{\man{F}}}}\le
C\sigma^{-m}\rho^{-(m+r-s)}(m+r-s)!
\end{equation*}
for\/ $m,r\in\integernn$\@,\/ $s\in\{0,1,\dots,m\}$\@, and\/
$x\in\nbhd{K}$\@.
\begin{proof}
We prove the lemma with a sort of meandering induction, covering various
special cases of $m$ and $s$ before giving a proof for the general case.

Before we embark on the proof, we organise some data that will arise in the
estimate that we prove.
\begin{compactenum}
\item We take $\nbhd{K}\subset\man{M}$ compact and define
$C_1,\sigma_1\in\realp$ as in the statement of the lemma.  We shall assume,
without loss of generality, that $C_1>1$ and $\sigma_1<1$\@.

\item Choose $\beta\in\realp$ sufficiently large that
\begin{equation*}
\sum_{k=0}^\infty\beta^{-k}<\infty,
\end{equation*}
and let $\alpha=\frac{\beta}{\beta-1}>1$ denote the value of this sum.  Let
$\gamma=4\alpha$\@.

\item \label{enum:(a+b)!/b!} We note that, for any\/ $a,b,c\in\integerp$ with
$b<c$\@, we have
\begin{equation*}
\frac{(a+b)!}{b!}<\frac{(a+c)!}{c!}.
\end{equation*}
This is a direct computation:
\begin{equation*}
\frac{(a+b)!}{b!}=(1+b)\cdots(a+b)<(1+c)\cdots(a+c)=\frac{(a+c)!}{c!}.
\end{equation*}

\item \label{enum:Cms} For $m\in\integernn$ and $s\in\{0,1,\dots,m\}$\@, we
denote
\begin{equation*}
C_{m,s}=\begin{cases}1,&m=0\ \textrm{or}\ s=0,\\
\binom{m-1}{s-1},&\textrm{otherwise}.\end{cases}
\end{equation*}
We note that
\begin{compactenum}
\item $C_{m,m}=1$\@, that
\item $C_{m,s}\le C_{m+1,s}$\@, that
\item $C_{m,s}\le C_{m+1,s+1}$\@, and that
\item $mC_{m,s}\le(m+1-s)C_{m+1,s}$\@.
\end{compactenum}
The first and second of these assertions is obvious.  For the third, for $m,s\in\integerp$ with
$s\le m$\@, we compute
\begin{equation*}
C_{m,s}=\frac{(m-1)!}{(s-1)!(m-s)!}\le
\frac{m}{s}\frac{(m-1)!}{(s-1)!(m-s)!}=C_{m+1,s+1}.
\end{equation*}
For the fourth, for $m\in\integerp$ and $s\in\integerp$ satisfying
$s\le m$\@, we compute
\begin{equation*}
mC_{m,s}=m\frac{(m-1)!}{(s-1)!(m-s)!}=(m-s+1)\frac{m!}{(s-1)!(m+1-s)!}=
(m-s+1)C_{m+1,s}.
\end{equation*}

\item We shall have occasion below, and also subsequently, to use a standard
multinomial estimate.  First let $\alpha_1,\dots,\alpha_n\in\realp$ and note
that
\begin{equation*}
(\alpha_1+\dots+\alpha_n)^m=\sum_{m_1+\dots+m_n=m}\frac{m!}{m_1!\cdots m_n!}
\alpha_1^{m_1}\cdots\alpha_n^{m_n}.
\end{equation*}
Taking $\alpha_1=\dots=\alpha_n=1$\@, we see that
\begin{equation}\label{eq:multinomialest}
\frac{m!}{m_1!\cdots m_n!}\le n^m
\end{equation}
whenever $m_1,\dots,m_n\in\integernn$ sum to $m$\@.
\end{compactenum}
Given all of this, we shall prove that
\begin{equation}\label{eq:DrAmsest}
\dnorm{D^r_{\nabla^{\man{M}},\nabla^{\pi_{\man{E}}\otimes\pi_{\man{F}}}}
A^m_s(x)}_{\metric_{\man{M},\pi_{\man{E}}\otimes\pi_{\man{F}}}}\le
C_1(C_1\sigma_1^{-1}\gamma)^mC_{m,s}
\left(\frac{\beta}{\sigma_1}\right)^{m+r-s}(m+r-s)!
\end{equation}
for $m,r\in\integernn$\@, $s\in\{0,1,\dots,m\}$\@, and $x\in\nbhd{K}$\@.

\paragraph{\mathversion{bold}Case $m=s=0$\@:}

Directly using the hypotheses, we have
\begin{align*}
\dnorm{D^r_{\nabla^{\man{M}},\nabla^{\pi_{\man{E}}\otimes\pi_{\man{F}}}}
A^0_0(x)}_{\metric_{\man{M},\pi_{\man{E}}\otimes\pi_{\man{F}}}}=&\;
\dnorm{D^r_{\nabla^{\man{M}},\nabla^{\pi_{\man{E}}\otimes\pi_{\man{F}}}}
\phi_0(x)}_{\metric_{\man{M},\pi_{\man{E}}\otimes\pi_{\man{F}}}}\\
\le&\;C_1\sigma_1^{-r}r!\le C_1(C_1\sigma_1^{-1}\gamma)^0C_{0,0}
\left(\frac{\beta}{\sigma_1}\right)^{0+r-0}(0+r-0)!
\end{align*}
for $r\in\integernn$ and $x\in\nbhd{K}$\@.  This gives~\eqref{eq:DrAmsest} in
this case.

\paragraph{\mathversion{bold}Case $m\in\integerp$ and $s=m$\@:}

By Lemma~\ref{lem:leibniz}\@, we have
\begin{equation*}
D^r_{\nabla^{\man{M}},\nabla^{\pi_{\man{E}}\otimes\pi_{\man{F}}}}A^m_m=
D^r_{\nabla^{\man{M}},\nabla^{\pi_{\man{E}}\otimes\pi_{\man{F}}}}
(\Phi^m_{m-1}\scirc\phi_m)
=\sum_{a=0}^r\binom{r}{a}D^a_{\nabla^{\man{M}},\nabla^{\pi_{\man{E}}}}
\Phi^m_{m-1}(D^{r-a}_{\nabla^{\man{M}},\nabla^{\pi_{\man{E}}\otimes\pi_{\man{F}}}}
\phi_m)
\end{equation*}
for $m,r\in\integernn$\@.  Therefore, by Lemma~\ref{lem:A(S)<AS}\@, using the
hypotheses, and by the preliminary observation~\ref{enum:(a+b)!/b!} above,
\begin{align*}
\dnorm{D^r_{\nabla^{\man{M}},\nabla^{\pi_{\man{E}}\otimes\pi_{\man{F}}}}
A^m_m}_{\metric_{\man{M},\pi_{\man{E}}\otimes\pi_{\man{F}}}}
\le&\;\sum_{a=0}^r\frac{r!}{a!(r-a)!}(C_1\sigma_1^{-a}a!)
(C_1\sigma_1^{-(m+r-a)}(r-a)!)\\
\le&\;C_1C_1\sigma_1^{-m}r!\sum_{a=0}^r\sigma_1^{-a}
\left(\frac{\beta}{\sigma_1}\right)^{r-a}
\le C_1C_1\sigma_1^{-m}
\left(\frac{\beta}{\sigma_1}\right)^rr!\sum_{a=0}^r\beta^{-a}\\
\le&\;C_1(C_1\sigma_1^{-1}\gamma)^mC_{m,m}
\left(\frac{\beta}{\sigma_1}\right)^{m+r-m}(m+r-m)!.
\end{align*}
As this holds for every $m\in\integerp$\@, $r\in\integernn$\@, and
$x\in\nbhd{K}$\@, this gives~\eqref{eq:DrAmsest} in this case.

\paragraph{\mathversion{bold}Case $m=1$ and $s=0$\@:}

By Lemma~\ref{lem:leibniz} we have
\begin{multline*}
D^r_{\nabla^{\man{M}},\nabla^{\pi_{\man{E}}\otimes\pi_{\man{F}}}}A^1_0=
\underbrace{\sum_{a=0}^r\binom{r}{a}
D^a_{\nabla^{\man{M}},\nabla^{\pi_{\man{E}}}}\Phi^0_0
(D^{r-a}_{\nabla^{\man{M}},\nabla^{\pi_{\man{E}}\otimes\pi_{\man{F}}}}
\nabla^{\man{M},\pi_{\man{E}}\otimes\pi_{\man{F}}}A^0_0)}_{\textrm{term~1}}\\
+\underbrace{\sum_{a=0}^r\binom{r}{a}
D^a_{\nabla^{\man{M}},\nabla^{\pi_{\man{E}}\otimes\pi_{\man{F}}}}\Psi^0_{00}
(D^{r-a}_{\nabla^{\man{M}},\nabla^{\pi_{\man{E}}\otimes\pi_{\man{F}}}}
A^0_0)}_{\textrm{term~2(a)}}
+\underbrace{\sum_{a=0}^r\binom{r}{a}
D^a_{\nabla^{\man{M}},\nabla^{\pi_{\man{E}}\otimes\pi_{\man{F}}}}\Psi^0_{10}
(D^{r-a}_{\nabla^{\man{M}},\nabla^{\pi_{\man{E}}\otimes\pi_{\man{F}}}}
A^0_0)}_{\textrm{term~2(b)}}.
\end{multline*}

As we showed in the proof of Lemma~\ref{lem:Jk+mdecomp}\@, we have
\begin{equation*}
D^{r-a}_{\nabla^{\man{M}},\nabla^{\pi_{\man{E}}\otimes\pi_{\man{F}}}}
(\nabla^{\man{M},\pi_{\man{E}}\otimes\pi_{\man{F}}}A_0^0)=
D^{r-a+1}_{\nabla^{\man{M},\pi_{\man{E}}\otimes\pi_{\man{F}}}}A_0^0.
\end{equation*}
Therefore, by Lemma~\ref{lem:A(S)<AS} and using the hypotheses,
\begin{align*}
\dnorm{\textrm{term~1}(x)}_{\metric_{\man{M},\pi_{\man{E}}\otimes\pi_{\man{F}}}}
\le&\;\sum_{a=0}^r\frac{r!}{a!(r-a)!}(C_1\sigma_1^{-a}a!)
(C_1\sigma_1^{-(r-a+1)}(r-a+1)!)\\
\le&\;C_1C_1(r+1)!\sum_{a=0}^r\sigma_1^{-a}
\left(\frac{\beta}{\sigma_1}\right)^{r-a+1}\le
C_1C_1\left(\frac{\beta}{\sigma_1}\right)^{r+1}(r+1)!\sum_{a=0}^r\beta^{-a}\\
\le&\;C_1(C_1\alpha)\left(\frac{\beta}{\sigma_1}\right)^{r+1}(r+1)!.
\end{align*}
In a similar manner,
\begin{equation*}
\dnorm{\textrm{term~2(a)}(x)}_{\metric_{\man{M},\pi_{\man{E}}\otimes\pi_{\man{F}}}},
\dnorm{\textrm{term~2(b)}(x)}_{\metric_{\man{M},\pi_{\man{E}}\otimes\pi_{\man{F}}}}
\le C_1(C_1\alpha)\left(\frac{\beta}{\sigma_1}\right)^rr!.
\end{equation*}
Therefore, for $r\in\integernn$ and $x\in\nbhd{K}$\@,
\begin{align*}
\dnorm{D^r_{\nabla^{\man{M}},\nabla^{\pi_{\man{E}}\otimes\pi_{\man{F}}}}
A^1_0(x)}_{\metric_{\man{M},\pi_{\man{E}}\otimes\pi_{\man{F}}}}\le&\;
\frac{1}{4}C_1(C_1\gamma)\left(\frac{\beta}{\sigma_1}\right)^{r+1}(r+1)!+
\frac{1}{2}C_1(C_1\gamma)\left(\frac{\beta}{\sigma_1}\right)^rr!\\
\le&\;C_1(C_1\sigma_1^{-1}\gamma)^1C_{1,0}
\left(\frac{\beta}{\sigma_1}\right)^{1+r-0}(1+r-0)!
\end{align*}
and this gives~\eqref{eq:DrAmsest} in this case.

\paragraph{\mathversion{bold}Case $m\in\integerp$ and $s=0$\@:}

We use induction on $m$\@, the desired estimate having been shown to be true
for $m=1$\@.  By Lemma~\ref{lem:leibniz} we have
\begin{multline*}
D^r_{\nabla^{\man{M}},\nabla^{\pi_{\man{E}}\otimes\pi_{\man{F}}}}A^{m+1}_0=
\underbrace{\sum_{a=0}^r\binom{r}{a}
D^a_{\nabla^{\man{M}},\nabla^{\pi_{\man{E}}}}\Phi^0_m
(D^{r-a}_{\nabla^{\man{M}},\nabla^{\pi_{\man{E}}\otimes\pi_{\man{F}}}}
\nabla^{\man{M},\pi_{\man{E}}\otimes\pi_{\man{F}}}A^m_0)}_{\textrm{term~1}}\\
+\underbrace{\sum_{a=0}^r\binom{r}{a}
D^a_{\nabla^{\man{M}},\nabla^{\pi_{\man{E}}\otimes\pi_{\man{F}}}}\Psi^0_{0m}
(D^{r-a}_{\nabla^{\man{M}},\nabla^{\pi_{\man{E}}\otimes\pi_{\man{F}}}}
A^m_0)}_{\textrm{term~2(a)}}\\
+\underbrace{\sum_{j=1}^m\sum_{a=0}^r\binom{r}{a}
D^a_{\nabla^{\man{M}},\nabla^{\pi_{\man{E}}\otimes\pi_{\man{F}}}}\Psi^0_{jm}
(D^{r-a}_{\nabla^{\man{M}},\nabla^{\pi_{\man{E}}\otimes\pi_{\man{F}}}}
A^m_0)}_{\textrm{term~2(b)}}.
\end{multline*}
For term~1, as above for the case $m=1$ and $s=0$\@, we have
\begin{equation*}
\dnorm{D^{r-a}_{\nabla^{\man{M}},\nabla^{\pi_{\man{E}}\otimes\pi_{\man{F}}}}
(\nabla^{\man{M},{\pi_{\man{E}}\otimes\pi_{\man{F}}}}
A^m_0)}_{\metric_{\man{M},\pi_{\man{E}}\otimes\pi_{\man{F}}}}\le
\dnorm{D^{r-a+1}_{\nabla^{\man{M}},\nabla^{\pi_{\man{E}}\otimes\pi_{\man{F}}}}
A^m_0}_{\metric_{\man{M},\pi_{\man{E}}\otimes\pi_{\man{F}}}}.
\end{equation*}
We now use Lemma~\ref{lem:A(S)<AS}\@, the hypotheses, the induction
hypotheses, and the preliminary observation~\ref{enum:(a+b)!/b!} above to
determine that
\begin{align*}
\dnorm{\textrm{term~1}(x)}_{\metric_{\man{M},\pi_{\man{E}}\otimes\pi_{\man{F}}}}
\le&\;\sum_{a=0}^r\frac{r!}{a!(r-a)!}(C_1\sigma_1^{-a}a!)\\
&\;\phantom{\sum_{a=0}^r}\times\left(C_1(C_1\sigma_1^{-1}\gamma)^mC_{m,0}
\left(\frac{\beta}{\sigma_1}\right)^{m+r-a+1}(m+r-a+1)!\right)\\
\le&\;C_1(C_1\sigma_1^{-1}\gamma)^mC_{m,0}
\left(\frac{\beta}{\sigma_1}\right)^{m+r+1}(m+r+1)!\sum_{a=0}^r\beta^{-a}\\
\le&\;C_1(C_1\sigma_1^{-1}\gamma)^mC_{m,0}\alpha
\left(\frac{\beta}{\sigma_1}\right)^{m+r+1}(m+r+1)!
\end{align*}
By a similar computation, we have
\begin{equation*}
\dnorm{\textrm{term~2(a)}(x)}_{\metric_{\man{M},\pi_{\man{E}}\otimes\pi_{\man{F}}}}
\le C_1(C_1\sigma_1^{-1}\gamma)^mC_{m,0}\alpha
\left(\frac{\beta}{\sigma_1}\right)^{m+r}(m+r)!.
\end{equation*}
We also have, making use of our observation~\ref{enum:Cms} from above,
\begin{align*}
\dnorm{\textrm{term~2(b)}(x)}_{\metric_{\man{M},\pi_{\man{E}}\otimes\pi_{\man{F}}}}
\le&\;C_1(C_1\sigma_1^{-1}\gamma)^mmC_{m,0}\alpha
\left(\frac{\beta}{\sigma_1}\right)^{m+r}(m+r)!\\
\le&\;C_1(C_1\sigma_1^{-1}\gamma)^m(m+1)C_{m+1,0}\alpha
\left(\frac{\beta}{\sigma_1}\right)^{m+r}(m+r)!\\
\le&\;C_1(C_1\sigma_1^{-1}\gamma)^mC_{m+1,0}\alpha
\left(\frac{\beta}{\sigma_1}\right)^{m+r}(m+r+1)!.
\end{align*}
Thus, for $x\in\nbhd{K}$\@,
\begin{align*}
\dnorm{D^r_{\nabla^{\man{M}},\nabla^{\pi_{\man{E}}\otimes\pi_{\man{F}}}}
A^{m+1}_0(x)}_{\metric_{\man{M},\pi_{\man{E}}\otimes\pi_{\man{F}}}}\le&\;
\frac{1}{4}C_1(C_1\sigma_1^{-1}\gamma)^mC_{m,0}\gamma
\left(\frac{\beta}{\sigma_1}\right)^{m+r+1}(m+r+1)!\\
&\;+\frac{1}{4}C_1(C_1\sigma_1^{-1}\gamma)^mC_{m,0}\gamma
\left(\frac{\beta}{\sigma_1}\right)^{m+r}(m+r)!\\
&\;+\frac{1}{4}C_1(C_1\sigma_1^{-1}\gamma)^mC_{m+1,0}\gamma
\left(\frac{\beta}{\sigma_1}\right)^{m+r}(m+r+1)!\\
\le&\;C_1(C_1\sigma_1^{-1}\gamma)^{m+1}C_{m+1,0}
\left(\frac{\beta}{\sigma_1}\right)^{m+1+r-0}
(m+1+r-0)!.
\end{align*}
This proves~\eqref{eq:DrAmsest} by induction in this case.

\paragraph{\mathversion{bold}Case $m\in\integerp$ and
$s\in\{1,\dots,m-1\}$\@:}

We use induction first on $m$ (the result having been proved for the case
$m=0$) and, for fixed $m$\@, by induction on $s$ (the result having been
proved for the case $s=0$).  By Lemma~\ref{lem:leibniz} we have
\begin{multline*}
D^r_{\nabla^{\man{M}},\nabla^{\pi_{\man{E}}\otimes\pi_{\man{F}}}}A^{m+1}_s=
\underbrace{\sum_{a=0}^r\binom{r}{a}
D^a_{\nabla^{\man{M}},\nabla^{\pi_{\man{E}}}}\Phi^s_m
(D^{r-a}_{\nabla^{\man{M}},\nabla^{\pi_{\man{E}}\otimes\pi_{\man{F}}}}
\nabla^{\man{M},\pi_{\man{E}}\otimes\pi_{\man{F}}}A^m_s)}_{\textrm{term~1}}\\
+\underbrace{\sum_{a=0}^r\binom{r}{a}
D^a_{\nabla^{\man{M}},\nabla^{\pi_{\man{E}}\otimes\pi_{\man{F}}}}\Psi^s_{0m}
(D^{r-a}_{\nabla^{\man{M}},\nabla^{\pi_{\man{E}}\otimes\pi_{\man{F}}}}
A^m_s)}_{\textrm{term~2(a)}}\\
+\underbrace{\sum_{j=1}^m\sum_{a=0}^r\binom{r}{a}
D^a_{\nabla^{\man{M}},\nabla^{\pi_{\man{E}}\otimes\pi_{\man{F}}}}\Psi^s_{jm}
(D^{r-a}_{\nabla^{\man{M}},\nabla^{\pi_{\man{E}}\otimes\pi_{\man{F}}}}
A^m_s)}_{\textrm{term~2(b)}}\\
+\underbrace{\sum_{a=0}^r\binom{r}{a}
D^a_{\nabla^{\man{M}},\nabla^{\pi_{\man{E}}\otimes\pi_{\man{F}}}}\Lambda^s_m
(D^{r-a}_{\nabla^{\man{M}},\nabla^{\pi_{\man{E}}\otimes\pi_{\man{F}}}}
A^m_{s-1})}_{\textrm{term~3}}.
\end{multline*}
We can argue just as in the preceding paragraph that, for $x\in\nbhd{K}$\@,
\begin{align*}
\dnorm{\textrm{term~1}(x)}_{\metric_{\man{M},\pi_{\man{E}}\otimes\pi_{\man{F}}}}
\le&\;C_1(C_1\sigma_1^{-1}\gamma)^mC_{m,s}
\alpha\left(\frac{\beta}{\sigma_1}\right)^{m+r+1-s}(m+r+1-s)!\\
\dnorm{\textrm{term~2(a)}(x)}_{\metric_{\man{M},\pi_{\man{E}}\otimes\pi_{\man{F}}}}
\le&\;C_1(C_1\sigma_1^{-1}\gamma)^mC_{m,s}\alpha
\left(\frac{\beta}{\sigma_1}\right)^{m+r-s}(m+r-s)!\\
\dnorm{\textrm{term~2(b)}(x)}_{\metric_{\man{M},\pi_{\man{E}}\otimes\pi_{\man{F}}}}
\le&\;C_1(C_1\sigma_1^{-1}\gamma)^mC_{m+1,s}\alpha
\left(\frac{\beta}{\sigma_1}\right)^{m+r-s}(m+r+1-s)!\\
\dnorm{\textrm{term~3}(x)}_{\metric_{\man{M},\pi_{\man{E}}\otimes\pi_{\man{F}}}}
\le&\;C_1(C_1\sigma_1^{-1}\gamma)^mC_{m,s-1}\alpha
\left(\frac{\beta}{\sigma_1}\right)^{m+r-s+1}(m+r-s+1)!.
\end{align*}
Adding these as in the previous case and using our
observation~\ref{enum:Cms} above, we have
\begin{equation*}
\dnorm{D^r_{\nabla^{\man{M}},\nabla^{\pi_{\man{E}}\otimes\pi_{\man{F}}}}
A^{m+1}_s(x)}_{\metric_{\man{M},\pi_{\man{E}}\otimes\pi_{\man{F}}}}\le
C_1(C_1\sigma_1^{-1}\gamma)^{m+1}C_{m+1,s}
\left(\frac{\beta}{\sigma_1}\right)^{m+1+r-s}
(m+1+r-s)!,
\end{equation*}
proving~\eqref{eq:DrAmsest} by induction in this case.

We now note that a standard binomial estimate via~\eqref{eq:multinomialest}
gives $C_{m,s}\le 2^m$\@.  The lemma now follows from~\eqref{eq:DrAmsest} by
taking
\begin{equation*}\eqqed
C=C_1,\enspace \sigma=2C_1\sigma_1^{-1}\gamma,
\enspace\rho=\frac{\beta}{\sigma_1}.
\end{equation*}
\end{proof}
\end{lemma}

We now apply the lemma to the recursion relations that we proved in
Lemmata~\ref{lem:PiterateddersI}\@, \ref{lem:PiterateddersII}\@,
\ref{lem:ViterateddersI}\@, \ref{lem:ViterateddersII}\@,
\ref{lem:HiterateddersI}\@, \ref{lem:HiterateddersII}\@,
\ref{lem:V*iterateddersI}\@, \ref{lem:V*iterateddersII}\@,
\ref{lem:LiterateddersI}\@, \ref{lem:LiterateddersII}\@,
\ref{lem:DiterateddersI}\@, \ref{lem:DiterateddersII}\@,
\ref{lem:CiterateddersI}\@, \ref{lem:CiterateddersII}\@,
\ref{lem:pbiterateddersI}\@, and~\ref{lem:pbiterateddersII}\@.  We first
provide the correspondence between the data from the preceding lemmata with
the data of Lemma~\ref{lem:pissynabla}\@.\label{page:pissy-translate}
\begin{compactenum}
\item Lemma~\ref{lem:PiterateddersI}\@: We have
\begin{compactenum}
\item $\man{M}=\man{E}$\@, $\man{E}=\man{F}=\real_{\man{E}}$\@,
\item $\phi_m(\beta_m)=\beta_m$\@, $\beta_m\in\tensor*[m]{\ctb{\man{M}}}$\@, $m\in\integernn$\@,
\item $\Phi^s_m(\alpha^{m+1}_s)=\alpha^{m+1}_s$\@,
$\alpha^{m+1}_s\in\tensor[s,m+1]{\ctb{\man{M}}\otimes\man{F}
\otimes\dual{\man{E}}})$\@, $m\in\integernn$\@, $s\in\{0,1,\dots,m+1\}$\@,
\item
$\Psi^s_{jm}(\alpha^m_s)(\beta_s)=-\alpha^m_s\otimes
\id_{\ctb{\man{M}}}(\Ins_j(\beta_s,B_{\man{E}}))$\@,
$\alpha^m_s\in\Hom(\tensor*[s]{\ctb{\man{M}}}\otimes\man{E};
\tensor*[m]{\ctb{\man{M}}}\otimes\man{F})$\@,
$\beta_s\in\tensor*[s]{\ctb{\man{M}}}\otimes\man{E}$\@, $m\in\integerp$\@,
$s\in\{1,\dots,m\}$\@, $j\in\{1,\dots,s\}$\@,
\item
$\Lambda^s_m(\alpha^m_{s-1})=\alpha^m_{s-1}\otimes\id_{\ctb{\man{M}}}$\@,
$\alpha^m_{s-1}\in\Hom(\tensor*[s-1]{\ctb{\man{M}}}\otimes\man{E};
\tensor*[m]{\ctb{\man{M}}}\otimes\man{F})$\@,
$m\in\integerp$\@, $s\in\{1,\dots,m\}$\@,
\item $\Psi^0_{jm}=0$\@, $m\in\integernn$\@, and
\item $\Lambda^0_m=0$\@, $m\in\integernn$\@.
\end{compactenum}

\item Lemma~\ref{lem:PiterateddersII}\@: We have
\begin{compactenum}
\item $\man{M}=\man{E}$\@, $\man{E}=\man{F}=\real_{\man{E}}$\@,
\item $\phi_m(\beta_m)=\beta_m$\@, $\beta_m\in\tensor*[m]{\ctb{\man{M}}}$\@, $m\in\integernn$\@,
\item $\Phi^s_m(\alpha^{m+1}_s)=\alpha^{m+1}_s$\@,
$\alpha^{m+1}_s\in\tensor[s,m+1]{\ctb{\man{M}}\otimes\man{F}
\otimes\dual{\man{E}}}$\@, $m\in\integernn$\@, $s\in\{0,1,\dots,m+1\}$\@,
\item
$\Psi^s_{jm}(\alpha^m_s)(\beta_s)=
\Ins_j(\alpha^m_s(\beta_s),B_{\man{E}})$\@,
$\alpha^m_s\in\Hom(\tensor*[s]{\ctb{\man{M}}}\otimes\man{E};
\tensor*[m]{\ctb{\man{M}}}\otimes\man{F})$\@,
$\beta_s\in\tensor*[s]{\ctb{\man{M}}}\otimes\man{E}$\@, $m\in\integerp$\@,
$s\in\{1,\dots,m\}$\@, $j\in\{1,\dots,m\}$\@, and
\item
$\Lambda^s_m(\alpha^m_{s-1})=\alpha^m_{s-1}\otimes\id_{\ctb{\man{M}}}$\@,
$\alpha^m_{s-1}\in\Hom(\tensor*[s-1]{\ctb{\man{M}}}\otimes\man{E};
\tensor*[m]{\ctb{\man{M}}}\otimes\man{F})$\@, $m\in\integerp$\@,
$s\in\{0,\dots,m\}$\@.
\end{compactenum}

\item Lemma~\ref{lem:ViterateddersI}\@: We have
\begin{equation*}
\man{M}=\man{E},\enspace\man{E}=\man{F}=\vb{\man{E}},
\end{equation*}
and all other data derived from Lemma~\ref{lem:ViterateddersI}\@, similarly
to the case of Lemma~\ref{lem:PiterateddersI}\@.

\item Lemma~\ref{lem:ViterateddersII}\@: We have
\begin{equation*}
\man{M}=\man{E},\enspace\man{E}=\man{F}=\vb{\man{E}},
\end{equation*}
and all other data derived from Lemma~\ref{lem:ViterateddersII}\@, similarly
to the case of Lemma~\ref{lem:PiterateddersII}\@.

\item Lemma~\ref{lem:HiterateddersI}\@: We have
\begin{equation*}
\man{M}=\man{E},\enspace\man{E}=\man{F}=\hb{\man{E}},
\end{equation*}
and all other data derived from Lemma~\ref{lem:HiterateddersI}\@, similarly
to the case of Lemma~\ref{lem:PiterateddersI}\@.

\item Lemma~\ref{lem:HiterateddersII}\@: We have
\begin{equation*}
\man{M}=\man{E},\enspace\man{E}=\man{F}=\hb{\man{E}},
\end{equation*}
and all other data derived from Lemma~\ref{lem:HiterateddersII}\@, similarly
to the case of Lemma~\ref{lem:PiterateddersII}\@.

\item Lemma~\ref{lem:V*iterateddersI}\@: We have
\begin{equation*}
\man{M}=\man{E},\enspace\man{E}=\man{F}=\cvb{\man{E}},
\end{equation*}
and all other data derived from Lemma~\ref{lem:V*iterateddersI}\@, similarly
to the case of Lemma~\ref{lem:PiterateddersI}\@.

\item Lemma~\ref{lem:V*iterateddersII}\@: We have
\begin{equation*}
\man{M}=\man{E},\enspace\man{E}=\man{F}=\cvb{\man{E}},
\end{equation*}
and all other data derived from Lemma~\ref{lem:V*iterateddersI}\@, similarly
to the case of Lemma~\ref{lem:PiterateddersII}\@.

\item Lemma~\ref{lem:LiterateddersI}\@: We have
\begin{equation*}
\man{M}=\man{E},\enspace\man{E}=\man{F}=\tensor[1,1]{\vb{\man{E}}},
\end{equation*}
and all other data derived from Lemma~\ref{lem:LiterateddersI}\@, similarly
to the case of Lemma~\ref{lem:PiterateddersI}\@.

\item Lemma~\ref{lem:LiterateddersII}\@: We have
\begin{equation*}
\man{M}=\man{E},\enspace\man{E}=\man{F}=\tensor[1,1]{\vb{\man{E}}},
\end{equation*}
and all other data derived from Lemma~\ref{lem:LiterateddersII}\@, similarly
to the case of Lemma~\ref{lem:PiterateddersII}\@.

\item Lemma~\ref{lem:DiterateddersI}\@: We have
\begin{compactenum}
\item $\man{M}=\man{E}$\@, $\man{E}=\real_{\man{E}}\oplus\cvb{\man{E}}$\@,
$\man{F}=\real_{\man{E}}\oplus\real_{\man{E}}$\@,

\item $\phi_m(\beta_m,\delta_m)=\beta_m$\@, $(\beta_m,\delta_m)\in\tensor*[m]{\ctb{\man{M}}}\otimes\man{E}$\@,
$m\in\integernn$\@,

\item
$\Phi^s_m(\alpha^{m+1}_s,\gamma^{m+1}_s)=(\alpha^{m+1}_m,\gamma^{m+1}_s)$\@,
$(\alpha^{m+1}_s,\gamma^{m+1}_s)\in\tensor[s,m+1]{\ctb{\man{M}}}\otimes
\man{F}\otimes\dual{\man{E}}$\@, $m\in\integernn$\@,
$s\in\{0,1,\dots,m-1\}$\@,

\item $\Phi^m_m(\alpha^{m+1}_m,\gamma^{m+1}_m)=(0,0)$\@,
$(\alpha^{m+1}_m,\gamma^{m+1}_m)\in\tensor[m,m+1]{\ctb{\man{M}}}\otimes
\man{F}\otimes\dual{\man{E}}$\@, $m\in\integernn$\@,

\item
$\Psi^s_{jm}(\alpha^m_s,\gamma^m_s)(\beta_s,\delta_s)=
(-\alpha^m_s\otimes\id_{\ctb{\man{M}}}(\Ins_j(\beta_s,B_{\man{E}})),$\\
$-\sum_{j=1}^{s+1}\gamma^m_s\otimes\id_{\ctb{\man{M}}}
(\Ins_j(\delta_s,B_{\man{E}})))$\@,
$(\alpha^m_s,\gamma^m_s)\in\tensor[s,m]{\ctb{\man{M}}}\otimes\man{F}
\otimes\dual{\man{E}}$\@,
$(\beta_s,\delta_s)\in\tensor*[s]{\ctb{\man{M}}}\otimes\man{E}$\@, $m\ge2$\@,
$s\in\{1,\dots,m-1\}$\@, $j\in\{1,\dots,s\}$\@,

\item
$\Psi^m_{jm}(\alpha^m_m,\gamma^m_m)(\beta_m,\delta_m)=
(-\Ins_j(\beta_m,B_{\man{E}}),\delta_m)$\@,
$(\alpha^m_m,\gamma^m_m)\in\tensor[m,m]{\ctb{\man{M}}}\otimes\man{F}
\otimes\dual{\man{E}}$\@,
$(\beta_m,\delta_m)\in\tensor*[m]{\ctb{\man{M}}}\otimes\man{E}$\@, $m\ge2$\@, $j\in\{1,\dots,m\}$\@,

\item
$\Psi^m_{00}(\alpha^m_0,\gamma^m_0)(\beta_0,\delta_0)=
(0,-\gamma^m_0\otimes\id_{\ctb{\man{M}}}(\Ins_1(\delta_0,B_{\man{E}})))$\@,
$(\alpha^m_0,\gamma^m_0)\in\tensor[m,0]{\ctb{\man{M}}}\otimes\man{F}
\otimes\dual{\man{E}}$\@, $(\beta_0,\delta_0)\in\man{E}$\@, $m\ge2$\@,

\item
$\Lambda^s_m(\alpha^m_s,\gamma^m_s)=(\alpha^m_{s-1}\otimes\id_{\ctb{\man{M}}},
\gamma^m_{s-1}\otimes\id_{\ctb{\man{M}}})$\@, $m\ge2$\@, $s\in\{1,\dots,m\}$\@.
\end{compactenum}

\item Lemma~\ref{lem:DiterateddersII}\@: We have
\begin{compactenum}
\item $\man{M}=\man{E}$\@, $\man{E}=\real_{\man{E}}\oplus\cvb{\man{E}}$\@,
$\man{F}=\real_{\man{E}}\oplus\real_{\man{E}}$\@,
\item $\phi_m(\beta_m,\delta_m)=\beta_m$\@, $(\beta_m,\delta_m)\in\tensor*[m]{\ctb{\man{M}}}\otimes\man{E}$\@,
$m\in\integernn$\@,

\item
$\Phi^s_m(\alpha^{m+1}_s,\gamma^{m+1}_s)=(\alpha^{m+1}_m,\gamma^{m+1}_s)$\@,
$(\alpha^{m+1}_s,\gamma^{m+1}_s)\in\tensor[s,m+1]{\ctb{\man{M}}}\otimes
\man{F}\otimes\dual{\man{E}}$\@, $m\in\integernn$\@,
$s\in\{0,1,\dots,m-1\}$\@,

\item $\Phi^m_m(\alpha^{m+1}_m,\gamma^{m+1}_m)=(0,0)$\@,
$(\alpha^{m+1}_m,\gamma^{m+1}_m)\in\tensor[m,m+1]{\ctb{\man{M}}}\otimes
\man{F}\otimes\dual{\man{E}}$\@, $m\in\integernn$\@,

\item
$\Psi^s_{jm}(\alpha^m_s,\gamma^m_s)(\beta_s,\delta_s)=
(\Ins_j(\alpha^m_s(\beta_s),B_{\man{E}})-
\Ins_{m+1}(\alpha^m_s(\beta_s),\dual{B_{\man{E}}}),-\ol{B}^m_s)$\@,\\
$(\alpha^m_s,\gamma^m_s)\in\tensor[s,m]{\ctb{\man{M}}}\otimes\man{F}
\otimes\dual{\man{E}}$\@,
$(\beta_s,\delta_s)\in\tensor*[s]{\ctb{\man{M}}}\otimes\man{E}$\@, $m\ge2$\@,
$s\in\{1,\dots,m-1\}$\@, $j\in\{1,\dots,m\}$\@,

\item
$\Psi^m_{jm}(\alpha^m_m,\gamma^m_m)(\beta_m,\delta_m)=
(\Ins_j(\beta_m,B_{\man{E}})-
\Ins_{m+1}(\beta_m,\dual{B_{\man{E}}}),-\delta_m)$\@,\\
$(\alpha^m_m,\gamma^m_m)\in\tensor[m,m]{\ctb{\man{M}}}\otimes\man{F}
\otimes\dual{\man{E}}$\@,
$(\beta_m,\delta_m)\in\tensor*[m]{\ctb{\man{M}}}\otimes\man{E}$\@, $m\ge2$\@, $j\in\{1,\dots,m\}$\@,

\item
$\Psi^m_{00}(\alpha^m_0,\gamma^m_0)(\beta_0,\delta_0)=
(0,-\gamma^m_0\otimes\id_{\ctb{\man{M}}}(\Ins_1(\delta_0,B_{\man{E}})))$\@,
$(\alpha^m_0,\gamma^m_0)\in\tensor[m,0]{\ctb{\man{M}}}\otimes\man{F}
\otimes\dual{\man{E}}$\@, $(\beta_0,\delta_0)\in\man{E}$\@, $m\ge2$\@,

\item
$\Lambda^s_m(\alpha^m_s,\gamma^m_s)=(\alpha^m_{s-1}\otimes\id_{\ctb{\man{M}}},
\gamma^m_{s-1}\otimes\id_{\ctb{\man{M}}})$\@, $m\ge2$\@, $s\in\{1,\dots,m\}$\@.
\end{compactenum}

\item Lemma~\ref{lem:CiterateddersI}\@: We have
\begin{equation*}
\man{M}=\man{E},\enspace\man{E}=\vb{\man{E}}\oplus\tensor[1,1]{\vb{\man{E}}},
\enspace\man{F}=\vb{\man{E}}\oplus\vb{\man{E}},
\end{equation*}
and all other data derived from Lemma~\ref{lem:CiterateddersI}\@, similarly
to the case of Lemma~\ref{lem:DiterateddersI}\@.

\item Lemma~\ref{lem:CiterateddersII}\@: We have
\begin{equation*}
\man{M}=\man{E},\enspace\man{E}=\vb{\man{E}}\oplus\tensor[1,1]{\vb{\man{E}}},
\enspace\man{F}=\vb{\man{E}}\oplus\vb{\man{E}},
\end{equation*}
and all other data derived from Lemma~\ref{lem:CiterateddersII}\@, similarly
to the case of Lemma\ref{lem:DiterateddersII}\@.

\item Lemma~\ref{lem:pbiterateddersI}\@: We have
\begin{equation*}
\man{M}=\man{M},\enspace\man{E}=\Phi^*\ctb{\man{N}},\enspace
\man{F}=\ctb{\man{M}},
\end{equation*}
and all other data derived from Lemma~\ref{lem:pbiterateddersI}\@, similarly
to the case of Lemma~\ref{lem:PiterateddersI}\@.

\item Lemma~\ref{lem:pbiterateddersII}\@: We have
\begin{equation*}
\man{M}=\man{M},\enspace\man{E}=\ctb{\man{M}},\enspace
\man{F}=\Phi^*\ctb{\man{N}},
\end{equation*}
and all other data derived from Lemma~\ref{lem:pbiterateddersII}\@, similarly
to the case of Lemma~\ref{lem:PiterateddersII}\@.
\end{compactenum}

Having now translated the lemmata of Section~\ref{sec:lift-isomorphisms} to
the general Lemma~\ref{lem:pissynabla}\@, we now need to show that the data
of the lemmata of Section~\ref{sec:lift-isomorphisms} satisfy the hypotheses
of Lemma~\ref{lem:pissynabla}\@.  As is easily seen, there are a few sorts of
expressions that appear repeatedly, and we shall simply give estimates for
these terms and leave to the reader the putting together of the pieces.

The following lemma gives the required bounds.
\begin{lemma}\label{lem:particular-bounds}
Let\/ $\map{\pi_{\man{E}}}{\man{E}}{\man{M}}$ be a real analytic vector
bundle, let\/ $\nabla^{\man{M}}$ be a real analytic affine connection on\/
$\man{M}$\@, and let\/ $\nabla^{\pi_{\man{E}}}$ be a real analytic vector
bundle connection in\/ $\man{E}$\@.  Let\/ $\metric_{\man{M}}$ be a real
analytic Riemannian metric on\/ $\man{M}$ and let\/ $\metric_{\man{E}}$ be a
real analytic fibre metrics for\/ $\man{E}$\@.  Let\/
$S\in\sections[\omega]{\tensor[1,2]{\tb{\man{M}}}}$\@.  Let\/
$\nbhd{K}\subset\man{M}$ be compact and let\/ $n$ be the larger of the
dimension of\/ $\man{M}$ and the fibre dimension of\/ $\man{E}$ and let\/
$\sigma_0=n^{-1}$\@.  Let\/ $m,r,a\in\integernn$ and\/
$s\in\{0,1,\dots,m\}$\@.  Then we have the following bounds for\/
$x\in\nbhd{K}$\@:
\begin{compactenum}[(i)]
\item \label{pl:particular-bounds1}
$\dnorm{D^r_{\nabla^{\man{M}},\nabla^{\pi_{\man{E}}}}
\id_{\tensor*[m]{\ctb{\man{M}}}\otimes\man{E}}(x)}_{\metric_{\man{M},\pi_{\man{E}}}}
\le\sigma_0^{m+r+1}$\@;

\item \label{pl:particular-bounds2}
$\dnorm{D^r_{\nabla^{\man{M}},\nabla^{\pi_{\man{E}}}}
\id_{\tensor[m,s]{\ctb{\man{M}}}\otimes\man{E}}
(x)}_{\metric_{\man{M},\pi_{\man{E}}}}\le\sigma_0^{2m+r+1}$\@;

\item \label{pl:particular-bounds3} if\/
$\Phi^s_m(\alpha^{m+1}_s)=\alpha^{m+1}_s$\@,\/
$\alpha^{m+1}_s\in\tensor[s,m+1]{\ctb{\man{M}}}\otimes\man{E}$\@, then there
exist\/ $C_1,\sigma_1\in\realp$ such that
\begin{equation*}
\dnorm{D^r_{\nabla^{\man{M}},\nabla^{\pi_{\man{E}}}}\Phi^s_m\scirc
D^a_{\nabla^{\man{M}},\nabla^{\pi_{\man{E}}}}
A^{m+1}_s(x)}_{\metric_{\man{M},\pi_{\man{E}}}}\le
\dnorm{D^a_{\nabla^{\man{M}},\nabla^{\pi_{\man{E}}}}
A^{m+1}_s(x)}_{\metric_{\man{M},\pi_{\man{E}}}};
\end{equation*}

\item \label{pl:particular-bounds4} if
\begin{multline*}
\Psi^s_{jm}(\alpha^m_s)(\beta_s)=(\alpha^m_s\otimes
\id_{\ctb{\man{M}}})(\Ins_j(\beta_s,S)),\\
\alpha^m_s\in\Hom(\tensor*[s]{\ctb{\man{M}}}\otimes\man{E};
\tensor*[m]{\ctb{\man{M}}}\otimes\man{E}),\
\beta_s\in\tensor*[s]{\ctb{\man{M}}}\otimes\man{E},
\end{multline*}
then there exist\/ $C_1,\sigma_1\in\realp$ such that
\begin{equation*}
\dnorm{D^r_{\nabla^{\man{M}},\nabla^{\pi_{\man{E}}}}\Psi^s_{jm}\scirc
D^a_{\nabla^{\man{M}},\nabla^{\pi_{\man{E}}}}
A^m_s(x)}_{\metric_{\man{M},\pi_{\man{E}}}}\le C_1\sigma_1^{-r}r!
\dnorm{D^a_{\nabla^{\man{M}},\nabla^{\pi_{\man{E}}}}
A^m_s(x)}_{\metric_{\man{M},\pi_{\man{E}}}};
\end{equation*}

\item \label{pl:particular-bounds5} if
\begin{multline*}
\Psi^s_{jm}(\alpha^m_s)(\beta_s)=\Ins_j(\alpha^m_s(\beta_s),S),\\
\alpha^m_s\in\Hom(\tensor*[s]{\ctb{\man{M}}}\otimes\man{E};
\tensor*[m]{\ctb{\man{M}}}\otimes\man{E}),\
\beta_s\in\tensor*[s]{\ctb{\man{M}}}\otimes\man{E},
\end{multline*}
then there exist\/ $C_1,\sigma_1\in\realp$ such that
\begin{equation*}
\dnorm{D^r_{\nabla^{\man{M}},\nabla^{\pi_{\man{E}}}}\Psi^s_{jm}\scirc
D^a_{\nabla^{\man{M}},\nabla^{\pi_{\man{E}}}}
A^m_s(x)}_{\metric_{\man{M},\pi_{\man{E}}}}\le C_1\sigma_1^{-r}r!
\dnorm{D^a_{\nabla^{\man{M}},\nabla^{\pi_{\man{E}}}}
A^m_s(x)}_{\metric_{\man{M},\pi_{\man{E}}}};
\end{equation*}

\item \label{pl:particular-bounds6} if
\begin{equation*}
\Lambda^s_m(\alpha^m_{s-1})=\alpha^m_{s-1}\otimes\id_{\ctb{\man{M}}},\qquad
\alpha^m_s\in\Hom(\tensor*[s-1]{\ctb{\man{M}}}\otimes\man{E};
\tensor*[m]{\ctb{\man{M}}}\otimes\man{E}),
\end{equation*}
then there exist\/ $C_1,\sigma_1\in\realp$ such that
\begin{equation*}
\dnorm{D^r_{\nabla^{\man{M}},\nabla^{\pi_{\man{E}}}}\Lambda^s_m\scirc
D^a_{\nabla^{\man{M}},\nabla^{\pi_{\man{E}}}}
A^m_{s-1}(x)}_{\metric_{\man{M},\pi_{\man{E}}}}\le 
\dnorm{D^a_{\nabla^{\man{M}},\nabla^{\pi_{\man{E}}}}
A^m_{s-1}(x)}_{\metric_{\man{M},\pi_{\man{E}}}}.
\end{equation*}
\end{compactenum}
\begin{proof}
Parts~\eqref{pl:particular-bounds1} and~\eqref{pl:particular-bounds2} follow
from Lemma~\ref{lem:idnorm} along with the fact that the covariant derivative
of the identity tensor is zero.  Part~\eqref{pl:particular-bounds3} is a
tautology, but one that arises in the lemmata of
Section~\ref{sec:lift-isomorphisms}\@.

For the next two parts of the proof, let $C_1,\sigma_1\in\realp$ be such that
\begin{equation}\label{eq:Sbound}
\dnorm{D^r_{\nabla^{\man{M}}}S(x)}_{\metric_{\man{M}}}\le C_1\sigma_1^{-r}r!,
\qquad x\in\nbhd{K},
\end{equation}
this being possible by Lemma~\ref{lem:Comegachar}\@, and recalling the r\^ole
of the factorials in the definition~\eqref{eq:olGpim} of the fibre norms.

\eqref{pl:particular-bounds4} Let us define
\begin{multline*}
\hat{\Psi}^s_{jm}(\beta^m_{s+1})(\alpha_s)=
(\beta^m_{s+1})(\Ins_j(\alpha_s,S)),\\
\beta^m_{s+1}\in\Hom(\tensor*[s+1]{\ctb{\man{M}}}\otimes\man{E};
\tensor*[m]{\ctb{\man{M}}}\otimes\man{E}),\
\alpha_s\in\tensor*[s]{\ctb{\man{M}}}\otimes\man{E}
\end{multline*}
and
\begin{equation*}
\tau^s_m(\alpha^m_s)=\alpha^m_s\otimes\id_{\ctb{\man{M}}},\qquad
\alpha^m_s\in\Hom(\tensor*[s]{\ctb{\man{M}}}\otimes\man{E};
\tensor*[m]{\ctb{\man{M}}}\otimes\man{E})
\end{equation*}
so that $\Psi^s_{jm}=\hat{\Psi}^s_m\scirc\tau^s_m$\@.  Note that
$\hat{\Psi}^s_{jm}=\Ins_{S,j}$ so that, by Lemma~\ref{lem:insder2}\@,
\begin{equation*}
D^r_{\nabla^{\man{M}},\nabla^{\pi_{\man{E}}}}\hat{\Psi}^s_{jm}
(D^a_{\nabla^{\man{M}},\nabla^{\pi_{\man{E}}}}B^m_{s+1})=
\Ins_{D^r_{\nabla^{\man{M}},\nabla^{\pi_{\man{E}}}}S,j}
(D^a_{\nabla^{\man{M}},\nabla^{\pi_{\man{E}}}}B^m_{s+1}).
\end{equation*}
Since the covariant derivative of the identity tensor is zero,
\begin{equation*}
D^a_{\nabla^{\man{M}},\nabla^{\pi_{\man{E}}}}(A^m_s\otimes\id_{\ctb{\man{M}}})=
(D^a_{\nabla^{\man{M}},\nabla^{\pi_{\man{E}}}}A^m_s)\otimes\id_{\ctb{\man{M}}}),
\end{equation*}
from which we deduce that
$D^a_{\nabla^{\man{M}},\nabla^{\pi_{\man{E}}}}\tau^s_m=\tau^s_{m+a}$\@.  Thus
\begin{align*}
D^r_{\nabla^{\man{M}},\nabla^{\pi_{\man{E}}}}\Psi^s_{jm}\scirc
D^a_{\nabla^{\man{M}},\nabla^{\pi_{\man{E}}}}A^m_s=&\;
D^r_{\nabla^{\man{M}},\nabla^{\pi_{\man{E}}}}
(\hat{\Psi}^s_{jm}\scirc\tau^s_m)\scirc
D^a_{\nabla^{\man{M}},\nabla^{\pi_{\man{E}}}}A^m_s\\
=&\;\Ins_{D^r_{\nabla^{\man{M}},\nabla^{\pi_{\man{E}}}}S,j}
(D^a_{\nabla^{\man{M}},\nabla^{\pi_{\man{E}}}}A^m_s\otimes\id_{\ctb{\man{M}}}).
\end{align*}
By Lemmata~\ref{lem:I1AS} and~\ref{lem:otimesnorm}\@, this part of the
lemma follows immediately.

\eqref{pl:particular-bounds5} Here we have
$\Psi^s_m(\alpha^m_s)=\Ins_{S,j}\scirc\alpha^m_s$ and, following the
arguments from the preceding part of the proof,
\begin{equation*}
D^r_{\nabla^{\man{M}},\nabla^{\pi_{\man{E}}}}\Psi^s_m\scirc
D^a_{\nabla^{\man{M}},\nabla^{\pi_{\man{E}}}}A^m_s=
\Ins_{D^r_{\nabla^{\man{M}},\nabla^{\pi_{\man{E}}}}S,j}
(D^a_{\nabla^{\man{M}},\nabla^{\pi_{\man{E}}}}A^m_s),
\end{equation*}
and so this part of the lemma follows from Lemma~\ref{lem:I2BS}\@.

\eqref{pl:particular-bounds6} This follows from Lemma~\ref{lem:otimesnorm}
and the fact that the covariant derivative of the identity tensor is zero.
\end{proof}
\end{lemma}

\section{Independence of topologies on connections and
metrics}\label{sec:independent}

The seminorms introduced in Section~\ref{subsec:seminorms} for defining
topologies for the space of real analytic sections of a vector bundle
$\map{\pi_{\man{E}}}{\man{E}}{\man{M}}$ are made upon a choice of various
objects, namely~(1)~an affine connection $\nabla^{\man{M}}$ on $\man{M}$\@,
(2)~a vector bundle connection $\nabla^{\pi_{\man{E}}}$ in $\man{E}$\@, (3)~a
Riemannian metric $\metric_{\man{M}}$ on $\man{M}$\@, and (4)~a fibre metric
$\metric_{\pi_{\man{E}}}$ for $\man{E}$\@.  In order for these topologies to
be useful, they should be independent of all of these choices.  This is made
more urgent by our very specific choice in Section~\ref{subsec:Esubmersion}
of a Riemannian metric $\metric_{\man{E}}$ on the total space $\man{E}$ and
its Levi-Civita connection.  These choices were made because they made
available to us the geometric constructions of
Section~\ref{sec:tensor-derivatives}\@, constructions of which substantial
use was made in Sections~\ref{sec:lift-isomorphisms}
and~\ref{sec:jet-estimates}\@, and of which will be made in
Section~\ref{sec:continuity}\@, as well as in the present section.

That the topologies are independent of choices of geometric objects is more
or less clear in the smooth case, but we will rather precisely point out why
this is so in our developments below.  In the real analytic case, one must
make use of all of the technical developments of Sections~\ref{sec:tensor-constructions}--\ref{sec:jet-estimates}\@.

\subsection{Comparison of iterated covariant derivatives for different
connections}

Our constructions start by comparing how covariant derivatives of high-order
differ when one changes connection.  The reader will see substantial
similarity between the results in this section and those in
Sections~\ref{subsec:formsvf} and~\ref{sec:lift-isomorphisms}\@.

We let $r\in\{\infty,\omega\}$ and let $\map{\pi_{\man{E}}}{\man{E}}{\man{M}}$
be a $\C^r$-vector bundle.  We consider $\C^r$-affine connections
$\nabla^{\man{M}}$ and $\ol{\nabla}^{\man{M}}$ on $\man{M}$\@, and vector
bundle connections $\nabla^{\pi_{\man{E}}}$ and $\ol{\nabla}^{\pi_{\man{E}}}$
in $\man{E}$\@.  It then holds that
\begin{equation*}
\ol{\nabla}^{\man{M}}_XY=\nabla^{\man{M}}_XY+S_{\man{M}}(Y,X),\quad
\ol{\nabla}^{\pi_{\man{E}}}_X\xi=\nabla^{\pi_{\man{E}}}_X\xi+
S_{\pi_{\man{E}}}(\xi,X)
\end{equation*}
for $S_{\man{M}}\in\sections[r]{\tensor[1,2]{\tb{\man{M}}}}$ and
$S_{\pi_{\man{E}}}\in\sections[r]{\dual{\man{E}}\otimes
\ctb{\man{M}}\otimes\man{E}}$\@.

First we relate covariant derivatives of higher-order tensors.
\begin{lemma}\label{lem:Tkdiffnabla}
Let\/ $r\in\{\infty,\omega\}$ and let\/
$\map{\pi_{\man{E}}}{\man{E}}{\man{M}}$ be a\/ $\C^r$-vector bundle.
Consider\/ $\C^r$-affine connections\/ $\nabla^{\man{M}}$ and\/
$\ol{\nabla}^{\man{M}}$ on\/ $\man{M}$\@, and\/ $\C^r$-vector bundle
connections\/ $\nabla^{\pi_{\man{E}}}$ and\/ $\ol{\nabla}^{\pi_{\man{E}}}$
in\/ $\man{E}$\@.  If\/ $k\in\integerp$ and if\/
$B\in\sections[1]{\tensor*[k]{\ctb{\man{M}}}\otimes\man{E}}$\@, then
\begin{equation*}
\ol{\nabla}^{\man{M},\pi_{\man{E}}}B=
\nabla^{\man{M},\pi_{\man{E}}}B-\sum_{j=1}^k\Ins_j(B,S_{\man{M}})-
\Ins_{k+1}(B,S_{\pi_{\man{E}}}).
\end{equation*}
\begin{proof}
We have
{\allowdisplaybreaks\begin{align*}
\lieder{X_{k+1}}{}&(B(X_1,\dots,X_k,\alpha))=
(\ol{\nabla}^{\man{M},\pi_{\man{E}}}_{X_{k+1}}B)(X_1,\dots,X_k,\alpha)\\
&\;+\sum_{j=1}^kB(X_1,\dots,\ol{\nabla}^{\man{M}}_{X_{k+1}}X_j,\dots,X_k,\alpha)
+B(X_1,\dots,X_k,\ol{\nabla}^{\pi_{\man{E}}}_{X_{k+1}}\alpha)\\
=&\;(\ol{\nabla}^{\man{M},\pi_{\man{E}}}_{X_{k+1}}B)(X_1,\dots,X_k,\alpha)+
\sum_{j=1}^kB(X_1,\dots,\nabla^{\man{M}}_{X_{k+1}}X_j,\dots,X_k,\alpha)\\
&\;+\sum_{j=1}^kB(X_1,\dots,S_{\man{M}}(X_j,X_{k+1}),\dots,X_k,\alpha)+
B(X_1,\dots,X_k,\nabla^{\pi_{\man{E}}}_{X_{k+1}}\alpha)\\
&\;+B(X_1,\dots,X_k,S_{\pi_{\man{E}}}(\alpha,X_{k+1})).
\end{align*}}%
This gives
\begin{equation*}
\ol{\nabla}^{\man{M},\pi_{\man{E}}}B=
\nabla^{\man{M},\pi_{\man{E}}}B-\sum_{j=1}^k\Ins_j(B,S_{\man{M}})-
\Ins_{k+1}(B,S_{\pi_{\man{E}}}),
\end{equation*}
as claimed.
\end{proof}
\end{lemma}

With this lemma, we can provide the following characterisation of iterated
covariant differentials of sections of $\man{E}$ with respect to different
connections.
\begin{lemma}\label{lem:nablaonablaiterI}
Let\/ $r\in\{\infty,\omega\}$ and let\/
$\map{\pi_{\man{E}}}{\man{E}}{\man{M}}$ be a\/ $\C^r$-vector bundle.
Consider\/ $\C^r$-affine connections\/ $\nabla^{\man{M}}$ and\/
$\ol{\nabla}^{\man{M}}$ on\/ $\man{M}$\@, and\/ $\C^r$-vector bundle
connections\/ $\nabla^{\pi_{\man{E}}}$ and\/ $\ol{\nabla}^{\pi_{\man{E}}}$
in\/ $\man{E}$\@.  For\/ $m\in\integernn$\@, there exist\/ $\C^r$-vector
bundle mappings
\begin{equation*}
(A^m_s,\id_{\man{E}})\in
\vbmappings[r]{\tensor*[s]{\ctb{\man{M}}\otimes\man{E}}}
{\tensor*[m]{\ctb{\man{M}}}\otimes\man{E}},\qquad s\in\{0,1,\dots,m\},
\end{equation*}
such that
\begin{equation*}
\ol{\nabla}^{\man{M},\pi_{\man{E}},m}\xi=
\sum_{s=0}^mA^m_s(\nabla^{\man{M},\pi_{\man{E}},s}\xi)
\end{equation*}
for all\/ $\xi\in\sections[m]{\man{E}}$\@.  Moreover, the vector bundle
mappings\/ $A^m_0,A^m_1,\dots,A^m_m$ satisfy the recursion relations
prescribed by\/ $A^0_0(\beta_0)=\beta_0$ and
\begin{align*}
A^{m+1}_{m+1}(\beta_{m+1})=&\;\beta_{m+1},\\
A^{m+1}_s(\beta_s)=&\;(\ol{\nabla}^{\man{M},\pi_{\man{E}}}A^m_s)(\beta_s)
+A^m_{s-1}\otimes\id_{\ctb{\man{M}}}(\beta_s)
-\sum_{j=1}^sA^m_s\otimes\id_{\ctb{\man{M}}}
(\Ins_j(\beta_s,S_{\man{M}}))\\
&\;-A^m_s\otimes\id_{\ctb{\man{M}}}(\Ins_{s+1}(\beta_s,S_{\pi_{\man{E}}})),
\enspace s\in\{1,\dots,m\},\\
A^{m+1}_0(\beta_0)=&\;(\ol{\nabla}^{\man{M},\pi_{\man{E}}}A^m_0)(\beta_0)-
A^m_0\otimes\id_{\ctb{\man{M}}}(\Ins_1(\beta_0,S_{\pi_{\man{E}}})),
\end{align*}
where\/ $\beta_s\in\tensor*[s]{\ctb{\man{M}}}\otimes\man{E}$\@,\/
$s\in\{0,1,\dots,m\}$\@.
\begin{proof}
The assertion clearly holds for $m=0$\@, so suppose it true for
$m\in\integerp$\@.  Thus
\begin{equation*}
\ol{\nabla}^{\man{M},\pi_{\man{E}},m}\xi=
\sum_{s=0}^mA^m_s(\nabla^{\man{M},\pi_{\man{E}},s}\xi),
\end{equation*}
where the vector bundle mappings $A^a_s$\@, $a\in\{0,1,\dots,m\}$\@,
$s\in\{0,1,\dots,a\}$\@, satisfy the recursion relations from the statement
of the lemma.  Then
\begin{align*}
\ol{\nabla}^{\man{M},\pi_{\man{E}},m+1}&\xi=
\sum_{s=0}^m(\ol{\nabla}^{\man{M},\pi_{\man{E}}}A^m_s)
(\nabla^{\man{M},\pi_{\man{E}},s}\xi)
+\sum_{s=0}^mA^m_s\otimes\id_{\ctb{\man{M}}}
(\ol{\nabla}^{\man{M},\pi_{\man{E}}}\nabla^{\man{M},\pi_{\man{E}},s}\xi)\\
=&\;\sum_{s=0}^m(\ol{\nabla}^{\man{M},\pi_{\man{E}}}A^m_s)
(\nabla^{\man{M},\pi_{\man{E}},s}\xi)+
\sum_{s=0}^mA^m_s\otimes\id_{\ctb{\man{M}}}
(\nabla^{\man{M},\pi_{\man{E}},s+1}\xi)\\
&\;-\sum_{s=1}^m\sum_{j=1}^sA^m_s\otimes\id_{\ctb{\man{M}}}
(\Ins_j(\nabla^{\man{M},\pi_{\man{E}},s}\xi,S_{\man{M}}))\\
&\;-\sum_{s=1}^mA^m_s\otimes\id_{\ctb{\man{M}}}
(\Ins_{s+1}(\nabla^{\man{M},\pi_{\man{E}},s}\xi,S_{\pi_{\man{E}}}))
-A^m_0\otimes\id_{\ctb{\man{M}}}(\Ins_1(\xi,S_{\pi_{\man{E}}}))\\
=&\;\nabla^{\man{M},\pi_{\man{E}},m+1}\xi+
\left.\sum_{s=1}^m\right((\ol{\nabla}^{\man{M},\pi_{\man{E}}}A^m_s)
(\nabla^{\man{M},\pi_{\man{E}},s}\xi)
+A^m_{s-1}\otimes\id_{\ctb{\man{M}}}(\nabla^{\man{M},\pi_{\man{E}},s}\xi)\\
&\;-\sum_{j=1}^sA^m_s\otimes\id_{\ctb{\man{M}}}
(\Ins_j(\nabla^{\man{M},\pi_{\man{E}},s}\xi,S_{\man{M}}))
-\left.\vphantom{\sum_{j=1}^s}A^m_s\otimes\id_{\ctb{\man{M}}}
(\Ins_{s+1}(\nabla^{\man{M},\pi_{\man{E}},s}\xi,S_{\pi_{\man{E}}}))\right)\\
&\;-(\ol{\nabla}^{\man{M},\pi_{\man{E}}}A^m_0)(\xi)-
A^m_0\otimes\id_{\ctb{\man{M}}}(\Ins_1(\xi,S_{\pi_{\man{E}}}))
\end{align*}
by Lemma~\ref{lem:Tkdiffnabla}\@.  From this, the lemma follows.
\end{proof}
\end{lemma}

The lemma has an ``inverse'' which we state next.
\begin{lemma}\label{lem:nablaonablaiterII}
Let\/ $r\in\{\infty,\omega\}$ and let\/
$\map{\pi_{\man{E}}}{\man{E}}{\man{M}}$ be a\/ $\C^r$-vector bundle.
Consider\/ $\C^r$-affine connections\/ $\nabla^{\man{M}}$ and\/
$\ol{\nabla}^{\man{M}}$ on\/ $\man{M}$\@, and\/ $\C^r$-vector bundle
connections\/ $\nabla^{\pi_{\man{E}}}$ and\/ $\ol{\nabla}^{\pi_{\man{E}}}$
in\/ $\man{E}$\@.  For\/ $m\in\integernn$\@, there exist\/ $\C^r$-vector
bundle mappings
\begin{equation*}
(B^m_s,\id_{\man{E}})\in
\vbmappings[r]{\tensor*[s]{\ctb{\man{M}}\otimes\man{E}}}
{\tensor*[m]{\ctb{\man{M}}}\otimes\man{E}},\qquad s\in\{0,1,\dots,m\},
\end{equation*}
such that
\begin{equation*}
\nabla^{\man{M},\pi_{\man{E}},m}\xi=
\sum_{s=0}^mB^m_s(\ol{\nabla}^{\man{M},\pi_{\man{E}},s}\xi)
\end{equation*}
for all\/ $\xi\in\sections[m]{\man{E}}$\@.  Moreover, the vector bundle
mappings\/ $B^m_0,B^m_1,\dots,B^m_m$ satisfy the recursion relations
prescribed by\/ $B^0_0(\alpha_0)=\beta_0$ and
\begin{align*}
B^{m+1}_{m+1}(\alpha_{m+1})=&\;\alpha_{m+1},\\
B^{m+1}_s(\alpha_s)=&\;(\ol{\nabla}^{\man{M},\pi_{\man{E}}}B^m_s)(\alpha_s)+
B^m_{s-1}\otimes\id_{\ctb{\man{M}}}(\alpha_s)+
\sum_{j=1}^m\Ins_j(B^m_s(\alpha_s),S_{\man{M}})\\
&\;+\Ins_{m+1}(B^m_s(\alpha_s\xi),S_{\pi_{\man{E}}}),\enspace s\in\{1,\dots,m\},\\
B^{m+1}_0(\alpha_0)=&\;(\ol{\nabla}^{\man{M},\pi_{\man{E}}}B^m_0)(\alpha_0)+
\sum_{j=1}^m\Ins_j(B^m_0(\alpha_0),S_{\man{M}})+
\Ins_{m+1}(B^m_0(\alpha_0),S_{\pi_{\man{E}}}),
\end{align*}
where\/ $\alpha_s\in\tensor*[s]{\ctb{\man{M}}}\otimes\man{E}$\@,\/
$s\in\{0,1,\dots,m\}$\@.
\begin{proof}
The lemma is clearly true for $m=0$\@, so suppose it true for
$m\in\integerp$\@.  Thus
\begin{equation}\label{eq:induct-nabla-onabla}
\nabla^{\man{M},\pi_{\man{E}},m}\xi=
\sum_{s=0}^mB^m_s(\ol{\nabla}^{\man{M},\pi_{\man{E}},s}\xi),
\end{equation}
where the vector bundle mappings $B^a_s$\@, $a\in\{0,1,\dots,m\}$\@,
$s\in\{0,1,\dots,a\}$\@, satisfy the recursion relations given in the lemma.
Then, working with the left-hand side of this relation,
\begin{align*}
\ol{\nabla}^{\man{M},\pi_{\man{E}}}\nabla^{\man{M},\pi_{\man{E}},m}\xi=&\;
\nabla^{\man{M},\pi_{\man{E}},m+1}\xi-
\sum_{j=1}^m\Ins_j(\nabla^{\man{M},\pi_{\man{E}},m}\xi,S_{\man{M}})-
\Ins_{m+1}(\nabla^{\man{M},\pi_{\man{E}},m}\xi,S_{\pi_{\man{E}}})\\
=&\;\nabla^{\man{M},\pi_{\man{E}},m+1}\xi-\sum_{s=0}^m\sum_{j=1}^m
\Ins_j(B^m_s(\ol{\nabla}^{\man{M},\pi_{\man{E}},s}\xi),S_{\man{M}})\\
&\;-\sum_{s=0}^m\Ins_{m+1}(B^m_s(\ol{\nabla}^{\man{M},\pi_{\man{E}},s}\xi),
S_{\pi_{\man{E}}}),
\end{align*}
by Lemma~\ref{lem:Tkdiffnabla}\@.  Now, working with the right-hand side
of~\eqref{eq:induct-nabla-onabla}\@,
\begin{equation*}
\ol{\nabla}^{\man{M},\pi_{\man{E}}}\nabla^{\man{M},\pi_{\man{E}},m}\xi=
\sum_{s=0}^m(\ol{\nabla}^{\man{M},\pi_{\man{E}}}B^m_s)
(\ol{\nabla}^{\man{M},\pi_{\man{E}},m}\xi)+
\sum_{s=0}^mB^m_s\otimes
\id_{\ctb{\man{M}}}(\ol{\nabla}^{\man{M},\pi_{\man{E}},m+1}\xi).
\end{equation*}
Combining the preceding two computations,
\begin{align*}
\nabla^{\man{M},\pi_{\man{E}},m+1}\xi=&\;
\ol{\nabla}^{\man{M},\pi_{\man{E}},m+1}\xi+
\left.\sum_{s=1}^m\right(
(\ol{\nabla}^{\man{M},\pi_{\man{E}}}B^m_s)
(\ol{\nabla}^{\man{M},\pi_{\man{E}},s}\xi)+
B^m_{s-1}\otimes\id_{\ctb{\man{M}}}
(\ol{\nabla}^{\man{M},\pi_{\man{E}},s}\xi)\\
&\;\left.+\sum_{j=1}^m\Ins_j(B^m_s(\ol{\nabla}^{\man{M},\pi_{\man{E}},s}\xi),
S_{\man{M}})+\Ins_{m+1}(B^m_s(\ol{\nabla}^{\man{M},\pi_{\man{E}},s}\xi),
S_{\pi_{\man{E}}})\right)\\
&\;+(\ol{\nabla}^{\man{M},\pi_{\man{E}}}B^m_0)(\xi)+
\sum_{j=1}^m\Ins_j(B^m_0(\xi),S_{\man{M}})+
\Ins_{m+1}(B^m_0(\xi),S_{\pi_{\man{E}}}),
\end{align*}
and from this the lemma follows.
\end{proof}
\end{lemma}

Now we give symmetrised versions of the preceding lemmata, since it is these
that are required for computations with jets.
\begin{lemma}\label{lem:symnablaonablaiterI}
Let\/ $r\in\{\infty,\omega\}$ and let\/
$\map{\pi_{\man{E}}}{\man{E}}{\man{M}}$ be a\/ $\C^r$-vector bundle.
Consider\/ $\C^r$-affine connections\/ $\nabla^{\man{M}}$ and\/
$\ol{\nabla}^{\man{M}}$ on\/ $\man{M}$\@, and\/ $\C^r$-vector bundle
connections\/ $\nabla^{\pi_{\man{E}}}$ and\/ $\ol{\nabla}^{\pi_{\man{E}}}$
in\/ $\man{E}$\@.  For\/ $m\in\integernn$\@, there exist\/ $\C^r$-vector
bundle mappings
\begin{equation*}
(\what{A}^m_s,\id_{\man{E}})\in
\vbmappings[r]{\tensor*[s]{\ctb{\man{M}}\otimes\man{E}}}
{\tensor*[m]{\ctb{\man{M}}}\otimes\man{E}},\qquad s\in\{0,1,\dots,m\},
\end{equation*}
such that
\begin{equation*}
(\Sym_m\otimes\id_{\man{E}})\scirc\ol{\nabla}^{\man{M},\pi_{\man{E}},m}\xi=
\sum_{s=0}^m\what{A}^m_s((\Sym_s\otimes\id_{\man{E}})
\scirc\nabla^{\man{M},\pi_{\man{E}},s}\xi)
\end{equation*}
for all\/ $\xi\in\sections[m]{\man{E}}$\@.
\begin{proof}
We define $\map{A^m}{\tensor*[\le m]{\ctb{\man{M}}}\otimes\man{E}}
{\tensor*[\le m]{\ctb{\man{M}}}\otimes\man{E}}$ by
\begin{equation*}
A^m(\xi,\nabla^{\pi_{\man{E}}}\xi,\dots,\nabla^{\man{M},\pi_{\man{E}},m}\xi)=
\left(A^0_0(\xi),\sum_{s=0}^1A^1_s(\nabla^{\man{M},\pi_{\man{E}},s}\xi),
\dots,\sum_{s=0}^mA^m_s(\nabla^{\man{M},\pi_{\man{E}},s}\xi)\right).
\end{equation*}
Let us organise the mappings we require into the following diagram:
\begin{equation}\label{eq:nablasymmetrise1}
\xymatrix{{\tensor*[\le m]{\ctb{\man{M}}}\otimes\man{E}}
\ar[rr]^{\Sym_{\le m}\otimes\id_{\man{E}}}\ar[d]_{A^m}&&
{\Symalg*[\le m]{\ctb{\man{M}}}\otimes\man{E}}
\ar[rr]^(0.55){S^m_{\ol{\nabla}^{\man{M}},\ol{\nabla}^{\pi_{\man{E}}}}}
\ar[d]^{\what{A}^m}&&{\jet{m}{\man{E}}}\ar@{=}[d]\\
{\tensor*[\le m]{\ctb{\man{M}}}\otimes\man{E}}
\ar[rr]^{\Sym_{\le m}\otimes\id_{\man{E}}}&&
{\Symalg*[\le m]{\ctb{\man{M}}}\otimes\man{E}}
\ar[rr]^(0.55){S^m_{\nabla^{\man{M}},\nabla^{\pi_{\man{E}}}}}&&
{\jet{m}{\man{E}}}}
\end{equation}
Here $\what{A}^m$ is defined so that the right square commutes.  We shall
show that the left square also commutes.  Indeed,
\begin{align*}
\what{A}^m\scirc\Sym_{\le m}\otimes\id_{\man{E}}&
(\xi,\nabla^{\pi_{\man{E}}}\xi,\dots,\nabla^{\man{M},\pi_{\man{E}},m}\xi)\\
=&\;(S^m_{\ol{\nabla}^{\man{M}},\ol{\nabla}^{\pi_{\man{E}}}})^{-1}\scirc
S^m_{\nabla^{\man{M}},\nabla^{\pi_{\man{E}}}}\scirc
(\Sym_{\le m}\otimes\id_{\man{E}})
(\xi,\ol{\nabla}^{\pi_{\man{E}}}\xi,\dots,\nabla^{\man{M},\pi_{\man{E}},m}\xi)\\
=&\;\Sym_{\le m}\otimes\id_{\vb{\man{E}}}
(\xi,\ol{\nabla}^{\pi_{\man{E}}}\xi,\dots,
\ol{\nabla}^{\man{M},\pi_{\man{E}},m}\xi)\\
=&\;(\Sym_{\le m}\otimes\id_{\man{E}})\scirc
A^m(\xi,\nabla^{\pi_{\man{E}}}\xi,\dots,\nabla^{\man{M},\pi_{\man{E}},m}\xi).
\end{align*}
Thus the diagram~\eqref{eq:nablasymmetrise1} commutes.  Thus, if we define
\begin{equation}\label{eq:AhatAnabla}
\what{A}^m_s((\Sym_s\otimes\id_{\man{E}})\scirc
\nabla^{\man{M},\pi_{\man{E}},s}\xi)=
(\Sym_m\otimes\id_{\man{E}})\scirc A^m_s(\nabla^{\man{M},\pi_{\man{E}},s}\xi),
\end{equation}
then we have
\begin{equation*}
(\Sym_m\otimes\id_{\man{E}})\scirc\ol{\nabla}^{\man{M},\pi_{\man{E}},m}\xi=
\sum_{s=0}^m\what{A}^m_s((\Sym_s\otimes\id_{\man{E}})
\scirc\nabla^{\man{M},\pi_{\man{E}},s}\xi),
\end{equation*}
as desired.
\end{proof}
\end{lemma}

The previous lemma has an ``inverse'' which we state next.
\begin{lemma}\label{lem:symnablaonablaiterII}
Let\/ $r\in\{\infty,\omega\}$ and let\/
$\map{\pi_{\man{E}}}{\man{E}}{\man{M}}$ be a\/ $\C^r$-vector bundle.
Consider\/ $\C^r$-affine connections\/ $\nabla^{\man{M}}$ and\/
$\ol{\nabla}^{\man{M}}$ on\/ $\man{M}$\@, and\/ $\C^r$-vector bundle
connections\/ $\nabla^{\pi_{\man{E}}}$ and\/ $\ol{\nabla}^{\pi_{\man{E}}}$
in\/ $\man{E}$\@.  For\/ $m\in\integernn$\@, there exist\/ $\C^r$-vector
bundle mappings
\begin{equation*}
(\what{B}^m_s,\id_{\man{E}})\in
\vbmappings[r]{\tensor*[s]{\ctb{\man{M}}\otimes\man{E}}}
{\tensor*[m]{\ctb{\man{M}}}\otimes\man{E}},\qquad s\in\{0,1,\dots,m\},
\end{equation*}
such that
\begin{equation*}
(\Sym_m\otimes\id_{\man{E}})\scirc\nabla^{\man{M},\pi_{\man{E}},m}\xi=
\sum_{s=0}^m\what{B}^m_s((\Sym_s\otimes\id_{\man{E}})
\scirc\ol{\nabla}^{\man{M},\pi_{\man{E}},s}\xi)
\end{equation*}
for all\/ $\xi\in\sections[m]{\man{E}}$\@.
\begin{proof}
The proof here is identical with the proof of
Lemma~\ref{lem:symnablaonablaiterI}\@, making the obvious notational
transpositions.
\end{proof}
\end{lemma}

The preceding four lemmata combine to give the following result.
\begin{lemma}
Let\/ $r\in\{\infty,\omega\}$ and let\/
$\map{\pi_{\man{E}}}{\man{E}}{\man{M}}$ be a\/ $\C^r$-vector bundle.
Consider\/ $\C^r$-affine connections\/ $\nabla^{\man{M}}$ and\/
$\ol{\nabla}^{\man{M}}$ on\/ $\man{M}$\@, and\/ $\C^r$-vector bundle
connections\/ $\nabla^{\pi_{\man{E}}}$ and\/ $\ol{\nabla}^{\pi_{\man{E}}}$
in\/ $\man{E}$\@.  For\/ $m\in\integernn$\@, there exist\/ $\C^r$-vector
bundle mappings
\begin{equation*}
A^m\in\vbmappings[r]{\jet{m}{\man{E}}}
{\Symalg*[\le m]{\ctb{\man{M}}}\otimes\man{E}},\quad
B^m\in\vbmappings[r]{\jet{m}{\man{E}}}
{\Symalg*[\le m]{\ctb{\man{M}}}\otimes\man{E}},
\end{equation*}
defined by
\begin{align*}
A^m(j_m\xi(x))=&\;
\Sym_{\le m}\otimes\id_{\man{E}}(\xi(x),\nabla^{\pi_{\man{E}}}\xi(x),\dots,
\nabla^{\man{M},\pi_{\man{E}},m}\xi(x)),\\
B^m(j_m\xi(x))=&\;\Sym_{\le m}\otimes\id_{\vb{\man{E}}}
(\xi(x),\ol{\nabla}^{\pi_{\man{E}}}\xi(x),\dots,
\ol{\nabla}^{\man{M},\pi_{\man{E}},m}\xi(x)).
\end{align*}
Moreover,\/ $A^m$ and\/ $B^m$ are isomorphisms, and
\begin{multline*}
B^m\scirc(A^m)^{-1}\scirc
(\Sym_{\le m}\otimes\id_{\man{E}})
(\xi(x),\nabla^{\pi_{\man{E}}}\xi(x),\dots,
\nabla^{\man{M},\pi_{\man{E}},m}\xi(x))\\
=\left(\xi(x),\sum_{s=0}^1\what{A}^1_s((\Sym_s\otimes\id_{\man{E}})
\scirc\nabla^{\man{M},\pi_{\man{E}},s}\xi(x)),\dots,
\sum_{s=0}^m\what{A}^m_s((\Sym_s\otimes\id_{\man{E}})
\scirc\nabla^{\man{M},\pi_{\man{E}},s}\xi(x))\right)
\end{multline*}
and
\begin{multline*}
A^m\scirc(B^m)^{-1}\scirc
(\Sym_{\le m}\otimes\id_{\man{E}})
(\xi(x),\ol{\nabla}^{\pi_{\man{E}}}\xi(x),\dots,
\ol{\nabla}^{\man{M},\pi_{\man{E}},m}\xi(x))\\
=\left(\xi(x),\sum_{s=0}^1\what{B}^1_s((\Sym_s\otimes\id_{\man{E}})
\scirc\ol{\nabla}^{\man{M},\pi_{\man{E}},s}\xi(e)),\dots,
\sum_{s=0}^m\what{B}^m_s((\Sym_s\otimes\id_{\man{E}})\scirc
\nabla^{\man{M},\pi_{\man{E}},s}\xi(x))\right),
\end{multline*}
where the vector bundle mappings\/ $\what{A}^m_s$ and\/ $\what{B}^m_s$\@,\/
$s\in\{0,1,\dots,m\}$\@, are as in Lemmata~\ref{lem:symnablaonablaiterI}
and~\ref{lem:symnablaonablaiterII}\@.
\end{lemma}

\subsection{Comparison of metric-related notions for different connections
and metrics}

We next consider how various constructions involving Riemannian metrics and
fibre metrics vary when one varies the fibre metrics.  The first result
concerns fibre norms for tensor products induced by a fibre metric.
\begin{lemma}\label{lem:G1G2comp}
Let\/ $\map{\pi_{\man{E}}}{\man{E}}{\man{M}}$ be a smooth vector bundle and
let\/ $\metric_1$ and\/ $\metric_2$ be smooth fibre metrics on\/ $\man{E}$\@.
Let\/ $\nbhd{K}\subset\man{M}$ be compact.  Then there exist\/
$C,\sigma\in\realp$ such that
\begin{equation*}
\frac{\sigma^{r+s}}{C}\dnorm{A(x)}_{\metric_2}\le\dnorm{A(x)}_{\metric_1}\le
\frac{C}{\sigma^{r+s}}\dnorm{A(x)}_{\metric_2}
\end{equation*}
for all\/ $A\in\sections[0]{\tensor[r,s]{\man{E}}}$\@,\/
$r,s\in\integernn$\@, and\/ $x\in\nbhd{K}$\@.
\begin{proof}
We begin by proving a linear algebra result.
\begin{proofsublemma}\label{psublem:metric-ind1}
If\/ $\metric_1$ and\/ $\metric_2$ are inner products on a
finite-dimensional\/ $\real$-vector space\/ $\alg{V}$\@, then there exists\/
$C\in\realp$ such that
\begin{equation*}
C^{-1}\metric_1(v,v)\le\metric_2(v,v)\le C\metric_1(v,v)
\end{equation*}
for all\/ $v\in\alg{V}$\@.
\begin{subproof}
Let $\metric^\flat_j\in\Hom_{\real}(\alg{V};\dual{\alg{V}})$ and
$\metric_j^\sharp\in\Hom_{\real}(\dual{\alg{V}};\alg{V})$\@, $j\in\{1,2\}$\@,
be the induced linear maps.  Note that
\begin{equation*}
\metric_1(\metric_1^\sharp\scirc\metric_2^\flat(v_1),v_2)=
\metric_2(v_1,v_2)=\metric_2(v_2,v_1)=
\metric_1(\metric_1^\sharp\scirc\metric_2^\flat(v_2),v_1),
\end{equation*}
showing that $\metric_1^\sharp\scirc\metric_2^\flat$ is
$\metric_1$-symmetric.  Let $\ifam{e_1,\dots,e_n}$ be a
$\metric_1$-orthonormal basis for $\alg{V}$ that is also a basis of
eigenvectors for $\metric_1^\sharp\scirc\metric_2^\flat$\@.  The matrix
representatives of $\metric_1$ and $\metric_2$ are then
\begin{equation*}
[\metric_1]=\begin{bmatrix}1&0&\cdots&0\\0&1&\cdots&0\\
\vdots&\vdots&\ddots&\vdots\\0&0&\cdots&1\end{bmatrix},\qquad{}
[\metric_2]=\begin{bmatrix}\lambda_1&0&\cdots&0\\0&\lambda_2&\cdots&0\\
\vdots&\vdots&\ddots&\vdots\\0&0&\cdots&\lambda_n\end{bmatrix},
\end{equation*}
where $\lambda_1,\dots,\lambda_n\in\realp$\@.  Let us assume, without loss of
generality, that
\begin{equation*}
\lambda_1\le\dots\le\lambda_n.
\end{equation*}
Then taking $C=\max\{\lambda_n,\lambda_1^{-1}\}$ gives the result, as one can
verify directly.
\end{subproof}
\end{proofsublemma}

Next we use the preceding sublemma to give the linear algebraic version of
the lemma.
\begin{proofsublemma}
Let\/ $\alg{V}$ be a finite-dimensional\/ $\real$-vector space and let\/
$\metric_1$ and\/ $\metric_2$ be inner products on\/ $\alg{V}$\@.  Then there
exist\/ $C,\sigma\in\realp$ such that
\begin{equation*}
\frac{\sigma^{r+s}}{C}\dnorm{A}_{\metric_2}\le\dnorm{A}_{\metric_1}\le
\frac{C}{\sigma^{r+s}}\dnorm{A}_{\metric_2}
\end{equation*}
for all\/ $A\in\tensor[r,s]{\alg{V}}$\@,\/ $r,s\in\integernn$\@.
\begin{subproof}
As in the proof of Sublemma~\ref{psublem:metric-ind1}\@, let
$\ifam{e_1,\dots,e_n}$ be a $\metric_1$-orthonormal basis for $\alg{V}$
consisting of eigenvectors for $\metric_1^\sharp\scirc\metric_2^\flat$\@.
Let $\lambda_1,\dots,\lambda_n\in\realp$ be the corresponding eigenvalues,
supposing that
\begin{equation*}
\lambda_1\le\dots\le\lambda_n.
\end{equation*}
Note that $\metric_2(e_j,e_k)=\delta_{jk}\lambda_j$\@, $j\in\{1,\dots,n\}$\@,
($\delta_{jk}$ being the Kronecker delta symbol) so that
$\ifam{\what{e}_1\eqdef\lambda_1^{-1}e_1,\dots,\what{e}_n\eqdef\lambda_n^{-1}e_n}$
is a $\metric_2$-orthonormal basis.  Denote by $\ifam{e^1,\dots,e^n}$ and
$\ifam{\what{e}^1,\dots,\what{e}^n}$ be the dual bases.  Note that
$\what{e}^j=\lambda_je^j$\@, $j\in\{1,\dots,n\}$\@.

Now let $A\in\tensor[r,s]{\alg{V}}$ and write
\begin{equation*}
A=\sum_{j_1,\dots,j_r=1}^n\sum_{k_1,\dots,k_s=1}^n
A^{j_1\cdots j_r}_{k_1\cdots k_s}e_{j_1}\otimes\dots\otimes e_{j_r}
\otimes e^{k_1}\otimes\dots\otimes e^{k_s}
\end{equation*}
and
\begin{equation*}
A=\sum_{j_1,\dots,j_r=1}^n\sum_{k_1,\dots,k_s=1}^n
\what{A}^{j_1\cdots j_r}_{k_1\cdots k_s}\what{e}_{j_1}\otimes\dots\otimes\what{e}_{j_r}
\otimes\what{e}^{k_1}\otimes\dots\otimes\what{e}^{k_s}.
\end{equation*}
We necessarily have
\begin{equation*}
\what{A}^{j_1\cdots j_r}_{k_1\cdots k_s}=
\lambda_{j_1}\cdots\lambda_{j_r}\lambda_{k_1}^{-1}\cdots\lambda_{k_s}^{-1}
A^{j_1\cdots j_r}_{k_1\cdots k_s},\qquad j_1,\dots,j_r,k_1,\dots,k_s\in\{1,\dots,n\}.
\end{equation*}
We have
\begin{equation*}
\dnorm{A}_{\metric_1}=\left(\sum_{j_1,\dots,j_r=1}^n\sum_{k_1,\dots,k_s=1}^n
\asnorm{A^{j_1\cdots j_r}_{k_1\cdots k_s}}^2\right)^{1/2},\quad
\dnorm{A}_{\metric_2}=\left(\sum_{j_1,\dots,j_r=1}^n\sum_{k_1,\dots,k_s=1}^n
\asnorm{\what{A}^{j_1\cdots j_r}_{k_1\cdots k_s}}^2\right)^{1/2}.
\end{equation*}
Therefore, if we let $\sigma=\min\{\lambda_1,\lambda_n^{-1}\}$\@, we have
\begin{equation*}
\dnorm{A}_{\metric_2}\le\sigma^{-(r+s)}\dnorm{A}_{\metric_1}.
\end{equation*}
This gives one half of the estimate in the sublemma, and the other is
established analogously.
\end{subproof}
\end{proofsublemma}

The lemma follows from the preceding sublemma since $C$ and $\sigma$ depend
only on $\metric_1$ and $\metric_2$ through the largest and smallest
eigenvalues of $\metric_1^\sharp\scirc\metric_2^\flat$\@, which are uniformly
bounded above and below on $\nbhd{K}$\@.
\end{proof}
\end{lemma}

Now we can compare fibre norms for jet bundles associated with different
metrics and connections.
\begin{lemma}\label{lem:normcompare}
Let\/ $r\in\{\infty,\omega\}$ and let\/
$\map{\pi_{\man{E}}}{\man{E}}{\man{M}}$ be a\/ $\C^r$-vector bundle.
Consider\/ $\C^r$-affine connections\/ $\nabla^{\man{M}}$ and\/
$\ol{\nabla}^{\man{M}}$ on\/ $\man{M}$\@, and\/ $\C^r$-vector bundle
connections\/ $\nabla^{\pi_{\man{E}}}$ and\/ $\ol{\nabla}^{\pi_{\man{E}}}$
in\/ $\man{E}$\@.  Consider\/ $\C^r$-Riemannian metrics\/ $\metric_{\man{M}}$
and\/ $\ol{\metric}_{\man{M}}$ for\/ $\man{M}$\@, and\/ $\C^r$-fibre
metrics\/ $\metric_{\pi_{\man{E}}}$ and\/ $\ol{\metric}_{\pi_{\man{E}}}$
for\/ $\man{E}$\@.  Let\/ $\nbhd{K}\subset\man{M}$ be compact.  Then there
exist\/ $C,\sigma\in\realp$ such that
\begin{equation*}
\frac{\sigma^m}{C}\dnorm{j_m\xi(x)}_{\ol{\metric}_{\man{M},\pi_{\man{E}},m}}
\le\dnorm{j_m\xi(x)}_{\metric_{\man{M},\pi_{\man{E}},m}}
\le\frac{C}{\sigma^m}\dnorm{j_m\xi(x)}_{\ol{\metric}_{\man{M},\pi_{\man{E}},m}}
\end{equation*}
for all\/ $\xi\in\sections[m]{\man{E}}$\@,\/ $m\in\integernn$\@, and\/
$x\in\nbhd{K}$\@.
\begin{proof}
We first make some preliminary constructions that will be useful.

By Lemma~\ref{lem:symnablaonablaiterI}\@, we have
\begin{equation}\label{eq:onabla-nabla}
\begin{aligned}
\xi(x)=&\;\what{A}^0_0\xi(x),\\
(\Sym_1\otimes\id_{\man{E}})\ol{\nabla}^{\pi_{\man{E}}}\xi(x)=&\;
\what{A}^1_1(\nabla^{\pi_{\man{E}}}\xi(x))+\what{A}^1_0(\xi(x)),\\
\vdots\,&\\
(\Sym_m\otimes\id_{\man{E}})\scirc\ol{\nabla}^{\man{M},\pi_{\man{E}},m}\xi(x)
=&\;\sum_{s=0}^m\what{A}^m_s((\Sym_s\otimes\id_{\man{E}})\scirc
\nabla^{\man{M},\pi_{\man{E}},s}\xi(x)).
\end{aligned}
\end{equation}
In like manner, by Lemma~\ref{lem:symnablaonablaiterII}\@, we have
\begin{equation}\label{eq:nabla-onabla}
\begin{aligned}
\xi(x)=&\;\what{B}^0_0\xi(x),\\
(\Sym_1\otimes\id_{\man{E}})\nabla^{\pi_{\man{E}}}\xi(x)=&\;
\what{B}^1_1(\ol{\nabla}^{\pi_{\man{E}}}\xi(x))+\what{B}^1_0(\xi(x)),\\
\vdots\,&\\
(\Sym_m\otimes\id_{\man{E}})\scirc\nabla^{\man{M},\pi_{\man{E}},m}\xi(x)=&\;
\sum_{s=0}^m\what{B}^m_s((\Sym_s\otimes\id_{\man{E}})\scirc
\ol{\nabla}^{\man{M},\pi_{\man{E}},s}\xi(x)).
\end{aligned}
\end{equation}

By Lemma~\ref{lem:A(S)<AS}\@, we have
\begin{equation*}
\dnorm{A^m_s(\beta_s)}_{\ol{\metric}_{\man{M},\pi_{\man{E}}}}\le
\dnorm{A^m_s}_{\ol{\metric}_{\man{M},\pi_{\man{E}}}}
\dnorm{\beta_s}_{\ol{\metric}_{\man{M},\pi_{\man{E}}}}
\end{equation*}
for $\beta_s\in\tensor*[s]{\ctb{\man{M}}}\otimes\man{E}$\@,
$m\in\integerp$\@, $s\in\{0,1,\dots,m\}$\@.  By Lemma~\ref{lem:Symnorm}\@,
\begin{equation*}
\dnorm{\Sym_s(A)}_{\ol{\metric}_{\man{M},\pi_{\man{E}}}}\le
\dnorm{A}_{\ol{\metric}_{\man{M},\pi_{\man{E}}}}
\end{equation*}
for $A\in\tensor*[s]{\ctb{\man{E}}}$ and $s\in\integerp$\@.  Thus,
recalling~\eqref{eq:AhatAnabla}\@,
\begin{equation*}
\dnorm{\what{A}^m_s(\Sym_s(\beta_s))}_{\ol{\metric}_{\man{M},\pi_{\man{E}}}}=
\dnorm{\Sym_m\scirc A^m_s(\beta_s)}_{\ol{\metric}_{\man{M},\pi_{\man{E}}}}
\le\dnorm{A^m_s}_{\ol{\metric}_{\man{E}}}
\dnorm{\beta_s}_{\ol{\metric}_{\man{M},\pi_{\man{E}}}},
\end{equation*}
for $\beta_s\in\tensor*[s]{\pi_{\man{E}}^*\ctb{\man{M}}}\otimes\man{E}$\@,
$m\in\integerp$\@, $s\in\{1,\dots,m\}$\@.

By Lemmata~\ref{lem:pissynabla} and~\ref{lem:particular-bounds} with
$r=0$\@, there exist $\sigma_1,\rho_1\in\realp$ such that
\begin{equation*}
\dnorm{A^k_s(x)}_{\ol{\metric}_{\man{M},\pi_{\man{E}}}}\le
\sigma_1^{-k}\rho_1^{-(k-s)}(k-s)!,\qquad
k\in\integernn,\ s\in\{0,1,\dots,k\},\ x\in\nbhd{K}.
\end{equation*}
Without loss of generality, we assume that $\sigma_1,\rho_1\le1$\@.  Thus,
abbreviating $\sigma_2=\sigma_1\rho_1$\@, we have
\begin{equation*}
\dnorm{\what{A}^k_s((\Sym_s\otimes\id_{\man{E}})\scirc
\ol{\nabla}^{\man{M},\pi_{\man{E}},s}\xi(x)}_{\ol{\metric}_{\man{M},\pi_{\man{E}}}}
\le C_1\sigma_2^{-k}(k-s)!
\dnorm{(\Sym_s\otimes\id_{\man{E}})\scirc\nabla^{\man{M},\pi_{\man{E}},s}
\xi(x)}_{\ol{\metric}_{\man{M},\pi_{\man{E}}}}
\end{equation*}
for $m\in\integernn$\@, $k\in\{0,1,\dots,m\}$\@, $s\in\{0,1,\dots,k\}$\@,
$x\in\nbhd{K}$\@.  Thus, by~\eqref{eq:12inftynorms}
and~\eqref{eq:nabla-onabla}\@,
\begin{align*}
\dnorm{j_m\xi(x)}_{\ol{\metric}_{\man{M},\pi_{\man{E}},m}}\le&\;
\sum_{k=0}^m\frac{1}{k!}\dnorm{(\Sym_k\otimes\id_{\man{E}})\scirc
\nabla^{\man{M},\pi_{\man{E}},k}\xi(x)}_{\metric_{\man{M},\pi_{\man{E}}}}\\
=&\;\sum_{k=0}^m\frac{1}{k!}\adnorm{\sum_{s=0}^k
\what{A}^k_s((\Sym_s\otimes\id_{\man{E}})\scirc\nabla^{\man{M},\pi_{\man{E}},s}
\xi(x))}_{\metric_{\man{M},\pi_{\man{E}}}}\\
\le&\;\sum_{k=0}^m\sum_{s=0}^kC_1\sigma_2^{-k}
\frac{s!(k-s)!}{k!}\frac{1}{s!}
\dnorm{(\Sym_s\otimes\id_{\man{E}})\scirc\nabla^{\man{M},\pi_{\man{E}},s}
\xi(x)}_{\metric_{\man{M},\pi_{\man{E}}}}
\end{align*}
for $x\in\nbhd{K}$ and $m\in\integernn$\@.  Now note that
\begin{equation*}
\frac{s!(k-s)!}{k!}\le1,\quad
C_1\sigma_2^{-k}\le C_1\sigma_2^{-m},
\end{equation*}
for $s\in\{0,1,\dots,m\}$\@, $k\in\{0,1,\dots,s\}$\@, since $\sigma_2\le1$\@.
Then
\begin{align*}
\dnorm{j_m\xi(x)}_{\ol{\metric}_{\man{M},\pi_{\man{E}},m}}\le&\;
C_1\sigma_2^{-m}\sum_{k=0}^m\sum_{s=0}^k\frac{1}{s!}
\dnorm{(\Sym_s\otimes\id_{\man{E}})\scirc
\nabla^{\man{M},\pi_{\man{E}},s}\xi(x)}_{\metric_{\man{M},\pi_{\man{E}}}}\\
\le&\;C_1\sigma_2^{-m}\sum_{k=0}^m\sum_{s=0}^m
\frac{1}{s!}\dnorm{(\Sym_s\otimes\id_{\man{E}})\scirc\nabla^{\man{M},s}
\xi(x)}_{\metric_{\man{M},\pi_{\man{E}}}}\\
=&\;(m+1)C_1\sigma_2^{-m}\sum_{s=0}^m\frac{1}{s!}
\dnorm{(\Sym_s\otimes\id_{\man{E}})\scirc\nabla^{\man{M},\pi_{\man{E}},s}
\xi(x)}_{\metric_{\man{M},\pi_{\man{E}}}}.
\end{align*}
Now let $\sigma<\sigma_2$ and note that
\begin{equation*}
\lim_{m\to\infty}(m+1)\frac{\sigma_2^{-m}}{\sigma^{-m}}=0.
\end{equation*}
Thus there exists $N\in\integerp$ such that
\begin{equation*}
(m+1)C_1\sigma_2^{-m}\le C_1\sigma^{-m},\qquad m\ge N.
\end{equation*}
Let
\begin{equation*}
C=\max\left\{C_1,2C_1\frac{\sigma}{\sigma_2},
3C_1\left(\frac{\sigma}{\sigma_2}\right)^2,
\dots,(N+1)C_1\left(\frac{\sigma}{\sigma_2}\right)^N\right\}.
\end{equation*}
We then immediately have $(m+1)C_1\sigma_2^{-m}\le C\sigma^{-m}$ for all
$m\in\integernn$\@.  We then have, using~\eqref{eq:12inftynorms}\@,
\begin{align*}
\dnorm{j_m\xi(x)}_{\ol{\metric}_{\man{M},\pi_{\man{E}},m}}\le&\;
C\sigma^{-m}\sum_{s=0}^m\frac{1}{s!}\dnorm{(\Sym_s\otimes\id_{\man{E}})
\scirc\nabla^{\man{M},\pi_{\man{E}},s}\xi(x)}_{\metric_{\man{M},\pi_{\man{E}}}}\\
=&\;C\sqrt{m}\sigma^{-m}\dnorm{j_m\xi(x)}_{\metric_{\man{M},\pi_{\man{E}},m}}.
\end{align*}
After modifying $C$ and $\sigma$ in the manner of the computations just
preceding, we have
\begin{equation*}
\dnorm{j_m\xi(x)}_{\ol{\metric}_{\man{M},\pi_{\man{E}},m}}\le
C\sigma^{-m}\dnorm{j_m\xi(x)}_{\metric_{\man{M},\pi_{\man{E}},m}}.
\end{equation*}
This gives one half of the desired pair of estimates.

For the other half of the estimate, we use~\eqref{eq:nabla-onabla}\@, and
Lemmata~\ref{lem:pissynabla} and~\ref{lem:particular-bounds} in the
computations above to arrive at the estimate
\begin{equation*}
\dnorm{j_m\xi(x)}_{\ol{\metric}_{\man{M},\pi_{\man{E}},m}}\le
C\sigma^{-m}\dnorm{j_m\xi(x)}_{\metric_{\man{M},\pi_{\man{E}},m}},
\end{equation*}
which gives the result.
\end{proof}
\end{lemma}

\begin{remark}
The preceding result holds in the smooth case, and with a much easier proof.
In the result, one can replace ``$C\sigma^{-m}$'' with a fixed constant
``$C$'' for each $m$\@.  For this reason, the proof is also far simpler, as
one need not keep track of all the factorial terms that give rise to the
exponential component in the estimates.\oprocend
\end{remark}

\subsection{Local descriptions of the real analytic topology}

We endeavour to make our presentation as unencumbered of coordinates as
possible.  While the intrinsic jet bundle characterisations of the seminorms
are useful for abstract definitions and proofs, concrete proofs often require
local descriptions of the topologies.  In this section we provide these local
descriptions of the topologies.  By proving that these local descriptions are
equivalent to the intrinsic descriptions above, we also prove that these
intrinsic descriptions of topologies do not depend on the choice of metrics
or connections.

Let us develop the notation for working with local descriptions of
topologies.  Let $\nbhd{U}\subset\real^n$ be an open set and let
$\vect{\Phi}\in\mappings[\omega]{\nbhd{U}}{\real^k}$\@.  We define local
seminorms as follows.  For $\nbhd{K}\subset\nbhd{U}$ compact and for
$\vect{a}\in\c_0(\integernn;\realp)$\@, denote
\begin{multline*}
p^{\prime\omega}_{\nbhd{K},\vect{a}}(\vect{\Phi})=
\sup\left\{\frac{a_0a_1\cdots a_m}{I!}
\snorm{\linder[I]{\Phi^a}(\vect{x})}\right|\\
\left.\vphantom{\frac{a_0}{I}}\vect{x}\in\nbhd{K},\ a\in\{1,\dots,k\},\
I\in\integernn^n,\ \snorm{I}\le m,\ m\in\integernn\right\}.
\end{multline*}
These seminorms, defined for all compact $\nbhd{K}\subset\nbhd{U}$ and
$\vect{a}\in\c_0(\integernn;\realp)$\@, define the \defn{local
$\C^\omega$-topology} for $\mappings[\omega]{\nbhd{U}}{\real^k}$\@.

There are many possible variations of the seminorms that one can use, and
these variations are equivalent to the seminorms above.  For example, rather
than using the $\infty$-vector norm, one might use the $2$-vector norm.  In
doing so, one uses~\eqref{eq:12inftynorms} to give
\begin{multline*}
\sup\setdef{\snorm{\linder[I]{\Phi^a}(\vect{x})}}{I\in\integernn^n,\
\snorm{I}=m,\ a\in\{1,\dots,k\}}\le
\dnorm{\linder[m]{\vect{\Phi}}(\vect{x})}\\
\le\sqrt{kn^m}\sup\setdef{\snorm{\linder[I]{\Phi^a}(\vect{x})}}
{I\in\integernn^n,\ \snorm{I}=m,\ a\in\{1,\dots,k\}}.
\end{multline*}
If we define
\begin{equation*}
b_0=2\sqrt{k}a_0,\ b_j=2\sqrt{n}a_j,\qquad j\in\integerp,
\end{equation*}
then, noting that $n^j\le n^m$ for $j\in\{0,1,\dots,m\}$ and that
$m+1\le 2^m$ for $m\in\integernn$\@, we have
\begin{multline*}
p^{\prime\omega}_{\nbhd{K},\vect{a}}(\vect{\Phi})\le
\sup\left\{\frac{a_0a_1\cdots a_m}{I!}
\dnorm{\linder[m]{\vect{\Phi}}(\vect{x})}\right|\\
\left.\vphantom{\frac{a_0}{I}}\vect{x}\in\nbhd{K},\
I\in\integernn^n,\ \snorm{I}\le m,\ m\in\integernn\right\}
\le p^{\prime\omega}_{\nbhd{K},\vect{b}}(\vect{\Phi}),
\end{multline*}
and this gives equivalence of the topologies using the $\infty$- and
$2$-norms.  Another variation in the seminorms is that one might scale the
derivatives by $\frac{1}{\snorm{I}!}$ rather than $\frac{1}{I!}$\@.  In this
case, we use the standard multinomial estimate~\eqref{eq:multinomialest} to
give
\begin{equation*}
\frac{\snorm{I}!}{I!}\le n^m.
\end{equation*}
Thus, if we take
\begin{equation*}
b_0=a_0,\ b_j=na_j,\qquad j\in\integerp,
\end{equation*}
we have
\begin{multline*}
p^{\prime\omega}_{\nbhd{K},\vect{b}}(\vect{\Phi})\le\\
\sup\asetdef{\frac{a_0a_1\cdots a_m}{\snorm{I}!}
\snorm{\linder[I]{\Phi^a}(\vect{x})}}
{\vect{x}\in\nbhd{K},\ a\in\{1,\dots,k\},\ I\in\integernn^n,\
\snorm{I}\le m,\ m\in\integernn}\\
\le p^{\prime\omega}_{\nbhd{K},\vect{a}}(\vect{\Phi}).
\end{multline*}
This gives the equivalence of the topologies defined using the scaling factor
$\frac{1}{\snorm{I}!}$ for derivatives in place of $\frac{1}{I!}$\@.  One can
also combine the previous modifications.  Indeed, if we use the $2$-norm and
the scaling factor $\frac{1}{\snorm{I!}}$\@, then one readily sees that we
recover the intrinsic seminorms on the trivial vector bundle
$\real^k_{\nbhd{U}}$ of Section~\ref{subsec:seminorms} using (1)~the
Euclidean inner product for the Riemannian metric on $\nbhd{U}$ and for the
fibre metric on $\real^k$ and (2)~standard differentiation as covariant
differentiation.  We shall use this observation in the proof of
Theorem~\ref{the:same-topology} below.

We wish to show that these local topologies can be used to define a topology
for $\sections[\omega]{\man{E}}$ that is equivalent to the intrinsic
topologies defined in Section~\ref{subsec:seminorms} using jet bundles,
connections, and metrics.  To state the result, let us indicate some
notation.  Let $(\nbhd{V},\psi)$ be a vector bundle chart for
$\map{\pi_{\man{E}}}{\man{E}}{\man{M}}$ with $(\nbhd{U},\phi)$ the induced
chart for $\man{M}$\@.  Suppose that
$\psi(\nbhd{V})=\phi(\nbhd{U})\times\real^k$\@.  Given a section $\xi$\@, we
define $\map{\psi_*(\xi)}{\phi(\nbhd{U})}{\real^k}$ by requiring that
\begin{equation*}
\psi\scirc\xi\scirc\phi^{-1}(\vect{x})=(\vect{x},\psi_*(\xi)(\vect{x})).
\end{equation*}
With this notation, we have the following result.
\begin{theorem}\label{the:same-topology}
Let\/ $\map{\pi_{\man{E}}}{\man{E}}{\man{M}}$ be a\/ $\C^\omega$-vector
bundle.  Let\/ $\metric_{\man{M}}$ be a Riemannian metric on\/ $\man{M}$\@,
let\/ $\metric_{\pi_{\man{E}}}$ be a vector bundle metric on\/ $\man{E}$\@,
let\/ $\nabla^{\man{M}}$ be an affine connection on\/ $\man{M}$\@, and let\/
$\nabla^{\pi_{\man{E}}}$ be a vector bundle connection on\/ $\man{E}$\@, with
all of these being of class\/ $\C^\omega$\@.  Then the following two
collections of seminorms for\/ $\sections[\omega]{\man{E}}$ define the same
topology:
\begin{compactenum}[(i)]
\item $p^\omega_{\nbhd{K},\vect{a}}$\@,\/
$\vect{a}\in\c_0(\integernn;\realp)$\@,\/ $\nbhd{K}\subset\man{M}$ compact;
\item $p^{\prime\omega}_{\nbhd{K},\vect{a}}\scirc\psi_*$\@,\/
$\vect{a}\in\c_0(\integernn;\realp)$\@,\/ $\nbhd{K}\subset\phi(\nbhd{U})$
compact,\/ $(\nbhd{V},\psi)$ is a vector bundle chart for\/ $\man{E}$ with\/
$(\nbhd{U},\phi)$ the induced chart for\/ $\man{M}$\@.
\end{compactenum}
\begin{proof}
As alluded to in the discussion above, it suffices to use the norm
\begin{equation*}
\dnorm{\linder[m]{\vect{\Phi}}(\vect{x})}_2=
\left(\sum_{\substack{I\in\integernn^n\\\snorm{I}=m}}\sum_{a=1}^k
\snorm{\linder[I]{\Phi^a}(\vect{x})}^2\right)^{1/2}
\end{equation*}
for derivatives of $\real^k$-valued functions on $\nbhd{U}\subset\real^n$\@.
If we denote
\begin{equation*}
j_m\vect{\Phi}(\vect{x})=(\vect{\Phi}(\vect{x}),
\linder{\vect{\Phi}}(\vect{x}),\dots,\linder[m]{\vect{\Phi}}(\vect{x})),
\end{equation*}
then we define
\begin{equation*}
\dnorm{j_m\vect{\Phi}(\vect{x})}^2_{2,m}=\sum_{j=0}^m
\frac{1}{(j!)^2}\dnorm{\linder[j]{\vect{\Phi}}(\vect{x})}^2_2,
\end{equation*}
this norm agreeing with the fibre norms used in
Section~\ref{subsec:seminorms} with the flat connections and with the
Euclidean inner products.  We use these norms to define seminorms that we
denote by $q^{\prime}$ in place of the local seminorms $p^{\prime}$ as above.

We might like to use Lemma~\ref{lem:normcompare} in this proof.  However, we
cannot do so.  The reason for this is that the proof of
Lemma~\ref{lem:normcompare} makes reference to Lemma~\ref{lem:pissynabla}\@.
The proof of this lemma relies on the bound~\eqref{eq:Sbound}\@, which is
deduced from Lemma~\ref{lem:Comegachar}\@.  The proof of
Lemma~\ref{lem:Comegachar} in~\cite{SJ/ADL:14a}\@, we note, relies on exactly
what we are now proving.  To intrude on the potential circular logic, we must
give a proof of this part of the theorem that does not rely on
Lemma~\ref{lem:normcompare}\@.  In fact, the only part of the chain of
results that we need to prove independently is the bound~\eqref{eq:Sbound}\@.
In particular, if we can show that Lemma~\ref{lem:normcompare} holds in the
current situation where
\begin{compactenum}
\item $\man{M}=\nbhd{U}\subset\real^n$ and $\man{E}=\real^k_{\nbhd{U}}$\@,
\item $\ol{\metric}_{\nbhd{U}}$ and $\ol{\metric}_{\pi_{\man{E}}}$ are the
Euclidean inner products, and
\item $\ol{\nabla}^{\man{M}}$ and $\ol{\nabla}^{\pi_{\man{E}}}$ are the flat connections,
\end{compactenum}
this will be enough to make use of this result.

Let $(\nbhd{V},\psi)$ be a vector bundle chart for $\man{E}$ with
$(\nbhd{U},\phi)$ the chart for $\man{M}$\@.  Standard estimates for real
analytic functions~\cite[\eg][Proposition~2.2.10]{SGK/HRP:02} give
$C_1,\sigma_1\in\realp$ such that
\begin{equation*}
\dnorm{D^r_{\ol{\nabla}_{\nbhd{U}},\ol{\nabla}^{\pi_{\man{E}}}}
S_{\nbhd{U}}(\vect{x})}_2,
\dnorm{D^r_{\ol{\nabla}_{\nbhd{U}},\ol{\nabla}^{\pi_{\man{E}}}}
S_{\pi_{\man{E}}}(\vect{x})}_2\le C_1\sigma_1^{-r}r!,\qquad\vect{x}\in\nbhd{K}.
\end{equation*}
This gives the bound~\eqref{eq:Sbound} in this case, and so we can use
Lemma~\ref{lem:particular-bounds}\@, and then Lemma~\ref{lem:pissynabla}\@,
and then the computation of Lemma~\ref{lem:normcompare} to give
\begin{equation*}
\frac{\sigma^m}{C}\dnorm{j_m\xi}_{\metric_{\man{M},\pi_{\man{E}},m}}\le
\dnorm{j_m(\psi_*(\xi))(\phi(x))}_{2,m}
\le\frac{C}{\sigma^m}\dnorm{j_m\xi}_{\metric_{\man{M},\pi_{\man{E}},m}}.
\end{equation*}

Now, having established Lemma~\ref{lem:normcompare} in the case of interest,
we proceed with the proof, making use of this fact.

Let $\nbhd{K}\subset\phi(\nbhd{U})$ be compact and let
$\vect{a}\in\c_0(\integernn;\realp)$\@.  As per our appropriate version of
Lemma~\ref{lem:normcompare}\@, there exist $C,\sigma\in\realp$ such that
\begin{equation*}
\dnorm{j_m(\psi_*(\xi))(\phi(x))}_{2,m}\le
\frac{C}{\sigma^m}\dnorm{j_m\xi(x)}_{\metric_{\man{M},\pi_{\man{E}},m}}
\end{equation*}
for every $\xi\in\sections[\omega]{\man{E}}$\@, $x\in\phi^{-1}(\nbhd{K})$\@,
and $m\in\integernn$\@.  Then
\begin{equation*}
a_0a_1\cdots a_m\dnorm{j_m(\psi_*(\xi))(\phi(x))}_{2,m}\le
\frac{Ca_0a_1\cdots a_m}{\sigma^m}
\dnorm{j_m\xi(x)}_{\metric_{\man{M},\pi_{\man{E}},m}}
\end{equation*}
for every $\xi\in\sections[\omega]{\man{E}}$\@, $x\in\phi^{-1}(\nbhd{K})$\@,
and $m\in\integernn$\@.  Define $\vect{b}\in\c_0(\integernn;\realp)$ by
\begin{equation*}
b_0=Ca_0,\ b_j=\frac{a_j}{\sigma},\qquad j\in\integerp.
\end{equation*}
Then, taking supremums of the preceding inequality gives
\begin{equation*}
q^{\prime\omega}_{\nbhd{K},\vect{a}}\scirc\psi_*(\xi)\le
p^\omega_{\phi^{-1}(\nbhd{K}),\vect{b}}(\xi)
\end{equation*}
for $\xi\in\sections[\omega]{\man{E}}$\@.

Now let $\nbhd{K}\subset\man{M}$ be compact and let
$\vect{a}\in\c_0(\integernn;\realp)$\@.  Let $x\in\nbhd{K}$ and let
$(\nbhd{V}_x,\psi_x)$ be a vector bundle chart for $\man{E}$ with
$(\nbhd{U}_x,\phi_x)$ the chart for $\man{M}$ with $x\in\nbhd{U}_x$\@.  We
suppose that $\nbhd{U}_x$ is relatively compact and that, by our appropriate
version of Lemma~\ref{lem:normcompare}\@, there exist $C_x,\sigma_x\in\realp$
such that
\begin{equation*}
\dnorm{j_m\xi(y)}_{\metric_{\man{M},\pi_{\man{E}},m}}\le
\frac{C_x}{\sigma_x^m}\dnorm{j_m(\psi_*\xi)(y)}_{m,2}
\end{equation*}
for $\xi\in\sections[\omega]{\man{E}}$\@, $y\in\closure(\nbhd{U}_x)$\@,
$m\in\integernn$\@.  Therefore,
\begin{equation*}
a_0a_1\cdots a_m\dnorm{j_m\xi(y)}_{\metric_{\pi_{\man{E}},m}}\le
\frac{C_xa_0a_1\cdots a_m}{\sigma_x^m}\dnorm{j_m(\psi_*\xi)(y)}_{m,2}
\end{equation*}
for $\xi\in\sections[\omega]{\man{E}}$\@, $y\in\closure(\nbhd{U}_x)$\@,
$m\in\integernn$\@.  Compactness of $\nbhd{K}$ gives
$x_1,\dots,x_s\in\nbhd{K}$ such that
$\nbhd{K}\subset\cup_{j=1}^s\nbhd{U}_{x_j}$ and we then take
\begin{equation*}
C=\max\{C_{x_1},\dots,C_{x_s}\},\quad
\sigma=\min\{\sigma_{x_1},\dots,\sigma_{x_s}\}.
\end{equation*}
We define $\vect{b}\in\c_0(\integernn;\realp)$ by
\begin{equation*}
b_0=Ca_0,\ b_j=\frac{a_j}{\sigma},\qquad j\in\integerp.
\end{equation*}
We then arrive at the inequality
\begin{equation*}
p^\omega_{\nbhd{K},\vect{a}}(\xi)\le
q^{\prime\omega}_{\closure(\nbhd{U}_{x_1}),\vect{b}}\scirc\psi_{x_1*}(\xi)+\dots+
q^{\prime\omega}_{\closure(\nbhd{U}_{x_s}),\vect{b}}\scirc\psi_{x_s*}(\xi)
\end{equation*}
which is valid for $\xi\in\sections[\omega]{\man{E}}$\@.
\end{proof}
\end{theorem}

\begin{remark}
The preceding theorem holds in the smooth case.  The proof is slightly
simpler in the smooth case, unlike in the proof of
Lemma~\ref{lem:normcompare} where the smooth case is significantly simpler
than the real analytic case.  Note also that, in the smooth case, one does
not need the local estimates for derivatives of real analytic functions, so
this also significantly simplifies the logic.\oprocend
\end{remark}

An immediate consequence of the theorem is that the topologies defined by the
seminorms of Section~\ref{subsec:seminorms} are independent of the choice of
connections $\nabla^{\man{M}}$ and $\nabla^{\pi_{\man{E}}}$\@, Riemannian
metric $\metric_{\man{M}}$\@, and vector bundle metric
$\metric_{\pi_{\man{E}}}$\@, since the preceding result shows that all such
topologies are the same as the one defined by local seminorms.

\section{Continuity of standard geometric operations}\label{sec:continuity}

In this section we put to use the somewhat complicated results of the
preceding sections to prove the continuity of standard algebraic and
differential operations on real analytic manifolds.  The reader will notice
as they go through the proofs that there are definite themes that emerge from
the various proofs of continuity.  Moreover, we take full advantage of the
results from Section~\ref{subsec:tensornorms} that were nominally developed
to prove the bounds of Lemma~\ref{lem:pissynabla}\@, so illustrating their
general utility.  We hope that a demonstration of the collection of
results\textemdash{}some easy, other less easy\textemdash{}will prove useful.

As a general comment on the results in this section, we shall prove in many
cases that certain linear mappings between spaces of sections of real
analytic vector bundles are continuous and open onto their
image,~\ie~homeomorphisms onto their image.  One might hope to do this with a
general Open Mapping Theorem.  Indeed, since the space of real analytic
sections of a vector bundle is both webbed and ultrabornological, one is in a
perhaps in a position to use the Open Mapping Theorem of \citet{MDW:67} (see
also \cite[Theorem~24.30]{RM/DV:97}).  However, since the images of our
mappings are not necessarily ultrabornological (even closed subspaces of
ultrabornological spaces may not be ultrabornological), we typically prove
the openness by a direct argument, by virtue of our having given in
Section~\ref{sec:lift-isomorphisms} relations between iterated covariant
derivatives going ``both ways.''  Moreover, the use of seminorms to prove
these results is in keeping with the general tenor of this work.

As we have indicated as we have been going along, the results in this section
are applicable to the smooth case.  We shall indicate the modifications
required in sample cases, with the general situation following easily from
these.

\subsection{Continuity of algebraic operations}

We begin with a consideration of continuity of standard algebraic operations
with vector bundles.
\begin{theorem}\label{the:+xcont}
Let\/ $\map{\pi_{\man{E}}}{\man{E}}{\man{M}}$ and\/
$\map{\pi_{\man{F}}}{\man{F}}{\man{M}}$ be\/ $\C^\omega$-vector bundles.
Then the following mappings are continuous:
\begin{compactenum}[(i)]
\item \label{pl:+cont} $\sections[\omega]{\man{E}}\oplus\sections[\omega]{\man{E}}
\ni(\xi,\eta)\mapsto\xi+\eta\in\sections[\omega]{\man{E}}$\@;
\item \label{pl:xcont} $\sections[\omega]{\man{F}\otimes\dual{\man{E}}}
\times\sections[\omega]{\man{E}}\ni(L,\xi)\mapsto L\scirc\xi\in
\sections[\omega]{\man{F}}$\@.\savenum
\end{compactenum}
Also, fixing an injective vector bundle mapping\/
$L\in\sections[\omega]{\man{F}\otimes\dual{\man{E}}}$\@, the following
mapping is open onto its image:
\begin{compactenum}[(i)]\resumenum
\item \label{pl:xopen}
$\sections[\omega]{\man{E}}\ni\xi\mapsto L\scirc\xi\in
\sections[\omega]{\man{F}}$\@.
\end{compactenum}
\begin{proof}
We suppose that we have a real analytic affine connection $\nabla^{\man{M}}$
on $\man{M}$\@, and real analytic vector bundle connections
$\nabla^{\pi_{\man{E}}}$ and $\nabla^{\pi_{\man{F}}}$ in $\man{E}$ and
$\man{F}$\@, respectively.  We suppose that we have a real analytic
Riemannian metric $\metric_{\man{M}}$ on $\man{M}$\@, and real analytic fibre
metrics $\metric_{\pi_{\man{E}}}$ and $\metric_{\pi_{\man{F}}}$ on $\man{E}$
and $\man{F}$\@, respectively.  This gives the seminorms
$p_{\nbhd{K},\vect{a}}^\omega$ and $q_{\nbhd{K},\vect{a}}^\omega$\@,
$\nbhd{K}\subset\man{M}$ compact, $\vect{a}\in\c_0(\integernn;\realp)$\@, for
$\sections[\omega]{\man{E}}$ and $\sections[\omega]{\man{F}}$\@,
respectively.  We denote the induced seminorms for
$\sections[\omega]{\man{F}\otimes\dual{\man{E}}}$ by
$q_{\nbhd{K},\vect{a}}^\omega\otimes p_{\nbhd{K},\vect{a}}^\omega$\@,
$\nbhd{K}\subset\man{M}$ compact, $\vect{a}\in\c_0(\integernn;\realp)$\@.

\eqref{pl:+cont} The fibre norms from Section~\ref{subsec:jet-norms} satisfy
the triangle inequality, and this readily gives
\begin{equation*}
p^\omega_{\nbhd{K},\vect{a}}(\xi+\eta)\le p^\omega_{\nbhd{K},\vect{a}}(\xi)+
p^\omega_{\nbhd{K},\vect{a}}(\eta),
\end{equation*}
which immediately gives this part of the result.

\eqref{pl:xcont} Let us make some preliminary computations from which this
part of the theorem will follow easily.

First, by Lemma~\ref{lem:A(S)<AS}\@, we have
\begin{equation}\label{eq:normLxi}
\dnorm{L\scirc\xi(x)}_{\metric_{\man{M},\pi_{\man{F}}}}\le
\dnorm{L(x)}_{\metric_{\man{M},\pi_{\man{F}}\otimes\pi_{\man{E}}}}
\dnorm{\xi(x)}_{\metric_{\man{M},\pi_{\man{E}}}}.
\end{equation}

Next, by Lemmata~\ref{lem:leibniz}\@,~\ref{lem:A(S)<AS}\@,
and~\ref{lem:Symnorm}\@, we have
\begin{equation*}
\dnorm{D^k_{\nabla^{\man{M}},\nabla^{\pi_{\man{F}}}}
(L\scirc\xi(x))}_{\metric_{\man{M},\pi_{\man{F}}}}
\le\sum_{j=0}^k\binom{k}{j}
\dnorm{D^j_{\nabla^{\man{M}},\nabla^{\pi_{\man{F}}\otimes\pi_{\man{E}}}}
L(x)}_{\metric_{\man{M},\pi_{\man{F}}\otimes\pi_{\man{E}}}}
\dnorm{D^{k-j}_{\nabla^{\man{M}},\nabla^{\pi_{\man{E}}}}
\xi(x)}_{\metric_{\man{M},\pi_{\man{E}}}}
\end{equation*}
for $k\in\integerp$\@.  By~\eqref{eq:12inftynorms} (twice) we have
\begin{align*}
\dnorm{j_m(L\scirc\xi)(x)}_{\metric_{\man{M},\pi_{\man{F}},m}}\le&\;
\sum_{k=0}^m\frac{1}{k!}\dnorm{D^k_{\nabla^{\man{M}},\nabla^{\pi_{\man{F}}}}
(L\scirc\xi(x))}_{\metric_{\man{M},\pi_{\man{F}}}}\\
\le&\;\sum_{k=0}^m\frac{1}{k!}\sum_{j=0}^k\binom{k}{j}
\dnorm{D^j_{\nabla^{\man{M}},\nabla^{\pi_{\man{F}}\otimes\pi_{\man{E}}}}
L(x)}_{\metric_{\man{M},\pi_{\man{F}}\otimes\pi_{\man{E}}}}
\dnorm{D^{k-j}_{\nabla^{\man{M}},\nabla^{\pi_{\man{E}}}}
\xi(x)}_{\metric_{\man{M},\pi_{\man{E}}}}\\
=&\;\sum_{k=0}^m\sum_{j=0}^k
\frac{\dnorm{D^j_{\nabla^{\man{M}},\nabla^{\pi_{\man{F}}\otimes\pi_{\man{E}}}}
L(x)}_{\metric_{\man{M},\pi_{\man{F}}\otimes\pi_{\man{E}}}}}{j!}
\frac{\dnorm{D^{k-j}_{\nabla^{\man{M}},\nabla^{\pi_{\man{E}}}}
\xi(x)}_{\metric_{\man{M},\pi_{\man{E}}}}}{(k-j)!}\\
\le&\;(m+1)^2\sup\asetdef{
\frac{\dnorm{D^j_{\nabla^{\man{M}},\nabla^{\pi_{\man{F}}\otimes\pi_{\man{E}}}}
L(x)}_{\metric_{\man{M},\pi_{\man{F}}\otimes\pi_{\man{E}}}}}{j!}}{j\le m}\\
&\;\times\sup\asetdef{\frac{
\dnorm{D^{k-j}_{\nabla^{\man{M}},\nabla^{\pi_{\man{E}}}}
\xi(x)}_{\metric_{\man{M},\pi_{\man{E}}}}}{(k-j)!}}{j\le m}\\
\le&\;(m+1)^{5/2}\dnorm{j_mL(x)}_{\metric_{\man{M},\pi_{\man{F}}
\otimes\pi_{\man{E}},m}}\dnorm{j_m\xi(x)}_{\metric_{\man{M},\pi_{\man{E}},m}}.
\end{align*}
Noting that $(m+1)^{5/2}\le3^{m+1}$\@, $m\in\integerp$\@, we finally get
\begin{equation}\label{eq:normjmLxi}
\dnorm{j_m(L\scirc\xi)(x)}_{\metric_{\man{M},\pi_{\man{F}},m}}\le
3^{m+1}\dnorm{j_mL(x)}_{\metric_{\man{M},\pi_{\man{F}}\otimes\pi_{\man{E}},m}}
\dnorm{j_m\xi(x)}_{\metric_{\man{M},\pi_{\man{E}},m}}.
\end{equation}
Let $\nbhd{K}\subset\man{M}$ be compact and let
$\vect{a}\in\c_0(\integernn;\realp)$\@.  Define define
$\vect{a}'\in\c_0(\integernn;\realp)$ by $a'_j=\sqrt{3}a_j$\@,
$j\in\integernn$\@.  We then have
\begin{equation*}
q^\omega_{\nbhd{K},\vect{a}}(L\scirc\xi)\le
Cq^\omega_{\nbhd{K},\vect{a}'}\otimes p^\omega_{\nbhd{K},\vect{a}'}(L)
p^\omega_{\nbhd{K},\vect{a}'}(\xi)=C(q^\omega_{\nbhd{K},\vect{a}'}\otimes p^\omega_{\nbhd{K},\vect{a}'})\otimes p^\omega_{\nbhd{K},\vect{a}'}(L\otimes\xi).
\end{equation*}
By~\cite[Theorem~15.1.2]{HJ:81}\@, this gives continuity of the bilinear map
$(L,\xi)\mapsto L\scirc\xi$\@.

\eqref{pl:xopen} We first prove a couple of technical lemmata.
\begin{prooflemma}\label{plem:open-mapping}
Let\/ $\alg{U}$ and\/ $\alg{V}$ be locally convex topological vector spaces,
and let\/ $L\in\lin{\alg{U}}{\alg{V}}$\@.  If, for every continuous
seminorm\/ $q$ for\/ $\alg{U}$\@, there exists a continuous seminorm\/ $p$
for\/ $\alg{V}$ such that
\begin{equation*}
q(u)\le p(L(u)),\qquad u\in\alg{U},
\end{equation*}
then\/ $L$ is an open mapping onto its image.
\begin{subproof}
First we prove that there are $0$-bases $\sB_{\alg{U}}$ for $\nbhd{U}$ and
$\sB_{\alg{V}}$ for $\alg{V}$ such that, for each
$\nbhd{B}\in\sB_{\alg{U}}$\@, there exists $\nbhd{C}\in\sB_{\alg{V}}$ such
that
\begin{equation*}
\nbhd{C}\cap\image(L)\subset L(\nbhd{B}).
\end{equation*}
To see this, first let $q$ be a continuous seminorm for $\alg{U}$ and let $p$
be a continuous seminorm for $\alg{V}$ satisfying
\begin{equation*}
q(u)\le p(L(u)),\qquad u\in\alg{U}.
\end{equation*}
Then
\begin{equation*}
p(L(u))<1\ \implies\ q(u)<1\ \implies\ L(u)\in L(q^{-1}(\interval[0,1)).
\end{equation*}
Thus
\begin{equation*}
p^{-1}(\interval[0,1))\cap\image(L)\subset L(q^{-1}(\interval[0,1))).
\end{equation*}
Now let $\sB_{\alg{U}}$ be the collection of all $0$-neighbourhoods of the
form
\begin{equation*}
\nbhd{B}=\cap_{j=1}^kq_j^{-1}(\interval[0,1)),\qquad k\in\integerp,\
q_j\ \textrm{a continuous seminorm}, j\in\{1,\dots,k\}.
\end{equation*}
This is a $0$-base for $\alg{U}$\@.  For each such $\nbhd{B}$\@, we let $p_j$
be continuous seminorm for $\alg{V}$ corresponding to $q_j$ by
\begin{equation*}
q_j(u)\le p_j(L(u)),\qquad u\in\alg{U},\ j\in\{1,\dots,k\}.
\end{equation*}
Then, by our above computations,
\begin{equation*}
\left(\cap_{j=1}^kp_j^{-1}(\interval[0,1))\right)\cap\image(L)\subset
L\left(\cap_{j=1}^kq_j^{-1}(\interval[0,1))\right).
\end{equation*}
Thus, the $0$-base
\begin{equation*}
\cap_{j=1}^kp_j^{-1}(\interval[0,1)),\qquad k\in\integerp,
p_j\ \textrm{a continuous seminorm}, j\in\{1,\dots,k\},
\end{equation*}
for $\alg{V}$ has the desired property.

Now let $\nbhd{O}\subset\alg{V}$ be open and let $u\in\nbhd{O}$\@.  Let
$\nbhd{B}\in\sB_{\alg{U}}$ be such that $u+\nbhd{B}\subset\nbhd{O}$ and let
$\nbhd{C}\in\sB_{\alg{V}}$ be such that $\nbhd{C}\cap\image(L)\subset
L(\nbhd{B})$\@.  Then
\begin{equation*}
L(u)+\nbhd{C}\cap\image(L)\subset L(u)+L(\nbhd{B})=L(u+\nbhd{B})\subset L(\nbhd{O}).
\end{equation*}
Thus $L(u)+\nbhd{C}\cap\image(L)$ is a neighbourhood of $L(u)$ in
$L(\nbhd{O})$ which shows that $L(\nbhd{O})$ is open in $\image(L)$\@.
\end{subproof}
\end{prooflemma}

\begin{prooflemma}\label{plem:vb-left-inverse}
If\/ $L$ is injective, then there exists a left-inverse\/
$L'\in\sections[\omega]{\man{E}\otimes\dual{\man{F}}}$\@.
\begin{subproof}
First we note that $\image(L)$ is a $\C^\omega$-subbundle of $\man{F}$\@ and
that $L$ is a $\C^\omega$-vector bundle isomorphism onto $\image(L)$\@.  Let
$\alg{G}\subset\alg{F}$ be the $\metric_{\pi_{\man{E}}}$-orthogonal
complement to $\image(L)$ which is then itself a $\C^\omega$-subbundle of
$\man{F}$\@.  Clearly, $\man{F}=\image(L)\oplus\man{G}$\@.  Let
\begin{equation*}
\mapdef{L'}{\image(L)\oplus\man{G}}{\man{E}}{(L(e),g)}{e,}
\end{equation*}
and note that $L'$ is obviously a left-inverse of $L$\@.  It is also of class
$\C^\omega$ since the projection from $\man{F}$ to the summand $\image(L)$ is
of class $\C^\omega$\@.
\end{subproof}
\end{prooflemma}

By Lemma~\ref{plem:vb-left-inverse} we suppose that there is a
$\C^\omega$-vector bundle mapping $L'$ that is a left-inverse for $L$\@.
Then, from the first part of the proof, for a compact
$\nbhd{K}\subset\man{M}$ and for $\vect{a}\in\c_0(\integernn;\realp)$\@, let
$C\in\realp$ be such that
\begin{equation*}
p^\omega_{\nbhd{K},\vect{a}}(L'\scirc\eta)\le Cq^\omega_{\nbhd{K},\vect{a}}(\eta),\qquad\eta\in\sections[\omega]{\man{F}}.
\end{equation*}
We then have, for $\xi\in\sections[\omega]{\man{E}}$\@,
\begin{equation*}
p_{\nbhd{K},\vect{a}}^\omega(\xi)=p_{\nbhd{K},\vect{a}}^\omega(L'\scirc L\scirc\xi)
\le Cq^\omega_{\nbhd{K},\vect{a}}(L\scirc\xi).
\end{equation*}
By Lemma~\ref{plem:open-mapping}\@, this suffices to establish that $L$ is
open onto its image.
\end{proof}
\end{theorem}

\begin{remark}
The preceding proof works equally well in the smooth case.  Indeed, the proof
is a little easier since one does not need to carefully keep track of the
growth in $m$ of the coefficient of the norm of the $m$-jet.\oprocend
\end{remark}

The following result is an important one, and is very much nontrivial in the
real analytic case.  It is established during the course of the proof of
their Lemma~2.5 by \citet{SJ/ADL:14a} using a local description of the real
analytic topology.  Here we use an intrinsic proof.
\begin{theorem}\label{the:compcont}
Let\/ $\man{M}$ and\/ $\man{N}$ be\/ $\C^\omega$-manifolds.  If\/
$\Phi\in\mappings[\omega]{\man{M}}{\man{N}}$\@, then the mapping
\begin{equation*}
\mapdef{\Phi^*}{\func[\omega]{\man{N}}}
{\func[\omega]{\man{M}}}{f}{f\scirc\Phi}
\end{equation*}
is continuous.  Moreover, if\/ $\Phi$ is a proper surjective submersion or a
proper embedding, then\/ $\Phi^*$ is open onto its image.  In
case\/ $\Phi$ is a proper embedding, for any compact\/
$\nbhd{K}\subset\man{M}$ and any\/ $\vect{a}\in\c_0(\integernn;\realp)$\@,
there exists\/ $\vect{a}'\in\c_0(\integernn;\realp)$ such that
\begin{equation*}
q^\omega_{\Phi(\nbhd{K}),\vect{a}}(f)\le
p^\omega_{\nbhd{K},\vect{a}'}(\Phi^*f),\qquad f\in\func[\omega]{\man{N}}.
\end{equation*}
\begin{proof}
We let $\nabla^{\man{M}}$ and $\nabla^{\man{N}}$ be $\C^\omega$-affine
connections on $\man{M}$ and $\man{N}$\@, respectively, and let
$\metric_{\man{M}}$ and $\metric_{\man{N}}$ be $\C^\omega$-Riemannian metrics
on $\man{M}$ and $\man{N}$\@, respectively.  For $\nbhd{K}\subset\man{M}$ and
$\nbhd{L}\subset\man{N}$ compact, and for
$\vect{a}\in\c_0(\integernn;\realp)$\@, we denote by
$p^\omega_{\nbhd{K},\vect{a}}$ and $q^\omega_{\nbhd{L},\vect{a}}$ the associated
seminorms for $\func[\omega]{\man{M}}$ and $\func[\omega]{\man{N}}$\@,
respectively.

From Lemma~\ref{lem:pbsymiterateddersI} we have the formula
\begin{equation}\label{eq:pbSMSN}
\Sym_m\scirc\nabla^{\man{M},m}\Phi^*f=
\sum_{s=0}^m\what{A}^m_s(\Sym_s\scirc\Phi^*\nabla^{\man{N},s}f).
\end{equation}
By Lemma~\ref{lem:A(S)<AS}\@, we have
\begin{equation*}
\dnorm{A^m_s(\beta_s)}_{\metric_{\man{M}}}\le
\dnorm{A^m_s}_{\metric_{\man{M}},\metric_{\man{N}}}
\dnorm{\beta_s}_{\metric_{\man{N}}}
\end{equation*}
for $\beta_s\in\tensor*[s]{\ctb[x]{\man{M}}}$\@, $m\in\integerp$\@, and
$s\in\{0,1,\dots,m\}$\@.  By Lemma~\ref{lem:Symnorm}\@,
\begin{equation*}
\dnorm{\Sym_s(A)}_{\metric_{\man{M}},\metric_{\man{N}}}\le
\dnorm{A}_{\metric_{\man{M}},\metric_{\man{N}}}
\end{equation*}
for $A\in\tensor*[s]{\ctb{\man{N}}}$ and $s\in\integerp$\@.  Thus,
recalling~\eqref{eq:AhatA} (and its analogue that would arise in a spelled
out proof of Lemma~\ref{lem:pbsymiterateddersI}),
\begin{equation*}
\dnorm{\what{A}^m_s(\Sym_s(\beta_s))}_{\metric_{\man{M}}}=
\dnorm{\Sym_m\scirc A^m_s(\beta_s)}_{\metric_{\man{M}}}
\le\dnorm{A^m_s}_{\metric_{\man{M}},\metric_{\man{N}}}
\dnorm{\beta_s}_{\metric_{\man{M}}},
\end{equation*}
for $\beta_s\in\tensor*[s]{\Phi^*\ctb{\man{N}}}$\@, $m\in\integernn$\@,
$s\in\{1,\dots,m\}$\@.

Let $\nbhd{K}\subset\man{M}$ be compact.  By Lemmata~\ref{lem:pissynabla}
and~\ref{lem:particular-bounds} with $r=0$\@, there exist
$C_1,\sigma_1,\rho_1\in\realp$ such that
\begin{equation*}
\dnorm{A^k_s(x)}_{\metric_{\man{M}},\metric_{\man{N}}}\le
C_1\sigma_1^{-k}\rho_1^{-(k-s)}(k-s)!,\qquad
k\in\integernn,\ s\in\{0,1,\dots,k\},\ x\in\nbhd{K}.
\end{equation*}
By Lemma~\ref{lem:pbnormbase}\@, let $C_2\in\realp$ be such that
\begin{equation*}
\dnorm{\Phi^*\nabla^{\man{N},m}f(x)}_{\metric_{\man{M}}}\le
C_2^m\dnorm{\nabla^{\man{N},m}f(\Phi(x))}_{\metric_{\man{N}}},\qquad
x\in\nbhd{K},\ m\in\integernn.
\end{equation*}
Without loss of generality, we assume that $C_1,C_2\ge1$ and
$\sigma_1,\rho_1\le1$\@.  Thus, abbreviating $\sigma_2=\sigma_1\rho_1$\@, we
have
\begin{equation*}
\dnorm{\what{A}^k_s(\Phi^*\Sym_s\scirc
\nabla^{\man{N},s}f(x))}_{\metric_{\man{M}}}\le C_1C_2^s\sigma_2^{-k}(k-s)!
\dnorm{\Sym_s\scirc\nabla^{\man{N},s}f(\Phi(x))}_{\metric_{\man{N}}}
\end{equation*}
for $k\in\integernn$\@, $s\in\{0,1,\dots,k\}$\@, $x\in\nbhd{K}$\@.  Thus,
by~\eqref{eq:12inftynorms} and~\eqref{eq:pbSMSN}\@,
\begin{align*}
\dlnorm j_m(\Phi^*f)(x)\drnorm_{\metric_{\man{M},m}}\le&\;
\sum_{k=0}^m\frac{1}{k!}\dnorm{\Sym_k\scirc
\nabla^{\man{M},k}\Phi^*f(x)}_{\metric_{\man{M}}}\\
=&\;\sum_{k=0}^m\frac{1}{k!}\adnorm{\sum_{s=0}^k
\what{A}^k_s(\Phi^*\Sym_s\scirc\nabla^{\man{N},s}
f(\Phi(x)))}_{\metric_{\man{M}}}\\
\le&\;\sum_{k=0}^m\sum_{s=0}^kC_1\sigma_2^{-k}
\frac{s!(k-s)!}{k!}\frac{C_2^s}{s!}
\dnorm{\Sym_s\scirc\nabla^{\man{N},s}f(\Phi(x))}_{\metric_{\man{N}}}
\end{align*}
for $x\in\nbhd{K}$ and $m\in\integernn$\@.  Now note that
\begin{equation*}
\frac{s!(k-s)!}{k!}\le1,\quad
C_1\sigma_2^{-k}C_2^s\le C_1C_2^m\sigma_2^{-m},
\end{equation*}
for $s\in\{0,1,\dots,m\}$\@, $k\in\{0,1,\dots,s\}$\@, since $\sigma_2\le1$\@.
Then
\begin{align*}
\dlnorm j_m(\Phi^*f)(x)\drnorm_{\metric_{\man{M},m}}
\le&\;C_1C_2^m\sigma_2^{-m}\sum_{k=0}^m\sum_{s=0}^k\frac{1}{s!}
\dnorm{\Sym_s\scirc\nabla^{\man{N},s}f(\Phi(x))}_{\metric_{\man{N}}}\\
\le&\;C_1C_2^m\sigma_2^{-m}\sum_{k=0}^m\sum_{s=0}^m
\frac{1}{s!}\dnorm{\Sym_s\scirc\nabla^{\man{N},s}
f(\Phi(x))}_{\metric_{\man{N}}}\\
=&\;(m+1)C_1C_2^m\sigma_2^{-m}\sum_{s=0}^m
\frac{1}{s!}\dnorm{\Sym_s\scirc\nabla^{\man{N},s}f(\Phi(x))}_{\metric_{\man{N}}}.
\end{align*}
Now let $\sigma<C_2^{-1}\sigma_2$ and note that
\begin{equation*}
\lim_{m\to\infty}(m+1)\frac{C_2^m\sigma_2^{-m}}{\sigma^{-m}}=0.
\end{equation*}
Thus there exists $N\in\integerp$ such that
\begin{equation*}
(m+1)C_1C_2^m\sigma_2^{-m}\le C_1\sigma^{-m},\qquad m\ge N.
\end{equation*}
Let
\begin{equation*}
C=\max\left\{C_1,2C_1C_2\frac{\sigma}{\sigma_2},
3C_1C_2^2\left(\frac{\sigma}{\sigma_2}\right)^2,
\dots,(N+1)C_1C_2^N\left(\frac{\sigma}{\sigma_2}\right)^N\right\}.
\end{equation*}
We then immediately have $(m+1)C_1C_2^m\sigma_2^{-m}\le C\sigma^{-m}$ for all
$m\in\integernn$\@.  We then have, using~\eqref{eq:12inftynorms}\@,
\begin{align*}
\dnorm{j_m(\Phi^*f)(x)}_{\metric_{\man{M},m}}\le&\;
C\sigma^{-m}\left(\sum_{s=0}^m\frac{1}{s!}
\dnorm{\Sym_s\scirc\nabla^{\man{N},s}
f(\Phi(x))}_{\metric_{\man{N}}}\right)\\
=&\;C\sqrt{m+1}\sigma^{-m}\dnorm{j_mf(\Phi(x))}_{\metric_{\man{N},m}}.
\end{align*}
By modifying $C$ and $\sigma$ guided by what we did just preceding, we get
\begin{equation*}
\dnorm{j_m(\Phi^*f)(x)}_{\metric_{\man{M},m}}\le
C\sigma^{-m}\dnorm{j_mf(\Phi(x))}_{\metric_{\man{N},m}}.
\end{equation*}
Now, for $\vect{a}\in\c_0(\integernn;\realp)$\@, let
$\vect{a}'\in\c_0(\integernn;\realp)$ be defined by $a'_0=Ca_0$ and
$a_j'=a_j\sigma^{-1}$\@, $j\in\integerp$\@.  Then we have
\begin{equation*}
p^\omega_{\nbhd{K},\vect{a}}(\Phi^*f)\le q^\omega_{\Phi(\nbhd{K}),\vect{a}'}(f),
\end{equation*}
and this gives continuity of $\Phi^*$\@.

Now we turn to the final assertion concerning the openness of $\Phi^*$ in
particular cases.  First we note that, by
Lemma~\ref{lem:pbsymiterateddersII}\@, we have
\begin{equation*}
\Sym_m\scirc\Phi^*\nabla^{\man{N},m}f(x)=
\sum_{s=0}^m\what{B}^m_s(\Sym_s\scirc
\nabla^{\man{M},s}\Phi^*f(x)).
\end{equation*}

First consider the case where $\Phi$ is a proper surjective submersion.  For
$\nbhd{L}\subset\man{N}$ compact and for $y\in\nbhd{L}$\@, since $\Phi$ is
surjective, there exists $x\in\man{M}$ such that $\Phi(x)=y$\@.  Also, since
$\Phi$ is proper, $\Phi^{-1}(\nbhd{L})$ is compact.  We can now reproduce the
steps from the proof above, now making use of the second part of
Lemma~\ref{lem:pbnormbase}\@, to prove that
\begin{equation*}
q^\omega_{\nbhd{L},\vect{a}}(f)\le
p^\omega_{\Phi^{-1}(\nbhd{L}),\vect{a}'}(\Phi^*f),
\end{equation*}
which suffices to prove the openness of $\Phi^*$ by
Lemma~\ref{plem:open-mapping} from the proof of Theorem~\ref{the:+xcont}\@.

Finally consider the case where $\Phi$ is a proper embedding.  Here we make
use of a lemma.
\begin{prooflemma}\label{plem:PB-embedding}
Let\/ $\man{M}$ be a\/ $\C^\omega$-manifold and let\/ $\man{S}\subset\man{M}$
be a\/ $\C^\omega$-embedded submanifold with\/
$\map{\iota_{\man{S}}}{\man{S}}{\man{M}}$ the inclusion.  Then
\begin{equation*}
\map{\iota_{\man{S}}^*}{\func[\omega]{\man{M}}}{\func[\omega]{\man{S}}}
\end{equation*}
is an epimorphism,~\ie~continuous, surjective, and open.
\begin{subproof}
First note that we can use Sublemma~\ref{psublem:Sextend} from the proof of
Lemma~\ref{lem:bundleform} to show that
$\map{\iota_{\man{S}}^*}{\func[\omega]{\man{M}}}{\func[\omega]{\man{S}}}$ is
surjective.  It, therefore, remains to show that $\iota^*_{\man{S}}$ is
continuous and open.  Continuity follows from Theorem~\ref{the:compcont}\@.
Since $\func[\omega]{\man{S}}$ and $\func[\omega]{\man{M}}$ are
ultrabornological webbed spaces, the De Wilde Open Mapping
Theorem~\cite[Theorem~24.30]{RM/DV:97} implies that $\iota^*_{\man{S}}$ is
open.
\end{subproof}
\end{prooflemma}

The lemma immediately gives openness of $\Phi^*$ in the case that $\Phi$ is a
proper embedding.  For the final assertion, we can follow the same argument
as was sketched for the openness of $\Phi^*$ when $\Phi$ is a proper
surjective submersion to give
\begin{equation*}
q^\omega_{\Phi(\nbhd{K}),\vect{a}}(f)\le
p^\omega_{\nbhd{K},\vect{a}'}(\Phi^*f),
\end{equation*}
as desired.
\end{proof}
\end{theorem}

The matter of determining general conditions under which $\Phi^*$ is an
homeomorphism onto its image or has closed image are taken up by
\citeauthor{PD/ML:03}~\cite*{PD/ML:03,PD/ML:06}\@.  The linear operator
$f\mapsto\Phi^*f$ is often called a composition operator.  Also of interest
is the nonlinear operator $\Phi\mapsto f\scirc\Phi$\@, which is variously
called a ``superposition operator,'' a ``nonlinear composition operator,'' or
the ``Nemytskii operator.''  Both operators are of substantial interest in
various areas of mathematics.  In our setting of real analytic analysis, we
work with the nonlinear operator in Section~\ref{subsec:nemytskii}\@.

\begin{remark}
The preceding proof can be adapted to the smooth case.  Indeed, much of the
elaborate work of the proof can be simplified by not having to pay attention
to the exponential growth of $m$-jet norms as $m\to\infty$\@.  In the smooth
case, one works with fixed orders of derivatives.  This comment applies to
all of our subsequent proofs in this section.\oprocend
\end{remark}

\subsection{Continuity of operations involving differentiation}

Next we consider a general version of the assertion that ``differentiation is continuous.''
\begin{theorem}\label{the:jetcont}
Let\/ $\map{\pi_{\man{E}}}{\man{E}}{\man{M}}$ be a\/ $\C^\omega$-vector
bundle.  If\/ $k\in\integernn$\@, then the map
\begin{equation*}
\mapdef{J_k^\omega}{\sections[\omega]{\man{E}}}
{\sections[\omega]{\jet{k}{\man{E}}}}{\xi}{j_k\xi}
\end{equation*}
is continuous.
\begin{proof}
We let $\nabla^{\man{M}}$ be a $\C^\omega$-affine connection on $\man{M}$\@,
$\nabla^\pi$ be a $\C^\omega$-vector bundle connection in $\man{E}$\@,
$\metric_{\man{M}}$ be a $\C^\omega$-Riemannian metric on $\man{M}$\@, and
$\metric_\pi$ be a $\C^\omega$-vector bundle connection in $\man{E}$\@.  We
denote the associated seminorms for $\sections[\omega]{\man{E}}$ by
$p^\omega_{\nbhd{K},\vect{a}}$ and for $\sections[\omega]{\jet{k}{\man{E}}}$
by $p^{k,\omega}_{\nbhd{K},\vect{a}}$\@, for $\nbhd{K}\subset\man{M}$ compact
and $\vect{a}\in\c_0(\integernn;\realp)$\@.

We recall from Section~\ref{subsec:prolongation} that we have the vector
bundle mapping
$\map{\pi_{m,k}}{\jet{k+m}{\man{E}}}{\jet{m}{\jet{k}{\man{E}}}}$ defined by
the requirement that $\pi_{m,k}\scirc j_{k+m}\xi=j_mj_k\xi$\@.  We begin the
proof by doing some computations that give norm estimates for this vector
bundle map.  To do this, we use the representation $\what{\Delta}_{m,k}$ of
$\pi_{m,k}$ relative to our decompositions of jet bundles, as in
Lemma~\ref{lem:Jk+mdecomp}\@.  We let $x\in\man{M}$ and let
$A_j\in\Symalg*[j]{\ctb[x]{\man{M}}}\otimes\man{E}_x$\@,
$j\in\{0,1,\dots,k+m\}$\@, and compute, using Lemmata~\ref{lem:Deltasym}
and~\ref{lem:Symnorm}\@,~\eqref{eq:Delathatrs}\@,
and~\eqref{eq:12inftynorms}\@,
\begin{align*}
\dnorm{\what{\Delta}_{m,k}\pi_{\man{E}}
(A_0,A_1,\dots,A_{m+k})}_{\metric_{\man{M},\pi_{\man{E}}}}
\le&\;\sum_{l=0}^m\frac{1}{l!}\sum_{j=0}^k
\dnorm{\Delta_{j,l}(A_{j+l})}_{\metric_{\man{M},\pi_{\man{E}}}}\\
\le&\;\sum_{l=0}^m\frac{(k+l)!}{l!}\sum_{j=0}^k\frac{1}{(j+l)!}
\dnorm{A_{j+l}}_{\metric_{\man{M},\pi_{\man{E}}}}\\
\le&\;\frac{(k+m)!}{m!}\sum_{l=0}^m\sum_{j=0}^k\frac{1}{(j+l)!}
\dnorm{A_{j+l}}_{\metric_{\man{M},\pi_{\man{E}}}}\\
\le&\;(m+k)^k(m+1)\sum_{j=0}^{k+m}\frac{1}{j!}
\dnorm{A_j}_{\metric_{\man{M},\pi_{\man{E}}}}.
\end{align*}
For $\sigma\in\interval(0,1)$\@,
\begin{equation*}
\lim_{m\to0}\sigma^{-m}(m+k)^k(m+1)=0.
\end{equation*}
Thus let $N\in\integerp$ be such that
\begin{equation*}
\sigma^{-m}(m+k)^k(m+1)<0,\qquad m\ge N.
\end{equation*}
Next let
\begin{equation*}
C=\max\left\{k^k,\frac{2(1+k)^k}{\sigma},\dots,
\frac{(N+1)(N+k)^k}{\sigma^N}\right\}.
\end{equation*}
Then, for any $m\in\integernn$\@,
\begin{equation*}
(m+k)^k(m+1)\le C\sigma^{-m},
\end{equation*}
and so, using~\eqref{eq:12inftynorms}\@,
\begin{equation*}
\dnorm{\what{\Delta}_{m,k}\pi_{\man{E}}
(A_0,A_1,\dots,A_{m+k})}_{\metric_{\man{M},\pi_{\man{E}}}}
\le\sqrt{k+m+1}C\sigma^{-m}
\dnorm{(A_0,A_1,\dots,A_{m+k})}_{\metric_{\man{M},\pi_{\man{E}}}}.
\end{equation*}
Modifying $C$ and $\sigma$ similarly to our constructions above shows that
\begin{equation*}
\dnorm{j_mj_k\xi(x)}_{\metric_{\man{M},\pi_k,m}}\le
C\sigma^{-m}\dnorm{j_{k+m}\xi(x)}_{\metric_{\man{M},\pi_{\man{E}},k+m}},
\qquad x\in\man{M}.
\end{equation*}

Let $\nbhd{K}\subset\man{M}$ be compact and let
$\vect{a}\in\c_0(\integernn;\realp)$\@.  Define
$\vect{a}'\in\c_0(\integernn;\realp)$ by $a'_0=a_0$\@, $a'_j=C$\@,
$j\in\{1,\dots,k\}$\@, and $a'_j=\sigma^{-1}a_{j-k}$\@,
$j\in\{k+1,k+2,\dots\}$\@.  The computations from the beginning of the proof
then give
\begin{align*}
a_0a_1\cdots a_m\dnorm{j_mj_k\xi(x)}\le&\;
C\sigma^{-m}a_0a_1\cdots a_m\dnorm{j_{k+m}\xi(x)}\\
\le&\;a_0C^k(\sigma^{-1}a_1)\cdots(\sigma^{-1}a_m)\dnorm{j_{k+m}\xi(x)}\\
=&\;a'_0a'_1\cdots a'_{k+m}\dnorm{j_{k+m}\xi(x)},
\end{align*}
since $C\ge1$\@.  We then immediately have
\begin{equation*}
p^{k,\omega}_{\nbhd{K},\vect{a}}(j_k\xi)\le p^\omega_{\nbhd{K},\vect{a}'}(\xi),
\end{equation*}
which gives the theorem.
\end{proof}
\end{theorem}

We can now prove a collection of results regarding standard operations of
differentiation, derived from the preceding result about basic prolongation.
\begin{corollary}\label{cor:dfcont}
Let\/ $\man{M}$ be a\/ $\C^\omega$-manifold.  Then the mapping
\begin{equation*}
\mapdef{\d}{\func[\omega]{\man{M}}}
{\sections[\omega]{\ctb{\man{M}}}}{f}{\d{f}}
\end{equation*}
is continuous.
\begin{proof}
Note that $\jet{1}{(\man{M};\real)}\simeq\real_{\man{M}}\oplus\ctb{\man{M}}$
and that, under this identification, $j_1f=f\oplus\d{f}$\@.  Thus
$\d{f}=\pr_2\scirc j_1f$\@, where
$\map{\pr_2}{\jet{1}{(\man{M};\real)}}{\ctb{\man{M}}}$ is the
$\C^\omega$-vector bundle mapping of projection onto the second factor.  The
result then immediately follows from Lemma~\pldblref{the:+xcont}{pl:xcont}
and Theorem~\ref{the:jetcont}\@.
\end{proof}
\end{corollary}

\begin{corollary}\label{cor:LXcont}
Let\/ $\man{M}$ be a\/ $\C^\omega$-manifold.  Then the map
\begin{equation*}
\mapdef{\lieder{}{}}{\sections[\omega]{\tb{\man{M}}}\times
\func[\omega]{\man{M}}}
{\func[\omega]{\man{M}}}{(X,f)}{\lieder{X}{f}}
\end{equation*}
is continuous.
\begin{proof}
We think of $X$ as being a $\C^\omega$-vector bundle mapping via
\begin{equation*}
\mapdef{X}{\ctb{\man{M}}}{\real_{\man{M}}}{\alpha_x}
{\natpair{\alpha_x}{X(x)}.}
\end{equation*}
Then the bilinear mapping of the lemma is given by the composition
\begin{equation*}
(X,f)\mapsto(X,\d{f})\mapsto X(\d{f}).
\end{equation*}
The left mapping is continuous since it is the product of the continuous
mappings $\id$ and $\d$\@.  The right mapping is continuous by
Theorem~\pldblref{the:+xcont}{pl:xcont}\@, and so the corollary follows.
\end{proof}
\end{corollary}

\begin{corollary}\label{cor:nablaXcont}
Let\/ $\map{\pi_{\man{E}}}{\man{E}}{\man{M}}$ be a\/ $\C^\omega$-vector
bundle with a\/ $\C^\omega$-vector bundle connection\/
$\nabla^{\pi_{\man{E}}}$\@.  Then the map
\begin{equation*}
\mapdef{\nabla^{\pi_{\man{E}}}}{\sections[\omega]{\tb{\man{M}}}\times
\sections[\omega]{\man{E}}}
{\sections[\omega]{\man{E}}}{(X,\xi)}{\nabla^{\pi_{\man{E}}}_X\xi}
\end{equation*}
is continuous.
\begin{proof}
As in the proof of Lemma~\ref{lem:Jmdecomp}\@, we have a $\C^\omega$-vector
bundle mapping $\map{S_{\nabla^{\pi_{\man{E}}}}}{\man{E}}{\jet{1}{\man{E}}}$
over $\id_{\man{M}}$ that determines the connection $\nabla^{\pi_{\man{E}}}$
by
\begin{equation*}
\nabla^{\pi_{\man{E}}}\xi(x)=j_1\xi(x)-S_{\nabla^{\pi_{\man{E}}}}(\xi(x)).
\end{equation*}
The mapping $\xi\mapsto\nabla^{\pi_{\man{E}}}\xi$ is continuous by
Theorems~\ref{the:jetcont} and~\ref{the:+xcont}\@.  We note that
$\nabla^{\pi_{\man{E}}}\xi$ is to be thought of as a $\C^\omega$-vector bundle
mapping by
\begin{equation*}
\mapdef{\nabla^{\pi_{\man{E}}}\xi}{\tb{\man{M}}}{\man{E}}{X}
{\nabla^{\pi_{\man{E}}}_X\xi.}
\end{equation*}
The bilinear mapping of the lemma is then given by the composition
\begin{equation*}
(X,\xi)\mapsto(X,\nabla^{\pi_{\man{E}}}\xi)\mapsto\nabla^{\pi_{\man{E}}}\xi(X).
\end{equation*}
The left mapping is continuous since it is the product of the continuous
mappings $\id$ and $\xi\mapsto\nabla^{\pi_{\man{E}}}\xi$\@.  The right
mapping is continuous by Theorem~\pldblref{the:+xcont}{pl:xcont}\@, and so
the lemma follows.
\end{proof}
\end{corollary}

\begin{corollary}\label{cor:[X,Y]cont}
Let\/ $\man{M}$ be a\/ $\C^\omega$-manifold.  Then the map
\begin{equation*}
\mapdef{{[\cdot,\cdot]}}{\sections[\omega]{\tb{\man{M}}}\times
\sections[\omega]{\tb{\man{M}}}}
{\sections[\omega]{\tb{\man{M}}}}{(X,Y)}{[X,Y]}
\end{equation*}
is continuous.
\begin{proof}
Let $\metric_{\man{M}}$ be a real analytic Riemannian metric on $\man{M}$ and
let $\nabla^{\man{M}}$ be the associated Levi-Civita connection.  Since
\begin{equation*}
[X,Y]=\nabla^{\man{M}}_XY-\nabla^{\man{M}}_XY,
\end{equation*}
the result follows from Corollary~\ref{cor:nablaXcont}\@.
\end{proof}
\end{corollary}

\begin{corollary}
Let\/ $\map{\pi_{\man{E}}}{\man{E}}{\man{M}}$ and\/
$\map{\pi_{\man{F}}}{\man{F}}{\man{M}}$ be\/ $\C^\omega$-vector bundles and
let\/ $\Phi\in\vbmappings[\omega]{\jet{k}{\man{E}}}{\man{F}}$\@.  Then the\/
$k$th-order linear partial differential operator\/
$\map{D_\Phi}{\sections[\omega]{\man{E}}}{\sections[\omega]{\man{F}}}$
defined by\/ $D_\Phi(\xi)(x)=\Phi(j_k\xi(x))$\@,\/ $x\in\man{M}$\@, is
continuous.
\begin{proof}
The operator $D_\Phi$ is the composition of the continuous mappings
$\xi\mapsto j_k\xi$ $\sections[\omega]{\man{E}}$ to $\sections[\omega]{\jet{k}{\man{E}}}$ and $\Xi\mapsto\Phi\scirc\Xi$ from
$\sections[\omega]{\jet{k}{\man{E}}}$ to $\sections[\omega]{\man{F}}$\@.
\end{proof}
\end{corollary}

The reader can no doubt imagine many extensions of results such as the ones
we give, and we leave these for the reader to figure out as they need them.

\subsection{Continuity of lifting operations}\label{subsec:liftcont}

In Sections~\ref{subsec:vbfunctions}--\ref{subsec:vbtensors} we introduced a
variety of constructions for lifting objects from the base space of a vector
bundle to the total space.  In Section~\ref{sec:lift-isomorphisms} we
considered how to differentiate these constructions in multiple ways, and how
relate these multiple differentiations.  In
Sections~\ref{subsec:Pjetnorm}--\ref{subsec:Cjetnorm} we described fibre
norms to give norms for these lifted objects.  In this section, we put this
all together to prove results that are the entire \emph{raison d'\^etre} for
all of these constructions, some of them quite elaborate.  That is, we show
that these lift operations are homeomorphisms onto their images.  Many of the
proofs are similar to one another, so we only give representative proofs.

We begin by considering horizontal lifts of functions.  We note that
continuity of the mapping in the next theorem follows from
Theorem~\ref{the:compcont}\@, but openness does not since the vector bundle
projection is not proper.  In any case, we give an independent proof of
continuity, as it is a model for the proof of subsequent statements for which
we will not give detailed proofs.
\begin{theorem}\label{the:pi*fcont}
Let\/ $\map{\pi_{\man{E}}}{\man{E}}{\man{M}}$ be a\/ $\C^\omega$-vector
bundle.  Then the mapping
\begin{equation*}
\func[\omega]{\man{M}}\ni f\mapsto\pi^*f\in\func[\omega]{\man{E}}
\end{equation*}
is an homeomorphism onto its image.
\begin{proof}
It is clear that the asserted map is injective, so we focus on its
topological attributes.

We let $\metric_{\man{M}}$ be a $\C^\omega$-Riemannian metric on $\man{M}$\@,
$\metric_\pi$ be a $\C^\omega$-vector bundle connection in $\man{E}$\@,
$\nabla^{\man{M}}$ be the Levi-Civita connection for $\metric_{\man{M}}$\@,
and $\nabla^\pi$ be a $\C^\omega$-vector bundle connection in $\man{E}$\@.
Corresponding to this, we have a Riemannian metric $\metric_{\man{E}}$ on
$\man{E}$ with its Levi-Civita connection $\nabla^{\man{E}}$\@, as in
Section~\ref{subsec:Esubmersion}\@.  We denote the associated seminorms for
$\func[\omega]{\man{M}}$ and $\func[\omega]{\man{E}}$ by
$p^\omega_{\nbhd{K},\vect{a}}$ and $q^\omega_{\nbhd{L},\vect{a}}$ for
$\nbhd{K}\subset\man{M}$ and $\nbhd{L}\subset\man{E}$ compact, and for
$\vect{a}\in\c_0(\integernn;\realp)$\@.

Let us make some preliminary computations.

By Lemma~\ref{lem:PsymiterateddersI}\@, we have
\begin{equation}\label{eq:PSESM}
\Sym_m\scirc\nabla^{\man{E},m}\pi_{\man{E}}^*f(e)=
\sum_{s=0}^m\what{A}^m_s(\Sym_{s+1}\scirc\pi_{\man{E}}^*\nabla^{\man{M},s}f(e)).
\end{equation}
By Lemma~\ref{lem:A(S)<AS}\@, we have
\begin{equation*}
\dnorm{A^m_s(\beta_s)}_{\metric_{\man{E}}}\le\dnorm{A^m_s}_{\metric_{\man{E}}}
\dnorm{\beta_s}_{\metric_{\man{E}}}
\end{equation*}
for $\beta_s\in\tensor*[s]{\ctb[e]{\man{E}}}$\@, $m\in\integerp$\@, and
$s\in\{0,1,\dots,m\}$\@.  By Lemma~\ref{lem:Symnorm}\@,
\begin{equation*}
\dnorm{\Sym_s(A)}_{\metric_{\man{E}}}\le\dnorm{A}_{\metric_{\man{E}}}
\end{equation*}
for $A\in\tensor*[s]{\ctb{\man{N}}}$ and $s\in\integerp$\@.  Thus,
recalling~\eqref{eq:AhatA}\@,
\begin{equation*}
\dnorm{\what{A}^m_s(\Sym_s(\beta_s))}_{\metric_{\man{E}}}=
\dnorm{\Sym_m\scirc A^m_s(\beta_s)}_{\metric_{\man{E}}}
\le\dnorm{A^m_s}_{\metric_{\man{E}}}\dnorm{\beta_s}_{\metric_{\man{E}}},
\end{equation*}
for $\beta_s\in\tensor*[s]{\pi_{\man{E}}^*\ctb{\man{M}}}$\@,
$m\in\integerp$\@, $s\in\{1,\dots,m\}$\@.

Let $\nbhd{L}\subset\man{E}$ be compact.  By Lemmata~\ref{lem:pissynabla}
and~\ref{lem:particular-bounds} with $r=0$\@, there exist
$C_1,\sigma_1,\rho_1\in\realp$ such that
\begin{equation*}
\dnorm{A^k_s(e)}_{\metric_{\man{E}}}\le
C_1\sigma_1^{-k}\rho_1^{-(k-s)}(k-s)!,\qquad
k\in\integerp,\ s\in\{0,1,\dots,k-1\},\ e\in\nbhd{L}.
\end{equation*}
Without loss of generality, we assume that $C_1\ge1$ and
$\sigma_1,\rho_1\le1$\@.  Thus, using Lemma~\ref{lem:Pnormbase} and
abbreviating $\sigma_2=\sigma_1\rho_1$\@, we have
\begin{equation*}
\dnorm{\what{A}^k_s(\pi_{\man{E}}^*\Sym_s\scirc\nabla^{\man{M},s}
f(e))}_{\metric_{\man{E}}}\le C_1\sigma_2^{-k}(k-s)!
\dnorm{\Sym_s\scirc\nabla^{\man{M},s}f(\pi_{\man{E}}(e))}_{\metric_{\man{M}}}
\end{equation*}
for $k\in\integernn$\@, $s\in\{0,1,\dots,k\}$\@, $e\in\nbhd{L}$\@.  Thus,
by~\eqref{eq:12inftynorms}\@,~\eqref{eq:PSESM}\@, and Lemma~\ref{lem:pbnormbase}\@,
\begin{align*}
\dlnorm j_m(\pi_{\man{E}}^*f)(e)\drnorm_{\metric_{\man{E},m}}\le&\;
\sum_{k=0}^m\frac{1}{k!}\dnorm{\Sym_k\scirc\nabla^{\man{E},k}
\pi_{\man{E}}^*f(e)}_{\metric_{\man{E}}}\\
=&\;\sum_{k=0}^m\frac{1}{k!}
\adnorm{\sum_{s=0}^k\what{A}^k_s(\pi_{\man{E}}^*\Sym_s\scirc
\nabla^{\man{M},s}f(e))}_{\metric_{\man{E}}}\\
\le&\;\sum_{k=0}^m\sum_{s=0}^kC_1\sigma_2^{-k}
\frac{s!(k-s)!}{k!}\frac{1}{s!}
\dnorm{\Sym_s\scirc\nabla^{\man{M},s}f(\pi_{\man{E}}(e))}_{\metric_{\man{M}}}
\end{align*}
for $e\in\nbhd{L}$ and $m\in\integernn$\@.  Now note that
\begin{equation*}
\frac{s!(k-s)!}{k!}\le1,\quad
C_1\sigma_2^{-k}\le C_1\sigma_2^{-m},
\end{equation*}
for $s\in\{0,1,\dots,m\}$\@, $k\in\{0,1,\dots,s\}$\@, since $\sigma_2\le1$\@.
Then
\begin{align*}
\dnorm{j_m(\pi_{\man{E}}^*f)(e)}_{\metric_{\man{E},m}}\le&\;
C_1\sigma_2^{-m}\sum_{k=0}^m\sum_{s=0}^k\frac{1}{s!}
\dnorm{\Sym_s\scirc\nabla^{\man{M},s}f(\pi_{\man{E}}(e))}_{\metric_{\man{M}}}\\
\le&\;C_1\sigma_2^{-m}\sum_{k=0}^m\sum_{s=0}^m
\frac{1}{s!}\dnorm{\Sym_s\scirc\nabla^{\man{M},s}
f(\pi_{\man{E}}(e))}_{\metric_{\man{M}}}\\
=&\;(m+1)C_1\sigma_2^{-m}\sum_{s=0}^m
\frac{1}{s!}\dnorm{\Sym_s\scirc\nabla^{\man{M},s}
f(\pi_{\man{E}}(e))}_{\metric_{\man{M}}}.
\end{align*}
Now let $\sigma<\sigma_2$ and note that
\begin{equation*}
\lim_{m\to\infty}(m+1)\frac{\sigma_2^{-m}}{\sigma^{-m}}=0.
\end{equation*}
Thus there exists $N\in\integerp$ such that
\begin{equation*}
(m+1)C_1\sigma_2^{-m}\le C_1\sigma^{-m},\qquad m\ge N.
\end{equation*}
Let
\begin{equation*}
C=\max\left\{C_1,C_1\frac{\sigma}{\sigma_2},
2C_1\left(\frac{\sigma}{\sigma_2}\right)^2,
\dots,NC_1\left(\frac{\sigma}{\sigma_2}\right)^N\right\}.
\end{equation*}
We then immediately have $(m+1)C_1\sigma_2^{-m}\le C\sigma^{-m}$ for all
$m\in\integernn$\@.  We then have, by~\eqref{eq:12inftynorms}\@,
\begin{align*}
\dnorm{j_m(\pi_{\man{E}}^*f)(e)}_{\metric_{\man{E},m}}\le&\;
C\sigma^{-m}\sum_{s=0}^m\frac{1}{s!}
\dnorm{\Sym_s\scirc\nabla^{\man{M},s}f(\pi_{\man{E}}(e))}_{\metric_{\man{M}}}\\
\le&\;C\sqrt{m+1}\sigma^{-m}\dnorm{j_mf(\pi_{\man{E}}(e))}_{\metric_{\man{M},m}}.
\end{align*}
By modifying $C$ and $\sigma$ guided by what we did just preceding, we get
\begin{equation*}
\dnorm{j_m(\pi_{\man{E}}^*f)(e)}_{\metric_{\man{E},m}}\le
C\sigma^{-m}\dnorm{j_mf(\pi_{\man{E}}(e))}_{\metric_{\man{M},m}}.
\end{equation*}
Now let $\vect{a}\in\c_0(\integernn;\realp)$ and define
$\vect{a}'\in\c_0(\integernn;\realp)$ be defined by $a'_0=Ca_0$ and
$a'_j=a_j\sigma^{-1}$\@, $j\in\integerp$\@.  Then we have
\begin{equation*}
q^\omega_{\nbhd{L},\vect{a}}(\pi_{\man{E}}^*f)\le
p^\omega_{\pi_{\man{E}}(\nbhd{L}),\vect{a}'}(f),
\end{equation*}
giving continuity in this case.

Now we show that $\pi_{\man{E}}^*$ is open onto its image.  The idea here is
to make some preliminary observations to put ourselves in a position to be
able to say, ``Now proceed as above.''

By Lemma~\ref{lem:PsymiterateddersII}\@, we have
\begin{equation}\label{eq:PSMSE}
\Sym_m\scirc\pi_{\man{E}}^*\nabla^{\man{M},m}f(e)=
\sum_{s=0}^m\what{B}^m_s(\Sym_s\scirc\nabla^{\man{E},s}\pi_{\man{E}}^*f(e)).
\end{equation}
For a compact $\nbhd{L}\subset\man{E}$ we can proceed as above to give a
bound
\begin{equation*}
\dnorm{j_mf(\pi_{\man{E}}(e))}_{\metric_{\man{M},m}}\le
C\sigma^{-m}\dnorm{j_m(\pi_{\man{E}}^*f)(e)}_{\metric_{\man{E},m}},
\qquad e\in\nbhd{L}.
\end{equation*}
We need to choose the compact set $\nbhd{L}$ in a specific way.  We let
$\nbhd{K}\subset\man{M}$ be compact and choose a continuous section
$\xi\in\sections[0]{\man{E}}$\@, and then take $\nbhd{L}=\xi(\nbhd{K})$\@.
Then we have the estimate
\begin{equation*}
\dnorm{j_mf(x)}_{\metric_{\man{M},m}}\le
C\sigma^{-m}\dnorm{j_m(\pi_{\man{E}}^*f)(\xi(x))}_{\metric_{\man{E},m}},
\qquad x\in\nbhd{K}.
\end{equation*}
Now we can mirror the arguments above for continuity to give the bound
\begin{equation*}
p^\omega_{\nbhd{K},\vect{a}}(f)\le
q^\omega_{\xi(\nbhd{K}),\vect{a}'}(\pi_{\man{E}}^*f),
\end{equation*}
and from this we conclude that $f\mapsto\pi_{\man{E}}^*f$ is indeed open onto
its image by Lemma~\ref{plem:open-mapping} from the proof of Theorem~\ref{the:+xcont}\@.
\end{proof}
\end{theorem}

Now we consider vertical lifts of sections.
\begin{theorem}
Let\/ $\map{\pi_{\man{E}}}{\man{E}}{\man{M}}$ be a\/ $\C^\omega$-vector
bundle.  Then the mapping
\begin{equation*}
\sections[\omega]{\man{E}}\ni\xi\mapsto\vl{\xi}\in
\sections[\omega]{\tb{\man{E}}}
\end{equation*}
is an homeomorphism onto its image.
\begin{proof}
This follows in the same manner as Theorem~\ref{the:pi*fcont}\@, using
Lemmata~\ref{lem:VsymiterateddersI}\@,~\ref{lem:VsymiterateddersII}\@, and~\ref{lem:Vnormbase}\@.
\end{proof}
\end{theorem}

One has the similar result for vertical lifts of endomorphisms.
\begin{theorem}\label{the:Lvcont}
Let\/ $\map{\pi_{\man{E}}}{\man{E}}{\man{M}}$ be a\/ $\C^\omega$-vector
bundle.  Then the mapping
\begin{equation*}
\sections[\omega]{\End(\man{E})}\ni L\mapsto\vl{L}\in
\sections[\omega]{\End(\tb{\man{E}})}
\end{equation*}
is an homeomorphism onto its image.
\begin{proof}
This follows in the same manner as Theorem~\ref{the:pi*fcont}\@, using
Lemmata~\ref{lem:LsymiterateddersI}\@,~\ref{lem:LsymiterateddersII}\@, and~\ref{lem:Lnormbase}\@.
\end{proof}
\end{theorem}

Now we consider horizontal lifts of vector fields.
\begin{theorem}\label{the:hliftX}
Let\/ $\map{\pi_{\man{E}}}{\man{E}}{\man{M}}$ be a\/ $\C^\omega$-vector
bundle.  Then the mapping
\begin{equation*}
\sections[\omega]{\tb{\man{M}}}\ni
X\mapsto\hl{X}\in\sections[\omega]{\tb{\man{E}}}
\end{equation*}
is an homeomorphism onto its image.
\begin{proof}
This follows in the same manner as Theorem~\ref{the:pi*fcont}\@, using
Lemmata~\ref{lem:HsymiterateddersI}\@,~\ref{lem:HsymiterateddersII}\@, and~\ref{lem:Hnormbase}\@.
\end{proof}
\end{theorem}

Now we consider vertical lifts of sections of the dual bundle.
\begin{theorem}\label{the:lambdavcont}
Let\/ $\map{\pi_{\man{E}}}{\man{E}}{\man{M}}$ be a\/ $\C^\omega$-vector
bundle.  Then the mapping
\begin{equation*}
\sections[\omega]{\dual{\man{E}}}\ni\lambda\mapsto\vl{\lambda}\in
\sections[\omega]{\ctb{\man{E}}}
\end{equation*}
is an homeomorphism onto its image.
\begin{proof}
This follows in the same manner as Theorem~\ref{the:pi*fcont}\@, using
Lemmata~\ref{lem:V*symiterateddersI}\@,~\ref{lem:V*symiterateddersII}\@, and~\ref{lem:V*normbase}\@.
\end{proof}
\end{theorem}

Next we consider vertical evaluations of sections of the dual bundle.
\begin{theorem}\label{the:lambdaecont}
Let\/ $\map{\pi_{\man{E}}}{\man{E}}{\man{M}}$ be a\/ $\C^\omega$-vector
bundle.  Then the mapping
\begin{equation*}
\sections[\omega]{\dual{\man{E}}}\ni\lambda\mapsto\ve{\lambda}\in
\func[\omega]{\man{E}}
\end{equation*}
is an homeomorphism onto its image.
\begin{proof}
Since the given map is clearly injective, we focus on its topological
properties.

We let $\metric_{\man{M}}$ be a $\C^\omega$-Riemannian metric on $\man{M}$\@,
$\metric_\pi$ be a $\C^\omega$-vector bundle connection in $\man{E}$\@,
$\nabla^{\man{M}}$ be the Levi-Civita connection for $\metric_{\man{M}}$\@,
and $\nabla^\pi$ be a $\C^\omega$-vector bundle connection in $\man{E}$\@.
Corresponding to this, we have a Riemannian metric $\metric_{\man{E}}$ on
$\man{E}$ with its Levi-Civita connection $\nabla^{\man{E}}$\@, as in
Section~\ref{subsec:Esubmersion}\@.  We denote the associated seminorms for
$\sections[\omega]{\dual{\man{E}}}$ and $\func[\omega]{\man{E}}$ by
$p^\omega_{\nbhd{K},\vect{a}}$ and $q^\omega_{\nbhd{L},\vect{a}}$ for
$\nbhd{K}\subset\man{M}$ and $\nbhd{L}\subset\man{E}$ compact, and for
$\vect{a}\in\c_0(\integernn;\realp)$\@.

Let us make some preliminary computations.

By Lemma~\ref{lem:DsymiterateddersI}\@, we have
\begin{equation}\label{eq:DSESM}
\begin{aligned}
\ve{\lambda}(e)=&\;\ve{\lambda}(e),\\
\nabla^{\man{E}}\ve{\lambda}(e)=&\;
\what{A}^1_1(\ve{(\nabla^{\pi_{\man{E}}}\lambda)}(e))+
\what{A}^1_0(\ve{\lambda}(e))+\what{C}^1_0(\vl{\lambda}(e)),\\
\vdots\,&\\
(\Sym_m\otimes\id_{\ctb{\man{E}}})\scirc\nabla^{\man{E},m}\ve{\lambda}(e)=&\;
\sum_{s=0}^m\what{A}^m_s((\Sym_s\otimes\id_{\ctb{\man{E}}})\scirc
\ve{(\nabla^{\man{M},\pi_{\man{E}},s}\lambda)}(e))\\
&\;+\sum_{s=0}^{m-1}\what{C}^m_s((\Sym_s\otimes\id_{\ctb{\man{E}}})
\scirc\vl{(\nabla^{\man{M},\pi_{\man{E}},s}\lambda)}(e)).
\end{aligned}
\end{equation}
Just as in the proof of Theorem~\ref{the:pi*fcont}\@, by
Lemmata~\ref{lem:A(S)<AS} and~\ref{lem:Symnorm}\@, and the appropriate
analogue of equation~\eqref{eq:Vsymmetrise1} that would appear in a fully
fleshed out proof of Lemma~\ref{lem:DiterateddersI}\@, we have bounds
\begin{equation*}
\dnorm{\what{A}^m_s(\Sym_s(\beta_s))}_{\metric_{\man{E}}}=
\dnorm{\Sym_m\scirc A^m_s(\beta_s)}_{\metric_{\man{E}}}
\le\dnorm{A^m_s}_{\metric_{\man{E}}}\dnorm{\beta_s}_{\metric_{\man{E}}}
\end{equation*}
and
\begin{equation*}
\dnorm{\what{C}^m_s(\Sym_s(\gamma_s))}_{\metric_{\man{E}}}=
\dnorm{\Sym_m\scirc C^m_s(\gamma_s)}_{\metric_{\man{E}}}
\le\dnorm{C^m_s}_{\metric_{\man{E}}}\dnorm{\gamma_s}_{\metric_{\man{E}}}.
\end{equation*}

Let $\nbhd{L}\subset\man{E}$ be compact.  By Lemmata~\ref{lem:pissynabla}
and~\ref{lem:particular-bounds} with $r=0$\@, there exist
$C_1,\sigma_1,\rho_1\in\realp$ such that
\begin{equation*}
\dnorm{A^k_s(e)}_{\metric_{\man{E}}}\le
C_1\sigma_1^{-k}\rho_1^{-(k-s)}(k-s)!,\qquad
k\in\integernn,\ s\in\{0,1,\dots,k\},\ e\in\nbhd{L},
\end{equation*}
and
\begin{equation*}
\dnorm{C^k_s(e)}_{\metric_{\man{E}}}\le
C_1\sigma_1^{-k}\rho_1^{-(k-s)}(k-s)!,\qquad
k\in\integernn,\ s\in\{0,1,\dots,k-1\},\ e\in\nbhd{L},
\end{equation*}
Without loss of generality, we assume that $C_1\ge1$ and
$\sigma_1,\rho_1\le1$\@.  Thus, using Lemma~\ref{lem:Dnormbase} and
abbreviating $\sigma_2=\sigma_1\rho_1$\@, we have
\begin{multline*}
\dnorm{\what{A}^k_s((\Sym_s\otimes\id_{\ctb{\man{E}}})\scirc
\ve{(\nabla^{\man{M},\pi_{\man{E}},s}
\lambda)}(e))}_{\metric_{\man{M},\pi_{\man{E}}}}\\
\le C_1\sigma_2^{-k}(k-s)!\dnorm{(\Sym_s\otimes\id_{\ctb{\man{E}}})\scirc
\ve{(\nabla^{\man{M},\pi_{\man{E}},s}
\lambda)}(e)}_{\metric_{\man{M},\pi_{\man{E}}}}
\end{multline*}
for $k\in\integernn$\@, $s\in\{0,1,\dots,k\}$\@, $e\in\nbhd{L}$\@, and
\begin{multline*}
\dnorm{\what{C}^k_s((\Sym_s\otimes\id_{\ctb{\man{E}}})\scirc
\ve{(\nabla^{\man{M},\pi_{\man{E}},s}
\lambda)}(e))}_{\metric_{\man{M},\pi_{\man{E}}}}\\
\le C_1\sigma_2^{-k}(k-s)!\dnorm{(\Sym_s\otimes\id_{\ctb{\man{E}}})\scirc
\ve{(\nabla^{\man{M},\pi_{\man{E}},s}
\lambda)}(e)}_{\metric_{\man{M},\pi_{\man{E}}}}
\end{multline*}
for $k\in\integernn$\@, $s\in\{0,1,\dots,k-1\}$\@, $e\in\nbhd{L}$\@.  Thus,
by~\eqref{eq:12inftynorms} and~\eqref{eq:DSESM}\@,
\begin{align*}
\dnorm{j_m\ve{\lambda}(e)}_{\metric_{\man{E},m}}\le&\;
\sum_{k=0}^m\frac{1}{k!}
\dnorm{\Sym_k\scirc\nabla^{\man{E},k}\ve{\lambda}(e)}_{\metric_{\man{E}}}\\
\le&\;\sum_{k=0}^m\frac{1}{k!}\adnorm{\sum_{s=0}^k
\what{A}^k_s((\Sym_s\scirc\id_{\ctb{\man{E}}})\scirc
\ve{(\nabla^{\man{M},\pi_{\man{E}},s}
\lambda)}(e))}_{\metric_{\man{M},\pi_{\man{E}}}}\\
&\;+\sum_{k=0}^{m-1}\frac{1}{k!}\adnorm{\sum_{s=0}^k
\what{C}^k_s((\Sym_s\scirc\id_{\ctb{\man{E}}})\scirc
\vl{(\nabla^{\man{M},\pi_{\man{E}},s}
\lambda)}(e))}_{\metric_{\man{M},\pi_{\man{E}}}}\\
\le&\;\sum_{k=0}^mC_1\sigma_2^{-k}\frac{s!(k-s)!}{k!}\frac{1}{s!}
\adnorm{\sum_{s=0}^k(\Sym_s\scirc\id_{\ctb{\man{E}}})\scirc
\ve{(\nabla^{\man{M},\pi_{\man{E}},s}
\lambda)}(e)}_{\metric_{\man{M},\pi_{\man{E}}}}\\
&\;+\sum_{k=0}^{m-1}C_1\sigma_2^{-k}\frac{s!(k-s)!}{k!}\frac{1}{s!}
\adnorm{\sum_{s=0}^k(\Sym_s\scirc\id_{\ctb{\man{E}}})\scirc
\vl{(\nabla^{\man{M},\pi_{\man{E}},s}
\lambda)}(e)}_{\metric_{\man{M},\pi_{\man{E}}}}
\end{align*}
for $e\in\nbhd{L}$ and $m\in\integernn$\@.  Now note that
\begin{equation*}
\frac{s!(k-s)!}{k!}\le1,\quad C_1\sigma_2^{-k}\le C_1\sigma_2^{-m},
\end{equation*}
for $s\in\{0,1,\dots,m-1\}$\@, $k\in\{0,1,\dots,s\}$\@, since
$\sigma_2\le1$\@.  Then
\begin{align*}
\dnorm{j_m\ve{\lambda}(e)}_{\metric_{\man{E},m}}\le&\;
C_1\sigma_2^{-m}\sum_{k=0}^m\sum_{s=0}^m\frac{1}{s!}
\adnorm{(\Sym_s\scirc\id_{\ctb{\man{E}}})\scirc
\ve{(\nabla^{\man{M},\pi_{\man{E}},s}
\lambda)}(e)}_{\metric_{\man{M},\pi_{\man{E}}}}\\
&\;+C_1\sigma_2^{-m}\sum_{k=0}^{m-1}\sum_{s=0}^{m-1}\frac{1}{s!}
\adnorm{(\Sym_s\scirc\id_{\ctb{\man{E}}})\scirc
\vl{(\nabla^{\man{M},\pi_{\man{E}},s}
\lambda)}(e)}_{\metric_{\man{M},\pi_{\man{E}}}}\\
=&\;(m+1)C_1\sigma_2^{-m}\sum_{s=0}^m\frac{1}{s!}
\adnorm{(\Sym_s\scirc\id_{\ctb{\man{E}}})\scirc
\ve{(\nabla^{\man{M},\pi_{\man{E}},s}
\lambda)}(e)}_{\metric_{\man{M},\pi_{\man{E}}}}\\
&\;+(m+1)C_1\sigma_2^{-m}\sum_{s=0}^{m-1}\frac{1}{s!}
\adnorm{(\Sym_s\scirc\id_{\ctb{\man{E}}})\scirc
\vl{(\nabla^{\man{M},\pi_{\man{E}},s}
\lambda)}(e)}_{\metric_{\man{M},\pi_{\man{E}}}}.
\end{align*}
Now let $\sigma<\sigma_2$ and note that
\begin{equation*}
\lim_{m\to\infty}(m+1)\frac{\sigma_2^{-m}}{\sigma^{-m}}=0.
\end{equation*}
Thus there exists $N\in\integerp$ such that
\begin{equation*}
(m+1)C_1\sigma_2^{-m}\le C_1\sigma^{-m},\qquad m\ge N.
\end{equation*}
Let
\begin{equation*}
C=\max\left\{C_1,2C_1\frac{\sigma}{\sigma_2},
3C_1\left(\frac{\sigma}{\sigma_2}\right)^2,
\dots,(N+1)C_1\left(\frac{\sigma}{\sigma_2}\right)^N\right\}.
\end{equation*}
We then immediately have $(m+1)C_1\sigma_2^{-m}\le C\sigma^{-m}$ for all
$m\in\integernn$\@.  We then have, using~\eqref{eq:12inftynorms}\@,
\begin{align*}
\dnorm{j_m\ve{\lambda}(e)}_{\metric_{\man{E},m}}\le&\;
C\sigma^{-m}\left(\sum_{s=0}^m\frac{1}{s!}
\adnorm{(\Sym_s\scirc\id_{\ctb{\man{E}}})\scirc
\ve{(\nabla^{\man{M},\pi_{\man{E}},s}
\lambda)}(e)}_{\metric_{\man{M},\pi_{\man{E}}}}\right.\\
&\;+\left.\sum_{s=0}^{m-1}\frac{1}{s!}
\adnorm{(\Sym_s\scirc\id_{\ctb{\man{E}}})\scirc
\vl{(\nabla^{\man{M},\pi_{\man{E}},s}
\lambda)}(e)}_{\metric_{\man{M},\pi_{\man{E}}}}\right)\\
\le&\;C\sqrt{m+1}\sigma^{-m}(
\dnorm{j_m\lambda(\pi_{\man{E}}(e))(e)}_{\metric_{\man{M},\pi_{\man{E}},m}}+
\dnorm{j_{m-1}\lambda(\pi_{\man{E}}(e))}_{\metric_{\man{M},\pi_{\man{E}},m-1}}).
\end{align*}
By modifying $C$ and $\sigma$ just as we did in the preceding, we get
\begin{equation*}
\dnorm{j_m\ve{\lambda}(e)}_{\metric_{\man{E},m}}\le
C\sigma^{-m}(\dnorm{j_m
\lambda(\pi_{\man{E}}(e))(e)}_{\metric_{\man{M},\pi_{\man{E}},m}}+
\dnorm{j_{m-1}\lambda(\pi_{\man{E}}(e))}_{\metric_{\man{M},\pi_{\man{E}},m-1}}).
\end{equation*}
We take
\begin{equation*}
\alpha=\max\{1,\sup\setdef{\dnorm{e}_{\metric_{\pi_{\man{E}}}}}{e\in\nbhd{L}}
\end{equation*}
and then use Lemma~\ref{lem:A(S)<AS} to arrive at
\begin{equation*}
\dnorm{j_m\ve{\lambda}(e)}_{\metric_{\man{E},m}}\le2\alpha C\sigma^{-m}
\dnorm{j_m\lambda(\pi_{\man{E}}(e))(e)}_{\metric_{\man{M},\pi_{\man{E}},m}}
\end{equation*}
Now, given $\vect{a}\in\c_0(\integernn;\realp)$\@, we define
$\vect{a}'\in\c_0(\integernn;\realp)$ by $a'_0=2\alpha Ca_0$ and
$a'_j=a_j\sigma^{-1}$\@, $j\in\integerp$\@, we then have
\begin{equation*}
q^\omega_{\nbhd{L},\vect{a}}(\ve{\lambda})\le p^\omega_{\pi_{\man{E}}(\nbhd{L}),\vect{a}'}(\lambda),
\end{equation*}
and this gives this part of the result.

Now we turn to showing that the mapping of the lemma is open onto its image.
By Lemma~\ref{lem:DsymiterateddersII}\@, we have
\begin{equation}\label{eq:DSMSE}
\begin{aligned}
\ve{\lambda}(e)=&\;\ve{\lambda}(e),\\
\ve{(\nabla^{\pi_{\man{E}}}\lambda)}(e)=&\;
\what{B}^1_1(\nabla^{\man{E}}\ve{\lambda}(e))+
\what{B}^1_0(\ve{\lambda}(e))+\what{D}^1_0(\vl{\lambda}(e)),\\
(\Sym_2\otimes\id_{\tb{\man{E}}})\scirc
\ve{(\nabla^{\man{M},\pi_{\man{E}},2}\lambda)}(e)=&\;
\what{B}^2_2(\nabla^{\man{E},2}\ve{\lambda}(e))+
\what{B}^2_1(\nabla^{\man{E}}\ve{\lambda}(e))+\what{B}^2_0(\ve{\lambda}(e))\\
&\;+\what{D}^2_1(\vl{(\nabla^{\man{M},\pi_{\man{E}}}\lambda)}(e))+
\what{D}^1_0(\vl{\lambda}(e)),\\
\vdots\,&\\
(\Sym_m\otimes\id_{\ctb{\man{E}}})\scirc
\ve{(\nabla^{\man{M},\pi_{\man{E}},m}\lambda)}(e)=&\;
\sum_{s=0}^m\what{B}^m_s((\Sym_s\otimes\id_{\ctb{\man{E}}})\scirc
\nabla^{\man{E},s}\ve{\lambda}(e))\\
&\;+\sum_{s=0}^{m-1}\what{D}^m_s((\Sym_s\otimes\id_{\ctb{\man{E}}})
\scirc\nabla^{\man{E},s}\vl{\lambda}(e)).
\end{aligned}
\end{equation}
Just as in the proof of Theorem~\ref{the:pi*fcont}\@, by
Lemmata~\ref{lem:A(S)<AS} and~\ref{lem:Symnorm}\@, and the appropriate
analogue of equation~\eqref{eq:Vsymmetrise1} that would appear in a fully
fleshed out proof of Lemma~\ref{lem:DiterateddersII}\@, we have bounds
\begin{gather*}
\dnorm{\what{B}^m_s(\Sym_s(\beta_s))}_{\metric_{\man{E}}}=
\dnorm{\Sym_m\scirc B^m_s(\beta_s)}_{\metric_{\man{E}}}
\le\dnorm{B^m_s}_{\metric_{\man{E}}}\dnorm{\beta_s}_{\metric_{\man{E}}},\\
\dnorm{\what{D}^m_s(\Sym_s(\gamma_s))}_{\metric_{\man{E}}}=
\dnorm{\Sym_m\scirc D^m_s(\gamma_s)}_{\metric_{\man{E}}}
\le\dnorm{D^m_s}_{\metric_{\man{E}}}\dnorm{\gamma_s}_{\metric_{\man{E}}}.
\end{gather*}
Proceeding analogously to the continuity proof above and using
Lemma~\ref{lem:Dnormbase}\@, we deduce that there exist
$C_1,\sigma_1\in\realp$ such that
\begin{equation}\label{eq:lambdaeopenest}
\dnorm{j_m\lambda(\pi_{\man{E}}(e))(e)}_{\metric_{\man{M},\pi_{\man{E}},m}}\le
C_1\sigma_1^{-m}(\dnorm{j_m\ve{\lambda}(e)}_{\metric_{\man{E},m}}+
\dnorm{j_{m-1}\vl{\lambda}(e)}_{\metric_{\man{E},m-1}}),\qquad
e\in\nbhd{L}.
\end{equation}
Now let $\nbhd{K}\subset\man{M}$ be compact and let
$\vect{a}\in\c_0(\integernn,\realp)$\@.  Define
\begin{equation*}
\nbhd{L}=\pi_{\man{E}}^{-1}(\nbhd{K})\cap
\setdef{e\in\man{E}}{\dnorm{e}_{\metric_{\pi_{\man{E}}}}=1},
\end{equation*}
noting that $\nbhd{L}$ is compact.  Let $n=\dim(\man{M})$ and let $k$ be the
fibre dimension of $\man{E}$\@.  By Lemma~\ref{lem:opnormupper}\@, and
equations~\eqref{eq:symalgdim} and~\eqref{eq:multinomialest}\@, we have
\begin{align*}
\dnorm{j_m\lambda(x)}_{\metric_{\man{M},\pi_{\man{E}},m}}\le&\;
\sum_{j=0}^m\sqrt{k\binom{n+j-1}{j}}\sup
\setdef{\dnorm{j_m\lambda(\pi_{\man{E}}(e))(e)}_{\metric_{\man{M},\pi_{\man{E}},m}}}
{e\in\nbhd{L}}\\
\le&\;\sum_{j=0}^mk\binom{n+j-1}{j}\sup
\setdef{\dnorm{j_m\lambda(\pi_{\man{E}}(e))(e)}_{\metric_{\man{M},\pi_{\man{E}},m}}}
{e\in\nbhd{L}}\\
\le&\;m^22^{n+m}\sup
\setdef{\dnorm{j_m\lambda(\pi_{\man{E}}(e))(e)}_{\metric_{\man{M},\pi_{\man{E}},m}}}
{e\in\nbhd{L}}
\end{align*}
for $x\in\nbhd{K}$\@.  For $\sigma_2<\frac{1}{2}$\@,
\begin{equation*}
\lim_{m\to\infty}m^2\frac{2^m}{\sigma_2^{-m}}=0.
\end{equation*}
By by now familiar arguments, one of which the reader can find in the first
part of the proof, we can combine this with~\eqref{eq:lambdaeopenest} to
arrive at $C,\sigma\in\realp$ for which
\begin{equation*}
\dnorm{j_m\lambda(x)}_{\metric_{\man{M},\pi_{\man{E}},m}}\le
C\sigma^{-m}(\sup
\setdef{\dnorm{j_m\ve{\lambda}(e)}_{\metric_{\man{M},\pi_{\man{E}},m}}}
{e\in\nbhd{L}}\\
+\sup\setdef{\dnorm{j_m\vl{\lambda}(e)}_{\metric_{\man{E},m}}}
{e\in\nbhd{L}})
\end{equation*}
for $x\in\nbhd{K}$\@.  Taking $\vect{a}'\in\c_0(\integernn;\realp)$ to be
defined by $a'_0=Ca_0$\@, $a'_j=\sigma^{-1}a_j$\@, $j\in\integerp$\@, we have
\begin{equation*}
q^\omega_{\nbhd{K},\vect{a}}(\lambda)\le
p^\omega_{\nbhd{L},\vect{a}'}(\ve{\lambda})+
p^\omega_{\nbhd{L},\vect{a}'}(\vl{\lambda}).
\end{equation*}
By Lemma~\ref{plem:open-mapping} from the proof of
Theorem~\ref{the:+xcont}\@, this shows that the mapping
\begin{equation*}
\sections[\omega]{\dual{\man{E}}}\ni\lambda\mapsto
(\ve{\lambda},\vl{\lambda})\in\func[\omega]{\man{E}}\oplus
\sections[\omega]{\tb{\man{E}}}
\end{equation*}
is open onto its image.  This part of the lemma now follows from the
following simple fact.
\begin{prooflemma}
Let\/ $\ts{S}$\@,\/ $\ts{T}_1$\@, and $\ts{T}_2$ be topological spaces and
let\/ $\map{\Phi}{\ts{S}}{\ts{T}_1\times\ts{T}_2}$ be an open mapping onto
its image.  Then the mappings\/ $\pr_1\scirc\Phi$ and\/ $\pr_2\scirc\Phi$ are
open onto their images.
\begin{subproof}
Let $\nbhd{O}\subset\ts{S}$ be open so that $\Phi(\nbhd{O})$ is open in
$\image(\Phi)$\@.  Then, for each $(y_1,y_2)\in\nbhd{O}$\@, there exists a
neighbourhood $\nbhd{N}_1\subset\image(\pr_1\scirc\Phi)$ of $y_1$ and a
neighbourhood $\nbhd{N}_2\subset\image(\pr_2\scirc\Phi)$ of $x_2$ such that
$\nbhd{N}_1\times\nbhd{N}_2\subset\Phi(\nbhd{O})$\@.  This immediately gives
the lemma.
\end{subproof}
\end{prooflemma}

Thus we arrive at the conclusion that the mapping
\begin{equation*}
\sections[\omega]{\dual{\man{E}}}\ni\lambda\mapsto
\ve{\lambda}\in\func[\omega]{\man{E}}
\end{equation*}
is open onto its image, as desired.
\end{proof}
\end{theorem}

Finally, we consider vertical evaluations of sections of the endomorphism
bundle.
\begin{theorem}
Let\/ $\map{\pi_{\man{E}}}{\man{E}}{\man{M}}$ be a\/ $\C^\omega$-vector bundle.
Then the mapping
\begin{equation*}
\sections[\omega]{\End(\man{E})}\ni L\mapsto\ve{L}\in
\sections[\omega]{\tb{\man{E}}}
\end{equation*}
is an homeomorphism onto its image.
\begin{proof}
This follows in the same manner as Theorem~\ref{the:lambdaecont}\@, using
Lemmata~\ref{lem:CsymiterateddersI}\@,~\ref{lem:CsymiterateddersII}\@, and~\ref{lem:Cnormbase}\@.
\end{proof}
\end{theorem}

As an illustration of how continuity of these lifts can be helpful, let us
consider the continuity of the map that assigns to a vector field on a
manifold the tangent lift of that vector field.  Precisely, let $\man{M}$ be
a real analytic manifold and let $X\in\sections[\omega]{\tb{\man{M}}}$ be a
real analytic vector field.  The \defn{tangent lift} of $X$ is the vector
field $\tlift{X}\in\sections[\omega]{\tb{\tb{\man{M}}}}$ on $\tb{\man{M}}$
whose flow is the derivative of the flow for $X$\@:
\begin{equation}\label{eq:XTflow}
\flow{\tlift{X}}{t}(v_x)=\tf[x]{\flow{X}{t}}(v_x)\enspace\implies\enspace
\tlift{X}=\derivatzero{}{t}\tf[x]{\flow{X}{t}}(v_x).
\end{equation}
Let us give a formula for the tangent lift that reduces the continuity of the
mapping $X\mapsto\tlift{X}$ to continuity of familiar operations.
\begin{lemma}\label{lem:XTdecomp}
Let\/ $r\in\{\infty,\omega\}$ and let\/ $\man{M}$ be a\/ $\C^r$-manifold with
a\/ $\C^r$-affine connection\/ $\nabla^{\man{M}}$\@.  Then
\begin{equation*}
\tlift{X}(v_x)=\horlift(v_x,X(x))+
\verlift(v_x,\nabla^{\man{M}}_{v_x}X+T^{\man{M}}(X(x),v_x)),
\end{equation*}
where\/ $T^{\man{M}}$ is the torsion of\/ $\nabla^{\man{M}}$\@.
\begin{proof}
Let $v_x\in\tb{\man{M}}$ and let $Y\in\sections[r]{\tb{\man{M}}}$ be such that
$Y(x)=v_x$\@.  Note that
\begin{equation*}
\derivatzero{}{s}\flow{X}{t}\scirc\flow{Y}{s}(x)=\tf[x]{\flow{X}{t}}(Y(x)).
\end{equation*}
Also compute
\begin{align*}
\derivatzero{}{s}\flow{X}{t}\scirc\flow{Y}{s}=&\;
\derivatzero{}{s}\flow{Y}{s}\scirc\flow{X}{t}\scirc\flow{X}{-t}
\flow{Y}{-s}\scirc\flow{X}{t}\scirc\flow{Y}{s}(x)\\
=&\;Y(\flow{X}{t}(x))+\tf[x]{\flow{X}{t}}\left(\derivatzero{}{s}
\flow{X}{-t}\scirc\flow{Y}{-s}\scirc\flow{X}{t}\scirc\flow{Y}{s}
(\flow{X}{t}(x))\right).
\end{align*}
Note that, for $f\in\func[r]{\man{M}}$\@,
\begin{equation*}
f\scirc\flow{X}{-t}\scirc\flow{Y}{-s}\scirc\flow{X}{t}\scirc\flow{Y}{s}(x)=
f(x)+st\lieder{[Y,X]}{f}(x)+o(\snorm{st}),
\end{equation*}
by~\cite[Proposition~4.2.34]{RA/JEM/TSR:88}\@.
Therefore,
\begin{align*}
\derivatzero{}{s}\flow{X}{-t}\scirc\flow{Y}{-s}\scirc\flow{X}{t}
\scirc\flow{Y}{s}(\flow{X}{t}(x))=t[Y,X](\flow{X}{t}(x)).
\end{align*}
Putting the above calculations together gives
\begin{equation*}
\tf[x]{\flow{X}{t}}(Y(x))=Y(\flow{X}{t}(x))-t[X,Y](\flow{X}{t}(x)).
\end{equation*}
Thus, making use of~\eqref{eq:XTflow}\@,
\begin{equation*}
\flow{\hl{X}}{t}\scirc\flow{\tlift{X}}{t}(Y(x))=
\tau^{(t,0)}_{\gamma_-}(Y(\flow{X}{t}(x))-t[X,Y](\flow{X}{t}(x))),
\end{equation*}
where $\gamma_-$ is the integral curve of $-X$ through $\flow{X}{t}(x)$ and
$\tau_{\gamma_-}$ is parallel translation along $\gamma_-$\@.  If $\gamma$ is
the integral curve of $X$ through $x$ note that
$\tau^{(t,0)}_{\gamma_-}=\tau^{(0,t)}_\gamma$\@.  Now we compute
\begin{align*}
\derivatzero{}{t}\flow{-\hl{X}}{t}\scirc\flow{\tlift{X}}{t}(Y(x))=&\;
\derivatzero{}{t}\tau^{(0,t)}_\gamma(Y(\flow{X}{t}(x))-
t[X,Y](\flow{X}{t}(x)))\\
=&\;\nabla_XY(x)-[X,Y](x)=\nabla_YX(x)+T(X(x),Y(x)).
\end{align*}
Note that, since $\tlift{X}$ and $\hl{X}$ are both vector fields over $X$\@,
it follows that
\begin{equation*}
t\mapsto\tau^{(0,t)}_\gamma(Y(\flow{X}{t}(x)))
\end{equation*}
is a curve in $\tb[x]{\man{M}}$\@.  Thus the derivative of this curve at
$t=0$ is in $\vb[Y(x)]{\tb{\man{M}}}$\@.  Thus we have shown that
\begin{equation*}
\derivatzero{}{t}\flow{-\hl{X}}{t}\scirc\flow{\tlift{X}}{t}(v_x)=
\verlift(v_x,\nabla_{v_x}X(x)+T(X(x),v_x)).
\end{equation*}

Finally, for $f\in\func[r]{\man{M}}$\@, by the BCH formula, we have
\begin{equation*}
f\scirc\flow{-\hl{X}}{t}\scirc\flow{\tlift{X}}{t}(v_x)=
f\scirc\flow{\tlift{X}-\hl{X}}{t}+o(\snorm{t}^2).
\end{equation*}
Differentiating with respect to $t$ and evaluating at $t=0$ gives the result.
\end{proof}
\end{lemma}

Now we can combine
Theorems~\pldblref{the:+xcont}{pl:+cont}\@,~\ref{the:Lvcont}\@,
and~\ref{the:hliftX}\@, and Corollary~\ref{cor:nablaXcont} to give the
following result.
\begin{corollary}
If\/ $\man{M}$ is a\/ $\C^\omega$-manifold, then the mapping
\begin{equation*}
\sections[\omega]{\tb{\man{M}}}\ni
X\mapsto\tlift{X}\in\sections[\omega]{\tb{\tb{\man{M}}}}
\end{equation*}
is continuous.
\end{corollary}

\section{Topology for the space of real analytic mappings}\label{sec:mappings}

As a final application of our analysis methods, we consider the problem of
topologising the space of real analytic mappings between real analytic
manifolds $\man{M}$ and $\man{N}$\@.  To do this, we shall connect our
constructions to a commonly used topology for finitely differentiable and
smooth mappings.

\subsection{The weak-PB topology for mappings}

We first consider a well-known topology for finitely differentiable and
smooth mappings.  We let $\man{M}$ and $\man{N}$ be $\C^\infty$-manifolds and
let $\nu\in\integernn\cup\{\infty\}$\@, with
$\mappings[\nu]{\man{M}}{\man{N}}$\@, therefore, the space of
$\C^\nu$-mappings from $\man{M}$ to $\man{N}$\@.  The topology we consider
for this space is called the ``weak topology'' by \cite[\S2.1]{MWH:76} and
the $\mathrm{CO}^\nu$-topology by \citet[\S4.3]{PWM:80}\@.  It is the
adaptation of the standard compact-open topology for spaces of continuous
mappings between topological spaces~\cite[\S43]{SW:70}\@.  Somewhat
precisely, it is the topology where a sequence
$\ifam{\Phi_j}_{j\in\integerp}$ converges to $\Phi$ if and only if:
\begin{compactenum}
\item $\nu\in\integerp$\@: the sequence of all derivatives of
$\ifam{\Phi_j}_{j\in\integerp}$ up to and including order $\nu$ converge
uniformly on compact subsets;
\item $\nu=\infty$\@: for all $k\in\integernn$\@, the sequence of all
derivatives of $\ifam{\Phi_j}_{j\in\integerp}$ up to and including order
$k$ converge uniformly on compact subsets.
\end{compactenum}
There is another characterisation of this topology, and it is the following.
\begin{theorem}
The topology described above is the same as the initial topology associated
with the family of mappings
\begin{equation*}
\mapdef{\Theta_f}{\mappings[\nu]{\man{M}}{\man{N}}}{\func[\nu]{\man{M}}}
{\Phi}{\Phi^*f,}\qquad f\in\func[\infty]{\man{N}},
\end{equation*}
where\/ $\func[\infty]{\man{N}}$ has the topology described above.
\end{theorem}

While we were not able to pinpoint a proof of this theorem anywhere, it will
not come as a surprise to those who understand this topology well; its
validity boils down to being able to find globally defined coordinate
functions about any point in $\man{M}$ or $\man{N}$\@.  Since this ability is
shared by real analytic manifolds, we use this as motivation for the
following definition.
\begin{definition}
Let $\man{M}$ and $\man{N}$ be real analytic manifolds.  The \defn{weak-PB
topology} for $\mappings[\omega]{\man{M}}{\man{N}}$ is the initial topology
associated with the family of mappings
\begin{equation*}
\mapdef{\Theta_f}{\mappings[\omega]{\man{M}}{\man{N}}}
{\func[\omega]{\man{M}}}{\Phi}{\Phi^*f,}\qquad f\in\func[\omega]{\man{N}},
\end{equation*}
where $\func[\omega]{\man{N}}$ has the usual topology.\oprocend
\end{definition}

The following lemma gives a few alternative means of characterising the
weak-PB topology, making use of the fact that a real analytic manifold can be
embedded in a Euclidean space of sufficiently high dimension~\cite[Theorem~3]{HG:58}\@.
\begin{lemma}\label{lem:embed-topology}
Let\/ $\man{M}$ and\/ $\man{N}$ be\/ $\C^\omega$-manifolds and let
\begin{equation*}
\mapdef{\vect{\iota}}{\man{N}}{\real^N}{y}{(\iota^1(y),\dots,\iota^N(y))}
\end{equation*}
be a proper\/ $\C^\omega$-embedding.  Then the following topologies for\/
$\mappings[\omega]{\man{M}}{\man{N}}$ agree:
\begin{compactenum}[(i)]
\item \label{pl:weakPB1} the initial topology associated with the family of
mappings
\begin{equation*}
\mapdef{\Psi_f}{\mappings[\omega]{\man{M}}{\man{N}}}{\func[\omega]{\man{M}}}
{\Phi}{\Phi^*f,}\qquad f\in\func[\omega]{\man{M}};
\end{equation*}
\item \label{pl:weakPB2} the initial topology associated with the family of
mappings
\begin{equation*}
\mapdef{\Psi_{\iota^j}}{\mappings[\omega]{\man{M}}{\man{N}}}
{\func[\omega]{\man{M}}}{\Phi}{\Phi^*\iota^j,}\qquad j\in\{1,\dots,N\};
\end{equation*}
\item \label{pl:weakPB3} the topology induced on\/
$\mappings[\omega]{\man{M}}{\man{N}}\subset
\mappings[\omega]{\man{M}}{\real^N}$ by the weak-PB topology for
\begin{equation*}
\mappings[\omega]{\man{M}}{\real^N}\simeq
\oplus_{j=1}^N\func[\omega]{\man{M}}.
\end{equation*}
\end{compactenum}
\begin{proof}
\eqref{pl:weakPB1}$\subset$\eqref{pl:weakPB2}\@: This follows from the
definition of initial topology.

\eqref{pl:weakPB2}$\subset$\eqref{pl:weakPB3}\@: Note that, because the
induced topology is the initial topology induced by the inclusion
$\mappings[\omega]{\man{M}}{\man{N}}\subset
\mappings[\omega]{\man{M}}{\real^N}$\@, the topology~\eqref{pl:weakPB3} is
the coarsest topology for which that inclusion is continuous.  By
Theorem~\ref{the:compcont}\@, the topology~\eqref{pl:weakPB2} has the
property that this inclusion is continuous.

\eqref{pl:weakPB3}$\subset$\eqref{pl:weakPB2}\@: The
topology~\eqref{pl:weakPB3} has the property that the mappings
$\Psi_{\iota^j}$ are continuous,
whence~\eqref{pl:weakPB3}$\subset$\eqref{pl:weakPB2}\@.

\eqref{pl:weakPB3}$\subset$\eqref{pl:weakPB1}\@: To prove this part of the
lemma, we shall show that $\Psi_f$ is continuous for all
$f\in\func[\omega]{\man{M}}$ for the topology~\eqref{pl:weakPB3}\@.  First of
all, this is obvious when $\man{N}=\real^N$ and we take the embedding
$\vect{\iota}=\id_{\real^N}$\@.

Now we use Lemma~\ref{plem:PB-embedding} from the proof of
Theorem~\ref{the:compcont}\@, taking $\man{M}=\real^N$ and
$\man{S}=\man{N}$\@.  Let $f\in\func[\omega]{\man{M}}$ and, by the lemma, let
$\ol{f}\in\func[\omega]{\real^N}$ be such that
$f=\ol{f}\scirc\vect{\iota}$\@.  Now consider the commutative diagram
\begin{equation*}
\xymatrix{{\mappings[\omega]{\man{M}}{\real^N}}\ar[r]^{\Psi_{\ol{f}}}&
{\func[\omega]{\man{M}}}\\\mappings[\omega]{\man{M}}{\man{N}}\ar[ru]_{\Psi_f}
\ar[u]^{\Phi\mapsto\vect{\iota}\scirc\Phi}&}
\end{equation*}
where $\mappings[\omega]{\man{M}}{\man{N}}$ has the
topology~\eqref{pl:weakPB3} and (as required by the definition of the
topology~\eqref{pl:weakPB3}) $\mappings[\omega]{\man{M}}{\real^N}$ has the
weak-PB topology.  Thus the vertical and horizontal arrows are continuous,
whence the diagonal arrow is continuous, as desired.

Finally, let us prove the assertion made in part~\eqref{pl:weakPB3} that the
weak-PB topology for $\mappings[\omega]{\man{M}}{\real^N}$ is isomorphic to
$\oplus_{j=1}^N\func[\omega]{\man{M}}$\@.  This is equivalent, given what we
have already proven, to the assertion that the initial topology for
$\mappings[\omega]{\man{M}}{\real^N}$ associated with the mappings
\begin{equation*}
\mapdef{\Psi_j}{\mappings[\omega]{\man{M}}{\real^N}}
{\func[\omega]{\man{M}}}{\vect{\Phi}}{\Phi^j,}
\end{equation*}
$j\in\{1,\dots,N\}$\@, is isomorphic to
$\oplus_{j=1}^N\func[\omega]{\man{M}}$\@.  Here we write
$\vect{\Phi}=(\Phi^1,\dots,\Phi^N)$\@.  However, this assertion follows from
a general assertion regarding initial
topologies~\cite[Theorem~8.12]{SW:70}\@.  Indeed, this general result asserts
that the mapping
\begin{equation*}
\mapdef{\Psi}{\mappings[\omega]{\man{M}}{\real^N}}
{\oplus_{j=1}^N\func[\omega]{\man{M}}}{\vect{\Phi}}{(\Phi^1,\dots,\Phi^N)}
\end{equation*}
is an homeomorphism onto its image.  Since it is surjective, it is,
therefore, an homeomorphism.  (Here we also make use of the fact that finite
direct sums are topologically isomorphic to products~\cite[Proposition~4.3.2]{HJ:81}\@.)
\end{proof}
\end{lemma}

Note that this topology can be characterised as being the uniform topology
associated with the family of semimetrics\footnote{It is not uncommon to call
this a ``pseudometric.''  We shall use ``semimetric,'' consistent with our
use of ``seminorm.''}
\begin{equation*}
\mapdef{\d^\omega_{\nbhd{K},\vect{a},f}}{\mappings[\omega]{\man{M}}{\man{N}}
\times\mappings[\omega]{\man{M}}{\man{N}}}{\real}{(\Phi_1,\Phi_2)}
{p^\omega_{\nbhd{K},\vect{a}}(\Phi^*_1f-\Phi^*_2f),}
\end{equation*}
for $\nbhd{K}\subset\man{M}$ compact, $\vect{a}\in\c_0(\integernn;\realp)$\@,
and $f\in\func[\omega]{\man{N}}$\@.  Moreover, the lemma ensures that it
suffices to restrict attention to the semimetrics of this form, and where $f\in\{\iota^1,\dots,\iota^N\}$\@.

Our objective is to come to a useful understanding of these semimetrics and
how to make use of them, in the same way as we have done for the seminorms
for the real analytic topology for sections of a vector bundle.

\subsection{Estimates for semimetrics for real analytic topology}

In this section we obtain a useful upper bound for the semimetrics
\begin{equation*}
\d^\omega_{\nbhd{K},\vect{a},f},\qquad\nbhd{K}\subset\man{M}\ \textrm{compact},\
\vect{a}\in\c_0(\integernn;\realp),\ f\in\func[\omega]{\man{N}},
\end{equation*}
for the weak-PB topology for $\mappings[\omega]{\man{M}}{\man{N}}$\@.  There
are a few difficulties that arise in obtaining such bounds.
\begin{compactenum}
\item First of all, for fixed $f\in\func[\omega]{\man{N}}$\@, the mapping
$\Phi\mapsto f\scirc\Phi$ is nonlinear and so requires a careful analysis.
Thankfully, the recursive analysis of Lemma~\ref{lem:pissynabla} can readily
be adapted to our needs.
\item Second, and not unrelated to the first, the expression
$\Phi^*_1f-\Phi^*_2f$ involves an evaluation of $f$ at points in $\man{N}$
that are possibly distant, and so the matter of evaluating the seminorm
$p^\omega_{\nbhd{K},\vect{a}}$ on this expression should involve comparing
the derivatives of $f$ at such distant points.  One way of doing
this\textemdash{}perhaps the way most faithful to our intrinsic
approach\textemdash{}would be to use the connection we assume on $\man{N}$ to
parallel transport the derivatives from one point to the other.  We elect to
\emph{not} carry this out here.  Instead we perform initial local estimates
which we then apply to the global case by standard compactness arguments.
\item A reader will not be surprised to learn that, in understanding the
weak-PB topology for $\mappings[\omega]{\man{M}}{\man{N}}$\@, we shall make
use of the tensors $\what{A}^m_{\Phi,s}$\@, $m\in\integernn$\@,
$s\in\{0,1,\dots,m\}$\@, defined in Lemma~\ref{lem:pbsymiterateddersI} and
associated with a mapping
$\Phi\in\mappings[\omega]{\man{M}}{\man{N}}$\@.\footnote{We now modify the
notation for these tensors to explicitly include their dependence on the
mapping.}  A difficulty that arises in working with these tensors is that
they are vector bundle mappings
\begin{equation*}
(\what{A}^m_{\Phi,s},\id_{\man{M}})\in
\vbmappings[r]{\Symalg*[s]{\Phi^*\ctb{\man{N}}}}
{\Symalg*[m]{\ctb{\man{M}}}},
\end{equation*}
\ie~vector bundle mappings whose domain depends on the mapping $\Phi$\@.
When working with two such mappings, it becomes necessary to ``canonicalise''
this vector bundle to eliminate its dependence on the mapping.
\end{compactenum}

Let us address the last issue first.  To do so, we make use of a proper real
analytic embedding $\map{\vect{\iota}}{\man{N}}{\real^N}$\@, as in
Lemma~\ref{lem:embed-topology}\@.  We also make use of the following
elementary fact.
\begin{lemma}
Let\/ $r\in\{\infty,\omega\}$\@, let\/ $\man{M}$ and\/ $\man{N}$ be\/
$\C^r$-manifolds, and let\/ $\Phi\in\mappings[r]{\man{M}}{\man{N}}$\@.  Then
the vector bundle\/ $\Phi^*{(\man{N}\times\real^m)}$ is naturally isomorphic
to\/ $\man{M}\times\real^m$\@.
\begin{proof}
One can show that the mapping
\begin{equation*}
\man{M}\times\real^m\ni(x,\vect{v})\mapsto((\Phi(x),\vect{v}),x)
\in\Phi^*(\man{N}\times\real^m)
\end{equation*}
furnishes the desired vector bundle isomorphism.
\end{proof}
\end{lemma}

By the lemma and other devices, we shall be able to reduce ourselves to
consideration of the vector bundle mappings
\begin{equation*}
(\what{A}^m_{\vect{\iota}\scirc\Phi,s},\id_{\man{M}})\in
\vbmappings[r]{\Symalg*[s]{(\real^N_{\man{M}})}}
{\Symalg*[m]{\ctb{\man{M}}}},\qquad m\in\integernn,\ s\in\{0,1,\dots,m\},
\end{equation*}
whose domain and codomain do not depend on~$\Phi$\@.

We begin by giving an upper bound for the semimetric
$\d^\omega_{\nbhd{K},\vect{a},f}$ that depends on the differences between the
tensors for pull-back introduced in Lemma~\ref{lem:pbsymiterateddersI}\@.  We
shall have need of seminorms for the compact-open topology for the space
$\sections[m]{\man{E}}$ of class $\C^m$-sections of a vector bundle:
\begin{equation*}
p^m_{\nbhd{K}}(\xi)=\sup\setdef{\dnorm{j_m\xi(x)}_{\metric_{\man{M},\pi_{\man{E}}}}}
{x\in\nbhd{K}},
\end{equation*}
and semimetrics for the compact-open topology for
$\mappings[0]{\man{M}}{\man{N}}$\@:
\begin{equation*}
\d^0_{\nbhd{K}}(\Phi,\Psi)=\sup\setdef{\d_{\metric_{\man{N}}}(\Phi(x),\Psi(x))}
{x\in\nbhd{K}},
\end{equation*}
for $\nbhd{K}\subset\man{M}$ compact, and where $\metric_{\man{N}}$ is a
Riemannian metric on $\man{N}$ with $\d_{\metric_{\man{N}}}$ the associated
metric.  As we will be performing a local analysis which we then use to give
a global result, we will find it useful to have the following comparison
result for the distance function associated to two Riemannian metrics.  This
is certainly a known result, although we could not locate a proof.
\begin{lemma}\label{lem:metric-ind}
If\/ $\metric_1$ and\/ $\metric_2$ are\/ $\C^\infty$-Riemannian metrics on
a\/ $\C^\infty$-manifold\/ $\man{M}$ with metrics\/ $\d_1$ and\/ $\d_2$\@,
respectively, and if\/ $\nbhd{K}\subset\man{M}$ is compact, then there
exists\/ $C\in\realp$ such that
\begin{equation*}
C^{-1}\d_1(x_1,x_2)\le\d_2(x_1,x_2)\le C\d_1(x_1,x_2)
\end{equation*}
for every\/ $x_1,x_2\in\nbhd{K}$\@.
\begin{proof}
First we give a local version of the result.  Let $x\in\man{M}$\@.  Let
$\nbhd{N}_1$ and $\nbhd{N}_2$ be geodesically convex neighbourhoods of $x$
with respect to the Riemannian metrics $\metric_1$ and $\metric_2$\@,
respectively~\cite[Proposition~IV.3.4]{SK/KN:63a}.  Thus every pair of points
in $\nbhd{N}_1$ can be connected by a unique distance-minimising geodesic for
$\metric_1$ that remains in $\nbhd{N}_1$\@, and similarly with $\nbhd{N}_2$
and $\metric_2$\@.  By Sublemma~\ref{psublem:metric-ind1} from the proof of
Lemma~\ref{lem:G1G2comp}\@, let $C_x\in\realp$ be such that
\begin{equation*}
C^{-2}_x\metric_1(v_x,v_x)<\metric_2(v_x,v_x)<C^2_x\metric_1(v_x,v_x),
\qquad v_x\in\tb[x]{\man{M}}.
\end{equation*}
By continuity of $\metric_1$ and $\metric_2$\@, we can choose $\nbhd{N}_1$
and $\nbhd{N}_2$ sufficiently small that
\begin{equation*}
C^{-2}_x\metric_1(v_y,v_y)<\metric_2(v_y,v_y)<C^2_x\metric_1(v_y,v_y),
\qquad y\in\nbhd{N}_1\cup\nbhd{N}_2.
\end{equation*}
Now define $\nbhd{U}_x=\nbhd{N}_1\cap\nbhd{N}_2$\@.  Then every pair of
points in $\nbhd{U}_x$ can be connected with a unique distance-minimising
geodesic of both $\metric_1$ and $\metric_2$ that remains in
$\nbhd{N}_1\cup\nbhd{N}_2$\@.  Now let $x_1,x_2\in\nbhd{U}_x$\@.  Let
$\map{\gamma}{\interval[0,1]}{\man{M}}$ be the unique distance-minimising
$\metric_1$-geodesic connecting $x_1$ and $x_2$\@.  Then
\begin{align*}
\d_2(x_1,x_2)\le&\;\ell_{\metric_2}(\gamma)
=\int_0^1\sqrt{\metric_2(\gamma'(t),\gamma'(t))}\,\d{t}\\
\le&\;C_x\int_0^1\sqrt{\metric_1(\gamma'(t),\gamma'(t))}\,]d{t}\\
=&\;C_x\ell_{\metric_1}(\gamma)=C_x\d_1(x_1,x_2).
\end{align*}
One similarly shows that $\d_1(x_1,x_2)\le C_x\d_2(x_1,x_2)$\@.

We now prove the assertion of the lemma.  Let $\nbhd{K}\subset\man{M}$ be
compact and, for each $x\in\nbhd{K}$\@, let $\nbhd{U}_x$ be a neighbourhood
of $x$ and let $C_x\in\realp$ be as in the preceding paragraph.  Let
$x_1,\dots,x_k\in\nbhd{K}$ be such that
$\nbhd{K}\subset\cup_{j=1}^k\nbhd{U}_{x_j}$\@.  Let
\begin{equation*}
D_a=\sup\setdef{\d_a(x,y)}{x,y\in\nbhd{K}},\qquad a\in\{1,2\}.
\end{equation*}
By the Lebesgue Number Lemma~\cite[Theorem~1.6.11]{DB/YB/SI:01}\@, let
$r_a\in\realp$ be such that, if $x_1,x_2\in\nbhd{K}$ satisfy $\d_a(x_1,x_2)<r_a$\@,
$a\in\{1,2\}$\@, then there exists $j\in\{1,\dots,k\}$ such that
$x_1,x_2\in\nbhd{U}_{x_j}$\@.  Let us denote
\begin{equation*}
C=\max\Bigl\{C_{x_1},\dots,C_{x_k},\frac{D_1}{r_2},\frac{D_2}{r_1}\Bigr\}.
\end{equation*}

Now let $x_1,x_2\in\nbhd{K}$\@.  If $\d_1(x_1,x_2)<r_1$\@, then let
$j\in\{1,\dots,k\}$ be such that $x_1,x_2\in\nbhd{U}_j$\@.  Then
\begin{equation*}
\d_2(x_1,x_2)\le C\d_1(x_1,x_2).
\end{equation*}
If $\d_1(x_1,x_2)\ge r_1$\@, then
\begin{equation*}
\frac{\d_2(x_1,x_2)r_1}{D_2}\le\frac{\d_2(x_1,x_2)r_1}{\d_2(x_1,x_2)}
\le d_1(x_1,x_2)
\end{equation*}
This gives $\d_2(x_1,x_2)\le C\d_1(x_1,x_2)$\@.  Swapping the r\^oles of
$\metric_1$ and $\metric_2$ gives $\d_1(x_1,x_2)\le C\d_2(x_1,x_2)$\@, giving
the lemma.
\end{proof}
\end{lemma}

With these developments at hand, we have the following result.
\begin{lemma}\label{lem:semimetric-bound}
Let\/ $\man{M}$ and\/ $\man{N}$ be\/ $\C^\omega$-manifolds and let\/
$\map{\vect{\iota}}{\man{N}}{\real^N}$ be a proper real analytic
embedding. Let\/ $\nbhd{K}\subset\man{M}$ be compact, let\/
$\vect{a}\in\c_0(\integernn;\realp)$\@, and let\/
$f\in\func[\omega]{\man{N}}$\@.  Let $\nbhd{V}\subset\man{N}$ be a relatively
compact open set for which\/
$\Phi(\nbhd{K})\cup\Psi(\nbhd{K})\subset\nbhd{V}$\@.  Then there exist\/
$C,\sigma\in\realp$ such that
\begin{multline*}
\dlnorm j_m(\ol{f}\scirc(\vect{\iota}\scirc\Phi)-
\ol{f}\scirc(\vect{\iota}\scirc\Psi))(x)\drnorm_{\metric_{\man{M}}}\\
\le C\sigma^{-m}p^{m+1}_{\closure(\nbhd{V})}(f)\left(
\sum_{s=0}^m\sum_{j=0}^s\frac{j!}{s!}
p^0_{\nbhd{K}}(\what{A}^s_{\vect{\iota}\scirc\Phi,j}-
\what{A}^s_{\vect{\iota}\scirc\Psi,j})+\d^0_{\nbhd{K}}(\Phi,\Psi)
\sum_{s=0}^m\sum_{j=0}^s\frac{j!}{s!}
p^0_{\nbhd{K}}(\what{A}^s_{\vect{\iota}\scirc\Phi,j})\right)
\end{multline*}
for\/ $x\in\nbhd{K}$ and\/ $m\in\integernn$\@.
\begin{proof}
First we prove a local version of the lemma.  We let $\nbhd{U}\subset\real^n$
and $\nbhd{V},\nbhd{W}\subset\real^k$ be open sets, let
$\vect{\Phi}\in\mappings[\omega]{\nbhd{U}}{\nbhd{V}}$ and
$\vect{\Psi}\in\mappings[\omega]{\nbhd{U}}{\nbhd{W}}$\@, and let\/
$f\in\func[\omega]{\nbhd{V}\cup\nbhd{W}}$\@.  Let\/ $\nbhd{K}\subset\nbhd{U}$
be compact and let\/ $\nbhd{L}\subset\nbhd{V}\cup\nbhd{W}$ be compact and
such that\/ $\vect{\Phi}(\nbhd{K}),\vect{\Psi}(\nbhd{K})\subset\nbhd{L}$\@.
We use the standard Euclidean metrics whose norms we simply denote by
$\dnorm{\cdot}$\@.  We also use the standard flat connections for which
covariant differentiation is the usual differentiation.

Let $m\in\integernn$\@.  For $s\in\{0,1,\dots,m\}$\@, we use
Lemma~\ref{lem:pbsymiterateddersI} and calculate
\begin{align*}
\dlnorm\linder[s]{(f\scirc&\vect{\Phi})}(\vect{x})-
\linder[s]{(f\scirc\vect{\Psi})}(\vect{x})\drnorm\\
=&\;\adnorm{\sum_{j=0}^s\left(
\what{\vect{A}}^s_{\vect{\Phi},j}(\vect{x})
(\linder[j]{f}(\vect{\Phi}(\vect{x})))-
\what{\vect{A}}^s_{\vect{\Psi},j}(\vect{x})
(\linder[j]{f}(\vect{\Psi}(\vect{x})))\right)}\\
\le&\;\sum_{j=0}^s\dnorm{\what{\vect{A}}^s_{\vect{\Phi},j}(\vect{x})-
\what{\vect{A}}^s_{\vect{\Psi},j}(\vect{x})}
\dnorm{\linder[j]{f}(\vect{\Psi}(\vect{x}))}\\
&\;+\sum_{j=0}^s\dnorm{\what{\vect{A}}^s_{\vect{\Phi},j}(\vect{x})}
\dnorm{\linder[j]{f}\scirc\vect{\Phi}(\vect{x)}-
\linder[j]{f}\scirc\vect{\Psi}(\vect{x})},
\end{align*}
using Lemma~\ref{lem:A(S)<AS}\@.  Let us give some estimates related to the
last two terms in the preceding expression.

For the first, by Cauchy\textendash{}Schwarz we have
\begin{align*}
\sum_{j=0}^s\frac{j!}{s!}\dnorm{\what{\vect{A}}^s_{\vect{\Phi},j}(\vect{x})&-
\what{\vect{A}}^s_{\vect{\Psi},j}(\vect{x})}\frac{1}{j!}
\dnorm{\linder[j]{f}(\vect{\Psi}(\vect{x}))}\\
\le&\;\left(\sum_{j=0}^s\left(\frac{j!}{s!}
\dnorm{\what{\vect{A}}^s_{\vect{\Phi},j}(\vect{x})-
\what{\vect{A}}^s_{\vect{\Psi},j}(\vect{x})}\right)^2\right)^{1/2}
\left(\sum_{j=0}^s\left(\frac{1}{j!}
\dnorm{\linder[j]{f}(\vect{\Psi}(\vect{x}))}\right)^2\right)^{1/2}\\
\le&\;\left(\sum_{j=0}^s\frac{j!}{s!}
\dnorm{\what{\vect{A}}^s_{\vect{\Phi},j}(\vect{x})-
\what{\vect{A}}^s_{\vect{\Psi},j}(\vect{x})}\right)p^m_{\nbhd{L}}(f)
\end{align*}
for $\vect{x}\in\nbhd{K}$ and $s\in\{0,1,\dots,m\}$\@.  In the last step, we
have used~\eqref{eq:12inftynorms}\@.

For the second, we consider two cases depending on specifying certain types
of open sets $\nbhd{U}$\@, $\nbhd{V}$\@, and $\nbhd{W}$\@.

In the first case, we suppose that $\nbhd{V}=\nbhd{W}$ are open balls and
that $\nbhd{L}$ is a closed ball.  Then, using convexity of $\nbhd{L}$ and
the Mean Value Theorem, we have
\begin{equation*}
\dnorm{\linder[j]{f}(\vect{y}_1)-\linder[j]{f}(\vect{y}_2)}\le
p^0_{\nbhd{L}}(\linder[j+1]{f})\dnorm{\vect{y}_1-\vect{y}_2},
\qquad\vect{y}_1,\vect{y}_2\in\nbhd{L},\ j\in\{0,1,\dots,m\}.
\end{equation*}
We can conclude, therefore, that
\begin{multline*}
\left(\sum_{j=0}^s\left(\frac{1}{j!}
\dnorm{\linder[j]{f}\scirc\vect{\Phi}(\vect{x})-
\linder[j]{f}\scirc\vect{\Psi}(\vect{x})}\right)^2\right)^{1/2}\\
\le p^{m+1}_{\nbhd{L}}(f)\dnorm{\vect{\Phi}(\vect{x})-\vect{\Psi}(\vect{x})},
\qquad\vect{x}\in\nbhd{K},\ s\in\{0,1,\dots,m\}.
\end{multline*}

In the second case, we suppose that $\nbhd{V}$ and $\nbhd{W}$ are open balls
with disjoint closures and let $\nbhd{A}\subset\nbhd{V}$ and
$\nbhd{B}\subset\nbhd{W}$ be closed balls.  We take
$\nbhd{L}=\nbhd{A}\cup\nbhd{B}$\@.  Let
\begin{equation*}
r=\inf\setdef{\dnorm{\vect{v}-\vect{w}}}
{\vect{v}\in\nbhd{V},\ \vect{w}\in\nbhd{W}},
\end{equation*}
noting that $r\in\realp$ since $\nbhd{V}$ and $\nbhd{W}$ are convex with
disjoint closures~\cite[Theorem~11.4]{RTR:70}\@.  Let $\vect{y}_1,\vect{y}_2\in\nbhd{L}$\@.  If
$\vect{y}_1,\vect{y}_2\in\nbhd{A}$ or $\vect{y}_1,\vect{y}_2\in\nbhd{B}$\@,
then, as in the preceding paragraph, we have
\begin{equation*}
\dnorm{\linder[j]{f}(\vect{y}_1)-\linder[j]{f}(\vect{y}_2)}\le
\max\{p^0_{\nbhd{A}}(\linder[j+1]{f}),p^0_{\nbhd{B}}(\linder[j+1]{f})\}
\dnorm{\vect{y}_1-\vect{y}_2},\qquad j\in\{0,1,\dots,m\}.
\end{equation*}
Otherwise we have
\begin{multline*}
\dnorm{\linder[j]{f}(\vect{y}_1)-\linder[j]{f}(\vect{y}_2)}\le
(p^0_{\nbhd{A}}(\linder[j]{f})+p^0_{\nbhd{B}}(\linder[j]{f}))\\
\le\frac{2}{r}p^0_{\nbhd{L}}(\linder[j]{f})\dnorm{\vect{y}_1-\vect{y}_2}
\qquad j\in\{0,1,\dots,m\}.
\end{multline*}
Therefore, in this case,
\begin{multline*}
\left(\sum_{j=0}^s\left(\frac{1}{j!}
\dnorm{\linder[j]{f}\scirc\vect{\Phi}(\vect{x})-
\linder[j]{f}\scirc\vect{\Psi}(\vect{x})}\right)^2\right)^{1/2}\\
\le\frac{2}{r}p^{m+1}_{\nbhd{L}}(f)
\dnorm{\vect{\Phi}(\vect{x})-\vect{\Psi}(\vect{x})},
\qquad\vect{x}\in\nbhd{K}, s\in\{0,1,\dots,m\}.
\end{multline*}

We see that, for both cases in the preceding paragraphs, there exists
$C\in\realp$ such that we have
\begin{multline*}
\left(\sum_{j=0}^sx\left(\frac{1}{j!}
\dnorm{\linder[j]{f}\scirc\vect{\Phi}(\vect{x})-
\linder[j]{f}\scirc\vect{\Psi}(\vect{x})}\right)^2\right)^{1/2}\\
\le Cp^{m+1}_{\nbhd{L}}(f)
\dnorm{\vect{\Phi}(\vect{x})-\vect{\Psi}(\vect{x})},
\qquad\vect{x}\in\nbhd{K},\ m\in\integernn,\ s\in\{0,1,\dots,m\}.
\end{multline*}
Therefore, again for both cases and using Cauchy\textendash{}Schwarz, we have
\begin{align*}
\sum_{j=0}^s\frac{j!}{s!}&\dnorm{\what{\vect{A}}^s_{\vect{\Phi},j}(\vect{x})}
\frac{1}{j!}\dnorm{\linder[j]{f}\scirc\vect{\Phi}(\vect{x)}-
\linder[j]{f}\scirc\vect{\Psi}(\vect{x})}\\
\le&\;\left(\sum_{j=0}^s\left(\frac{j!}{s!}
\dnorm{\what{\vect{A}}^s_{\vect{\Phi},j}(\vect{x})}\right)^2\right)^{1/2}
\left(\sum_{j=0}^s\left(\frac{1}{j!}
\dnorm{\linder[j]{f}\scirc\vect{\Phi}(\vect{x)}-
\linder[j]{f}\scirc\vect{\Psi}(\vect{x})}\right)^2\right)^{1/2}\\
\le&\;C\left(\sum_{j=0}^s\frac{j!}{s!}
\dnorm{\what{\vect{A}}^s_{\vect{\Phi},j}(\vect{x})}\right)p^{m+1}_{\nbhd{L}}(f)
\dnorm{\vect{\Phi}(\vect{x})-\vect{\Psi}(\vect{x})}
\end{align*}
for $\vect{x}\in\nbhd{K}$\@, $m\in\integernn$\@, and $s\in\{0,1,\dots,m\}$\@,
making use of~\eqref{eq:12inftynorms}\@.

Since
\begin{equation*}
\dlnorm j_m(f\scirc\vect{\Phi}-f\scirc\vect{\Psi})(\vect{x})\drnorm\\
\le\sum_{s=0}^m\frac{1}{s!}
\dnorm{\linder[s]{(f\scirc\vect{\Phi})}(\vect{x})-
\linder[s]{(f\scirc\vect{\Psi})}(\vect{x})}
\end{equation*}
by~\eqref{eq:12inftynorms}\@, we can conclude that there exists $C\in\realp$
such that
\begin{multline*}
\dnorm{j_m(f\scirc\vect{\Phi}-f\scirc\vect{\Psi})(\vect{x})}\\\le
Cp^{m+1}_{\nbhd{L}}(f)\left(\sum_{s=0}^m\sum_{j=0}^s\frac{j!}{s!}
\dnorm{\what{\vect{A}}^s_{\vect{\Phi},j}(\vect{x})-
\what{\vect{A}}^s_{\vect{\Psi},j}(\vect{x})}+
\dnorm{\vect{\Phi}(\vect{x})-\vect{\Psi}(\vect{x})}
\sum_{s=0}^m\sum_{j=0}^s\frac{j!}{s!}
\dnorm{\what{\vect{A}}^s_{\vect{\Phi},j}(\vect{x})}\right)
\end{multline*}
for\/ $\vect{x}\in\nbhd{K}$ and $m\in\integernn$\@, this formula being valid
when
\begin{compactenum}
\item $\nbhd{V}=\nbhd{W}$ are open balls and
$\nbhd{L}\subset\nbhd{V}=\nbhd{W}$ is a closed ball or
\item $\nbhd{V}$ and $\nbhd{W}$ are open balls with disjoint closures,
$\nbhd{A}\subset\nbhd{V}$ and $\nbhd{B}\subset\nbhd{W}$ are closed balls, and
$\nbhd{L}=\nbhd{A}\cup\nbhd{B}$\@.
\end{compactenum}

Now we prove the global result using the local result.  For
$\nbhd{K}\subset\man{M}$ compact, let $\nbhd{V}\subset\man{N}$ be open,
precompact, and with the property that
\begin{equation*}
\Phi(\nbhd{K}),\Psi(\nbhd{K})\subset\nbhd{V}.
\end{equation*}
According to Lemma~\ref{plem:PB-embedding} from the proof of
Theorem~\ref{the:compcont}\@, we let $\ol{f}\in\func[\omega]{\real^N}$ be
such that $\vect{\iota}^*\ol{f}=f$\@.

To make use of the local computations, for $x\in\nbhd{K}$ we consider two
cases.

First we consider the case when $\Phi(x)=\Psi(x)$\@.  In this case we have
the following data:
\begin{compactenum}
\item $\chi^1,\dots,\chi^n\in\func[\omega]{\man{M}}$ that comprise a
coordinate chart $\vect{\chi}$ on a neighbourhood $\nbhd{U}_x$ of $x$\@;
\item an open ball $\nbhd{V}_x\subset\real^N$ about
$\vect{\iota}\scirc\Phi(x)=\vect{\iota}\scirc\Psi(x)$\@;
\item a precompact neighbourhood $\nbhd{U}'_x\subset\nbhd{U}_x$ of $x$ and a
precompact neighbourhood $\nbhd{V}'_x\subset\nbhd{V}_x$ of
$\vect{\iota}\scirc\Phi(x)=\vect{\iota}\scirc\Psi(x)$ for which
\begin{compactenum}
\item $\nbhd{V}'_x\subset\real^N$ is an open ball,
\item $\nbhd{V}'_x\cap\man{N}\subset\nbhd{V}$ is contained in a
$\metric_{\man{N}}$-geodesically convex neighbourhood of
$\vect{\iota}\scirc\Phi(x)=\vect{\iota}\scirc\Psi(x)$\@, and
\item $\vect{\iota}\scirc\Phi(\closure(\nbhd{U}'_x))\subset\nbhd{V}'_x$\@.
\end{compactenum}
\end{compactenum}
In case $\Phi(x)\not=\Psi(x)$\@, we have the following data:
\begin{compactenum}
\item $\chi^1,\dots,\chi^n\in\func[\omega]{\man{M}}$ that comprise a
coordinate chart $\vect{\chi}$ on a neighbourhood $\nbhd{U}_x$ of $x$\@;
\item an open ball $\nbhd{V}_x\subset\real^N$ about
$\vect{\iota}\scirc\Phi(x)$\@;
\item an open ball $\nbhd{W}_x\subset\real^N$ about
$\vect{\iota}\scirc\Psi(x)$ for which $\closure(\nbhd{V}_x)\cap\closure(\nbhd{W}_x)=\emptyset$\@;
\item a precompact neighbourhood $\nbhd{U}'_x\subset\nbhd{U}_x$ of $x$\@, a
precompact neighbourhood $\nbhd{V}'_x\subset\nbhd{V}_x$ of
$\vect{\iota}\scirc\Phi(x)$\@, and a precompact neighbourhood
$\nbhd{W}'_x\subset\nbhd{W}_x$ of $\vect{\iota}\scirc\Psi(x)$ for which
\begin{compactenum}
\item $\nbhd{V}'_x$ is an open ball,
\item $\nbhd{W}'_x$ is an open ball,
\item $\nbhd{V}'_x\cap\man{N}\subset\nbhd{V}$ is contained in a
$\metric_{\man{N}}$-geodesically convex neighbourhood of $\vect{\iota}\scirc\Phi(x)$\@, and
\item $\nbhd{W}'_x\cap\man{N}\subset\nbhd{V}$ is contained in a
$\metric_{\man{N}}$-geodesically convex neighbourhood of
$\vect{\iota}\scirc\Psi(x)$\@, and
\item $\vect{\iota}\scirc\Phi(\closure(\nbhd{U}'_x))\subset\nbhd{V}'_x$\@,
and
\item $\vect{\iota}\scirc\Psi(\closure(\nbhd{U}'_x))\subset\nbhd{W}'_x$\@.
\end{compactenum}
\end{compactenum}
To save having to consider multiple versions of the same estimates, if we are
in the first of the two cases above, let us denote $\nbhd{W}_x=\nbhd{V}_x$
and $\nbhd{W}'_x=\nbhd{V}'_x$\@.

By Lemma~\ref{lem:metric-ind}\@, let $C_{1,x}\in\realp$ be such that
\begin{equation*}
C^{-1}_{1,x}\dnorm{\vect{\iota}(y)-\vect{\iota}(z)}\le
\d_{\metric_{\man{N}}}(y,z)\le
C_{1,x}\dnorm{\vect{\iota}(y)-\vect{\iota}(z)},
\qquad y\in\closure(\nbhd{V}'_x\cap\man{N}),\
z\in\closure(\nbhd{W}'_x\cap\man{N}).
\end{equation*}
By Lemma~\ref{lem:G1G2comp}\@, let $C_{2,x},\sigma_{2,x}\in\realp$ be such
that
\begin{align*}
&\frac{\sigma^{m+s}_{2,x}}{C_{2,x}}\dnorm{\vect{\alpha}^m_s(\vect{\chi}(x'))}
\le\dnorm{\alpha^m_s(x')}_{\metric_{\man{M}},\metric_{\real^N}}\le
\frac{C_{2,x}}{\sigma^{m+s}_{2,x}}\dnorm{\vect{\alpha}^m_s(\vect{\chi}(x'))},\\
&\frac{\sigma^{m+s}_{2,x}}{C_{2,x}}\dnorm{\vect{\beta}^m_s(\vect{\chi}(x'))}
\le\dnorm{\beta^m_s(x')}_{\metric_{\man{M}},\metric_{\real^N}}\le
\frac{C_{2,x}}{\sigma^{m+s}_{2,x}}\dnorm{\vect{\beta}^m_s(\vect{\chi}(x'))},
\end{align*}
for
\begin{equation*}
(\alpha^m_s,\id_{\man{M}}),(\beta^m_s,\id_{\man{M}})\in
\vbmappings[r]{\tensor*[s]{(\real^N_{\man{M}})}}
{\tensor*[m]{\ctb{\man{M}}}},
\end{equation*}
for
\begin{equation*}
x'\in\closure(\nbhd{U}'_x),\ m\in\integernn,\ s\in\{0,1,\dots,m\},
\end{equation*}
and where $\vect{\alpha}^m_s$ and $\vect{\beta}^m_s$ are the local
representatives for $\alpha^m_s$ and $\beta^m_s$\@.  By
Lemma~\ref{lem:normcompare}\@, let $C_{3,x},\sigma_{3,x}\in\realp$ be such
that
\begin{multline*}
\frac{\sigma_{3,x}^m}{C_{3,x}}
\dnorm{j_m(\ol{f}\scirc(\vect{\iota}\scirc\vect{\Phi})-
\ol{f}\scirc(\vect{\iota}\scirc\vect{\Psi}))(\vect{\chi}(x'))}\le
\dnorm{j_m(\ol{f}\scirc(\vect{\iota}\scirc\Phi)-
\ol{f}\scirc(\vect{\iota}\scirc\Psi))(x')}_{\metric_{\man{M}}}\\
\le\frac{C_{3,x}}{\sigma^m_{3,x}}
\dnorm{j_m(\ol{f}\scirc(\vect{\iota}\scirc\vect{\Phi})-
\ol{f}\scirc(\vect{\iota}\scirc\vect{\Psi}))(\vect{\chi}(x'))},\qquad
m\in\integernn,\ x'\in\closure(\nbhd{U}'_x),
\end{multline*}
where $\vect{\Phi}$ and $\vect{\Psi}$ are the ``semi'' local representatives
of $\Phi$ and $\Psi$\@:
\begin{equation*}
\Phi=\vect{\Phi}\scirc\vect{\chi},\quad\Psi=\vect{\Psi}\scirc\vect{\chi}.
\end{equation*}
Finally, by our local computations above, there exists $C_{4,x}\in\realp$
such that
\begin{multline*}
\dnorm{j_m(\ol{f}\scirc(\vect{\iota}\scirc\vect{\Phi})-
\ol{f}\scirc(\vect{\iota}\scirc\vect{\Psi}))(\vect{\chi}(x'))}\\
\le C_{4,x}p^{m+1}_{\closure((\nbhd{V}'_x\cup\nbhd{W}'_x)\cap\man{N})}(f)
\left(\sum_{s=0}^m\sum_{j=0}^s\frac{j!}{s!}
\dnorm{\what{\vect{A}}^s_{\vect{\iota}\scirc\vect{\Phi},j}(\vect{\chi}(x'))-
\what{\vect{A}}^s_{\vect{\iota}\scirc\vect{\Psi},j}(\vect{\chi}(x'))}\right.\\
\left.+\dnorm{\vect{\iota}\scirc\Phi(x')-\vect{\iota}\scirc\Psi(x')}
\sum_{s=0}^m\sum_{j=0}^s\frac{j!}{s!}
\dnorm{\what{\vect{A}}^s_{\vect{\iota}\scirc\vect{\Phi},j}(\vect{\chi}(x'))}\right)
\end{multline*}
for $x'\in\closure(\nbhd{U}'_x)$ and $m\in\integernn$\@.  In the above
expressions, we can without loss of generality assume that
\begin{equation*}
\sigma_{2,x},\sigma_{3,x}\le1.
\end{equation*}
A useful consequence of this is that
$\sigma_{2,x}^{-(m+s)}\le\sigma_{2,x}^{-2m}$ if $s\in\{0,1,\dots,m\}$\@.

Now we assemble all of the above together into a single estimate.  To do so,
let $x_1,\dots,x_r\in\nbhd{K}$ be such that
$\nbhd{K}\subset\cup_{l=1}^r\nbhd{U}'_{x_l}$\@.  Let $x\in\nbhd{K}$ and let
$l\in\{1,\dots,r\}$ be such that $x\in\nbhd{U}'_{x_l}$\@.  Then, using the
formulae above, we directly compute
\begin{align*}
\dlnorm j_m(&\ol{f}\scirc(\vect{\iota}\scirc\Phi)-
\ol{f}\scirc(\vect{\iota}\scirc\Psi))(x)\drnorm_{\metric_{\man{M}}}\\
\le&\;C_{3,x_l}C_{4,x_l}\sigma^{-m}_{3,x_l}
p^{m+1}_{\closure(\nbhd{V})}(f)\left(
\sum_{s=0}^m\sum_{j=0}^s\frac{j!}{s!}C_{2,x_l}\sigma^{-(m+s)}_{2,x_l}
\dnorm{\what{A}^s_{\vect{\iota}\scirc\Phi,j}(x)-
\what{A}^s_{\vect{\iota}\scirc\Psi,j}(x)}_{\metric_{\man{M}},\metric_{\man{N}}}
\right.\\
&\left.+C_{1,x_l}\d_{\metric_{\man{N}}}(\Phi(x),\Psi(x))
\sum_{s=0}^m\sum_{j=0}^s\frac{j!}{s!}C_{2,x_l}\sigma^{-(m+s)}_{2,x_l}
\dnorm{\what{A}^s_{\vect{\iota}\scirc\Phi,j}}_{\metric_{\man{M}},\metric_{\man{N}}}
\right)\\
\le&\;C\sigma^{-m}p^{m+1}_{\closure(\nbhd{V})}(f)\left(
\sum_{s=0}^m\sum_{j=0}^s\frac{j!}{s!}
p^0_{\nbhd{K}}(\what{A}^s_{\vect{\iota}\scirc\Phi,j}-
\what{A}^s_{\vect{\iota}\scirc\Psi,j})+\d^0_{\nbhd{K}}(\Phi,\Psi)
\sum_{s=0}^m\sum_{j=0}^s\frac{j!}{s!}
p^0_{\nbhd{K}}(\what{A}^s_{\vect{\iota}\scirc\Phi,j})\right),
\end{align*}
provided we take
\begin{equation*}
C=\max\setdef{\max\{C_{2,x_l}C_{3,x_l}C_{4,x_l},
C_{1,x_l}C_{2,x_l}C_{3,x_l}C_{4,x_l}\}}{l\in\{1,\dots,r\}}
\end{equation*}
and
\begin{equation*}\eqqed
\sigma=\min\setdef{\sigma_{2,x_l}^2\sigma_{3,x_l}}{l\in\{1,\dots,r\}}.
\end{equation*}
\end{proof}
\end{lemma}

\begin{remark}
In the smooth case, one can apply the same methods to arrive at a bound where
the term ``$C\sigma^{-m}$'' can be replaced by a constant depending on $m$\@,
and where the sums of the form
$\sum_{s=0}^m\sum_{j=0}^s\frac{j!}{s!}\alpha^s_j$ can be replaced with
\begin{equation*}\eqoprocend
(m+1)^2\max\setdef{\alpha^s_j}{s\in\{0,1,\dots,m\},\ j\in\{0,1,\dots,s\}}
\end{equation*}
\end{remark}

From the preceding lemma, we see that there are two components critically
involved in estimating $\d^\omega_{\nbhd{K},\vect{a},f}(\Phi,\Psi)$\@;
namely, to make this semimetric expression small, one must appropriately
bound each of
\begin{equation*}
\sum_{s=0}^m\sum_{j=0}^sp^0_{\nbhd{K}}(\widehat{A}^s_{\vect{\iota}\scirc\Phi,j}-
\widehat{A}^s_{\vect{\iota}\scirc\Psi,j})\enspace\textrm{and}\enspace
\d^0_{\nbhd{K}}(\Phi,\Psi).
\end{equation*}
The second of these is easy to understand.  For the first, we shall relate
the differences between these tensors and the tensors $A_\Phi$ and $A_\Psi$
defined in Lemma~\ref{lem:APhidef}\@, or, more precisely, the modified
tensors $B_\Phi$ and $B_\Psi$ introduced in the statement of
Lemma~\ref{lem:pullbackform}\@.  Since we will be taking differences, of
course we work with $\vect{\iota}\scirc\Phi$ and $\vect{\iota}\scirc\Psi$ in
place of $\Phi$ and $\Psi$\@, where $\map{\vect{\iota}}{\man{N}}{\real^N}$ is
a real analytic proper embedding.
\begin{lemma}\label{lem:Amscont}
Let\/ $\man{M}$ be a real analytic manifold and let\/
$\vect{\Phi},\vect{\Psi}\in\mappings[\omega]{\man{M}}{\real^N}$\@.  Then,
for\/ $\nbhd{K}\subset\man{M}$ compact, there exists\/
$C,\sigma,\rho\in\realp$ such that
\begin{equation*}
\dnorm{D^r_{\nabla^{\man{M}}}(A^m_{\vect{\Phi},s}-
A^m_{\vect{\Psi},s})(x)}\le C\sigma^{-m}\rho^{-(m+r-s)}(m+r-s)!
p^{m+r-s-1}_{\nbhd{K}}(B_{\vect{\Phi}}-B_{\vect{\Psi}}),
\end{equation*}
for\/ $m,r\in\integernn$\@,\/ $s\in\{0,1,\dots,m\}$\@, and\/
$x\in\nbhd{K}$\@.
\begin{proof}
The reader will not be surprised to learn that a proof can be made modelled
closely upon that of Lemma~\ref{lem:pissynabla}\@, applied to the particular
case of Lemma~\ref{lem:pbiterateddersI}\@.  This is precisely what we shall
do.  In particular, we shall make use of the notation and observations
introduced at the beginning of the proof of Lemma~\ref{lem:pissynabla}\@,
specifically, we introduce $\beta,\alpha,\gamma,C_{m,s}\in\realp$ as in that
proof.  We shall also, sometimes without explicit mention, use some of the
computations from the proof of Lemma~\ref{lem:pissynabla}\@.  Thus a reader
feeling as if they are missing something may wish to refer back to that proof
as an attempt to locate the absent steps.

First we obtain the desired recursion relations for the tensors
\begin{equation*}
A^m_{\vect{\Phi},s}-A^m_{\vect{\Psi},s},\qquad m\in\integernn,\
s\in\{0,1,\dots,m\},
\end{equation*}
recalling from Lemma~\ref{lem:pullbackform} the tensors
\begin{equation*}
B_{\vect{\Phi}},B_{\vect{\Psi}}\in
\sections[\omega]{\tensor*[2]{\ctb{\man{M}}}\otimes\real^N_{\man{M}}}.
\end{equation*}
We do this simply by subtracting the recursion relations for
$A^m_{\vect{\Phi},s}$ and $A^m_{\vect{\Psi},s}$ from
Lemma~\ref{lem:pbiterateddersI}\@, and using the equality
\begin{multline*}
A^m_{\vect{\Phi},s}\otimes\id_{\ctb{\man{M}}}(\Ins_j(\beta_s,B_{\vect{\Phi}}))-
A^m_{\vect{\Psi},s}\otimes\id_{\ctb{\man{M}}}(\Ins_j(\beta_s,B_{\vect{\Psi}}))\\
=(A^m_{\vect{\Phi},s}-A^m_{\vect{\Psi},s})\otimes
\id_{\ctb{\man{M}}}(\Ins_j(\beta_s,B_{\vect{\Phi}}))+
A^m_{\vect{\Psi},s}\otimes\id_{\ctb{\man{M}}}
(\Ins_j(\beta_s,B_{\vect{\Phi}}-B_{\vect{\Psi}}))
\end{multline*}
to obtain
\begin{align*}
A^{m+1}_{\vect{\Phi},m+1}-A^{m+1}_{\vect{\Psi},m+1}=&\;0,\\
A^{m+1}_{\vect{\Phi},s}-A^{m+1}_{\vect{\Psi},s}=&\;
\nabla^{\man{M}}(A^m_{\vect{\Phi},s}-A^m_{\vect{\Psi},s})+
\Lambda^s_m\scirc(A^m_{\vect{\Phi},s-1}-A^m_{\vect{\Psi},s-1})\\
&+\sum_{j=1}^s\Theta^s_{jm}(A^m_{\vect{\Phi},s}-A^m_{\vect{\Psi},s})
+\sum_{j=1}^s\Omega^s_{jm}(A^m_{\vect{\Psi},s}),\quad s\in\{1,\dots,m\},\\
A^{m+1}_{\vect{\Phi},0}-A^{m+1}_{\vect{\Psi},0}=&\;
\nabla^{\man{M}}(A^m_{\vect{\Phi},0}-A^m_{\vect{\Psi},0}),
\end{align*}
for $m\in\integernn$\@, where
\begin{align*}
\Lambda^s_m(\alpha^m_{s-1})=&\;\alpha^m_{s-1}\otimes\id_{\ctb{\man{M}}},\\
\Theta^s_{jm}(\alpha^m_s)(\beta_s)=&\;-\alpha^m_s\otimes\id_{\ctb{\man{M}}}
(\Ins_j(\beta_s,B_{\vect{\Phi}})),\\
\Omega^s_{jm}(\alpha^m_s)(\beta_a)=&\;-\alpha^m_s\otimes\id_{\ctb{\man{M}}}
(\Ins_j(\beta_s,B_{\vect{\Phi}}-B_{\vect{\Psi}})),
\end{align*}
\cf~the constructions from page~\pageref{page:pissy-translate}\@.  Let us
provide some bounds for some of the terms in the recursion, all coming from
Lemma~\ref{lem:particular-bounds}\@, and with a contribution from
Lemma~\ref{lem:pissynabla} in the final estimate:
\begin{align*}
&\dnorm{D^r_{\nabla^{\man{M}}}\Lambda^s_m\scirc D^a_{\nabla^{\man{M}}}
(A^m_{\vect{\Phi},s-1}-A^m_{\vect{\Psi},s-1})(x)}_{\metric_{\man{M}}}
\le\dnorm{D^a_{\nabla^{\man{M}}}
(A^m_{\vect{\Phi},s-1}-A^m_{\vect{\Psi},s-1})(x)}_{\metric_{\man{M}}},\\
&\dnorm{D^r_{\nabla^{\man{M}}}\Theta^s_{jm}\scirc D^a_{\nabla^{\man{M}}}
(A^m_{\vect{\Phi},s}-A^m_{\vect{\Psi},s})(x)}_{\metric_{\man{M}}}
\le C_1\sigma_1^{-r}r!\dnorm{D^a_{\nabla^{\man{M}}}
(A^m_{\vect{\Phi},s-1}-A^m_{\vect{\Psi},s-1})(x)}_{\metric_{\man{M}}},\\
&\dnorm{D^r_{\nabla^{\man{M}}}\Omega^s_{jm}\scirc D^a_{\nabla^{\man{M}}}
A^m_{\vect{\Psi},s}(x)}_{\metric_{\man{M}}}\le
C_1(C_1\sigma_1^{-1}\gamma)^mC_{m,s}\left(\frac{\beta}{\sigma_1}\right)^{m+a-s}
(m+a-s)!\\
&\phantom{\dnorm{D^r_{\nabla^{\man{M}}}\Omega^s_{jm}\scirc D^a_{\nabla^{\man{M}}}
A^m_{\vect{\Psi},s}(x)}_{\metric_{\man{M}}}\le}\times
\dnorm{D^r_{\nabla^{\man{M}}}(B_{\vect{\Phi}}-
B_{\vect{\Psi}})(x)}_{\metric_{\man{M}}},
\end{align*}
for some suitable $C_1,\sigma_1,\gamma\in\realp$\@, and we can assume,
without loss of generality, that $C_1,\gamma\ge1$ and $\sigma_1\le1$\@.  For
the last of these inequalities, we are using the detailed
bound~\eqref{eq:DrAmsest} for
$\dnorm{D^a_{\nabla^{\man{M}}}A^m_{\vect{\Psi},s}(x)}_{\metric_{\man{M}}}$
that was proved during the course of the proof of
Lemma~\ref{lem:pissynabla}\@.

We note that the recursion relations immediately imply that
\begin{equation*}
A^m_{\vect{\Phi},0}-A^m_{\vect{\Psi},0}=0,\
A^m_{\vect{\Phi},m}-A^m_{\vect{\Psi},m}=0,\qquad m\in\integernn.
\end{equation*}
We shall show that, for $m\ge2$ and $s\in\{1,\dots,m-1\}$\@, we have
\begin{multline*}
\dnorm{D^r_{\nabla^{\man{M}}}(A^m_{\vect{\Phi},s}-
A^m_{\vect{\Psi},s})}_{\metric_{\man{M}}}\\
\le C_1(C_1\sigma_1^{-1}\gamma)^mC_{m,s}
\left(\frac{\beta}{\sigma_1}\right)^{m+r-s}
(m+r-s)!p^{m+r-s-1}_{\nbhd{K}}(B_{\vect{\Phi}}-B_{\vect{\Psi}}).
\end{multline*}

Thus, to prove the result, we prove the desired estimate by a double
induction, first on $m$ and, for fixed $m$\@, by induction on $s$\@.

By Lemma~\ref{lem:leibniz} we have
\begin{multline*}
D^r_{\nabla^{\man{M}}}(A^{m+1}_{\vect{\Phi},s}-A^{m+1}_{\vect{\Psi},s})=
\underbrace{D^r_{\nabla^{\man{M}}}\nabla^{\man{M}}(A^m_{\vect{\Phi},s}-
A^m_{\vect{\Psi},s})}_{\textrm{term~1}}\\
+\underbrace{\sum_{j=1}^s\sum_{a=0}^r\binom{r}{a}D^a_{\nabla^{\man{M}}}
\Theta^s_{jm}(D^{r-a}_{\nabla^{\man{M}}}
(A^m_{\vect{\Phi},s}-A^m_{\vect{\Psi},s}))}_{\textrm{term~2}}
+\underbrace{\sum_{j=1}^s\sum_{a=0}^r\binom{r}{a}
D^a_{\nabla^{\man{M}}}\Omega^s_{jm}
(D^{r-a}_{\nabla^{\man{M}}}A^m_{\vect{\Psi},s})}_{\textrm{term~3}}\\
+\underbrace{\sum_{a=0}^r\binom{r}{a}D^a_{\nabla^{\man{M}}}\Lambda^s_m
(D^{r-a}_{\nabla^{\man{M}}}
(A^m_{\vect{\Phi},s-1}-A^m_{\vect{\Psi},s-1}))}_{\textrm{term~4}}.
\end{multline*}
We shall separately estimate the four terms in the preceding equation.

As we showed in the proof of Lemma~\ref{lem:Jk+mdecomp}\@, we have
\begin{equation*}
D^r_{\nabla^{\man{M}}}
(\nabla^{\man{M}}(A^m_{\vect{\Phi},s}-A^m_{\vect{\Psi},s}))=
D^{r+1}_{\nabla^{\man{M}}}(A^m_{\vect{\Phi},s}-A^m_{\vect{\Psi},s}).
\end{equation*}
Therefore, by Lemma~\ref{lem:A(S)<AS}\@, we have
\begin{align*}
\dnorm{\textrm{term~1(x)}}_{\metric_{\man{M}}}\le&\;
\dnorm{D^{r+1}_{\nabla^{\man{M}}}(A^m_{\vect{\Phi},s}-
A^m_{\vect{\Psi},s})(x)}_{\metric_{\man{M}}}\\
\le&\;C_1(C_1\sigma_1^{-1}\gamma)^mC_{m,s}
\left(\frac{\beta}{\sigma_1}\right)^{m+r+1-s}(m+r+1-s)!\\
&\;\times p^{m+r-s}_{\nbhd{K}}(B_{\vect{\Phi}}-B_{\vect{\Psi}}).
\end{align*}
In like manner, using the bounds above for the tensors involved,
\begin{align*}
\dnorm{\textrm{term~2(x)}}_{\metric_{\man{M}}}\le&\;
\sum_{j=1}^s\sum_{a=0}^r\frac{r!}{a!(r-a)!}C_1\sigma_1^{-a}a!
\dnorm{D^{r-a}_{\nabla^{\man{M}}}(A^m_{\vect{\Phi},s}-
A^m_{\vect{\Psi},s})(x)}_{\metric_{\man{M}}}\\
\le&\;\sum_{a=0}^r\frac{r!}{a!(r-a)!}
mC_1\sigma_1^{-a}a!(C_1(C_1\sigma_1^{-1}\gamma)^mC_{m,s}\\
&\;\times\left(\frac{\beta}{\sigma_1}\right)^{m+r-a-s}(m+r-a-s)!
p^{m+r-a-s-1}_{\nbhd{K}}(B_{\vect{\Phi}}-B_{\vect{\Psi}})\\
\le&\;C_1(C_1\sigma_1^{-1}\gamma)^mmC_{m,s}
\left(\frac{\beta}{\sigma_1}\right)^{m+r-s}(m+r-s)!\\
&\;\times\left(\sum_{a=0}^r\beta^{-a}\right) p^{m+r-s-1}_{\nbhd{K}}(B_{\vect{\Phi}}-B_{\vect{\Psi}})\\
\le&\;C_1(C_1\sigma_1^{-1}\gamma)C_{m+1,s}\alpha
\left(\frac{\beta}{\sigma_1}\right)^{m+r-s}
p^{m+r-s-1}_{\nbhd{K}}(B_{\vect{\Phi}}-B_{\vect{\Psi}}).
\end{align*}
For the third term, we begin by using Cauchy\textendash{}Schwarz, and then
proceed as above to get
\begin{align*}
\dlnorm\textrm{term~3(x)}&\drnorm_{\metric_{\man{M}}}\le\sum_{a=0}^r
C_1(C_1\sigma_1^{-1}\gamma)^mmC_{m,s}
\left(\frac{\beta}{\sigma_1}\right)^{m+r-a-s}(m+r-a-s)!\\
&\;\phantom{\drnorm_{\metric_{\man{M}}}\le}
\times\dnorm{D^a_{\nabla^{\man{M}}}(B_{\vect{\Phi}}-
B_{\vect{\Psi}})(x)}_{\metric_{\man{M}}}\\
\le&\;\left(\sum_{a=0}^r\left(C_1(C_1\sigma_1^{-1}\gamma)^mmC_{m,s}
\left(\frac{\beta}{\sigma_1}\right)^{m+r-a-s}\frac{(m+r-a-s)!}{a!}
\right)^2\right)^{1/2}\\
&\;\times\left(\sum_{a=0}^r\left(\frac{1}{a!}
\dnorm{D^a_{\nabla^{\man{M}}}(B_{\vect{\Phi}}-
B_{\vect{\Psi}})(x)}_{\metric_{\man{M}}}\right)^2\right)^{1/2}\\
\le&\;\sum_{a=0}^rC_1(C_1\sigma_1^{-1}\gamma)^mmC_{m,s}
\left(\frac{\beta}{\sigma_1}\right)^{m+r-a-s}(m+r-a-s)!
p^r_{\nbhd{K}}(B_{\vect{\Phi}}-B_{\vect{\Psi}})\\
\le&\;C_1(C_1\sigma_1^{-1}\gamma)^mC_{m+1,s}\alpha
\left(\frac{\beta}{\sigma_1}\right)^{m+r-s}(m+1+r-s)!
p^r_{\nbhd{K}}(B_{\vect{\Phi}}-B_{\vect{\Psi}}).
\end{align*}
Finally, we have
\begin{align*}
\dlnorm\textrm{term~4(x)}&\drnorm_{\metric_{\man{M}}}\le
\sum_{a=0}^r\frac{r!}{a!(r-a)!}
\dnorm{D^{r-a}_{\nabla^{\man{M}}}(A^m_{\vect{\Phi},s-1}-
A^m_{\vect{\Psi},s-1})(x)}_{\metric_{\man{M}}}\\
\le&\;\sum_{a=0}^r\frac{r!}{a!(r-a)!}C_1(C_1\sigma_1^{-1}\gamma)^m
C_{m,s-1}\left(\frac{\beta}{\sigma_1}\right)^{m+r-a-s+1}(m+r-a-s+1)!\\
&\;\times p^{m+r-a-s}_{\nbhd{K}}(B_{\vect{\Phi}}-B_{\vect{\Psi}})\\
\le&\;C_1(C_1\sigma_1^{-1}\gamma)^mC_{m,s-1}\alpha
\left(\frac{\beta}{\sigma_1}\right)^{m+r-s+1}(m+r-s+1)!\\
&\;\times p^{m+r-s}_{\nbhd{K}}(B_{\vect{\Phi}}-B_{\vect{\Psi}}).
\end{align*}

Assembling the four components of the estimate and taking
\begin{equation*}
C=C_1,\ \sigma^{-1}=2C_1\sigma_1^{-1}\gamma,\ \rho=\frac{\beta}{\sigma_1}
\end{equation*}
gives the desired result.
\end{proof}
\end{lemma}

\begin{remark}
In the smooth case, as is usual, a similar but simpler analysis can be
carried out, resulting in a bound where the complicated exponential and
factorial terms on the right in the conclusion of the lemma can be replaced
by a single constant, one that will depend on $m$\@, $s$\@, and
$r$\@.\oprocend
\end{remark}

The final ingredient in arriving at a complete understanding of the first
component of the estimate in Lemma~\ref{lem:semimetric-bound} is to
understand how the tensor $B_\Phi$ depends on $\Phi$\@, or, more properly,
how the tensor $B_{\vect{\iota}\scirc\Phi}$ depends on $\Phi$\@.  First of
all, let us give the explicit form of $A_{\vect{\Phi}}$ for
$\vect{\Phi}\in\mappings[\omega]{\man{M}}{\real^N}$\@.  In doing so, of
course we use the standard Euclidean metric for the vector bundle
$\real^N_{\man{M}}$ along with its attendant flat connection for which
covariant differentiation is ordinary differentiation.  Upon doing so, we see
that
\begin{multline*}
A_{\vect{\Phi}}(X(x),Y(x))\\
=(x,(\natpair{\d{\Phi^1}(x)}{\nabla^{\man{M}}_XY(x)}-
X(\natpair{\d{\Phi^1}}{Y})(x),
\dots,\natpair{\d{\Phi^N}(x)}{\nabla^{\man{M}}_XY(x)}-
X(\natpair{\d{\Phi^N}}{Y})(x))),
\end{multline*}
for $X,Y\in\sections[\omega]{\tb{\man{M}}}$\@.

With this formula at hand, we can now assert the nature of the dependence of
$B_{\vect{\iota}\scirc\Phi}$ on $\Phi$\@.  We shall show, indeed, that the
mapping $\Phi\mapsto B_{\vect{\iota}\scirc\Phi}$ is continuous and bounded.
This requires that we have an appropriate notion of boundedness for subsets
of $\mappings[\omega]{\man{M}}{\man{N}}$\@.  The notion of boundedness we
give here is not a completely standard one, but is adapted from a notion of
boundedness for uniform spaces that can be found
in~\cite[\eg][Theorem~1.12]{JH:59}\@.
\begin{definition}\label{def:CwMN-bdd}
For real analytic manifolds $\man{M}$ and $\man{N}$\@, a subset
$\nbhd{B}\subset\mappings[\omega]{\man{M}}{\man{N}}$ is \defn{bounded} if
$\d^\omega_{\nbhd{K},\vect{a},f}|\nbhd{B}\times\nbhd{B}$ is bounded for each
compact $\nbhd{K}\subset\man{M}$\@, $\vect{a}\in\c_0(\integernn;\realp)$\@,
and $f\in\func[\omega]{\man{N}}$\@.\oprocend
\end{definition}

We now have the following result.
\begin{lemma}\label{lem:BPhicont}
If\/ $\man{M}$ and\/ $\man{N}$ are real analytic manifolds and if\/
$\map{\vect{\iota}}{\man{N}}{\real^N}$ is a proper real analytic embedding,
then the mapping
\begin{equation*}
\mappings[\omega]{\man{M}}{\man{N}}\ni\Phi\mapsto B_{\vect{\iota}\scirc\Phi}
\in\sections[\omega]{\tensor*[2]{\ctb{\man{M}}}\otimes\real^N_{\man{M}}}
\end{equation*}
is continuous and bounded.
\begin{proof}
We claim, first, that the mapping
\begin{equation*}
\mappings[\omega]{\man{M}}{\real^N}\ni\vect{\Phi}\mapsto B_{\vect{\Phi}}
\in\sections[\omega]{\tensor*[2]{\ctb{\man{M}}}\otimes\real^N_{\man{M}}}
\end{equation*}
is continuous.  To see this, note that we have an homeomorphism
\begin{equation*}
\mappings[\omega]{\man{M}}{\real^N}\simeq\sections[\omega]{\real^N_{\man{M}}},
\end{equation*}
\ie~because of the vector space structure of the codomain, the space of
mappings is topologically equivalent to the space of sections of the trivial
bundle.  Now we note that $\vect{\Phi}\mapsto A_{\vect{\Phi}}$ is continuous
since $A_{\vect{\Phi}}$ is obtained from $\vect{\Phi}$ by differentiations
and algebraic operations; here we are making use of
Theorems~\ref{the:+xcont} and~\ref{the:jetcont}\@.  Since $B_{\vect{\Phi}}$
is obtained from $A_{\vect{\Phi}}$ by algebraic operations, the claim is
indeed correct.

Now consider the commutative diagram
\begin{equation*}
\xymatrix{{\mappings[\omega]{\man{M}}{\man{N}}}
\ar[r]^(0.37){\Phi\mapsto B_{\vect{\iota}\scirc\Phi}}
\ar[d]_{\Phi\mapsto\vect{\iota}\scirc\Phi}&
{\sections[\omega]{\tensor*[2]{\ctb{\man{M}}}\otimes\real^N_{\man{M}}}}\\
{\mappings[\omega]{\man{M}}{\real^N}}\ar[ru]_{\vect{\Phi}\mapsto
B_{\vect{\Phi}}}&}
\end{equation*}
The vertical arrow is continuous by Lemma~\ref{lem:embed-topology} and we
have just shown that the diagonal arrow is continuous.  This proves the lemma
as concerns continuity.

As concerns boundedness, let
$\nbhd{B}\subset\mappings[\omega]{\man{M}}{\man{N}}$ be bounded.  Let
$\nbhd{K}\subset\man{M}$ be compact and let
$\vect{a}\in\c_0(\integernn;\realp)$\@.  Boundedness of $\nbhd{B}$ implies
that $\d^\omega_{\nbhd{K},\vect{a},\iota^j}|\nbhd{B}\times\nbhd{B}$ is
bounded for $j\in\{1,\dots,N\}$\@.  Now let $\Phi_0\in\nbhd{B}$ be chosen
arbitrarily and note that
\begin{equation*}
p^\omega_{\nbhd{K},\vect{a}}(\vect{\iota}\scirc\Phi)\le
p^\omega_{\nbhd{K},\vect{a}}(\vect{\iota}\scirc\Phi_0)+
p^\omega_{\nbhd{K},\vect{a}}(\vect{\iota}\scirc\Phi-\vect{\iota}\scirc\Phi_0).
\end{equation*}
As the first term on the right in the preceding expression is constant as a
function of $\Phi$\@, we will show that the second term on the right is a
bounded function of $\Phi\in\nbhd{B}$\@.  To see this, we let
$\vect{e}_j\in\real^N$\@, $j\in\{1,\dots,N\}$\@, be the standard basis for
$\real^N$ and calculate
\begin{align*}
p^\omega_{\nbhd{K},\vect{a}}(\vect{\iota}\scirc\Phi-\vect{\iota}\scirc\Phi_0)
=&\;\sum_{j=1}^Np^\omega_{\nbhd{K},\vect{a}}((\iota^j\scirc\Phi)\vect{e}_j-
(\iota^j\scirc\Phi_0)\vect{e}_j)\\
\le&\;\sum_{j=1}^N\d^\omega_{\nbhd{K},\vect{a},\iota^j}(\Phi,\Phi_0),
\end{align*}
and we see that this shows that
\begin{equation*}
\nbhd{B}'\eqdef\setdef{\vect{\iota}\scirc\Phi}{\Phi\subset\nbhd{B}}
\end{equation*}
is bounded in
$\sections[\omega]{\real^N_{\man{M}}}\simeq
\mappings[\omega]{\man{M}}{\real^N}$\@.  Since $\vect{\Phi}\mapsto
B_{\vect{\Phi}}$ is a continuous linear map,
\begin{equation*}
\nbhd{B}''=\setdef{B_{\vect{\iota}\scirc\Phi}}{\Phi\in\nbhd{B}}
\end{equation*}
is bounded, which completes the proof of the lemma.
\end{proof}
\end{lemma}

\subsection{The real analytic superposition operator}\label{subsec:nemytskii}

In this section we complete the story begun with Theorem~\ref{the:compcont}
of composition of real analytic mappings.  For real analytic manifolds
$\man{M}$ and $\man{N}$\@, there are three sorts of composition operators we
can consider:
\begin{gather*}
\mapdef{C_\Phi}{\func[\omega]{\man{N}}}{\func[\omega]{\man{M}}}
{f}{\Phi^*f,}\qquad
\mapdef{S_f}{\mappings[\omega]{\man{M}}{\man{N}}}{\func[\omega]{\man{M}}}
{\Phi}{f\scirc\Phi,}\\
\mapdef{C_{\man{M},\man{N}}}{\mappings[\omega]{\man{M}}{\man{N}}\times
\func[\omega]{\man{N}}}{\func[\omega]{\man{M}}}{(\Phi,f)}{f\scirc\Phi,}
\end{gather*}
the first being defined for fixed
$\Phi\in\mappings[\omega]{\man{M}}{\man{N}}$ and the second for fixed
$f\in\func[\omega]{\man{N}}$\@.  We call these the \defn{composition
operator} associated with $\Phi$\@, the \defn{superposition operator}
associated with $f$\@, and the \defn{joint composition operator}\@,
respectively.  The superposition operator is also known as the ``nonlinear
composition operator'' or the ``Nemytskii operator.''  In general, one
studies these mappings for classes of function spaces,~\eg~Lebesgue spaces or
Hardy spaces.  The questions one can ask for such operators include the
following.
\begin{compactenum}
\item \emph{Well-definedness}\@: Here, for instance, one wishes to know for
which $f$'s or $\Phi$'s do the operators $S_f$ or $C_\Phi$ maps one function
space into another.
\item \emph{Continuity}\@: The continuity of the linear composition operator
$C_\Phi$ is often fairly easily established, and also often coincides with
the well-definedness of the operator.  The continuity of the nonlinear
superposition operator $S_f$\@, however, is often quite difficult to
establish.  Moreover, there are important cases where continuity of this
operator does not coincide with its well-definedness, a well-known example
of this being in the Lipschitz class~\cite{PD:75}\@.
\item \emph{Boundedness}\@: Of course, in the linear case, continuity of
$C_\Phi$ implies boundedness, although not necessarily the converse if the
function spaces in question are not metrisable.  For example, in the real
analytic case in which we are interested here, these spaces are not
metrisable.  In the nonlinear case, the boundedness and continuity of $S_f$
are not generally logically comparable.
\item \emph{Real analyticity}\@: It is sometimes the case that the
superposition operator, though nonlinear, admits a convergent power series
expansion, in which case it is said to be ``real analytic.''  We mention this
here mostly because this is \emph{not} what we are considering here; here we
are considering operators on spaces of real analytic mappings, not operators
which are themselves real analytic.
\end{compactenum}
A reader interested in a detailed discussion of superposition operators for
various classes of function spaces is referred to~\cite{JA/PPZ:90}\@.

For our purposes, we wish to prove continuity and boundedness of the joint
composition operator.  The space of real analytic mappings has the bornology
of Definition~\ref{def:CwMN-bdd}\@.  Of course, $\func[\omega]{\man{M}}$ and
$\func[\omega]{\man{N}}$ have their natural bornologies as topological vector
spaces, and this notion of bounded sets corresponds to continuous seminorms
being bounded.  We shall also use the product bornology on products, that
having a base consisting of products of bounded sets.

We can now state the desired result.
\begin{theorem}
If\/ $\man{M}$ and\/ $\man{N}$ are real analytic manifolds, then the joint
composition operator\/ $C_{\man{M},\man{N}}$ is continuous and bounded.
\begin{proof}
We let $\map{\vect{\iota}}{\man{N}}{\real^N}$ be a proper real analytic
embedding.

We first prove continuity of the joint composition operator.  Let
$\Phi_0\in\mappings[\omega]{\man{M}}{\man{N}}$\@, let
$f_0\in\func[\omega]{\man{N}}$\@, let $\nbhd{K}\subset\man{M}$ be compact,
let $\vect{a}\in\c_0(\integernn;\realp)$\@, and let $\epsilon\in\realp$\@.
Note that
\begin{equation}\label{eq:jointcomp1}
p^\omega_{\nbhd{K},\vect{a}}(f\scirc\Phi-f_0\scirc\Phi_0)=
p^\omega_{\nbhd{K},\vect{a}}(f_0\scirc\Phi-f_0\scirc\Phi_0)+
p^\omega_{\nbhd{K},\vect{a}}((f-f_0)\scirc\Phi).
\end{equation}
We shall separately estimate the two terms on the right.

First, by Lemmata~\ref{lem:pissynabla}\@,~\ref{lem:semimetric-bound}\@,
and~\ref{lem:Amscont}\@, let $C,\sigma\in\realp$ be such that
\begin{compactenum}
\item for $m\in\integernn$ and $s\in\{0,1,\dots,m\}$\@, we
have \begin{equation*}
p^0_{\nbhd{K}}(\widehat{A}^m_{\vect{\iota}\scirc\Phi_0,s})\le
C\sigma^{-m}(m-s)!,
\end{equation*}
\item for $m\in\integernn$ and $x\in\nbhd{K}$\@, we have
\begin{multline*}
\dnorm{j_m(f\scirc\Phi-f\scirc\Psi_0)(x)}_{\metric_{\man{M}}}\le\\
C\sigma^{-m}p^{m+1}_{\closure(\nbhd{V})}(f)\left(
\sum_{s=0}^m\sum_{j=0}^s\frac{j!}{s!}
p^0_{\nbhd{K}}(\what{A}^s_{\vect{\iota}\scirc\Phi,j}-
\what{A}^s_{\vect{\iota}\scirc\Psi,j})+\d^0_{\nbhd{K}}(\Phi,\Psi)
\sum_{s=0}^m\sum_{j=0}^s\frac{j!}{s!}
p^0_{\nbhd{K}}(\what{A}^s_{\vect{\iota}\scirc\Phi,j})\right),
\end{multline*}
and
\item for $m\in\integernn$ and $s\in\{0,1,\dots,m\}$\@, we have
\begin{equation*}
p^0_{\nbhd{K}}(\widehat{A}^m_{\vect{\iota}\scirc\Phi,s}-
\widehat{A}^m_{\vect{\iota}\scirc\Phi_0,s})\le
C\sigma^{-m}(m-s)!p^{m-s}(B_{\vect{\iota}\scirc\Phi}-B_{\vect{\iota}\scirc\Phi_0}).
\end{equation*}
\end{compactenum}
We then have
\begin{align*}
\sum_{s=0}^m\sum_{j=0}^s\frac{j!}{s!}
p^0_{\nbhd{K}}(\what{A}^s_{\vect{\iota}\scirc\Phi,j}-
\what{A}^s_{\vect{\iota}\scirc\Phi_0,j})\le&\;
C\sum_{s=0}^m\sum_{j=0}^s\frac{j!}{s!}\sigma^{-s}(s-j)!
p^{s-j-1}_{\nbhd{K}}(B_{\vect{\iota}\scirc\Phi}-B_{\vect{\iota}\scirc\Phi_0})\\
\le&\;Cp^{m-1}_{\nbhd{K}}(B_{\vect{\iota}\scirc\Phi}-B_{\vect{\iota}\scirc\Phi_0})
\sigma^{-m}\sum_{s=0}^m\sum_{j=0}^s\frac{j!(s-j)!}{s!}\\
\le&\;Cp^{m-1}_{\nbhd{K}}(B_{\vect{\iota}\scirc\Phi}-B_{\vect{\iota}\scirc\Phi_0})
\sigma^{-m}\underbrace{\sum_{s=0}^m2^s}_{2^{m+1}-1}\\
\le&\;4C(2\sigma^{-1})^m
p^{m-1}_{\nbhd{K}}(B_{\vect{\iota}\scirc\Phi}-B_{\vect{\iota}\scirc\Phi_0})
\end{align*}
and
\begin{equation*}
\sum_{s=0}^m\sum_{j=0}^s\frac{j!}{s!}
p^0_{\nbhd{K}}(\what{A}^s_{\vect{\iota}\scirc\Phi_0,j})\le
\sum_{s=0}^m\sum_{j=0}^s\frac{j!}{s!}C\sigma^{-s}(s-j)!
\le4C(2\sigma^{-1})^m
\end{equation*}
for every $\Phi,\Psi\in\mappings[\omega]{\man{M}}{\man{N}}$ and
$m\in\integernn$\@.  Combining these we have
\begin{align*}
a_0a_1\cdots a_m\dlnorm&j_m(f_0\scirc\Phi-f_0\scirc\Phi_0)(x)
\drnorm_{\metric_{\man{M}}}\\
\le&\;a_0a_1\cdots a_mC\sigma^{-m}p^{m+1}_{\closure(\nbhd{V})}(f_0)\\
&\;\times\left(4C(2\sigma^{-1})^m
p^{m-1}_{\nbhd{K}}(B_{\vect{\iota}\scirc\Phi}-B_{\vect{\iota}\scirc\Phi_0})+
4C(2\sigma^{-1})^m\d^0_{\nbhd{K}}(\Phi,\Phi_0)\right)\\
\le&\;b_{1,0}b_{1,1}\cdots b_{1,m+1}p^{m+1}_{\closure(\nbhd{V})}(f_0)
b_{2,0}b_{2,1}\cdots b_{2,m+1}p^{m+1}_{\nbhd{K}}
(B_{\vect{\iota}\scirc\Phi}-B_{\vect{\iota}\scirc\Phi_0})\\
&\;+b_{3,0}b_{3,1}\cdots b_{3,m+1}p^{m+1}_{\closure(\nbhd{V})}(f_0)
\d^0_{\nbhd{K}}(\Phi,\Phi_0),
\end{align*}
where
\begin{align*}
&b_{1,0}=C,\ b_{1,j}=\frac{\sqrt{a_{j-1}}}{\sigma},\qquad
j\in\{1,\dots,m+1\},\\
&b_{2,0}=4C,\ b_{2,j}=\frac{2\sqrt{a_{j-1}}}{\sigma},\qquad
j\in\{1,\dots,m+1\},\\
&b_{3,0}=4C^2,\ b_{3,j}=\frac{2a_{j-1}}{\sigma^2},\qquad
j\in\{1,\dots,m+1\}.
\end{align*}
Therefore,
\begin{equation}\label{eq:jointcomp2}
p^\omega_{\nbhd{K},\vect{a}}(f_0\scirc\Phi-f_0\scirc\Phi_0)\le
p^\omega_{\closure(\nbhd{V}),\vect{b}_1}(f_0)
p^\omega_{\nbhd{K},\vect{b}_2}(B_{\vect{\iota}\scirc\Phi}-
B_{\vect{\iota}\scirc\Phi_0})+p^\omega_{\closure(\nbhd{V}),\vect{b}_3}(f_0)
\d^0_{\nbhd{K}}(\Phi,\Phi_0),
\end{equation}
giving a useful bound for the first term on the right in~\eqref{eq:jointcomp1}\@.

For the second, we first compute
\begin{align*}
p^0_{\nbhd{K}}(\what{A}^m_{\vect{\iota}\scirc\Phi,s})\le&\;
p^0_{\nbhd{K}}(\what{A}^m_{\vect{\iota}\scirc\Phi,s}-
\what{A}^m_{\vect{\iota}\scirc\Phi_0,s})+
p^0_{\nbhd{K}}(\what{A}^m_{\vect{\iota}\scirc\Phi_0,s})\\
\le&\;C\sigma^{-m}(m-s)!p^{m-s-1}_{\nbhd{K}}(B_{\vect{\iota}\scirc\Phi}-
B_{\vect{\iota}\scirc\Phi_0})+C\sigma^{-m}(m-s)!
\end{align*}
and
\begin{align*}
a_0a_1\cdots a_m\dlnorm&j_m(f\scirc\Phi-f_0\scirc\Phi)(x)
\drnorm_{\metric_{\man{M}}}\\
=&\;a_0a_1\cdots a_m\left\dlnorm\bigoplus_{s=0}^m
D^s_{\nabla_{\man{M}}}(\ol{f}\scirc\vect{\iota}\scirc\Phi)(x)-
\bigoplus_{s=0}^mD^s_{\nabla_{\man{M}}}(\ol{f}_0\scirc\vect{\iota}\scirc\Phi)(x)
\right\drnorm_{\metric_{\man{M}}}\\
=&\;a_0a_1\cdots a_m\left\dlnorm\bigoplus_{s=0}^m\sum_{j=0}^s
\widehat{A}^s_{\vect{\iota}\scirc\Phi,j}(x)(\linder[j]{\ol{f}}
(\vect{\iota}\scirc\Phi(x))-
\linder[j]{f_0}(\vect{\iota}\scirc\Phi(x)))\right\drnorm_{\metric_{\man{M}}}\\
\le&\;a_0a_1\cdots a_m\left(\sum_{s=0}^m\sum_{j=0}^s\frac{j!}{s!}
p^0_{\nbhd{K}}(\widehat{A}^s_{\vect{\iota}\scirc\Phi,j})\right)
p^m_{\closure(\nbhd{V})}(f-f_0)\\
\le&b_{4,0}b_{4,1}\cdots b_{4,m}p^m_{\closure(\nbhd{V})}(f-f_0)
b_{5,0}b_{5,1}\cdots b_{6,m}p^m_{\nbhd{K}}(B_{\vect{\iota}\scirc\Phi}-
B_{\vect{\iota}\scirc\Phi_0})\\
&\;+b_{6,0}b_{6,1}\cdots b_{6,m}
p^m_{\closure(\nbhd{V})}(f-f_0),
\end{align*}
where
\begin{align*}
&b_{4,j}=\sqrt{a_j},\qquad j\in\{0,1,\dots,m\},\\
&b_{5,0}=4C\sqrt{a_0},\ b_{6,j}=\frac{2\sqrt{a_j}}{\sigma},\qquad
j\in\{1,\dots,m\},\\
&b_{6,0}=4Ca_0,\ b_{6,j}=\frac{2a_j}{\sigma},\qquad j\in\{1,\dots,m\}.
\end{align*}
Thus
\begin{equation}\label{eq:jointcomp3}
p^\omega_{\nbhd{K},\vect{a}}((f-f_0)\scirc\Phi)\le
p^\omega_{\closure(\nbhd{V}),\vect{b}_4}(f-f_0)
p^\omega_{\nbhd{K},\vect{b}_5}(B_{\vect{\iota}\scirc\Phi}-B_{\vect{\iota}\scirc\Phi_0})
+p^\omega_{\closure(\nbhd{V}),\vect{b}_6}(f-f_0),
\end{equation}
giving a bound for the second term on the right in~\eqref{eq:jointcomp1}\@.

Now we combine the bounds~\eqref{eq:jointcomp2} and~\eqref{eq:jointcomp3}\@.
We take
\begin{equation*}
b_j=\max\{b_{1,j},\dots,b_{6,j}\},\qquad j\in\integernn,
\end{equation*}
and let $\nbhd{O}\subset\mappings[\omega]{\man{M}}{\man{N}}$ and
$\nbhd{N}\subset\func[\omega]{\man{N}}$ be neighbourhoods of $\Phi_0$ and
$f_0$ such that
\begin{align*}
&\d^0_{\nbhd{K}}(\Phi,\Phi_0)<\frac{\epsilon}
{4p^\omega_{\closure(\nbhd{V}),\vect{b}_1}(f_0)},\\
&p^\omega_{\nbhd{K},\vect{b}}(B_{\vect{\iota}\scirc\Phi}-B_{\vect{\iota}\scirc\Phi_0})
<\frac{\epsilon}{4p^\omega_{\closure(\nbhd{V}),\vect{b}_1}(f_0)},\\
&p^\omega_{\closure(\nbhd{V}),\vect{b}}(f-f_0)<\frac{\sqrt{\epsilon}}{2},\\
&p^\omega_{\nbhd{K},\vect{b}}(B_{\vect{\iota}\scirc\Phi}-B_{\vect{\iota}\scirc\Phi_0})
<\frac{\sqrt{\epsilon}}{2},\\
&p^\omega_{\closure(\nbhd{V}),\vect{b}}(f-f_0)<\frac{\epsilon}{4},
\end{align*}
for $\Phi\in\nbhd{O}$ and $f\in\nbhd{N}$\@; the special case where
$p^\omega_{\closure(\nbhd{V}),\vect{b}_1}(f_0)=0$ is accommodated by mere
omission of the first two conditions.  We note that, in the second and fifth
of these conditions, we have made use of Lemma~\ref{lem:BPhicont}\@.  We then
immediately have
\begin{equation*}
p^\omega_{\nbhd{K},\vect{a}}(f\scirc\Phi-f_0\scirc\Phi_0)<\epsilon,\qquad
(\Phi,f)\in\nbhd{O}\times\nbhd{N},
\end{equation*}
which gives the result as concerns continuity.

Now let $\nbhd{B}\subset\mappings[\omega]{\man{M}}{\man{N}}$ and
$\nbhd{C}\subset\func[\omega]{\man{N}}$ be bounded.  Let
$\nbhd{K}\subset\man{M}$ and let $\vect{a}\in\c_0(\integernn;\realp)$\@.  Fix
$\Phi_0\in\nbhd{B}$\@.  Let $(\Phi,f)\in\nbhd{B}\times\nbhd{C}$ and note that
\begin{equation}\label{eq:jointcomp4}
p^\omega_{\nbhd{K},\vect{a}}(f\scirc\Phi)\le
p^\omega_{\nbhd{K},\vect{a}}(f\scirc\Phi-f\scirc\Phi_0)+
p^\omega_{\nbhd{K},\vect{a}}(f\scirc\Phi_0).
\end{equation}
By our analysis for continuity above,
\begin{equation*}
p^\omega_{\nbhd{K},\vect{a}}(f\scirc\Phi-f\scirc\Phi_0)\le
p^\omega_{\closure(\nbhd{V}),\vect{b}_1}(f)
p^\omega_{\nbhd{K},\vect{b}_2}(B_{\vect{\iota}\scirc\Phi}-
B_{\vect{\iota}\scirc\Phi_0})+p^\omega_{\closure(\nbhd{V}),\vect{b}_3}(f)
\d^0_{\nbhd{K}}(\Phi,\Phi_0)
\end{equation*}
for suitable $\vect{b}_j\in\c_0(\integernn;\realp)$\@, $j\in\{1,2,3\}$\@.
Because the mapping $\Phi\mapsto B_{\vect{\iota}\scirc\Phi}$ is bounded (by
Lemma~\ref{lem:BPhicont}) and because the mapping
$\Phi\mapsto\d^0_{\nbhd{K}}(\Phi,\Phi_0)$ is bounded on $\nbhd{B}$ (because
the real analytic uniformity is finer than the compact-open uniformity), we
conclude that the first term of the right in~\eqref{eq:jointcomp4} defines a
bounded function on $\nbhd{B}\times\nbhd{C}$\@.  Since the second term on the
right in~\eqref{eq:jointcomp4} is obviously bounded as a function on
$\nbhd{C}$\@, the boundedness part of the assertion follows.
\end{proof}
\end{theorem}


\printbibliography[heading=bibintoc]

\end{document}